\documentclass[draft,singlespacing]{dcthesis}

\usepackage{amsmath,amssymb,amsthm,amscd}
\usepackage[final]{graphicx}
\usepackage{index}

\allowdisplaybreaks

\newtheorem{prop}{Proposition}[chapter]
\newtheorem{theorem}[prop]{Theorem}
\newtheorem{lemma}[prop]{Lemma}
\newtheorem{corr}[prop]{Corollary}
\theoremstyle{definition}
\newtheorem{definition}[prop]{Definition}
\theoremstyle{remark}
\newtheorem{remark}[prop]{Remark}
\newtheorem{example}[prop]{Example}
\newtheorem*{claim}{Claim}

\makeindex
\newindex{not}{ndx}{nnd}{Notation and Symbol Index}

\DeclareMathOperator{\Res}{Res}
\DeclareMathOperator{\Prim}{Prim}

\DeclareMathOperator{\Aut}{Aut}
\DeclareMathOperator{\lt}{lt}
\DeclareMathOperator{\supp}{supp}
\DeclareMathOperator{\ran}{ran}
\DeclareMathOperator{\Ad}{Ad}
\DeclareMathOperator{\Ind}{Ind}
\DeclareMathOperator{\Ex}{Ex}
\DeclareMathOperator{\id}{id}

\DeclareMathOperator{\ev}{ev}
\DeclareMathOperator{\spn}{span}
\DeclareMathOperator{\isom}{Iso}
\DeclareMathOperator{\ess}{ess}

\newcommand{\mcal}[1]{\mathcal{#1}}
\newcommand{\mfrk}[1]{\mathfrak{#1}}
\newcommand{\lset}[1]{{}_{#1}}
\newcommand{\inv}{^{-1}}
\newcommand{\llangle}{\langle\!\langle}
\newcommand{\rrangle}{\rangle\!\rangle}
\newcommand{\cspn}{\overline{\spn}}
\newcommand{\sidehat}{^{\wedge}}
\newcommand{\doubledual}[1]{\,\widehat{\!\widehat{#1}}}
\newcommand{\unit}{^{(0)}}
\newcommand{\neghalf}{^{-\frac{1}{2}}}
\newcommand{\poshalf}{^{\frac{1}{2}}}
\newcommand{\N}{\mathbb{N}}
\newcommand{\Z}{\mathbb{Z}}
\newcommand{\R}{\mathbb{R}}
\newcommand{\C}{\mathbb{C}}
\newcommand{\Q}{\mathbb{Q}}
\newcommand{\T}{\mathbb{T}}
\newcommand{\erune}{{\raisebox{-.04in}
{\includegraphics[width=.11in]{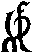}}}}
\newcommand{\suberune}
{_{\text{\includegraphics[width=.07in]{humpf.png}}}}

\title{Groupoid Crossed Products}
\author{Geoff Goehle}
\date{May 29, 2009}
\field{Mathematics}
\degree{Doctor of Philosophy}
\committee{Dana Williams}{Paul Muhly}{John Trout}{David Webb}

\begin{document}

\frontmatter

\maketitle

\chapter*{Abstract}
\addcontentsline{toc}{section}{Abstract}
We present a number of findings concerning groupoid dynamical systems
and groupoid crossed products.  The primary result is an
identification of the spectrum of the groupoid crossed product when
the groupoid has continuously varying abelian stabilizers and a well behaved
orbit space.  In this case, 
the spectrum of the crossed product is homeomorphic, via
an induction map, to a quotient of the spectrum of the crossed product
by the stabilizer group bundle.  The main theorem is also
generalized in the 
groupoid algebra case to an identification of the primitive ideal
space.  This generalization replaces the assumption that the orbit
space is well behaved with an amenability hypothesis.  We then use induction
to show that the primitive ideal space of the groupoid algebra is
homeomorphic to a quotient of the dual of the stabilizer group bundle.
In both cases the identification is topological.  We then apply these
theorems in a number of examples, and examine when a groupoid algebra
has Hausdorff spectrum.  As a separate result, we also develop a
theory of principal groupoid group bundles and 
locally unitary groupoid actions.  We prove that such 
actions are characterized, up to exterior equivalence, by a cohomology
class which arises from a principal bundle.  Furthermore, we also
demonstrate how to construct a locally unitary action from a given
principal bundle.  This last result uses a duality theorem for abelian
group bundles which is also included as part of this thesis.   

\chapter*{Preface}
\addcontentsline{toc}{section}{Preface}
This thesis contains a number of results concerning groupoid crossed
products and groupoid $C^*$-algebras which have been developed through
the author's graduate studies at Dartmouth College.  Groupoid crossed
products inherit the generality of group\-oids.  In particular, they
simultaneously generalize group crossed products, 
transformation group algebras, and groupoid algebras.  
\begin{figure}[ht]
\centering
\includegraphics[height=2in]{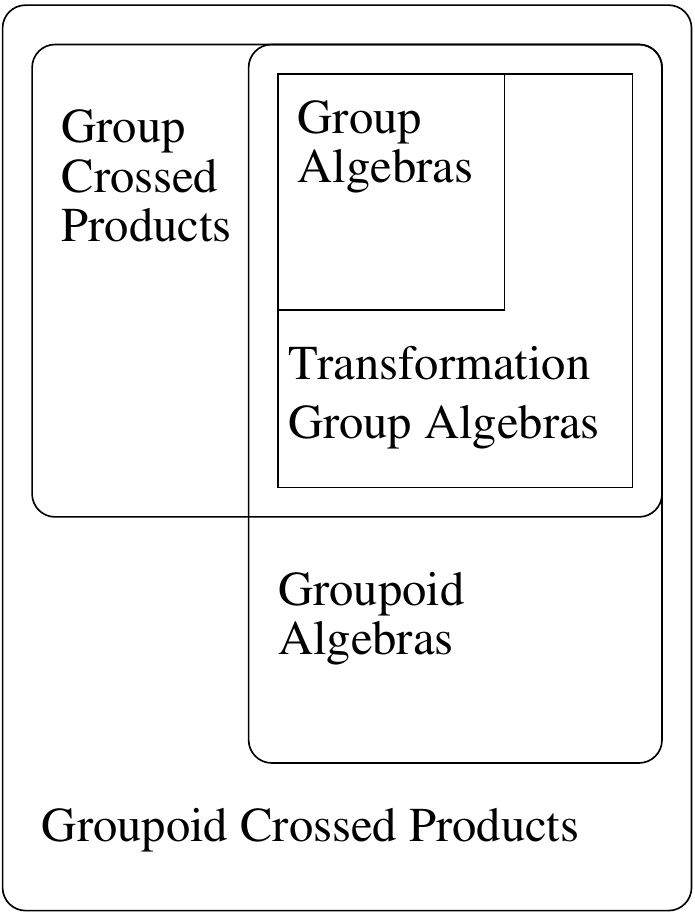}
\caption{Groupoid crossed products are very general objects.}
\label{fig:algdiag}
\end{figure}

Before we describe the structure of the thesis we should give
some idea of what the prerequisites are for understanding its
contents.  The reader who is familiar with groupoids and groupoid
crossed products will encounter no difficulties.  Because the field is
so new, the author has taken some trouble to keep the presentation as
self-contained as possible.  Someone with a basic knowledge of
$C^*$-algebras and functional analysis could expect to understand much
of this thesis, particularly the first two chapters, but may eventually
run into trouble.  A reasonable set of required reading is the
following list of references:
\begin{itemize}
\item {\em An Invitation to $C^*$-algebras}, William Arveson, Chapter
  1, \cite{invitation}
\item {\em Morita Equivalence and Continuous-Trace $C^*$-algebras}, 
  Iain Raeburn and Dana P. Williams, Chapters 1, 2, 3, Sections
  4.1, 4.2, and Appendix A, \cite{tfb}
\item {\em Crossed Products of $C^*$-algebras}, Dana P. Williams,
  Chapters 1, 2 and Appendix~C. \cite{tfb2}
\end{itemize}
These references come highly recommended by the author and are each
worth reading.  In particular, 
we will be citing these books frequently.   Some other works that we
will cite frequently are:
\begin{itemize}
\item {\em A Groupoid Approach to $C^*$-algebras}, Jean Renault,
  \cite{groupoidapproach}
\item {\em Coordinates In Operator Algebras}, Paul Muhly, \cite{coords}
\item {\em Continuous-Trace Groupoid $C^*$-algebras III}, Paul Muhly,
  Jean Renault, and Dana Williams. \cite{ctgIII}.
\item {\em Renault's Equivalence Theorem for Groupoid Crossed
    Products}, Paul Muhly and Dana P. Williams, \cite{renaultequiv}
\item {\em The Ideal Structure of Groupoid Crossed Product
    $C^*$-algebras}, Jean Renault, \cite{renaultgcp}
\end{itemize}
Those readers interested in chasing down citations will find that much
of the groupoid theory in this work is inspired by the first three
references listed above 
and much of the crossed product theory is inspired by the
last two.  What's more, readers are encouraged to look up references.
There has been some effort made to cite results as they appear in
their original context, or if that is not possible, to include a
remark which explains how to extract the given statement from the
statement in the reference.  

As for the structure of this thesis, 
because groupoid crossed products rely heavily on groupoid theory, and
because groupoid theory is a relatively new field in and of itself, we
begin with an introduction to groupoid basics in Chapter
\ref{cha:groupoids}.  This includes
definitions and elementary properties of groupoids in Section
\ref{sec:basics} and actions of
groupoids on topological spaces in Section \ref{sec:actions}.  
Also included in this chapter are
more advanced results concerning groupoid equivalence in Section
\ref{sec:equivalence} and groupoid
amenability in Section \ref{sec:amenable}.  
Next, in Chapter \ref{cha:bundles}
we explore the structure
of groupoid group bundles.  In Section \ref{sec:principal}
we develop the notion of a
principal $S$-bundle and show that they are characterized, up to
isomorphism, by an associated cohomology class.  In Section
\ref{sec:duality} we demonstrate a generalization of Pontryagin
duality for abelian, continuously varying group bundles.  In Section
\ref{sec:opencounter} we describe a counterexample which, in addition
to being interesting in its own right, shows that the work done in Section
\ref{sec:duality} is necessary. In Chapter
\ref{cha:crossed} we introduce the basics of groupoid
dynamical systems and groupoid crossed products.  We start by giving a
brief overview of upper-semicontinuous bundle theory in Section
\ref{sec:cstarbundles} and then define a groupoid dynamical system in
Section \ref{sec:dynamical}.  In Section \ref{sec:covariant} we
develop the theory of covariant representations of groupoid dynamical
systems.  We then use these covariant representations in Section
\ref{sec:crossedprod} to define the groupoid crossed product.  In this
section we also introduce Renault's Disintegration Theorem, which will
be an important tool.  
In Chapter \ref{cha:special-cases} we describe a number of special
cases of groupoid crossed products, and show that
groupoid crossed products generalize groupoid algebras and group
crossed products as in Figure \ref{fig:algdiag}.  Not only
does this connect the theory to existing mathematics, but these
constructions will prove essential in later chapters.  A modest
result, that is nonetheless interesting, is a generalization of the
Stone-von Neumann theorem to groupoids, presented in Section
\ref{sec:neumann}.  Next, in
Chapter \ref{cha:basic} we present some useful and
interesting properties of groupoid crossed products.  In Sections 
\ref{sec:transitive} and \ref{sec:ideals} we mainly deal
with technical results.  In particular, since it is the first really
high level portion of the text, Section \ref{sec:transitive} contains
restatements of quite a few theorems which are too complicated to
prove here.  On the other hand, this section also presents results
concerning transitive groupoid crossed products which, while basic, are
new.  In Section \ref{sec:unitary} and \ref{sec:locally-unitary} we
define the notion of unitary and locally unitary actions.
In particular, we show that for unitary actions the crossed product
reduces to a tensor product.  For locally unitary actions the results
are more interesting.  We show that these actions are characterized, up to
exterior equivalence, by a principal bundle.  We also show that any
principal bundle can be used to construct a locally unitary action.
Moving on, Chapter \ref{cha:fine-structure} contains the primary
results of the thesis.  Section \ref{sec:indreps} describes a
technique for inducing representations from a closed subgroupoid up to
the whole crossed product.  We will eventually this induction technique to
identify the spectrum of certain crossed product algebras.  In Section
\ref{sec:regularity} we show, as long as the orbit space of the
groupoid is $T_0$, that every irreducible representation of the
crossed product is equivalent to the induction of an irreducible
representation of a fibre.  Then in Section \ref{sec:crossedstab}
we show, whenever the stabilizers are abelian and continuously
varying and the orbit space is $T_0$, 
that this induction map factors to a homeomorphism of a
quotient of the spectrum of $A\rtimes S$ onto the spectrum of
$A\rtimes G$. In Chapter \ref{cha:examples} we apply these
results to various examples and special cases.  In particular, Section
\ref{sec:groupstab} contains a strengthening of the results of
Section \ref{sec:crossedstab} in the groupoid $C^*$-algebra case.
Section \ref{sec:redux} describes how these results can
be applied to transformation groupoids.  Moreover, examples are
given and the theory is connected back to similar results for
transformation group algebras.  The last portion of the thesis is Section
\ref{sec:haussdorff}, which contains an analysis of when groupoid
$C^*$-algebras have Hausdorff spectrum.  While no conclusive answer is
given, this section has some intriguing constructions as well as
several nice counterexamples.  

We have outlined the logical structure of the thesis in the following
diagram. In particular, those readers interested only in induction and the fine
structure result may skip the branch containing 
Chapter \ref{cha:bundles} while whose readers only
interested in locally unitary actions can ignore Chapters
\ref{cha:fine-structure} and \ref{cha:examples}.  
\begin{figure}[ht]
\begin{center}
\includegraphics[width=5in]{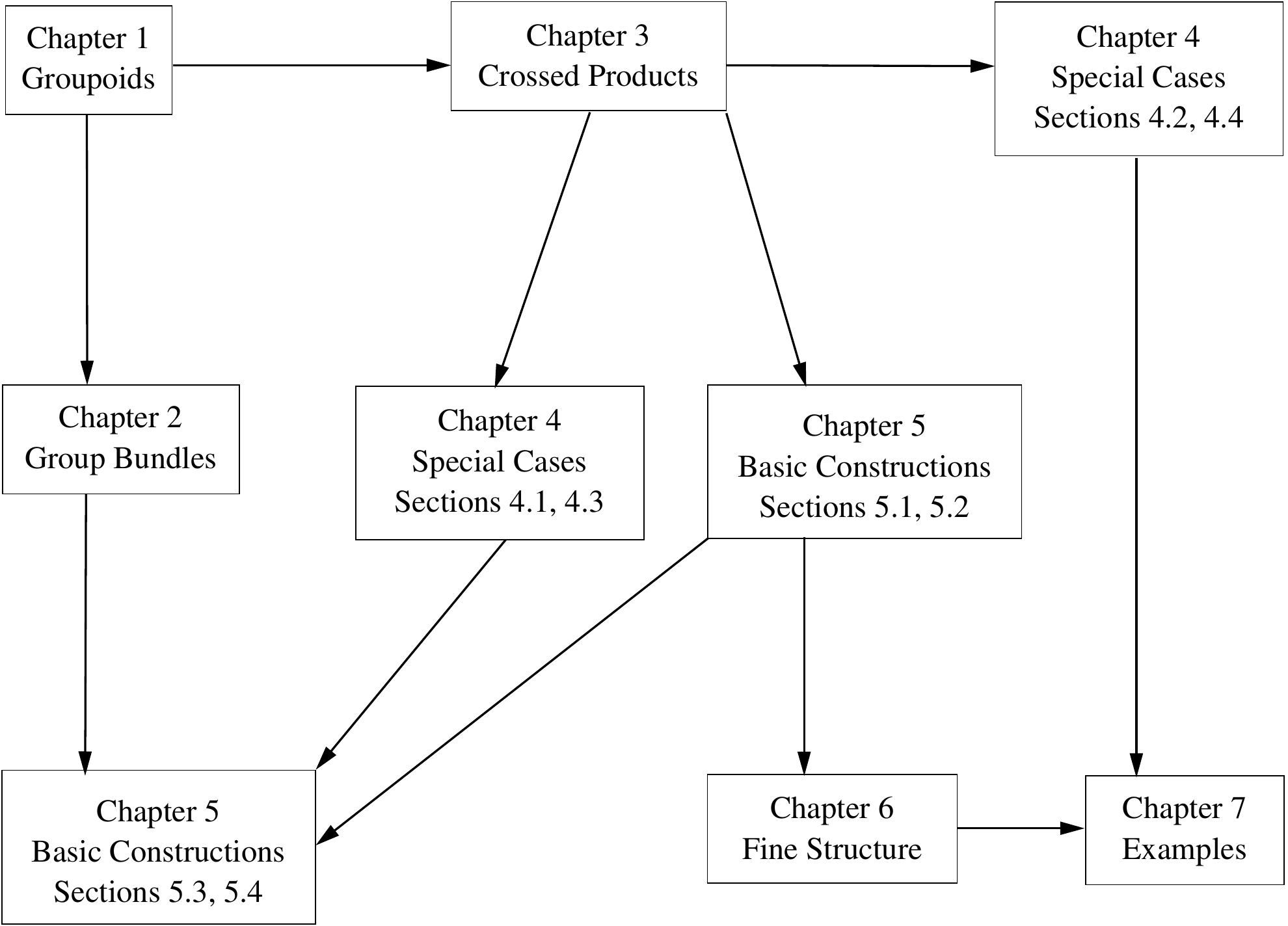}
\caption{The logical structure of the thesis.}
\label{fig:chapters}
\end{center}
\end{figure}

\subsection*{Acknowledgments}
\addcontentsline{toc}{subsection}{Acknowledgments}
First and foremost I would like to thank my adviser, Dana
Williams.\index{Dana Williams}
He has answered all of my incessant questions patiently ever since I got here
five years ago.  I would also like to thank the host of other people,
including Jonathan Brown and Paul Muhly, 
who have helped me with proofs, supplied clever arguments, or just
made their knowledge available.  Lastly, I would like to thank Naomi
for being awesome.


\tableofcontents

\mainmatter

\chapter{Groupoids}
\label{cha:groupoids}
In this chapter we present an overview of basic groupoid theory.  In
Section \ref{sec:basics} we define what a groupoid is and
outline some of the elementary facts and notation.  We also discuss
the notion of a stabilizer groupoid and an orbit groupoid.  The
stabilizer groupoid will play an important role through out.   Section
\ref{sec:actions} contains the basic constructions concerning
groupoid actions.  It is notable that we will separate the condition that
the structure map be open from the usual definition of a groupoid
action, see Remark \ref{rem:3}. 
The remainder of the section is basically review.  We define
groupoid equivalence and give a construction of the imprimitivity
groupoid and then in Section \ref{sec:amenable} we give the
briefest description of groupoid amenability.  

\section{Groupoid Basics}
\label{sec:basics}

Groupoids are essentially groups with a partially defined
multiplication.  While it may not seem like much, this has a
tremendous impact on their structure.  
This section will introduce some of the basic properties of groupoids, 
but we must start with their definition. The following draws heavily from
\cite{coords}.  

\begin{definition}
\label{def:1}
\index{groupoid}
Suppose $G$ is a set and $G^{(2)}\subseteq G\times G$.  Then $G$ is a
{\em groupoid} if there are maps $(\gamma,\eta)\mapsto \gamma\eta$
from $G^{(2)}$ into $G$ and $\gamma\mapsto \gamma\inv$ from $G$ into
$G$ such that:
\begin{enumerate}
\item {\em (associativity)} If $(\gamma,\eta)$ and $(\eta,\xi)$ are in $G^{(2)}$ then so are
  $(\gamma\eta,\xi)$ and $(\gamma,\eta\xi)$, and we have
  $(\gamma\eta)\xi = \gamma(\eta\xi)$.  
\item {\em(involution)} For all $\gamma \in G$ we have $(\gamma\inv)\inv = \gamma$.
\item {\em (cancellation)} For all $\gamma \in G$ we have $(\gamma\inv,\gamma)\in G^{(2)}$
  and if $(\gamma,\eta)\in G^{(2)}$ then $\gamma\inv(\gamma\eta) =
  \eta$ and similarly $(\gamma\eta)\eta\inv = \gamma$. 
\end{enumerate}
The set $G^{(2)}$ is called the set of {\em composable pairs}, when
$(\gamma,\eta)\in G^{(2)}$ we say $\gamma$ and $\eta$ are {\em composable}, and
$\gamma\inv$ is called the {\em inverse} of $\gamma$.  
\end{definition}

Some of the formulas in Definition \ref{def:1} are reminiscent of the
usual group axioms.  However, unlike the group case, the partially defined
multiplication implies that many different elements of $G$ act like
units.

\begin{definition}
\label{def:2}
\index[not]{$G^u$} \index[not]{$G_u$}
\index[not]{$G\unit$} 
\index{unit space $G\unit$}
\index{range map $r$}
\index{source map $s$}
Suppose $G$ is a groupoid.  Then the set of elements of $G$ such that
$\gamma = \gamma\inv = \gamma^2$ is denoted $G\unit$ and is called
the {\em unit space}.  The map $r:G\rightarrow G\unit$ such that
$r(\gamma) = \gamma\gamma\inv$ is called the {\em range map} and the
map $s:G\rightarrow G\unit$ such that $s(\gamma) = \gamma\inv\gamma$
is called the {\em source map}.  Given $u\in G\unit$ we will use the
notation $G_u := s\inv(u)$ and $G^u := r\inv(u)$.  
\end{definition}

\begin{remark}
\label{rem:4}
Suppose $G$ is a groupoid with $A$ and $B$ subsets of $G$.  We will use
the notation 
\begin{align*}
AB=A\cdot B &:=\{\gamma\eta: \gamma\in A, \eta\in B, (\gamma,\eta)\in G^{(2)}\},
& A\inv &:= \{\gamma\inv : \gamma\in A\}.
\end{align*}
It's important to realize that $AB$ may be badly behaved.  For
instance, $AB$ may not contain either $A$ or $B$, and is actually
empty if $A\times B\cap G^{(2)}=\emptyset$.  
\end{remark}

Of course, given any new class of objects, there is also a new class of
homomorphisms.  

\begin{definition}
\label{def:61}
\index{groupoid homomorphism}
Suppose $G$ and $H$ are groupoids.  A map $\phi:G\rightarrow H$ is a
{\em groupoid homomorphism} if and only if whenever $(\gamma,\eta)\in
G^{(2)}$ then $(\phi(\gamma),\phi(\eta))\in H^{(2)}$ and in this case
$\phi(\gamma\eta) = \phi(\gamma)\phi(\eta)$.  If $\phi$ is also
bijective then its called a {\em groupoid isomorphism}.   
\end{definition}

The next proposition outlines some of the basic properties of the
range map, source map, and the elements of $G\unit$.  

\begin{prop}
\label{prop:1}
Suppose $G$ is a groupoid.  
\begin{enumerate}
\item Given $\gamma,\eta\in G$ we have $(\gamma,\eta)\in G^{(2)}$ if and
  only if $s(\gamma) = r(\eta)$.  
\item If $(\gamma,\eta)\in G^{(2)}$  then $r(\gamma\eta) = r(\gamma)$
  and $s(\gamma\eta) = s(\eta)$.  
\item If $\gamma\in G$ then $r(\gamma) = s(\gamma\inv)$ and $s(\gamma)
  = r(\gamma\inv)$.  
\item If $(\gamma,\eta)\in G^{(2)}$ then $(\eta\inv,\gamma\inv)\in
  G^{(2)}$ and $(\gamma\eta)\inv = \eta\inv\gamma\inv$.  
\item If $\gamma\in G$ then $r(\gamma),s(\gamma)\in G\unit$. 
      Furthermore, $r$ and $s$ are retractions onto
      $G\unit$.\footnote{Given a set $X$ and a subset $A$ of $X$ a map
        $f:X\rightarrow A$ is a retraction if $f$ restricted to $A$ is
        the identity.} 
\item If $\gamma\in G$ then $(r(\gamma),\gamma),(\gamma,s(\gamma))\in
  G^{(2)}$, $r(\gamma)\gamma = \gamma$, and $\gamma
  s(\gamma) = \gamma$.  
\end{enumerate}
\end{prop}

\begin{proof}
Part {\bf(a)}:  Suppose $(\gamma,\eta)\in G^{(2)}$.  We know from the
cancellation condition 
of Definition \ref{def:1} that $(\gamma\inv,\gamma)\in G^{(2)}$.
Using associativity we have $(\gamma\inv\gamma)\eta =
\gamma\inv(\gamma\eta)$ which, after applying cancellation, gives us
$(\gamma\inv\gamma)\eta = \eta$.  Using the first part of cancellation
again,
along with involution, we get $(\eta,\eta\inv)\in G^{(2)}$, allowing
us to multiply the previous equality by $\eta$.  This yields
\[
((\gamma\inv\gamma)\eta)\eta\inv = \eta\eta\inv.
\]
Finally, using cancellation once again, we conclude
\[
s(\gamma) = \gamma\inv\gamma = \eta\eta\inv = r(\eta). 
\]
Next, suppose $s(\gamma) = r(\eta)$.  Once more, condition (c) of Definition
\ref{def:1} tells us that $(\gamma\inv,\gamma)\in G^{(2)}$, and if we
use involution we can similarly conclude that $(\gamma,\gamma\inv) \in G^{(2)}$.
Associativity implies that $(\gamma,\gamma\inv\gamma)\in
G^{(2)}$.  Since $\gamma\inv\gamma = s(\gamma) = r(\eta) =
\eta\eta\inv$ it follows that $(\gamma,\eta\eta\inv)\in G^{(2)}$.
Next, using Definition \ref{def:1} on $\eta$ in a similar fashion, we have
$(\eta\eta\inv, \eta)\in G^{(2)}$ and therefore, by associativity, 
$(\gamma,(\eta\eta\inv)\eta)\in G^{(2)}$.  However, cancellation
implies that $(\eta\eta\inv)\eta = \eta$ so that $(\gamma,\eta)\in
G^{(2)}$.  

Part {\bf (c)}:  Suppose $\gamma\in G$.  Then using involution
\[
r(\gamma\inv) = \gamma\inv (\gamma\inv)\inv = \gamma\inv\gamma =
s(\gamma). 
\]
The calculation that $s(\gamma\inv) = r(\gamma)$ is similar.  

Part {\bf (b)}:  Suppose $(\gamma,\eta)\in G^{(2)}$.  Then
$((\gamma\eta)\inv,\gamma\eta)\in G^{(2)}$, as well as
$(\eta,\eta\inv)$.  Applying associativity to $(\gamma,\eta)$ and
$(\eta,\eta\inv)$ gives us $(\gamma\eta,\eta\inv)\in G^{(2)}$.  Applying
associativity again to $(\gamma\eta, \eta\inv)$ 
and $((\gamma\eta)\inv,\gamma\eta)$ implies
$((\gamma\eta)\inv, (\gamma\eta)\eta\inv)\in G^{(2)}$.  Using
cancellation we conclude that $((\gamma\eta)\inv,\gamma)\in G^{(2)}$.
It follows from (c) and (a) that
\[
r(\gamma\eta) = s((\gamma\eta)\inv) = r(\gamma). 
\]
The calculation which shows $s(\gamma\eta) = s(\eta)$ is similar.  

Part {\bf (d)}:  Suppose $(\gamma,\eta)\in G^{(2)}$.  Then
$s(\gamma)=r(\eta)$ and using part (c) we have $r(\gamma\inv) =
s(\eta\inv)$ so that $(\eta\inv,\gamma\inv)\in G^{(2)}$. Applying
cancellation we get
\[
\eta\inv\gamma\inv = ((\eta\inv\gamma\inv)(\gamma\eta))(\gamma\eta)\inv.
\]
However, applying associativity and cancellation, we have
\begin{align*}
\eta\inv\gamma\inv &= (\eta\inv\gamma\inv(\gamma\eta))(\gamma\eta)\inv \\
&= (\eta\inv(\gamma\inv(\gamma\eta)))(\gamma\eta)\inv \\ 
&= (\eta\inv\eta)(\gamma\eta)\inv.
\end{align*}
Technically, we have to know that $\eta$ and $(\gamma\eta)\inv$ are
composable
before we can apply associativity.  However part (b) implies $s(\eta)
= s(\gamma\eta) = r((\gamma\eta)\inv)$ so that
$(\eta,(\gamma\eta)\inv) \in G^{(2)}$.  Now we can use associativity
and conclude
\begin{align*}
\eta\inv\gamma\inv &= (\eta\inv\eta)(\gamma\eta)\inv \\
&= \eta\inv(\eta(\gamma\eta)\inv) = (\gamma\eta)\inv.
\end{align*}

Part {\bf (e)}:  If $\gamma\in G$ then $s(\gamma)\inv =
(\gamma\inv\gamma)\inv = \gamma\inv\gamma = s(\gamma)$ by involution
and part (d).  Since $s(\gamma)\inv = s(\gamma)$ we conclude from part
(c) that the range and source of $s(\gamma)$ are equal and therefore 
$(s(\gamma),s(\gamma))\in G^{(2)}$.  Finally, 
\[
s(\gamma)s(\gamma) = \gamma\inv\gamma\gamma\inv\gamma = \gamma\inv
\gamma = s(\gamma)
\]
by cancellation.  Thus $s(\gamma)\in G\unit$.  Since $r(\gamma) =
s(\gamma\inv)$ this also shows $r(\gamma)\in G\unit$. Next, if $u\in
G\unit$ then $s(u) = u\inv u = u^2 = u$ and $r(u) = u
u\inv = u^2 = u$.  Hence $r$ and $s$ are retractions onto
$G\unit$.  

Part {\bf (f)}: Suppose $\gamma\in G$.  Then $s(r(\gamma)) =
s(\gamma\gamma\inv) = r(\gamma)$ by parts (b) and (c).  Thus $r(\gamma)$ and
$\gamma$ are composable and by cancellation 
\[
r(\gamma)\gamma = \gamma\gamma\inv\gamma = \gamma.
\]
The proof that $\gamma$ and $s(\gamma)$ are composable and that $\gamma
s(\gamma) = \gamma$ is similar.  
\end{proof}

\begin{remark}
\index{range map $r$}
\index{source map $s$}
Properties (a) - (b) of Proposition \ref{prop:1} help explain the
terminology behind the range and source maps and why groupoid
elements are sometimes called arrows.  Colloquially, every arrow in a
groupoid has a range and a source given by $r$ and $s$.  
Arrows are composable if and only
if the range of one matches the source of the second, and the range
and source of the composition are exactly what you would expect them
to be.  Properties (c) - (d) show that the inverse of an arrow goes in
the ``opposite direction'' and that
composition and inverses get along nicely.  Property (f) explains why
elements of $G\unit$ are called units; they act like identities on
elements with which they are composable. 

An alternative source for intuition regarding the partially defined
multiplication, and another reason for the arrow terminology, is to
think of a groupoid as a (small) category where every morphism is an
isomorphism.  This is actually equivalent to Definition \ref{def:1}.  
\end{remark}

We can also describe some basic properties of groupoid homomorphisms. 

\begin{prop}
\index{groupoid homomorphism}
Suppose $G$ and $H$ are groupoids and $\phi:G\rightarrow H$ is a
groupoid homomorphism.
\begin{enumerate}
\item Given $u\in G\unit$ we have $\phi(u)\in H\unit$.  
\item Given $\gamma\in G$ we have $\phi(\gamma\inv) =
  \phi(\gamma)\inv$.
\item For all $\gamma\in G$ we have $r(\phi(\gamma))=\phi(r(\gamma))$
  and $s(\phi(\gamma))=\phi(s(\gamma))$.  
\end{enumerate}
\end{prop}
\begin{proof}
Part {\bf (a)}:  Suppose $u\in G\unit$.  Then 
\[
\phi(u) = \phi(u^2) = \phi(u)\phi(u).
\]
Composing both sides with $\phi(u)\inv$ and using cancellation yields
\[
\phi(u)\phi(u)\inv = (\phi(u)\phi(u))\phi(u)\inv = \phi(u).
\]
It follows that $\phi(u) = r(\phi(u)) \in H\unit$.  

Part {\bf (b)}:  Suppose $\gamma\in G$.  Since $\gamma$ and $\gamma\inv$ are
composable we have $\phi(\gamma\gamma\inv) =
\phi(\gamma)\phi(\gamma\inv)$.  Next, we can compose both sides with
$\phi(\gamma)\inv$ and use cancellation to obtain
\[
\phi(\gamma)\inv \phi(\gamma\gamma\inv) = \phi(\gamma)\inv
\phi(\gamma) \phi(\gamma\inv) = \phi(\gamma\inv).
\]
However $\gamma\gamma\inv =
r(\gamma)\in G\unit$.  By part (a) we know $\phi(r(\gamma))\in
H\unit$ and using Proposition \ref{prop:1} to view $\phi(r(\gamma))$
as a right identity we have 
\[
\phi(\gamma)\inv = \phi(\gamma\inv). 
\]

Part {\bf (c)}:  Using part (b), we have
\[
\phi(r(\gamma)) = \phi(\gamma\gamma\inv) =
\phi(\gamma)\phi(\gamma)\inv = r(\phi(\gamma)).
\]
The proof for the source map is exactly the same.  
\end{proof}

In order to do any interesting functional analysis using groupoids you
have to assume that there is a topology floating around, or at least a
Borel structure. 

\begin{definition}
\label{def:3}
\index{groupoid}
Suppose $G$ is a groupoid with a topology and $G^{(2)}$ is endowed
with the relative product topology.  Then $G$ is a {\em
  topological groupoid} if the maps $(\gamma,\eta)\mapsto \gamma\eta$
from $G^{(2)}$ to $G$ and $\gamma\mapsto\gamma\inv$ from $G$ to $G$
are continuous.  If $G$ has a Borel structure such that $G^{(2)}$ is a
Borel subset of $G\times G$ and the above maps are Borel then we call
$G$ a {\em Borel groupoid}.  Furthermore, if $G$ is a topological
groupoid then we view $G$ as a Borel groupoid with the Borel structure
coming from the topology, after we give $G^{(2)}$ the relative product
Borel structure.   
\end{definition}

\begin{remark}
Almost without exception we will only be interested in topological
groupoids where the topology is locally compact Hausdorff.
Oftentimes, we will also assume that the topology is second
countable.
\end{remark}

\begin{prop}
If $G$ is a locally compact Hausdorff 
groupoid then 
\begin{enumerate}
\item the range and source maps are continuous,
\item the unit space $G\unit$ is closed in $G$, and
\item the set $G^{(2)}$ is closed in $G\times G$.  
\end{enumerate}
\end{prop}
\begin{proof}
It is clear that $r$ and $s$ are continuous since the 
composition and inversion operations are continuous.  
Now, suppose $\{u_i\}\in G$ is a net and $u_i \rightarrow u$.  Using
the fact that $r$ is continuous we have $r(u_i) \rightarrow r(u)$.
However, $r$ is a retraction onto $G\unit$ by Proposition
\ref{prop:1} so $r(u_i) = u_i$.  It follows that $u_i\rightarrow
r(u)$.  Since $G$ is Hausdorff $u = r(u)\in G\unit$ and $G\unit$ is
closed.  Finally, suppose $(\gamma_i,\eta_i)\in G^{(2)}$ and
$(\gamma_i,\eta_i)\rightarrow (\gamma,\eta)$.  Then we have
$\gamma_i\rightarrow \gamma$ and $\eta_i\rightarrow \eta$.  It follows
that $s(\gamma_i)\rightarrow s(\gamma)$ and $r(\eta_i)\rightarrow
r(\eta)$.  However, $s(\gamma_i)=r(\eta_i)$ for all $i$ and $G$ is
Hausdorff so $s(\gamma)=r(\eta)$ and $(\gamma,\eta)\in G^{(2)}$. 
\end{proof}

An important class of groupoids are those for which the unit space is
also open. We will see later that they are very rigid objects with
some nice properties. 

\begin{definition}
\index{r-discrete@$r$-discrete groupoid}
Suppose $G$ is a locally compact Hausdorff groupoid.  If $G\unit$ is
open in $G$ then we say that $G$ is an {\em $r$-discrete} groupoid. 
\end{definition}

\begin{remark}
We are using the older definition of $r$-discrete as given in
\cite{groupoidapproach}.  However, this definition has fallen out of
favor.  Currently $r$-discrete groupoids are those for which the unit
space is open and the range map is a local homeomorphism.  These
groupoids are also called etal\'e groupoids.  We will see in
Proposition \ref{prop:111}
that this is equivalent to assuming that the groupoid is
$r$-discrete, in the classical sense, and has a Haar system. 
\end{remark}

\begin{figure}[hb]
\centering
\includegraphics[height=2in]{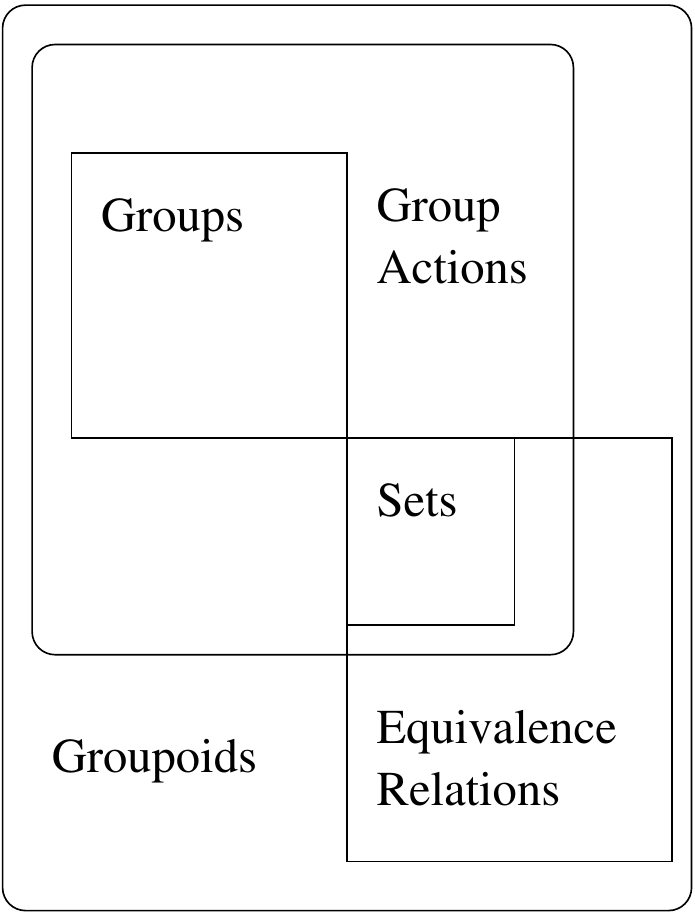}
\caption{Groupoids generalize many different objects.}
\label{fig:groupdiag}
\end{figure}

Groupoids are very general objects and extend a number of
well understood structures.  The following examples show how groupoids
generalize groups, sets, equivalence relations, and transformation
groups as in Figure \ref{fig:groupdiag}.

\begin{example}
\label{ex:1}
Suppose $H$ is a locally compact Hausdorff group.  If we let $H^{(2)}
= H\times H$ and give $H$ its group operations then $H$ is a locally
compact Hausdorff groupoid.  In this case $r(g) = s(g) = e$ for all
$g\in H$, where $e$ is the identity of $H$. 
\end{example}

\begin{example}
\index{cotrivial groupoid}
\label{ex:2}
Suppose $X$ is a locally compact Hausdorff space.  If we let 
$X^{(2)}$ be the diagonal in $X\times X$ then, 
with the trivial operations $(x,x)\mapsto x$ and
$x\inv \mapsto x$, it is easy to see that 
$X$ is a locally compact Hausdorff groupoid.  
In this case $X\unit = X$, every element is a
unit, and $X$ is known as a ``cotrivial'' groupoid.
\end{example}

\begin{example}
\label{ex:3}
\index{trivial groupoid}
Let $X$ be a locally compact Hausdorff space.  Suppose $R\subset X\times
X$ is locally compact Hausdorff in the relative topology and defines
an equivalence relation on $X$ by $x\sim y$ if and only if $(x,y)\in
R$.  
We let 
\[
R^{(2)} = \{((x,y),(w,z))\in R\times R : y=w \}
\]
and define 
\begin{align*}
(x,y)(y,z) &:= (x,z), &
(x,y)\inv &:= (y,x).
\end{align*}
With these operations $R$ is a locally compact Hausdorff groupoid.  
The unit space is $R\unit =
\{(x,x): x\in X\}$, and we usually make the obvious identification of $R\unit$ with $X$.
Under this identification $r(x,y) = x$ and $s(x,y) = y$.  
If $R =
X\times X$ then $R$ is called the ``trivial'' groupoid.  
If $R =
R\unit\cong X$ then $X$ is the ``cotrivial'' groupoid from Example \ref{ex:2} 
\end{example}

\begin{example}
\label{ex:4}
\index{transformation groupoid}
\index{transformation group groupoid|see{transformation groupoid}}
Suppose $H$ is a locally compact Hausdorff group acting on a locally
compact Hausdorff space $X$.  Let $G = H\times X$ and 
\[
G^{(2)} = \{((g,x),(h,y))\in G\times G : y = g\inv \cdot x\}.
\]
Given $((g,x),(h,y))\in G^{(2)}$ we define  
\begin{align*}
(g,x)(h,y) &:= (gh,x), &
(g,x)\inv &:= (g\inv, g\inv \cdot x).
\end{align*}
It's not hard to see that with these operations $G$ is a locally
compact Hausdorff groupoid.  The range and source maps are
\begin{align*}
s(g,x) &= (g\inv,g\inv\cdot x)(g,x) = (g\inv g,g\inv\cdot x) =
(e,g\inv \cdot x),\quad\text{and} \\
r(g,x) &= (g,x)(g\inv,g\inv\cdot x) = (gg\inv,x) = (e,x)
\end{align*} 
where $e$ is the unit of $H$.  In this case $G\unit = \{(e,x):x\in
X\}$ and we 
will usually identify $G\unit$ with $X$.  Under this identification
$s(g,x)= g\inv\cdot x$ and $r(g,x)=x$.  This type of
groupoid is called a ``transformation group groupoid'' or just a
``transformation groupoid'' for short.  Transformation group groupoids
generalize group actions in the sense that the group action is
completely determined by the associated transformation groupoid.  
\end{example}

\begin{remark}
Suppose $G= H\times X$ is a transformation groupoid.  If $H$ is
discrete then $X$ is open in $G$ and $G$ is $r$-discrete.  Similarly
if $X$ is open in $G$ then the identity must be open (as a singleton)
in $H$, and therefore $H$ is discrete.  Thus $r$-discrete
transformation groupoids correspond to discrete group actions.  
One of the reasons that $r$-discrete groupoids are important is because
they generalize discrete group actions in this way.  
\end{remark}

Examples \ref{ex:1} through \ref{ex:3}  
are all slightly degenerate in some sense.
Example \ref{ex:4} describes a class of groupoids which is much more
general.  In fact, much of the inspiration for groupoids can be traced
back to the transformation group case.  

\begin{remark}
\index{unit space $G\unit$}
In Examples \ref{ex:3} and \ref{ex:4} we were able to identify the
unit space of $G$ with an associated space not contained in $G$.  This
kind of identification happens fairly frequently, and we will often
treat $G\unit$ as if it exists ``outside'' $G$.
\end{remark}

There are also many groupoids which do not come from
transformation groups, equivalence relations, or one of the examples
presented above.  

\begin{example}
\label{ex:9}
\index{Deaconu-Renault groupoid}
Suppose $X$ is a {\em compact} Hausdorff space and $\sigma:
X\rightarrow X$ is a covering map.  
Let
\[
G = \{(x,n,y)\in X\times \Z \times X : \exists\ k,l \geq 0\
\text{s.t.}\ n = l-k,
\sigma^k x = \sigma^l y\}. 
\] 
Then define 
\[
G^{(2)} = \{((x,n,y),(w,m,z))\in G\times G: y=w\}
\]
and give $G$ the operations
\begin{align*}
(x,n,y)(y,m,z)&:= (x,n+m,z), &
(x,n,y)\inv &:= (y,-n,x).
\end{align*}
With these operations $G$ is a groupoid with unit space $G\unit =
\{(x,0,x)\in G:x\in X\}$.  We usually make the obvious identification of 
$G\unit$ with $X$.  Under this identification $r(x,n,y) = x$ and
$s(x,n,y) = y$.  Furthermore,
in these circumstances $G$ carries a topology making it into a
locally compact Hausdorff $r$-discrete groupoid \cite[Theorem
1]{deaconuendo}.  
This is known as the ``Deaconu-Renault groupoid'' associated to $(X,\sigma)$.  
\end{example}

\begin{example}
\label{ex:6}
\index{graph groupoid}
Suppose $E = (E^0,E^1,r,s)$ is a row-finite\footnote{A directed graph
  is {\em row-finite} if each vertex emits at most finitely many
  edges.} directed graph without
sources.  Let
$E^\infty$ denote the infinite path space of $E$.  Two paths
$\alpha,\beta\in E^\infty$ are shift equivalent with lag $n\in \Z$,
denoted $\alpha\sim_n \beta$, if there exists $N\in \N$ such that
$\alpha_i = \beta_{i+n}$ for all $i\geq N$.  Let
\[
G = \{(\alpha,n,\beta)\in E^\infty\times \Z\times E^\infty: \alpha
\sim_n \beta\}. 
\]
Next, define 
\[
G^{(2)} = \{((\alpha,n,\beta),(\gamma,m,\delta))\in G\times G : \beta=\gamma\}
\]
and let 
\begin{align*}
(\alpha, n,\beta)(\beta, m,\delta) &:= (\alpha, n+m, \delta), &
(\alpha, n,\beta)\inv &:= (\beta,-n,\alpha).
\end{align*}
Then $G$ is a groupoid.  The unit space $G\unit = \{(\alpha,
0,\alpha)\in G: \alpha \in E^\infty\}$ can be naturally identified with $E^\infty$ and the range and source maps are given by ${r(\alpha,n,\beta) = \alpha}$ and
$s(\alpha, n,\beta) = \beta$.  It is shown in 
\cite[Proposition 2.6]{graphgroupoid}
that $G$ carries a topology making it into a locally compact Hausdorff
$r$-discrete groupoid, called the ``graph groupoid'' associated to $E$. 
\end{example}

\begin{remark}
The reason that the groupoids in Examples \ref{ex:9} and \ref{ex:6}
look so similar is that they are both 
associated to generalizations of Cuntz-Krieger algebras \cite{deconeucuntz,graphgroupoid}. 
\end{remark}

Haar measure is essential to the study of locally compact groups
because it allows one to integrate.
We will also want to integrate on groupoids and to do that we will
need the following generalization of Haar measure.  

\begin{definition} 
\label{def:6}
\index{Haar system}
\index[not]{$\{\lambda^u\}$}
\index[not]{$\{\lambda_u\}$}
A (left) {\em Haar system} on a locally compact
  Hausdorff groupoid $G$ is a family $\lambda = \{\lambda^u\}_{u\in G\unit}$ of
    non-negative Radon measures on $G$ such that 
\begin{enumerate}
\item $\supp(\lambda^u) = G^u$ for all $u\in G\unit$,\footnote{For a
    Borel measure $\mu$ on a topological space $X$ the support of
    $\mu$, denoted $\supp\mu$, is defined to be the largest (closed)
    subset of $X$ for which every open neighborhood of every point of
    the set has positive measure.}
\item for $f\in C_c(G)$ the function 
\[
u \mapsto \int_G f(\gamma) d\lambda^u(\gamma)
\]
on $G\unit$ is in $C_c(G\unit)$; and
\item for $\gamma\in G$ we have
  $\gamma\lambda^{s(\gamma)}= \lambda^{r(\gamma)}$.  In other words,
  given $f\in C_c(G)$, 
\[
\int_G f(\gamma\eta)d\lambda^{s(\gamma)}(\eta) = \int_G
f(\eta)d\lambda^{r(\gamma)}(\eta).
\]
\end{enumerate}
Given $u\in G\unit$ we will use $\lambda_u$ to denote
$(\lambda^u)\inv$.  In other words, given $f\in C_c(G)$,
\[
\int_G f(\gamma)\lambda_u(\gamma) = \int_G f(\gamma\inv) \lambda^u(\gamma).
\]
\end{definition}

The following is a technical lemma which we will use ever so often. 

\begin{lemma}
\label{lem:10}
Given a groupoid $G$ with Haar system $\{\lambda^u\}$ and a compact
set $K\subset G$ then the set $\{\lambda^u(K)\}$ is bounded. 
\end{lemma}

\begin{proof}
Choose a compact neighborhood $L$ of $K$ and a positive function $f$ 
which is one on $K$ and zero off $L$.  Then $\lambda^u(K)\leq\int_G
f\lambda^u$ for all $u$ and the function $u\mapsto \int_G f\lambda^u$ is
continuous and compactly supported.  The result follows. 
\end{proof}

Unlike Haar measure, Haar systems are not always guaranteed to exist
and may not be unique in any reasonable sense.  What's more, only
groupoids with open range and source maps can have Haar systems.  The following
is asserted in \cite{groupoidapproach} and proved in \cite{seda}.

\begin{prop}
\label{prop:7}
If $G$ is a locally compact Hausdorff groupoid with a Haar
system then the range and source maps are open. 
\end{prop}

This is a good opportunity to mention a characterization of
surjective open maps which we will use constantly and is stated and
proved in \cite[Proposition 1.15]{tfb2}. 

\begin{prop}
\label{prop:9}
Let $p:X\rightarrow Y$ be a continuous surjection between two
topological spaces.  Then $p$ is an open map if
and only if given a net $\{y_i\}_{i\in I}$ converging to $p(x)$ in
$Y$, there is a subnet $\{y_{i_j}\}_{j\in J}$ and a net $\{x_j\}_{j\in
  J}$ indexed by the same set which converges to $x$ in $X$, and which
also satisfies $p(x_j) = y_{i_j}$.  
\end{prop}

Because Haar systems will be necessary to build groupoid
$C^*$-algebras we will usually assume that they exist.  Luckily, for
most of the groupoids that we are interested in there is a reasonable
Haar system.  

\begin{example}
\label{ex:24}
Suppose $X$ is a locally compact Hausdorff space and $G=X\times X$ is
the associated trivial groupoid.  Let $\lambda$ be any measure on
$X$ with full support and define $\lambda^x = \delta_x\times
\lambda$.  Then it is straightforward to show that 
$\{\lambda^x\}$ is a Haar system for $G$. Now suppose
we view $X$ as the cotrivial groupoid.  Then the collection of Dirac
delta measures $\{\delta_x\}$ forms a Haar system for $X$.  Since
integration against $\delta_x$ is just evaluation it's easy to see that
the continuity condition is satisfied and all of the other conditions
follow from the fact that the operations are ``cotrivial.''
\end{example}

\begin{example}
\label{ex:5}
\index{transformation groupoid!Haar system}
Suppose $H$ is a locally compact Hausdorff group acting on a locally
compact Hausdorff space $X$ and $G=H\times X$ is the associated
transformation groupoid.  Let $\lambda$ be a Haar measure for $H$ and
define $\lambda^x = \lambda\times \delta_x$.  Then $\{\lambda^x\}$ is
a Haar system for $G$.  We will always give transformation group
groupoids this Haar system.  
\end{example}

\begin{example}
\label{ex:7}
The groupoids in both Example \ref{ex:9} and Example \ref{ex:6} can be
given a Haar system by letting $\lambda^u$ be counting measure on
$G^u$ for all $u\in G^{(0)}$ \cite{deaconuendo,graphgroupoid}.  
\end{example}

The situation from Example \ref{ex:7} is actually much more generic.
The following proposition is proved in
\cite[Propositions 2.7,2.8]{groupoidapproach}. 

\begin{prop}
\label{prop:111}
\index{r-discrete@$r$-discrete groupoid}
Suppose $G$ is an $r$-discrete groupoid.
\begin{enumerate}
\item For any $u\in G^{(0)}$, $G^u$ and $G_u$ are discrete spaces. 
\item If $\{\lambda^u\}$ is a Haar system on $G$ then each $\lambda^u$
  is a multiple of the counting measure. 
\item The following are equivalent: 
\begin{enumerate}
\item $G$ admits a Haar system, 
\item $r$ and $s$ are local homeomorphisms, 
\item the product map $G^{(2)}\rightarrow G$ is a local
  homeomorphism. 
\end{enumerate}
\end{enumerate}
\end{prop}

\subsection{The Stabilizer Subgroupoid}

One slightly surprising fact is that a groupoid (potentially) contains
many different groups. 

\begin{prop}
\label{prop:2}
Suppose $G$ is a locally compact Hausdorff groupoid and $u \in
G\unit$.  Then $S_u=G_u\cap G^u = \{\gamma\in G: r(\gamma) = s(\gamma) = u\}$,
with the operations inherited from $G$, is a locally compact Hausdorff
group which is closed in $G$.  
\end{prop}

\begin{proof}
First, it's clear that $u\in S_u$ so that $S_u$ is not empty.  Now,
every element in $S_u$ has range and source $u$ so that any two
elements are composable.  Thus the groupoid operation is everywhere
defined on $S_u\times S_u$ and its associative because of the
associativity condition in Definition \ref{def:1}.  Given $\gamma \in
S_u$, since $s(\gamma)=r(\gamma) = u$, we know from Proposition
\ref{prop:1} that $\gamma u = u \gamma = \gamma$.  Finally, given
$\gamma \in S_u$ we have $\gamma\inv \gamma = s(\gamma) =u$ and
$\gamma\gamma\inv = r(\gamma) = u$.  Thus, $S_u$ is a group.  

Next suppose we have a net 
$\gamma_i \rightarrow \gamma$ in $G$ such that $\gamma_i\in S_u$
for all $i$.  The fact that the range and source maps are continuous
implies $r(\gamma_i) \rightarrow r(\gamma)$ and
$s(\gamma_i)\rightarrow s(\gamma)$.  However,
$r(\gamma_i)=s(\gamma_i)=u$ for all $i$ so clearly, because $G$ is
Hausdorff, ${r(\gamma)=s(\gamma)=u}$.  Thus $S_u$ is closed and it
follows that the relative topology on $S_u$ is locally compact
Hausdorff \cite[Lemma 1.26]{tfb2}.  Finally, since the operations are
continuous on $G$ they are continuous on $S_u$.  Thus $S_u$ is
a locally compact Hausdorff group.  
\end{proof}

These groups will play an important role and are given their own
special name.  

\begin{definition}
\label{def:4}
\index{stabilizer subgroup}
\index{stabilizer subgroupoid}
\index[not]{$S_u$}
Suppose $G$ is a groupoid and $u\in G\unit$ then the group 
\[
S_u = \{\gamma\in G: s(\gamma)=r(\gamma) = u\}
\]
is known as the {\em stabilizer subgroup} of $G$ at $u$.  The set 
\[
S := \{\gamma\in G : s(\gamma) = r(\gamma)\} = \bigcup_{u\in G\unit}
S_u
\]
is called the {\em stabilizer subgroupoid} of $G$.  Well use $p$ to
denote the restriction of the range (and source) map to $S$.  Oftentimes the
word {\em isotropy} is used interchangeably with {\em stabilizer}.  
\end{definition}

\begin{remark}
Since $S_u$ is a group we will generally denote elements of $S_u$ and
$S$ by lowercase Roman letters, instead of the Greek letters used to
denote generic elements of $G$. 
\end{remark}

\begin{prop}
\label{prop:3}
Suppose $G$ is a locally compact Hausdorff groupoid and let $S$ be the
stabilizer subgroupoid of $G$.  Then $S$ is a locally compact
Hausdorff subgroupoid which is closed in $G$.  Furthermore, $S$ 
can be viewed as a
bundle over $G\unit$ with bundle map $p$ whose fibres are the
isotropy subgroups.   
\end{prop}

\begin{proof}
Suppose $g,h\in S$.  If $g$ and $h$ are composable then there exists
$u\in G\unit$ such that $g,h\in S_u$ and it follows that $gh\in S_u\subset S$.
Similarly, we find that $S$ is closed under the inverse operation.
Since $S$ is closed under the operations inherited from $G$ it's
clear that $S$ is a groupoid in its own right, where $S^{(2)} =
\{(g,h)\in G^{(2)}: g,h\in S\}$.  Furthermore, since $s$ and $r$ are
continuous, if $\gamma_i\rightarrow \gamma$ in $G$
and $s(\gamma_i) = r(\gamma_i)$ for all $i$ then $s(\gamma) =
r(\gamma)$.  Thus $S$ is closed in $G$, and it follows that the
relative topology makes $S$ into a locally compact Hausdorff
groupoid. The remaining statements of the proposition are clear.  
\end{proof}

The stabilizer subgroupoid is a very important object.  It can
oftentimes tell us a lot about its parent groupoid.  An even better
situation is when the stabilizer subgroupoid is everything.  

\begin{definition}
\label{def:24}
\index{group bundle}
A {\em groupoid group bundle}, or {\em group bundle} for short, is a
locally compact Hausdorff groupoid $S$ such that the range and source
maps are equal.  We will denote the range (and source) map by
$p$.  We view $S$ as a bundle over $S\unit$ with bundle map $p$ and
denote the fibres by $S_u = p\inv(u)$ for all $u\in S\unit$.  We say
that $S$ is an {\em abelian group bundle} if $S_u$ is an abelian group
for all $u\in S\unit$.  
\end{definition}

\begin{remark}
Groupoid group bundles are different
than the kinds of group bundles one usually encounters. For instance,
groupoid group bundles carry no kind of local triviality condition and
the fibres can and will vary over the base space.  Furthermore, the
injection of the unit space into $S$ always 
gives a continuous section of the bundle map.  
\end{remark}

\begin{example}
Suppose we have a locally compact
Hausdorff group $H$ and a topological space $X$.  Then we can view
$S=X\times H$ as a group bundle where the bundle map is
just the projection onto the first factor.  In this case the unit space
can be identified with $X$ and the groupoid operations are obvious.  
\end{example}

It turns out that the existence of a Haar system on a group bundle $S$
is equivalent to requiring the fibres of $S$ vary ``continuously.''
In order to make this notion precise we recall the following from
\cite[Section H.1]{tfb2}. 

\begin{definition}
\label{def:8}
\index{Fell Topology}
Let $X$ be an arbitrary topological space and $\mcal{C}(X)$ the
collection of all closed subsets of $X$ (including the empty set).
Given a finite collection $\mcal{F}$ of open sets of $X$ and a
compact subset $K$ of $X$ we define
\[
\mcal{U}(K;\mcal{F}) := \{F\in \mcal{C}(X):\text{$F\cap K = \emptyset$
  and $F\cap U \ne \emptyset$ for all $U\in\mcal{F}$}\}. 
\]
The collection $\{\mcal{U}(K;\mcal{F})\}$ forms a basis for a compact
Hausdorff topology on $\mcal{C}(X)$ called the {\em Fell Topology}. 
\end{definition}

\begin{prop}
\label{prop:8}
Suppose $X$ is a locally compact space and let $\{F_i\}_i\in I$ be a net
in $\mcal{C}(X)$.  Then $F_i\rightarrow F$ in $\mcal{C}(X)$ if and
only if 
\begin{enumerate}
\item given $t_i\in F_i$ such that $t_i\rightarrow t$, then $t\in F$,
  and 
\item if $t\in F$, then there is a subnet $\{F_{i_j}\}$ and
  $t_{i_j}\in F_{i_j}$ such that $t_{i_j}\rightarrow t$.  
\end{enumerate}
\end{prop}

Using Definition \ref{def:8} to pin down the appropriate notion of
continuity we can make the following

\begin{definition}
\label{def:9}
\index{continuously varying stabilizers|see{stabilizer subgroupoid}}
\index{stabilizer subgroupoid!continuously varying}
Suppose $S$ is a locally compact Hausdorff group bundle.  We say that
$S$ is {\em continuously varying} if given a net $\{u_i\}$ in
$S\unit$ such that $u_i\rightarrow u$ in $S\unit$ then
$S_{u_i}\rightarrow S_u$ with respect to the Fell topology.  If $G$ is
a locally compact Hausdorff groupoid we say that $G$ has {\em
  continuously varying stabilizers} if its stabilizer subgroupoid $S$
is continuously varying. 
\end{definition}

At this point we can give some reasonable conditions for the existence
of a Haar system for a group bundle.  The following also provides a partial
converse to Proposition \ref{prop:7}.

\begin{prop}
\label{prop:10}
\index{group bundle!Haar system}
\index{group bundle!continuously varying}
Suppose $S$ is a locally compact Hausdorff groupoid group bundle with
bundle map $p$.  The following are equivalent:
\begin{enumerate}
\item $S$ has a Haar system, 
\item $p$ is open, 
\item $S$ is continuously varying. 
\end{enumerate}
\end{prop}
\begin{proof}
It is shown in \cite[Lemma 1.3]{renaultgcp} that (a) and (b) are
equivalent.  Now suppose $p$ is open and $u_i\rightarrow u$ in
$S\unit$.  We will show $S_{u_i}\rightarrow S_u$ using Proposition
\ref{prop:8}.  If $s_i\in S_{u_i}$ and $s_i\rightarrow s$ then
$p(s_i)=u_i\rightarrow p(s)$.  It follows that $p(s)=u$ and $s\in
S_u$.  Now suppose $s\in S_u$.  Because $p$ is open we can use
Proposition \ref{prop:9} to pass to a subnet and find $s_{i_j}\in
S_{u_{i_j}}$ such that $s_{i_j}\rightarrow s$.  Thus
$S_{u_i}\rightarrow S_{u}$.  

Next suppose $S$ varies continuously.  We will show $p$ is open using
the characterization in Proposition \ref{prop:9}.  Let $u_i\rightarrow
u$ be a net which converges in $S\unit$ and suppose $s\in S_u$.
Using Proposition \ref{prop:8} we can pass to a subnet and find
$s_{i_j}\in S_{u_{i_j}}$ such that $s_{i_j}\rightarrow s$.   Since
$p(s_{i_j}) = u_{i_j}$ this proves $p$ is open. 
\end{proof}

\begin{remark}
We will pass to subnets frequently, and may even make use of
sub-subnets.  
In order to avoid notational clutter we will usually relabel so that only one
index is shown.  
\end{remark}

\begin{remark}
\label{rem:5}
Given a locally compact groupoid group bundle $S$ with a Haar system
$\{\lambda^u\}$ each $\lambda^u$ is a measure supported on the group $S_u =
p\inv(u)$.  It follows from the left invariance condition of
Definition \ref{def:6} that $\lambda^u$ is a Haar measure on $S_u$.
Thus a Haar system on $S$ is just a ``continuously varying''
collection of
Haar measures.  It is straightforward 
to see that given another Haar system $\{\mu^u\}$ on $S$
the unicity of Haar measure guarantees the existence of a
continuous function $f:G\unit\rightarrow \C$ such that $\lambda^u =
f(u)\mu^u$ for all $u\in G\unit$.  
\end{remark}

\begin{example}
\label{ex:33}
Let $T$ be the cone of length 1 and maximum radius 1 aligned on the
positive $x$-axis as in Figure \ref{fig:bundle} 
and let $p_T$ be the projection onto the
$x$-axis.  It is fairly clear that $T$ is a continuously varying group
bundle with each fibre equal to the circle group, after a scaling.
Similarly let $Z$ be the collection of line segments given by $y =
nx$ for all $n \in \mathbb{Z}$ and $x \in [0,1]$ as in Figure
\ref{fig:bundle} and let $p_Z$ be
the projection onto the $x$-axis.  Then $Z$ is also a continuously
varying group bundle and each fibre is equal to the integers, again
after a scaling. 
\begin{figure}[h]
\begin{center}
\includegraphics[width=2.5in]{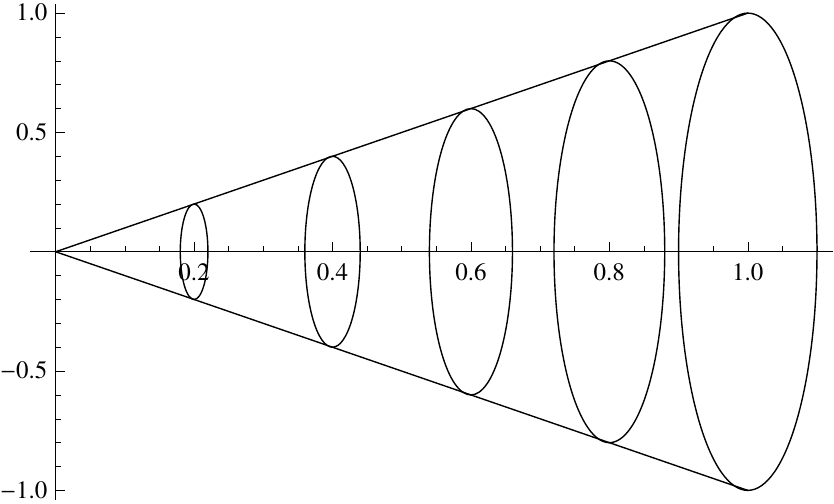}\hspace{.5in}
\includegraphics[width=2.5in]{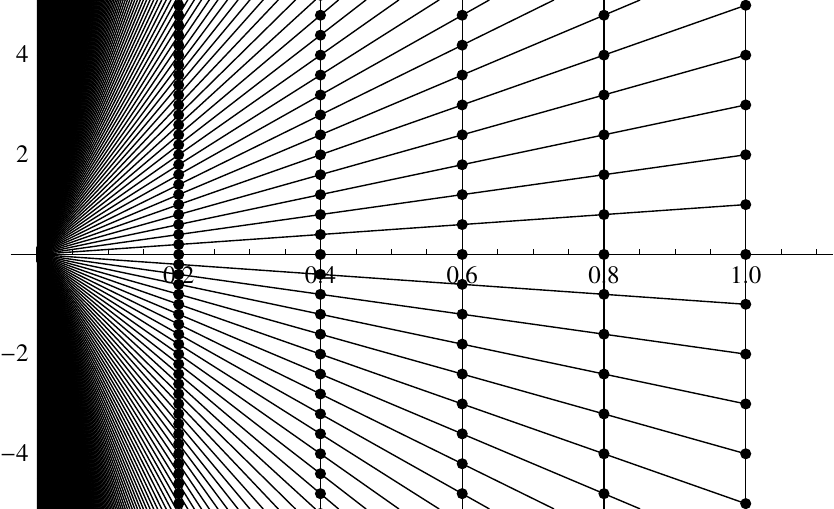}
\end{center}
\caption{Examples of continuously varying group bundles.}
\label{fig:bundle}
\end{figure}
\end{example}

\begin{remark}
If $G$ is a locally compact Hausdorff groupoid then Proposition
\ref{prop:10} implies that $G$ has continuously varying stabilizers if
and only if the stabilizer subgroupoid has a Haar system.  This is a
natural hypothesis and will be frequently invoked.  However, it's
important to understand that this is a very strong condition.  For
instance, if $G = H\times X$ is a transformation groupoid associated
to an action of $H$ on $X$ then it's easy to see that the stabilizer
subgroup $S_x$ is isomorphic to $\{s\in H : s\cdot x = x\}$.  In other
words, the stabilizer subgroup $S_x$ is exactly the stabilizer
subgroup of $H$ at $x$.  It's also straightforward to show that 
the stabilizers vary
continuously in $G$ if and only if they vary continuously in $H$.
However, most group actions do not have continuously varying
stabilizers; even really nice actions.  
For example, if $\mathbb{T}$ acts on $\R^2$ by rotation
then the stabilizers are discontinuous at the origin.  
\end{remark}

\subsection{The Orbit Groupoid}

Yet another important class of groupoid are those that are ``opposite''
of Definition \ref{def:24}.  In other words, groupoids whose isotropy
subgroupoid is as trivial as possible.  

\begin{definition}
\index{principal groupoid}
\index{groupoid!principal}
A groupoid $G$ is called {\em principal} if $r(\gamma)=s(\gamma)$
implies $\gamma\in G\unit$ for all $\gamma\in G$. 
\end{definition}

Any groupoid gives rise to a canonical principal groupoid.  

\begin{definition}
\label{def:7}
\index{orbit equivalence relation}
\index{orbit groupoid}
\index[not]{$\pi=(r,s)$}
Suppose $G$ is a groupoid.  We define the {\em orbit equivalence
  relation} on $G\unit$ to be given by
$u\sim v$ if and only if there exists
$\gamma\in G$ such that $u=r(\gamma)$ and $v=s(\gamma)$.  We define
the {\em orbit groupoid} to be $R = \{(u,v)\in
G\times G: u\sim v\}$ where $R^{(2)} = \{((u,v),(w,z))\in R\times R: v=w\}$
and the operations are given by
\begin{align*}
(u,v)(v,w) &:= (u,w), &
(u,v)\inv &:= (v,u).
\end{align*}
Finally, we call the map $\pi=(r,s):G\rightarrow R$ defined by
$\pi(\gamma) = (r(\gamma),s(\gamma))$ the {\em canonical homomorphism}.
\end{definition}

\begin{prop}
\label{prop:4}
Suppose $G$ is a groupoid.  Let $\sim$ be the orbit equivalence relation,
$R$ the orbit groupoid, and $\pi=(r,s)$ the canonical
homomorphism from Definition \ref{def:7}.  Then $\sim$ is an
equivalence relation, and $R$ is a principal groupoid.  The unit space
of $R$ can be identified with $G\unit$ and under this identification
$r(u,v) = u$ and $s(u,v)= v$.  Furthermore,
$\pi$ is a surjective groupoid homomorphism. 
\end{prop}

\begin{proof}
Since $u = r(u) = s(u)$ we see that $\sim$ is reflexive.  We
know $\sim$ is symmetric because if $u = r(\gamma)$ and $v=s(\gamma)$
then $u = s(\gamma\inv)$ and $v = r(\gamma\inv)$.  Finally if $u=
r(\gamma)$, $v = s(\gamma) = r(\eta)$ and $w=s(\eta)$ then, citing
Proposition \ref{prop:1}, $u = r(\gamma\eta)$ and $v = s(\gamma\eta)$.
Thus $\sim$ is transitive. 

As for the groupoid $R$, this is exactly the situation in Example
\ref{ex:3}.  Observe that the composition operation is well defined, 
i.e. it maps
into $R$, because of the transitivity of $\sim$.  Similarly the inverse
is well defined because of symmetry.  
Now, suppose 
\[
((x,y),(y,w)),((y,w),(w,z))\in R^{(2)},
\]
then it's clear that 
\[
((x,w),(w,z)),((x,y),(y,z))\in R^{(2)},
\] 
and that 
\[
((x,y)(y,w))(w,z) = (x,z) = (x,y)((y,w)(w,z)).
\]
Therefore associativity is satisfied, and it's obvious that involution
is also satisfied.  Finally, suppose $(x,y)\in R$.   Clearly
$((x,y),(y,x))\in R^{(2))}$ and given $(y,w)\in R$
\begin{align*}
(y,x)((x,y)(y,w)) &= (y,x)(x,w) = (y,w) \\
((x,y)(y,w))(w,y) &= (x,w)(w,y) = (x,y).
\end{align*}
Thus cancellation also holds and $R$ is a groupoid.  Given $(x,y)\in
R$ we have $r((x,y)) = (x,y)(y,x) = (x,x)$.  From here it is clear
that $G\unit = r(G) = \{(x,x): x\in G\unit\}$.  This set is
trivially identifiable with $G\unit$ and under this identification
$r(x,y) = x$.  Similarly $s(x,y) = y$. Lastly, if $r(x,y) = s(x,y)$
then $x=y$ and $(x,y)\in R\unit$, making $R$ principal. 

Observe that $(r(\gamma),s(\gamma))\in R$ for all $\gamma\in G$ so
$\pi$ is well defined.  Furthermore if $(\gamma,\eta)\in G^{(2)}$ then
$s(\gamma) = r(\eta)$ so that
$((r(\gamma),s(\gamma)),(r(\eta),s(\eta)))\in R^{(2)}$ and 
\begin{align*}
\phi(\gamma\eta) &= (r(\gamma\eta),s(\gamma\eta))
= (r(\gamma),s(\eta)) \\
&= (r(\gamma),s(\gamma))(r(\eta),s(\eta)) = \phi(\gamma)\phi(\eta). 
\end{align*}
Finally, it is clear from the definition of $R$ that $\pi$ is
surjective. 
\end{proof}

\begin{prop}
\label{prop:5}
If $G$ is a principal groupoid then the canonical homomorphism
$\pi:G\rightarrow R$ is an isomorphism. 
\end{prop}

\begin{proof}
We already know that $\pi$ is a surjective homomorphism. Suppose
$\pi(\gamma) = \pi(\eta)$ for $\gamma,\eta\in G$.  Then $r(\gamma)=
r(\eta)$ and $s(\gamma)= s(\eta)$.  This means that $\gamma$ and
$\eta\inv$ are composable and that 
\[
s(\gamma\eta\inv) = r(\eta) = r(\gamma) = r(\gamma\eta\inv).
\]
Therefore $\gamma\eta\inv\in G\unit$.  If we compose
$\gamma\eta\inv$ on both sides of the equation $\eta = \eta$ we get 
\[
\gamma = \gamma\eta\inv\eta = \eta
\]
where the left hand equality comes from cancellation and the right
hand equality comes from the fact that $\gamma\eta\inv$ is a unit.
Thus $\pi$ is injective and therefore an isomorphism. 
\end{proof}

\begin{remark}
Proposition \ref{prop:5} states that every principal groupoid is
isomorphic to its orbit groupoid.  However, the orbit
groupoid is determined by the orbit equivalence relation.  Thus,
disregarding topology for a moment, every
principal groupoid is (isomorphic to) one of the groupoids defined in
Example \ref{ex:3}. 
\end{remark}

The isotropy groupoid of a locally compact Hausdorff groupoid is
naturally a locally compact Hausdorff groupoid.  Unfortunately the
situation is not so simple for the orbit groupoid.  

\begin{definition}
\label{def:5}
\index[not]{$R_Q$}
\index[not]{$R_P$}
Suppose $G$ is a locally compact Hausdorff groupoid and $R$ is the
orbit groupoid determined by $G$. We denote $R$ with the relative
topology of ${G\unit\times G\unit}$ by $R_P$. We denote $R$ with the
quotient topology with respect to the canonical homomorphism $\pi=(r,s)$ by
$R_Q$.  When $R$ is used as a topological groupoid it will always mean
$R_Q$.
\end{definition}

\begin{prop}
\label{prop:6}
Suppose $G$ is a locally compact Hausdorff groupoid.  Then the
topology on $R_Q$ is finer than the topology on $R_P$.  Furthermore $R_P$ and
$R_Q$ are both Hausdorff and the map $\pi$ is continuous as
a function into both $R_Q$ and $R_P$.  Next, $R_P$ is a topological
groupoid and, if either $G$ is second countable or $R_Q$ is locally
compact, then $R_Q$ is a topological groupoid.  Finally, if $G$ has
open range and source then the range and source maps are open as maps
on $R_P$ and $R_Q$.
\end{prop}

\begin{proof}
First we will show $\pi:G\rightarrow R_P$ is continuous.  Suppose
$\gamma_i\rightarrow \gamma\in G$. Since the range and source maps are
continuous we have
$(r(\gamma_i),s(\gamma_i))\rightarrow(r(\gamma),s(\gamma))$ in
$G\unit\times G\unit$ and hence in $R_P$.  Since $R_Q$ has the
quotient topology determined by $\pi$, clearly $\pi:G\rightarrow R_Q$
must be continuous.  

Next, suppose $O$ is open in $R_P$, then $\pi\inv(O)$ is open in $G$,
but this implies that $O$ is open in $R_Q$.  Thus $R_Q$ has a finer
topology than $R_P$.  Furthermore, since any subset of a Hausdorff
space inherits a Hausdorff topology, $R_P$ is Hausdorff.  This implies
$R_Q$ is Hausdorff as well, since it carries a finer topology. 

It's pretty easy to see that the operations on $R_P$ are
continuous with respect to the product topology.  After all
$(u_i,v_i)\rightarrow (u,v)$  in $R_P$ if and only if $u_i\rightarrow u$ and
$u_i\rightarrow u $.  Proving that the operations are continuous on
$R_Q$ takes more work, and more hypotheses.  Let $I:G\rightarrow G$ be
defined by $I(\gamma) = \gamma\inv$ and consider $\pi\circ I:
G\rightarrow R_Q$.  It's easy to see that $I$ is a continuous map and
that if $\pi(\gamma)=\pi(\eta)$ then $\pi\circ I(\gamma) = \pi\circ
I(\eta)$.  Thus $I\circ \pi$ factors to a continuous map from $R_Q$ into
$R_Q$ and clearly this factorization is nothing more than the inversion
operation on $R_Q$.  We would like to use the same argument with the
multiplication.  The main issue is the following

\begin{claim}  The map $\pi\times \pi: G\times G \rightarrow R_Q\times
R_Q$ is a quotient map.  
\end{claim}

\begin{proof}[Proof of Claim.]
It is known \cite[Section 8]{biquotient} 
that the product of a quotient map with
itself need not be quotient.  However, there are rather minimal
conditions on $R_Q$ which will guarantee that $\pi\times \pi$ is a
quotient map.  We know from \cite[Theorem 1.5]{biquotient} that if 
$G$ and $R_Q\times R_Q$ are Hausdorff $k$-spaces\footnote{A
  topological space $X$ is called a {\em $k$-space} if a set $A\subset
  X$ is closed whenever $A\cap K$ is closed in $K$ for every compact
  $K\subset X$.  Such a space is also called {\em compactly
    generated}.} 
then $\pi\times \pi$ will be a quotient map. 
It's clear that $G$ and $R_Q\times R_Q$ are Hausdorff
and, since locally compact spaces are always $k$-spaces, all that is left
is to show that $R_Q\times R_Q$ is a $k$-space.  If, on one hand,
$R_Q$ is locally compact then $R_Q\times R_Q$ is locally compact and
we are done.  On the other hand, suppose $G$ is second countable.  Then
by choosing a countable basis of compact neighborhoods we can find a
countable collection $\{K_n\}$ of compact sets which cover $G$ and have the
property that a set $A\subset G$ is closed if and only if $A\cap K_n$ is
closed for all $n$.  Such a space is known as a $k_\omega$-space.  It
follows from \cite[Remark 7.5]{biquotient} that quotients and products
of $k_\omega$-spaces are $k_\omega$-spaces so that we can conclude $R_Q\times
R_Q$ is a $k_\omega$-space. Thus assuming either $G$ is second
countable or $R_Q$ is locally compact proves our claim. 
\end{proof}

Suppose the claim is satisfied.  It's
straightforward to see that $R_Q^{(2)}$ is closed in $R_Q\times R_Q$
and, since $G^{(2)}$ is just the set of those elements whose range and
sources match up, that $G^{(2)} = (\pi\times\pi)\inv(R^{(2)})$.  It
is also straightforward to see that the restriction of a quotient map
to the inverse image of a closed set results in a quotient map.  Thus
the restriction $\pi\times \pi: G^{(2)}\rightarrow R_Q^{(2)}$ is a
quotient map.  Let $M:G^{(2)}\rightarrow G$ be given by
$M(\gamma,\eta) = \gamma\eta$ and consider $\pi\circ M:G^{(2)}\rightarrow
R_Q$.  Then it's easy to see that if $\pi\times\pi(\gamma,\eta) =
\pi\times \pi(\gamma',\eta')$ then $\pi\circ M(\gamma,\eta) = \pi\circ
M(\gamma',\eta')$ so that $\pi\circ M$ factors to a continuous map
from $R^{(2)}_Q$ into $R_Q$.  Furthermore, this map is clearly the composition
operation on $R_Q$ implying that composition is continuous.    

Finally, suppose $G$ has open range and source maps, that
$u_i\rightarrow u$ in $G\unit$ and that $(u,v)\in R$.  First, choose
any $\gamma\in G$ such that $\pi(\gamma) = (u,v)$.  Since the range
on $G$ is open we can pass to a subnet and find
$\gamma_i\rightarrow\gamma$ such that $r(\gamma_i)= u_i$.  Then
$\pi(\gamma_i)\rightarrow \pi(\gamma)=(u,v)$ in both $R_P$ and $R_Q$ 
and $r(\pi(\gamma_i)) = u_i$.  Thus
$r$ is open on $R_P$ and $R_Q$.  The proof that $s$ is open is similar.  
\end{proof}

\begin{remark}
\label{rem:1}
The astute reader will have noticed that there are several things
missing from Proposition \ref{prop:6}.  For instance, neither $R_P$
nor $R_Q$ are necessarily locally compact. What's more, in extreme cases
it is not clear if $R_Q$ is even a topological groupoid.  On the other
hand, the operations on $R_Q$ are always continuous if $G$ is second
countable, which includes almost all of the examples that we will care
about.  It's also nice to note that if the topology on $R_Q$ is well
behaved (i.e. locally compact) then the operations are continuous then
as well.  However, even in this case $R_Q$ may not have a Haar system
when $G$ does.
\end{remark}

Interestingly enough, given a groupoid $G$ there is a duality between
its orbit groupoid and its stabilizer subgroupoid.  The
following is stated in \cite[Remark 1.2]{renaultgcp}. 

\begin{prop}
\label{prop:11}
Suppose $G$ is a locally compact Hausdorff groupoid.  Then $G$ has 
continuously varying stabilizers if and only if the canonical map
$\pi = (r,s)$ is open onto $R_Q$.  
Furthermore, under these conditions $R_Q$ is a
locally compact Hausdorff groupoid, and if $G$ is second countable
then $R_Q$ is also.
\end{prop}

\begin{proof}
Let $G$ be a locally compact Hausdorff groupoid, $S$ the stabilizer
subgroupoid, and $p$ the bundle map for $S$.  
First, suppose $G$ has continuously varying stabilizers so that
Proposition \ref{prop:10} implies that $p$ is open.  Let $O$ be an
open set in $G$.  We want to show that $\pi\inv(\pi(O))$ is open.
Well
\[
\pi\inv(\pi(O)) = \{\gamma\in G:\exists\ \eta\in O\ \text{s.t.}\
\pi(\gamma)=\pi(\eta)\}.
\]
However, if $\pi(\eta)=\pi(\gamma)$ then $\gamma\eta\inv \in S$.
Therefore, if $\gamma\in \pi\inv(\pi(O))$ then $\gamma\in S\cdot O =
\{s\eta: s\in S, \eta\in O, p(s)=r(\eta)\}$.  What's more, if
$s\eta\in S\cdot O$ then $\pi(s\eta)=\pi(\eta)$ so that $s\eta\in
\pi\inv(\pi(O))$.  It follows that $\pi\inv(\pi(O)) = S\cdot O$.  Now
suppose $S\cdot O$ is not open so that there exists $s\in S,\gamma\in
O$ and a net $\{\gamma_i\}$ such that $\gamma_i\rightarrow s\gamma$ and 
for all $i$ we have
$\gamma_i\notin S\cdot O$. Now, $r(\gamma_i)\rightarrow r(\gamma)$
and $p(s) = r(\gamma)$ so that, using Proposition \ref{prop:9}, we can
pass to a subnet, relabel, and find $s_i\in S_{r(\gamma_i)}$ such that
$s_i\rightarrow s$.  It follows that $s_i\inv \gamma_i \rightarrow
s\inv s\gamma = \gamma$.  Since $O$ is open this implies
$s_i\inv\gamma_i\in O$ eventually and this is a contradiction since
$\gamma_i = s_i(s_i\inv\gamma_i)$.  Hence $S\cdot O$ is open
so that $\pi$ is open onto $R_Q$.

Next, let $\pi:G\rightarrow R_Q$ be open.  Suppose $u_i\rightarrow u$
is a convergent net in $G\unit$ and that $p(s) = u$.  Observe that
$\pi(u_i)=(u_i,u_i)\rightarrow \pi(u)=(u,u)$ and that $\pi(s) =
(u,u)$.  Since $\pi$ is open we can use Proposition \ref{prop:9} to
pass to a subnet, relabel, and find $\gamma_i \in G$ such that
$\gamma_i\rightarrow s$ and $\pi(\gamma_i)=(u_i,u_i)$.  However, this
implies $\gamma_i \in S_{u_i}$ for all $i$.  Thus, using Proposition
\ref{prop:9} again, $p$ is an open map.  

Suppose $\pi$ is open onto $R_Q$, or equivalently that the stabilizers vary
continuously.  It follows that $R_Q$ is locally compact, since the image of a
basis of compact neighborhoods under $\pi$ will be a basis of compact
neighborhoods for $R_Q$.  In this case, Proposition
\ref{prop:6} implies that the operations on $R_Q$ are continuous.
Furthermore, since the image of a countable basis under $\pi$ will be a
countable basis, if $G$ is second countable then so is $R_Q$.
\end{proof}

\begin{remark}
It is not necessary for the stabilizers to vary continuously for $R_Q$
to be a locally compact Hausdorff (topological) groupoid.  
For instance, consider $\mathbb{T}$ acting on the closed
unit ball in $\R^2$ by rotation.  Then 
$R_Q$ is compact since it's the continuous image of a compact space,
and is therefore a topological groupoid.  However, the stabilizers are
clearly discontinuous at the origin. 
\end{remark}


\section{Groupoid Spaces}
\label{sec:actions}

The notion of a groupoid action on a space is a straightforward
generalization of group actions.  The only caveat is
that the action is only ``partially defined'' in the same sense that
the groupoid multiplication is only partially defined. Once again,
much of this section is inspired by \cite{coords}.  

\begin{definition}
\label{def:10}
\index{G-space@$G$-space}
\index{groupoid action|see{$G$-space}}
\index[not]{$G*X$}
Suppose $G$ is a groupoid and $X$ is a set.  We say that $G$ acts (on
the left) of $X$, and that $X$ is a {\em left $G$-space}, if there is
a surjection $r_X:X\rightarrow G^{(0)}$ and a map $(\gamma,x)\mapsto
\gamma\cdot x$ from $G*X:=\{(\gamma,x)\in G\times X:
s(\gamma)=r_X(x)\}$ to $X$ such that 
\begin{enumerate}
\item if $(\eta,x)\in G*X$ and $(\gamma,\eta)\in G^{(2)}$, then 
$(\gamma\eta,x),(\gamma,\eta\cdot x)\in G*X$ and 
\[
\gamma\cdot(\eta\cdot x) = \gamma\eta \cdot x,
\]
\item and $r_X(x) \cdot x = x$ for all $x\in X$. 
\end{enumerate}
Right actions and right $G$-spaces are defined similarly except we use
$s_X$ to denote the map from $X$ to $G^{(0)}$ and we define the action
on the set $X*G := \{(x,\gamma)\in X\times G: s_X(x) = r(\gamma)\}$.
\end{definition}

\begin{remark}
\index{range map $r$}
\index{source map $s$}
Given a left $G$-space we call $r_X(x)$ the range of $x$ and in a
right $G$-space $s_X(x)$ is the source of $x$.  
In order to avoid notational clutter we will almost always
drop the subscripts on $r_X$ and $s_X$.  When an action is not
specified to act on the right or the left we will assume that it acts
on the left.  (Unless it acts on the right, of course.)
\end{remark}

\begin{definition}
\label{def:19}
Let $G$ be a groupoid acting on both $X$ and $Y$.  A map
$\phi:X\rightarrow Y$ is {\em $G$-equivariant} if and only if $r_X(x) =
r_Y(\phi(x))$ and $\phi(\gamma\cdot x) = \gamma\cdot \phi(x)$ for all
$x\in X$ and $\gamma\in G_{r(x)}$.  
\end{definition}

\begin{remark}
Suppose $G$ is a groupoid, $X$ is a $G$-space, $H\subset G$ and 
$A\subset X$.  We
will use the notation 
\[
H\cdot A := \{\gamma\cdot x: \gamma\in H, x\in A, s(\gamma) = r(x)\}. 
\]
As in Remark \ref{rem:4}, it's important to realize that $H\cdot A$ may
be poorly behaved.  For instance if $s(H)\cap r(A) = \emptyset$ then
$H\cdot A$ is empty. 
\end{remark}

\begin{prop}
\label{prop:12}
Suppose $G$ is a groupoid and $X$ is a left $G$-space. 
\begin{enumerate}
\item Given $\gamma\in G$ and $x\in X$ such that $s(\gamma) = r(x)$ we
  have $r(\gamma\cdot x) = r(\gamma)$.  
\item Given $\gamma\in G$ and $x\in X$ such that $s(\gamma) = r(x)$ we
  have $\gamma\inv \cdot (\gamma\cdot x) = x$.  
\item Given $\gamma\in G$ and $x,y\in X$ such that $s(\gamma) = r(x)$
  and $y = \gamma\cdot x$ we have $\gamma\inv \cdot y = x$.  
\end{enumerate}
Similar statements hold if $X$ is a right $G$-space.  
\end{prop}

\begin{proof}
Part {\bf (a)}: If $s(\gamma) = r(x)$ then $(\gamma,x)\in G*X$.  Furthermore,
$(r(\gamma),\gamma)\in G^{(2)}$ so that by Definition \ref{def:10} we
have $(r(\gamma),\gamma\cdot x)\in G*X$.  However, this implies
\[
s(r(\gamma)) = r(\gamma) = r(\gamma\cdot x). 
\]

Part {\bf (b)}:  Given $(\gamma,x)\in G*X$ we have 
\begin{align*}
\gamma\inv\cdot(\gamma\cdot x) &= (\gamma\inv\gamma)\cdot x =
s(\gamma)\cdot x \\
&= r(x)\cdot x = x.
\end{align*}

Part {\bf (c)}: Given $(\gamma,x)\in G*X$ and $y= \gamma\cdot x$
observe that $r(y) = r(\gamma\cdot x) = r(\gamma)$.  Therefore we can
act on both sides by $\gamma\inv$ and use part (b) to obtain 
\[
\gamma\inv \cdot y = \gamma\inv\cdot(\gamma\cdot x) = x. 
\]

The corresponding statements for a right $G$ action are proved
similarly. 
\end{proof}

The appropriate topological notion of a $G$-space is almost exactly
what one would think.  

\begin{definition}
\label{def:11}
\index{G-space@$G$-space!continuous}
\index{G-space@$G$-space!strong/strongly continuous}
Suppose $G$ is a topological groupoid and $X$ is a
left $G$-space equipped with a topology.  We say that $G$ acts
continuously on the left of $X$, and call $X$ a {\em continuous left $G$-space},
if the maps $r_X:X\rightarrow G\unit$ and 
$(\gamma,x)\mapsto \gamma\cdot x$ from $G*X$ to $X$ are
continuous.  Continuous right $G$-spaces are defined similarly.  We
will call $X$ a {\em strongly continuous left $G$-space} (resp. {\em
  strongly continuous right
  $G$-space}) if $r_X$ (resp. $s_X$) is an open map.  Furthermore, we
will often refer to a strongly continuous $G$-space as a {\em
  strong $G$-space}.  
\end{definition}

\begin{remark}
\label{rem:3}
Definition \ref{def:11} is different from the usual definition found
in the literature.  It is standard to make the requirement that $r_X$
(or $s_X$) be open part of the definition of a $G$-space and to forgo
the notion of a ``strongly continuous $G$-space'' altogether.  
The author has chosen to introduce some new terminology and break from the
literature for (at least) three reasons.  
The first is that the material in Section 
\ref{sec:principal} does not require $r_X$ to be open.  The
second is that if $r_X$ is required to be open then technically
groupoids without Haar systems might not act ``continuously'' on
themselves.  Lastly, and most importantly, given a groupoid
dynamical system $(A,G,\alpha)$ as in Section \ref{sec:dynamical}, 
we will be interested in the induced action
of $G$ on $\Prim A$.  However, unless $A$ is a {\em continuous}
$C^*$-bundle the range map on $\Prim A$ will not be open.  Since
it will be necessary to deal with these actions at some point, we have
decided to make this distinction part of the definition.  
\end{remark}

\begin{example}
Suppose the locally compact group $H$ acts on a locally compact space
$X$.  Then, if we view $H$ as a groupoid with only one unit and let
the range map for $X$ map onto this point, it's easy
to see that $H$ is a strongly continuous groupoid action on $X$.  
\end{example}

\begin{example}
\label{ex:10}
The following is a particularly important example of a groupoid
action.  Suppose $G$ is a locally compact Hausdorff groupoid and $H$
is a closed subgroupoid.  Let $X=s\inv(H\unit)$ and give $X$ the
relative topology.  Since $H$ is closed in $G$, $H\unit$ is closed in
$G\unit$ so that $X$ is closed in $G$ and must be locally compact
Hausdorff.  Let $s:X\rightarrow H\unit$ be the restriction of the
source map on $G$ to $X$.  Given $(\xi,\eta)\in X*H$ we have $s(\xi) =
r(\eta)$ so that we can define $\xi\cdot \eta = \xi\eta$.  Observe
that $s(\xi\eta)=s(\eta)\in H\unit$ so that the action is well defined
on $X$.  Since this action is defined via the groupoid operation it is
straightforward to see that $X$ is a continuous right $G$-space.
Furthermore, if $G$ has open range and source then, because $X$ is a
saturated\footnote{Given a surjective map $f:X\rightarrow Y$ a set
  $A\subset X$ is saturated if $A=f\inv(f(A))$.}
closed set of $G$ with respect to $s$, it follows that the restriction
of $s$ to $X$ is open and in this case $X$ is a strongly continuous right
$G$-space.  We can also define an analogous left $G$-space in a
similar fashion. 
\end{example}

\begin{example}
\label{ex:11}
This is a special case of Example \ref{ex:10} but
deserves to be singled out.  Suppose $G$ is a locally compact
Hausdorff groupoid.  If we treat $G$ as a closed subgroupoid of itself
and use the construction from Example \ref{ex:10} we find that $X=G$
and that there is a continuous right action of $G$ on itself defined
by multiplication.  Furthermore, if $G$ has open range and source then clearly
this action is strongly continuous.  
\end{example}

\begin{example}
\label{ex:12}
Suppose $G$ is a locally compact Hausdorff groupoid.  Let
$r:G\unit\rightarrow G\unit$ be the identity map (or the restriction
of $r$ to $G\unit$, whichever you prefer).  Given
$(\gamma,u)\in G*G\unit$ define 
\[
\gamma\cdot u := r(\gamma) = \gamma u\gamma\inv.
\]
It is straightforward to check that this defines a groupoid action and
that with this action $G\unit$ is a strong $G$-space.  As usual, we can
let $G$ act on the right of $G\unit$ in a similar fashion. 
\end{example}

\begin{example}
Suppose $G$ is a locally compact Hausdorff groupoid and $S$ is the
stabilizer subgroupoid.  Let $p:S\rightarrow G\unit$ be the
restriction of the range map to $S$ (or the bundle map associated to
$S$, whichever you prefer).  Given $(\gamma,s)\in G*S$ define 
\[
\gamma\cdot s := \gamma s \gamma\inv.
\]
This action is well defined since $(\gamma,s)\in G*S$ implies
$s(\gamma) = p(s)$.  Furthermore, it's easy to see that $\gamma\cdot
\eta\cdot s = \gamma\eta \cdot s$ and that $p(s) \cdot s = s$.  It's
fairly clear that the action is continuous so that $S$ is a continuous
$G$-space.  Furthermore, if the stabilizers of $G$ vary continuously
then $p$ is open and $S$ is a strong $G$-space.  As usual, we can let
$G$ act on the right of $S$ in a similar fashion.  
\end{example}

As with group actions, the action of a groupoid on a space defines an
equivalence relation, and (fortunately) the quotient map associated to
this equivalence relation is open whenever the range map on $G$ is
open.   This is true even if $X$ is not a
strongly continuous $G$-space.

\begin{definition}
\index{orbit equivalence relation}
\index[not]{$X/G$, $G\backslash X$}
\index[not]{$G\cdot x$, $x\cdot G$}
\index[not]{$[x]$}
\label{def:12}
Suppose $X$ is a (left) $G$-space.  Define the {\em orbit equivalence
  relation} on $X$ determined by $G$ to be $x\sim y$ if
and only if there exists $\gamma\in G$ such that $\gamma\cdot x =y$.  
The quotient space with respect to this relation is denoted $X/G$, the
elements of $X/G$ are denoted $G\cdot x$, and
the canonical quotient map is (often) denoted by $q$.  When $X$ is a
continuous $G$-space we will give $X/G$ the quotient topology with
respect to $q$.  When $X$ is a right $G$-space the orbit equivalence
relation is defined similarly and we will use exactly the same
notation. 

In some cases $X$ will be both a left $G$-space and
a right $H$-space and in these situations we will denote the orbit
space with respect to the $G$ action by $G\backslash X$ and the orbit space with
respect to the $H$ action by $X/H$ and we will denote elements of the
orbit space by $G\cdot x$ and $x\cdot H$, respectively.   
It will occasionally be useful to denote the orbit $G\cdot x$ or
$x\cdot G$ by $[x]$ to conserve notation. 
\end{definition}

\begin{remark}
Since the orbit equivalence relation on $G\unit$ with respect to $G$
as defined in Definition \ref{def:7} is exactly the orbit equivalence
relation on $G\unit$ with respect to $G$ when we view $G\unit$ as a
$G$-space via Example \ref{ex:12}, there is no problem reusing the
orbit equivalence relation terminology. 
\end{remark}

\begin{prop}
\label{prop:13}
Let $G$ be a locally compact Hausdorff groupoid and 
suppose $X$ is a continuous $G$-space.  Then the orbit equivalence
relation defined in Definition \ref{def:12} is an equivalence
relation.  If the range and source maps for $G$ are open
then the canonical quotient map $q:X\rightarrow X/G$ is open and $X/G$
is locally compact. In particular, this is true if $G$ has a Haar
system. Furthermore, in this case if $X$ is second countable then $X/G$ is
second countable.  
\end{prop}

\begin{proof}
We will assume that $G$ is a left $G$-space, but the proof is entirely
analogous when $G$ acts on the right.  
Since $r(x)\cdot x = x$ it is clear that $\sim$ is reflexive.  Given
$x,y\in X$ if $x\sim y$ then there exists $\gamma\in G$ such that
$y = \gamma\cdot x$.  However it follows from Proposition
\ref{prop:12} that $\gamma\inv \cdot y = x$ and $y\sim x$.  Finally if
$x\sim y$ and $y\sim z$ then find $\gamma,\eta\in G$ such that $y =
\gamma\cdot x$ and $z = \eta\cdot y$.  Then $z = (\eta\gamma)\cdot x$
and $x\sim z$.  Thus $\sim$ is an equivalence relation.  

Now suppose $G$ has open range and source maps.  By Proposition
\ref{prop:7} this is true whenever $G$ has a Haar system.  Since $X/G$
has the quotient topology it suffices to see that $q\inv(q(O))=G\cdot
O$ is open whenever $O$ is. Suppose $x_i$ is a net in $X$ such that
$x_i\rightarrow \gamma\cdot x$ where $\gamma\in G$ and $x\in O$.  Let
$u_i = r(x_i)$, $u=r(\gamma\cdot x)=r(\gamma)$ and observe that 
$u_i\rightarrow u$.  Since the range on $G$ is an open map we can
pass to a subnet, reindex, and find $\gamma_i\in G$ such that
$\gamma_i \rightarrow \gamma$ and $r(\gamma_i) = u_i$ for all $i$.  
It follows then that $\gamma_i\inv
\rightarrow \gamma\inv$ and, since the group action is continuous
\[
\gamma_i\inv \cdot x_i \rightarrow \gamma\inv\cdot(\gamma\cdot x) = x.
\]
However, $O$ is open so eventually $\gamma_i\inv \cdot x_i \in O$.
Thus, eventually, $x_i = \gamma_i \cdot (\gamma_i\inv \cdot x_i) \in
G\cdot O$.  This suffices to show that $G\cdot O$, and hence $q$,
is open.  Finally, the open image of a locally compact space is
locally compact, since a basis of compact neighborhoods will map to a
basis of compact neighborhoods.  The same argument shows that if
$X$ is second countable then so is $X/G$.  
\end{proof}

Just as in the group case we can associate a groupoid to a groupoid
action. 

\begin{definition}
\index{transformation groupoid}
\index[not]{$G\ltimes X$}
Suppose $G$ is a locally compact Hausdorff groupoid which acts
continuously on the left of a locally compact Hausdorff space $X$.
The {\em transformation groupoid} associated to $G$ and $X$ is the
space $G\ltimes X = \{(\gamma,x)\in G\times X : r(\gamma) = r(x)\}$
with 
\[
G\ltimes X^{(2)} = \{((\gamma,x),(\eta,y))\in (G\ltimes X)\times
(G\ltimes X) : y = \gamma\inv \cdot x\}
\] 
and, when $((\gamma,x),(\eta,y))\in G\rtimes X^{(2)}$, the operations
\begin{align*}
(\gamma,x)(\eta,y) &= (\gamma\eta,x), &
(\gamma,x)\inv &= (\gamma\inv, \gamma\inv\cdot x).
\end{align*}
When $G$ acts on the right of $X$ the transformation groupoid is
defined in an analogous fashion and is denoted $X\rtimes G$.  
\end{definition}

This groupoid is generally well behaved, as we will see after we prove
the following utility lemma.  

\begin{lemma}
\label{lem:8}
Suppose that $X$ is a locally compact space and $C$ is a closed
subset of $X$.  
Given $f\in C_c(C)$ we can extend $f$ to a function in
$C_c(X)$.  
\end{lemma}

\begin{proof}
Now, if everything is second countable then $X$ is normal so
we can use the usual Tietze Extension Theorem.  If we want to work with the
nonseparable case we will need to make the following local argument. 
First, let $U$ be some neighborhood of $K=\supp
f$ in $X$ with compact closure.  
Since $C$ is closed in $X$ it follows that 
$C\cap \overline{U}$ is a closed set 
in $\overline{U}$ and is therefore compact.  We can now use the Tietze
Extension Theorem for compact sets \cite[Lemma 1.42]{tfb2} to find 
a function $g\in C_c(X)$ such that $g$ is equal to $f$ on
$C\cap\overline{U}$.  Observe that $U$ is an open
neighborhood of the compact set $K$ so that we may use
Urysohn's Lemma for compact sets \cite[Lemma 1.41]{tfb2} to find $h\in
C_c(X)$ such that  $h$ is one on $K$
and $h$ is zero off $U$.  
Consider $hg$.  Given $x\in K\subset C\cap\overline{U}$ we have 
$h(x) = 1$, and $g(x)=f(x)$.  Thus $hg = f$ on
$K$.  If $x\in
C\setminus K$ then $f(x)=0$ and there are two 
cases to consider.  In the first
case $x\in U$ so that $x\in C\cap \overline{U}$ and $g(x)=f(x)=0.$
Otherwise $x\not\in U$ so that $h(x)=0$.   In
either case $hg(x)=0$ so that $hg$ is an
extension of $f\in C_c(C)$ to $C_c(X)$.  
\end{proof}

Transformation groupoids are very similar to their group action
analogues.  For instance, they have a Haar system as long as $G$ does,
even if the action of $G$ is not strongly continuous.  

\begin{prop}
\label{prop:14}
\index{transformation groupoid!Haar system}
Suppose $G$ is a locally compact Hausdorff groupoid acting
continuously on a locally compact Hausdorff space $X$.  Then the
transformation groupoid $G\ltimes X$ is a locally compact Hausdorff
groupoid which is second countable if $G$ and $X$ are.  
The unit space $G\ltimes X\unit$ can be naturally
identified with $X$ and under this identification 
\begin{align*}
r(\gamma,x) &= x, &
s(\gamma,x) &= \gamma\inv \cdot x. 
\end{align*}
The range and source maps are open if the range and source maps of
$G$ are open.  Furthermore, if $G$ has a Haar system $\{\lambda^u\}$ 
then $\mu^x = \lambda^{r(x)}\times \delta_x$ is a
Haar system for $G\ltimes X$.  
\end{prop}

\begin{proof}
Suppose $((\gamma,x),(\eta,y)),((\eta,y),(\xi,z))\in G\ltimes
X^{(2)}$.  Then 
\[
(\gamma\eta)\inv\cdot x = \eta\inv \cdot (\gamma\inv \cdot x) =
\eta\inv \cdot y = z
\]
and we have assumed $\gamma\inv \cdot x = y$ so that
$((\gamma\eta,x),(\xi,z)),((\gamma,x),(\eta\xi,y))\in G\ltimes
X^{(2)}$.  Furthermore we clearly have
\[
(\gamma,x)((\eta,y)(\xi,z)) = (\gamma\eta\xi,x) =
((\gamma,x)(\eta,y))(\xi,z).
\]
Next we calculate that
\[
((\gamma,x)\inv)\inv = (\gamma\inv,\gamma\inv\cdot x)\inv 
= (\gamma,\gamma\cdot(\gamma\inv\cdot x)) = (\gamma, x).
\]
Finally suppose $(\gamma,x)\in G\ltimes X$.  Then it's easy to see 
$((\gamma\inv,\gamma\inv\cdot x),(\gamma,x))\in G\ltimes X^{(2)}$ and
if $((\gamma,x),(\eta,y))\in G\ltimes X^{(2)}$ then 
\begin{align*}
(\gamma\inv,\gamma\inv\cdot x)((\gamma,x)(\eta,y)) &= 
(\gamma\inv,\gamma\inv\cdot x)(\gamma\eta,x) \\
&= (\gamma\inv\gamma\eta,\gamma\inv\cdot x) \\
&= (\eta,y).
\end{align*}
Similarly $((\gamma,x)(\eta,y))(\eta,y)\inv=(\gamma,x)$.  Thus
$G\ltimes X$ is a groupoid.  We can now calculate
\begin{align*}
(\gamma\inv,\gamma\inv \cdot x)(\gamma,x) & = (s(\gamma),\gamma\inv\cdot
x) = (r(\gamma\inv\cdot x),\gamma\inv \cdot x)\\
(\gamma,x)(\gamma\inv,\gamma\inv\cdot x) & = (r(\gamma),x)
= (r(x),x).
\end{align*}
It follows that $G\ltimes X\unit = r(G\ltimes X) = \{(r(x),x): x\in
X\}$ and the identification with $X$ is obvious.  Furthermore, under
this identification, the previous calculations show that $r(\gamma,x) =
x$ and $s(\gamma,x) = \gamma\inv \cdot x$.  

Observe that, since the range maps for both $G$ and $X$ are
continuous, the set $G\ltimes X$ is closed in $G\times X$ and is
therefore locally compact Hausdorff.  Clearly $G\ltimes X$ is second
countable if $G$ and $X$ are.  Furthermore, since the topology
on $G\ltimes X$ is inherited from the product topology and because
the operations on $G$ and the action on $X$ are continuous, it's easy to
see that $G\ltimes X$ is a topological groupoid.  

Now we show that the range and source maps of $G\ltimes X$ are open if
the range and source maps of $G$ are open.  
Suppose $x_i\rightarrow x$ and $(\gamma,x)\in G\ltimes X$.  Since
$r(\gamma) = r(x)$ and $r(x_i)\rightarrow r(x)$ we can
pass to a subnet, reindex, and find $\gamma_i\in G$ such that
$\gamma_i\rightarrow \gamma$ and $r(\gamma_i)=r(x)$.  Since
$(\gamma_i,x_i)\in G\ltimes X$ and $(\gamma_i,x_i)\rightarrow
(\gamma,x)$ we see that $r$ is open.  Since $s$ is the composition of
$r$ with the inversion map $s$ must be open as well. 

Next, suppose $\{\lambda^u\}$ is a Haar system for $G$ and for all
$x\in X$ define
$\mu^x = \lambda^{r(x)}\times\delta_x$ where $\delta_x$ is the Dirac
delta measure at $x$.  Then clearly $\mu^x$ is a
non-negative Radon measure and it's fairly easy to see that 
$\supp \mu^x = \supp\lambda^{r(x)}\times\supp \delta_x =
G^{r(x)}\times \{x\}$.  However, it's also clear that $G\ltimes X^x =
G^{r(x)}\times \{x\}$.  Thus condition (a) of Definition \ref{def:6}
is satisfied.  Now given $f\in C_c(G\ltimes X)$ we have 
\[
\int_{G\ltimes X} f(\eta,y) d\mu^x(\eta,y) = \int_G f(\eta,x)
d\lambda^{r(x)}(\eta)
\]
and we would like to show that the function 
\begin{equation}
\label{eq:1}
x\mapsto \int_G f(\eta,x)d\lambda^{r(x)}(\eta)
\end{equation}
is continuous.

Given $f\in C_c(G\ltimes X)$ we can extend $f$
to a function in $C_c(G\times X)$, also denoted $f$, using Lemma
\ref{lem:8}.  
For $h\in
C_0(G)$ and $g\in C_0(X)$ define $h\otimes g$ by
$h\otimes g(\gamma,x) = h(\gamma)g(x)$.  It follows
from \cite[Corollary B.17]{tfb} that sums of functions of the form $h\otimes
g$ are dense in $C_0(G\times X)$.  Therefore, we
can find sets $\{h_i^j\}\subset C_0(G)$ and $\{g_i^j\}\subset C_0(X)$ such that 
$k_i = \sum_j h_i^j\otimes g_i^j$ is a net and 
$k_i \rightarrow f$ uniformly.  Let $L_1$ be the
projection of $\supp f$ to $G$ and let $g\in C_c(G)$ be one on $L_1$ and
zero off some compact neighborhood of $L_1$.  Similarly let $L_2$ be
the projection of $\supp f$ to $X$ and let $h\in C_c(X)$ be one on
$L_2$ and zero off some compact neighborhood of $L_2$.  By replacing
$g_i^j$ with $gg_i^j$ and $h_i^j$ with $hh_i^j$
we can assume, without loss of generality, that there exists a
compact set $K$ such that $\supp(k_i)\subset K$ for all $i$.  
Observe that for each $i$ and $j$ the function 
\[
x \mapsto g_i^j(x)\int_G h_i^j(\eta)d\lambda^{r(x)}(\eta) = \int_G
h_i^j\otimes g_i^j(\eta,x) d\lambda^{r(x)}(\eta)
\]
is continuous since it's built from continuous functions using
composition and multiplication.  Hence 
\[
x \mapsto \int_G k_i(\eta,x) d\lambda^{r(x)}(\eta)
\]
is the sum of continuous functions and is continuous. 
Next, let $K$ be the compact set
given above and $K'$ be
its projection onto $G$.  Since $K'$ is compact it follows from Lemma
\ref{lem:10} that  
$\{\lambda^u(K')\}$ is bounded by some number $M$.  Thus given $y\in X$
\begin{align*}
\left| \int_G f(\eta,y) d\lambda^{r(y)}(\eta) - \int_G
  k_i(\eta,y)d\lambda^{r(y)}(\eta) \right| &= 
\left| \int_G (f-k_i)(\eta,y) d\lambda^{r(y)}(\eta)\right| \\
&\leq \int_G |(f-k_i)(\eta,y)| d\lambda^{r(y)}(\eta) \\
&\leq M \|f-k_i\|_\infty.
\end{align*}
Now suppose $x_j\rightarrow x$ is a net converging in $X$ and
$\epsilon > 0$.  Choose $i$ large enough so that $\|f-k_i\|_\infty <
\epsilon/4M$ and $J$ such that for all $j\geq J$ we have 
\[
\left| \int_G k_i(\eta,x_j)d\lambda^{r(x_j)}(\eta) - 
\int_G k_i(\eta,x)d\lambda^{r(x)}(\eta)\right| < \epsilon/2
\]
Then for all $j\geq J$ we have 
\begin{align*}
\Biggl| \int_G f(\eta,x_j)  d\lambda^{r(x_j)}(\eta) &- \int_G f(\eta,x)
  d\lambda^{r(x)}(\eta) \Biggr| \\
\leq & \Biggl| \int_G 
f(\eta,x_j) d\lambda^{r(x_j)}(\eta) - \int_G
  k_i(\eta,x_j)d\lambda^{r(x_j)}(\eta) \Biggr| \\
&+\Biggl| \int_G k_i(\eta,x_j)d\lambda^{r(x_j)}(\eta) - 
 \int_G k_i(\eta,x)d\lambda^{r(x)}(\eta)\Biggr| \\
&+\Biggl| \int_G f(\eta,x) d\lambda^{r(x)}(\eta) - \int_G
  k_i(\eta,x)d\lambda^{r(x)}(\eta) \Biggr| \\
\leq & 2M\|f-k_i\|_\infty + \epsilon/2 < \epsilon
\end{align*}
This proves that the function in \eqref{eq:1} is continuous.  The
support of this function must be contained in $r(\supp f)$ and this
implies that condition (b) of Definition \ref{def:6} is satisfied.  

Finally, suppose
$(\gamma,x)\in G\ltimes X$ and $f\in C_c(G\ltimes X)$.  Then 
\begin{align*}
\int_{G\ltimes X} f((\gamma,x)(\eta,y)) d\mu^{\gamma\inv\cdot x}(\eta,y) &=
\int_G f((\gamma,x)(\eta,\gamma\inv\cdot x)) \lambda^{r(\gamma\inv\cdot
  x)}(\eta) \\
&= \int_G f(\gamma\eta,x) \lambda^{s(\gamma)}(\eta) \\
&= \int_G f(\eta,x) \lambda^{r(\gamma)}(\eta) \\
&= \int_{G\ltimes X} f(\eta,y) \mu^x(\eta,y)
\end{align*}
This proves the left invariance condition and that $\{\mu^x\}$ is a
Haar system for $G\ltimes X$.  
\end{proof}

\begin{remark}
\index{transformation group groupoid}
When we view a group action of $H$ on $X$ as a groupoid action we can
use Proposition \ref{prop:14} to form the transformation groupoid
$H\ltimes X$.  In this case $H\ltimes X$ is exactly the transformation
group groupoid from Example \ref{ex:4} and the Haar system on
$H\ltimes X$ is the one constructed in Example \ref{ex:5}.  Therefore
Proposition \ref{prop:14} justifies the statements made in those two
examples and we can use the phrase ``transformation groupoid'' without
causing too much confusion.
\end{remark}

\begin{remark}
In \cite{coords} the groupoid associated to a groupoid action is
defined to be $G*X=\{(\gamma,x)\in G\times X: s(\gamma)=r(x)\}$ 
with the operations 
\begin{align*}
(\gamma,\eta\cdot x)(\eta,x) &= (\gamma\eta, x), &
(\gamma,x)\inv &= (\gamma\inv, \gamma\cdot x).
\end{align*}
It's easy to see that the map $\phi:G\ltimes X \rightarrow G*X$ defined
by $\phi(\gamma,x) = (\gamma,\gamma\inv\cdot x)$ is an isomorphism
between these two groupoids.  The principal difference between them 
is that given
$(\gamma,x) \in G\ltimes X$ we have $r(\gamma,x) = x$ and $s(\gamma,x)
= \gamma\inv\cdot x$ and for $(\gamma,x)\in G*X$ we have
$r(\gamma,x) = \gamma\cdot x$ and $s(\gamma,x) = x$.  The reason that
we are using the former groupoid is that it interacts nicely with
crossed products.  
\end{remark}

Statements about groupoid actions often translate into similar statements about
the transformation groupoid.  For instance

\begin{definition}
\index{groupoid!transitive}
Suppose $G$ is a groupoid.  We say $G$ is {\em transitive} if given $u,v\in
G\unit$ there exists $\gamma\in G$ such that $s(\gamma) = u$ and
$r(\gamma) = v$.  
\end{definition}

\begin{remark}
Note that a groupoid is transitive if and only if its associated orbit
groupoid is trivial.  Also, in this case it is easy to see, using
conjugation by elements in $G$, that all the stabilizer subgroups are
isomorphic.
\end{remark}

\begin{definition}
\index{transitive action}
\index{orbit transitive action}
Suppose $G$ is a groupoid and $X$ is a $G$-space.  We say $X$ is {\em
  transitive} if given $x,y\in X$ there exists $\gamma\in G$ such that
$\gamma\cdot x =y$.  We say $X$ is {\em orbit transitive} if given
$x,y\in X$ such that $r(x)$ is orbit equivalent to $r(y)$ in $G\unit$
then there exists $\gamma\in G$ such that $\gamma\cdot x =y$. 
\end{definition}

\begin{remark}
Using Proposition \ref{prop:12}, it is straightforward to show 
that if $X$ is a transitive
$G$-space then $G$ must be a transitive groupoid.  If a non-transitive
groupoid acts on $X$ then
the most that one could hope for is that $X$ is orbit transitive.  
\end{remark}

\begin{prop}
Suppose $G$ is a groupoid, $X$ is a $G$-space and $G\ltimes X$ is the
associated transformation groupoid.  Then $X$ is transitive if and
only if $G\ltimes X$ is.  Furthermore $X$ is orbit transitive if and
only if the range map on $X$ factors to a bijection from $X/G$ onto
$G\unit / G$.  
\end{prop}
\begin{proof}

Well $X$ is transitive if and only if given $x,y\in X$ we have
$\gamma\in G$ such that $x = \gamma\cdot y$.  However this occurs if
and only if there exits $\gamma\in G$ such that 
\[
y = \gamma\inv\cdot x = s(\gamma,x),\quad\text{and}\quad x = r(\gamma,x).
\]
Finishing the chain, this is equivalent to requiring that $G\ltimes X$ be
transitive.  

Next, $r:X\rightarrow G\unit$ is always equivariant by Proposition
\ref{prop:12} and as such factors to a surjection from $X/G$ onto
$G\unit/G$.  If $X$ is orbit transitive then $G\cdot r(x) = G\cdot
r(y)$ implies that $r(x)$ and $r(y)$ are orbit equivalent so that
there exits $\gamma\in G$ such that $\gamma\cdot x =y$.  Thus $G\cdot
x = G\cdot y$ and $r$ is injective.  The reverse direction is just as
simple. 
\end{proof}

One of the important things about transformation groupoids is that
the original groupoid action can be recovered from the action of the
transformation groupoid on its unit space.  What this will allow us to
do is extend theorems about groupoids to theorems about groupoid
actions with relatively little effort.  The following proposition gives
some indication of how this works because it will allow us to
translate statements concerning the stabilizers and the orbit space of
the action of $G$ on $X$ into statements concerning the stabilizers
and orbit space of $G\ltimes X$.  

\begin{definition}
\index{stabilizer subgroup}
\index{stabilizer subgroupoid!continuously varying}
Suppose $G$ is a groupoid acting on a set $X$.  Given $x\in X$ the
{\em stabilizer subgroup} of $G$ at $x$ is
\[
G_x := \{ \gamma\in G : \gamma\cdot x = x\}
\]
If $G$ is a locally compact Hausdorff groupoid and $X$ is a continuous
$G$-space then we say the stabilizers {\em vary continuously} in $G$
if $x_i\rightarrow x$ in $X$ implies $G_{x_i}\rightarrow G_x$ in $G$ with
respect to the Fell topology. 
\end{definition}

\begin{prop}
\label{prop:15}
Suppose $G$ is a locally compact groupoid which acts continuously on
the locally compact space $X$.  Then $G_x$ is a closed subgroup of the 
stabilizer subgroup $S_{r(x)}$ for all $x\in X$.  
Furthermore, $G_x$ is naturally isomorphic to the stabilizer
subgroup of $G\ltimes X$ at $x$ and the stabilizers $G_x$ vary
continuously in $G$ if and only if
${G\ltimes X}$ has continuously varying stabilizers.  
Finally, the orbit space $X/G$ is naturally
homeomorphic to $G\ltimes X\unit/G\ltimes X$.
\end{prop}

\begin{proof}
If $\gamma \in G_x$ then $s(\gamma) = r(x)$ and 
$r(\gamma) = r(\gamma\cdot x) = r(x)$.  Thus $G_x \subseteq
S_{r(x)}$.  Furthermore, it's straightforward to show that, because the action is
continuous, $G_x$ is closed.  If $s,t\in G_x$ then $(st)\cdot
x = s\cdot (t\cdot x) = x$ so that $st\in G_x$.  Lastly if $s\in G_x$
then $s\cdot x = x$ so that, using Proposition \ref{prop:12}, $x = s\inv
\cdot x$.  Thus $G_x$ is a closed subgroup of $S_{r(x)}$.  Next, let
$T_x$ the be the stabilizer subgroup of $G\ltimes X$ at $x$ and 
define $\phi:G_x \rightarrow T_x$ by $\phi(t) = (t,x)$.  It
is clear that $r(t,x) = s(t,x) = x$ so that $\phi$ is well defined.
What's more, $\phi$ is a homomorphism because 
\[
\phi(st) = (st,x) = (s,x)(t,x) = \phi(s)\phi(t).
\]
Moving on, it's easy to see that $\phi$ is continuous.  If we let $\psi :
T_x \rightarrow G_x$ be defined by $\psi(t,x) = t$ then, since 
\[
t\inv \cdot x  = s(t,x) = r(t,x) = x
\]
we can conclude that $t\inv$, and hence $t$, are elements of $G_x$,
making $\psi$ well defined.  Furthermore $\psi$ is clearly a
continuous inverse to $\psi$, making $\psi$ an isomorphism of locally
compact groups.  

The following will make heavy use of Proposition \ref{prop:8}. 
Suppose the stabilizers $G_x$ vary continuously with respect to the
Fell topology and suppose $x_i\rightarrow x$.  Next, let
$(s_i,x_i)\in T_{x_i}$ and suppose $(s_i,x_i)\rightarrow (s,x)$. Then,
because the range and source maps are continuous and $x_i\rightarrow
x$, we have $(s,x)\in T_x$.  Now suppose $(s,x)\in T_x$.  We know
$s\in G_x$ so that we can pass to a subnet, relabel, and find $s_i\in
G_{x_i}$ such that $s_i\rightarrow s$.  Thus $(s_i,x_i)\in T_{x_i}$
and $(s_i,x_i)\rightarrow (s,x)$.  It follows that $T_{x_i}\rightarrow
T_x$ in the Fell topology.  The opposite directly is proved in an
entirely similar fashion.  

Next, recall that we identify the unit space of $G\ltimes X$ with
$X$ via the map $(r(x),x)\leftrightarrow x$.  In order to show that $X/G$ and $X/G\ltimes X$ are naturally
isomorphic it suffices to show that the actions induce the same orbit
equivalence relation.  However $x\sim y$ with respect to the action of
$G$ if and only if there exists $\gamma\in G$ such that $\gamma\cdot x
= y$.  This is true if and only if $x = \gamma\inv \cdot y$ which is
true if and only if there exists $(\gamma,x)\in G\ltimes X$ such that $y =
s(\gamma,x)$ and $x=r(\gamma,x)$.  Continuing the string of
implications, this is true if and only if $x \sim y$ with respect to
the action of $G\ltimes X$.  Since the actions induce the same orbit
equivalence relation, the quotient of $X$ with respect to each action
is the same.
\end{proof}

\begin{remark}
\label{rem:2}
\index{stabilizer subgroupoid}
Proposition \ref{prop:15} helps answer another question.  Given a
continuous $G$-space $X$ we would like to form a bundle out of the
stabilizer groups $G_x$.  Since
Proposition \ref{prop:15} tells us that the stabilizers of the action
are exactly the stabilizers of the transformation groupoid, it is
natural to bundle the $G_x$ together inside $G\ltimes X$.  In other
words, the stabilizer group bundle $T$ for the action of $G$ on $X$ is
nothing more than the stabilizer subgroupoid of $G\ltimes X$.  Of
course, once you work out all the definitions this just boils down to
defining the stabilizer bundle to be the set
\[
T:= \{(s,x)\in G\times X:s\in G_x\}.
\]
\end{remark}

\subsection{Groupoid Equivalence}
\label{sec:equivalence}
One very important application of groupoid actions is the the notion
of groupoid equivalence.  This is a fundamental idea and will be a
key tool in Section \ref{sec:indreps}.  The following material is
taken, and expanded, from \cite{groupoidequiv}.

\begin{definition}
\label{def:13}
\index{G-space@$G$-space!principal}
Suppose $G$ is a groupoid and $X$ is a $G$-space.  We say the action
of $G$ on $X$ is {\em free} if $\gamma\cdot x = x$ implies $\gamma\in
G\unit$.  If $G$ is a locally compact Hausdorff groupoid and $X$ is a
continuous $G$-space then we say the action is {\em proper} if the map
from $G*X$ to $X\times X$ given by $(\gamma,x)\mapsto (\gamma\cdot x,
x)$ is proper.\footnote{A continuous map is proper if the inverse
  image of every compact set is compact.}  The action is {\em
  principal}, and $X$ is called a {\em principal $G$-space}, if it 
is both free and proper.
\end{definition}

The following proposition gives a useful characterization of 
proper actions that helps explain why we are interested in them. The
proof is lifted from \cite[Lemma 3.24]{tfb2}.

\begin{prop}
\label{prop:18}
Suppose $G$ is a locally compact Hausdorff groupoid acting
continuously on a
locally compact Hausdorff space $X$.  The action is
proper if and only if given nets $\{x_i\}_{i\in I}\in X$ and
$\{\gamma_i\}_{i\in I}\in G$ 
such that $x_i\rightarrow x$ and 
$\gamma_i\cdot x_i\rightarrow y$ then $\{\gamma_i\}$ has a convergent
subnet.  
\end{prop}

\begin{proof}
Suppose the action is proper so that the map $\phi:G*X\rightarrow
X\times X$ given by $\phi(\gamma,x)=(\gamma\cdot x, x)$ is a proper
map.  Given elements $x_i,x,y\in X$ and $\gamma_i,\gamma\in G$ as in the 
statement of the proposition, let $K$ be a
compact neighborhood of $x$ and $y$.  We eventually have
$x_i,\gamma_i\cdot x_i\in K$ so that eventually 
$(\gamma_i\cdot x_i,x_i)\in K\times K$.  Hence, eventually,
$(\gamma_i,x_i)\in \phi\inv(K\times K)$.  Since $\phi$ is proper,
$\phi\inv (K\times K)$ is compact and $\{\gamma_i\}$ must have a
convergent subnet.  

Now we will prove the reverse direction.  Suppose $K$ is a compact
subset of $X\times X$ and $\{(\gamma_i,x_i)\}$ is a net in
$\phi\inv(K)$.  Then $\{(\gamma_i\cdot x_i,x_i)\}$ is a net in
$K$ so that we can pass to a subnet, relabel, and find $(y,x)\in
K$ such that $\gamma_i\cdot x_i\rightarrow y$ and
$x_i\rightarrow x$.  However, we can now use our hypothesis to pass to
a subnet, relabel, and find $\gamma$ such that $\gamma_i\rightarrow
\gamma$.  Since the action is continuous $\gamma_i\cdot x_i \rightarrow
\gamma\cdot x$ and we have $\gamma\cdot x = y$.  
Finally, $(\gamma_i,x_i)\rightarrow (\gamma,x)$
and $\phi(\gamma,x) = (y,x)$ so that $\{(\gamma_i,x_i)\}$
has a convergent subnet in $\phi\inv(K)$ and we are done. 
\end{proof}

One good thing about proper actions of well behaved groupoids 
is that they have nice orbit spaces. 

\begin{prop}
\label{prop:17}
If $G$ is a locally compact Hausdorff groupoid with open range and
source maps and $X$ is a proper
$G$-space then the orbit space $X/G$ is locally compact Hausdorff.  
\end{prop}
\begin{proof}
We proved in Proposition \ref{prop:13} that the quotient map $q$ is
open and that $X/G$ is locally compact.  Suppose we have a net $G\cdot
x_i$ in $X/G$ such that $G\cdot x_i\rightarrow G\cdot x$ and $G\cdot
x_i\rightarrow G\cdot y$.  It will suffice to show $G\cdot x = G\cdot
y$.  Using the fact that $q$ is open, we can pass to a subnet,
relabel, and choose new representatives $x_i$ so that $x_i\rightarrow
x$.  Then we pass to another subnet, relabel, and this time find
$\gamma_i\in G$ such that $\gamma_i\cdot x_i \rightarrow y$.  Since
the action is proper we can use Proposition \ref{prop:18} to pass to
yet another subnet, relabel, and find $\gamma$ such that
$\gamma_i\rightarrow \gamma$.  It now follows from the continuity of
the action that $\gamma_i\cdot x_i\rightarrow \gamma\cdot x$.  Since
$X$ is Hausdorff, $\gamma\cdot x = y$ and $G\cdot x = G\cdot y$.  
\end{proof}

\begin{example}
The most basic example of a principal $G$-space is the action from
Example \ref{ex:11} where we let $G$ act on itself by multiplication.
It is easy enough to show that this action is free and proper.  
\end{example}

\begin{definition}
\label{def:14}
\index{G,H-equivalence@$(G,H)$-equivalence}
Suppose $G$ and $H$ are locally compact Hausdorff groupoids which have
open range and source maps.  We say
that a locally compact space $Z$ is a {\em $(G,H)$-equivalence} if 
\begin{enumerate}
\item $Z$ is a left {\em strong} principal $G$-space,
\item $Z$ is a right {\em strong} principal $H$-space, 
\item the $G$ and $H$ actions commute,
\item the map $s_X$ induces a bijection of $Z/H$ onto $G\unit$, and
\item the map $r_X$ induces a bijection of $G\backslash Z$ onto
  $H\unit$. 
\end{enumerate}
\end{definition}

\begin{remark}
Recall that a $G$-space $X$ is ``strong'' if the structure map from
$X$ to $G\unit$ is open.  This is a necessary condition when dealing
with $(G,H)$-equivalences.  Colloquially, when dealing with groupoid
equivalence every range and source map in sight has to be open.  
\end{remark}

\begin{remark}
It is remarked in \cite{groupoidequiv} that $(G,H)$-equivalence
induces an equivalence relation on locally compact groupoids.  It is
easy to see that the left and right action of $G$ on $G$ makes $G$ a
$(G,G)$-equivalence.  Given a $(G,H)$-equivalence one can just swap
the left and right actions to form a $(H,G)$-equivalence.  Finally,
given a $(G,H)$-equivalence $Z$ and a $(H,K)$-equivalence $Y$ we
define $Z*_H Y$ to be the quotient of $Z*Y=\{(z,y)\in Z\times Y:
s(z)=r(y)\}$ where we identify $(z,y)$ with $(z\cdot \gamma,
\gamma\cdot y)$ 
for all $z\in Z$, $y\in Y$ and $\gamma\in H$.  It's not hard to see
that $G$ and $K$ act naturally on $Z*_HY$ and that $Z*_HY$ is a
$(G,K)$-equivalence.  
\end{remark}

The definition of a $(G,H)$-equivalence is a little complicated so one
might expect that they are relatively rare.  In fact, any strong
principal $G$-space gives rise to an equivalence. 

\begin{definition}
\label{def:15}
\index{imprimitivity groupoid}
\index[not]{$X^H$}
Suppose $H$ is a locally compact Hausdorff groupoid with open range
and source maps and $X$ is a
strong principal right $H$-space.  Let $H$ act on $X*X = \{(x,y)\in
X\times X: s(x)=s(y)\}$ via the diagonal action 
\[
(x,y)\cdot\gamma := (x\cdot \gamma,y\cdot \gamma) 
\]
and define $X^H := (X*X)/H$.  Denoting the image of $(x,y)$ in $X^H$ as
$[x,y]$ we let $(X^H)^{(2)} = \{([x,y],[w,z])\in X^H\times X^H :
y=w\}$ and we define groupoid operations on $X^H$ via
\begin{align*}
[x,y][y,z]&=[x,z], &
[x,y]\inv &= [y,x].
\end{align*}
Then equipped with these operations, $X^H$ is called the {\em imprimitivity
  groupoid} associated to $X$ and $H$. The imprimitivity groupoid
associated to a left action is defined analogously.  
\end{definition}

\begin{remark}
Propositions \ref{prop:16}, \ref{prop:19} and \ref{prop:20} all have
corresponding statements for imprimitivity groupoids built from strong
principal {\em left} $G$-spaces and the proofs are all exactly the
same.  We will endeavor to always build
imprimitivity groupoids from right actions, but there isn't any real
difference one way or the other.  
\end{remark}

Of course, we made a number of claims in the definition that need to
be verified. 

\begin{prop}
\label{prop:16}
Let $H$ be a locally compact Hausdorff groupoid with open range and
source maps and $X$ a strong
principal right $H$-space.  Then the imprimitivity groupoid $X^H$ is a
locally compact Hausdorff groupoid which is second countable if $X$
is.  The unit space of $X^H$ can be identified with $X/H$
and under this identification $r([x,y])=[x]$ and
$s([x,y])=[y]$.  Furthermore, the range and source maps of
$X^H$ are open.  
\end{prop}

\begin{proof}
First we show that $X^H$ is a locally compact Hausdorff space.  It is
straightforward to show 
that the diagonal action of $H$ on $X*X=\{(x,y):s(x)=s(y)\}$ is a
continuous groupoid action.  Using Proposition \ref{prop:13} we
conclude that the quotient map $r:X*X\rightarrow X^H$ is open, 
$X^H=X*X/H$ is locally compact, and $X^H$ is
second countable if $X$, and hence $X*X$, is second countable.  
Next we will show that
the action of $H$ on $X*X$ is proper.  Suppose $\{(x_i,y_i)\}$ is a net in
$X*X$ and $\{\gamma_i\}$ is a net in $H$ such that
$(x_i,y_i)\rightarrow (x,y)$ and $(x_i\cdot \gamma_i, y_i\cdot
\gamma_i) \rightarrow (w,z)$.  Then, using Proposition \ref{prop:18}, and the
fact that the action of $H$ on $X$ is proper, we can find a convergent
subnet of $\gamma_i$.  Thus the action of $H$ on $X*X$ is proper and
$X^H$ must be a Hausdorff space by Proposition \ref{prop:17}.  

We must spend some time showing that the operations in Definition 
\ref{def:15} are well defined.  For instance, if $[x,y]=[x',y']$ and
$[y,z]=[y',z']$ then there exists $\gamma,\eta\in H$ such that
$y=y'\cdot\gamma$ and $y=y'\cdot\eta$.  However, it's easy to see
that, because the action of $H$ on $X$ is free, we must have $\gamma =
\eta$.  Hence $(x,z) = (x',z')\cdot\gamma$ and multiplication is well
defined.  We can also do a similar calculation to show that the
inverse is well defined.  At this point it is straightforward to show
that these operations make
$X^H$ into a groupoid. In order to see that the actions are continuous
one must use the fact that the quotient map $Q:X*X\rightarrow X*X/H$ 
is open.  For example,
if $[x_i,y_i]\rightarrow [x,y]$ then by passing to a subnet,
relabeling, and possibly choosing new representatives $x_i$ and
$y_i$, we can assume $x_i\rightarrow x$ and $y_i\rightarrow y$.
However, it's now clear that $[y_i,x_i]\rightarrow[y,x]$ so that the
inverse is continuous.  Similar considerations show that the
multiplication is also continuous and that $X^H$ is a topological
groupoid.  

Now $r([x,y])=[x,y][y,x]=[x,x]$ and similarly $s([x,y])=[y,y]$.  The
map $\phi:X\rightarrow X*X$ such that $\phi(x) = (x,x)$ is a
homeomorphism of $X$ onto the diagonal in $X*X$.  Since 
$\phi(x\cdot\gamma) = \phi(x)\cdot\gamma$, it is straightforward to
show that this homeomorphism factors to a homeomorphism $\bar{\phi}$
from $X/H$ onto $\{[x,x]\in X^H: x\in X\}$.  Under this identification
we clearly have $r([x,y])=[x]$ and $s([x,y])=[x]$.  Finally, we show
that the range and source maps are open.  Suppose $[x_i]\rightarrow
[x]$ and $r([x,y])=[x]$.  Since $q:X\rightarrow X/H$ is open we can
pass to a subnet, relabel, and choose new representatives so that
$x_i\rightarrow x$.  Now $s(x_i)\rightarrow s(x)$ and $s(y)=s(x)$.
Using the fact that the action of $H$ on $X$ is strongly continuous 
we can pass to
another subnet, reindex, and find $y_i$ in $X$ such that
$y_i\rightarrow y$ and $s(y_i)=s(x_i)$.  Thus $(x_i,y_i)\in X*X$ and
clearly $(x_i,y_i)\rightarrow (x,y)$.  It follows that
$[x_i,y_i]\rightarrow [x,y]$ and that $r([x_i,y_i])=[x_i]$.  Hence the
range map is open on $X^H$.  We could run through a similar argument
to show that the source map is open, or we could observe that the
source map is equal to the range map composed with the inverse map and
that the inverse map is a homeomorphism.  
\end{proof}

The reason we care about the imprimitivity groupoid is that it turns
out to be naturally equivalent to $H$. 

\begin{prop}
\label{prop:19}
Suppose $H$ is a locally compact Hausdorff groupoid with open range
and source maps and $X$ is a
strong principal right $H$-space.  Then the imprimitivity groupoid $X^H$ has
a strong principal action on the left of $X$.  The range map
$r_X:X\rightarrow X/H$ is the quotient map and the action is defined by 
\[
[x,y]\cdot z = x\cdot\gamma
\]
where $\gamma$ is the unique element of $H$ such that $z =
y\cdot\gamma$.  
Furthermore, with these two actions $X$ is a $(X^H,H)$-equivalence. 
\end{prop}

\begin{proof}
First, observe that $r_X$ is open by Proposition
\ref{prop:13}.  If $s([x,y]) =
[x]= [z] = r(z)$ then there exists $\gamma\in H$ such that
$y\cdot\gamma = z$.  Furthermore if $\eta\in G$ such that $y\cdot\eta
=z$ then $y\cdot \gamma\eta\inv = y$ so that $\gamma\eta\inv \in
G\unit$ and $\gamma = \eta$.  Thus $\gamma$ is unique and under our
definition $[x,y]\cdot z = x\cdot\gamma$.  Now, if
$[x',y']=[x,y]$ then there exists $\zeta \in G$ such that 
$(x,y) = (x',y')\cdot\zeta$.  It follows that $y' \cdot \zeta\gamma =
y\cdot \gamma = z$ so that $[x',y']\cdot z = x'\cdot \zeta\gamma =
x\cdot\gamma$.  This shows that the action of $X^H$ on $X$ is well
defined.  

Next we will show it is a groupoid action.  Suppose $z\in X$, and
$[x,y],[y,w]\in X^H$ such that $[w]=[z]$.  Let $\gamma$ be
the unique element such that $w\cdot\gamma = z$.   Then $[y,w]\cdot z
=  y\cdot \gamma$.  Clearly $[y\cdot\gamma] = [y]$ so that
$[x,y]$ acts on $y\cdot\gamma$.  The unique element of $G$ we are
looking for is again $\gamma$ so that 
\[
[x,y]\cdot([y,w]\cdot z) = [x,y]\cdot (y\cdot\gamma) = x\cdot\gamma. 
\]
On the other hand we clearly have $[x,w]\cdot z= x\cdot \gamma$ so
that condition (a) of Definition \ref{def:10} is satisfied.
Since $x\cdot s(x) = x$ we have $[x,x]\cdot x= x$ for all
$x\in X$ and condition (b) is satisfied as well.  Now suppose
$z_i\rightarrow z$ in $X$ and $[x_i,y_i]\rightarrow[x,y]$ in $X^H$
such that $[y_i]=[z_i]$ for all $i$.
First, pass to a subnet.  We will show that there is a sub-subnet such
that $[x_i,y_i]\cdot z_i \rightarrow [x,y]\cdot z$.  Now, 
pass to a subnet, relabel, and choose new representatives so that
$x_i\rightarrow x$ and $y_i\rightarrow y$.  
For each $i$ let $\gamma_i$ be such that $y_i\cdot\gamma_i = z_i$
and let $\gamma$ be such that $y\cdot\gamma = z$.  Using the fact
that the action of $G$ is proper we can conclude that by passing to a
subnet we can find $\eta$ such that $\gamma_i\rightarrow \eta$.
However, we then have $y\cdot\eta = z$ so that $\eta = \gamma$ by
freeness.  We can use this trick to show that every subnet of
$\gamma_i$ has a subnet which converges to $\gamma$.  This implies
$\gamma_i\rightarrow \gamma$.  Hence $x_i\cdot\gamma_i \rightarrow
x\cdot\gamma$, but this is exactly what we needed to show.  

We have shown that $X$ is a strong left $X^H$-space.  We will now
show that it is principal.  Suppose $[x,y]\cdot z = z$.  Let $\gamma$
be such that $y\cdot\gamma = z$.  Since $[x,y]\cdot z = x\cdot\gamma = z$
we can conclude that $y = z\cdot\gamma\inv = x$ and $[x,y]\in
(X^H)\unit$.  Thus the action is free. 
Now suppose we have $z_i,w,z\in X$ and $[x_i,y_i]\in
X^H$ such that $z_i\rightarrow z$ and $[x_i,y_i]\cdot z_i \rightarrow
w$.  Choose $\gamma_i\in G$ so that $y_i\cdot\gamma_i = z_i$ and
$[x_i,y_i]\cdot z_i = x_i\cdot \gamma_i$.  Now we have $x_i \cdot
\gamma_i \rightarrow w$ and $y_i\cdot \gamma_i\rightarrow z$.  However,
this implies $(x_i,y_i)\cdot \gamma_i\rightarrow (w,z)$ in $X*X$ and
that $[x_i,y_i]\rightarrow [w,z]$ in $X^H$.  Thus the action is proper
and we are done.  

Now we will show that $X$ is a $(X^H,H)$-equivalence.  We have already
shown, or assumed, that $H$ and $X^H$ have open range and source maps,
that $X$ is a strong principal right $H$-space, and that $X$ is a 
strong principal left $X^H$-space.  Furthermore we have already seen
that the range map on $X$ is nothing more than the quotient map from
$X$ to $X/H\cong (X^H)\unit$ which clearly factors to a bijection of
$X/H$ onto itself.  Now consider the source map $s:X\rightarrow
G^{(0)}$.  All that we need to do to show $s$ factors to a bijection 
is show that if $s(x) = s(y)$ then
there exists $[w,z]\in X^H$ such that $[w,z]\cdot x = y$.  However if
$s(x) = s(y)$ then $[x,y]\in X^H$ and it's easy to see that $[y,x]\cdot
x = y$.  All that's left is to show that the actions commute.
Well, suppose $\gamma \in G$, $z\in X$ and $[x,y]\in X^H$ with the
appropriate ranges and sources.  Let $\eta\in G$ be such that
$y\cdot \eta = z$.  Then $([x,y]\cdot z)\cdot\gamma = x\cdot
\eta\gamma$.  Since $y\cdot\eta\gamma = z\cdot\gamma$ it follows that 
$[x,y]\cdot (z\cdot\gamma) = x\cdot\eta\gamma$.  Thus the actions
commute and $X$ is a $(X^H,H)$-equivalence. 
\end{proof}

It turns out that this situation is the most general one.  

\begin{prop}
\label{prop:20}
Suppose $G$ and $H$ are locally compact Hausdorff groupoids and 
$Z$ is a $(G,H)$-equivalence.  There is a natural isomorphism $\phi$
between $G$ and the imprimitivity groupoid $Z^H$ defined as follows.
Given $[x,y]\in Z^H$ we have $s(x) = s(y)$ so that there exists a
unique $\gamma\in G$ such that $x = \gamma\cdot y$ and we define
$\phi([x,y]) := \gamma$.  
\end{prop}

\begin{proof}
First we show that $\phi$ is well defined.  Suppose $[x,y]=[x',y']$.
Then there exists $\eta\in H$ such that $x\cdot \eta = x'$ and $y \cdot
\eta = y'$.  Let $\gamma$ be the unique element of $G$ such that
$x= \gamma\cdot y$.  Since the actions commute $x\cdot \eta =
(\gamma\cdot y)\cdot\eta = \gamma\cdot(y\cdot\eta)$.  Therefore $\phi([x,y]) =
\gamma = \phi([x',y'])$.  Next we will show that $\phi$ is a
homomorphism. 
Now suppose $[x,y],[y,z]\in Z^H$ and let $\gamma,\eta\in G$ such that
$x= \gamma\cdot y$ and $y=\eta\cdot z$.  Then $s(\gamma) = r(y) =
r(\eta\cdot z) = r(\eta)$ and $x = \gamma\cdot y =
\gamma\eta\cdot z$ so that $\phi([x,z]) = \gamma\eta =
\phi([x,y])\phi([y,z])$.  

Now suppose $[x_i,y_i]\rightarrow [x,y]$, $\phi([x_i,y_i])=\gamma_i$,
and $\phi([x,y])=\gamma$.  Pass to a subnet.  We will show that there
is a sub-subnet such that $\gamma_i\rightarrow \gamma$.  Use the
fact that the quotient map $Z*Z\rightarrow Z^H$ is open to pass to a
subnet, reindex, and possibly choose new representatives $x_i$ and
$y_i$ so that $x_i\rightarrow x$ and $y_i\rightarrow y$.  Since $x_i =
\gamma_i \cdot y_i$ for all $i$ the fact that the action of $G$ is
proper implies that we can pass to a subnet and assume
$\gamma_i\rightarrow \eta$.  The continuity of the actions now implies
that $\gamma_i\cdot y_i\rightarrow \eta\cdot y$.  Since $Z$ is
Hausdorff we have $\eta\cdot y=\gamma\cdot y$ and freeness implies
$\eta=\gamma$ and we are done.  

Now suppose $\gamma\in G$.  Then $s(\gamma)\in G\unit$ we can choose an
$x\in Z$ such that $r(x)=s(\gamma)$.  Furthermore, since the range map
on $Z$ factors to a bijection, if $y\in Z$ such that $r(y) = s(\gamma)$
then there exists $\eta\in H$ such that $x = y\cdot \eta$.  Now define
$\psi:G\rightarrow Z^H$ by $\psi(\gamma) = [x,\gamma\cdot x]$.  We need
to show $\psi$ is well defined.  Given $\eta\in H$ we have
$(x,\gamma\cdot x)\cdot \eta = (x\cdot\eta,\gamma\cdot(x\cdot\eta))$ and
therefore $[x,\gamma\cdot x] = [x\cdot\eta,\gamma\cdot(x\cdot\eta)]$.  Thus
$\psi(\gamma)$ is independent of the choice of $x$.  Next, it's clear
that $\phi(\psi(\gamma)) = \phi([x,\gamma\cdot x]) = \gamma$.
Furthermore, if $\phi([x,y])=\gamma$ then $y=\gamma\cdot x$ so that
$\psi(\phi([x,y])) = \psi(\gamma) = [x,\gamma\cdot x] = [x,y]$.
Therefore $\phi$ and $\psi$ are inverses and $\phi$ is a bijection.
The last thing we need to show is that $\psi$ is continuous.  Suppose
$\gamma_i\rightarrow \gamma$.  Pass to a subnet.  Our goal is to show
a sub-subnet of $\psi(\gamma_i)$ converges to $\psi(\gamma)$.  Well
$s(\gamma_i)\rightarrow s(\gamma)$ and we can use the fact that the
range map on $Z$ is open to pass to a subnet and find $x_i\rightarrow
x$ such that $s(x_i)=r(\gamma_i)$.  It follows that $\gamma_i\cdot
x_i\rightarrow \gamma\cdot x$ and $[x_i,\gamma_i\cdot x_i]\rightarrow
[x,\gamma\cdot x]$. Thus $\psi$ is continuous and $\phi$ is an
isomorphism of locally compact Hausdorff groupoids. 
\end{proof}

The following proposition describes the imprimitivity groupoid that
will be of interest in Chapter \ref{cha:fine-structure}.  
The principal action is exactly the one
described in Example \ref{ex:10}.

\begin{prop}
\label{prop:21}
\index{imprimitivity groupoid}
Suppose $G$ is a locally compact Hausdorff groupoid with a Haar system
and that $H$ is a closed subgroupoid such that the range and source
maps restricted to $H$ are open.  In particular, this is true if $H$
has its own Haar system.  If $X=s\inv(H\unit)$ then $H$ acts on the
right of $X$ via multiplication and with this action $X$ is a strong
principal $H$-space.  The associated imprimitivity groupoid, in this
case denoted $G^H$, has
a Haar system $\{\mu^{[\xi]}\}$ defined for $f\in C_c(G^H)$ by 
\begin{equation}
\label{eq:2}
\int_{G^H} f([\xi,\eta])d\mu^{[\zeta]}([\xi,\eta]) = \int_G f([\zeta,\eta])
d\lambda_{s(\zeta)}(\eta).
\end{equation}
\end{prop}

\begin{proof}
The action of $H$ on $X$ is exactly the action described in 
Example \ref{ex:10} 
and it was shown there that if $s:X\rightarrow H\unit$ is the restriction of the
source map to $X$ and $\xi\cdot \eta = \xi\eta$ then $X$ is a strong
right $H$-space.  If $\xi\cdot \eta = \xi\eta = \xi$ for $\xi\in X$
and $\eta\in H$ then we can multiply both sides by $\xi\inv$ and get 
$\eta = \xi\inv\xi = s(\xi)\in H\unit$.  Thus the action of $H$ on $X$
is free.  Now suppose $\xi_i\rightarrow \xi$ and
$\xi_i\cdot\eta_i\rightarrow \zeta$.  Using the fact that the groupoid
operations are continuous we have 
\[
\xi_i\inv(\xi_i\eta_i) = \eta_i\rightarrow \xi\inv\zeta.
\]
Since $H$ is closed, $\xi\inv\zeta\in H$ and $\eta_i$ converges in $H$.
Thus the action of $H$ on $X$ is proper and therefore principal. Since
$H$ has open range and source maps, by assumption, and $X$ is a strong
principal right $H$-space, we can use Proposition \ref{prop:19} to
construct the imprimitivity groupoid $G^H$.  

We would like to show that the $\mu^{[\xi]}$ defined in the statement
of the proposition form a Haar system for $G^H$.  This is actually
proved in \cite{irredreps} and has been expanded here for
reference.  First we show that \eqref{eq:2} is well defined.  Suppose
$\gamma \in H$ such that $r(\gamma) = s(\zeta)$.  Then, noting that
applying left invariance to $\lambda_u=(\lambda^u)\inv$ gives us ``right invariance'',
we get
\[
\int_G f([\zeta\gamma,\eta])d\lambda_{s(\gamma)}(\eta) = 
\int_G f([\zeta\gamma,\eta\gamma]) d\lambda_{s(\zeta)}(\eta)
= \int_G f([\zeta,\eta])d\lambda_{s(\zeta)}(\eta).
\]
Thus \eqref{eq:2} is independent of the representative $\zeta$.  It is
clear that this defines a positive linear functional on $C_c(G^H)$ so
that $\mu^{[\zeta]}$ is a non-negative Radon measure for all $[\zeta]\in
G^H$.  Now suppose $[\zeta,\gamma]\in (G^H)^{[\zeta]}$ and $U$ is an open
neighborhood of $[\zeta,\gamma]$.  We may as well assume $U$ is
relatively compact and choose $f\in C_c(G^H)^+$ such that, $0\leq f
\leq 1$, $f$ is one on
$[\zeta,\gamma]$, and $f$ is zero off $U$.  It follows that 
\[
\mu^{[\zeta]}(U) \geq \int_{G^H} f([\xi,\eta])d\mu^{[\zeta]}([\xi,\eta]) = 
\int_G f([\zeta,\eta])d\lambda_{s(\zeta)}(\eta) > 0.
\]
Now if $[\xi,\gamma]\not\in G^{[\zeta]}$ we can pick a relatively
compact open neighborhood $U$ of $[\xi,\gamma]$ which is disjoint from
$(G^H)^{[\zeta]}$.  Since $U$ is relatively compact we can choose another
relatively compact neighborhood $V$ such that $U\subset V$ and $V$ is
disjoint from $(G^H)^{[\zeta]}$ as well.  
Choose $f\in C_c(G^H)^+$ so that, $0\leq f \leq 1$, $f$ is one on $U$,
and $f$ is zero off $V$.  Then
\[
\mu^{[\zeta]}(U) \leq
\int_{G^H} f([\xi,\eta])d\mu^{[\zeta]}([\xi,\eta]) = 
\int_G f([\zeta,\eta])d\lambda_{s(\zeta)}(\eta) = 0
\]
where the last equality holds because $f$ is supported off
$(G^H)^{[\zeta]}$.  It follows that $\supp\mu^{[\zeta]} =
(G^H)^{[\zeta]}$.  

Next we need to prove the continuity condition.  We start by showing
that given $f\in C_c(X* X)$ the function 
\begin{equation}
\label{eq:39}
\gamma \mapsto \int_G f(\gamma,\eta)d\lambda_{s(\gamma)}(\eta)
\end{equation}
is continuous on $X$.  Use Lemma \ref{lem:8} to extend $f$ to 
$C_c(X\times X)$.  Just as in the proof of Proposition \ref{prop:14}
we find a net of sums of elementary tensors $k_i = \sum_j g_i^j\otimes
h_i^j$ which converge to
$f$ uniformly 
and are all supported in some compact set.  We can then use the fact that 
\[
\gamma\mapsto \int_G g_i^j\otimes h_i^j(\gamma,\eta)d\lambda_{s(\gamma)}(\eta)
= g_i^j(\gamma) \int_G h_i^j(\eta)d\lambda_{s(\gamma)}(\eta)
\]
is clearly a continuous function to prove that \eqref{eq:39} is
continuous.  The computation is exactly the same as in the proof of
Proposition \ref{prop:14} and won't be reproduced here.  Next, we
need to show that the ``factorization'' of \eqref{eq:39} is
continuous.  This result is proved in \cite[Lemme 1.3]{frenchrenault}
but that paper is in French so the proof is included for reference.  
Define 
\[
Y = G^H*X := \{([\gamma,\eta],\zeta)\in G^H\times X:[\gamma]=[\zeta]\}.
\]
and let $\phi:X*X\rightarrow Y$ be given by $\phi(\gamma,\eta) =
([\gamma,\eta],\gamma)$.  We claim that $\phi$ is a homeomorphism.  It
is clear that $\phi$ is continuous and it is straightforward to show
that $\phi$ is bijective.  In particular, the inverse of
$([\gamma,\eta],\zeta)$ is $(\zeta,\eta\delta)$ where $\delta$ is the
unique element of $H$ such that $\zeta=\gamma\delta$.  
We will restrict ourselves to showing that
it has a continuous inverse.  Suppose
$([\gamma_i,\eta_i],\zeta_i)\rightarrow ([\gamma,\eta],\zeta)$
and let $\delta_i$ and $\delta$ be as above.  Pass to a subnet.  It
will suffice to show that there is a sub-subnet such that
$\eta_i\delta_i\rightarrow \eta\delta$.  By passing to a subnet we may
assume that $\eta_i\rightarrow \eta$ and $\gamma_i\rightarrow\gamma$.
However $\zeta_i = \gamma_i\delta_i$ and $\zeta=\gamma\delta$ so that
we can use the fact that the action of $H$ is principal to pass to
another subnet and assume $\delta_i\rightarrow \delta$.  The result
follows.  

Next, observe that given $f\in C_c(Y)$ we have
\begin{align*}
\int_G f(\phi(\zeta,\eta))d\lambda_{s(\zeta)}(\eta) &= 
\int_G f([\zeta,\eta],\zeta)d\lambda_{s(\zeta)}(\eta)\\ 
&= \int_Y f([\gamma,\eta],\xi)d(\mu^{[\zeta]}\times
\delta_\zeta)([\gamma,\eta],\xi).
\end{align*}
In other words $\phi$ identifies the measures $\lambda_{s(\gamma)}$
and $\mu^{[\zeta]}\times \delta_\zeta$.  Fix $f\in C_c(G^H)$ and
suppose $g\in C_c(X)$.  Define $F\in C_c(Y)$ by $F = f\otimes g$ so that
$F([\gamma,\eta],\zeta) = f([\gamma,\eta])g(\zeta)$ and observe that
\[
\zeta \mapsto \int_Y
F([\gamma,\eta],\xi)d(\mu^{[\zeta]}\times\delta_\zeta)([\gamma,\eta],\xi) = 
\int_G F(\phi(\zeta,\eta))d\lambda_{s(\zeta)}(\eta)
\]
is continuous since \eqref{eq:39} is continuous.  However this implies
that 
\begin{equation}
\label{eq:40}
\zeta \mapsto \int_Y F([\gamma,\eta],\xi) d(\mu^{[\zeta]}\times \delta_\zeta)([\gamma,\eta],\xi) 
= g(\zeta) \int_{G^H} f([\gamma,\eta])d\mu^{[\zeta]}([\gamma,\eta])
\end{equation}
is continuous.  Suppose $[\zeta_i]\rightarrow[\zeta]$ and, by possibly
passing to a subnet and choosing new representatives, assume that
$\zeta_i\rightarrow \zeta$.  We can, without loss of generality
suppose that this sequence converges inside some compact neighborhood
$K$.  Choose $g\in C_c(X)$ so that $g$ is one on $K$.  Then the fact
that \eqref{eq:40} is continuous implies that 
\[
\int_{G^H} f([\gamma,\eta])d\mu^{[\zeta_i]}([\gamma,\eta]) \rightarrow
\int_{G^H} f([\gamma,\eta])d\mu^{[\zeta]}([\gamma,\eta]).
\]
This is enough to show that $[\zeta]\mapsto \int_{G^H} fd\mu^{[\zeta]}$ is
continuous and it's easy to see that this function is compactly
supported.  

All we have to do now is verify the left invariance condition.
Suppose $[\gamma,\eta]\in G^H$ and $f\in C_c(G^H)$.  Then 
\begin{align*}
\int_{G^H} f([\gamma,\eta][\zeta,\xi])d\mu^{[\eta]}([\zeta,\xi]) &= 
\int_G f([\gamma,\eta][\eta,\xi])d\lambda_{s(\eta)}(\xi) = 
\int_G f([\gamma,\xi])d\lambda_{s(\eta)}(\xi) \\
&= \int_G f([\gamma,\xi])d\lambda_{s(\gamma)}(\xi)\quad\text{since
  $s(\gamma) = s(\eta)$} \\
&= \int_{G^H} f([\zeta,\xi])d\mu^{[\gamma]}([\zeta,\xi]).
\end{align*}
Thus $\mu^{[\zeta]}$ is left invariant and the collection
$\{\mu^{[\zeta]}\}$ forms a Haar system for $G^H$. 
\end{proof}

\begin{remark}
\index{imprimitivity groupoid}
Observe that in Proposition
\ref{prop:21} we did not have to 
assume that $H$ has a Haar system for $G^H$ to have one, 
only that it has open range and
source maps.  It seems to be an open, and difficult question, to ask
if there are two equivalent groupoids such that one has a Haar system
and the other doesn't.
\end{remark}


\section{Amenable Groupoids}
\label{sec:amenable}

Groupoid amenability will only play a role as an 
assumption in one of the major results of the thesis, and will never
be used in a technical way.  However, since groupoid amenability is a very
new subject, it seems appropriate to include some of the basic
definitions and theorems.  We will use the notion of amenability that
Anantharaman-Delaroche and Renault use in \cite{amenable}.  

\begin{definition}
Suppose $X$ and $Y$ are locally compact Hausdorff spaces and
$\pi:X\rightarrow Y$ is a continuous surjection.  A continuous system
of measures for $\pi$ is a set of positive Radon measures 
$\alpha=\{\alpha^y\}_{y\in Y}$ on $X$ such that 
\begin{enumerate}
\item $\supp \alpha^y \subseteq \pi\inv(y)$, 
\item and for every $f\in C_c(X)$ the function $\alpha(f):Y\rightarrow \C$
such that
\[
\alpha(f)(y) := \int_X f d\alpha^y
\]
is continuous and compactly supported. 
\end{enumerate}
We say $\alpha$ is {\em proper} if $\alpha^y \ne 0$ for all $y\in Y$
and we say $\alpha$ is {\em full} if $\supp \alpha^y = \pi\inv(y)$ for all
$y\in Y$.  
\end{definition}

\begin{definition}
Let $G$ be a locally compact groupoid, $X$ and $Y$ locally compact
strong $G$-spaces, and $\pi:X\rightarrow Y$ a continuous
$G$-equivariant surjection.  An {\em invariant continuous} $\pi$-system
is a continuous system of measures $\alpha$ for $\pi$ such
that given $(\gamma,y)\in G*Y$ we have $\gamma \alpha^y =
\alpha^{\gamma\cdot y}$.  Here $\gamma \alpha^{y}$ denotes the measure 
defined on $C_c(X)$ via 
\[
\gamma \alpha^{y}(f) = \int_X f(\gamma\cdot x) d\alpha^{y}(x).
\]
\end{definition}

\begin{example}
\index{Haar system}
Given a locally compact groupoid $G$ with Haar system $\lambda$, it is
clear that $\lambda$ is a full, invariant, continuous $r$-system.  In
fact, ``full, invariant, continuous $r$-system'' is just a very short,
alliterative way to define a Haar system.  
\end{example}

\begin{definition}
\index{acim@a.c.i.m.}
Let $G$ be a locally compact groupoid, $X$ and $Y$ locally compact
strong $G$-spaces, and $\pi:X\rightarrow Y$ a continuous
$G$-equivariant surjection.  An {\em approximate continuous invariant
  mean} (a.c.i.m.) for $\pi$ is a net $\{m_i\}$ of continuous systems
of probability measures for $\pi$ such that $\| \gamma m_i^y -
m_i^{\gamma\cdot y}\|_1\rightarrow 0$ uniformly on compact
subsets of $G*Y$.\footnote{For a finite measure $m$ on $Y$ we define
  $\|m\|_1 := |m|(Y)$.}  
\end{definition}

It may be helpful to observe that the invariance condition for an
invariant $\pi$-system $\alpha$ can be written $\|\gamma\alpha^{y} -
\alpha^{\gamma\cdot y}\|_1 = 0$.  In this light, an a.c.i.m.\ is
a net of systems of {\em probability} measures that, in the limit, behave like
an invariant system.  Next, we define amenability for maps,
groupoids and $G$-spaces.  It's notable that amenability for
groupoid actions is defined in terms of the transformation groupoid.  

\begin{definition}
\label{def:16}
\index{amenable}
\index{groupoid!amenable}
We say a continuous $G$-equivariant surjection $\pi:X\rightarrow
Y$ between strong $G$-spaces $X$ and $Y$ is {\em amenable} if it
admits an approximate continuous invariant mean.  In particular, we
say that a locally compact Hausdorff groupoid with open range
and source maps is (topologically) amenable if the range map
$r:G\rightarrow G\unit$ is amenable, where we view $G$ and $G\unit$ as
$G$-spaces in the usual way.  We say that a locally compact
$G$-space $X$ is amenable if the groupoid $G\ltimes X$ is.  
\end{definition}

\begin{remark}
Suppose $G$ is a locally compact Hausdorff groupoid with open range
and source maps.  It's clear that the action of
$G$ on itself is strongly continuous.  Furthermore, the action of $G$ on
$G\unit$ is also strongly continuous.  Finally, given $\gamma,\eta$ in
$G$ such that $s(\gamma) = r(\eta)$ we have
\[
r(\gamma\eta) = r(\gamma) = \gamma\cdot r(\eta)
\]
so that $r$ is a continuous $G$-equivariant surjection.  Thus,
Definition \ref{def:16} makes sense.
\end{remark}

\begin{example}
Suppose $G$ is a locally compact Hausdorff groupoid which admits a
Haar system $\{m^u\}$ of probability measures.  Then $m$ is,
by definition, a continuous invariant $r$-system.  We can form an
a.c.i.m.\ using the constant sequence $m_i^u = m^u$ for all $i$. In
some sense, amenability is meant to be a generalization of this
situation.  
\end{example}

Next we give a couple of useful characterizations of amenability.  As
usual, we have to start out with a definition. 

\begin{definition}
\label{def:55}
\index{conditionally compact}
Given a continuous surjection $\pi:X\rightarrow Y$ between locally
compact Hausdorff spaces a set $A\subset X$ will be called
{\em $\pi$-compact} if, for every compact $K\subset Y$, its intersection
with $\pi\inv(K)$ is compact.  We will use $C_{c,\pi}(X)$ to denote
the space of continuous functions on $X$ with $\pi$-compact support.
When $G$ is a locally compact Hausdorff groupoid and 
$r:G\rightarrow G\unit$ is the range map then $r$-compact sets will be
called {\em conditionally compact}.
\end{definition}

\begin{definition}
\label{def:21}
Suppose the locally compact Hausdorff groupoid $G$ acts on a locally 
compact Hausdorff space $X$.  A function $e$ on $G\ltimes X$ is said
to be of {\em positive type} if, for every $x\in X$, every $n\in\N$,
and every $\gamma_1,\ldots,\gamma_n\in G^{r(x)}$, and
$z_1,\ldots,z_n\in\C$ we have
\[
\sum_{i=1}^n \sum_{j=1}^n e(\gamma_i\inv\gamma_j,\gamma_i\inv\cdot
x)\bar{z}_i z_j \geq 0.
\]
We denote by $e\unit$ the restriction of $e$ to $G\unit\ltimes X$
\end{definition}

The following is a restatement of \cite[Propositions
2.2.5,2.2.6]{amenable}.

\begin{prop}
\label{prop:22}
Suppose $G$ is a locally compact Hausdorff groupoid and let $G$ have a
strongly continuous action on $X$ and $Y$.  Now suppose we have a
continuous $G$-equivariant surjection $\pi:X\rightarrow Y$ such that
there exists an invariant continuous $\pi$-system $\alpha$ of
measures.  Then properties {\em (b)},{\em (c)} and {\em (d)} below are all equivalent
and each of them implies {\em (a)}.  Furthermore, if $X$ is a proper
$G$-space then the following are equivalent.  
\begin{enumerate}
\item $\pi$ is amenable. 
\item There exists a net of positive functions 
$\{g_i\}\subset C_{c,\pi}(X)$ such that 
\begin{enumerate}
\item $\int g_i d\alpha^y =1$ for all $y\in Y$ and
\item $\int |g_i(\gamma\cdot x)-g_i(x)|d\alpha^y(x)$ converges to zero
  uniformly on the compact subsets of $G*Y$. 
\end{enumerate}
\item There exists a net $\{\xi_i\}\subset C_{c,\pi}(X)$ such that 
\begin{enumerate}
\item $\int |\xi_i|^2 d\alpha^y = 1$ for all $y\in Y$
\item $\int |\xi_i(\gamma\cdot x)-\xi_i(x)|^2 d\alpha^y(x)$ converges
  to zero uniformly on the compact subsets of $G*Y$.  
\end{enumerate}
\item There exists a net $\{\xi_i\}\subset C_{c,\pi}(X)$ such that the
  net $\{e_i\}$ of positive type functions on the groupoid $G\ltimes
  Y$ defined by 
\[
e_i(\gamma,y) = \int_X \overline{\xi_i(x)} \xi_i(\gamma\inv\cdot x)
d\alpha^y(x)
\]
satisfies: 
\begin{enumerate}
\item $e_i\unit=1$ for all $i$, 
\item $\lim_i e_i = 1$ uniformly on compact subsets of $G\ltimes Y$.  
\end{enumerate}
\end{enumerate}
\end{prop}

\begin{remark}
Proposition \ref{prop:22} is especially useful in the case of groupoid
amenability.  If the locally compact Hausdorff groupoid $G$ has a Haar
system then, by definition, there exists a continuous invariant
$r$-system.  Furthermore, the action of $G$ on itself is always
proper.  Thus, in this situation, conditions (a)--(d) in Proposition
\ref{prop:22} are all equivalent. 
\end{remark}

There are three more results in \cite{amenable} that we will state here
without proof.  
In \cite{amenable} they are Propositions 2.2.13, 5.1.1, and
5.1.2 respectively. 

\begin{prop}
\index{G,H-equivalence@$(G,H)$-equivalence}
For locally compact Hausdorff groupoids, amenability is invariant under
groupoid equivalence. 
\end{prop}

\begin{prop}
A locally closed subgroupoid $H$ of an amenable
locally compact Hausdorff groupoid $G$ is amenable.\footnote{A set is locally closed if it is the
  intersection of a closed set and an open set.  In particular, every
  closed set is locally closed.} In particular, the stabilizer subgroupoid
and all of the stabilizers subgroups are amenable if $G$ is.  
\end{prop}

\begin{prop}
\label{prop:23}
Let $G$ and $H$ be locally compact Hausdorff groupoids with open range and
source maps and $\pi:G\rightarrow H$ a continuous open surjective
homomorphism such that $\pi\unit:G\unit\rightarrow H\unit$ is a
homeomorphism.  We denote by $N=\{\gamma\in G: \pi(\gamma)\in
H\unit\}$ the kernel of $\pi$.  Then $G$ is amenable if and only if
$N$ and $H$ are both amenable.  
\end{prop}

\begin{corr}
\index{amenable}
Let $G$ be a locally compact Hausdorff groupoid with open range and
source and continuously
varying stabilizer.  Then $G$ is amenable
if and only if both the stabilizer subgroupoid $S$ and the orbit
groupoid $R_Q$ are. 
\end{corr}

\begin{proof}
Let $\pi:G\rightarrow R_Q$ be the canonical homomorphism.  Since the
stabilizers vary continuously we know from Proposition \ref{prop:11}
that $\pi$ is a continuous, surjective, open homomorphism and that
$R_Q$ is a locally compact Hausdorff groupoid with open range and source.
Furthermore, $\pi$ clearly restricts to a homeomorphism of $G\unit$
with $R_Q\unit$.  (We usually identify $R_Q\unit$ with $G\unit$ via
this map.)  It is clear that the kernel of $\pi$ is the stabilizer
subgroupoid $S$.  It follows from Proposition \ref{prop:23} that $G$
is amenable if and only if both $R_Q$ and $S$ are.
\end{proof}

\begin{remark}
We have given the briefest of introductions to topological
amenability.  None of the proofs of the above theorems are particularly
difficult and, as stated before, can all be found in \cite{amenable}.
There is also a measure theoretic notion of amenability.  
In \cite{amenable} Anantharaman-Delaroche and Renault go into
extensive detail about measurable amenability and the properties of
amenable measured groupoids.  They also go into some depth about the
relationship between amenability and the operator algebras associated
to groupoids.  
\end{remark}


\chapter{Groupoid Group Bundles}
\label{cha:bundles}
This chapter contains two results concerning abelian group bundles.
The first is the definition of a principal $S$-bundle where $S$ is an
abelian group bundle.  We construct a sheaf cohomology theory for
$S$ and show that the isomorphism classes of principal $S$-bundles are
in one-to-one correspondence with elements of the first cohomology
group.  We also connect this material back to existing work done for
``locally $\sigma$-trivial'' bundles.  In Section
\ref{sec:duality} we prove a Pontryagin duality theorem for
continuously varying, abelian group bundles and in Section
\ref{sec:opencounter} we present an interesting counterexample which
shows that not every continuous bijective homomorphism between second
countable, continuously varying, abelian group bundles is an
isomorphism.  This material is (mostly) low level and has been
included before Chapter \ref{cha:crossed} to emphasize that.  

\section{Principal Group Bundles}
\label{sec:principal}

The goal of this section is to generalize principal group bundles to
situations where the bundle may not be locally trivial.  In
particular, we would like to be able to deal with group bundles as
defined in Definition \ref{def:24}. However, as in the group case,
commutativity will be an essential assumption.  

\begin{remark} 
For the rest of this section $S$ will always denote an abelian locally compact
Hausdorff group bundle with bundle map $p$.  We will endeavor to
state these hypothesis on all of the important theorems.  
\end{remark}

We will develop a theory of principal
$S$-bundles that mirrors the classic theory of principal $H$-bundles where
$H$ is a locally compact Hausdorff group.  As in the group case, there
is a nice one-to-one correspondence between $S$-bundles and principal
$S$-spaces with local sections.  However, in order to be consistent we
must start with the bundle definition and will develop the
correspondence later.  The material in this section is modeled off \cite[Section
4.2]{tfb}.  This first definition is really a matter of notation and
terminology.   Any surjection can be viewed as a bundle
map, although it's not always useful to do so.  

\begin{definition}
Suppose $Y$ and $X$ are topological spaces and 
$q:X\rightarrow Y$ is a continuous surjection.  Then $X$ is called a
{\em (topological) bundle} over $Y$ with bundle map $q$ 
and the fibres $q\inv(y)$ are denoted by $X_y$ for all $y\in Y$.  
\end{definition}

Principal $S$-bundles are bundles which are locally isomorphic to
$S$.  When dealing with ``local triviality'' conditions one must consider
the fact that there may be different trivializations for the same
bundle.  Our first
definition of principal bundle will depend on the trivialization.  

\begin{definition}
\index{principal S-bundle@principal $S$-bundle}
\index{group bundle}
Let $S$ be an abelian locally compact Hausdorff group bundle with
bundle map $p$.  Suppose
$X$ is a locally compact Hausdorff bundle over $S\unit$ with bundle
map $q$.  Furthermore, 
suppose there is an open cover $\mcal{U}=\{U_i\}_{i\in I}$ of $S\unit$ such
that for each $i\in I$ there is a homeomorphism $\phi_i:
q\inv(U_i)\rightarrow p\inv(U_i)$ with $p\circ \phi_i = q$.
Finally, suppose that for all $i,j\in I$ there is a section
$\gamma_{ij}$ of $S|_{U_i\cap U_j} = p\inv(U_i\cap U_j)$ such that 
\[
\phi_i\circ \phi_j\inv(s) = \gamma_{ij}(p(s))s
\]
for all $s\in S|_{U_i\cap U_j}$.  Such a bundle is called a {\em
  principal $S$-bundle with trivialization $(\mcal{U},\phi,\gamma)$.}  The
maps $\phi=\{\phi_i\}$ are referred to as {\em trivializing maps} and
the sections $\gamma= \{\gamma_{ij}\}$ are referred to as {\em
  transition maps.}
\end{definition}

\begin{remark}
\index[not]{$U_{ij}$}
If $\mcal{U} = \{U_i\}$ is an open cover then we use the
usual notation $U_{ij} := U_i\cap U_j$.  
\end{remark}

We would like to see that there is some freedom with respect to the
neighborhoods in a trivialization.  In particular, we would
like to be able to refine them at will.  

\begin{prop}
\label{prop:24}
If $X$ is a principal $S$-bundle with trivialization
$(\mcal{U},\phi,\gamma)$ 
then
for any refinement $\mcal{V}$ of the open cover $\mcal{U}$ we can make $X$ 
a principal $S$-bundle with trivialization
$(\mcal{V},\phi_{\mcal{V}},\gamma_{\mcal{V}})$ where $\phi_{\mcal{V}}$
and $\gamma_{\mcal{V}}$ denote the
natural restriction of the trivializing and transition 
maps to $\mcal{V}$. 
\end{prop}

\begin{proof}
Suppose $\mcal{V} = \{V_j\}_{j\in J}$ is a refinement of
$\mcal{U}=\{U_i\}_{i\in I}$ with refining map ${r:I\rightarrow J}$ so
that $V_j\subseteq U_{r(j)}$ for all $j\in J$.  We can define
trivializing maps on $V_j$ by $\phi'_j = \phi_{r(j)}|_{V_j}$ for all
$j\in J$ and transition maps $\gamma'_{ij} =
\gamma_{r(i)r(j)}|_{V_{ij}}$ for all $i,j\in J$.  It is easy to
see that $\phi'_i$ is a homeomorphism and that $p\circ \phi'_i = q$;
all that you have to check is that $\phi'_i$ is surjective, which
follows from the fact that $\phi$ preserves fibres.  Furthermore
\[
\phi'_i\circ (\phi'_j)\inv(s) = \phi_{r(i)}\circ \phi_{r(j)}\inv(s) 
= \gamma_{r(i)r(j)}(p(s))s = \gamma_{ij}'(p(s))s. 
\]
Thus the maps $\phi_{\mcal{V}}= \{\phi'_j\}$ and $\gamma_{\mcal{V}} =
\{\gamma_{ij}\}$ make $X$ into a principal $S$-bundle.  
\end{proof}

\begin{remark}
We will usually drop the $\mcal{V}$ from $\phi_{\mcal{V}}$ and
$\gamma_{\mcal{V}}$ to avoid notational clutter. 
\end{remark}

Now we define the notion of an $S$-bundle isomorphism.  
The basic idea is that an isomorphism is
locally given by sections. 

\begin{definition}
\label{def:18}
Suppose $q:X\rightarrow S\unit$ and $r:Y\rightarrow S\unit$ are both
principal $S$-bundles with trivializations $(\mcal{U},\phi,\gamma)$
and $(\mcal{V},\psi,\eta)$ respectively.  Let $\mcal{W}$ be some
common refinement of $\mcal{U}$ and $\mcal{V}$ with refining maps $\sigma$
and $\rho$ respectively. Furthermore, suppose 
$\Omega:X\rightarrow Y$ is a homeomorphism such that $r\circ\Omega = q$ and
that for all $W_i\in \mcal{W}$
$\beta_i:W_i\rightarrow S$ is a section of $p$ 
such that for all $s\in p\inv(W_i)$
\[
\psi_{\sigma(i)} \circ \Omega \circ \phi_{\rho(i)}\inv(s) = \beta_i(p(s))s.
\]
Then $(\mcal{W},\Omega,\beta)$ is an $S$-bundle {\em isomorphism} of
$X$ onto $Y$.  
\end{definition}

A special case of a principal $S$-bundle isomorphism is when some
principal $S$-bundle $X$ comes with two trivializations which
``agree'' on a common refinement.  

\begin{definition}
\label{def:25}
Suppose $X$ is a principal $S$-bundle and that
$(\mcal{U},\phi,\gamma)$ and $(\mcal{V},\psi,\eta)$ are two
trivializations of $X$.  We say that $(\mcal{U},\phi,\gamma)$ and
$(\mcal{V},\psi,\eta)$ are {\em equivalent} if there exists some
common refinement $\mcal{W}$ of $\mcal{U}$ and $\mcal{V}$ such that
the maps
\begin{align*}
\id:X\rightarrow X&: x\mapsto x, &
\iota_i:W_i\rightarrow S&: u\mapsto u,
\end{align*}
define an $S$-bundle isomorphism $(\mcal{W},\id,\iota)$ from $X$ with
trivialization $(\mcal{U},\phi,\gamma)$ onto $X$ with trivialization
$(\mcal{V},\psi,\eta)$.  
\end{definition}

\begin{remark}
We won't need to use this, but it's easy to see that Definition
\ref{def:25} defines an equivalence relation on the set of
trivializations of some bundle $X$.  
\end{remark}

\begin{example}
\label{ex:8}
Suppose $X$ is a principal $S$-bundle with trivialization
$(\mcal{U},\phi,\gamma)$ and $\mcal{V}$ is a refinement of $\mcal{U}$.
Then $(\mcal{U},\phi,\gamma)$ and
$(\mcal{V},\phi_{\mcal{V}},\gamma_\mcal{V})$ are equivalent.  
\end{example}

We can now write down a better definition of principal $S$-bundle.  

\begin{definition}
\label{def:17}
\index{principal S-bundle@principal $S$-bundle}
Suppose $S$ is an abelian locally compact Hausdorff group bundle and
$X$ is a locally compact Hausdorff bundle over $S\unit$.  Suppose
$\mcal{A}$ is a set of triplets of the form $(\mcal{U},\phi,\gamma)$
such that $X$ is a principal
$S$-bundle with trivialization $(\mcal{U},\phi,\gamma)$ for all
$(\mcal{U},\phi,\gamma)\in \mcal{A}$.  Then $\mcal{A}$ is called
{\em pairwise equivalent} if any two trivializations in $\mcal{A}$ are
equivalent.  A pairwise equivalent collection $\mcal{A}$ is called 
{\em maximal} if it is
not properly contained in any other pairwise equivalent collection.  
A locally compact Hausdorff bundle over
$S\unit$ equipped with a maximal pairwise equivalent collection
of trivializations is called a {\em principal $S$-bundle}.  
\end{definition}

\begin{remark}
Every principal $S$-bundle $X$ with fixed trivialization
$(\mcal{U},\phi,\gamma)$ can be turned into a principal
$S$-bundle by using Zorn's Lemma to find the maximal pairwise equivalent
collection $\mcal{A}$ containing
$(\mcal{U},\phi,\gamma)$.  It follows from Example \ref{ex:8} 
that any refinement of a trivialization in $\mcal{A}$ is in
$\mcal{A}$. 
\end{remark}

\begin{definition}
Given two principal $S$-bundles $X$ and $Y$ we say $X$ is {\em
  isomorphic} to $Y$ if there is an isomorphism between the two
bundles with respect to some pair of trivializations. 
\end{definition} 

\begin{remark}
It is clear from Definition \ref{def:17} that if two principal
$S$-bundles $X$ and $Y$ are isomorphic with respect to some given pair of
trivializations then they are isomorphic with respect to all
trivializations.  Furthermore, since any two trivializations in the
maximal collections associated to $X$ and $Y$ are isomorphic with respect
to the identity map, the isomorphism between $X$ and $Y$ is basically
the same for all trivializations.  
\end{remark}

Next, we would like to mimic the group case and characterize the set
of all principal $S$-bundles by classes in some cohomology group.  We
will be using sheaf cohomology as defined and developed in
\cite[Section 4.1]{tfb}.  

\begin{prop}
\label{prop:26}
\index{cohomology}
\index[not]{$H^n(S)$}
Let $S$ be an abelian locally compact Hausdorff group bundle and for
$U$ open in $S\unit$ let $\mcal{S}(U) = \Gamma(U,S)$ be the set of
continuous sections from $U$ into $S$.  Then
$\mcal{S}$ is a sheaf and as such gives rise to a sheaf cohomology 
$H^n(S\unit;\mcal{S})$ which we shall denote by $H^n(S)$.  
\end{prop}

\begin{proof}
First we show that $\mcal{S}$ is a presheaf.  Observe that
$\Gamma(U,S)$ is non-empty for all $U$ because the inclusion of $U$
into $S$ is always a section for $p$.  It follows from the continuity of the
operations on $S$ that the 
pointwise multiplication of two continuous sections, or the inverse of a
continuous section, is still a continuous section.  It is easy to see that the
group axioms hold with respect to these operations, with the inclusion
of $U$ into $S$ acting as the identity.  

Next, observe that $\Gamma(\emptyset,S) = \{0\}$ since there is only
one ``function'' with the empty domain.  Finally, given open subsets
$U$ and $V$ of $S\unit$ such that $U\subset V$ we can define 
$\rho_{V,U}:\Gamma(V,S)\rightarrow \Gamma(U,S)$ to be given by
restriction to $U$.  It is clear that $\rho_{U,U} = \id$ and that for
$U\subset V\subset W$ we have 
\[
\rho_{W,V}\circ\rho_{V,U} = \rho_{W,U}.
\]
Thus $\mcal{S}$ defines a presheaf on $S\unit$.

Now suppose we have an open set $U\subset S\unit$ and a decomposition
$U = \bigcup_{i\in I} U_i$ of $U$ into open sets $U_i$.   Furthermore, suppose
we have $\gamma_i\in \Gamma(U_i,S)$ for all $i\in I$ and for all
$i,j\in I$ 
\[
\rho_{U_i,U_{ij}}(\gamma_i) = \rho_{U_j,U_{ij}}(\gamma_j).
\]
Tracing through the definitions we see that each $\gamma_i$ is a
continuous section on $U_i$ such that the $\gamma_i$ agree on
overlaps.  Therefore, we can define a continuous section $\gamma$ on $U$
in a piecewise fashion so that $\rho_{U,U_i}(\gamma) = \gamma_i$.
Furthermore, it is clear that $\gamma$ is uniquely determined by the
$\gamma_i$.  

Thus $\mcal{S}$ is a sheaf of groups on $S\unit$.   Furthermore,
because $S$ is an abelian group bundle it's easy to see that each
$\Gamma(U,S)$ is an abelian group and that $\mcal{S}$ is an abelian
sheaf.  As such it has an associated cohomology described in
\cite[Section 4.1]{tfb} which we use to define $H^n(S\unit;\mcal{S}) =
H^n(S)$.  
\end{proof}

At this point we can build the desired correspondence between principal
$S$-bundles and elements of $H^1(S)$.  The proof uses the details of
the construction of $H^1(S)$.  These details can be found in
\cite[Section 4.1]{tfb} and anyone unfamiliar with sheaf cohomology
should at least look through this section before working through the
following proof.  

\begin{theorem}
\index{New Result}
\label{prop:principcohom}
Suppose $S$ is an abelian locally compact Hausdorff group bundle.
There is a one-to-one correspondence between the isomorphism classes
of principal $S$-bundles and elements of the sheaf cohomology group
$H^1(S)$.  Given a principal bundle $X$ with trivialization 
$(\mcal{U},\phi,\gamma)$, the cohomology class in $H^1(S)$ associated
to $X$ is realized by the cocycle $\gamma$.  
\end{theorem}

\begin{proof}
First, suppose $X$ is a principal $S$-bundle and pick a trivialization
$(\mcal{U},\phi,\gamma)$.  Let $\mcal{S}$ be the sheaf defined in
Proposition \ref{prop:26}.  We will use the sheaf cohomology 
notation and definitions from \cite[Section 4.1]{tfb}.  Observe
that $\gamma_{ij}\in \mcal{S}(U_{ij}) =\Gamma(U_{ij},S)$ by
definition.  This implies that $\gamma = \{\gamma_{ij}\}$ forms a
chain in $C^1(\mcal{U},\mcal{S})$.  We need to show that this chain is
a cocycle; we need to show that $d(\gamma)$ is the trivial section.
This amounts to proving that 
\[
\gamma_{ij}(u)\gamma_{jk}(u) = \gamma_{ik}(u)
\]
for all $u\in U_{ijk}$.  Well, given $u\in U_{ijk}$ we have 
\begin{align*}
\gamma_{ij}(u)\gamma_{jk}(u) &= \gamma_{ij}(u)\gamma_{jk}(u)u = 
\phi_i\circ\phi_j\inv(\phi_j\circ\phi_k\inv(u)) \\
&= \phi_i\circ\phi_k\inv(u) = \gamma_{ik}(u) u \\
&= \gamma_{ik}(u).
\end{align*}
Thus $\gamma$ is a cocycle and as such we can use $\gamma$ to define a
class $[\gamma]\in H^1(S\unit;\mcal{S})$.  

We need to see that if we take a different trivialization for $X$ then
we get the same cohomology class.  We start by showing that if we take
a refinement of $\mcal{U}$ then we don't change
$[\gamma]$.  Suppose $\mcal{V}$ is a refinement of $\mcal{U}$
with refining map $r$. Then the trivialization of $X$ associated to
$\mcal{V}$ is, according to Proposition \ref{prop:24}, given by
$(\mcal{V},\phi_{\mcal{V}},\gamma_{\mcal{V}})$ where, in particular, 
$(\gamma_\mcal{V})_{ij} = \gamma_{r(i)r(j)}|_{V_{ij}}$.  It follows,
by definition, that $\gamma_{\mcal{V}} = r_*(\gamma)$ where $r_* :
H^n(\mcal{U},\mcal{S}) \rightarrow H^n(\mcal{V},\mcal{S})$ is the
homomorphism defined by restriction.  Since $H^n(S\unit;\mcal{S})$ is
a direct limit of the $H^n(\mcal{U};\mcal{S})$ with respect to the
restriction maps $r_*$ it is clear that $[\gamma] =
[\gamma_{\mcal{V}}]$ so that we can safely pass to refinements without
changing the associated cohomology class. 

Now suppose that $Y$ is another bundle isomorphic to $X$.  Suppose
$(\mcal{U},\phi,\gamma)$ is a trivialization for $X$ and 
$(\mcal{V},\psi,\eta)$ a trivialization for $Y$.  In
particular, $Y$ could be $X$ with an equivalent trivialization.  
The goal is to show that the cohomology classes associated to $X$ and
$Y$ with these trivializations are equal.  Let
$(\mcal{W},\Omega,\beta)$ be an isomorphism from $X$ to $Y$.  
By passing to a common
refinement and using the previous paragraph 
we can assume, without loss of generality, that $\mcal{U} =
\mcal{V}=\mcal{W}$.  Then, for all
$u\in U_{ij}$, 
\begin{align*}
\eta_{ij}(u)\beta_j(u) &= \eta_{ij}(u)\beta_j(u) u \\
&=\psi_i\circ\psi_j\inv(\psi_j\circ \Omega \circ \phi_j\inv(u)) \\
&= \psi_i\circ\Omega \circ\phi_i\inv(\phi_i\circ \phi_j\inv(u)) \\
&= \beta_i(u)\gamma_{ij}(u)u \\
&= \beta_i(u)\gamma_{ij}(u).
\end{align*}
Hence $\gamma\inv\eta$ is a boundary and
therefore $[\gamma]=[\eta]$ in $H^1(S\unit;\mcal{S})$.  This proves a
number of things.  First, it shows that no matter which trivialization
for $X$ we choose we get the same cohomology class in
$H^1(S\unit,\mcal{S})$.  Therefore we can define a map $X\mapsto [X]$
from the set of principal $S$-bundles to $H^1(S\unit;\mcal{S})$ by
$[X] = [\gamma]$ where $\gamma$ is the cocycle defined by any set of
transition maps for $X$. Furthermore, since we were working with an
arbitrary isomorphic bundle $Y$, this also shows that the map
$X\mapsto [X]$ is a well defined function from the set of isomorphism
classes of principal $S$-bundles into $H^1(S)$.  

Now we are going to construct an inverse map.  Suppose $c\in
H^1(S\unit;\mcal{S})$ is realized by $\gamma\in Z^1(\mcal{U},\mcal{S})$
for some open cover $\mcal{U}$.  Let $C = \coprod_i p\inv(U_i)$ be the
disjoint union of the $p\inv(U_i)$ and denote elements of $C$ by
$(s,i)$ where $s\in p\inv(U_i)$.  Observe that we can define
$\tilde{p}:C \rightarrow S\unit$ by 
by $\tilde{p}(s,i) = p(s)$, and that this map will be continuous since
\[
\tilde{p}\inv(O) = \bigcup_i p\inv(U_i\cap O)\times \{i\}
\] 
for any open set $O\subset S\unit$. From here on we will often denote
the clopen subset $p\inv(U_i)\times \{i\}$ of $C$ by $p\inv (U_i)$.
Since $p\inv(U_i)\times \{i\}$ and $p\inv (U_i)$ are clearly
homeomorphic there should be relatively little confusion.  
Define a relation on $C$
by $(s,i)\equiv (t,j)$ if and only if 
$p(s)=p(t)=u$ and $s = \gamma_{ij}(u)t$.  We need
to show that this is an equivalence relation.
Since $\gamma$ is a cocycle we know that 
\begin{equation}
\label{eq:3}
\gamma_{ij}\gamma_{jk} = \gamma_{ik}
\end{equation}
for all $i,j,k$.
In particular  \eqref{eq:3} implies 
$\gamma_{ii}\gamma_{ii} = \gamma_{ii}$ so that
for all $u\in U_i$ we have
$\gamma_{ii}(u) = \gamma_{ii}(u)\inv\gamma_{ii}(u)$ and
$\gamma_{ii}(u)\in S\unit$.  However, because $\gamma_{ii}$ is a
section, this implies that $\gamma_{ii}(u) = u$ for all $u\in U_i$ and  that
$\gamma_{ii}$ is just the inclusion of $U_i$ into $S$.   It follows
that $(s,i)\equiv (s,i)$ for all $(s,i)\in C$.  Since
$\gamma_{ii}$ is inclusion we have 
$\gamma_{ij}(u)\gamma_{ji}(u) = \gamma_{ii}(u) = u$.  
This implies $\gamma_{ij}(u)\inv = \gamma_{ji}(u)$ for all $u\in U_{ij}$ 
and this is enough to show that $\equiv$ is symmetric.  Finally, it
is easy enough to use the cocycle
identity \eqref{eq:3} to show that $\equiv$ is transitive.  Let
$X_\gamma$ be the quotient of $C$ by $\equiv$ with
equivalence classes denoted by $[s,i]$ for $(s,i)\in C$
and associated quotient map $Q$.  Since
$\tilde{p}$ is constant on $[s,i]$, we can factor $\tilde{p}$ through $Q$
to obtain a continuous
surjection $q:X_\gamma\rightarrow S\unit$ such that $q([s,i])=p(s)$.

Now, it is clear that $Q(p\inv(U_i)) \subset q\inv(U_i)$.  If
$[s,j]\in q\inv(U_i)$ then $u=p(s)\in U_{ij}$ and $(s,j)\equiv
(\gamma_{ji}(u)s, i)$.  However, $(\gamma_{ji}(u)s,i)\in p\inv(U_i)$
and therefore $Q(p\inv(U_i))=q\inv(U_i)$.    
Next, observe that if $(s,i)\equiv (t,i)$ then $s =
\gamma_{ii}(p(s)) t = t$.  Since the equivalence relation is trivial
on $p\inv(U_i)$ the restriction $Q|_{p\inv(U_i)} : p\inv(U_i)
\rightarrow q\inv(U_i)$ is a continuous bijection, which we will
denote $Q_i$.  Suppose $O$ is
open in $p\inv(U_i)$.  In order to show that $Q_i(O)$ is open we must
show $Q\inv(Q_i(O)) = Q\inv(Q(O))$ is open.   Suppose
$(t_l,j_l)\rightarrow (t,j)$ in $C$ and $(t,j)\in Q\inv(Q(O))$.  Since
$C$ is a disjoint union, $j_l = j$ eventually.  So, disregarding some
initial segment, we can assume $j_l = j$ for all $j$.  
Now $(t,j)\in Q\inv(Q(O))$
so there exits $(s,i)\in O$ such that $(t,j)\equiv (s,i)$.  In
particular this implies that $u = p(s)=p(t) \in U_{ij}$.  Let
$u_l = p(t_l)$ and observe that $u_l \rightarrow u$ so that eventually
$u_l \in U_{ij}$.  As before we can assume without loss of generality 
that this is always
true.  Since $u_l \in U_{ij}$ we can define $s_l = \gamma_{ij}(u_l)
t_l$.  Because $\gamma_{ij}(u_l)\rightarrow \gamma_{ij}(u)$ and
$t_l\rightarrow t$ we have $s_l\rightarrow s$ so that eventually
$s_l\in O$.  Then $(s_l,i)\equiv (t_l,j)$ implies that,
eventually, $(t_l,j)\in Q\inv(Q(O))$ so that $Q\inv(Q(O))$ is open.  
Hence $Q_i$ is a homeomorphism from $p\inv (U_i)$, which we now view as a subset
of $S$, onto $q\inv (U_i)$.  

First, note that because $X_\gamma$ is locally homeomorphic to $S$ it
is straightforward to show that $X_\gamma$ is locally compact
Hausdorff, and that we can view $X_\gamma$ as a bundle over $S\unit$ with
bundle map $q$.  We will show that $X_\gamma$ is a principal
$S$-bundle.  Let $\phi_i = Q_i\inv$ and observe that since 
$q\circ Q = p$ 
we have $q = p\circ \phi_i$.  Furthermore, given $s\in p\inv(U_{ij})$
\begin{align*}
\phi_i \circ \phi_j\inv(s) &= \phi_i([s,j]) \\
& = Q_i\inv([\gamma_{ij}(p(s))s, i]) = \gamma_{ij}(p(s)) s
\end{align*}
where the second equality holds because, by definition,
$(s,j) \equiv (\gamma_{ij}(p(s)) s,i)$.  Therefore, the $\phi_i$ define
local trivializations of $X_\gamma$ and the transition maps are the
$\gamma_{ij}$.  Hence $X_\gamma$ is a principal $S$-bundle with
trivialization $(\mcal{U},\phi,\gamma)$.
Furthermore, it is clear from our construction that the cohomology
class associated to $X_\gamma$ is $[X_\gamma] = [\gamma] = c$.  

Next, we must see that our map is well defined in the sense that if we
choose two different realizations of $c$ we end up with isomorphic
principal bundles.  Let $\eta = \{\eta_{ij}\}$ be some other cocycle
which implements $c$ on an open cover $\mcal{V}$.  Since
$[\eta]=[\gamma]=c$ we can pass to some common refinement of
$\mcal{U}$ and $\mcal{V}$, say $\mcal{W}$ with refining maps $r$ and
$\rho$ respectively, and find continuous
sections $\beta_i\in \Gamma(W_i,S)$ such that 
\[
\eta_{\rho(i)\rho(j)} \beta_j = \beta_i \gamma_{r(i)r(j)}.
\]
We define $\Omega : X_\gamma\rightarrow X_\eta$ locally by
$\Omega([s,r(i)]) = [\beta_i(p(s)) s, \rho(i)]$.  This is well defined
because if $(s,r(i))\equiv (t,r(j))$ then $p(s) = p(t) = u$ and 
\[
\beta_i(u)s = \beta_i(u) \gamma_{r(i)r(j)}(u) t = \eta_{\rho(i)\rho(j)}
\beta_j(u)t
\]
so that $(\beta_i(u)s,\rho(i)) \equiv (\beta_j(u)t,\rho(j))$.  Since each
$\beta_i$ is continuous, $\Omega$ is locally continuous and therefore
continuous.  It's easy to see that we can construct a
continuous inverse for $\Omega$ by using $\beta_i\inv$ and that
$\Omega$ is a homeomorphism.  Furthermore, it is clear that $\Omega$
preserves fibres.  Finally, given $s\in p\inv(W_i)$ we have
\[
\psi_{\rho(i)} \circ \Omega \circ \phi_{r(i)}\inv (s) = 
\psi_{\rho(i)} \circ \Omega ([s,r(i)]) = 
\psi_{\rho(i)}([\beta_i(p(s))s,\rho(i)]) = \beta_i(p(s))s.
\]
It follows that $(\mcal{W},\Omega,\beta)$ is an isomorphism from $X_\gamma$ to
$X_\eta$.  Thus we have constructed a well defined map
$[\gamma]\mapsto X_\gamma$ from
$H^1(S)$ into the set of isomorphism classes of principal
$S$-bundles.  Furthermore, it is clear that this map is a right
inverse for $X\mapsto [X]$.  

The last thing we must show is that this is a left inverse.  Suppose
$X$ is a principal $S$-bundle with trivialization
$(\mcal{U},\phi,\gamma)$ and $X_\gamma$ is constructed as above.
Define $\Omega: X\rightarrow X_\gamma$ by $\Omega(x) = [\phi_i(x),i]$
when $q(x)\in U_i$.  Suppose we have $q(x)\in U_j$ as well.  Then
$\phi_i\inv(\phi_j(s)) = \gamma_{ij}(p(s)) s$ for all $s\in
p\inv(U_{ij})$ so that for all $x\in q\inv(U_{ij})$ we get
\[
\phi_i(x) = \gamma_{ij}(q(x)) \phi_j(x).
\]
However, this shows that $(\phi_i(x),i)\equiv (\phi_j(x),j)$ and that
$\Omega$ is well defined.  Furthermore, it is easy to see by construction
that $\Omega$ preserves the fibres.  Next, if $x_l\rightarrow x$ and $q(x)\in
U_i$ then, eventually, $q(x_l)\in U_i$ and
$\Omega(x_l)=[\phi_i(x_l),i]$.  Now,  $\phi_i(x_l)\rightarrow
\phi_i(x)$ and it follows then that $\Omega(x_l)\rightarrow \Omega(x)$
so that $\Omega$ is continuous.  Furthermore, it is easy to see that we
can construct a continuous inverse via the map $[s,i]\mapsto
\phi_i\inv(s)$ so that $\Omega$ is a homeomorphism.  Now let $\psi_i$
be the local trivializations for $X_\gamma$ and observe that for $s\in
p\inv(U_i)$
\[
\psi_i\circ\Omega\circ \phi_i\inv(s) = \psi_i([s,i]) = s.
\]
Therefore letting $\iota_i$ be inclusion turns $(\mcal{U},\Omega,\iota)$ into an
isomorphism of principal $S$-bundles.  It follows that $X$ and
$X_\gamma$ have the same isomorphism class and that we have
constructed a bijection between these classes and $H^1(S)$.  
\end{proof}

\begin{remark}
Given an abelian locally compact Hausdorff group bundle $S$ we
can view $S$ as a principal $S$-bundle with respect to any cover
$\mcal{U}$ by letting $\phi$ be the identity map and $\gamma_{ij}$ be
the inclusion of $U_{ij}$ into $S$.  Since the
trivial class in $H^1(S)$ is realized on any cover $\mcal{U}$ by the
inclusion maps of $U_{ij}$ into $S$ we see that the identity is the
cohomological invariant associated to the class of trivial
$S$-bundles. 
\end{remark}

We continue our exploration of principal $S$-bundles by showing that
they are equivalent to a certain class of principal $S$-spaces.  

\begin{prop}
\label{prop:27}
Suppose $X$ is a principal $S$-bundle with trivialization
$(\mcal{U},\phi,\gamma)$.  Define the source map on $X$ to be its bundle
map.  Then, for $s\in S$, $x\in X$ such that $p(s)=q(x)\in U_i$, 
\[
s\cdot x = \phi_i\inv(s\phi_i(x))
\]
defines a continuous action of $S$ on $X$.  Furthermore the following
hold:
\begin{enumerate}
\item The action of $S$ on $X$ is principal.
\item For all $i$ the map $\phi_i$ is equivariant with respect to this
  action and the action of $S$ on itself by left multiplication. 
\item The action of $S$ on $X$ is orbit transitive.  
\end{enumerate}
\end{prop}

\begin{proof}
First, we need to make sure the action is well defined on overlaps.
Suppose $u = q(x) = p(s)$ and $u\in U_{ij}$. Then, using the fact that
$S$ is abelian,  
\begin{align*}
\phi_i\inv(s\phi_i(x)) &=
\phi_j\inv\circ\phi_j\circ\phi_i\inv(s\phi_i(x)) = 
\phi_j\inv(\gamma_{ji}(u)s\phi_i(x)) \\
&= \phi_j\inv(s \phi_j(\phi_i\inv(\phi_i(x)))) =
\phi_j\inv(s\phi_j(x)).
\end{align*}
Hence the action is well defined.  Now suppose $s,t\in S$ and
$x\in X$ such that $p(s)=p(t)=q(x)=u\in U_i$.  Then 
\begin{align*}
s\cdot(t\cdot x) &= s\cdot (\phi_i\inv(t\phi_i(x))) = 
\phi_i\inv(s\phi_i(\phi_i\inv(t\phi_i(x)))) \\
&= \phi_i\inv(st\phi_i(x)) = st\cdot x.
\end{align*}
It is also easy to see that $q(x)\cdot x = x$ for all $x\in X$.  Next,
suppose $s_l\rightarrow s$ and $x_l\rightarrow x$ such that
$p(s_l)=q(x_l)=u_l$ for all $l$ and $p(s)=q(x)=u\in U_i$.  Eventually
$u_l\in U_i$ and over $U_i$ the action is clearly continuous.
Therefore the action of $S$ on $X$ is a continuous groupoid action.  

Part {\bf (a)}: Suppose $s\cdot x = x$ for $s\in S$ and $x\in X$.  
Then for some $i$ we have
$\phi_i\inv(s\phi_i(x)) = x$ so that $s \phi_i(x) = \phi_i(x)$.  It
follows that $s\in S\unit$ and that the action is free.  Now suppose
$x_l$ and $s_l$ are nets in $X$ and $S$ respectively so that
$x_l\rightarrow x$ and $s_l\cdot x_l\rightarrow y$.  We can pass to a
subnet and assume that $p(s)=q(x)=q(y)=u \in U_i$ and 
$p(s_l)=q(x_l)=u_l\in U_i$ for all $l$.  In this case $s_l\phi_i(x_l)
\rightarrow \phi_i(y)$ and, combining this with the fact that
$\phi_i(x_l)\rightarrow \phi_i(x)$, we have $s_l\rightarrow
\phi_i(y)\phi_i(x)\inv$.  It follows that the action of $S$ on $X$ is
proper, and therefore principal.  

Part {\bf (b)}:  Suppose $s\in S$ and $x\in X$ such that
$p(s)=q(x)=u\in U_i$.  Then 
\[
\phi_i(s\cdot x) = \phi_i(\phi_i\inv(s\phi_i(x))) = s\phi_i(x)
\]
so that $\phi$ is equivariant with respect to the action of $S$ on $X$
and the action of $S$ on itself by left multiplication.  

Part {\bf (c)}:  Suppose $x,y\in X$ such that $q(x) = q(y)=u$.  We need
to find $s\in S$ such that $s\cdot x = y$.  Choose $U_i$ so that $u\in
U_i$ and let $s = \phi_i(y)\phi_i(x)\inv$.  Then we are done since
\[
s\cdot x = \phi_i\inv(\phi_i(y)\phi_i(x)\inv\phi_i(x)) = y. \qedhere
\]
\end{proof}

This next proposition shows that we can view principal $S$-bundles as
particularly nice $S$-spaces.  From now on we will think of principal
$S$-bundles in this manner.  

\begin{theorem}
\index{New Result}
\label{thm:principalspace}
\index{principal S-bundle@principal $S$-bundle}
\index{G-space@$G$-space!principal}
Suppose $S$ is an abelian locally compact Hausdorff group bundle and
$X$ is a locally compact Hausdorff space.
Then $X$ is a principal $S$-bundle if and only if $X$ is a principal,
orbit transitive, $S$-space such that the range map on $X$
has local sections. 
\end{theorem}

\begin{proof}
If $X$ is a principal $X$-bundle then let $(\mcal{U},\phi,\gamma)$ be
a trivialization of $X$ and let $S$ act on $X$ as in Proposition
\ref{prop:27}.  All we need to show is that $q$ has local sections. 
On $U_i$ define $\sigma_i: U_i\rightarrow X$
by $\sigma_i(u) = \phi_i\inv(u)$.  It is easy to see that $\sigma_i$
is a continuous section of $q$ on $U_i$ and we are done.  

Next, suppose $S$ acts on $X$ as in the statement of the theorem and
that $\mcal{U}$ is an open cover of $S\unit$ such that there
are local sections $\sigma_i:U_i\rightarrow X$ of $q$.  We define
$\psi_i:p\inv(U_i)\rightarrow q\inv(U_i)$ by $\psi_i(s) = s\cdot
\sigma_i(p(s))$.  Since everything in sight is continuous, it is clear
that $\psi_i$ is continuous.  Furthermore, using Proposition \ref{prop:12},
\[
q(\psi_i(s)) = q(s\cdot \sigma_i(p(s))) = q(\sigma_i(p(s))) = p(s).
\]
Now, if $\psi_i(s) = \psi_i(t)$ then we have $p(s) = p(t)$ and, after
multiplying by $t\inv$,
\[
t\inv s\cdot \sigma_i(p(s)) = \sigma_i(p(s)).
\]
Since the action is free this implies $t\inv s\in S\unit$ and $t=s$ so
$\phi_i$ is injective.  Next, if $y \in q\inv(U_i)$ then $q(y) =
q(\sigma_i(q(y)))$ so that, by orbit transitivity, there exists $s\in
S$ such that $y = s\cdot \sigma_i(q(y))$.  It is clear that $\psi_i(s)
= y$ and that $\psi_i$ is surjective.  Now suppose
$\psi_i(s_l)\rightarrow \psi_i(s)$.  Pass to a subnet.  We will show
that there is a sub-subnet such that $s_l\rightarrow s$.  By
definition we have $s_l \cdot \sigma_i(p(s_l)) \rightarrow s\cdot
\sigma_i(p(s))$.  Furthermore $q(s_l\cdot \sigma_i(p(s_l))) = p(s_l)$
for all $l$ and $q(s\cdot \sigma_i(p(s)))=p(s)$.  Since $q$ is continuous
we have $p(s_l)\rightarrow p(s)$ and therefore
$\sigma_i(p(s_l))\rightarrow \sigma_i(p(s))$.  Since the action of $S$
is proper this implies that we can pass to a subnet, relabel, and find
$t$ such that $s_l\rightarrow t$.  However, using the continuity of
the action, this implies $s_l\cdot\sigma_i(p(s_l))\rightarrow t\cdot
\sigma_i(p(s))$.  Using the fact that $X$ is Hausdorff and the action
is free we have $s=t$ and we are done.  Therefore $\psi_i$ is a
homeomorphism and we define the trivializing maps to be $\phi_i =
\psi_i\inv$.  

Next, we need to compute the transition functions.  Suppose $s\in
p\inv(U_{ij})$.  Then 
\begin{align*}
\phi_i\circ \phi_j\inv(s) &= \psi_i\inv\circ \psi_j(s) 
= \psi_i\inv(s\cdot\sigma_j(p(s))) \\
&= \psi_i\inv(\gamma_{ij}(p(s))s\cdot \sigma_i(p(s))) =
\gamma_{ij}(p(s))s
\end{align*}
where $\gamma_{ij}(u)$ is the unique element of $S$ such that 
\[
\gamma_{ij}(u) \cdot \sigma_i(u) = \sigma_j(u).
\]
We know $\gamma_{ij}(u)$ is 
guaranteed to exist because the action is orbit transitive and
that $\gamma_{ij}(u)$ 
is unique because the action is free.  Since $\gamma_{ij}$ is clearly
a section of $S$ on $U_{ij}$ and since it satisfies the right
algebraic properties, all that is left is to show it is continuous.
However, if $u_l\rightarrow u$ in $U_{ij}$ then
$\sigma_i(u_l)\rightarrow \sigma_i(u)$ and 
\[
\sigma_j(u_l) = \gamma_{ij}(u_l)\cdot \sigma_i(u_l)\rightarrow 
\sigma_j(u) = \gamma_{ij}(u) \cdot \sigma_i(u).
\]
Using properness as we did before, it is easy to show that
$\gamma_{ij}(u_l)\rightarrow \gamma_{ij}(u)$.  It follows that $X$ is
a principal $S$-bundle with trivialization $(\mcal{U},\phi,\gamma)$.  
\end{proof}

\begin{remark}
As in the classical case, trivial principal
$S$-bundles are exactly principal, orbit transitive
$S$-spaces whose range maps have a {\em global} section.
\end{remark}

This next proposition is nice because it frees our idea of principal
bundle isomorphism from the hassle of having to keep track of local
trivializations.  It is also mildly remarkable that $\Omega$ is not
required to be a homeomorphism, or even a bijection.  

\begin{prop}
Suppose $X$ and $Y$ are principal $S$-bundles. Then $X$ and $Y$ are
isomorphic if and only if there exists a continuous map $\Omega: X\rightarrow
Y$ which is $S$-equivariant with respect to the actions of
$S$ on $X$ and $Y$. 
\end{prop}

\begin{proof}
Let $X$ and $Y$ be as above with bundle maps $q$ and $r$, and
trivializations $(\mcal{U},\phi,\gamma)$ and $(\mcal{V},\psi,\eta)$
respectively.
Suppose $(\mcal{W},\Omega,\beta)$ is a 
principal bundle isomorphism from $X$ to $Y$.
Definition \ref{def:18} requires that $q = r\circ\Omega$ so that
$\Omega$ respects the range maps on $X$ and $Y$.  It is obvious that
$\Omega$ is continuous.  Now if $s\in S$ and $x\in X$ such that $p(s)
= q(x) \in U_i$ then 
\begin{align*}
\Omega(s\cdot x) &= \psi_i\inv \circ \psi_i\circ \Omega
\circ\phi_i\inv(s\phi_i(x)) \\
&= \psi_i\inv(\beta_i(u)s\phi_i(x)) \\
&= \psi_i\inv(s\psi_i\circ\Omega\circ\phi_i\inv(\phi_i(x))) \\
&= \psi_i\inv(s\psi_i(\Omega(x))) = s\cdot \Omega(x).
\end{align*}
Hence $\Omega$ is a continuous equivariant map.  

Now suppose $\Omega:X\rightarrow Y$ is a continuous equivariant map.
Since the range maps of $X$ and $Y$ are precisely their bundle maps,
it is part of Definition \ref{def:19}
that $\Omega$ preserves the fibres.  Let $(\mcal{U},\phi,\gamma)$ and
$(\mcal{V},\phi,\eta)$ be trivializations for $X$ and $Y$
respectively.  By passing to a common refinement we may assume without
loss of generality that $\mcal{U}=\mcal{V}$.  Given $U_i$ let $\Omega_i
= \phi_i \circ \Omega \circ \phi_i\inv$.  Since each of its component
maps preserves fibres $\Omega_i$ does as well, and therefore
$\Omega_i|_{S_u}$ maps $S_u$ into $S_u$ for $u\in U_i$.  Suppose
$s,t\in S_u$, then 
\begin{align*}
\Omega_i(st)&=
\psi_i\circ\Omega\circ\phi_i\inv(s\phi_i(\phi_i\inv(t))) 
= \psi_i\circ\Omega(s\cdot\phi_i\inv(t)) \\
&= \psi_i(s\cdot \Omega(\phi_i\inv(t)) 
= \psi_i(\psi_i\inv(s\cdot \psi_i\circ\Omega\circ\phi_i\inv(t))) \\
&= s\Omega_i(t).
\end{align*}

Now, the following general nonsense implies any map $h$ from a
group $H$ into itself such that $h$ is equivariant with respect to the
action of $H$ on itself by left multiplication is  given by left
multiplication.  That is, given map $h:H\rightarrow H$ such that 
$h(st)=sh(t)$ for all
$s,t\in H$, we have 
\[
h(s) = h(es) = h(e)s.
\]
Hence, $h$ is actually just left multiplication by $h(e)$.  Applying
this to the current situation we find that $\Omega_i|_{S_u}$ is given
on $S_u$ by left multiplication by $\Omega_i(u)$.  Define $\beta_i$ on
$U_i$ by $\beta_i(u) = \Omega_i(u)$.  The function $\beta_i$ is a section
of $S$ on $U_i$ which is continuous because $\Omega_i$ is.  Since
$\Omega_i$ is defined by left multiplication of the continuous section
$\beta_i$ it follows immediately that $\Omega_i$ has a continuous
inverse given by left multiplication by $\beta_i\inv$.  Thus
$\Omega_i$ is a homeomorphism.  But then $\Omega|_{q\inv(U_i)} =
\psi_i\inv\circ \Omega_i \circ \phi_i$ is a homeomorphism for all $i$.  It is
straightforward to show that this implies that $\Omega$ is a
homeomorphism.  Furthermore, we know that for $s\in q\inv(U_i)$ 
\[
\psi_i \circ\Omega\circ\phi_i\inv(s) = \Omega_i(s) = \beta_i(p(s)) s.
\]
Hence $(\mcal{U},\Omega,\beta)$ is a principal bundle isomorphism of $X$ onto $Y$.
\end{proof}

It is philosophically important to see that the theory of principal
$S$-bundles is an extension of the classical theory of principal group bundles.
Since we will not use the classical theory of principal group bundles,
we will not reproduce those definitions here and will instead refer the
reader to \cite[Section 4.2]{tfb}.  As with Theorem \ref{prop:principcohom},
the following will draw heavily from \cite{tfb}.  

\begin{prop}
\label{prop:91}
\index{principal S-bundle@principal $S$-bundle}
Suppose $H$ is an abelian locally compact Hausdorff group and $X$ and $Y$ are
locally compact Hausdorff spaces.  Let $S = Y\times H$ be the trivial group
bundle. Then
$X$ is a principal $H$-bundle over $Y$ if and only if $X$ is a principal
$S$-bundle.
\end{prop}

\begin{proof}
According to \cite{tfb} a principal $H$-bundle over $Y$ is just a
fibre bundle $X$ with fibres $H$ and structure group $H$ where $H$
acts on itself by left multiplication.  In particular, if
$q:X\rightarrow Y$ is a principal $H$-bundle then there exists an open
cover $\mcal{U}$, homeomorphisms  $\phi_i : q\inv(U_i)\rightarrow
U_i\times H$ such that $\phi_i(x)$ has the form $(q(x),s)$ for
some $s\in H$, and continuous maps $\gamma_{ij}:U_{ij}\rightarrow H$ such that 
\begin{equation}
\label{eq:4}
\phi_i \circ \phi_j\inv(x,s) = (x,\gamma_{ij}(x) s)
\end{equation}
for all $x\in U_{ij}$ and $s\in H$.  However, if we view $S=Y\times H$
as a group bundle over $Y$ with bundle map $p$ then $q\inv(U_i) =
U_i\times H$ and $\phi_i$ is nothing more than a bundle homeomorphism
from $q\inv(U_i)$ onto $p\inv(U_i)$.  It is then clear from
\eqref{eq:4} that the $\phi_i$ form trivializing maps and the maps
$\tilde{\gamma}_{ij}$ defined in the natural way by
$\tilde{\gamma}_{ij}(y) = (y,\gamma_{ij}(y))$ form transition
maps so that $X$ is a principal $S$-bundle with trivialization
$(\mcal{U},\phi,\tilde{\gamma})$.  This exact same argument holds in reverse,
and it's straightforward to show that if $X$ is a principal $S$-bundle
then the trivialization $(\mcal{U},\phi,\gamma)$ makes $X$ into a
principal $H$-bundle.  
\end{proof}

In particular, it's nice to observe that the sheaf cohomology theory
associated to principal $S$-bundles is just an extension of the group
sheaf cohomology theory associated to principal $H$-bundles. 

\begin{prop}
\label{prop:29}
\index{cohomology}
Suppose $H$ is an abelian locally compact Hausdorff group, $Y$ a
locally compact Hausdorff space, and $S$ is the trivial bundle
$Y\times H$.  Then $H^n(S)\cong H^n(Y;H)$ and under this
identification $X$ generates the same cohomology class when viewed as
either a principal $S$-bundle or a principal $H$-bundle.
\end{prop}

\begin{proof}
By Proposition \ref{prop:26} we know $H^n(S)$ is the sheaf cohomology
generated by the sheaf $\mcal{S}$ where $S(U) =\Gamma(U,S)$ is the set
of continuous sections on $U$ for a given open set $U\in Y$.  The
sheaf cohomology of $Y$ with coefficients $H$, denoted $H^n(Y;H)$, is
the sheaf cohomology associated to the sheaf $\mcal{T}$ where
$\mcal{T}(U) = C(U,H)$ is the set of continuous functions from $U$
into $H$.  Given $f\in C(U,H)$ we can define $\tilde{f}\in
\Gamma(U,S)$ by $\tilde{f}(y) = (y,f(y))$.  It's easy to see that this
map is an isomorphism of $\mcal{T}(U)$ onto $\mcal{S}(U)$ and extends
to a sheaf isomorphism of $\mcal{T}$ onto $\mcal{S}$.  This allows us
to identify the cohomology groups $H^n(S)=H^n(Y,\mcal{S})$ and
$H^n(Y;H)=H^n(Y;\mcal{T})$ via $[\gamma]\mapsto [\tilde{\gamma}]$.  

Now suppose $q:X\rightarrow Y$ is a topological bundle.  Proposition
\ref{prop:29} tells us that $X$ is a principal $H$-bundle if and only
if it's a principal $Y$-bundle.  Furthermore, it follows from the proof
of Proposition \ref{prop:29} that if $\gamma_{ij}$ are the transition
maps for $X$ as a principal $H$-bundle then $\tilde{\gamma}_{ij}$ are
the transition maps for $X$ as a principal $S$-bundle.  We showed in
Theorem \ref{prop:principcohom} that the cohomology invariant associated to
$X$ as a principal $S$-bundle is $[\tilde{\gamma}_{ij}]$.  It is shown
in \cite[Proposition 4.53]{tfb} that the cohomology invariant
associated to $X$ as a principal $H$-bundle is $[\gamma_{ij}]$.  Since
these two classes are precisely the ones identified under our
isomorphism, we are done. 
\end{proof}

\subsection{Locally $\sigma$-trivial Spaces}

Principal ``group bundle bundle'' theory is a natural extension of
classical principal group bundle theory.  
The real question is if there are principal
$S$-bundles which are not generated by principal $H$-bundles.
Fortunately, principal $S$-bundles are also an extension of the
notion of $\sigma$-trivial spaces as defined in \cite{locunitarystab}
and there are nontrivial examples given there.  We will include these
examples as part of this section for completeness.  However, first we
have to define $\sigma$-trivial spaces.  

\begin{definition}
\label{def:20}
\index{locally sigma trivial@locally $\sigma$-trivial}
Suppose the abelian locally compact Hausdorff group $H$ acts on the
locally compact Hausdorff space $X$ and the stabilizers vary
continuously in $H$.  We
shall say that $X$ is a locally $\sigma$-trivial space if $X/H$ is
Hausdorff and if every $x\in X$ has a $G$-invariant neighborhood $U$
such that there exists $\phi:U\rightarrow (U/H\times H)/\cong$ where 
\[
(H\cdot x,s)\cong(H\cdot y,t)\ \text{if and only if $H\cdot x=H\cdot y$
  and $st\inv\in H_x.$}
\]
Furthermore we require that  
\begin{enumerate}
\item If $x\in U$ then $\phi(x) = [H\cdot x, s]$ for some $s\in H$
  and, 
\item If $x\in U$ and $s\in H$ and $\phi(x) = [H\cdot x, t]$ then 
\[
\phi(s\cdot x) = [H\cdot x, st].
\]
\end{enumerate}
\end{definition}

Our goal will be to construct a continuously varying abelian group
bundle $S$ associated
to $G$ and $X$.  We start with a useful lemma.  

\begin{lemma}
\label{lem:1}
Suppose $H$, $X$ and $\cong$ are as in Definition \ref{def:20}.  Then
$\cong$ is an equivalence relation, the orbit space $S=(X/H\times
H)/\cong$ is locally compact Hausdorff, and the quotient map from
$X/H\times H$ onto $S$ is open. 
\end{lemma}

\begin{proof}
First, given $s\in H$ and $x\in X$ suppose $t\in H_x$.  Since $H$ is
abelian $t\cdot(s\cdot x) = st\cdot x = s\cdot x$.  Hence $t\in
H_{s\cdot x}$ and $H_x\subset H_{s\cdot x}$.  We can be tricky and
apply this argument to $s\cdot x$ and $s\inv$ to see that the
inclusion is actually an equality.  

Now, we must show that $\cong$ is an equivalence relation.  It is
clear that $(H\cdot x,s)\cong (H\cdot x,s)$.  Suppose $(H\cdot x,
s)\cong (H\cdot y,t)$.  Then, since $H\cdot x = H\cdot y$, 
by the previous paragraph, $H_x = H_y$
and since $H_y$ is a group $st\inv\in H_y$ implies $ts\inv\in H_y$.
Hence $(H\cdot y, t)\cong (H\cdot x, s)$.  Finally, suppose $(H\cdot
x, s) \cong (H\cdot y,t)$ and $(H\cdot y, t)\cong (H\cdot z, r)$.
Then, $H\cdot x = H\cdot y = H\cdot z$, and therefore $H_x =H_y
=H_z$.  By assumption $st\inv$ and $tr\inv$ are elements of $H_x$, but
of course this implies $sr\inv\in H_x$ and we are done. 

Let $Q:X/H\times H\rightarrow S$ be the quotient map and denote
elements of $S$ by $[H\cdot x, s]$.  Suppose $U$ is open in
$X/H$ and $V$ is open in $H$.  We must show that ${O = Q\inv(Q(U\times
V))}$ is open.  Notice that a generic element of $O$ has the form
$(H\cdot x, st)$ where $H\cdot x \in U$, $s\in V$, and $t\in H_x$.
Suppose to the contrary that $O$ is not open.  This implies that the
complement of $O$ is not closed so that we can find a net such that
$(H\cdot x_i, s_i)\not\in O$ for all $i$, $(H\cdot x, st)\in O$, and 
such that $(H\cdot x_i, s_i)\rightarrow (H\cdot x, st)\in O$.  In
particular, this implies that $H\cdot x_i\rightarrow H\cdot x$.  Since
the quotient map from $X\rightarrow X/H$ is always open we can pass to
a subnet, reindex, and possibly choose new representatives, so that
$x_i\rightarrow x$.  Since the stabilizers of $H$ vary continuously
with respect to the Fell topology this implies $H_{x_i}\rightarrow
H_x$.  However, $t\in H_x$ so, by Proposition \ref{prop:8}, we can
pass to another subnet, reindex, and find $t_i\in H_{x_i}$ such that
$t_i\rightarrow t$.  Since $s_i\rightarrow st$ it follows that
$s_it_i\inv \rightarrow s$ and since $s\in O$ this implies
$s_it_i\in O$ eventually.  It follows that, eventually,
$(H_{x_i},s_i) \in O$, which is a contradiction.  Therefore $O$,
and hence $Q$, must be open. 

Since $Q$ is open and $X/H\times H$ is locally compact, it follows that
$S$ is locally compact.  All that's left is to show that it is
Hausdorff.  Suppose $[H\cdot x_i, s_i]$ is a net in $S$
which converges to both $[H\cdot x,s]$ and $[H\cdot y, t]$.  Using the
fact that $Q$ is open we can pass to a subnet, twice actually,
reindex, choose new representatives $s_i$, and find $t_i\in H_{x_i}$
such that 
\begin{align*}
H\cdot x_i &\rightarrow H\cdot x  & s_i &\rightarrow s \\
H\cdot x_i &\rightarrow H\cdot y  & s_it_i & \rightarrow t.
\end{align*}

Now, $X/H$ is assumed to be Hausdorff so that $H\cdot x = H\cdot y$.  We
pass to yet another subnet, relabel, and choose new $x_i$ so that
$x_i\rightarrow x$.  Since the stabilizers vary continuously this
implies $H_{x_i}\rightarrow H_x$.  Observe that 
\[
s_i\inv(s_it_i) = t_i \rightarrow s\inv t
\]
and that, via Proposition \ref{prop:8} because $t_i\in H_{x_i}$ for
all $i$, we must have $s\inv t \in H_x$.  This implies that $(H\cdot
x, s) \cong (H\cdot x, t)$ and we are done. 
\end{proof}

Now that we know the quotient space $(X/H\times H)/\cong$ is well
behaved topologically we can prove a more interesting proposition. 

\begin{prop}
\label{prop:30}
\index{group bundle}
Suppose the abelian locally compact Hausdorff group $H$ acts on the
locally compact Hausdorff space $X$.  Furthermore, suppose the
stabilizers vary continuously in $H$ and that $X/H$ is Hausdorff.
Define $S_{(X,H)} := (X/H\times H)/\cong$, often denoted $S$, 
where $\cong$ is defined as in Definition \ref{def:20}.  Then $S$ is
an abelian, continuously varying, locally compact Hausdorff group
bundle whose unit space can be identified with $X/H$.  
The bundle map is given by $p([H\cdot x, s]) = H\cdot x$ and the
operations are 
\begin{align*}
[H\cdot x, s][H\cdot x, t] &:= [H\cdot x, st], &
[H\cdot x, s]\inv &:= [H\cdot x, s\inv].
\end{align*}
The fibre $S_{H\cdot x}$ over $H\cdot x$ is (isomorphic to) $H/H_x$.  
\end{prop}

\begin{proof}
Define $S$ as in the statement of the proposition and observe that it
follows from Lemma \ref{lem:1} that $S$ is locally compact Hausdorff.
Furthermore, let $Q:X/H\times H\rightarrow S$ be the quotient map and
recall that we showed that $Q$ is open.  Define 
\[
S^{(2)} = \{([H\cdot x, s],[H\cdot y,t])\in S\times S: H\cdot x =
H\cdot y\}
\]
and observe that, since all representatives of the class $[H\cdot x
,s]$ are of the form $(H\cdot x, st)$ where $t\in H_x$, $S^{(2)}$ is
well defined.  We would like to show that the operations given above
are well defined.  Suppose $(H\cdot x, s), (H\cdot x,t),(H\cdot
y,s'),(H\cdot y,t')\in X/H\times H$ and that $(H\cdot x,s)\cong
(H\cdot y,s')$ and $(H\cdot x,t)\cong (H\cdot y,t')$.  Then $H\cdot x
= H\cdot y$ and there exists $u,v\in H_x$ such that $s' = su$ and $t'
= tv$.  Since $s't' = (su)(tv) = st(uv)$ it follows that $(H\cdot
x,st)\cong (H\cdot y,s't')$ and that the multiplication is well
defined.  The proof that the inverse operation is well defined is
similar.  It is easy to use the fact that $H$ is a group to prove that
$S$, with these operations, is a groupoid.  Furthermore, the range and
source maps on $S$ are equal and are both given by 
\[
p([H\cdot x, s])= [H\cdot x, e]
\]
where $e\in H$ is the identity.  Thus $S$ is a group bundle and it is
straightforward to see that 
$S\unit$ can be identified with
$X/H$ via the map $[H\cdot x, e]\mapsto H\cdot x$.  Furthermore, under this
identification the bundle map $p$ has the required form. 

Next, we have to show that the operations are continuous.  Suppose
$[H\cdot x_i,s_i]\rightarrow [H\cdot x,s]$ and $[H\cdot
x_i,t_i]\rightarrow [H\cdot x,t]$ in $S$.  We can pass to a subnet,
twice, relabel, and choose new representatives $s_i$ and $t_i$ such that
$H\cdot x_i\rightarrow H_x$, $s_i\rightarrow s$ and
$t_i\rightarrow t$.  It follows immediately that 
\[
(H\cdot x_i, s_it_i)\rightarrow (H\cdot x, st).
\]
Showing that the inverse operation is continuous is a similar
exercise.  Suppose $H\cdot x_i\rightarrow H\cdot x$ in $X/H$ and that
$[H\cdot x,s]\in S$.  Then $(H\cdot x_i,s)\rightarrow (H\cdot x,s)$ in
$X/H\times H$.  It follows immediately that the bundle map $p$ is
open so that $S$ is a continuously varying locally compact Hausdorff
group bundle.  

Given $H\cdot x\in X/H$ we have $S_{H\cdot x} = \{[H\cdot x,s ]: s\in
H\}$.  We can define a continuous surjection $\phi:H\rightarrow S_{H\cdot x}$ by
$\phi(s) = [H\cdot x, s]$ and it is clear from the definition of the
operations on $S$ that this is a homomorphism.  Next, if 
$[H\cdot x, s_i]\rightarrow [H\cdot x, s]$ in $S_{H\cdot x}$ then we
can use the fact that $Q$ is open to pass to a subnet, relabel, and
choose new $s_i$ such that $s_i\rightarrow s$.  However, this implies
that $\phi$ is an open map.  Finally, it is clear from the
definition of $\cong$ that $\phi(s) = \phi(t)$ if and only if
$st\inv\in H_x$.  It follows that $\phi$ factors to an isomorphism of
$H/H_x$ with $S_{H\cdot x}$.  Since $H/H_x$ is clearly abelian this
proves that $S$ has abelian fibres and we are done.  
\end{proof}

The reason we went through all of this rigmarole is that given a
locally $\sigma$-trivial system $(H,X)$ we would like to show that $X$ is a
principal $S_{(H,X)}$ bundle.  

\begin{prop}
\label{prop:31}
\index{principal S-bundle@principal $S$-bundle}
\index{locally sigma trivial@locally $\sigma$-trivial}
Suppose $H$ is an abelian locally compact Hausdorff group acting on a
locally compact Hausdorff space $X$ with continuously varying
stabilizers such that $X/H$ is Hausdorff.  If $X$ is locally 
$\sigma$-trivial then
$X$ is a principal $S_{(H,X)}$-bundle. 
\end{prop}

\begin{proof}
Let $q:X\rightarrow X/H$ be the quotient map.  We know from Definition
\ref{def:20} that if $x\in X$ then there is an $H$-invariant
neighborhood $U$ that is homeomorphic to $U/H\times H/\cong$.  If we
let $V = q(U)= U/H$ then $V$ is an open neighborhood of $H\cdot x$ and
$q\inv (V) = U$. Furthermore, let $p$ be the bundle map for $S$ and
observe that 
\[
p\inv(V) = \{[H\cdot x,s]\in S:H\cdot x\in V\} =  U/H\times H/\cong. 
\]
Thus, as per Definition \ref{def:20}, there is a homeomorphism
$\phi_{V}:q\inv(V)\rightarrow p\inv(V)$.  Furthermore, since
$\phi_{V}(x) = [H\cdot x, s]$ for some $s\in H$ it is clear that $q =
p\circ \phi_{V}$.  Now, find one of these neighborhoods for every
$x\in X$ and use them to build an open cover $\mcal{V}$ of $X/H$.  

For $V_i$ in this open cover let $\phi_i= \phi_{V_i}$ and note that we
have already shown that each $\phi_i$ is a bundle map.  Given
$V_{ij}$ define $\gamma_{ij}:U_{ij}\rightarrow S$ by
$\gamma_{ij}(H\cdot x) = \phi_i\circ\phi_j\inv([H\cdot x, e])$.  It is
clear that $\gamma_{ij}$ is a continuous section on $V_{ij}$.  Since
$\gamma_{ij}$ is a section we can find a function 
$\tilde{\gamma}_{ij}$ from $V_{ij}$ into $H$ such that
$\gamma_{ij}(H\cdot x) = [H\cdot x, \tilde{\gamma}_{ij}(H\cdot x)]$.  
Suppose
$[G\cdot x, s]\in p\inv(V_{ij})$.  Then, using the equivariance
condition of Definition \ref{def:20}, we have 
\begin{align*}
\phi_i\circ\phi_j\inv([H\cdot x, s]) &= 
\phi_i(s\cdot \phi_j\inv([H\cdot x,e])) \\
&= [H\cdot x, \tilde{\gamma}_{ij}(H\cdot x)s] \\
&= \gamma_{ij}(H\cdot x) [H\cdot x,s].
\end{align*}
Thus $X$ is a principal $S$-bundle with trivialization
$(\mcal{V},\phi,\gamma)$.  
\end{proof}

\begin{remark}
\index{cohomology}
In \cite{locunitarystab} $\sigma$-trivial space is said to be 
locally liftable if every
continuous section $c:X/H\rightarrow S$ is given locally by a
continuous map $\tilde{c}:X/H \rightarrow H$ such that $c(H\cdot x) = [H\cdot x,
\tilde{c}(H\cdot x)]$.  Locally $\sigma$-trivial {\em bundles} are defined to
be locally $\sigma$-trivial spaces which are also locally liftable.
The reason for this extra requirement has to do with finding a
cohomological invariant for $X$. Let $\mcal{T}$ be the sheaf
defined for $U\subset X/H$ by $\mcal{T}(U) = C(U,H)$ and 
$\mcal{R}$ be the subsheaf of $\mcal{T}$ where $\mcal{R}(U)$ is the
subset of $\mcal{T}(U)$ such that $f(H\cdot x)\in H_x$ for all $x$.  Then sheaf
cohomological considerations will show that we can construct a
quotient sheaf $\mcal{T}/\mcal{R}$ and an associated cohomology
$H^n(X/H; \mcal{T}/\mcal{R})$.  Given a $\sigma$-trivial space one
would like to use the transition maps $\gamma_{ij}$, as defined in
the proof of Proposition \ref{prop:31}, to construct an element of
$H^1(X/H;\mcal{T}/\mcal{R})$.  The problem is that while $\gamma_{ij}$
is a continuous section of $U_{ij}$ into $S$ the associated map
$\tilde{\gamma}_{ij}:U_{ij}\rightarrow H$ may not be continuous.  If
$\tilde{\gamma}_{ij}$ is not continuous then it doesn't define an
element of $\mcal{T}(U_{ij})$ and we cannot construct the appropriate
cohomology element.  However, if $X$ is required to be locally
liftable then, by passing to a smaller open set, we can guarantee that
$\tilde{\gamma}_{ij}$ is a continuous function.  As such it defines an
element of $\mcal{T}(U_{ij})$ and hence a cohomology element in
$H^1(X/H;\mcal{T}/\mcal{R})$.  In fact, it is shown in
\cite{locunitarystab} that this construction leads to a one-to-one
correspondence between locally $\sigma$-trivial {\em bundles} with a
fixed orbit space $X/H$ and $H^1(X/H;\mcal{T}/\mcal{R})$.  

This is an artificial restriction in our setting.  The
$\gamma_{ij}$ can always be used to define an element of
$H^1(S_{(X,H)})$, regardless of whether $\sigma$ is locally liftable or
not.  It is comforting to observe the following however.  Let
$\mcal{S}$ be the sheaf of local sections of $S$ so that $H^n(S) =
H^n(X/H;\mcal{S})$ by definition.  It is straightforward to show that if
$\sigma$ is locally liftable then we get a short exact sequence of
sheaves 
\[
0\rightarrow \mcal{R}\rightarrow \mcal{T}\rightarrow
\mcal{S}\rightarrow 0
\]
and that $H^n(X/H;\mcal{S})$ is naturally isomorphic to
$H^n(X/H,\mcal{T}/\mcal{R})$.  Furthermore, once one sorts out all of
the various constructions, it is clear that the different
cohomological invariants of a $\sigma$-trivial bundle are identified
under this isomorphism. 

One reason for preferring the locally liftable case is that the
cohomology group $H^1(X/H;\mcal{T}/\mcal{R})$ is easier to deal with.
This is because it is a quotient of $H^1(X/H;\mcal{T})$ and $H^1(X/H;\mcal{T})$
is just the classical sheaf cohomology of $X/H$ with coefficients in
$H$.  In general the cohomology group $H^1(S_{(X,H)})$ is much more
mysterious.  For instance, it is not immediately clear that there are group
bundles $S$ such that $H^1(S)$ is non-trivial.  However, any nontrivial
example of a locally $\sigma$-trivial space will give rise to a
nontrivial principal $S$-bundle.  This will in turn guarantee that the
cohomology group $H^1(S)$ is nontrivial.
\end{remark}

The following examples are lifted straight from \cite{locunitarystab}
and are included for completeness.  

\begin{example}[{\cite[Example 4.6]{locunitarystab}}]
\label{ex:13}
Choose a complex line bundle $p:L\rightarrow Y$, and give it a
Hermitian structure.  Then we can define an action of $H=\R$ on $L$ by 
\[
r\cdot z = \begin{cases} e^{2\pi i r/|z|} \cdot z & \text{if $|z|\ne
    0$,}\\  z & \text{if $|z|=0$.}
\end{cases}
\]
Then we have 
\[
H_z = \begin{cases} |z|\Z &\text{if $|z|\ne 0$,} \\
\R &\text{if $|z|=0$.}
\end{cases}
\]
It is fairly easy to see that the stabilizers vary continuously. 
The local triviality of $L$
as a bundle implies that there are local cross sections of $p$, so
that by \cite[Proposition 4.3]{locunitarystab} $L$ is a locally
$\sigma$-trivial space.  In general, $L$ is not globally
$\sigma$-trivial.  In fact, it is easy to see that it is globally
$\sigma$-trivial exactly when it is trivial as a line bundle.  Thus,
if $L$ is a non-trivial line bundle then $L$ is a non-trivial
principal $S$-bundle where $S = (L/\R\times \R)/\cong$.
\end{example}

\begin{example}[{\cite[Example 4.15]{locunitarystab}}]
\label{ex:26}
Let $H$ be a locally compact abelian group, $q:Y\rightarrow Z$ a
locally trivial principal $H$-bundle, and $\tau:Z\rightarrow
\mathbb{RP}^n$ a continuous map onto $n$-dimensional real projective
space.  Now define $X = Y\times \R^{n+1}/\sim$ where $(x,a)\sim (y,b)$
if and only if $x=y$ and $a-b\in \tau(x)$.  Let $K=H\times \R^{n+1}$
act on $X$ by $(s,a)\cdot (y,b):= (s\cdot y, a+b)$.  
Then it is shown in \cite[Example 4.15]{locunitarystab} that $X$ is a
locally $\sigma$-trivial bundle over $Z$ which is
globally $\sigma$-trivial if and only if $Y$ is trivial.  It follows
that $X$ is a principal $S_{(K,X)}$-bundle.  For an interesting
concrete example of such a space we can take 
$H=\Z_2$, $Y=S^n$, $Z=\mathbb{RP}^n$, $q:Y\rightarrow Z$ the canonical
map, and $\tau$ the identity map.
\end{example}


\section{Group Bundle Duality}
\label{sec:duality}

The goal of this section is to show that when 
$S$ is an abelian, continuously varying, locally compact
Hausdorff group bundle then $S$ has a Pontryagin duality theory which
parallels the duality theory of abelian locally compact Hausdorff
groups.  We start by recalling the basic facts concerning the group
case.  

\begin{remark}
\label{rem:6}
\index[not]{$\hat{f}$}
\index{Fourier transform}
If $H$ is an abelian locally compact Hausdorff group then the
Pontryagin dual of $H$, denoted $\widehat{H}$, is defined to be the set of continuous $\T$-valued homomorphisms on $H$.  Elements of $\widehat{H}$ are called
characters and $\widehat{H}$ is an abelian group under the
operations of pointwise multiplication and conjugation.  Furthermore,
the topology of uniform convergence on compact sets 
makes $\widehat{H}$ into an abelian
locally compact Hausdorff group \cite{rudinfourier}.  
Recall that for second countable
spaces, the topology of uniform convergence on compacta is
characterized by $\omega_i\rightarrow \omega$ if and only if
$\omega_i(s_i)\rightarrow \omega(s)$ for all $s_i\rightarrow s$
\cite[Lemma 1.30]{tfb2}.

We can also realize $\widehat{H}$ as the spectrum
of the group $C^*$-algebra in the following way.  Those readers
unfamiliar with group $C^*$-algebras are referred to \cite[Appendix
C.3]{tfb}.  First, observe that
$C^*(H)$ is abelian since $H$ is.  Given a character
$\omega\in \widehat{H}$ we define a function on $C_c(H)$, also denoted
$\omega$, by 
\begin{equation}
\label{eq:5}
\omega(f) = \hat{f}(\omega) := \int_H f(s)\omega(s)ds.
\end{equation}
Then $\omega$ extends to a character on $C^*(G)$.  
Furthermore this character is uniquely determined
by $\omega$ and every element of the spectrum can be obtained in this
fashion.  Whats more, the topology of uniform convergence on compacta
on $\widehat{H}$ is exactly the Gelfand topology when $\widehat{H}$ is
viewed as the spectrum of $C^*(H)$.  

Given a function $f\in C_c(H)$ the function $\hat{f}$ defined in
\eqref{eq:5} by 
\[
\omega \mapsto \omega(f) = \int_H f(s)\omega(s)ds
\]
is called the Fourier transform of $f$.  The driving result behind
harmonic analysis is that the map defined on $C_c(H)$
via $f\mapsto \hat{f}$ extends to an isomorphism from
$C^*(H)$ onto $C_0(\widehat{H})$.  Of course, this is a special case
of the usual
Gelfand-Naimark theorem for abelian $C^*$-algebras.  This discussion
has been a short version of \cite[Section 3.1]{tfb2}.
\end{remark}

The basic idea is that given an abelian group bundle $S$ we just ``bundle
together'' the duals of the fibres of $S$.  The question is what to use
as the topology on the total space.  

\begin{definition}
\label{def:22}
\index{dual bundle}
\index{group bundle}
\index[not]{$\widehat{S}$}
\index[not]{$\hat{p}$}
Suppose $S$ is an abelian, continuously varying, locally compact
Hausdorff group bundle with bundle map $p$.  We define the {\em dual
  bundle} of $S$ to be the disjoint union
\[
\widehat{S} := \coprod_{u\in S\unit} \widehat{S}_u
\]
where $\widehat{S}_u$ denotes the Pontryagin dual of $S_u$.  We define
the {\em dual bundle map} $\hat{p}:\widehat{S}\rightarrow S\unit$ to
be given for $\omega\in \widehat{S}$ by $\hat{p}(\omega)=u$ if 
$\omega\in \widehat{S}_u$.  If $\beta$ is a Haar system for $S$ then
given $\omega\in \widehat{S}$ and $f\in C_c(S)$ we define
\begin{equation}
\label{eq:6}
\omega(f) := \int_S f(s)\omega(s) d\beta^{\hat{p}(\omega)}(s).
\end{equation}
Finally, we define a topology on
$\widehat{S}$ as follows.  Given a net $\{\omega_i\}\subset \widehat{S}$ and
$\omega\in \widehat{S}$ then $\omega_i\rightarrow \omega$ if and only
if $\omega_i(f)\rightarrow \omega(f)$ for all $f\in C_c(S)$.  
\end{definition}

\begin{remark}
Defining a topology by specifying the convergent sequences can be a
subtle process.  We will show that in this case these
convergent sequences characterize the Gelfand topology on
$\widehat{S}$ as the spectrum of the groupoid $C^*$-algebra of $S$.
\end{remark}

\begin{remark}
Since $\hat{p}:\widehat{S}\rightarrow S\unit$ is a surjection we can
view $\widehat{S}$ as a bundle over $S\unit$ with fibres
$\hat{p}\inv(u)$.  However, by construction, 
$\hat{p}\inv(u)$ is clearly equal to the
dual of $S_u$ which lies inside $\widehat{S}$.  
Therefore the notation
$\widehat{S}_u$ can be used to denote both the dual of $S_u$ and the
fibre of $\widehat{S}$ over $u$ without confusion.  
\end{remark}

\begin{example}
\label{ex:34}
Before checking the details of Definition \ref{def:22} below, it is
worth pointing out that if $Z$ and $T$ are the bundles described in
Example \ref{ex:33} then $Z$ is isomorphic to $\widehat{T}$ and $T$ is
isomorphic to $\widehat{Z}$.  
\end{example}

There are a lot of assertions in Definition \ref{def:22} that deserve
to be checked and unfortunately we lack the technology to do them all
justice at the moment.  Fortunately, these issues are worked through
in \cite{ctgIII}.  We will be developing the necessary technology
later on, however, and for those readers not opposed to a little
forward referencing a proof to the following lemma is provided.
Basically, we need a little operator algebra theory to show that the
topology on $\widehat{S}$ is well behaved.  

\begin{lemma}
\label{lem:2}
Suppose $S$ is an abelian, continuously varying, second countable, 
locally compact
Hausdorff groupoid.  Then the dual bundle $\widehat{S}$ is a second
countable locally compact
Hausdorff space.  Furthermore, the
map $\hat{p}:\widehat{S}\rightarrow S\unit$ is a continuous
surjection, and $\hat{p}\inv(u)$ may be identified with the dual of
$S_u$ topologically.   
\end{lemma}

\begin{proof}
These statements are proved in the discussion at the beginning of
\cite[Section 3]{ctgIII}.   An explicit proof, that unfortunately
relies on some material in Section \ref{sec:bundleprod}, is given below. 

Since $S$ is continuously varying it has a Haar system $\beta$.
Therefore we can construct the groupoid $C^*$-algebra $C^*(S)$.
Furthermore, it is easy to see that, because $S$ is abelian, $C^*(S)$ is
abelian also.  It then follows from the Gelfand-Naimark theorem
\cite{invitation} that given the Gelfand topology the spectrum of 
$C^*(S)$ is a locally compact
Hausdorff space.  Next, Proposition \ref{prop:groupbundalg} implies
that $C^*(S)$ is a $C_0(S\unit)$-algebra and that given $u\in S\unit$ the
fibre $C^*(S)(u)$ is (isomorphic to) $C^*(S_u)$.  It is a result of
general $C_0(X)$-algebra theory, 
reproduced in Section \ref{sec:c0x-algebras} as Proposition \ref{prop:36}, 
that there is a continuous surjection
$\hat{p}:C^*(S)\sidehat\rightarrow S\unit$ such that if $\pi\in
C^*(S)\sidehat$ then $\pi$ factors to an irreducible representation of
$C^*(S)(\hat{p}(\pi))=C^*(S_{\hat{p}(\pi)})$.  
In this way we can identify $C^*(S)\sidehat$ with the disjoint union
$\coprod_{u\in S\unit} C^*(S_u)\sidehat$.  However, we know from
Remark \ref{rem:6} that $C^*(S_u)\sidehat = \widehat{S}_u$.
It follows that $\widehat{S}$ can be identified with
$C^*(S)\sidehat$ as a set and that $\hat{p}$ has exactly the right definition
under this identification.  

Now, it is also a result of Proposition \ref{prop:groupbundalg} that
the quotient map from $C^*(S)$ to $C^*(S)(u) = C^*(S_u)$ is given on
$C_c(S)$ by restriction to $S_u$.  Suppose $\omega\in \widehat{S}$ and
let $u = \hat{p}(\omega)$.  Then we can lift $\omega$ to a
representation of $C^*(S)$ and for $f\in C_c(S)$ this is given by 
\[
\omega(f) = \int_{S_u} f|_{S_u}(s)\omega(s) d\beta^u(s) = \int_S
f(s)\omega(s) d\beta^u(s).
\]
Therefore the action of $\omega$ on $C_c(S)$ defined in
\eqref{eq:6} is precisely the action of $\omega$ on $C_c(S)$ as an
element of $C^*(S)\sidehat$.  Since the Gelfand topology on
$C^*(S)\sidehat$ is characterized by pointwise convergence
\cite{invitation}, we conclude that the topology on $\widehat{S}$
defined in Definition \ref{def:22} is exactly the Gelfand topology
when $\widehat{S}$ is identified with $C^*(S)\sidehat$.  Thus
$\widehat{S}$ is locally compact Hausdorff.  Furthermore, 
it immediately follows that the map $\hat{p}$ is continuous as a function on
$\widehat{S}$, and that the topology on $\widehat{p}\inv(u)$ is the
Gelfand topology on $C^*(S_u)\sidehat$, which in turn is the topology
on the dual of $S_u$.  Finally, it follows from Corollary \ref{cor:23} that 
$C^*(S)$ is separable.  Hence $C^*(S)\sidehat$ is
second countable \cite[Proposition 3.3.4]{dixmiercstar} and we are done
\end{proof}

\begin{remark}
It follows from the proof of Proposition \ref{lem:2} that
$\widehat{S}$ can be identified with the spectrum of $C^*(S)$ so that
characters act on elements of $C_c(S)$ as in \eqref{eq:6}.  However,
we can then use the Gelfand-Naimark theorem to conclude that
the Fourier transform induces an isomorphism of $C^*(S)$ with
$C_0(\widehat{S})$.  In this way bundle duality is a generalization of
the usual Pontryagin duality for groups.  It is notable that the dual
of $S$ is often {\em defined} to be the spectrum of the
$C^*$-algebra of $S$.  
\end{remark}

The following proposition is important because it gives an alternative 
sequential characterization of the topology on $\widehat{S}$.  This
proposition, and its proof, can be found in \cite[Proposition
3.3]{ctgIII} and are only reproduced here for ease of reference. In
particular, the argument used in this proof will be used again and
again later on. 

\begin{prop}
\label{prop:33}
Suppose $S$ is an abelian, continuously varying, second countable,
locally compact Hausdorff group bundle.  A sequence
$\{\omega_i\}$ in $\widehat{S}$ converges to $\omega_0$ in
$\widehat{S}$ if and only if 
\begin{enumerate}
\item $\hat{p}(\omega_i)\rightarrow \hat{p}(\omega_0)$ in $S\unit$,
  and
\item if $s_i\in S_{\hat{p}(\omega_i)}$ for all $i\geq 0$ and
  $s_i\rightarrow s$ in $S$, then $\omega_i(s_i)\rightarrow
  \omega_0(s_0)$.  
\end{enumerate}
\end{prop}

\begin{proof}
First, suppose that $\omega_i$ converges to
$\omega_0$ and let $u_i = \hat{p}(\omega_i)$ for all $i\geq 0$.  
The continuity of $\hat{p}$ implies that $u_i\rightarrow u_0$.  If
condition (b) fails then there are $s_i\in S_{u_i}$ converging to
$s_0\in S_{u_0}$ with $\omega_i(s_i)$ not converging to
$\omega_0(s_0)$.  Clearly, we may assume that no subsequence converges
to $\omega_0(s_0)$ either.  Next, we observe that we may assume
$u_i \ne u_0$ for all $i$; otherwise we obtain an immediate
contradiction by passing to a subsequence and relabeling so that $u_i
= u_0$ for all $i$, and $\omega_i\rightarrow\omega_0$ in
$\widehat{S}_{u_0}$.  Furthermore,
again passing to a subsequence and relabeling, we can assume that
$u_i\ne u_j$ if $i\ne j$.  In particular, we can define an integer
valued function on $C = p\inv(\{u_i\}_{i\geq0})$ by $\iota(s) = i$
when $p(s) = u_i$.  Now fix $f\in C_c(S)$ with $\omega_0(f) = 1$.
Notice that $C$ is closed, and $g_0:C\rightarrow\C$, defined by
\[
g_0(t) = f(s_{\iota(t)}\inv t)\quad\text{for $t\in C$}, 
\]
is continuous and compactly supported.  The Tietze Extension Theorem
implies that there is a $g\in C_c(S)$ extending $g_0$.  But
\[
\omega_i(g) = \omega_i(s_i) \omega_i(f)
\]
for $i\geq 0$.  We obtain the desired contradiction by noting that
$\omega_i(f)\rightarrow 1$ and $\omega_i(g)\rightarrow \omega_0(s_0)$.  

Conversely, now assume that $\omega_i$ satisfies conditions (a) and
(b).  Let $u_i = \hat{p}(\omega_i)$ for all $i\geq 0$.  Suppose
there is a $f\in C_c(S)$ such that $\omega_i(f)$ fails to converge to
$\omega_0(f)$.  As above we can reduce to the case that $u_i\ne u_j$ if
$i\ne j$.  This time we define $g_0:C\rightarrow \C$ by 
\[
g_0(t) = \omega_{\iota(t)}(t)f(t)\quad \text{for $t\in C.$}
\]
Again a few moments of reflection reveal that $g_0$ is continuous and
compactly supported so that there is a $g\in C_c(S)$ extending $g_0$.
The continuity of the Haar system on $S$ implies that 
\[
\int_S g(s) d\beta^{u_i}(s) \rightarrow \int_S g(s)d\beta^{u_0}(s).
\]  
Since $\int_S gd\beta^{u_i} = \omega_i(f)$ for all $i\geq 0$ 
we obtain the necessary contradiction.
\end{proof}

This characterization of the topology on $\widehat{S}$ is so nice that
we will restrict ourselves to the second countable case for the rest
of the section.  The following proposition is, more or less, 
\cite[Corollary 3.4]{ctgIII}.  

\begin{prop}
\label{prop:34}
\index{dual bundle}
Suppose $S$ is an abelian, continuously varying, second countable
locally compact Hausdorff group bundle.  Then the dual bundle
$\widehat{S}$ is an abelian, second countable, locally compact Hausdorff group
bundle where $\widehat{S}^{(2)} = \{(\omega,\chi)\in \widehat{S}\times
\widehat{S} : \hat{p}(\omega)=\hat{p}(\chi)\}$ and the operations are
given pointwise by 
\begin{align*}
\omega\chi(s) &:= \omega(s)\chi(s) &
\omega\inv(s) &:= \overline{\omega(s)}
\end{align*}
The unit space $\widehat{S}\unit$ can be identified with
$S\unit$ and under this identification the bundle map for
$\widehat{S}$ is $\hat{p}$.  Furthermore, the fibres of $\widehat{S}$ are
the Pontryagin duals of the fibres of $S$.  
Finally, the topology on
$\widehat{S}$ is independent of the choice of Haar measure.  
\end{prop}

\begin{proof}
We proved in Lemma \ref{lem:2} that $\widehat{S}$ is a second
countable locally compact Hausdorff space.  Furthermore, it is clear
from the definition of the dual bundle that the fibres of
$\widehat{S}$ are the Pontryagin duals of the fibres of $S$ set
theoretically.  It follows from Lemma \ref{lem:2} 
that the fibres of $\widehat{S}$ can be
identified with the Pontryagin duals of the fibres of $S$
topologically.   Finally, since we have characterized the
topology on $\widehat{S}$ independently of  Haar measure in
Proposition \ref{prop:33}, it follows that the topology is independent
of Haar measure.  

Next, we have defined the operations on $\widehat{S}$ fibrewise by the
usual pointwise operations on the Pontryagin duals.  Since the dual of
a group is again a group it is easy to see that $\widehat{S}$ is a
groupoid.  We would like to see that the operations are continuous.
Suppose $\omega_i\rightarrow \omega$ in $\widehat{S}$.  Let
$s_i\rightarrow s$ such that $p(s_i) = \hat{p}(\omega_i)$ and $p(s) =
\hat{p}(\omega)$.  Then, by Proposition \ref{prop:33}
$\omega_i(s_i)\rightarrow \omega(s)$ and therefore
\[
\omega_i\inv(s_i) = \overline{\omega_i(s_i)}\rightarrow
\overline{\omega(s)} = \omega\inv(s).
\]
It follows from Proposition \ref{prop:33} that this implies
$\omega_i\inv\rightarrow \omega$.  The proof that multiplication is
continuous is similar.  

We constructed $\widehat{S}$ so that the range and source maps are
both given by $\omega\mapsto \omega\inv\omega$.  However,
$\omega\inv\omega = 1_u$ where $1_u$ denotes the trivial character on
$S_u$.  We would like to identify $S\unit$ and $\widehat{S}\unit$ via
the map $u\mapsto 1_u$. It is easy to see that this map is a
bijection.  Proposition \ref{prop:33} implies that if
$1_{u_i}\rightarrow 1_u$ then $u_i\rightarrow u$ and it is easy to see
that Proposition \ref{prop:33} also implies that $1_{u_i}\rightarrow
1_u$ if $u_i\rightarrow u$.  Therefore we can topologically identify
$S\unit$ and $\widehat{S}\unit$ and clearly the bundle map on
$\widehat{S}$ is $\hat{p}$ under this identification.  Finally,
since the dual of an abelian group is abelian it follows that
$\widehat{S}$ is an abelian group bundle. 
\end{proof}

At this point we need to recall some basic facts about abelian
harmonic analysis.  These theorems can, for the most part be found in
\cite[Chapter 1]{rudinfourier}.  Those portions of the theorems which
are not explicitly proved in \cite{rudinfourier} are proved here.  It
is notable that \cite{rudinfourier} uses the opposite conjugation
convention and that the following theorems have been modified
accordingly.  

\begin{definition}
\label{def:23}
Let $H$ be a locally compact Hausdorff abelian group.  A function $f$
defined on $H$ is said to be {\em positive definite} if the inequality
\[
\sum_{i,j=1}^N c_i\overline{c_j} f(s_n s_m\inv) \geq 0
\]
holds for every choice of $s_1,\ldots, s_N\in H$ and for every
choice of $c_1,\ldots, c_N\in \C$.  
\end{definition}

\begin{remark}
Definition \ref{def:23} is a special case of Definition \ref{def:21}
in that if we let $H$ act on the trivial space consisting of one point
then the positive definite functions are exactly those of positive
type. 
\end{remark}

\begin{example}
\label{ex:14}
Suppose $H$ is a locally compact abelian group and 
$g\in L^2(H)$ where $H$ is given Haar measure.  For $s\in H$ define 
\[
f(s) = g^* * g(s) := \int \overline{g(t\inv)} g(t\inv s) dt.
\]
It follows from \cite[1.4.2]{rudinfourier} that $f$ is continuous and
positive definite. 
\end{example}

This the premiere example of a positive definite function and are the
only positive definite functions that we will be using.  The next
theorem says that positive definite functions are all given by
integration with respect to some measure on $\widehat{H}$.  

\begin{theorem}[Bochner's Theorem {\cite[1.4.3]{rudinfourier}}]  
\label{thm:bochner}
Let
  $H$ be an abelian locally compact Hausdorff group.  A continuous
  function $f$ on $H$ is positive definite if and only if there is a
  finite non-negative measure $\mu$ on $\widehat{H}$ such that for $s\in H$
\begin{equation}
\label{eq:7}
f(s) = \int_{\widehat{H}} \overline{\omega(s)}d\mu(\omega).
\end{equation}
Furthermore $\mu(\widehat{H}) = \|f\|_\infty$.  
\end{theorem}

\begin{proof}
The first part of the theorem is Bochner's Theorem as stated in
\cite{rudinfourier}.  For the last statement, since $f$ is positive
definite, $|f(s)|\leq f(0)$ for all $s\in H$ \cite[1.4.1]{rudinfourier}.  Hence 
\[
\|f\|_\infty = f(0) = \int_{\widehat{H}} \overline{\omega(0)}
d\mu(\omega) = \mu(\widehat{H}).\qedhere
\]
\end{proof}

It is not particularly hard to see that the span of the functions of
positive type are exactly those functions defined via \eqref{eq:7}
except where $\mu$ is a complex measure.  The following theorem says
that the Fourier transform is very well behaved for this kind of
function. 

\begin{theorem}[Inversion Theorem {\cite[1.5.1]{rudinfourier}}]  
\label{thm:inversion}
Let
  $H$ be an abelian locally compact Hausdorff group with Haar measure
  $\lambda$.  Suppose $f\in L^1(H)$ is such that for all $s\in H$
\[
f(s) = \int_{\widehat{H}} \overline{\omega(s)} d\mu(\omega)
\]
for some complex measure $\mu$ on $\widehat{H}$.  Then the Fourier
transform $\hat{f}$ (defined in Remark \ref{rem:6}) 
is in $L^1(\widehat{H})$.  Furthermore, there is a
Haar measure $\widehat{\lambda}$ on $\widehat{H}$ such that for all
functions of this form 
\[
f(s) = \int_{\widehat{H}} \hat{f}(\omega)\overline{\omega(s)}
d\widehat{\lambda}(\omega).
\]
\end{theorem}

\begin{remark}
Given $f\in L^1(\widehat{H})$ the function 
\[
\check{f}(s) = \int_{\widehat{H}} f(\omega)\overline{\omega(s)}
d\widehat{\lambda}(\omega)
\]
is called the ``inverse Fourier transform'' of $f$.  Colloquially
Theorem \ref{thm:inversion} says that for a certain class of functions
on $H$ the inverse Fourier transform of the Fourier
transform of $f$ is equal to $f$. 
\end{remark}

\begin{definition} 
\index[not]{$\widehat{\lambda}$}
Given an abelian locally compact Hausdorff group $H$ with Haar measure
$\lambda$ we call the measure $\widehat{\lambda}$ coming from Theorem
\ref{thm:inversion} the {\em dual Haar measure}.  We generally denote
integration with respect to $\lambda$ by $ds$ and integration with
respect to $\widehat{\lambda}$ by $d\omega$.  
\end{definition}

\begin{theorem}[Plancharel's Theorem {\cite[1.6.1]{rudinfourier}}]
\label{thm:plancharel}
Let $H$ be an abelian locally compact Hausdorff group with Haar
measure $\lambda$ and dual measure $\widehat{\lambda}$.  The Fourier
transform, restricted to $L^1(H)\cap L^2(H)$, is an isometry with
respect to the $L^2$-norms on $L^2(H,\lambda)$ and
$L^2(\widehat{H},\widehat{\lambda})$.  Furthermore it maps onto a
dense subspace of $L^2(\widehat{H})$ and can be extended to an
isometry of $L^2(H)$ onto $L^2(\widehat{H})$.  
\end{theorem}

It is notable that the isomorphism between $L^2(H)$ and
$L^2(\widehat{H})$ given by Theorem \ref{thm:plancharel} cannot be
explicitly defined off of $L^1(H)\cap L^2(H)$.  

\begin{lemma}
\label{lem:3}
If $H$ and $f$ are as in Theorem \ref{thm:inversion} then $d\mu =
\hat{f}d\omega$.  
\end{lemma}

\begin{proof}
For $g\in C_c(H)$ we have, by Theorem \ref{thm:plancharel}, 
\begin{align}
\label{eq:8}
\int_{\widehat{H}} \overline{\hat{g}}\hat{f}d\omega &=
\int_H \overline{g} f ds = \int_H\int_{\widehat{H}}
\overline{g(s)\omega(s)} d\mu(\omega) ds \\
&= \int_{\widehat{H}}\overline{\int_G g(s)\omega(s)ds}d\mu(\omega) =
\int_{\widehat{H}} \overline{\hat{g}}d\mu(\omega). \nonumber
\end{align}
It is well known \cite{rudinfourier} that the image of 
$C_c(H)$ under the Fourier transform is dense in
$C_0(\widehat{H})$.  Given 
$\epsilon > 0$ and $\phi\in C_0(\widehat{H})$ choose $g\in
C_c(H)$ such that 
\[
\| \hat{g} - \phi\|_\infty <
\frac{\epsilon}{2\max\{|\mu|(\widehat{H}),\|\hat{f}\|_1\}}.
\]
Then, using \eqref{eq:8}, we have
\begin{align*}
\left| \int_{\widehat{H}} \overline{\phi} \hat{f} d\omega - \int_{\widehat{H}}
  \overline{\phi} d\mu \right| &\leq \left|
  \int_{\widehat{H}}(\overline{\phi}-\overline{\hat{g}})\hat{f}d\omega\right| +
\left|\int_{\widehat{H}} (\overline{\phi}-\overline{\hat{g}})d\mu\right| \\
&\leq \int_{\widehat{H}} |\phi-\hat{g}||\hat{f}|d\omega +
\int_{\widehat{H}} |\phi-\hat{g}|d|\mu|  \\
&\leq \|\phi-\hat{g}\|_\infty(\|\hat{f}\|_1 + |\mu|(\widehat{H})) <
\epsilon.
\end{align*}
Since this is true for arbitrary $\epsilon > 0$ and $\phi\in
C_0(\widehat{H})$, after replacing $\phi$ by $\overline{\phi}$ we may
conclude that  $\hat{f}d\omega = d\mu$.  
\end{proof}

At this point we have all the classical harmonic analysis that we
need.  It is time to return to the groupoid case.  First, we would
like to see that given a continuously varying abelian group bundle the
dual bundle is also continuously varying. Fortunately, this very
thing is shown in \cite{ctgIII}, once we recall that a group bundle
has a Haar system if and only if it is continuously varying.   We have
reproduced the proof here for completeness.  

\begin{prop}[{\cite[Proposition 3.6]{ctgIII}}]
\index{Haar system}
If $S$ is an abelian, second countable, locally compact Hausdorff group
bundle with Haar system $\beta = \{\beta^u\}$ then the dual measures
$\hat{\beta}^u$ form a Haar system for the dual bundle $\widehat{S}$.  
\end{prop}

\begin{proof}
Suppose that $K$ is compact in $\widehat{S}$.  We claim that $u\mapsto
\hat{\beta}^u(K)$ is bounded on $S\unit$.  Of course, it suffices to
consider only $u\in \hat{p}(K)$.  Let $f\in C_c(S)$ be a non-negative
function such that 
\begin{equation}
\label{eq:16}
\int_{S_u} f(s)^2 d\beta^u(s) = 1\  \text{for all $u\in
  \hat{p}(K)$.}
\end{equation}
Since $\hat{p}(K)$ is compact, there is an $\epsilon > 0$ so that 
\begin{equation}
\label{eq:19}
\int_{S_u} f(s) d\beta^u(s) > \epsilon\ \text{for all $u\in
  \hat{p}(K)$.}
\end{equation}
Using Theorem \ref{thm:plancharel} \eqref{eq:16} implies that
\begin{equation}
\label{eq:20}
\int_{\widehat{S}_u} |\hat{f}(\omega)|^2 d\hat{\beta}^u(\omega) = 1
\ \text{for all $u\in \hat{p}(K)$}.
\end{equation}
Moreover the continuity of $\hat{f}$ and \eqref{eq:19} imply that $U =
\{\omega \in \widehat{S} : |\hat{f}(\omega)|^2 > \epsilon^2 \}$ is an
open neighborhood of $\hat{p}(K)$. 

If $\omega\in K$, then $\omega\inv\omega \in \hat{p}(K)$.  The
continuity of multiplication implies that there is a neighborhood $V$
of $\omega$ such that $V\inv V\subset U$.  Therefore there is a cover
$V_1,\ldots,V_m$ of $K$ such that $V_j\inv V_j \subset U$ for each
$1\leq j \leq m$.  In view of \eqref{eq:20}, $\hat{\beta}^u(U) \leq
\epsilon^{-2}$ if $u\in \hat{p}(K)$.  Furthermore, if
$\hat{\beta}^u(V_j)\ne 0$, then there is an $\omega \in V_j$ with
$\hat{p}(\omega) = u$.  Then 
\[
\hat{\beta}^u(V_j) = \hat{\beta}^u(\omega\omega\inv V_j) \leq
\hat{\beta}^u(\omega V_j\inv V_j) \leq \hat{\beta}^u(\omega U) = 
\hat{\beta}^u(U) \leq 1/\epsilon^2.
\]
It follows that $\hat{\beta}^u(K)\leq m/\epsilon^2$ for all $u \in
S\unit$.  This proves the claim.  

Now let $\{u_i\}_{i\in I}$ be a net in $S\unit$ converging to $u\in
S\unit$.  If $\phi\in C_c(\widehat{S})$, then let
$\hat{\beta}(\phi)(v) =
\int_{\hat{S}_v}\phi(\omega)d\hat{\beta}^v(\omega)$.  The above
argument implies that $\{\hat{\beta}(\phi)(u_i)\}_{i\in I}$ is
bounded.  Thus if $\omega$ is a generalized limit\footnote{A
  generalized limit is a norm one extension of the ordinary limit
  functional on the subspace $c_0$ of $\ell^\infty(I)$ consisting of those nets
  $\{a_i\}$ such that $\lim_i a_i$ exists.} on $\ell^\infty(I)$,
then we obtain a positive linear functional $\mu$ on
$C_c(\widehat{S})$ by $\mu(\phi) =
\omega(\{\hat{\beta}(\phi)(u_i)\})$.  Suppose that
$\phi,\psi\in C_c(\widehat{S})$ agree on $\widehat{S}_u$.  Then if $K$
is a compact set containing the supports of $\phi$ and $\psi$, 
\begin{align*}
|\hat{\beta}(\phi)(u_i) - \hat{\beta}(\psi)(u_i)| &\leq
\int_{\widehat{S}_{u_i}}
|\phi(\omega)-\psi(\omega)|d\hat{\beta}^{u_i}(\omega) \\
&\leq \sup_{\omega\in \widehat{S}_{u_i}} |\phi(\omega)-\psi(\omega)|
\hat{\beta}^{u_i}(K).
\end{align*}
Since $\sup_{\omega\in\widehat{S}_v}|\phi(\omega)-\psi(\omega)|$ tends
to zero as $v$ tends to $u$ and since $v\mapsto \hat{\beta}^v(K)$ is
bounded, it follows that $\mu(\phi) = \mu(\psi)$.  Since every
function in $C_c(\widehat{S}_u)$ has an extension to an element of
$C_c(\widehat{S})$, we can view $\mu$ as a Radon measure on
$\widehat{S}_u$.  

However, if $f\in C_c(S)$ then, by the Plancharel Theorem,
$\hat{\beta}(|\hat{f}|^2)(u_i) = \beta(|f|^2)(u_i)$
which converges to $\beta(|f|^2)(u) = \hat{\beta}(|\hat{f}|^2)(u)$.  It
follows that 
\begin{equation}
\label{eq:21}
\hat{\beta}(|\hat{f}|^2)(u) = \mu(|\hat{f}|^2)\ \text{for all $f\in
  C_c(S_u)$.}
\end{equation}
By density, \eqref{eq:21} holds for all $\hat{f}\in
L^2(\widehat{S}_u,\hat{\beta}^u)$.  In particular $\mu =
\hat{\beta}^u$ on $C_c(\widehat{S}_u)$.  

We have shown that if $\{u_i\}$ is any net converging to $u$ in
$S\unit$, then $\omega(\{\hat{\beta}(\phi)(u_i)\}) =
\hat{\beta}(\phi)(u)$.  Therefore $\lim_i \hat{\beta}(\phi)(u_i) =
\hat{\beta}(\phi)(u)$, and it follows that $\{\hat{\beta}^u\}$ is a
Haar system. 
\end{proof}

\begin{definition}
Given an abelian, second countable locally compact Hausdorff group
bundle with Haar system $\beta = \{\beta^u\}$ then the Haar system
formed by the collection of dual measures $\hat{\beta}
=\{\hat{\beta}^u\}$ is called the {\em dual Haar system}.
\end{definition}

This is interesting because it means that given an abelian continuously
varying group bundle $S$ we can form the double dual $\doubledual{S}$
by taking the dual of $\widehat{S}$.  It is natural to ask whether or
not this is isomorphic to the original group bundle.  The following
lemma gets us most of the way there.  

\begin{lemma}
\label{lem:4}
\index[not]{$\hat{s}$}
Given an abelian continuously varying group bundle $S$ there is a
continuous bijective groupoid homomorphism $\Phi:S\rightarrow
\doubledual{S}$ given for $s\in S$ and $\omega\in \widehat{S}$ by 
\begin{equation}
\label{eq:9}
\Phi(s)(\omega) = \hat{s}(\omega) := \omega(s).
\end{equation}
\end{lemma}
\begin{proof}
It follows from Proposition \ref{prop:34} that $\doubledual{S}_u$ is
the double dual of $S_u$ for all $u\in S\unit$.  Furthermore,
classical Pontryagin duality \cite{rudinfourier} 
states that $s\mapsto \hat{s}$ is a group
isomorphism from $S_u$ onto $\doubledual{S}_u$ for all $u\in S\unit$.
Since $\Phi$ is formed by gluing all of these fibre isomorphisms
together it is clear that $\Phi$ is at least a bijective groupoid
homomorphism.  Next, we show that it is continuous.  Suppose
$s_i\rightarrow s$ in $S$ and let $u_i = p(s_i)$ and $u=p(s)$.  
By Proposition \ref{prop:33} it will
suffice to show, one, that $\hat{\hat{p}}(\Phi(s_i))\rightarrow
\hat{\hat{p}}(\Phi(s))$ and, two, that given $\omega_i\in \widehat{S}_{u_i}$
and $\omega \in \widehat{S}_u$ such that $\omega_i\rightarrow\omega$
then $\Phi(s_i)(\omega_i)\rightarrow\Phi(s)(\omega)$.  
Since $\Phi$ preserves the bundle maps (it's a groupoid homomorphism)
we know that $\hat{\hat{p}}(\Phi(s_i)) = u_i$ and
$\hat{\hat{p}}(\Phi(s))=u$ so that clearly the first condition
holds. Now suppose $\omega_i\in \widehat{S}_{u_i}$ and
$\omega\in\widehat{S}$ are such that $\omega_i\rightarrow \omega$.
All we have to do is cite Proposition \ref{prop:33} yet again to see
that 
\[
\Phi(s_i)(\omega_i) = \omega_i(s_i)\rightarrow \omega(s) =
\Phi(s)(\omega).\qedhere
\]
\end{proof}

If we were working with groups we would be done since continuous
bijective group homomorphisms 
between second countable locally compact groups are
automatically  bicontinuous.  This follows from Souslin's Theorem \cite[Theorem
3.2.3]{invitation}, which states that continuous bijections have a
Borel inverse, and the fact that measurable homomorphisms between
second countable locally compact groups are automatically continuous 
\cite[Theorem D.3]{tfb2}.  However, we will
show in Section \ref{sec:opencounter} that this is ``open mapping
theorem'' is not true for group bundles.  Fortunately, it
turns out that in this specific case $\Phi$ is a homeomorphism.  It
should be noted that the the following theorem is also stated, without
proof, in \cite[Proposition 1.3.7]{ramazan}.

\begin{theorem}[{\cite[Theorem 16]{bundleduality}}]
\label{thm:duality}
\index{dual bundle}
\index{New Result}
If $S$ is an abelian, continuously varying, second countable, locally compact 
Hausdorff group bundle then the map $\Phi:S\rightarrow \doubledual{S}$ such that
$\Phi(s) = \hat{s}$ is a groupoid isomorphism.  
\end{theorem}

\begin{proof}
It follows from Lemma \ref{lem:4} that all we need to do is show that if 
$\{s_i\}_{i\geq0}\subset S$ such that $\hat{s}_i\rightarrow \hat{s}_0$ 
then $s_i\rightarrow s_0$.  Let $u_i = p(s_i)$ for all $i\geq 0$ and 
recall from \eqref{eq:9} that
$\hat{s}(\omega) := \omega(s)$.  It follows from Definition
\ref{def:22} that $\hat{s}_i(\phi)\rightarrow \hat{s}_0(\phi)$ for all
$\phi\in C_c(\widehat{S})$.  Using
\eqref{eq:5} we have, for all $\phi\in C_c(\widehat{S})$, 
\begin{equation}
\label{eq:10}
\int_{\widehat{S}} \phi(\omega)\omega(s_i)d\hat{\beta}^{u_i}(\omega)
\rightarrow
\int_{\widehat{S}} \phi(\omega)\omega(s_0)d\hat{\beta}^{u_0}(\omega).
\end{equation}

Now suppose we have a relatively compact open neighborhood $V$ of
$u_0$ in $S$.  Then, using the continuity of the operations, there
exists a relatively compact open neighborhood $U$ of $u_0$ in $S$ such
that $U = U\inv$ and 
$U^2\subset V$.  Choose $h\in C_c(S)$ such that $h(u_0)=1$ and
$\supp h\subset U$.  Define $f$ on $S$ by
\[
f(s) := h^* * h(s) = \int_S \overline{h(t\inv)} h(t\inv s) d\beta^{p(s)}(t).
\]
It is straightforward to check that $f$ is continuous.  Furthermore, $f(s) = 0$
unless there exists $t\in S$ such that $t\inv \in U$ and $t\inv s\in
U$.  In other words, unless $s\in U^2\subset V$.  Therefore $\supp
f\subset V$ and $f\in C_c(S)$.  From now on let $f^u$ denote the
restriction of $f$ to $S_u$.  It is clear from the definition of $f$, and
Example \ref{ex:14}, that $f^u$ is positive definite. Therefore $f^u$
satisfies the conditions of Bochner's theorem and the inversion
theorem for all $u\in S\unit$.  In particular, this implies that for
each $u\in S\unit$ there exists a finite positive measure $\mu^u$ on
$\widehat{S}_uS$, which we extend to $\widehat{S}$ by giving everything else
measure zero, such that for all $s\in S$
\[
f(s) = \int_{\widehat{S}} \overline{\omega(s)} \mu^{p(s)}(\omega).
\]
Theorem \ref{thm:bochner} also implies that $\mu^u(\widehat{S}) =
\mu^u(\widehat{S}_u) \leq \|f^u\|_\infty \leq \|f\|_\infty$ for all
$u\in S\unit$.  It follows from Lemma \ref{lem:3} that for all $u\in S\unit$.
\begin{equation}
\label{eq:11}
\hat{f} d\hat{\beta}^u = d\mu^u
\end{equation}
as measures on $\widehat{S}_u$.  However, since both $\hat{\beta}^u$
and $\mu^u$ have support contained in $\widehat{S}_u$ equation
\eqref{eq:11} holds for $\hat{\beta}^u$ and $\mu^u$ as measures on all
of $\widehat{S}$.  

Now, we don't know that $\hat{f}$ is compactly supported.  In fact,
it's probably not.  However if $\phi\in C_c(\widehat{S})$ then
$\phi\hat{f}$ is compactly supported.  It follows from \eqref{eq:10}
that 
\begin{equation}
\label{eq:12}
\int_{\widehat{S}}
\phi(\omega)\hat{f}(\omega)\omega(s_i)d\hat{\beta}^{u_i}(\omega)
\rightarrow \int_{\widehat{S}}
\phi(\omega)\hat{f}(\omega)\omega(s_0)d\hat{\beta}^{u_0}(\omega). 
\end{equation}
Using \eqref{eq:11} we can rewrite \eqref{eq:12} as 
\begin{equation}
\label{eq:13}
\int_{\widehat{S}}
\phi(\omega)\omega(s_i)d\mu^{u_i}(\omega)\rightarrow
\int_{\widehat{S}} \phi(\omega)\omega(s_0) d\mu^{u_0}(\omega).
\end{equation}
In order to make the notation a little more palatable we will
temporarily define, for all $i\geq 0$,
\[
a_i(\phi) := \int_{\widehat{S}}
\phi(\omega)\omega(s_i)d\mu^{u_i}(\omega).
\]

We would like to extend \eqref{eq:13} to functions in
$C_0(\widehat{S})$.  Suppose we have $\psi \in C_0(\widehat{S})$ and
are given $\epsilon > 0$.  Choose $\phi\in C_c(\widehat{S})$ such that
$\| \psi - \phi\|_\infty < \frac{\epsilon}{4\|f\|_\infty}$.  Now
choose $I$ such that $|a_i(\phi)-a_0(\phi)|< \epsilon/2$ for all
$i\geq I$.  We then compute for $i\geq I$ that
\begin{align*}
|a_i(\psi) - a_0(\psi)| =& \left| \int_{\widehat{S}}
  \psi(\omega)\omega(s_i) d\mu^{u_i} - \int_{\widehat{S}}
  \psi(\omega)\omega(s_0) d\mu^{u_0} \right| \\
\leq& \left| \int_{\widehat{S}} (\psi(\omega)-\phi(\omega))\omega(s_i)
  d\mu^{u_i}\right| \\ 
&+\left| \int_{\widehat{S}}\phi(\omega)\omega(s_i)d\mu^{u_i} -
  \int_{\widehat{S}}\phi(\omega)\omega(s_0)d\mu^{u_0} \right| \\
&+\left| \int_{\widehat{S}}(\psi(\omega)-\phi(\omega))\omega(s_0)d\mu^{u_0} \right| \\
\leq& |a_i(\phi)-a_0(\phi)| + \int_{\widehat{S}}|\psi-\phi| d\mu^{u_i} + \int_{\widehat{S}}
|\psi-\phi| d\mu^{u_0} \\
\leq& |a_i(\phi)-a_0(\phi)| +
(\mu^{u_i}(\widehat{S})+\mu^{u_0}(\widehat{S}))\|\psi-\phi\|_\infty \\
<& \frac{\epsilon}{2} + \frac{2\|f\|_\infty\epsilon}{4\|f\|_\infty} =
\epsilon.
\end{align*}
It follows that if $\phi\in C_0(\widehat{S})$ then
$a_i(\phi)\rightarrow a_0(\phi)$, or equivalently, that \eqref{eq:12}
holds for all $\phi\in C_0(\widehat{S})$.  

Next, if $g\in C_c(S)$ then, since both $g$ and $f$ are compactly
supported and $\beta^x$ is regular for all $u\in S\unit$, we have
$f^u,g^u\in L^1(S_u,\beta^u)\cap L^2(S_u,\beta^u)$ for all $u\in S\unit$.
The Plancharel Theorem states that the Fourier transform is an
isometry on ${L^1(S_u)\cap L^2(S_u)}$.   Therefore we have, for all
$u\in S\unit$,
\[
\int_{S_u} \overline{g^u}f^u d\beta^u = \int_{\widehat{S}_u}
\overline{\widehat{g^u}} \widehat{f^u} d\hat{\beta}^u.
\]
Observe that, for all $i\geq 0$ and $\omega\in\widehat{S}_{u_i}$,
\begin{align*}
\overline{\widehat{g^{u_i}}(\omega)}\omega(s_i) &= \int_{S_{u_i}}
  \overline{g^{u_i}(s)\omega(s)}\omega(s_i) d\beta^{u_i}(s) = 
\int_{S_{u_i}} \overline{g^{u_i}(s)}\omega(s\inv s_i) d\beta^{u_i}(s)
\\
&= \int_{S_{u_i}} \overline{g^{u_i}(s_is)}\omega(s\inv) d\beta^{u_i}(s)
= \overline{(\lt_{s_i\inv} g^{u_i})\sidehat}(\omega).
\end{align*}
where we make the usual definition $\lt_s f(t) := f(s\inv t)$.
Therefore, for all $i\geq 0$, we can compute
\begin{align*}
\int_{\widehat{S}} \overline{\hat{g}(\omega)}
\hat{f}(\omega) \omega(s_i) d\hat{\beta}^{u_i}(\omega) &=
\int_{\widehat{S}_{u_i}}
\overline{\widehat{g^{u_i}}(\omega)}\widehat{f^{u_i}}(\omega)\omega(s_i)
d\hat{\beta}^{u_i}(\omega) \\
&= \int_{\widehat{S}_{u_i}} \overline{(\lt_{s_i\inv}
  g^{u_i})\sidehat}\widehat{f^{u_i}} d\hat{\beta}^{u_i} =
\int_{S_{u_i}} \overline{\lt_{s_i\inv}g^{u_i}} f^{u_i}
d\beta^{u_i}.
\end{align*}
Using \eqref{eq:12} it follows that for all $g\in C_c(S)$ 
\begin{equation}
\label{eq:14}
\int_{S_{u_i}} \overline{\lt_{s_i\inv}g^{u_i}}f^{u_i}d\beta^{u_i}
\rightarrow
\int_{S_{u_0}} \overline{\lt_{s_0\inv}g^{u_0}}f^{u_0}d\beta^{u_0}.
\end{equation}

Now we are finally ready to attack the convergence of the $s_i$.
Choose an open neighborhood $O$ of $s_0$.  Using the continuity of
multiplication we can find relatively compact open neighborhoods $V$
and $W$ in $S$ such that $u_0\in V$, $s_0\in W$, and $VW\subset O$.
Furthermore, by intersecting $V$ and $V\inv$ we can assume that $V\inv
= V$.  Construct $f$ for $V$ as in the beginning of the proof.  Now
choose $g\in C_c(S)$ so that $0\leq g\leq 1$, $g(s_0)=1$, and $g$ is
zero off $W$.  Then $\overline{g} = g$ so that by \eqref{eq:14} we
have 
\begin{equation}
\label{eq:15}
\int_{S_{u_i}} g(s_it)f(t) d\beta^{u_i}(t) \rightarrow \int_{S_{u_0}}
g(s_0t)f(t)d\beta^{u_0}(t).
\end{equation}
Given $i\geq 0$, $\int_{S_{u_i}} g(s_it)f(t)d\beta^{u_i}(t) = 0$ unless
there exists $t$ such that $s_it \in \supp g\subset W$ and $t\in \supp
f \subset V$.  This implies that the integral is zero unless $s_i \in
WV\inv = WV \subset O$.  Furthermore, both $g(s_0 u_0)$ and $f(u_0)$
are nonzero by construction and, since both are positive continuous
functions, this implies $ \int_{S_{u_0}} g(s_0 t) f(t) d\beta^{u_0}(t) \ne 0.$
It follows from \eqref{eq:15} that eventually $\int_{S_{u_i}} g(s_i t)
f(t) d\beta^{u_i}(t) \ne 0$ so that, eventually, $s_i \in O$.  Of
course, it follows that $s_i\rightarrow s_0$ and we are done. 
\end{proof}


\section{Open Mapping Counterexample}
\label{sec:opencounter}

As we noted in Section \ref{sec:duality}, given a second
countable continuously
varying abelian group bundle $S$ it is easy to see that the natural
map from $S$ to its double dual is a continuous,
bijective, groupoid homomorphism.
Furthermore, in the second countable locally compact Hausdorff 
group case this map would automatically have a
continuous inverse.  If this kind of ``open
mapping theorem'' were true for second countable group bundles then
Theorem \ref{thm:duality} would be trivial.  In this section we will
exhibit an example which shows that not every continuous bijective
groupoid homomorphism between second countable, abelian, continuously
varying,
group bundles is necessarily a homeomorphism.  Specifically we will
prove the following 

\begin{theorem}
\label{thm:counter}
\index{group bundle}
\index{New Result}
There exists second countable, locally compact Hausdorff, abelian
group bundles $S$ and $T$ and a map $\phi:T\rightarrow S$ such that
$\phi$ is a continuous, bijective, groupoid
homomorphism which is not a homeomorphism.  
\end{theorem}

Finding an example of such a homomorphism is especially tricky.
Fibrewise such a map is a continuous bijective group
homomorphism of second countable groups, and as such must be a
homeomorphism when restricted to the fibres.  Of course, we have to
start by defining $S$ and $T$.  

\begin{remark}
In this section it will be convenient to define 
\[
\Z_{2n+1} =
\{-n,{-n+1},\ldots,{-1},0,1,\ldots,n-1,n\}
\] 
for all $n\in\N$.
Furthermore, when appropriate we will still give $\Z_{2n+1}$ the group
operation of addition modulo $2n+1$.  For the remainder of this
section we are also going to define $x_n = 1/n$ for $n>0$ and $x_0=0$
and will let $X = \{x_n\}_{n=0}^\infty$.  
\end{remark}

The following technical lemma will see a lot of use in this section. 

\begin{lemma}
\label{lem:5}
Suppose $x_{n_i}\rightarrow x_N$ converges in $X$.  Then we can pass
to a subsequence, relabel, and assume either
\begin{enumerate}
\item $n_i=N$ for all $i$, or 
\item $N=0$, $n_i\rightarrow\infty$, and $n_i\ne 0$ for all $i$.  
\end{enumerate}
\end{lemma}
\begin{proof}
There are two cases to consider.
First, if $N\ne 0$ then, because $x_{n_i}\rightarrow x_N$ and $x_N$ is
open as a point in $X$, eventually $x_{n_i} = x_N$.  Therefore, 
$n_i = N$ eventually, and by
passing to a subsequence we can assume that $n_i=N$ for all $i$. On
the other hand, 
suppose $N=0$.  If $n_i = 0$ infinitely often then we can pass to a
subsequence and assume that $n_i = N = 0$ for all $i$.  If it is not
true that $n_i = 0$ infinitely often then, because
$x_{n_i}\rightarrow x_0$, we know that $n_i$ does not equal any
$N\in\N$ infinitely often.  However, this implies that
$n_i\rightarrow \infty$.  We could then pass to a subsequence again
and assume that $n_i \ne 0$ for all $i$. 
\end{proof}

Next we define the total space of one of our bundles and
show that it is topologically well behaved.  

\begin{prop}
\label{prop:25}
Let 
\[
T = \{(m,x_n)\in \Z\times X : \text{$m\in \Z_{6n+3}$ if $n>0$ and
  $m\in\Z$ otherwise}\}.
\] 
Then $T$ is closed in $\Z\times X$, and hence second countable locally compact
Hausdorff.  Furthermore the map $p_T:T\rightarrow X$ defined by
projection onto the second factor is continuous and open.  This makes
$T$ into a bundle over $X$ with fibres denoted $T_n$.  
\end{prop}

\begin{proof}
Suppose $(m_i,x_{n_i})\rightarrow (m,x_N)$ in $\Z\times X$.  We would
like to show that $(m,x_N)\in T$.  Apply Lemma \ref{lem:5} and
pass to a subsequence.  We know
that there are two cases to consider.  First, suppose that $n_i =N$
for all $i$.  Then for all $i$ we have either $N\ne 0$ and
$m = m_i \in \Z_{6n_i+3} = \Z_{6N+3}$ or $N=0$ and $m = m_i \in \Z$.
It follows that $(m,x_N)\in T$.  In the second case $N=0$, but then
$(m,x_0)\in T$ automatically.  Thus $T$ is closed and since it is a
closed subset of a second countable locally compact Hausdorff space it
is also second countable locally compact Hausdorff. 

Let $p_T:T\rightarrow X$, denoted by $p$ when convenient, be given by
$p(m,x) = x$ for $(m,x)\in T$.  Since $p$ is the restriction of a
continuous map it is continuous.  We will show $p$ is open.  Suppose
$x_{n_i}\rightarrow x_N$ and that $(m,x_N)\in T$.  We will show that
we can pass to a subsequence, relabel, and find $m_i$ such that
$(m_i,x_{n_i})\rightarrow (m,x_N)$ and $(m_i,x_{n_i})\in T$ for all
$i$.  By applying Lemma \ref{lem:5}, and possibly passing to a subsequence
and relabeling, we
can either assume that $n_i = N$ eventually or that $N=0$, $n_i \rightarrow
\infty$ and $n_i \ne 0$ for all $i$.  Suppose the former is true.
Then either $N\ne 0$ and 
$m\in \Z_{6n_i +3} = \Z_{6N+3}$ or $N=0$ and $m\in\Z$ for all $i$ so that
$(m,x_{n_i})\in T$ for all $i$ and clearly $(m,x_{n_i})\rightarrow
(m,x_N)$.  Suppose the latter is true.  Then eventually $|m|\leq
3n_i+1$ and by passing to a subsequence we can assume that this is
true for all $i$.  Hence $m\in \Z_{6n_i+1}$ and $(m,x_{n_i})\in T$ for
all $i$, and clearly $(m,x_{n_i})\rightarrow (m,x_0)$.  It now follows
from Proposition \ref{prop:9} that $p$ is open. 
\end{proof}

Next, we add the groupoid structure on $T$. 

\begin{prop}
If we endow $T_n$ with the operations of addition modulo $6n+3$ for
$n> 0$ and $T_0$ with the usual integer addition, then $T$ is a second
countable, continuously varying, abelian, locally compact Hausdorff
group bundle and the bundle map for $T$ can be identified with $p$. 
\end{prop}

\begin{proof}
It follows from Proposition \ref{prop:25} that $T$ has the appropriate
topological conditions.  Furthermore, it is easy to see that, using
the operations defined in the statement of the proposition, $T$ is an
abelian group bundle with unit space $X$ and that the bundle map is
exactly $p$.  After all, algebraically $T$ is just the disjoint union
of the groups $T_n$.  

All that is left is to show that the operations on $T$ are
continuous.  Suppose $(m_i,x_{n_i})\rightarrow (x_N)$ and
$(m_i',x_{n_i})\rightarrow (m',x_N)$ in $T$.  We would like to show
that 
\begin{align*}
(m_i,x_{n_i}) + (m_i',x_{n_i}) &\rightarrow (m,x_N) + (m',x_N),\quad\text{and} \\
-(m_i,x_{n_i}) &\rightarrow -(m,x_N).
\end{align*}
It will suffice to show that for each subsequence, a sub-subsequence
converges in the above fashion.  So pass to a subsequence.  Using
Lemma \ref{lem:5} we can pass to another subsequence and either assume
that $n_i = N$ for all $i$ or that $N=0$ and $n_i\rightarrow\infty$.
Suppose the former is true.  Then, eventually $m_i = m$ and $m_i'=m'$
so that eventually $(m_i,x_{n_i})+(m_i',x_{n_i}) = (m,x_N)+(m',x_N)$.
Similarly in this case $-(m_i,x_{n_i})=-(m,x_N)$ and at this point
convergence is clear.  Now suppose we are in the second case so that
$n_i\rightarrow \infty$, $n_i\ne 0$ for all $i$ and $N=0$.  
As before we have $m_i=m$ and $m_i'=m'$
eventually so that, for large $i$, $|m_i+m_i'| = |m+m'| \leq n_i$.
However, when this is true $m+m'\mod 6n_i+3 = m+m'$.  It follows that
for large enough $i$
\[
(m_i,x_{n_i})+(m_i',x_{n_i}) = (m+m',x_{n_i})
\]
and it is clear that $(m+m',x_{n_i})\rightarrow (m+m',x_0) =
(m,x_0)+(m',x_0)$.  Since $-(m_i,x_{n_i}) = (-m_i,x_{n_i})$ for all $i$
it is easier to see that, assuming $i$ is large enough,
\[
-(m_i,x_{n_i}) = (-m,x_{n_i})\rightarrow (-m,x_0) = -(m,x_0).
\]
It follows that the operations on $T$ are continuous and that $T$ is a
topological groupoid.  However, we can now conclude that $T$ is
continuously varying, since the bundle map $p$ is open.
\end{proof}

\begin{remark}
Since $T$ is clearly $r$-discrete the Haar system on $T$ is
given by $\lambda^{x_n} = \mu\times\delta_{x_n}$ where $\mu$ is
counting measure and $\delta_{x_n}$ is the Dirac delta measure at
$x_n$.  
\end{remark}

Next we define the total space of our other bundle and show that it is
well behaved topologically.  

\begin{prop}
\label{prop:28}
Let $A_n = \{-1/n,0,1/n\}$ and $S_n = A_n \times \Z_{2n+1} \times
\{x_n\}$ for $n> 1$.  Let $A_0 = \{0\}$ and $S_0 =
A_0\times\Z\times\{x_0\}$.  Define 
\[
S = \bigcup_{n=0}^\infty S_n
\]
and give $S$ the relative topology as a subset of $\R\times\Z\times
X$.  Then $S$ is closed and is therefore a second countable locally
compact Hausdorff space.  Furthermore, the map $p_S:S\rightarrow X$
defined by projection onto $X$ is continuous and open.  
\end{prop}

\begin{proof}
Suppose $(a_i,m_i,x_{n_i})\rightarrow (a,m,x_N)$ in $S$.  Using Lemma
\ref{lem:5}, pass to a subsequence, relabel, and assume that either
$n_i = N$ for all $i$ or that $n_i\rightarrow \infty$, $n_i\ne 0$ for
all $i$ and $N=0$.  Consider the first case.  
Observe that eventually either $N\ne 0$ and $m_i=m \in
\Z_{2n_i+1} = \Z_{2N+1}$ or $N=0$ and $m_i = m \in \Z$. 
Furthermore, we also have, for large $i$,
$a_i \in A_{n_i} = A_N$.  Since $A_N$ is a closed set (it consists of
at most three distinct points), $a = \lim a_i \in A_N$.  Thus $(a,m,x_N)\in
S_N \subset S$.  On the other hand suppose that $n_i\rightarrow\infty$
and $n_i\ne 0$ for all $i$.   Since $-1/n_i \leq a_i \leq
1/n_i$ for all $i$ we can conclude that $0 = \lim a_i = a$.  Hence
$(a,m,x_N)\in S_0\subset S$.  It follows that $S$ is closed.  

Let $p_S:S\rightarrow X$, often denoted by $p$, be defined by
$p_S(a,m,x) = x$.  Since $p$ is the restriction to $S$ of a continuous
map it must be continuous.  We will show that it is open.  Suppose
$x_{n_i}\rightarrow x_N$ and that $(a,m,x_N)\in S$.  Then, using Lemma
\ref{lem:5}, we pass to a subsequence, relabel, and assume that either
$n_i = N$ for all $i$ or that $n_i\rightarrow\infty$, $n_i\ne 0$ for
all $i$, and $N=0$.  Suppose the former is true.  Then either $N\ne
0$, $m \in \Z_{2N+1} = \Z_{2n_i+1}$, and $a\in A_N = A_{n_i}$ for all
$i$ or $N=0$, $m\in \Z$, and $a = 0$.  It follows that
$(a,m,x_{n_i})\in S_{n_i}$ for all $i$ and we clearly have
$(a,m,x_{n_i})\rightarrow (a,m,x_N)$.  On the other hand suppose the
latter case is true.  Since $N=0$ we know $a=0$.  Furthermore $|m|<
n_i$ eventually, so we might as well pass to a subsequence and assume
that this is always true.  It follows that $m\in Z_{2n_i+1}$ for all
$i$ and that $(0,m,x_{n_i})\in S_{n_i}$ for all $i$.  Furthermore it
is clear that $(a,m,x_{n_i})\rightarrow (a,m,x_N)$.  It follows that
$p$ is open. 
\end{proof}

\begin{lemma}
\label{lem:6}
Given $m\in\Z$ and $n>0$ there exists a unique $d\in\Z$ and $r\in\Z$
such that $m = d(2n+1)+r$ and $-n\leq r\leq n$.  Furthermore if $|m|
\leq 6n+3$ then $-1\leq d \leq 1$.  
\end{lemma}

\begin{proof}
This is a straightforward modification of the division algorithm where
we divide by $2n+1$ and allow the remainder to take on values between $-n$ and
$n$.  Furthermore if $|d|> 1$ then, as long as $-n\leq r \leq n$, we
have $|m|> 6n+3$.
\end{proof}

\begin{lemma}
\label{lem:7}
For $d\in\Z_3$ let $a^n(d) = d/n$ for all $n> 0$ and let $a^0(d) = 0$.  
Then every element of $A_n$ is of the form $a^n(d)$ for
$d\in \Z_3$.  Furthermore, 
given $n > 0$ we can define a bijective map $\phi:T_n\rightarrow S_n$
such that $\phi_n(m,x_n) = (a^n(d),r,x_n)$ where $d$ and $r$ are as in
Lemma \ref{lem:6}.  Finally, the map $\phi_0:T_0\rightarrow S_0$
such that $\phi(m,x_0) = (0,m,x_0)$ is a bijection. 
\end{lemma}

\begin{proof}
First, it is clear from the definition of $A_n$ that every element is
of the form $a^n(d)$ for $d=-1,0,1$.  
Let $n> 0$.  Given $(m,x_n)\in S_n$ we know $|m|\leq 6n+3$ by definition. Let
$d,r\in \Z$ be as in Lemma \ref{lem:6}.  Then $r\in \Z_{2n+1}$ and,
since $|d|\leq 1$, we know $a^n(d)\in A_n$.  Hence $\phi_n(m,x_n) =
(a^n(d),r,x_n)\in T_n$ and $\phi_n$ is well defined. Furthermore, if
$\phi_n(m,x_n) = (a^n(d),r,x_n) = \phi(m',x_n)$ then $m = d(2n+1)+r =
m'$ so that $(m,x_n) = (m',x_n)$.  Next, given $(a,r,x_n)\in S_n$
choose $d\in \Z_3$ so that $a = a^n(d)$.  Then $m = d(2n+1)+r \in
\Z_{6n+3}$ and, by the uniqueness of the factorization, $\phi(m,x_n) =
(a,r,x_n)$.  It follows that $\phi_n$ is a bijection.  Finally, it is
clear that $\phi_0$ is a bijection.  
\end{proof}

Next, we use the above maps to define the group structure on the
fibres.  

\begin{prop}
\label{prop:32}
Endow $S_0$ with the operations of integer addition and negation
in the second factor.  For $n>0$ define operations on $S_n$ via 
\begin{align}
\label{eq:17}
&\quad\quad\ \! -(a,r,x_n) = (-a,-r,x_n) \\
\label{eq:18}
(a^n(d),r,x_n)&+(a^n(d'),r',x_n) = \\ \nonumber
& \begin{cases}
(a^n(d+d'+1\mod 3), r+r' - (2n+1), x_n) & r+r' > n \\
(a^n(d+d'\mod 3), r+r', x_n) & -n \leq r+r' \leq n \\
(a^n(d+d'-1\mod 3), r+r' + (2n+1), x_n) & r+r' < -n.
\end{cases}
\end{align}
With these operations $S_n$ is an abelian group for all $n\geq 0$ and $\phi_n$
is a group isomorphism for all $n \geq 0$. 
\end{prop}

\begin{proof}
First observe that the topologies on $S_n$ and $T_n$ are discrete for
all $n$ so that $\phi_n$ is a homeomorphism for all $n$.  All we need
to do is show that each $\phi_n$ preserves the operations. Then the
operations on $S_n$ will automatically make $S_n$ into an abelian group and
$\phi_n$ will be an isomorphism.  This is trivial for $\phi_0$.  

Let $n>0$.  Verifying that $\phi_n$ is a homomorphism is
straightforward but tedious.  For example, if
$(a^n(d),r,x_n),(a^n(d'),r',x_n)\in S_n$ such that $r+r' >n$ then let
$m= d(2n+1)+r$ and $m' = d'(2n+1)+r'$.  Observe that $2n\geq r+r'>n$ so
that $-n \leq r+r' -(2n+1) \leq n$. Next we compute
\[
m+m' = (d+d')(2n+1)+r+r' = (d+d'+1)(2n+1) + (r+r'-(2n+1)).
\]
Suppose $d' = d = 1$.  Then $m+m' = 3(2n+1) + (r+r'-(2n+1)) > 3n+1$
since $r+r' - (2n+1)\geq -n$.  It follows that, using our version of
modulo addition on $\Z_{6n+3}$,  
\begin{align*}
m+m'\mod 6n+3 &= m+m' - (6n+3) \\
&= (d+d'+1-3)(2n+1) + (r+r'-(2n+1)) \\
&= (d+d+1 \mod 3)(2n+1) +(r+r'-(2n+1)).
\end{align*}
The cases for the other possible values of $d$ and $d'$ are similar
and in general we have 
\[
m+m' \mod 6n+3 = (d+d'+1\mod 3)(2n+1) + (r+r'-(2n+1)).
\]
Now, observe that in $T_n$
\[
(m,x_n)+(m',x_n) = (m+m'\mod 6n+3, x_n).
\]
By our
construction $-n\leq r+r'-(2n+1)\leq n$ and $-1\leq d+d'+1\mod 3 \leq
1$ so that, by the definition of $\phi_n$ and \eqref{eq:18}, 
\begin{align*}
\phi_n(m+m'\mod 6n+3,x_n) &= (a^n(d+d'+1\mod 3),r+r'-(2n+1),x_n) \\
&=(a^n(d),r,x_n)+(a^n(d'),r',x_n).
\end{align*}
Of course, we have only shown that $\phi_n$ preserves addition when
$r+r' > n$.  However,  the computations for the other cases are
analogous and it is straightforward, but tedious, to see that $\phi_n$
respects the multiplication operations on $S_n$ and $T_n$.    

Proving that $\phi_n$ preserves inversion is much easier.
Suppose $(m,x_n)\in T_n$.  Then if $m = d(2n+1)+r$ is the
decomposition given by Lemma \ref{lem:6} we have $-m = -d(2n+1)-r$.
Since $-a^n(d) =a^n(-d)$ we have
\[
\phi_n(-m,x_n) = (a^n(-d),-r,x_n) = -(a^n(d),r,x_n).
\]
Thus $\phi_n$ respects the inverse operation.  Since $\phi$ is a
bijection which respects the operations on $S_n$ and $T_n$ the fact
that $T_n$ is an abelian group implies that $S_n$ is also an  abelian 
group and that
$\phi_n$ is an isomorphism.  
\end{proof}

We have been stepping around it for the last few propositions, but it
is time to put everything together and show that we can form a group
bundle out of the $S_n$.  

\begin{prop}
With the operations on $S_n$ defined in Proposition \ref{prop:32} $S$
is an abelian, continuously varying, second countable, locally compact Hausdorff group bundle with
bundle map $p_S$.  
\end{prop}

\begin{proof}
First, it follows from Proposition \ref{prop:28} that $S$ is second
countable, locally compact Hausdorff.  It's easy to see from the way
we defined $S$ as the disjoint union of groups over $X$ that $S$ is a
group bundle with unit space $X$ and bundle map $p_S$.  
Furthermore,  each $S_n$ is abelian so that $S$
is abelian.  All that is left is to show that the operations are
continuous.  Suppose $(a^{n_i}(d_i),r_i,x_{n_i})\rightarrow
(a^N(d),r,x_N)$ and $(a^{n_i}(d'_i),r'_i,x_{n_i})\rightarrow
(a^N(d'),r',x_N)$ in $S$.  We would like to show that 
\begin{align*}
(a^{n_i}(d_i),r_i,x_{n_i})+(a^{n_i}(d_i'),r_i',x_{n_i})&\rightarrow
(a^N(d),r,x_N)+(a^N(d'),r',x_N),\quad\text{and} \\
-(a^{n_i}(d_i),r_i,x_{n_i})&\rightarrow -(a^N(d),r,x_N).
\end{align*}
It will suffice to show that for every subsequence, we can pass to a
sub-subsequence and obtain the required convergence.  So, pass to a
subsequence and use Lemma \ref{lem:5} to pass to another subsequence
and assume that either $n_i = N$ for all $i$ or that
$n_i\rightarrow\infty$, $n_i\ne 0$ for all $i$, and $N=0$.  Consider
the first case.  Eventually $r_i = r$ and $r_i' =r'$.  Furthermore,
since $A_N$ is a discrete space containing
$\{a^{n_i}(d_i)\}$ and $\{a^{n_i}(d'_i)\}$ we also eventually have
$a^{n_i}(d_i)=a^N(d)$ and $a^{n_i}(d'_i) = a^N(d')$.  Hence,
eventually, $d_i=d$ and $d_i' = d'$.  However, this implies that for
very large $i$
\begin{align*}
(a^{n_i}(d_i),r_i,x_{n_i})+(a^{n_i}(d'_i),r'_i,x_{n_i}) &=
(a^N(d),r,x_N)+(a^N(d'),r',x_N),\quad\text{and} \\
-(a^{n_i}(d_i),r_i,x_{n_i}) &= -(a^N(d),r,x_N)
\end{align*}
and at this point convergence is clear.  

Next, consider the second case so that $n_i\rightarrow \infty$ and
$n_i \ne 0$ for all $i$.  Given any arbitrary sequence $\{c_i\}\subset \Z_3$ we
know that $-1/n_i\leq a^{n_i}(c_i) \leq 1/n_i$ and this implies that
$a^{n_i}(c_i)\rightarrow 0$.  Since $r_i=r$ eventually it follows from
\eqref{eq:17} that 
\[
-(a^{n_i}(d_i),r_i,x_{n_i}) = (a^{n_i}(-d_i),-r_i,x_{n_i}) \rightarrow
(0,-r,x_0) = -(a^0(d),r,x_0).
\]
Next, eventually $r_i = r$ and $r_i'=r'$ and $|r+r'|<n_i$ so that we
can pass to a subsequence and assume this always holds. Furthermore,
let $c_i = d_i + d_i'\mod 3$ for all $i$.  Then, by \eqref{eq:18}, we
have
\[
(a^{n_i}(d_i),r_i,x_{n_i})+(a^{n_i}(d_i'),r_i',x_{n_i}) =
(a^{n_i}(c_i),r+r',x_{n_i}).
\]
As before, whatever the $c_i$, we know that $a^{n_i}(c_i)\rightarrow
0$.  This implies that 
\[
(a^{n_i}(c_i),r+r',x_{n_i})\rightarrow (0,r+r',x_0) = (a^0(d),r,x_0) +
(a^0(d'),r',x_0).
\]
Thus, both of the required sequences converge and it follows that the
operations on $S$ are continuous.  This makes $S$ into a topological 
groupoid, and since we have shown that $p_S$
is open, $S$ is continuously varying. 
\end{proof}

\begin{remark}
It is easy enough to see that $S$ is $r$-discrete so that 
the Haar system for $S$ is given on
$S_n$ by $\mu\times \delta_{x_n}$ where $\mu$ is counting measure on
$A_n\times \Z_{2n+1}$ (or $A_0\times\Z$) and $\delta_{x_n}$ is the
Dirac delta measure at $x_n$. 
\end{remark}

We now prove the main result of the section. 

\begin{proof}[Proof of Theorem \ref{thm:counter}]
Let $S$ be as in Proposition \ref{prop:28} and $T$ be as in
Proposition \ref{prop:25}.  Define $\phi:T \rightarrow S$ such that
$\phi(m,x_n) = \phi_n(m,x_n)$ for all $(m,x_n)\in T$.  
Algebraically everything is straightforward.  It is clear that we can
define such a map on $T$ and that, since each $\phi_n$ is a group
homomorphism, the resulting $\phi$ will be a groupoid homomorphism.
Furthermore, each $\phi_n$ is bijective so $\phi$ is bijective.  Let
us show that it is continuous.  Suppose $(m_i,x_{n_i})\rightarrow
(m,x_N)$.  As before it suffices to show that given a subsequence of
$\phi(m_i,x_{n_i})$ we can find a sub-subsequence converging to
$\phi(m,x_N)$.  So let us pass to a subsequence, and then do so again
to assume that $m_i = m$ for all $i$.  Now use Lemma
\ref{lem:5} to pass to yet another subsequence and assume that either
$n_i = N$ for all $i$ or that $n_i\rightarrow\infty$, $n_i\ne 0$ for
all $i$, and that $N=0$.  When the former is true $(m_i,x_{n_i})$ is a
constant sequence so of course $\phi(m_i,x_{n_i})$ converges to
$\phi(m,x_N)$.  Suppose the latter is true.  Eventually $|m|\leq n_i$
and when this happens $\phi(m,x_{n_i}) = (0,m,x_{n_i})$.  Clearly
$(0,m,x_{n_i})\rightarrow (0,m,x_0)$ and therefore $\phi$ is
continuous.  

Now consider the sequence $(2n+1,x_n)$.  It is easy to see that
$\phi(2n+1,x_n) = (1/n,0,x_n)$ so that clearly 
\[
\phi(2n+1,x_n) = (1/n,0,x_n)\rightarrow (0,0,x_0) = \phi(0,x_0).
\]
However, the sequence $(2n+1,x_n)$ does not converge to anything in
$T$.  It follows that $\phi$ is not a homeomorphism. 
\end{proof}


\chapter{Groupoid Crossed Products}
\label{cha:crossed}
In this chapter we give the definition of a groupoid dynamical system
and construct the groupoid crossed product.  Unfortunately, many
elements of the construction are rather technical and we will need to
draw upon a wealth of existing mathematics.  In Section
\ref{sec:cstarbundles} we give a brief overview of
upper-semicontinuous bundles and their relation to $C_0(X)$-algebras.
In Section \ref{sec:dynamical} we define a groupoid dynamical
system and construct the function algebra from which we will build the
crossed product.  Section \ref{sec:covariant} concerns covariant
representations.  In order to properly define a covariant
representation we will need to deal with both groupoid representations
and decompositions of representations of $C^*$-algebras.  After
developing these tools we then construct the crossed product
in Section \ref{sec:crossedprod}.  The most important result in this
section is Renault's Disintegration Theorem which will free us from
having to deal with covariant representations directly. In particular,
the beginner should not be discouraged if they don't understand all of
Section \ref{sec:covariant} on their first read through.  

\section{Upper-semicontinuous Bundles}
\label{sec:cstarbundles}
This section is essentially a collection of the important
results concerning upper-semicontinuous bundles that we will need
for the study of groupoid crossed products.  
Those readers unfamiliar with $C_0(X)$-algebras and their
related bundles are referred to \cite[Appendix C]{tfb2}.  This
reference is self-contained and does a very good job of covering the
basics of $C_0(X)$-algebra theory.  In fact, for the most part, the
definitions and theorems in this section are lifted from
\cite{tfb2} and we will cite a number of results from this reference
without proof.  All of this theory has its roots in \cite{dupregillette}.

Our main concern will be to develop a theory of bundles of
$C^*$-algebras.  However, in order to define our induction techniques
in Section \ref{sec:indreps} we will need to start with something a
little more general. 

\begin{definition}
\label{def:26}
\index{upper-semicontinuous!Banach bundle}
An {\em upper-semicontinuous Banach bundle} over a locally compact
Hausdorff space
$X$ is a topological space $\mcal{A}$ together with a continuous, open
surjection $p = p_\mcal{A}:\mcal{A}\rightarrow X$ and complex Banach
space structures on each fibre $\mcal{A}_x := p\inv(x)$ satisfying the
following axioms. 
\begin{enumerate}
\item The map $a\mapsto \|a\|$ is upper-semicontinuous from
  $\mcal{A}$ to $\R^+$.  (That is, for all $\epsilon > 0$, the set $\{
  a\in\mcal{A} : \|a\| \geq \epsilon\}$ is closed.)
\item If $\mcal{A}*\mcal{A} := \{(a,b)\in \mcal{A}\times\mcal{A} :
  p(a) = p(b)\}$, then $(a,b)\mapsto a+b$ is continuous from
  $\mcal{A}*\mcal{A}$ to $\mcal{A}$. 
\item For each $\lambda\in\C$, $a\mapsto \lambda a$ is continuous from
  $\mcal{A}$ to $\mcal{A}$.  
\item If $\{a_i\}$ is a net in $\mcal{A}$ such that $p(a_i)\rightarrow
  x$ and such that $\|a_i\|\rightarrow 0$, then $a_i\rightarrow 0_x$
  (where $0_x$ is the zero element of $\mcal{A}_x$).  
\end{enumerate}
\end{definition}

The following proposition is something of a utility belt for dealing with 
upper-semicontinuous bundles.  In particular, the fourth part gives us
a handle on the topology of the total space, which can be difficult to
deal with.

\begin{prop}
\label{prop:35}
Suppose $\mcal{A}$ is an upper-semicontinuous Banach bundle over $X$
with bundle map $p$.  Then the following statements hold. 
\begin{enumerate}
\item If $a_i\rightarrow 0_x$ in $\mcal{A}$ then $\|a_i\|\rightarrow
  0$.  
\item For all $x\in X$ the topology of $A_x$ as a subset of $\mcal{A}$
  is exactly its norm topology as a Banach space.  
\item The map $(\lambda,a)\mapsto \lambda a$ is continuous from
  $\C\times\mcal{A}$ into $\mcal{A}$. 
\item Let $\{a_i\}$ be a net in $\mcal{A}$ such that
  $p(a_i)\rightarrow p(a)$ for some $a\in \mcal{A}$.  Suppose that for
  all $\epsilon > 0$ there is a net $\{u_i\}$ in $\mcal{A}$ and
  $u\in\mcal{A}$ such that 
\begin{enumerate}
\item $u_i\rightarrow u$ in $\mcal{A}$, 
\item $p(u_i) = p(a_i)$ for all $i$, 
\item $\|a-u\|<\epsilon$, and 
\item eventually $\|a_i-u_i\|<\epsilon$.   
\end{enumerate}
Then $a_i\rightarrow a$. 
\end{enumerate}
\end{prop}

\begin{proof}
Part {\bf (a)}:  Since the norm is upper-semicontinuous on
$\mcal{A}$ the set $\{a\in\mcal{A}:{\|a\|<\epsilon}\}$ is open for all
$\epsilon > 0$.  Thus we eventually have $\|a_i\|<\epsilon$ for all
$\epsilon > 0$ and the result is proved. 

Part {\bf (b)}:  Suppose that $a_i\rightarrow a$ in $\mcal{A}$ with
$p(a_i) = p(a)$ for all $i$.  Then $a_i-a \rightarrow 0_{p(a)}$ by the
continuity of addition and $\|a_i-a\|\rightarrow 0$ by part (a).
Conversely, if $\|a_i-a\|\rightarrow 0$ then $a_i-a\rightarrow
0_{p(a)}$ by the last axiom of Definition \ref{def:26}, and
$a_i\rightarrow a$ by the continuity of addition.  

Part {\bf (d)}:  Since $X$ is Hausdorff we must have $p(u) = p(a)$ so
that condition (iii) makes sense.  Pass to a subnet of $\{a_i\}$.  It
will suffice to show that there is sub-subnet converging to $a$.
Since $p$ is open, we can pass to a subnet, relabel, and find $c_i\in
\mcal{A}_{p(a_i)}$ such that $c_i\rightarrow a$.  Fix $\epsilon > 0$
and choose $u_i$ as in part (d).  Since addition is continuous,
$c_i-u_i \rightarrow a-u$ in $\mcal{A}$.  Since $\|a-u\|<\epsilon$ by
assumption, and since $\{b\in\mcal{A}:\|b\|<\epsilon\}$ is open, we
eventually have $\|c_i-u_i\|<\epsilon$.  The triangle inequality then
implies that we eventually have $\|a_i-c_i\|<2\epsilon$.  As $\epsilon$
was arbitrary, we've shown that $\|a_i-c_i\|\rightarrow 0$.  Therefore
axiom (d) implies that $a_i-c_i\rightarrow 0_{p(a)}$.  Thus
\[
a_i = (a_i-c_i)+c_i\rightarrow 0_{p(a)}+a = a.
\]

Part {\bf (c)}:  Suppose $a_i\rightarrow a$ in $\mcal{A}$ and
$\lambda_i\rightarrow \lambda$ in $\C$.  We will apply part (d) with
$u_i = \lambda a_i$ and $u=\lambda a$.  It is clear that
$p(\lambda_ia_i) = p(a_i)\rightarrow p(a) = p(\lambda a)$.  Suppose
$\epsilon > 0$.  By axiom (c) we have $\lambda a_i \rightarrow \lambda
a$.  Conditions (ii) and (iii) are both trivial.  For condition (iv),
first observe that the set $\{b\in\mcal{A} : \|b\| < \|a\|+1\}$ is
open by axiom (a).  Since this set contains $a$, it eventually
contains $a_i$.  It follows immediately that the set $\{\|a_i\|\}$ is
bounded.  Hence, because $\lambda_i\rightarrow \lambda$, we must have,
eventually, 
\[
\|\lambda_i a_i - \lambda a_i\| = |\lambda_i-\lambda|\|a_i\| <
\epsilon
\]
and therefore condition (iv) holds.  Thus $\lambda_ia_i\rightarrow
\lambda a$. 
\end{proof}

\begin{remark}
What makes the proof of part (c) so complicated is that formulas
like $\|\lambda_i a_i - \lambda a\|$ don't make sense because $a_i$
and $a$ could possibly live in different fibres. 
\end{remark}

Next, we can add structure to an upper-semicontinuous Banach bundle to
make it a bundle of $C^*$-algebras in the obvious way. 

\begin{definition}
\label{def:27}
\index{upper-semicontinuous!cstar-bundle@$C^*$-bundle}
An {\em upper-semicontinuous $C^*$-bundle} is an upper-semicontinuous
Banach bundle $p:\mcal{A}\rightarrow X$ such that each fibre is a
$C^*$-algebra and such that the following additional axioms hold. 
\begin{enumerate}
\item[{\bf(e)}] The map $(a,b)\mapsto ab$ is continuous from
  $\mcal{A}*\mcal{A}$ to $\mcal{A}$. 
\item[{\bf(f)}] The map $a\mapsto a^*$ is continuous from $\mcal{A}$
  to $\mcal{A}$.  
\end{enumerate}
\end{definition}

There is also the more restrictive notion of continuous bundles which
deserves to be mentioned.  

\begin{definition}
An  upper-semicontinuous Banach bundle
(resp. $C^*$-bundle) $\mcal{A}$ is a Banach bundle
(resp. $C^*$-bundle) if the map $a\mapsto \|a\|$ is continuous. 
\end{definition}

It may seem strange that we are working with upper-semicontinuous
bundles as opposed to continuous bundles.  However, we will see that,
at least in the $C^*$-algebraic case, upper-semicontinuous bundles are the more
natural object.  

\begin{definition}
Suppose $\mcal{A}$ and $\mcal{B}$ are upper-semicontinuous Banach
bundles 
(resp. $C^*$-bundles) over $X$ with bundle maps $p$ and $q$
respectively.  A continuous map $\phi:\mcal{A}\rightarrow\mcal{B}$ is called a
Banach bundle (resp. $C^*$-bundle) homomorphism if $q\circ \phi = p$ and for
each $x\in X$ the restriction $\phi_x :\mcal{A}_x\rightarrow
\mcal{B}_x$ is a Banach space (resp. $C^*$-algebra) homomorphism. 
A Banach bundle (resp. $C^*$-bundle) {\em isomorphism} is a 
bijective, bicontinuous, Banach bundle (resp. $C^*$-bundle) homomorphism.  
\end{definition}

Given an upper-semicontinuous bundle the primary object of interest
will be the space of sections.  

\begin{definition}
\label{def:28}
\index[not]{$\Gamma_0(X,\mcal{A}),\Gamma_c(X,\mcal{A})$}
Suppose $\mcal{A}$ is an upper-semicontinuous Banach bundle.  Then we
will denote the space of sections of the bundle map by
$\Gamma(X,\mcal{A})$.  Given $f\in \Gamma(X,\mcal{A})$ we say that $f$
{\em vanishes at infinity} if the set $\{x\in X: \|f(x)\| \geq
\epsilon \}$ is compact for all $\epsilon \geq 0$.  We 
will denote the subspace of sections which
vanish at infinity by $\Gamma_0(X,\mcal{A})$.  Furthermore, we will
let $\Gamma_c(X,\mcal{A})$ be the subspace of sections which have
compact support.  

We endow $\Gamma(X,\mcal{A})$ with the operations of
pointwise addition and pointwise scalar multiplication.  Furthermore we equip
$\Gamma_0(X,\mcal{A})$ with the uniform norm $\|f\|_\infty = \sup_{x\in
  X}\|f(x)\|$.  If $\mcal{A}$ is an upper-semicontinuous
$C^*$-bundle then we give $\Gamma(X,\mcal{A})$ the operations of
pointwise multiplication and involution. Finally, given $\phi\in C(X)$
and $f\in\Gamma(X,\mcal{A})$ we define the section $\phi\cdot f$ via
$\phi\cdot f(x) := \phi(x)f(x)$ for all $x\in X$.  
\end{definition}

\begin{remark}
It is not clear at the outset that there are any nontrivial sections
in $\Gamma(X,\mcal{A})$.  A bundle $\mcal{A}$ is said to have {\em
  enough sections} if given $x\in X$ and $a\in \mcal{A}_x$ there exists
$f\in \Gamma(X,\mcal{A})$ such that $f(x) = a$.  If we are dealing
with a Banach bundle then it is a result of Douady and
Soglio-H\'erault that there are enough sections \cite[Appendix
C]{felldoran}.  Hoffman has noted that the same is true for
upper-semicontinuous Banach bundles \cite{hoff1}, although the details remain
unpublished \cite{hoff2}.  
We will not need to worry about this because, as we show in
Remark \ref{rem:30}, 
in all of our examples there will obviously be enough sections. 
\end{remark}

The point of all this is that the objects in Definition \ref{def:28}
have fairly nice algebraic properties and will fill roles analogous to
$C(X)$, $C_0(X)$ and $C_c(X)$.  

\begin{prop}
Suppose $\mcal{A}$ is an upper-semicontinuous Banach bundle.  Then the
following hold. 
\begin{enumerate}
\item $\Gamma(X,\mcal{A})$ is a vector space with respect to the
  natural pointwise operations.  If $\mcal{A}$ is a $C^*$-bundle then
  $\Gamma(X,\mcal{A})$ is a $*$-algebra. 
\item $\Gamma_0(X,\mcal{A})$ is complete with respect to the uniform
  norm.  Furthermore, $\Gamma_0(X,\mcal{A})$ is closed under the
  pointwise operations so that it is a Banach space.  If $\mcal{A}$ is
  a $C^*$-bundle then $\Gamma_0(X,\mcal{A})$ is a $C^*$-algebra. 
\item Given $\phi\in C_0(X)$ and $f\in \Gamma_0(X,\mcal{A})$ we have
  $\phi\cdot\gamma \in \Gamma_0(X,\mcal{A})$ and in particular
  $\Gamma_0(X,\mcal{A})$ is a $C_0(X)$-module. 
\end{enumerate}
\end{prop}

\begin{proof}
The algebraic statements are all straightforward to check.  We will
content ourselves with showing that $\Gamma_0(X,\mcal{A})$ is
complete.  Suppose $f_i$ is a Cauchy sequence in
$\Gamma_0(X,\mcal{A})$.  Since each $\mcal{A}_x$ is complete we can at
least define a section $f:X\rightarrow \mcal{A}$ by $f(x) = \lim_i
f_i(x)$.  Now suppose $\epsilon > 0$ and choose $N$ such that $i,j\geq
N$ implies $\|f_i - f_j\|_\infty < \epsilon$.  Given $x\in X$ pick
$i_x\geq N$ so that $\|f_{i_x}(x) - f(x)\| < \epsilon$.  Then for all
$x\in X$ and $i \geq N$ we have 
\[
\|f_i(x) - f(x) \| \leq \|f_i(x) - f_{i_x}(x)\| + \|f_{i_x}(x) -
f(x)\| < 2\epsilon.
\]
It follows that $\|f_i-f\|_\infty \rightarrow 0$.  
We need to show $f$ is continuous.  Suppose $x_i\rightarrow x $ and
fix $\epsilon > 0$.  Choose $N$ so that $\|f_N-f\|<\epsilon$.  Since
$f_N(x_i)\rightarrow f_N(x)$ we can let $a_i = f(x_i)$ and $u_i =
f_N(x_i)$ and then part (d) of Proposition \ref{prop:35} implies that
$f(x_i)\rightarrow f(x)$.  Next, since $x\mapsto
\|f(x)\|$ is the uniform limit of functions which vanish at infinity, it's easy
to see that $x\mapsto \|f(x)\|$ vanishes at infinity and therefore
$f\in \Gamma_0(X,\mcal{A})$.  
\end{proof}

The following proposition gives us another nice tool for dealing with
the topology on the total space. It also shows that the topology on
$\mcal{A}$ is determined by its space of sections.  

\begin{prop}
\label{prop:37}
Let $p:\mcal{A}\rightarrow X$ be an upper-semicontinuous Banach
bundle.  Suppose that $\{a_i\}$ is a net in $\mcal{A}$, that $a\in
\mcal{A}$, and that $f\in \Gamma_0(X,\mcal{A})$ is such that $f(p(a))
= a$.  If $p(a_i)\rightarrow p(a)$ and if
$\|a_i-f(p(a_i))\|\rightarrow 0$ then $a_i\rightarrow a$. 
\end{prop}

\begin{proof}
We have $a_i - f(p(a_i)) \rightarrow 0$ by axiom (d) of Definition
\ref{def:26}.  However, since $f$ is continuous we also have
$f(p(a_i))\rightarrow f(p(a))=a$.  Hence
\[
a_i = (a_i-f(p(a_i)))+f(p(a_i)) \rightarrow 0_{p(a)} + a = a. \qedhere
\]
\end{proof}

The following proposition is important because it gives a very
convenient criterion for a subspace of $\Gamma_0(X,\mcal{A})$ to be
dense.  This will be useful because we will often want to use some
dense subspace of particularly simple functions.   This is proved for
upper-semicontinuous $C^*$-bundles in \cite[Proposition C.24]{tfb2}
and the extension to Banach bundles is basically the same. 

\begin{prop}
\label{prop:42}
Suppose $p:\mcal{A}\rightarrow X$ is an upper-semicontinuous Banach
bundle over $X$ and $\Gamma$ is a subspace of $\Gamma_0(X,\mcal{A})$
such that 
\begin{enumerate}
\item $f\in \Gamma$ and $\phi\in C_0(X)$ implies $\phi\cdot f\in
  \Gamma$, and 
\item for each $x\in X$ the set $\{f(x):f\in\Gamma\}$ is dense in
  $\mcal{A}_x$.  
\end{enumerate}
Then $\Gamma$ is dense in $\Gamma_0(X,\mcal{A})$.  
\end{prop}

\begin{proof}
Fix $f\in \Gamma_0(\mcal{A})$ and $\epsilon > 0$.  We need to find
$g\in \Gamma$ such that $\|f-g\|_\infty < \epsilon$.  Let $K$ be the
compact set $\{x\in X:\|f(x)\| \geq \epsilon/3\}$.  Given $x\in K$,
there is a $g\in \Gamma$ such that $\|f(x)-g(x)\|< \epsilon/3$.  Using
upper-semicontinuity, there is a neighborhood $U$ of $x$ such that 
\[
\|f(y)-g(y)\| < \epsilon/3\ \text{if $y\in U$}.
\]
Since $K$ is compact, there is a cover $U_1,\ldots,U_n$ of $K$ by open
sets with compact closure, and $g_i\in \Gamma$ such that 
\[
\|f(y)-g_i(y)\| < \epsilon/3\ \text{if $y\in U_i$}.
\]
Using \cite[Lemma 1.43]{tfb2} we can find a partition of unity
$\{\phi_i\}_{i=1}^n\subset C_c(X)$ such that $0\leq\phi_i(x)\leq 1$
for all $x\in X$, $\supp\phi_i\subset U_i$, if $x\in K$ then $\sum
\phi_i(x) = 1$, and if $x\not\in K$ then $\sum \phi_i(x) \leq 1$.  By
assumption, $\sum \phi_i\cdot g_i \in \Gamma$.  Now, if
$x\in K$ then 
\begin{align*}
\left\|f(x) - \sum_{i=1}^n \phi_i(x)g_i(x)\right\| &= 
\left\|\sum_{i=1}^n \phi_i(x)(f(x)-g_i(x))\right\| \\
&\leq \sum_{i=1}^n \phi_i(x)\|f(x)-g_i(x)\| \\
&\leq \epsilon/3 \leq \epsilon.
\end{align*}
But if $x\in U_i\setminus K$, then $\|g_i(x)\|< 2\epsilon/3$.  Since
$\supp \phi_i \subset U_i$, for any $x\not\in K$ we have
$\phi_i(x)\|g_i(x)\| \leq \frac{2\epsilon}{3} \phi_i(x)$.  Thus if $x\not\in
K$, we still have 
\begin{align*}
\left\| f(x) - \sum_i \phi_i(x)g_i(x)\right\| &\leq \|f(x)\| + \sum_i
\phi_i(x)\|g_i(x)\| \\
&\leq \frac{\epsilon}{3} + \frac{2\epsilon}{3} = \epsilon.
\end{align*}
Therefore $\sup_{x\in X}\|f(x)-(\sum \phi_i\cdot g_i)(x)\|< \epsilon$
as required. 
\end{proof}

\subsection{$C_0(X)$-algebras}
\label{sec:c0x-algebras}

The following objects play the same role for groupoid crossed products
that $C^*$-algebras do for group crossed products.  They 
will eventually explain our preference for
upper-semicontinuous bundles over continuous bundles.  

\begin{definition}
\index{czeroofxalgebra@$C_0(X)$-algebra}
Suppose that $A$ is a $C^*$-algebra and that $X$ is a locally compact
Hausdorff space.  Then $A$ is a {\em $C_0(X)$-algebra} if there is a
homomorphism $\Phi_A$ from $C_0(X)$ into the center of the
multiplier algebra $ZM(A)$ which is nondegenerate in that the set
\[
\Phi_A(C_0(X))\cdot A := \spn\{\Phi_A(f)a: f\in C_0(X), a\in A\}
\]
is dense in $A$.  
\end{definition}

\begin{remark}
Suppose $A$ is a $C_0(X)$-algebra, $B\subset A$, and $C\subset C_0(X)$.
We will use the notation 
\[
C\cdot B = \Phi_A(C)\cdot B := \spn\{\Phi_A(f)a : f\in C, a\in B\}.
\]
\end{remark}

Our eventual goal is to show that
there is a one-to-one correspondence between $C_0(X)$-algebras and 
upper-semicontinuous $C^*$-bundles.  For starters, what follows
next  shows how we can view $C_0(X)$-algebras as ``fibred''
objects. 

\begin{prop}
Suppose $A$ is a $C_0(X)$-algebra and $J$ is an ideal in $C_0(X)$.
Then the closure of $\Phi_A(J)\cdot A$ is an ideal in $A$.  
\end{prop}
\begin{proof}
Let $I$ be the closure of $\Phi_A(J)\cdot A$ and observe that $I$ is
just the closed linear span of $I_0 = \{\Phi(f)a : f\in J, a\in A\}$.
Therefore it will suffice to show that given $a\in A$ and $\Phi(f)b\in
I_0$ then $a(\Phi(f)b), (\Phi(f)b)a \in I_0$.  However, $\Phi(f)$ is
in the center of $M(A)$ so that $a(\Phi(f)b) = \Phi(f)(ab)$ and
$(\Phi(f)b)a = \Phi(f)(ba)$.  The result follows. 
\end{proof}

\begin{definition}
\label{def:29}
Suppose $A$ is a $C_0(X)$-algebra.  Given $x\in X$ let $J_x$ be the
ideal of functions in $C_0(X)$ which vanish at $x$.   Then we will
denote the ideal $\overline{\Phi_A(J_x)\cdot A}$ by $I_x$ and the
quotient $A/I_x$ by $A(x)$.  We think of $A(x)$ as the {\em fibre of
  $A$ over $x$} and given $a\in A$ we write $a(x)$ for the image of
$a$ in $A(X)$.  In this way we think of $a$ as a function from $X$
onto the disjoint union $\coprod_{x\in X} A(x)$.  
\end{definition}

The following are some particularly nice examples of
$C_0(X)$-algebras.  

\begin{example}
If $D$ is any $C^*$-algebra and $X$ is a locally compact Hausdorff
space then $A = C_0(X,D)$ is a $C_0(X)$-algebra in a natural way: 
\[
\Phi_A(f)(a)(x) := f(x)a(x)
\]
for all $f\in C_0(X)$ and $a\in A$.  In this case each fibre $A(x)$ is
easily identified with $D$ and the identification of the elements of
$A$ with functions on $X$ is the obvious one.  
\end{example}

\begin{example}
\label{ex:15}
Suppose that $X$ and $Y$ are locally compact Hausdorff spaces and that
$\phi:Y\rightarrow X$ is a continuous surjection.  Then $C_0(Y)$ becomes a
$C_0(X)$-algebra with respect to the map defined by
\[
\Phi_{C_0(Y)}(f)g(y) := f(\phi(y))g(y).
\]
The only issue is to see that $\Phi_{C_0(Y)}$ is nondegenerate, but
this is easy enough to do using the Stone-Weierstrass theorem.  In
this example, the fibres $C_0(Y)(x)$ are isomorphic to
$C_0(\phi\inv(x))$.  If $f\in C_0(Y)$ then $f(x)$ is just the
restriction of $f$ to $\phi\inv(x)$.  
\end{example}

\begin{example}
\label{ex:16}
Let $\mcal{A}$ be an upper-semicontinuous $C^*$-bundle and $A=
\Gamma_0(X,\mcal{A})$.  Then $A$ is a $C_0(X)$-algebra with respect to
the map defined by 
\[
\Phi_{A}(\phi)f(x) := \phi\cdot f(x) = \phi(x) f(x)
\]
for $\phi\in C_0(X)$ and $f\in A$.  This is really just
\cite[Proposition C.23]{tfb2}, however everything is fairly
straightforward to prove.  The only part that could be difficult is
the nondegeneracy but this is taken care of by Proposition
\ref{prop:42}.  It is also easy enough to show that in this
case $A(x) \cong \mcal{A}_x$ for all $x\in X$.  The isomorphism is
given on $A$ by evaluation at $x$ so that if $f\in A$ then
$f(x)$ as an element of $\mcal{A}_x$ is identified with $f(x)$ as an
element of the quotient $A(x)$.  
\end{example}

Next, we define the homomorphisms associated to $C_0(X)$-algebras.  In
particular we will show that they preserve the ``fibering'' process.  

\begin{definition}
\index{czeroofxlinear@$C_0(X)$-linear}
Suppose $A$ and $B$ are $C_0(X)$-algebras.  A map $\phi:A\rightarrow
B$ is called {\em $C_0(X)$-linear} 
if $\phi(\Phi_A(f)a) = \Phi_B(f)\phi(a)$ for all
$f\in C_0(X)$ and $a\in A$. 
\end{definition}

\begin{prop}
\label{prop:43}
Suppose $A$ and $B$ are $C_0(X)$-algebras and $\phi:A\rightarrow B$ is
a $C_0(X)$-linear homomorphism.  Then, for all $x\in X$, $\phi$ 
factors to a homomorphism
$\phi_x:A(x)\rightarrow B(x)$ such that $\phi_x(a(x))=\phi(a)(x)$. 
Furthermore, if $\phi$ is an isomorphism then each $\phi_x$
is as well. 
\end{prop}

\begin{proof}
Given $x\in X$ let $J_x$ be the ideal of functions on $C_0(X)$
vanishing at $x$.  Furthermore, 
let $I_x^A$ and $I_x^B$ be the ideals in $A$ and $B$,
respectively, such that $A(x) = A/I_x^A$ and $B(x) = B/I_x^B$.  We
would like to show that $\phi(I_x^A)\subset \phi(I_x^B)$.  Since
$\phi$ is a homomorphism and $I_x^A$ is the closure of the
set $J_x\cdot A$ it suffices to show that $\phi(f\cdot a)\in I_x^B$
for all $a\in A$ and $f\in J_x$.  However $\phi(f\cdot a) =
f\cdot\phi(a)$ and the result follows.  At this point we can compose
$\phi$ with the quotient map $b\mapsto b(x)$ and this will factor to a
homomorphism $\phi_x:A(x)\rightarrow B(x)$ defined via
$\phi_x(a(x))=\phi(a)(x)$.  Furthermore, if $\phi$ is an isomorphism
then $\phi\inv$ is $C_0(X)$-linear and we can construct $(\phi_x)\inv$.
However, it is straightforward to check that $(\phi_x)\inv =
\phi_x\inv$ so that in this case each $\phi_x$ is an isomorphism. 
\end{proof}

An essential fact about $C_0(X)$-algebras is that their primitive
ideal spaces are fibred over $X$, and that there is a tight
relationship between the action and this fibration.  

\begin{remark}
\index[not]{$f\cdot a$}
Recall from the Dauns-Hofmann Theorem \cite[Theorem A.34]{tfb} that
given a $C^*$-algebra $A$ there is an isomorphism $\Psi:C^b(\Prim
A)\rightarrow ZM(A)$ given as follows.  For $f\in C^b(\Prim A)$ and
$P\in \Prim A$ let
\begin{equation}
\label{eq:22}
(\Psi(f)(a))(P) := f(P)a(P)
\end{equation}
where $a(P)$ denotes the image of $a$ in the quotient $A/P$.  Then
\eqref{eq:22} defines a unique element of $A$ which we denote
$\Psi(f)(a)$.  In this way $\Psi(f)$ defines an element in the center
of the multiplier algebra $M(A)$. 
\end{remark}

\begin{prop}[{\cite[Proposition C.5]{tfb2}}]
\label{prop:36}
Suppose $A$ is a $C^*$-algebra and that $X$ is a locally compact Hausdorff
space.  If there is a continuous map $\sigma_A:\Prim A\rightarrow X$
then $A$ is a $C_0(X)$-algebra with 
\begin{equation}
\label{eq:23}
\Phi_A(f)a := \Psi(f\circ \sigma_A)a
\end{equation}
for all $f\in C_0(X)$ and $a\in A$.  Conversely, if $A$ is a $C_0(X)$-algebra
then there is a continuous map $\sigma_A:\Prim A \rightarrow X$ such
that \eqref{eq:23} holds. 

In particular, every irreducible representation of $A$ is lifted from
a fibre $A(x)$ for some $x\in X$.  More precisely, if
$\pi\in\widehat{A}$ then the ideal $I_{\sigma_A(\ker\pi)}$ is
contained in $\ker\pi$ and $\pi$ is lifted from an irreducible
representation of $A(\sigma_A(\ker\pi))$.  In this way we can identify
$\widehat{A}$ with the disjoint union $\coprod_{x\in X}
A(x)\sidehat$.  
\end{prop}

Thus, the map $\sigma_A$ gives us our fibration of $\Prim A$. 

\begin{corr}
Suppose $A$ is a $C_0(X)$-algebra and $\sigma_A:\Prim A\rightarrow X$
is the map given in Proposition \ref{prop:36}.  We can view
$\Prim A$ as a bundle over $X$ and the fibre $\sigma_A\inv(x)$ can be
identified with $\Prim A(x)$ for all $x\in X$. 
\end{corr}
\begin{proof}
This is nothing more than a restatement of the second part of 
Proposition \ref{prop:36} in terms of primitive ideals.  In
particular, given $P\in \Prim A$ choose any $\pi\in \widehat{A}$ such
that $P = \ker \pi$ and it follows that $I_{\sigma(P)}\subset P$ and
that $P$ is lifted from an element of $\Prim A(x)$.  
\end{proof}

Proposition \ref{prop:36} allows us to present another example of a
$C_0(X)$-algebra that will be particularly important in Section
\ref{sec:locally-unitary}. 

\begin{example}
\label{ex:21}
Suppose $A$ is a $C^*$-algebra with Hausdorff spectrum $\widehat{A}$.  
Since the spectrum is always locally compact it follows that
$\widehat{A}$ 
is locally compact Hausdorff.  It is straightforward to show \cite[Lemma
5.1]{tfb} that the map $\pi\mapsto \ker\pi$ induces a homeomorphism of
$\widehat{A}$ onto $\Prim A$.  Therefore, if we identify $\Prim A$
with $\widehat{A}$ via this map, then
$\sigma_A = \id$ allows us to view $A$ as a $C_0(\widehat{A})$-algebra.  Given
$f\in C_0(\widehat{A})$ we combine \eqref{eq:22} and \eqref{eq:23} to get
\[
\Phi_A(f)a(\pi) = f(\pi)a(\pi)
\]
where $a(\pi)$ is the image of $a$ in $A/\ker\pi$.  From here it is
straightforward to identify the fibres as $A(\pi) = A/\ker\pi$.  It also
follows from \cite[Lemma 5.1]{tfb} (and is easy to show directly) that
each fibre $A(x)$ is simple and has, up to equivalence, a unique faithful
irreducible representation.  Moreover, in the separable case each
$A(x)$ is elementary
\end{example}

Moving on, the ``fibration'' of $A$ given by a $C_0(X)$-action is much more
rigorous than one might think.  The key link between $C_0(X)$-algebras
and upper-semicontinuous $C^*$-bundles is given by the following
theorem, which is, more or less, a summary of the results in
\cite[Appendix C]{tfb2}.  It also justifies our preference for
upper-semicontinuous bundles since there are many well behaved 
$C_0(X)$-algebras for which the map
$\sigma:\Prim A\rightarrow X$ is not open.

\begin{theorem}[{\cite[Theorem C.26]{tfb2}}]
\label{thm:c0xalgs}
\index{czeroofxalgebra@$C_0(X)$-algebra}
\index{upper-semicontinuous!cstar-bundle@$C^*$-bundle}
Suppose $A$ is a $C^*$-algebra.  Then the following statements are
equivalent.  
\begin{enumerate}
\item $A$ is a $C_0(X)$-algebra. 
\item There is a continuous map $\sigma_A : \Prim A\rightarrow X$. 
\item There is an upper-semicontinuous $C^*$-bundle
  $p:\mcal{A}\rightarrow X$ over $X$ and a $C_0(X)$-linear isomorphism
  of $A$ onto $\Gamma_0(X,\mcal{A})$.
\end{enumerate}
Moreover, $\mcal{A}$ is a $C^*$-bundle over $X$ if and only
if $\sigma_A$ is open.  
\end{theorem}

The following corollary is nothing more than a basic rehashing of
Theorem \ref{thm:c0xalgs}.  It is important, however, because it
presents the view of $C_0(X)$-algebras and upper-semicontinuous
$C^*$-bundles that we will use from now on.  

\begin{corr}
\label{cor:1}
Suppose $A$ is a $C_0(X)$-algebra.  Then we can endow the disjoint
union $\mcal{A} = \coprod_{x\in X} A(x)$ with a unique topology which
makes it
into an upper-semi\-continuous $C^*$-bundle such that the map which
sends $a\in A$ to the section $x\mapsto a(x)$ is a $C_0(X)$-linear
isomorphism of $A$ onto $\Gamma_0(X,\mcal{A})$.  Moreover, every
upper-semi\-continuous $C^*$-bundle can be obtained in this fashion. 
\end{corr}

\begin{proof}
Suppose $A$ is a $C_0(X)$-algebra and $\mcal{A}$ is defined as above.
Let $\mcal{B}$ be an upper-semicontinuous $C^*$-bundle such that there
is a $C_0(X)$-linear isomorphism $\phi:A\rightarrow
\Gamma_0(X,\mcal{B})$.  
First, we use the canonical action of $C_0(X)$ on $B = \Gamma_0(X,\mcal{B})$
to view $B$ as a $C_0(X)$-algebra.  Given $x\in X$ let $I_x^A$ be the
ideal in $A$ generated by $\Phi_A(J_x)\cdot A$ and $I_x^B$ the ideal
in $B$ generated by $\Phi_B(J_x)\cdot B$.  For a little while we will
use the notation $a+I_x^A := a(x)$ since we don't want to confuse
elements of quotients with function evaluation.  It follows from
Proposition \ref{prop:43}
that $\phi$ factors to an isomorphism $\phi_x:A(x)\rightarrow B(x)$
which is defined via $\phi_x(a+I_x^A) = \phi(a)+I_x^B$ for all $a+I_x^A\in A(x)$.
Next, is not hard to 
check that $I_x^B: =\{ f\in \Gamma_0(X,\mcal{B}): f(x) = 0\}$ and that
$\psi_x(f+I_x^B) = f(x)$ defines an isomorphism of $B(x)$ onto
$\mcal{B}_x$.  

Now we can define $\Omega: \mcal{A}\rightarrow\mcal{B}$
by $\Omega(a) = \psi_{p(a)}(\phi_{p(a)}(a))$ for all $a\in\mcal{A}$.  Once we
sort out all of the definitions it is easy to see that $\Omega$ is a
bijection and that restricted to a fibre $\Omega_x = \psi_x\circ
\phi_x$ is a $C^*$-algebra isomorphism.  It follows that we can pull
back the topology from $\mcal{B}$ to $\mcal{A}$ and, with this topology,
$\mcal{A}$ will be an upper-semicontinuous $C^*$-bundle.
Furthermore, in this situation $\Omega$ will be a $C^*$-bundle isomorphism.  

Now we have to see that sections have the right form.  It is
straightforward to show that the bundle
isomorphism $\Omega$ induces an isomorphism $\omega :
\Gamma_0(X,\mcal{B})\rightarrow \Gamma_0(X,\mcal{A})$ by $\omega(f) =
\Omega\inv\circ f$ for all $f\in \Gamma_0(X,\mcal{B})$.  We can compose
$\omega$ with $\phi$ to conclude that $A$ is isomorphic to 
$\Gamma_0(X,\mcal{A})$ and, furthermore, we can calculate
\begin{align*}
\omega\circ\phi(a)(x) &= \omega(\phi(a))(x) =
\Omega\inv(\phi(a)(x)) \\
&= \phi_x\inv\circ\psi_x\inv(\phi(a)(x)) 
= \phi_x\inv(\phi(a)+I^B_x) \\
&= a+I^A_x.
\end{align*}
However, reverting back to our former notation, this implies that 
$\omega\circ\phi(a)$ is exactly the section which sends $x$ to
$a(x)$.  

Next, let $\Upsilon: A\rightarrow \Gamma_0(X,\mcal{A})$ be given by
$\Upsilon(a)(x) = a(x)$ for all $a\in A$ and $x\in X$.  
Suppose that $\mcal{A}'$ is equal to $\mcal{A}$ as a set but has
a different topology such that $\Upsilon$ is a $C_0(X)$-isomorphism onto
$\Gamma_0(X,\mcal{A}')$.  We will use $\Upsilon'$ to denote this new
isomorphism.  It follows from Proposition \ref{prop:43} that for each $x\in X$ 
both $\Upsilon$ and $\Upsilon'$
factor to isomorphisms from $A(x)$ to $\mcal{A}_x$ and $\mcal{A}'_x$,
respectively.  Thus, fibrewise $\mcal{A}$ and $\mcal{A}'$ have the same
norm.  Suppose $a_i\rightarrow a$ in
$\mcal{A}$ and let $u_i = p(a_i)$ and $u= p(a)$.  Choose $b\in A$ such
that $b(u) = a$.  Observe that $u_i \rightarrow u$ and that, by
viewing $b$ as a continuous section of $\mcal{A}$ we have $a_i -
b(u_i) \rightarrow 0$.  It follows from Proposition \ref{prop:35} that
$\|a_i-b(u_i)\| \rightarrow 0$.  However, by using $\Upsilon'$ to view
$b$ as a section of $\mcal{A}'$, it follows from Proposition
\ref{prop:37} that $a_i\rightarrow a$ in $\mcal{A}'$.  Since the
situation is entirely symmetric this implies that the topology on
$\mcal{A}$ is unique.  

Finally, the fact that every upper-semicontinuous
$C^*$-bundle can be obtained in this fashion is an implication of the
equivalence in Theorem \ref{thm:c0xalgs}.
\end{proof}

\begin{definition}
\label{def:31}
\index{upper-semicontinuous!associated to a $C_0(X)$-algebra}
Given a $C_0(X)$-algebra $A$ we define the upper-semicontinuous
$C^*$-bundle {\em associated to $A$} to be $\mcal{A} = \coprod_{x\in
  X} A(x)$ with the topology from Corollary \ref{cor:1}.
\end{definition}

\begin{remark}
\label{rem:30}
Observe that if $A$ is a $C_0(X)$-algebra and $\mcal{A}$ is the
upper-semicontinuous bundle associated to $A$ then $\mcal{A}$ has
enough sections.  Indeed, if $a\in A(x)$ then we can view $A(x)$ as a
quotient of $A$ to find $b\in A$ such that $b(x)=a$.  However, we can
also view $b$ as a section in $\Gamma_0(X,\mcal{A})$ which takes on
the value $a$ at $x$. 
\end{remark}

\begin{remark}
We will need to make sure we don't confuse the $C_0(X)$-algebra
$A$ with its associated bundle $\mcal{A}$.  One reason we must
do this is because the topology on $\mcal{A}$ is not at all
straightforward and we will need to be extra careful in dealing with
it. For instance, $\mcal{A}$ may not even be Hausdorff \cite[Example
C.27]{tfb2}.  (Although it turns out that $\mcal{A}$ has to be
Hausdorff if it is a continuous bundle.)
\end{remark}

This duality between upper-semicontinuous bundles and
$C_0(X)$-algebras allows us to construct a similar duality between the
homomorphisms of these two categories. 

\begin{prop}
\label{prop:38}
Suppose $A$ and $B$ are $C_0(X)$-algebras and $\mcal{A}$ and
$\mcal{B}$ are the associated upper-semicontinuous bundles.  Then a
$C_0(X)$-linear homomorphism $\phi:A\rightarrow B$ induces a
$C^*$-bundle homomorphism $\hat{\phi}:\mcal{A}\rightarrow\mcal{B}$ via
$\hat{\phi}(a(x))= \phi(a)(x)$ for all $a(x)\in \mcal{A}$.

Conversely, a $C^*$-bundle homomorphism
$\phi:\mcal{A}\rightarrow\mcal{B}$ induces a $C_0(X)$-linear
homomorphism $\check{\phi}:A\rightarrow B$ where $\check{\phi}(a)$ is
uniquely determined by the relation $\check{\phi}(a)(x) =
\phi(a(x))$ for all $x\in X$.  
\end{prop}

\begin{proof}
This is really a matter of sorting out definitions. Given a
$C_0(X)$-linear map $\phi:A\rightarrow B$ it follows from Proposition
\ref{prop:43} that, for
each $x\in X$, there is a well defined homomorphism $\phi_x:A(x)\rightarrow
B(x)$ defined by $\phi_x(a(x)) = \phi_x(a)(x)$.  We can glue each of
these homomorphisms together to get the map
$\hat{\phi}:\mcal{A}\rightarrow \mcal{B}$.  It is clear that
$\hat{\phi}$ preserves fibres and that restricted to fibres
$\hat{\phi}$ is a homomorphism.  All we need to do is show
that $\hat{\phi}$ is continuous. Suppose $b_i\rightarrow b$ in
$\mcal{A}$ and let $x_i = p(b_i)$ and $x = p(b)$.  
Lift $b$ from the quotient $A(x)$ to find $a\in A$ such that $a(x) =
b$.  First, observe that because $p$ is continuous $x_i\rightarrow
x$.  Next, observe that $b_i - a(x_i) \rightarrow 0_x$ so that by
Proposition \ref{prop:35} $\|b_i - a(x_i)\|\rightarrow 0$.  Since
$\phi_x$ is contractive for all $x$, we have $\|\phi_{x_i}(b_i) -
\phi_{x_i}(a(x_i))\| \leq \|b_i - a(x_i)\|$ so that, using the
definition of $\phi_{x_i}$, 
\[
\| \phi_{x_i}(b_i) - \phi(a)(x_i) \| \rightarrow 0.
\]
However, $\phi(a)$ is a section of $\mcal{B}$ such that $\phi(a)(x) =
\phi_x(a(x))= \phi_x(b)$ so that it follows from
Proposition \ref{prop:37} that $\phi_{x_i}(b_i)\rightarrow \phi_{x}(b)$.  
For the reverse direction, identify $A$ and $B$ as the section
algebras of $\mcal{A}$ and $\mcal{B}$ respectively and 
define $\check{\phi}:A\rightarrow B$ by $\check{\phi}(a) = \phi\circ
a$.  The result follows without too much difficulty. 
\end{proof}

\begin{corr}
\label{cor:2}
Suppose $A$ and $B$ are $C_0(X)$-algebras and $\mcal{A}$ and
$\mcal{B}$ are the associated upper-semicontinuous bundles.  If
$\phi:A\rightarrow B$ is a $C_0(X)$-linear isomorphism then
$\hat{\phi}$ is a $C^*$-bundle isomorphism.  Conversely, if $\phi:
\mcal{A}\rightarrow \mcal{B}$ is a $C^*$-bundle isomorphism then
$\check{\phi}$ is a $C_0(X)$-linear isomorphism.  
\end{corr}

\begin{proof}
For the first direction, use Proposition \ref{prop:38} on both $\phi$
and $\phi\inv$.  Then use the characterization of  $\hat{\phi}$ and
$\hat{\phi}\inv$ to show that these maps are inverses.  The other
direction is exactly the same. 
\end{proof}

\begin{remark}
It follows from Corollary \ref{cor:2} that two 
$C_0(X)$-algebras are $C_0(X)$-isomorphic if and only if their associated
bundles $\mcal{A}$ and $\mcal{B}$ are isomorphic. Thus, citing Theorem
\ref{thm:c0xalgs}, the map sending $A$ to its associated bundle is a 
bijection between isomorphism classes of
$C_0(X)$-algebras and upper-semicontinuous $C^*$-bundles. 
\end{remark}

\subsection{Pull Back Bundles}

The last bit of $C_0(X)$-algebra theory that we need is the notion of
a pull back.  

\begin{definition}
\label{def:30}
\index{pull back}
Suppose $X$ and $Y$ are locally compact Hausdorff spaces, 
$\mcal{A}$ is an upper-semicontinuous Banach bundle over $X$,
and that $\tau:Y\rightarrow X$ is continuous.   The {\em
  pull back} of $\mcal{A}$ is defined to be the set 
\[
\tau^*\mcal{A} = \{(y,a)\in Y\times\mcal{A} : \tau(y) = p(a)\}.
\]
In this case $\tau^*\mcal{A}$ is equipped with the relative topology and
the bundle map $q:\tau^*\mcal{A}\rightarrow Y$ defined by
$q(y,a) = y$.
\end{definition}

Of course, we made this definition with every intention of proving the
following 

\begin{prop}
\label{prop:39}
Suppose $X$ and $Y$ are locally compact Hausdorff spaces, $\mcal{A}$ is
an upper-semicontinuous Banach bundle, and $\tau:Y\rightarrow X$ is
continuous.  Then the pull back $\tau^*\mcal{A}$ is an
upper-semicontinuous Banach bundle.  What's more, $\tau^*\mcal{A}$ is
an upper-semicontinuous $C^*$-bundle if $\mcal{A}$ is, and
if $\mcal{A}$ is a continuous bundle then $\tau^*\mcal{A}$ is as well.  
\end{prop}
\begin{proof}
First, observe that $\tau^*\mcal{A}_y$ can be easily identified with
$\mcal{A}_{\tau(y)}$ so that we can give $\tau^*\mcal{A}_y$ whatever
structure $\mcal{A}_{\tau(y)}$ has.  Next, note that the bundle map
$q:\tau^*\mcal{A}\rightarrow Y$ is continuous since it's the
restriction of a continuous map.  Let us show that it is open.
Suppose $y_i\rightarrow y$ in $Y$ and $a\in \mcal{A}_{\tau(y)}$.  Then
$\tau(y_i)\rightarrow \tau(y)$ and we can use the fact that 
the bundle map for $\mcal{A}$ is
open to pass to a subnet, relabel, and find $a_i \in
\mcal{A}_{\tau(y_i)}$ such that $a_i\rightarrow a$.  It follows that
$(y_i,a_i)\rightarrow (y,a)$ so that $q$ is open. 

All that is left is to verify the various bundle axioms.  The axioms
concerning the continuity of the operations are straightforward, as is
axiom (d).  We will content ourselves with showing that axiom (a)
holds.  Suppose $\epsilon > 0$.  We would like to show that the set
$C=\{(y,a)\in \tau^*\mcal{A} : \|a\|\geq \epsilon\}$ is closed.  Suppose
$\{(y_i,a_i)\}$ is a net in $C$ and that $(y_i,a_i)\rightarrow
(y,a)$.  Since $\mcal{A}$ is an upper-semicontinuous bundle $\|a\|\geq
\epsilon$ and we are done.  Finally, if $a\mapsto\|a\|$ is continuous
then clearly its composition with $(y,a)\mapsto a$ is continuous.  
\end{proof}

\begin{prop}
\label{prop:40}
Suppose $X$ and $Y$ are locally compact Hausdorff spaces, $\mcal{A}$
is an upper-semicontinuous Banach bundle, and $\tau:Y\rightarrow X$ is
continuous.  Then $f\in \Gamma(Y,\tau^*\mcal{A})$ if and
only if there exists a continuous function $\tilde{f}:Y\rightarrow
\mcal{A}$ such that $p(\tilde{f}(y)) = \tau(y)$ and $f(y) =
(y,\tilde{f}(y))$ for all $y\in Y$.  
Furthermore, $\tilde{f}$ is compactly supported if
and only if $f$ is as well.
\end{prop}

\begin{proof}
Given $f\in \Gamma(Y,\tau^*\mcal{A})$ we define $\tilde{f}$ to be the
composition of $f$ with the projection from $\tau^*\mcal{A}$ onto
$\mcal{A}$.  Given a continuous $\tilde{f}:Y\rightarrow \mcal{A}$ such
that $p(\tilde{f}(y)) = \tau(y)$ we define $f$ by $f(y) =
(y,\tilde{f}(y))$.  Given $f$ and $\tilde{f}$ as in the statement of
the proposition it is clear that $\|f(y)\| = \|\tilde{f}(y)\|$.  It
follows immediately that $f$ is compactly supported if and only if
$\tilde{f}$ is as well.  
\end{proof}

\begin{remark}
\label{rem:8}
We
will often times denote the element $(y,a)\in\tau^*\mcal{A}$ by just
$a$ and will usually not distinguish between the maps $f$ and
$\tilde{f}$.
\end{remark}

\begin{definition}
\label{def:32}
\index{pull back}
Let $X$ and $Y$ be locally compact Hausdorff spaces, $A$ be a
$C_0(X)$-algebra, $\mcal{A}$ its associated
upper-semicontinuous $C^*$-bundle, and $\tau:X\rightarrow Y$ a
continuous map.  We define the {\em pull back} of
$A$ to be $\tau^* A := \Gamma_0(Y,\tau^*\mcal{A})$.  
\end{definition}

\begin{prop}
\label{prop:45}
Let $X$ and $Y$ be locally compact Hausdorff spaces, $A$ be a
$C_0(X)$-algebra, and $\tau:Y\rightarrow X$ a
continuous map.  Then there is a natural identification of
$\tau^*A(y)$ with $A(\tau(y))$ for all $y\in Y$. 
\end{prop}

\begin{proof}
This is really just working out the definitions.  Let $\mcal{A}$ be
the bundle associated to $A$ so that $\tau^*A =
\Gamma_0(Y,\tau^*\mcal{A})$.  Then, as we have seen in Example
\ref{ex:16}, 
$\tau^*A(y)= \tau^*\mcal{A}_y$.  
It follows, almost by definition, that $\tau^*\mcal{A}_y
= \mcal{A}_{\tau(y)} = A(\tau(y))$ and we are done. 
\end{proof}

\begin{remark}
When $\tau$ is a surjection, $\tau^* A$ is usually 
defined to be the balanced tensor product
$C_0(Y)\otimes_{C_0(X)} A$ where we view $C_0(Y)$ as a
$C_0(X)$-algebra as in Example \ref{ex:15}.  We will show that this
is equivalent to our definition in Section \ref{sec:tensor}.  
However, the following
proposition captures an important aspect of this identification. 
\end{remark}

\begin{prop}
\label{prop:46}
\index[not]{$f\otimes a$}
\index[not]{$C_c(Y)\odot A$}
Suppose $X$ and $Y$ are locally compact Hausdorff spaces, $A$ is a
$C_0(X)$-algebra, $\mcal{A}$ is its associated upper-semicontinuous
bundle, and $\tau:X\rightarrow Y$ is continuous.  
Given $f\in C_c(Y)$ and $a\in A$ define $f\otimes a(y)=
f(y) a(\tau(y))$ for all $y\in Y$.  Then $f\otimes a\in
\Gamma_c(Y,\tau^*\mcal{A})$ and
\[
C_c(Y)\odot A := \spn\{f\otimes a : f\in C_c(Y), a\in A\}.
\]
is dense in $\tau^*A$. 
\end{prop}

\begin{remark}
\index{elementary tensors}
We will often refer to elements of the form $f\otimes a$ as
{\em elementary tensors}, because, as we will see in Section  
\ref{sec:tensor}, they correspond to elementary tensors in a tensor product.
\end{remark}

\begin{proof}
Given $f$ and $a$ as above view $a$ as a section of the associated
bundle. Now define $g(y) = (y,f(y)a(\tau(y)))$.  Since everything in
sight is continuous, it is clear that $g\in\Gamma(Y,\tau^*\mcal{A})$.
Furthermore, given $y\in Y$ we have $\|g(y)\| =
|f(y)|\|a(\tau(y))\|$ so that $\supp g\subset \supp f$.  Thus
$g\in\Gamma_c(Y,\tau^*\mcal{A})$.  Once we make the
identification mentioned in Remark \ref{rem:8}, this shows
$f\otimes a \in \Gamma_c(Y,\tau^*\mcal{A})$.  

We would like to see that the set 
$C_c(Y)\odot A = \spn\{f\otimes a: f\in C_c(Y), a\in
A\}$ is dense in $\tau^*A = \Gamma_0(Y,\tau^*\mcal{A})$.  First
observe that if $g\in C_0(Y)$, $f\in C_c(Y)$, and $a\in A$ 
then $g\cdot (f\otimes a) = gf\otimes a$.  It follows that
$C_c(Y)\odot A$ is closed
under the $C_0(Y)$ action.  Now suppose $b\in\tau^*\mcal{A}$.  Choose
$a\in\mcal{A}$ so so that $a(\tau(y)) = b$ and $f\in C_c(Y)$ so that
$f(y) = 1$.  Then $f\otimes a(y) = b$.  We can now conclude from
Proposition \ref{prop:42} that $C_c(Y)\odot A$ is dense in $\tau^*A$. 
\end{proof}

Of course, we don't need to be working with pull backs for Proposition
\ref{prop:46} to hold. 

\begin{corr}
Suppose $A$ is a $C_0(X)$-algebra, and let $\mcal{A}$ be is its
associated  upper-semicontinuous bundle.  Given $f\in C_c(X)$ and
$a\in A$ define $f\otimes a(x)=  f(x) a(x)$ for all $x\in X$.
Then $f\otimes a\in \Gamma_c(X,\mcal{A})$ and
\[
C_c(X)\odot A := \spn\{f\otimes a : f\in C_c(X), a\in A\}.
\]
is dense in $A$. 
\end{corr}

\begin{proof}
This result follows immediately from Proposition \ref{prop:46} with
$\tau = \id$.  
\end{proof}

This is a good opportunity to introduce something that will be
fundamental to our study of crossed products.

\begin{definition}
\label{def:34}
\index{inductive limit topology}
Suppose $X$ is a locally compact Hausdorff space and $\mcal{A}$ is an
upper-semicontinuous Banach bundle over $X$.  Given a net $\{f_i\}_{i\in
  I}\subset \Gamma(X,\mcal{A})$ and $f\in \Gamma(X,\mcal{A})$ we say
that $f_i\rightarrow f$ with respect to the {\em inductive limit
  topology} if and only if $f_i\rightarrow f$ uniformly and there
exists a compact set $K$ in $X$ such that, eventually, all the $f_i$
and $f$ vanish off $K$.   Furthermore, we will say that a function
$F:\Gamma(X,\mcal{A})\rightarrow Y$ 
is continuous in the inductive limit topology if
$F(f_i)\rightarrow F(f)$ whenever $f_i\rightarrow f$ with respect to
the inductive limit topology. 
\end{definition}

\begin{remark}
First, we will often use Definition \ref{def:34} in the degenerate
situation where $X$ is locally compact Hausdorff and $\mcal{A}$ is the
trivial Banach bundle $X\times \C$.  In this case there is actually a
topology $\mcal{T}$ on $C_c(X)$ such that a function from $C_c(X)$ into a convex
space is continuous with respect to $\mcal{T}$ if and only if
it respects nets which converge in the inductive limit topology
\cite[Lemma D.10]{tfb}.  
However, we are not claiming in general that there is actually a topology on
$\Gamma(X,\mcal{A})$ which is characterized by these convergent
nets and even in the scalar case there may be nets which
converge in $C_c(X)$ with respect to $\mcal{T}$ and do not satisfy 
Definition \ref{def:34}.
\end{remark}

\begin{corr}
\label{cor:3}
Suppose $X$ and $Y$ are locally compact Hausdorff spaces, $A$ is a
$C_0(X)$-algebra, $\mcal{A}$ is its associated upper-semicontinuous
bundle, and $\tau:X\rightarrow Y$ is continuous.  Then $C_c(Y)\odot A$
is dense in $\Gamma_c(Y,\tau^*\mcal{A})$ with respect to the inductive
limit topology.  
\end{corr}

\begin{proof}
Suppose $g\in \Gamma_c(Y,\tau^*\mcal{A})$.  We know from Proposition
\ref{prop:46} that there exists a net $f_i\in C_c(Y)\odot A$
such that $f_i\rightarrow g$ uniformly. Let $K$ be a compact
neighborhood of $\supp g$ and choose $\phi\in C_c(Y)$ such that $\phi$
is one on $\supp g$ and $\phi$ is zero off $K$.  We showed in the proof
of Proposition \ref{prop:46} that $C_c(Y)\odot A$ is closed
under the $C_0(Y)$ action so that $\phi\cdot f_i\in C_c(Y)\odot
A$ for all $i$.  Furthermore, it follows immediately from the
fact that $\phi=1$ on $\supp g$ that we still have
$\phi\cdot f_i\rightarrow g$ uniformly.  Since clearly $\supp
\phi\cdot f_i \subset K$ we conclude that $\phi\cdot f_i\rightarrow g$
with respect to the inductive limit topology. 
\end{proof}


\section{Groupoid Dynamical Systems}

\label{sec:dynamical}

We are finally at a point where we can define what it means to be a
groupoid dynamical system. The reason we needed to introduce
$C_0(X)$-algebras in the last section is because, just like groupoid
actions on spaces, groupoids must act on fibred $C^*$-algebras. 

\begin{definition}
\label{def:33}
\index{dynamical system}
Suppose $G$ is a locally compact Hausdorff groupoid with a Haar
system.  
Let $A$ be a
$C_0(G\unit)$-algebra and $\mcal{A}$ its associated
upper-semicontinuous bundle.  An action $\alpha$ of $G$ on $A$ is a
family of functions $\{\alpha_\gamma\}_{\gamma\in G}$ such that, 
\begin{enumerate}
\item for each $\gamma\in G$ the map
  $\alpha_\gamma:A(s(\gamma))\rightarrow A(r(\gamma))$ is an
  isomorphism,  
\item $\alpha_{\gamma\eta} = \alpha_\gamma\circ\alpha_\eta$ for all
  $(\gamma,\eta)\in G^{(2)}$, and  
\item $\gamma\cdot a:= \alpha_\gamma(a)$ defines a (strongly) continuous action
  of $G$ on $\mcal{A}$.  
\end{enumerate}
The triple $(A,G,\alpha)$ is called a groupoid dynamical system.  We
say that $(A,G,\alpha)$ is separable if $A$ is separable and $G$ is
second countable.  
\end{definition}

\begin{remark}
The bundle map associated to an upper-semicontinuous bundle $\mcal{A}$
is assumed to be open and this is exactly the structure map for the
action in Definition \ref{def:33}.   Thus, as long as $G$ acts on
$\mcal{A}$ continuously, condition (c) above will be satisfied. 
\end{remark}

\begin{remark}
\index{Haar system}
The assumption that $G$ has a Haar system is not really an integral
part of Definition \ref{def:33}.  However, we don't care about
dynamical systems with no Haar system so it's useful to include it as
part of the definition.
\end{remark}

We will use the following frequently and will not bother to reference
it. 

\begin{prop}
Suppose $(A,G,\alpha)$ is a dynamical system.  Then 
\begin{enumerate}
\item $\alpha_u = \id_{A(u)}$ for all $u\in G\unit$, and 
\item $\alpha_{\gamma\inv} = \alpha_\gamma\inv$ for all $\gamma\in
  G$. 
\end{enumerate}
\end{prop}

\begin{proof}
Given $u\in G\unit$ we have $\alpha_{u} = \alpha_{u^2} = \alpha_u
\circ \alpha_u$.  Since $\alpha_u$ is an automorphism of $A(u)$, the
result holds.  However, we can now conclude from the fact that
$\gamma\inv\gamma \in G\unit$ that $\id = \alpha_{\gamma\inv\gamma} =
\alpha_{\gamma\inv}\circ \alpha_\gamma$ and the second half of the
proposition follows. 
\end{proof}

This next proposition gives an alternate characterization of groupoid
dynamical systems which is often easier to work with given that the
topology on $\mcal{A}$ can be poorly behaved.  However, there is no
way to completely dodge this fact and continuity is almost always the
hardest condition to verify.  

\begin{prop}
\label{prop:44}
Suppose $(A,G,\alpha)$ is a groupoid dynamical system.  Then
\[
\alpha(f)(\gamma) = \alpha_\gamma(f(\gamma))
\]
defines a $C_0(G)$-linear isomorphism of $s^*A$ onto $r^*A$.  

Conversely, if $G$ is a groupoid, $A$ is a $C_0(G\unit)$-algebra,
and there is a $C_0(G)$-linear isomorphism $\alpha:s^*A\rightarrow
r^*A$ then there are isomorphisms
$\alpha_\gamma:A(s(\gamma))\rightarrow A(r(\gamma))$ for all
$\gamma\in G$.  Furthermore, if $\alpha_{\gamma\eta} =
\alpha_\gamma\circ\alpha_\eta$ for all $(\gamma,\eta)\in G^{(2)}$ then
$(A,G,\alpha)$ is a dynamical system. 
\end{prop}

\begin{proof}
Given $f\in s^*A$ observe that $\alpha_\gamma(f(\gamma))\in
A(r(\gamma))$ for all $\gamma\in G$.  Therefore, if we define
$\alpha(f)$ as in the statement of the proposition, it follows that
$\alpha(f)$ is a section of $r^*A$.  If $\gamma_i\rightarrow
\gamma$ then $f(\gamma_i)\rightarrow f(\gamma)$.  It follows from
condition (c) of Definition \ref{def:33} that
$\alpha_{\gamma_i}(f(\gamma_i))\rightarrow \alpha_\gamma(f(\gamma))$.
Thus $\alpha(f)\in \Gamma(G,r^*\mcal{A})$.  Furthermore, each
$\alpha_\gamma$ is an isomorphism so that 
\begin{equation}
\label{eq:24}
\|\alpha(f)(\gamma)\| =
\|\alpha_\gamma(f(\gamma))\| = \|f(\gamma)\|.  
\end{equation}
It immediately follows that
$\alpha(f)$ vanishes at infinity if $f$ does.  Thus $\alpha(f)\in r^*A
= \Gamma_0(G,r^*\mcal{A})$.  
Now we will show that $\alpha$ is a $*$-isomorphism.  Given
$f,g\in s^*A$ and $\gamma\in G$ we have 
\[
\alpha(f+g)(\gamma) = \alpha_\gamma(f(\gamma)+g(\gamma)) 
= \alpha(f)(\gamma) + \alpha(g)(\gamma)
\]
where the second equality follows from the fact that $\alpha_\gamma$
is linear.  It is just as easy to show that $\alpha$ preserves the
rest of the operations.  Furthermore, it follows from \eqref{eq:24}
that $\|\alpha(f)\|_\infty = \|f\|_\infty$ and that $\alpha$ is
isometric.  Lastly, given $f\in r^*A$ define $g\in s^*A$ by $g(\gamma)
= \alpha_{\gamma\inv}(f(\gamma))$.  It is easy enough to see that $g$
is a continuous section which vanishes at infinity and that
$\alpha(g)=f$.  Thus $\alpha$ is a $*$-isomorphism.  Finally, we compute,
for $f\in s^*A$ and $\phi\in C_0(G)$,
\[
\alpha(\phi\cdot f)(\gamma) = \alpha_\gamma(\phi\cdot f(\gamma)) 
=\alpha_\gamma(\phi(\gamma)f(\gamma)) =
\phi(\gamma)\alpha_\gamma(f(\gamma)) = \phi\cdot \alpha(f)(\gamma).
\]

Now we prove the opposite direction.  Suppose that
$\alpha:s^*A\rightarrow r^*A$ is a
$C_0(G)$-isomorphism. It follows from Proposition \ref{prop:43} that for
each $\gamma\in G$ there is a $*$-isomorphism
$\alpha_\gamma:s^*A(\gamma) \rightarrow r^*A(\gamma)$ such that
$\alpha_\gamma(f(\gamma)) = \alpha(f)(\gamma)$ for all $f\in s^*A$.
It then follows from Proposition \ref{prop:45} that we can make the 
identification $s^*A(\gamma) =
A(s(\gamma))$ and $r^*A(\gamma) = A(r(\gamma))$.  Thus we
have satisfied condition (a) of Definition \ref{def:33}.  Condition
(b) is satisfied by assumption.  Suppose $\gamma_i\rightarrow
\gamma$ in $G$ and $a_i\rightarrow a$ in $\mcal{A}$ such that $p(a_i)
= s(\gamma_i)$ for all $i$ and $p(a) = s(\gamma)$.  Choose $g\in s^*A$
such that $g(\gamma) = a$.  Then $\alpha(g)\in r^*A$ and we have,
using the fact that both $a_i$ and $g(\gamma_i)$ converge to $a$, 
\[
\|\alpha_{\gamma_i}(a_i) - \alpha(g)(\gamma_i)\| = 
\|\alpha_{\gamma_i}(a_i-g(\gamma_i))\| = 
\|a_i-g(\gamma_i)\|\rightarrow 0.
\]
Since $p(\alpha_{\gamma_i}(a_i)) = r(\gamma_i)\rightarrow r(\gamma) =
p(\alpha_\gamma(a))$ it follows from Proposition \ref{prop:37} that
$\gamma_i\cdot a_i\rightarrow \gamma\cdot a$ and that $(A,G,\alpha)$
is a groupoid dynamical system.  
\end{proof}

\begin{remark}
\label{rem:9}
Given a dynamical system we will construct a $*$-algebra structure on $\Gamma_c(G,r^*,\mcal{A})$.  This
algebra will eventually be completed into the groupoid crossed
product.  However, in order to define the convolution operation we will need to
use vector valued integration.  Most of the time vector valued
integrals ``just work'' and can be treated like scalar valued
integrals.  Those readers looking for a good reference are referred to
\cite[Appendix B]{tfb2}.  The short version is that given a Radon
measure $\mu$ on a locally compact Hausdorff space $X$ and a
separable Banach
algebra $B$ we define $\mcal{L}^1(X,B)$ to be the set of measurable
functions $f:X\rightarrow B$ such that\footnote{Since $B$ is separable,
  this is just the usual definition of measurability.}
\[
\|f\|_1 := \int_X\|f(x)\|d\mu(x) < \infty.
\]
Then there is a linear function
\[
f\mapsto \int_X f(x) d\mu(x)
\]
from $\mcal{L}^1(X,B)$ into $B$ satisfying
\[
\left\|\int_X f(x) d\mu(x)\right\| \leq \int_X \|f(x)\| d\mu(x).
\] \index{elementary tensors}
Given $b\in B$ and $f\in C_c(X)$ we can define the {\em elementary
  tensor} $f\otimes b\in C_c(X,B)$ by $(f\otimes b)(x) = f(x)b$ for all $x$.  We
then have 
\begin{equation}
\label{eq:25}
\int_X f\otimes b(x)d\mu(x) = \int_X f(x)d\mu(x) b.
\end{equation}
Furthermore, if $L:B\rightarrow B'$ is a bounded linear map onto
another Banach space $B'$ then
\begin{equation}
\label{eq:37}
L\left( \int_X f(x)d\mu(x)\right) = \int_X L(f(x))d\mu(x)
\end{equation}
for all $f\in \mcal{L}^1(X,B)$.  In fact \eqref{eq:37} characterizes the
integral when you allow $L$ to be any bounded linear functional.  
Lastly, there is a Fubini's Theorem for
vector valued integrals that is analogous to the scalar one.  In
particular, when integrating $L^1$ functions we will
freely reorder the integrals.  
\end{remark}

\begin{remark}
\label{rem:10}
Suppose $G$ is a locally compact Hausdorff groupoid with Haar system
$\lambda$ and $A$ is a $C_0(G^{(0)})$-algebra with associated bundle
$\mcal{A}$.  Given $f\in \Gamma_c(G,r^*\mcal{A})$ it is clear that
$f(\gamma)\in A(u)$ for all $\gamma\in G^u$.
Therefore, by Remark \ref{rem:9} we can form the integral
$\int_{G^u}f|_{G^u}d\lambda^u$.  Since the support of $\lambda^u$ is equal to
$G^u$ there is no harm denoting this integral by $\int_G f\lambda^u$.
\end{remark}

\begin{prop}
\label{prop:47}
Suppose $(A,G,\alpha)$ is a dynamical system.  
Given $f\in \Gamma_c(G,r^*\mcal{A})$ the function 
\begin{equation}
\label{eq:26}
u \mapsto \int_G f(\gamma)d\lambda^u(\gamma)
\end{equation}
is continuous.  Furthermore if $g\otimes a\in C_c(G)\odot A$ then 
\[
\int_G (g\otimes a)(\gamma) d\lambda^u(\gamma) = 
\int_G g(\gamma)d\lambda^u(\gamma) a(u)
\]
for all $u\in G\unit$.  
\end{prop}

\begin{proof}
It follows from Remark \ref{rem:10} that \eqref{eq:26} is well
defined.  We will address the second half of the proposition first.  Let
$g\otimes a$ be an elementary tensor with $g\in C_c(G)$ and $a\in A$.
When we restrict $g\otimes a$ to $G^u$ we get a new elementary tensor
$g|_{G^u}\otimes a(u)\in C_c(G^u,A(u))$ in the sense of Remark
\ref{rem:9}.  However, the result then follows from \eqref{eq:25}.  

Now suppose $f\in \Gamma_c(G,r^*\mcal{A})$, that $u_i\rightarrow u$ in
$G\unit$, and fix $\epsilon > 0$.  
Use Corollary \ref{cor:3}
to find a collection of elementary tensors $\{g_j^i\otimes a_j^i\}$ such
that the net $k_j = \sum_i g_j^i \otimes a_j^i$ converges to $f$
with respect to the inductive limit topology.  Let $K$ be a compact
set which eventually contains the supports of the $k_j$ and
$f$.  Next, since $K$ is compact and the $\lambda^u$ vary continuously, we
can find an upper bound $M$ for $\lambda^u(K)$. Now choose $J$ so that 
$\supp k_J\subset K$ and $\|f-k_J\|_\infty <
\epsilon/M$.  Then for all $v\in G\unit$ we have 
\begin{align*}
\left\| \int_G fd\lambda^v - \int_G k_Jd\lambda^v\right\|
&= \left\| \int_G f-k_J d\lambda^v\right\| \\
&\leq \int_G \|f-k_J\| d\lambda^v \\
&\leq \|f-k_J\|_\infty \lambda^v(K) < \epsilon.
\end{align*}
However, it is clear enough that 
\[
v\mapsto  \int_G g_j^i\otimes a_j^i(\gamma) d\lambda^v(\gamma) = \int_G g_j^i(\gamma)d\lambda^v(\gamma)a_j^i(v)
\]
is continuous for all $j$ and $i$.  Since sums of continuous functions
are continuous we conclude that 
\[
v\mapsto \int_G k_J(\gamma)d\lambda^v(\gamma)
\]
is continuous as well.  It now follows from the previous paragraph, 
and using the last part of Proposition \ref{prop:35}, that $\int_G
fd\lambda^{u_i}\rightarrow \int_G fd\lambda^u$.  
\end{proof}

We are now ready to turn $\Gamma_c(G,r^*\mcal{A})$ into a
$*$-algebra.  This material is all worked out, in greater generality, in
\cite[Section 4]{renaultequiv} and many of these proofs are copied from
there.  

\begin{prop}
\label{prop:48}
\index{inductive limit topology}
Let $G$ be a locally compact Hausdorff groupoid with Haar
system $\{\lambda^u\}$, $A$ a $C_0(G\unit)$-algebra with associated
bundle $\mcal{A}$, and $\alpha$ an action of $G$ on $A$.  Then 
$\Gamma_c(G,r^*\mcal{A})$ becomes a $*$-algebra with respect to the
operations
\[
f*g(\gamma) = \int_G f(\eta)\alpha_\eta(g(\eta\inv\gamma))
d\lambda^{r(\gamma)}(\eta)\quad\text{and}\quad
f^*(\gamma)= \alpha_\gamma(f(\gamma\inv)^*).
\]
Furthermore, these operations are continuous with respect to the
inductive limit topology.
\end{prop}

We are going to need the following lemma, which is quite similar in
nature to Proposition \ref{prop:47}.  It can also be proved in the
same fashion, but we have presented a different way of approaching the
situation.  

\begin{lemma}
\label{lem:9}
Let $G* G = \{(\gamma,\eta)\in G\times G : r(\gamma)=r(\eta)\}$ and
let $r^*\mcal{A}$ be the pull back of $\mcal{A}$ to $G* G$.  Given
$F\in \Gamma_c(G*G,r^*\mcal{A})$ then 
\[
f(\gamma) = \int_G F(\eta,\gamma) d\lambda^{r(\gamma)}(\eta)
\]
defines a section in $\Gamma_c(G,r^*\mcal{A})$.  
\end{lemma}

\begin{proof}
First, observe that $f$ is clearly a section. 
Now, it is straightforward to check that if $F_i\rightarrow F$ with respect
to the inductive limit topology in $\Gamma_c(G* G,r^*\mcal{A})$ then
$f_i\rightarrow f$ with respect to the inductive limit topology.
Furthermore, observe that in this case if we show each $f_i$ is continuous and
compactly supported then $f$ must be as well.  It
follows then, that it suffices to show that the lemma holds for
elementary tensors $g\otimes a$ where $g\in C_c(G* G)$ and $a\in
A$.  We can use Lemma \ref{lem:8} 
to extend $g$ to all of $G\times G$.  Now, sums of functions
of the form $(\gamma,\eta)\mapsto h_1(\gamma)h_2(\eta)$ are dense in
$C_c(G\times G)$ with respect to the inductive limit topology
\cite[Corollary B.17]{tfb}, and it is
not hard to see that if $f_i\rightarrow f$ with respect to the
inductive limit topology then $f_i\otimes a\rightarrow f\otimes a$
with respect to the inductive limit topology.  Thus, using the above
argument, we can assume without loss of generality that
$g(\gamma,\eta) = h_1(\gamma) h_2(\eta)$ where $h_1,h_2\in C_c(G)$.
However, we now have
\[
f(\gamma) = \int_G g(\eta,\gamma) a(r(\gamma))
d\lambda^{r(\gamma)}(\eta) = 
h_2(\gamma)a(r(\gamma)) \int_G h_1(\eta) d\lambda^{r(\gamma)}(\eta).
\]
It is clear that in this case $f$ is a continuous compactly supported
section, so we are done.  
\end{proof}

\begin{proof}[Proof of Proposition \ref{prop:48}]
First, we have to check that the operations are well defined.  This is
straightforward for $f^*$ since everything in sight is continuous.
The fact that convolution produces a continuous section is exactly
Lemma \ref{lem:9} once you realize that $(\eta,\gamma)\mapsto
f(\eta)\alpha_\eta(\eta\inv\gamma)$ is a section in $\Gamma_c(G*
G,r^*\mcal{A})$.  

At this point we have to show that the operations are well behaved
algebraically.  For the most part these computations are omitted.
However, we will verify two of them as examples of how the rest should
work.  First, we will show that the convolution is associative.
Suppose $f,g,h\in \Gamma_c(G,r^*\mcal{A})$ and $\gamma\in G$, then 
\begin{align*}
(f*g)*h(\gamma) &= \int_G
f*g(\eta)\alpha_\eta(h(\eta\inv\gamma))d\lambda^{r(\gamma)}(\eta) \\
&= \int_G\int_G
f(\zeta)\alpha_\zeta(g(\zeta\inv\eta))\alpha_\eta(h(\eta\inv\gamma))
d\lambda^{r(\eta)}(\zeta)d\lambda^{r(\gamma)}(\eta) \\
&= \int_G\int_G
f(\zeta)\alpha_\zeta(g(\eta))\alpha_{\zeta\eta}(h(\eta\inv\zeta\inv\gamma))
d\lambda^{s(\zeta)}(\eta)d\lambda^{r(\gamma)}(\zeta)
\end{align*}
where we switched the order of the integrals and used the left
invariance of the Haar measure to get the last equality.  Continuing
the computation by using the fact that $\alpha$ respects the groupoid
operations we get
\begin{align*}
(f*g)*h(\gamma)&= \int_G\int_G
f(\zeta)\alpha_\zeta(g(\eta)\alpha_\eta(h(\eta\inv\zeta\inv\gamma)))
d\lambda^{s(\zeta)}(\eta)d\lambda^{r(\gamma)}(\zeta) \\
&= \int_G f(\zeta)\alpha_\zeta(g*h(\zeta\inv\gamma))
d\lambda^{r(\gamma)}(\zeta) \\
&= f*(g*h)(\gamma).
\end{align*}
Notice that we used the fact that vector valued integrals are
preserved by bounded linear maps to pass the integral through $\alpha$.

We will also show that $(f*g)^* = g^**f^*$.  Suppose $f,g\in
\Gamma_c(G,r^*\mcal{A})$ and $\gamma\in G$.  We have 
\begin{align*}
g^**f^*(\gamma) &= \int_G g^*(\eta)\alpha_\eta(f^*(\eta\inv\gamma))
d\lambda^{r(\gamma)}(\eta) \\
&= \int_G
\alpha_\eta(g(\eta\inv)^*)\alpha_\eta(\alpha_{\eta\inv\gamma}(f(\gamma\inv\eta)^*))
d\lambda^{r(\gamma)}(\eta) \\
&= \int_G \alpha_\eta(g(\eta\inv))^*
\alpha_\gamma(f(\gamma\inv\eta))^*d\lambda^{r(\gamma)}(\eta) \\
&= \alpha_\gamma\left(\int_G
  f(\gamma\inv\eta)\alpha_{\gamma\inv\eta}(g(\eta\inv))
  d\lambda^{r(\gamma)}(\eta)\right)^*
\end{align*}
To get the last equality we pulled both the $\alpha_\gamma$ and the
$^*$ operation out of the integral using the fact that they are bounded
linear maps.\footnote{Technically the $^*$ operation is conjugate
  linear but it is linear when viewed as a map into the conjugate
  algebra.}  Now, using left invariance, we get
\begin{align*}
g^**f^*(\gamma) &= \alpha_\gamma\left(\int_G
  f(\eta)\alpha_{\eta}(g(\eta\inv\gamma\inv))d\lambda^{s(\gamma)}(\eta)\right)^*
\\
&= \alpha_\gamma(f*g(\gamma\inv)^*) = (f*g)^*(\gamma).
\end{align*}
The rest of the algebraic computations are similar and it is apparent
now why we would want to skip them.  

The only thing that remains to be verified is that the operations are
continuous with respect to the inductive limit topology.  This is
clearly true for the involution since each $\alpha_\gamma$ is
isometric.  The convolution is only slightly more complicated.  Suppose
$f_i\rightarrow f$ and $g_i\rightarrow g$ with respect to the
inductive limit topology in $\Gamma_c(G,r^*\mcal{A})$.  Let
$F_i(\eta,\gamma) = f_i(\eta)\alpha_\eta(g_i(\eta\inv\gamma))$ for
each $i$ and let $F(\eta,\gamma) =
f(\eta)\alpha_\eta(g(\eta\inv\gamma))$.  Then it is easy to show that
$F_i\rightarrow F$ with respect to the inductive limit topology in
$\Gamma_c(G* G,r^*\mcal{A})$.  For instance, if the support of $f_i$
is eventually contained in the compact set $K$ and the support of
$g_i$ is eventually contained in the compact set $L$ then the support
of $F_i$ will eventually be contained in $K\times KL$, which is
compact.  But then $f_i*g_i$ and $f*g$ are defined via integration as in Lemma
\ref{lem:9}.  It is now straightforward to show that
$f_i*g_i\rightarrow f*g$ with respect to the inductive limit
topology.
\end{proof}

If $f\in \Gamma_c(G,r^*\mcal{A})$ then $\gamma\mapsto \|f(\gamma)\|$
is upper-semicontinuous and compactly supported.  Therefore this
function is integrable on $G$ with respect to any Radon measure.  This
allows us to define the following norm. 

\begin{definition}
\label{def:35}
\index[not]{$\|\cdot\|_I$}
\index{I-norm@$I$-norm}
If $(G,A,\alpha)$ is a dynamical system we define the {\em $I$-norm}
of $f\in\Gamma_c(G,r^*\mcal{A})$ to be 
\[
\|f\|_I = \max\left\{\sup_{u\in G\unit} \int_G
  \|f(\gamma)\|d\lambda^u(\gamma), 
\sup_{u\in G\unit} \int_G \|f(\gamma)\|d\lambda_u(\gamma) \right\}.
\]
Recall that we define $\lambda_u := (\lambda^u)\inv$.  
\end{definition}

This norm structure interacts nicely with the existing structure on
$\Gamma_c(G,r^*\mcal{A})$.  Actually, the $I$-norm was defined to play
along, so to speak.  For example, we have to use both
supremums in the definition of $\|\cdot\|_I$ in order to make the
involution an isometry.  

\begin{prop}
\label{prop:60}
\index{inductive limit topology}
Suppose $(A,G,\alpha)$ is a dynamical system.  Then the $I$-norm is a
norm on $\Gamma_c(G,r^*\mcal{A})$.  Furthermore, for $f,g\in
\Gamma_c(G,r^*\mcal{A})$ we have $\|f*g\|_I \leq \|f\|_I \|g\|_I$ and
$\|f^*\|_I = \|f\|_I$.  
Finally, if $f_i\rightarrow f$ in $\Gamma_c(G,r^*\mcal{A})$ with
respect to the inductive limit topology then $f_i\rightarrow f$ with
respect to the $I$-norm.
\end{prop}

\begin{proof}
First we will show that $\|f\|_I<\infty$ for all
$f\in\Gamma_c(G,r^*\mcal{A})$.  Since $\supp f$ is compact, and since
the $\lambda^u$ vary continuously, we can find an upper bound $M$ for
the set $\{\lambda^u(\supp f)\}$.  We can also increase $M$,
if necessary, so that it is also an upper bound for $\{\lambda_u(\supp
f)\}$.  However, it is now clear that $\|f\|_I \leq M\|f\|_\infty$.  
Showing that $\|\cdot\|_I$ is a norm is straightforward.  We will
restrict ourselves to showing that it is positive definite.  Suppose
$f\in \Gamma_c(G,r^*\mcal{A})$ and $\|f\|_I = 0$.  Then, in
particular, $\int_G\|f(\gamma)\|d\lambda^u(\gamma) = 0$ for all $u\in
G\unit$.  Suppose
$f(\gamma)\ne 0_{r(\gamma)}$ for some $\gamma\in G$.  When restricted
to $G^{r(\gamma)}$ the function $\eta\mapsto \|f(\eta)\|$ is
positive, continuous, and nonzero at $\gamma$.  
Since $\supp \lambda^{r(\gamma)} = G^{r(\gamma)}$ this
implies that $\int_G\|f(\gamma)\|d\lambda^u(\gamma) \ne 0$, which is a
contradiction.  

Now, suppose $f,g\in \Gamma_c(G,r^*\mcal{A})$.  Then 
\begin{align*}
\int_G \|f*g(\gamma)\|d\lambda^u(\gamma) &\leq 
\int_G\int_G \|f(\eta)\|\|\alpha_\eta(g(\eta\inv\gamma))\|
d\lambda^u(\eta)d\lambda^u(\gamma) \\
&= \int_G \|f(\eta)\| \int_G \|g(\eta\inv\gamma)\|d\lambda^u(\gamma)
d\lambda^u(\eta) \\
&= \int_G \|f(\eta)\| \int_G \|g(\gamma)\|d\lambda^{s(\eta)}(\gamma)
d\lambda^u(\eta) \\
&\leq \|g\|_I \int_G \|f(\eta)\| d\lambda^u(\eta) \leq \|g\|_I \|f\|_I
\end{align*}
Similar considerations show that $\int_G
\|f*g(\gamma)\|d\lambda_u(\gamma) \leq \|g\|_I\|f\|_I$.  It follows
that the $I$-norm is submultiplicative.  Next, we compute
\begin{align*}
\int_G \|f^*(\gamma)\| d\lambda^u(\gamma) &=
\int_G \|\alpha_\gamma(f(\gamma\inv)^*)\| d\lambda^u(\gamma) \\
&=  \int_G \|f(\gamma\inv)\|d\lambda^u(\gamma) = \int_G\|f(\gamma)\|
d\lambda_u(\gamma).
\end{align*}
Therefore, it is clear from the definition of the $I$-norm that
$\|f^*\|_I = \|f\|_I$.  

Lastly suppose that $f_i\rightarrow f$ with respect to the inductive
limit topology in $\Gamma_c(G,r^*\mcal{A})$.  Let $K$ be a compact set
which eventually contains $\supp f_i$ and $\supp f$.
Furthermore, let $M$ be an upper bound for both $\{\lambda^u(K)\}$ and
$\{\lambda_u(K)\}$.  Fix $\epsilon > 0$ and
observe that eventually 
$\|f-f_i\|_\infty<\epsilon /M$.  It follows that, eventually, 
\begin{align*}
\int_G \|f(\gamma)-f_i(\gamma)\|d\lambda^u(\gamma) &\leq
\|f-f_i\|_\infty \lambda^u(K) < \epsilon\quad\text{and,} \\
\int_G \|f(\gamma)-f_i(\gamma)\|d\lambda_u(\gamma) &\leq
\|f-f_i\|_\infty \lambda_u(K) < \epsilon. \\
\end{align*}
for all $u\in G\unit$.  Thus $\|f_i-f\|_I\rightarrow 0$ and we
are done. 
\end{proof}

\begin{remark}
While the $I$-norm does make $\Gamma_c(G,r^*\mcal{A})$ into a
$*$-algebra it does not, however, satisfy the $C^*$-identity.  This
means that the completion of $\Gamma_c(G,r^*\mcal{A})$ is not a
$C^*$-algebra.  We will instead use the ``universal norm'' to
construct the crossed product in Section \ref{sec:crossedprod}.
\end{remark}

\begin{corr}
\label{cor:6}
Suppose $(A,G,\alpha)$ is a dynamical system.  Then $C_c(G)\odot A$ 
is dense in $\Gamma_c(G,r^*\mcal{A})$ with respect to the
$I$-norm. 
\end{corr}

\begin{proof}
This follows immediately from Corollary \ref{cor:3} and the last
statement in Proposition \ref{prop:60}.
\end{proof}


\section{Covariant Representations}
\label{sec:covariant}
Our goal is to define the notion of a covariant representation of a groupoid
dynamical system because these are the representations we will use to
define the universal norm.  However, in order to do that we first have to
discuss the notion of a groupoid representation. 

\subsection{Groupoid Representations}

Because groupoids are fibred objects they must be represented on
fibred objects.  The appropriate bundle in this case is a Borel
Hilbert bundle.  These objects are relatively classical and only those
proofs that seem relevant will be included.  All of the necessary technology
can be found in \cite[Appendix F.1]{tfb2}.  In particular, the
following definition and remarks are lifted straight from there.
Another good reference for this material is \cite[Section 7]{renaultequiv}.  

\begin{remark}
The reader may also wish to consider \cite[Chapter 3]{invitation} where you will
learn about analytic and standard spaces.  
While it is important to understand the the difference
between an analytic space and a standard space at some point,
beginners are advised to just ignore these Borel considerations on
their first pass through the material.
\end{remark}

\begin{definition}
\label{def:36}
\index{Borel Hilbert bundle}
\index[not]{$X*\mfrk{H}$}
Suppose $\mfrk{H} = \{\mcal{H}(x)\}_{x\in X}$ is a collection of
separable (non-zero) complex Hilbert spaces indexed by an analytic
Borel space $X$.  We define the total space to be the disjoint union 
\[
X*\mfrk{H} := \{(x,h) : h\in\mcal{H}(x)\}
\]
and let $\pi:X*\mfrk{H}\rightarrow X$ be the obvious projection map.
Then $X*\mfrk{H}$ is an analytic (resp. standard) 
{\em Borel Hilbert bundle} if
$X*\mfrk{H}$ has an analytic (resp. standard) Borel structure such that 
\begin{enumerate}
\item $\pi$ is a Borel map and 
\item there is a sequence $\{f_n\}$ of sections such that 
\begin{enumerate}
\item the maps $\bar{f}_n : X*\mfrk{H} \rightarrow \C$ defined by
\[
\bar{f}_n(x,h) := (f_n(x),h), 
\]
are Borel for each $n$, 
\item for each $n$ and $m$, 
\[ 
x\mapsto (f_n(x),f_m(x))
\]
is Borel, and 
\item the functions $\{\bar{f}_n\}$ and $\pi$ separate points of
  $X*\mfrk{H}$.  
\end{enumerate}
\end{enumerate}
The sequence $\{f_n\}$ is called a {\em fundamental sequence} for 
$X*\mfrk{H}$.  We let $B(X*\mfrk{H})$ be the set of Borel sections
of $X*\mfrk{H}$.
\end{definition}

\begin{remark}
We are using the same notational trickery that is described in Remark
\ref{rem:8}.  In particular, a section $f$ of $X*\mfrk{H}$ is of the
form $f(x) = (x,\tilde{f}(x))$ where $\tilde{f}$ maps $X$ into the
disjoint union of the $\mcal{H}(x)$ and $\tilde{f}(x)\in \mcal{H}(x)$
for all $x\in X$.  Of course, $f$ is completely determined by
$\tilde{f}$ and just as in Remark \ref{rem:8} we will not distinguish
between the two functions.  
\end{remark}

We collect some useful facts from \cite[Appendix F]{tfb2} into the
next proposition. However, we must first prove a useful lemma, which is an
immediate consequence of the Unique Structure Theorem \cite[Theorem
3.3.5]{invitation} and modeled off \cite[Lemma D.20]{tfb2}.  

\begin{lemma}
\label{lem:11}
Suppose that $(X,\mcal{B})$ is an analytic Borel space and that
$f_n:X\rightarrow Y_n$ is a sequence of Borel functions on $X$ which map into
countably generated Borel spaces  $Y_n$, and 
which separate points.  Then $\mcal{B}$ is the smallest
$\sigma$-algebra in $X$ such that each $f_n$ is Borel.  In particular,
$g:Y\rightarrow X$ is Borel if and only if $f_n\circ g$ is Borel for
all $n$.  
\end{lemma}
\begin{proof}
Let $\mcal{B}_0$ be the smallest $\sigma$-algebra such that each $f_n$
is Borel.  Let $U^n_k$ be a countable generating set for the Borel
structure on $Y_n$.  Then $\{f_n\inv(U_k^n)\}$ is a countable family
that generates $\mcal{B}_0$ and which separates points.  Thus
$\mcal{B}=\mcal{B}_0$ by the Unique Structure Theorem.  The rest of
the proposition is straightforward.  
\end{proof}

The following is little more than a Swiss army knife for dealing with
Borel Hilbert bundles, although the fact that all analytic Borel
Hilbert bundles are, in some sense, trivial is interesting in its own
right.  

\begin{prop}
\label{prop:49}
Let $X*\mfrk{H}$ be an analytic Borel Hilbert bundle with fundamental sequence
$f_n$.  Then the following are true. 
\begin{enumerate}
\item We have $f\in B(X*\mfrk{H})$ if and only if $x\mapsto (f(x),f_n(x))$ is
  Borel for all $n$. 
\item If $f,g\in B(X*\mfrk{H})$ then $x\mapsto (f(x),g(x))$ is
  Borel. It follows that $x\mapsto \|f(x)\|$ is Borel.  
\item There exists a fundamental sequence $\{e_k\}$ in $B(X*\mfrk{H})$
  such that 
\begin{enumerate}
\item for each $x\in X$ the set $\{e_k(x)\}_k$, minus any possible zero
vectors, is an orthonormal basis for $\mcal{H}(x)$ and 
\item for each $k$ there is a Borel partition $X = \bigcup_i B_i^k$ and
  for each $(i,k)$ finitely many Borel functions $\phi_j^{i,k}$ where
  $1\leq j \leq l(i,k)$ such that 
\[
e_k(x) = \sum_{j=1}^{l(i,k)} \phi_j^{i,k}(x)f_j(x)
\]
for all $x\in B_i^k$. 
\end{enumerate}
\index{special orthogonal fundamental sequence}
Such a sequence is called a {\em special orthogonal fundamental sequence}. 

\item There exists a Borel partition $X=X_\infty\cup X_1\cup
  X_2\cup\ldots$ of $X$ such that, if $\mcal{H}_d$ is a fixed Hilbert
  space of dimension $1\leq d\leq \aleph_0$, then $X*\mfrk{H}$ is
  Borel isomorphic to the disjoint union $\coprod_{d=1}^{d=\infty}
  X_d\times \mcal{H}_d$.  
\end{enumerate}
\end{prop}
\begin{proof}[Remark]
These are all results in \cite{tfb2} and we will just reference their
locations here.  Part (a) is demonstrated in \cite[Remark F.3]{tfb2}.
In particular it follows from Lemma \ref{lem:11} and the fact
that the $\bar{f}_n$ and $\pi$ separate points.  (Notice that $\pi\circ f =
\id$ is always Borel and that $X$ is countably generated since it's
analytic.) 
Parts (b) and (c) are proved in
\cite[Proposition F.6]{tfb2} and Part (d) is \cite[Corollary
F.9]{tfb2}.  
\end{proof}

It is worth pointing out that Borel Hilbert bundles do not usually
come equipped with an existing Borel structure.  Usually we give them
one using the following

\begin{prop}[{\cite[Proposition F.8]{tfb2}}]
\label{prop:113}
Suppose that $X$ is an analytic Borel space and that $\mfrk{H} =
\{\mcal{H}(x)\}_{x\in X}$ is a family of separable Hilbert spaces.
Suppose that $\{f_n\}$ is a countable family of sections of
$X*\mfrk{H}$ such that conditions (ii) and (iii) of axiom (b) in
Definition \ref{def:36} are satisfied.  Then there is a unique analytic
Borel structure on $X*\mfrk{H}$ such that $X*\mfrk{H}$ becomes an
analytic Borel Hilbert bundle and $\{f_n\}$ is a fundamental sequence.
\end{prop}

\begin{example}
\label{ex:17}
Suppose $G$ is a locally compact Hausdorff groupoid with Haar system
$\lambda$.  Let $L^2(\lambda):=\{L^2(G^u,\lambda^u)\}_{u\in G\unit}$ 
and form the bundle $G\unit*L^2(\lambda)$. 
The Borel structure on $G\unit* L^2(\lambda)$ is determined by a
sequence of sections $\{\xi_n\}_{n=1}^\infty$ defined as follows.
Choose a point separating sequence of functions $\{f_n\}_{n=1}^\infty$
in $C_c(G)$ and define $\xi_n:G\unit\rightarrow G\unit*L^2(\lambda)$
by the formula 
\[
\xi_n(u)(\gamma) = f_n(\gamma)\quad\text{for $\gamma\in G^u$.}
\]
It follows easily from Proposition \ref{prop:113} that there is a unique
Borel structure on $G\unit*L^2(\lambda)$ such that
$G\unit*L^2(\lambda)$ is an analytic Borel Hilbert bundle and $\xi_n$
is a fundamental sequence. 
\end{example}

Since Borel Hilbert bundles are 
fibred objects, they give rise to the following isomorphism groupoid.
This groupoid will eventually form the range of a unitary groupoid
representation. 

\begin{definition}
\label{def:38}
\index{groupoid}
\index{isomorphism groupoid}
\index[not]{$\isom(X*\mfrk{H},\mu)$}
If $X*\mfrk{H}$ is an analytic Borel Hilbert bundle then its {\em
  isomorphism groupoid} is defined to be 
\[
\isom(X*\mfrk{H}):= \{(x,V,y):V\in U(\mcal{H}(y),\mcal{H}(x))\}
\]
equipped with the weakest Borel structure such that $(x,V,y)\mapsto
(Vf(y),g(x))$ is Borel for all $f,g\in B(X*\mfrk{H})$.  We define 
the set of composable pairs to be
\[
(\isom(X*\mfrk{H}))^{(2)} :=\{((x,V,y),(w,U,z))\in\isom(X*\mfrk{H})\times\isom(X*\mfrk{H}) : y = w\}
\]
and the operations to be 
\begin{align*}
(x,V,y)(y,U,z) &:= (x,VU,z), &
(x,V,y)\inv &:= (y,V^*,z).
\end{align*}
\end{definition}

Of course, we made a number of claims in Definition \ref{def:38} 
which need to be verified.  

\begin{prop}
\label{prop:52}
Suppose $X*\mfrk{H}$ is an analytic Borel Hilbert bundle.  Then
$\isom(X*\mfrk{H})$ is a Borel groupoid.  Furthermore,
$\isom(X*\mfrk{H})$ is an analytic Borel space and is standard if
$X*\mfrk{H}$ is.  Finally, the unit space of $\isom(X*\mfrk{H})$ can be
identified with $X$ and under this identification the range and source
maps are given by $r(x,V,y)=x$ and $s(x,V,y)=y$.  
\end{prop}

\begin{proof}
It is clear from the definition of the groupoid operations that
$\isom(X*\mfrk{H})$ is a groupoid, that $(x,\id,x)\mapsto x$ is an
identification of the unit space with $X$, and that the range and
source maps have the appropriate form.  All that remains is to
demonstrate the statements concerning the Borel structure.  

Recall from Proposition \ref{prop:49} that there is a Borel
partition $\{X_d\}_{d=0}^{d=\infty}$ of $X$ such that $X*\mfrk{H}$ is
Borel 
isomorphic to the disjoint union 
$C = \coprod_{d=0}^{d=\infty} X_d\times \mcal{H}_d$ where $\mcal{H}_d$ is a
Hilbert space of dimension $d$.  It suffices to see that $\isom(C)$
has the required properties.  However, it is easy to check that
$\isom(C)$ is Borel isomorphic to the disjoint union
$\coprod_d \isom(X_d\times \mcal{H}_d)$.  It is also easy to check that
the trivial bundle $\isom(X_d\times \mcal{H}_d)$ is Borel
isomorphic to the space $X_d\times U(\mcal{H}_d)\times X_d$ where
$U(\mcal{H}_d)$ is given the (standard) Borel structure coming from the weak
operator topology.  Thus,
$\isom(X*\mfrk{H})$ is Borel isomorphic to the disjoint union
\[
\coprod_{d=0}^{d=\infty} X_d\times U(\mcal{H}_d)\times X_d.
\]
Now, $U(\mcal{H}_d)$ has a standard Borel structure.  
It takes an application of
\cite[Theorem 3.34]{invitation}, but is otherwise straightforward, to
show that a Borel subset of a standard space is standard and a Borel
subset of an analytic space is analytic.  The upshot is that the $X_d$
are at least analytic and are standard if $X$ is.  Thus the product
$X_d\times U(\mcal{H}_d)\times X_d$ is at least analytic and is
standard if $X$ is.  As a result $\isom(X*\mfrk{H})$ is analytic and
is standard if $X$ is.  Finally, it is easy to see that the groupoid
operations are Borel on $X_d\times U(\mcal{H}_d)\times X_d$ and it
follows quickly that they are Borel on the disjoint union.  
\end{proof}

The following lemma is useful because it allows us to use a
fundamental sequence to check when a map into $\isom(X*\mfrk{H})$ is
Borel, as opposed to using every section in $B(X*\mfrk{H})$.  It will
be particularly useful for groupoid representations because condition
(a) will turn out to be trivial.  

\begin{lemma}
\label{lem:51}
Suppose $X*\mfrk{H}$ is an analytic Borel Hilbert bundle with
fundamental sequence $\{f_n\}$, that $Z$ is
some Borel space, and that 
$U:Z\rightarrow
\isom(X*\mfrk{H})$.  Then $U$ is Borel if only if 
\begin{enumerate}
\item $r\circ U$ is Borel and, 
\item $\psi_{n,m}\circ U$ is Borel for all $n,m$ where
  $\psi_{n,m}(x,V,y) = (Vf_n(y),f_m(x))$.
\end{enumerate}
\end{lemma}
\begin{proof}
First observe that if $r\circ U$ is Borel then $s\circ U$ is Borel
since $s$ is just the composition of $r$ with inversion.  Using Lemma
\ref{lem:11}, and noting that $X$ is countably generated since it is
analytic, it will suffice to show that $\{\psi_{n,m}\}$, $r$, and $s$
separate points.  Suppose $(x,V,y),(w,U,z)\in\isom(X*\mfrk{H})$ such
that $r(x,V,y)=r(w,U,z)$,\ $s(x,V,y)=s(w,U,z)$ and
$\psi_{n,m}(x,V,y)=\psi_{n,m}(w,U,z)$ for all $n$ and $m$.  
This implies that $x=w$,
$y=z$, and that $(Vf_n(y),f_m(x)) = (Uf_n(y),f_m(x))$ for all $n$ and
$m$.  It follows by part (iii) of Definition \ref{def:36} that
$V^*f_m(x) = U^*f_m(x)$ for all $m$.  However, given $h\in \mcal{H}(y)$ we
have
\[
(Vh,f_m(x)) = (h,V^*f_m(x)) = (h,U^*f_m(x)) = (Uh,f_m(x)).
\]
As above, it follows that $Vh = Uh$.  This is true for all $h\in
\mcal{H}(y)$ so that $U=V$ and we are done.  
\end{proof}

We are almost at the point were we can define a groupoid
representation, but first we need to deal with quasi-invariant
measures.  This will also be our introduction to the modular function which
has so far been missing in our treatment of groupoid dynamical systems.  

\begin{definition}
\label{def:39}
\index{quasi-invariant}
\index{modular function}
\index{induced measures}
\index[not]{$\Delta(\gamma)$}
\index[not]{$\nu,\nu\inv$}
Suppose $G$ is a locally compact Hausdorff groupoid with Haar system
$\lambda$.  Given a Radon measure $\mu$ on $G\unit$ we define
the {\em induced measures} $\nu$ and $\nu\inv$ to be the Radon
measures on $G$ defined by the equations
\begin{align*}
\nu(f) &:= \int_{G\unit}\int_G
f(\gamma)d\lambda^u(\gamma)d\mu(u),\quad\text{and} \\ 
\nu\inv(f) &:= \int_{G\unit}\int_G f(\gamma)d\lambda_u(\gamma)d\mu(u)
\end{align*}
for all $f\in C_c(G)$.  We call the measure $\mu$ {\em
  quasi-invariant} if $\nu$ and $\nu\inv$ are mutually absolutely
continuous.  In this case we write $\Delta$ for the Radon-Nikodym
derivative $d\nu/d\nu\inv$ and call it the {\em modular function} of
$\mu$.  If $\Delta \equiv 1$ $\nu$-almost everywhere then $\mu$ is
said to be {\em invariant}.
\end{definition}

\begin{remark}
Suppose $G$ is a groupoid with a Haar system $\lambda$. 
Given a measure $\mu$ on $G\unit$ we will often denote the measures on
$G$ induced from $\mu$ by the suggestive notation
\begin{align*}
\nu &:= \int_G \lambda^u d\mu(u), & \nu\inv &:= \int_G \lambda_ud\mu(u).
\end{align*}
\end{remark}

The reason for the terminology ``modular function'' is that $\Delta$
behaves like the modular function for a locally compact group.  We
will learn more about this relationship in Section
\ref{sec:groupprod}.  Moving on, we can use the following
theorem to choose $\Delta$ so that it is a groupoid homomorphism from
$G$ into $\R_+^\times$.  

\begin{theorem}[{\cite[Corollary 3.14]{hahnhaar}}]
\label{thm:quasi}
Given a quasi-invariant measure $\mu$ on the unit space of a groupoid
$G$ with Haar system $\lambda$ it is possible to choose the
modular function of $\mu$, $\Delta$, to be a Borel homomorphism from
$G$ to $\R_+^\times$.  Moreover, if $\mu'$ is another quasi-invariant
measure on $G\unit$ that is equivalent to $\mu$ so that $\mu' = g\mu$
for a suitable non-negative function $g$, and if $\Delta'$ is the
modular function of $\mu'$, then $\Delta'(\gamma) =
g(r(\gamma))\Delta(\gamma)g(s(\gamma))\inv$ $\nu$-almost everywhere, 
where $\nu$ is
the measure induced by $\mu$.
\end{theorem}

It is not immediately clear that there are such things as
quasi-invariant measures.  However, the following proposition shows
that they are easily constructed.  The details are thanks to Dana
Williams\index{Dana Williams}.  
It is also proved in
\cite[Pages 24,25]{groupoidapproach} and detailed in
\cite{coords}.  

\begin{prop}
\label{prop:50}
\index{saturation of a measure}
Suppose $G$ is a second countable locally compact groupoid
with Haar system $\lambda$.  Given a Radon
measure $\mu_0$ on $G\unit$ let $\nu_0=\int_G\lambda^ud\mu_0(u)$ 
be the measure
induced by $\mu_0$.  Then $\nu_0$ is a $\sigma$-finite measure on $G$
(but not necessarily finite).  Let $\nu$ be a probability measure on $G$
which is equivalent to $\nu_0$ and define $\mu$
to be the image of $\nu$ under the source map, i.e. $\mu = s_*\nu$.  
Then $\mu$ is quasi-invariant and, if $\mu_0$ was
quasi-invariant to begin with then $\mu_0$ is equivalent to $\mu$.
\end{prop}

\begin{remark}
\index[not]{$[\mu]$}
The measure $\mu$ is called the {\em saturation} of $\mu_0$ and is
denoted by $[\mu_0]$.  
\end{remark}

\begin{proof}
To see that $\nu_0$ is $\sigma$-finite it suffices to produce $f\in
\mcal{L}^1(G,\nu_0)$ such that $f(\gamma)>0$ for all $\gamma\in G$.
However, since $G$ is second countable, $G$ is $\sigma$-compact.  Let
$G = \bigcup K_n$ with each $K_n$ compact, and let $f_n\in C_c^+(G)$ be
such that $f_n(\gamma) >0$ for all $\gamma\in K_n$.  Define 
\[
\lambda(f_n)(u) = \int_G f_n(\gamma) d\lambda^u(\gamma)
\]
for all $n$.  Then $\lambda(f_n)\in C_c^+(G\unit)$ and, because 
$\mu_0$ is a Radon measure, 
\[
\nu_0(f_n) = \mu(\lambda(f_n)) = \alpha_n < \infty
\]
for all $n$.  We may as well assume that we have chosen $f_n$ such
that $\|f_n\|_\infty \leq 2^{-n}$ and $\alpha_n \leq 2^{-n}$.  Then $f
:= \sum_n f_n$ will do.  

Now let $\nu$ and $\mu$ be as in the statement of the proposition and
note that, by definition, $\mu(E) = \nu(s\inv(E))$.  Looking at the
characteristic functions, we see that 
\[
\int_{G\unit} f(u)d\mu(u) = \int_G f(s(\gamma))d\nu(\gamma).
\]
For convenience we let $\phi$ be the Radon-Nikodym derivative of
$d\nu/d\nu_0$ and assume, as we can, that $\phi(\gamma)\in (0,\infty)$
for all $\gamma\in G$.

Next, we let $\nu'=\int_G \lambda^u d\mu(u)$.  Since $(\nu')\inv$ is
nothing more than the composition of $\nu$ with the inversion map,
for the first part of
the proof it will suffice to show that if $f$ is a bounded
non-negative Borel function on $G$ such that 
\[
\int_G f(\gamma)d\nu'(\gamma) = 0
\]
then 
\[
\int_G \tilde{f}(\gamma) d\nu'(\gamma) = 0
\]
where $\tilde{f}(\gamma) = f(\gamma\inv)$.  We now compute as
follows:\footnote{Technically we need Proposition \ref{prop:57} to do
  these calculations.  However, in the interest of not getting
  sidetracked we will postpone these considerations until the next section.}
\[
0= \int_G f(\gamma)d\nu'(\gamma) = \int_{G\unit}\int_G
f(\gamma)d\lambda^u(\gamma) d\mu(u)
= \int_G \int_G f(\gamma)d\lambda^{s(\eta)}(\gamma)d\nu(\eta)
\]
which, using the definition of $\phi$, is
\begin{align*}
&= \int_G \int_G
f(\gamma)d\lambda^{s(\eta)}(\gamma)\phi(\eta)d\nu_0(\gamma) 
= \int_{G\unit}\int_G\int_G f(\gamma)\phi(\eta)
d\lambda^{s(\eta)}(\gamma)d\lambda^u(\eta) d\mu_0(u) \\
&= \int_{G\unit}\int_G\int_G f(\eta\inv\gamma)\phi(\eta)
d\lambda^u(\gamma)d\lambda^u(\eta)d\mu_0(u),
\end{align*}
where we used left invariance to get the last equality.  Now, switch
the order of integration, and use left invariance again to get
\[
0= \int_{G\unit}\int_G\int_G f(\eta\inv)\phi(\gamma\eta)
d\lambda^{s(\gamma)}(\eta) d\lambda^u(\gamma)d\mu_0(u) = \int_G \int_G \tilde{f}(\eta)
\phi(\gamma\eta)d\lambda^{s(\gamma)}(\eta) d\nu_0(\gamma)
\]
which, using the fact that $\phi\inv = d\nu_0/d\nu$, is 
\begin{equation}
\label{eq:27}
0=\int_G \int_G \tilde{f}(\eta) \phi(\gamma\eta) \phi(\gamma)\inv
d\lambda^{s(\gamma)}(\eta)d\nu(\gamma).
\end{equation}
Notice that, since $f$ is non-negative and 
$\phi(\gamma)>0$ for all $\gamma$, we have
\begin{equation}
\label{eq:28}
\int_G \tilde{f}(\eta)\phi(\eta\gamma)\phi(\gamma)\inv
d\lambda^{s(\gamma)}(\eta) = 0
\end{equation}
if and only if 
\begin{equation}
\label{eq:29}
\int_G \tilde{f}(\eta)d\lambda^{s(\gamma)}(\eta) = 0.
\end{equation}

At this point it follows from \eqref{eq:27} that there exists a
$\nu$-null set $N$ such that $\gamma\ne N$ implies that \eqref{eq:28}
holds.  Clearly $N$ is $s$-saturated:  $s\inv(s(N))=N$.  Thus
$s(N)$ is $\mu$-null and, since \eqref{eq:28} holds if and only if
\eqref{eq:29} does, we've established that 
\[
\int_G \tilde{f}(\eta)d\lambda^u(\eta) = 0
\]
for $\mu$-almost all $u$.  But then 
\[
\int_G \tilde{f}(\eta)d\nu'(\eta) = \int_{G\unit}\int_G
\tilde{f}(\eta)d\lambda^u(\eta)d\mu(u) = 0.
\]

For the second assertion suppose that $\mu_0$ is quasi-invariant.  But
then if $f$ is a bounded non-negative Borel function on $G\unit$
\[
\mu(f) = \int_G f(s(\gamma)) d\nu(\gamma) = \int_G
f(s(\gamma))\phi(\gamma)d\nu_0(\gamma)
\]
which, using the quasi-invariance of $\mu_0$, is 
\begin{align*}
&= \int_G f(r(\gamma))\phi(\gamma\inv)\Delta(\gamma) d\nu_0(\gamma) 
= \int_{G\unit}\int_G f(r(\gamma))\phi(\gamma\inv)\Delta(\gamma)
d\lambda^u(\gamma)d\mu_0(u) \\
&= \int_{G\unit} f(u) \left(\int_G \phi(\gamma\inv)\Delta(\gamma)
  d\lambda^u(\gamma)\right)d\mu_0(u).
\end{align*}
Since 
\[
\alpha(u) := \int_G \phi(\gamma\inv) \Delta(\gamma) d\lambda^u(\gamma)
\]
is a non-negative, (extended) real-valued function it follows that
$\mu\ll \mu_0$.  The argument is symmetric in $\mu$ and $\mu_0$ so
$\mu_0 \ll \mu$ and $\mu$ is equivalent to $\mu_0$. 
\end{proof}

Being able to restrict ourselves to finite quasi-invariant measures
will be useful later on.  

\begin{corr}
\label{cor:4}
\index{quasi-invariant}
Suppose $G$ is a second countable locally compact Hausdorff groupoid
with a Haar system.  Then every quasi-invariant measure on $G\unit$ is
equivalent to a finite quasi-invariant measure. 
\end{corr}
\begin{proof}
This follows immediately from Proposition \ref{prop:50} since the
measure $\mu$ constructed there is finite.  
\end{proof}

We are finally ready to define the notion of a unitary
representation of a groupoid.  It is perhaps notable that, unlike the
group case, groupoid representations are Borel creatures. 

\begin{definition}
\label{def:40}
\index{groupoid representation}
Let $G$ be a locally compact Hausdorff groupoid with a Haar system.  A
{\em groupoid representation} of $G$ is a triple $(\mu,G\unit*\mfrk{H},L)$ where
$\mu$ is a finite quasi-invariant measure on $G\unit$, $G\unit*\mfrk{H}$ is an
analytic Borel Hilbert bundle, and $L:G\rightarrow
\isom(G\unit*\mfrk{H})$ is a Borel homomorphism such that $L(\gamma) =
(r(\gamma),L_\gamma,s(\gamma))$ for some unitary
$L_\gamma:\mcal{H}(s(\gamma))\rightarrow \mcal{H}(r(\gamma))$.    
\end{definition}

\begin{remark}
\label{rem:12}
We have taken the quasi-invariant measure $\mu$ to be finite.  The
reason we have done this is because it is convenient.  However,
there is no harm in working with $\sigma$-finite quasi-invariant
measures either.  Any such measure is, at least in the second
countable case, equivalent to a finite one by
Corollary \ref{cor:4} and
Remarks \ref{rem:11} and \ref{rem:13} address how the the theory could
be expanded to include the $\sigma$-finite case using this fact.  
\end{remark}

\begin{example}
\index{left regular representation!groupoid}
\label{ex:18}
Let $\mu$ be a finite quasi-invariant measure on the unit space of a
second countable locally compact Hausdorff groupoid $G$ equipped with
a Haar system $\lambda$.  Let $G\unit*L^2(\lambda)$ be as in
Example \ref{ex:17}.  Define $L:G\rightarrow
\isom(G\unit*L^2(\lambda))$ by $L(\gamma) =
(r(\gamma),L_\gamma,s(\gamma))$ where 
$L_\gamma : L^2(G^{s(\gamma)},\lambda^{s(\gamma)}) \rightarrow
L^2(G^{r(\gamma)},\lambda^{r(\gamma)})$ is defined by the formula
\[
(L_\gamma\xi)(\eta) = \xi(\gamma\inv\eta).
\]
All we really need to show is that $L$ is a Borel map.  We will use
Lemma \ref{lem:51}.  First, observe that $r\circ L = r$ is clearly
Borel.   We must show that the function 
\[
\psi_{n,m}(\gamma) = (L_\gamma \xi_n(s(\gamma)),\xi_m(r(\gamma))) =
\int_G \overline{\xi_n(s(\gamma))(\gamma\inv\eta)}\xi_m(r(\gamma))(\eta)d\lambda^{r(\gamma)}(\eta)
\]
is Borel for all $n,m$ where $\xi_n$ is the 
fundamental sequence constructed in Example
\ref{ex:17}.  However, the $\xi_n$ are given by $\xi_n(u)(\gamma) =
f_n(\gamma)$ 
where $f_n$ is a continuous compactly supported function for each $n$. 
So in this case $\psi_{n,m}$
is actually continuous for all $n,m$.  
Thus the triple $(\mu,G\unit*L^2(\lambda),L)$ forms a representation
of $G$.  This type of
representation is called the {\em left regular representation}
associated to $\mu$.  
\end{example}

\subsection{Decomposable Representations}

Now that we have covered the groupoid half of a covariant
representation we have to deal with the $C^*$-algebraic half.
Unfortunately this is not entirely straightforward 
since we have to deal with the fibred structure on the algebras and this
means working with decomposable operators and representations.  
We will not present a self contained exposition of 
decomposable operators here.  Instead we will only provide 
those proofs which seem relevant.  For a more complete treatment 
the reader is referred to \cite[Section F.3]{tfb2} and 
\cite[Chapter 4]{invitation}.

\begin{remark}
We will always assume, unless explicitly stated otherwise, that all of
our representations are nondegenerate.
\end{remark}

We begin by forming a Hilbert space out of a
given Borel Hilbert bundle.  The definition is fairly
straightforward. 

\begin{definition}
\label{def:37}
\index{direct integral}
\index[not]{$L^2(X*\mfrk{H},\mu)$}
Suppose $X*\mfrk{H}$ is an analytic Borel Hilbert bundle and
$\mu$ is a measure on $X$.  Let 
\[
\mcal{L}^2(X*\mfrk{H},\mu) = \{f\in B(X*\mfrk{H}):\text{$x\mapsto
  \|f(x)\|^2$ is integrable}\}
\]
and give $\mcal{L}^2(X*\mfrk{H},\mu)$ the operations of pointwise
addition and scalar multiplication.  
Let $L^2(X*\mfrk{H},\mu)$ be the quotient of
$\mcal{L}^2(X*\mfrk{H},\mu)$ where functions which agree $\mu$-almost
everywhere are identified.  Equipped with the operations coming from
$\mcal{L}^2(X*\mfrk{H},\mu)$ and the inner product
\[
(f,g) = \int_X (f(x),g(x)) d\mu(x)
\]
$L^2(X*\mfrk{H},\mu)$ becomes a Hilbert space known as the {\em direct
  integral} of $\mfrk{H}$ with respect to $\mu$.  
\end{definition}

\begin{remark}
The fact that the integral even makes sense follows from Proposition
\ref{prop:49} and an
application of the Cauchy-Schwartz inequality.  The rest of the
assertions made in Definition \ref{def:37} are straightforward to
verify.  In any case, they are all addressed in \cite[Appendix
F.2]{tfb2}.  It is worthwhile to point out that the direct integral is
classically denoted 
\[
\int_X^\oplus \mcal{H}(x)d\mu(x).
\]
\end{remark}

The following proposition, which we cite without proof, guarantees that
we will not have to deal with any nonseparable weirdness.  

\begin{prop}[{\cite[Lemma F.17]{tfb2}}]
\label{prop:54}
If $X*\mfrk{H}$ is an analytic Borel Hilbert bundle and 
$\mu$ is a finite Borel measure on $X$ then $L^2(X*\mfrk{H},\mu)$ is a
separable Hilbert space. 
\end{prop}

\begin{remark}
\label{rem:11}
Proposition \ref{prop:54} is one reason why we chose
to restrict ourselves to finite quasi-invariant measures in Definition
\ref{def:40}.  However, as noted in Remark \ref{rem:12}, if $\mu$ is
$\sigma$-finite then it is equivalent to a finite measure.  It
is easy to see that, up to unitary equivalence, 
$L^2(X*\mfrk{H},\mu)$ only depends on the equivalence class of $\mu$.
Thus, $L^2(X*\mfrk{H},\mu)$ is separable in the $\sigma$-finite case as well.
\end{remark}

We take a moment to discuss pull back Borel Hilbert bundles. 

\begin{example}[{[\cite[Example F.18]{tfb2}}]
\label{ex:35}
\index{pull back}
Suppose that $X*\mfrk{H}$ is an analytic Borel Hilbert bundle with
fundamental sequence $\{f_n\}$ and that $\sigma:Y\rightarrow X$ is a
Borel map.  Then we can form the pull back Borel Hilbert bundle
\[
\sigma^*(X*\mfrk{H}) := \{(y,h):h\in \mcal{H}(\sigma(y))\}.
\]
We use Proposition \ref{prop:113} to give $\sigma^*(X*\mfrk{H})$ the
Borel structure coming from the fundamental sequence $\{f_n\circ
\sigma\}$.  It follows that $f\in B(X*\mfrk{H})$ implies that $f\circ
\sigma\in B(\sigma^*(X*\mfrk{H}))$.  If $\nu$ is a finite Borel
measure on $Y$ and if $\sigma_*\nu$ is the push forward measure on $X$
then it turns out that $W(f)(y)= f(\sigma(y))$ defines an isometry 
\[
W:L^2(X*\mfrk{H},\sigma_*\nu) \rightarrow L^2(\sigma^*(X*\mfrk{H}),\nu)
\]
which is an isomorphism if $\sigma$ is a Borel isomorphism. 
\end{example}

Moving on, the basic idea will be that certain representations of $C^*$-algebras
on the direct integral $L^2(X*\mfrk{H},\mu)$ can be decomposed into
representations on the fibres $\mcal{H}(x)$.  In order to make sense
out of this we will eventually need the following

\begin{definition}
Suppose that $X*\mfrk{H}$ is an analytic Borel Hilbert bundle and $\mu$
is a finite measure on $X$.  An
operator $T$ on $L^2(X*\mfrk{H},\mu)$ is called {\em diagonal} if
there is a bounded Borel function $\phi\in \mcal{B}^b(X)$ such that 
\[
Th(x) = \phi(x)h(x)
\]
for $\mu$-almost every $x$.  The collection of diagonal operators on
$L^2(X*\mfrk{H},\mu)$ is denoted by $\Delta(X*\mfrk{H},\mu)$.  If
$\phi\in \mcal{B}^b(X)$ then the associated diagonal operator is
denoted by $T_\phi$.  
\end{definition}

\begin{prop}{{\cite[Lemma F.15]{tfb2}}}
Suppose that $X*\mfrk{H}$ is an analytic Borel Hilbert bundle and that
$\mu$ is a finite Borel measure on $X$.  Then $\Delta(X*\mfrk{H},\mu)$
is an abelian von Neumann subalgebra of $B(L^2(X*\mfrk{H},\mu))$, and
the map $\phi\mapsto T_\phi$ induces an isomorphism of
$L^\infty(X,\mu)$ onto $\Delta(X*\mfrk{H},\mu)$. 
\end{prop}

\begin{remark}
\label{rem:13}
As discussed, we are assuming that the measure $\mu$ is finite, but these
definitions and their theory can be extended to
$\sigma$-finite measures by again using the fact that every $\sigma$-finite
measure is equivalent to a finite one.  We actually do this for the theory
of direct integrals of operators and representations.\footnote{For the
  curious reader,   this was actually added in after much of the
  thesis was finished  because it is needed in Section
  \ref{sec:crossedstab}.} 
The key fact is that the unitary induced by the
Radon-Nikodym derivative of two equivalent measures will intertwine
the diagonal operators, and as such will respect the theory of
decomposable operators.  In short, the finiteness of $\mu$ will not be an
essential part of the theory of covariant representations and is, as
stated in Remark \ref{rem:12}, a convenience.  
\end{remark}

Now, before we can decompose representations on $L^2(X*\mfrk{H})$ we
have to be able to decompose operators.  In order to do this we must
have some idea of what happens after such a decomposition.  

\begin{definition}
\label{def:41}
\index{Borel field of operators}
Suppose $X*\mfrk{H}$ is an analytic Borel Hilbert
bundle with fundamental sequence $\{f_n\}$.  A family of bounded linear maps
$T(x):\mcal{H}(x)\rightarrow \mcal{H}(x)$ is a {\em Borel field of
  operators} if
\[
x\mapsto (T(x)f_n(x) , f_m(x))
\]
is Borel for all $n$ and $m$.  
\end{definition}

Of course, we would like to see how this relates to the Borel
structure on $X*\mfrk{H}$. 

\begin{prop}
Suppose $X*\mfrk{H}$ is an analytic Borel Hilbert
bundle and that we have a family of bounded 
linear maps $T(x):\mcal{H}(x)\rightarrow
\mcal{H}(x)$.  We can define a bundle map
$\widehat{T}:X*\mfrk{H}\rightarrow X*\mfrk{H}$ such that
$\widehat{T}(x,h) = (x,T(x)h)$.  Then $\widehat{T}$ is
Borel if and only if $\{T(x)\}$ is a Borel field of operators.  
\end{prop}

\begin{remark}
Thus, a
Borel field of operators is essentially nothing more than an
endomorphism of a Borel Hilbert bundle.  
\end{remark}

\begin{proof}
Let $e_l$ be a special orthogonal fundamental sequence for
$X*\mfrk{H}$.  Suppose $\widehat{T}$ is Borel.  
Then $\overline{\bar{e}_k\circ
  \widehat{T}\circ e_l}$ is clearly Borel and, tracing through the
  definitions,
\begin{equation}
\label{eq:30}
\overline{\bar{e}_k\circ\widehat{T}\circ e_l(x)} = 
\overline{\bar{e}_k(x,T(x)e_l(x))} = 
(T(x)e_l(x),e_k(x)).
\end{equation}

For the reverse direction, suppose $x\mapsto (T(x)e_l(x),e_k(x))$ is
Borel for all $l$ and $k$. We want to show that $\bar{e}_k \circ
\widehat{T}$ is Borel for all $k$.  Well, using the Fourier Identity,
\[
\bar{e}_k\circ \widehat{T}(x,h) = (e_k(x),T(x)h) = 
\sum_{l} (e_l(x),h)\overline{(T(x)e_l(x),e_k(x))}.
\]
However, $(x,h)\mapsto (e_l(x),h)$ is Borel since $e_l$ is a
fundamental sequence and $x\mapsto (T(x)e_l(x),e_k(x))$ is Borel by
assumption.  This suffices to show that 
$\bar{e}_k\circ\widehat{T}$ is Borel and we are done.  
\end{proof}

Given a Borel field of bounded linear 
operators we would like to be able to glue them
together to form a bounded linear operator on the direct integral of
the Hilbert spaces.  

\begin{prop}
\label{prop:51}
Suppose $X*\mfrk{H}$ is an analytic Borel Hilbert bundle and that $\mu$ is a
$\sigma$-finite measure on $X$.  Let $\{T(x)\}$ be a Borel field of bounded
linear operators
such that 
\[
\Lambda := \ess\sup_{x\in X}\|T(x)\| < \infty.
\]
Then there exists a bounded linear operator
$T\in B(L^2(X*\mfrk{H},\mu))$ defined by
$Tf(x) = T(x)f(x)$ for all $f\in \mcal{L}^2(X*\mfrk{H},\mu)$ such that 
$\|T\| = \Lambda$.  
\end{prop}

\begin{remark}
\index{direct integral}
The classical notation for the operator defined in Proposition
\ref{prop:51} is
\[
\int_X^\oplus T(x) d\mu(x).
\]
This, of course, is in line with the direct integral notation for Borel Hilbert
bundles. We will also occasionally refer to $T$ as the {\em direct
  integral} of the $T(x)$.  
\end{remark}

\begin{proof}
Let us start by supposing that $\mu$ is actually a finite measure.  
The first thing we have to do is make sure that everything is
well defined. It is not too difficult to show that if $S$ is the
countable family of rational linear combinations of some special
orthogonal fundamental
sequence for $X*\mfrk{H}$ then $\{h(x):h\in S\}$ is dense in
$\mcal{H}(x)$ for all $x\in X$.  This implies that 
\[
\|T(x)\| = \sup_{h\in S, h(x)\ne 0} \|T(x)h(x)\|\|h(x)\|\inv.
\]
It then follows that the function $x\mapsto \|T(x)\|$ is Borel and
that taking the essential supremum $\Lambda$ makes sense.

Given a Borel field of operators and $f\in B(X*\mfrk{H})$ we have 
$\widehat{T}\circ f(x) = T(x)f(x)$ for all $x\in X$.  Since
$\widehat{T}$ is Borel, we have $\widehat{T}\circ f\in B(X*\mfrk{H})$.  
Now, define $T$ on $\mcal{L}^2(X*\mfrk{H},\mu)$ as in the statement of the
proposition and observe that $Tf = \widehat{T}\circ f$ so that $Tf$ is
a Borel section.  Furthermore $\|Tf(x)\| \leq \Lambda \|f(x)\|$
$\mu$-almost everywhere so, because $x\mapsto \|f(x)\|^2$ was integrable,
$x\mapsto \|Tf(x)\|^2$ must be integrable as well.  Hence $Tf \in
\mcal{L}^2(X*\mfrk{H},\mu)$.  It is straightforward to show that $T$ is a
linear map.  Finally 
\[
\|Tf\|^2 = \int_X \|T(x)f(x)\|^2 d\mu(x) \leq \Lambda^2 \int_X \|f(x)\|^2
d\mu(x) = \Lambda^2\|f\|^2
\]
so that $T$ factors to an operator on $L^2(X*\mfrk{H})$ and
is bounded with norm at most $\Lambda$.  

Next we will show the reverse inequality.  First, because of the last
part of Proposition \ref{prop:49} we can assume without loss of
generality that $X*\mfrk{H}$ is trivial.  If not, decompose
$X*\mfrk{H}$ as a disjoint union of trivial bundles $X*\mfrk{H} =
\coprod_{d=0}^\infty X_d\times \mcal{H}_d$.  Then if we show
that $\ess\sup_{x\in X_d} \|T(x)\| \leq \|T\|$ for all $d$ it will
follow that $\Lambda \leq \|T\|$. So assume that $X*\mfrk{H} = X\times
\mcal{H}$ is trivial and note that $L^2(X*\mfrk{H},\mu) =
L^2(X,\mcal{H},\mu)$.  Suppose $h\in\mcal{H}$ is a unit vector 
and observe that, since $\mu$ is finite, we can view $h$ as a constant
function in $L^2(X,\mcal{H},\mu)$.  It follows that $x\mapsto
\|T(x)h\|$ is a scalar function in $L^2(X,\mu)$.  Now, take $f\in
\mcal{L}^2(X,\mu)$ such that $\|f\|_2 =1$.  Observe that we can define
the function $f\otimes h(x) = f(x)h$ for all $x\in X$, that $f\otimes h\in
L^2(X,\mcal{H})$, and $\|f\otimes h\|_2 = 1$.  It follows that 
\begin{equation}
\label{eq:34}
\int_X \|T(x)h\|^2|f(x)|^2 d\mu(x) = \int_X \|T(f\otimes h)(x)\|^2d\mu(x) =
\|T(f\otimes h)\|^2 \leq \|T\|^2.
\end{equation}
We then conclude from the following general nonsense that $\|T(x)h\|\leq \|T\|$
everywhere off a $\mu$-null set $N_h$.  Let $\phi(x) = \|T(x)h\|^2$
and define $M_\phi$ to be the multiplication operator on $L^1(X)$.
Now, given $k\in L^1(X)$ such that $\|k\|_1 = 1$ let $f(x) =
\sqrt{|k(x)|}$ and observe that $f\in L^2(X)$ with $\|f\|_2 = 1$.  
It follows from \eqref{eq:34} that $\|M_\phi  k\|_1\leq \|T\|^2$.
Hence $\|\phi\|_\infty = \|M_\phi\| \leq \|T\|^2$.  Thus, as required, there
exists a $\mu$-null set $N_h$ such that $x\notin N_h$ implies
$\|T(x)h\|\leq \|T\|$.  Next, let $S$ be the
set of all rational linear combinations of a countable basis for
$\mcal{H}$ such that $\|h\|=1$ for all $h\in S$.  
We can find a new $\mu$-null set $N=\bigcup_{h\in S} N_h$
such that given $x\notin N$ we have $\|T(x)h\|\leq \|T\|$ for all
$h\in S$.  
It follows that $\|T(x)\|\leq \|T\|$ $\mu$-almost everywhere and we are done.  

Now suppose $\mu$ is $\sigma$-finite.  Then in the usual fashion we
can find a finite measure $\nu$ which is equivalent to $\mu$.  Let
$d\nu/d\mu$ be the Radon-Nikodym derivative and define
$U:L^2(X*\mfrk{H},\mu) \rightarrow L^2(X*\mfrk{H},\nu)$ by $Uf(x)=
(d\mu/d\nu)^{1/2}(x)f(x)$.  It is straightforward to show that $U$ is a
unitary.  Since $\nu$ is finite we can form the direct integral 
\[
T' = \int_X^\oplus T(x) d\nu(x)
\]
on $L^2(X*\mfrk{H},\nu)$.  This allows us to define an operator $T$ on
$L^2(X*\mfrk{H},\mu)$ by $T = U^* T' U$.  We immediately have $\|T\|
= \|T'\| = \Lambda$.  Furthermore, we can compute
\begin{align*}
Tf(x) &= U^* T' Uf(x) = \left(\frac{d\nu}{d\mu}(x)\right)^{1/2} T' Uf(x)
= \left(\frac{d\nu}{d\mu}(x)\right)^{1/2} T(x)U f(x) \\
&= \left(\frac{d\nu}{d\mu}(x)\right)^{1/2} T(x)
\left(\frac{d\mu}{d\nu}(x)\right)^{1/2} f(x) = T(x)f(x). \qedhere
\end{align*}
\end{proof}

We can now ``integrate'' operators and will shortly describe how to
integrate representations.  Our eventual goal will be
to show that $C_0(X)$-linear representations all have this form. 

\begin{prop}
\label{prop:41}\index{direct integral}
Suppose $X$ is a second countable locally compact Hausdorff space, $A$
is a separable $C_0(X)$-algebra, $X*\mfrk{H}$ is an analytic Borel
Hilbert bundle and $\mu$ is a $\sigma$-finite measure on $X$.  Given a
collection of representations $\pi_x:A(x)\rightarrow B(\mcal{H}(x))$
for all $x\in X$ such that given $a\in A$ the set $\{\pi_x(a(x))\}$ is
a Borel field of operators then we can form a representation 
\[
\pi = \int_X^\oplus \pi_x d\mu(x)
\]
of $A$ on $L^2(X*\mfrk{H},\mu)$ 
called the {\em direct integral} and defined for $a\in A$ by 
\[
\pi(a) = \int_X^\oplus \pi_x(a(x))d\mu(x).
\]
\end{prop}

\begin{remark}
Of course, $\pi = \int_X^\oplus \pi_x d\mu(x)$ is just the classic
direct integral notation of representations.  In fact, as long as you
are willing to confuse $\pi_x$ with its lift to $A$, Proposition
\ref{prop:41} and the upcoming Proposition \ref{prop:53} are just the
classic theory of decomposable representations. 
\end{remark}

\begin{proof}
Let $\{\pi_x\}$ be as above. Given $a\in A$ we assumed that
$\{\pi_x(a(x))\}$ is a Borel field of operators.  Since
$\|\pi_x(a(x))\| \leq \|a(x)\| \leq \|a\|$ for all $x \in X$, it
follows that $\{\pi_x(a(x))\}$ is bounded by $\|a\|$.  Therefore we
can form the direct integral 
\[
\pi(a) = \int_X^\oplus \pi_x(a(x))d\mu(x).
\]
All we need to do is show that $\pi$ preserves the algebraic
operations.  However, given $f\in L^2(X*\mfrk{H},\mu)$ we have 
\[
\pi(a+b)f(x) = \pi_x((a+b)(x))f(x) = \pi_x(a(x))f(x)+\pi_x(b(x))f(x) =
(\pi(a)+\pi(b))f(x).
\]
It is just as easy to show that the rest of the operations are
preserved.  
\end{proof}

Now that we know what ``integrated'' operators and representations 
look like we can define the decomposable ones.  

\begin{definition}
\index{decomposable operator}
Given an analytic Borel Hilbert bundle $X*\mfrk{H}$ and a finite measure $\mu$
on $X$, an operator $T\in B(L^2(X*\mfrk{H},\mu))$ is called
{\em decomposable} if there exists an essentially bounded Borel field
of operators $\{T(x)\}_{x\in X}$ such that
\[
T = \int_X^\oplus T(x) d\mu(x).
\]
\end{definition}

This lemma shows that the decomposition of an operator is unique
almost everywhere.  
\begin{lemma}[{\cite[Lemma F.20]{tfb2}}]
\label{lem:12}
Suppose that $X*\mfrk{H}$ is a Borel Hilbert bundle and that $\mu$ is
a finite measure on $X$.  Let $\{T(x)\}$ be an essentially bounded
Borel field of operators and let $T$ be the direct integral operator
on $L^2(X*\mfrk{H},\mu)$.  
\begin{enumerate}
\item If $T=0$ then $T(x)=0$ for $\mu$-almost all $x$. 
\item If $\{T'(x)\}$ is another essentially bounded Borel field of
  operators such
  that $T = \int_X^\oplus T'(x) d\mu(x)$ then
  $T'(x) = T(x)$ for $\mu$-almost all $x$. 
\end{enumerate}
\end{lemma}
\begin{proof}
Clearly it suffices to prove part (a).  Let $\{e_l\}$ be a
special orthogonal fundamental sequence.  Then $e_l\in
L^2(X*\mfrk{H},\mu)$ for all $l$ and for all $l$ and $k$ we have 
\[
(T(x)e_l(x)|e_k(x)) = 0\quad\text{for $\mu$-almost all $x$.}
\]
It follows that $T(x) = 0$ almost everywhere. 
\end{proof}

Now, if $T$ is the direct integral of a Borel
field of operators $\{T(x)\}$ then, for all $\phi\in L^\infty(X,\mu)$, 
\[
TT_\phi(f)(x) = \phi(x)T(x)f(x) = T_\phi T f(x).
\]
Thus, $T$ is in the commutant of the von Neumann algebra
$\Delta(X*\mfrk{H},\mu)$.  It is a deep result that this characterizes
the decomposable operators.  

\begin{theorem}[{\cite[II.2.5 Theorem 1]{dixmier}}]
\label{thm:decomposable}
Suppose that $X*\mfrk{H}$ is an analytic Borel Hilbert bundle and that
$\mu$ is a finite measure on $X$.  Let $T\in B(L^2(X*\mfrk{H},\mu))$.
Then $T$ is decomposable if and only if $T$ is in the commutant of the
diagonal operators $\Delta(X*\mfrk{H},\mu)$.  
\end{theorem}

The reason we went through all of this is that we would like to use
this decomposition theorem to decompose certain representations of
$C_0(X)$-algebras so that they are given by direct integrals as in
Proposition \ref{prop:41}.  

\begin{definition}
\index{czeroofxlinear@$C_0(X)$-linear}
Suppose $X$ is a second countable locally compact Hausdorff space, 
$A$ is a $C_0(X)$-algebra, $X*\mfrk{H}$ is an analytic
Borel Hilbert bundle and $\mu$ is a finite Borel measure on $X$.  
We say a representation $\pi:A\rightarrow
B(L^2(X*\mfrk{H},\mu))$ is {\em $C_0(X)$-linear} if 
\[
\pi(\phi\cdot a) = T_\phi \pi(a).
\]
for all  $a\in A$ and $\phi\in C_0(X)$.  
\end{definition}

This definition tells us what kind of representations are
decomposable, so now lets decompose them. 

\begin{prop}
\label{prop:53}
Suppose $X$ is a second countable locally compact Hausdorff space, 
$A$ is a separable $C_0(X)$-algebra, $X*\mfrk{H}$ is an analytic
Borel Hilbert bundle and $\mu$ is a finite Borel measure on $X$.  
Given a $C_0(X)$-linear representation $\pi$ of $A$ on
$L^2(X*\mfrk{H},\mu)$ there exists
possibly degenerate 
representations $\pi_x : A(x)\rightarrow B(\mcal{H}(x))$ such that
given $a\in A$ the set $\{\pi_x(a(x))\}$ is an essentially bounded
Borel field of operators and 
\[
\pi = \int_X^\oplus \pi_xd\mu(x)
\]
Furthermore, the representations $\pi_x$ are nondegenerate
$\mu$-almost everywhere and are uniquely
determined up to a $\mu$-null set.  
\end{prop}

\begin{remark}
Given $\pi$ and $A$ as above we will generally refer to $\{\pi_x\}$ as a 
decomposition of $\pi$ and call $\pi$ the direct integral of $\{\pi_x\}$.
\end{remark}

\begin{proof}
Suppose $a\in A$ and let $\Phi_A$ implement the $C_0(X)$-action on
$A$.  Then 
\[
T_\phi \pi(a) = \pi(\Phi_A(\phi) a) =
\pi(\Phi_A(\overline{\phi})a^*)^* 
= (T_{\overline{\phi}} \pi(a^*))^* = \pi(a)T_\phi
\]
where we used the fact that $\Phi_A(\phi)$ is in the center of the
multiplier algebra.  Thus $\pi(a)\in (\Delta(X*\mfrk{H}))'$ and it
follows from Theorem \ref{thm:decomposable} that there exists a Borel
field of bounded linear operators $\{\pi'_x(a)\}$ such that $\pi(a) =
\int_X^\oplus\pi'_x(a)d\mu(x)$.  

Define $\pi'_x:A\rightarrow B(\mcal{H}(x))$ such that $\pi'_x(a)$ is given
by the above decomposition.  We would be done if not for the fact 
that each $\pi'_x$ is
only a representation ``almost everywhere.''  
Now, it uses the separability of $A$ in a
fundamental way, but we can actually make our choices so that each
$\pi'_x$ is a $*$-homomorphism.  The details are worked out for a
slightly different context in
\cite[Section 4.2]{invitation} and have also been included here.
Suppose $a,b\in A$.  It is clear that
$\{\pi'_x(a)+\pi'_x(b)\}$ is an essentially bounded 
Borel field of operators.  As such we have the direct integral 
\[
T = \int_X^\oplus \pi'_x(a)+\pi'_x(b) d\mu(x).
\]
Furthermore, it is clear that for all $f\in L^2(X*\mfrk{H})$ we have
\begin{align*}
Tf(x) &= \pi'_x(a)f(x) + \pi'_x(b)f(x) =\pi(a)f(x) + \pi(b)f(x).
\end{align*}
Thus $T = \pi(a)+\pi(b) = \pi(a+b)$ and  
$\{\pi'_x(a)+\pi'_x(b)\}$ is another decomposition for $\pi(a+b)$.  
It follows from Lemma \ref{lem:12}
that there exists a $\mu$-null set $N_{a,b}$ such that for all
$x\not\in N_{a,b}$ we have $\pi_x'(a+b) = \pi_x'(a)+\pi_x'(b)$.  Now,
by enlarging $N_{a,b}$ we can use similar arguments to 
assume that for $x\not\in N_{a,b}$ we have
\begin{align}
\label{eq:32}
\pi_x'(a+b) &= \pi_x'(a)+\pi_x'(b), \\
\pi_x'(ab) &= \pi_x'(a)\pi_x'(b), \\
\pi_x'(a^*) &= \pi_x'(a)^*,\quad \text{and}\\
\label{eq:33}
\pi_x'(\lambda a) &= \lambda\pi_x'(a)\quad \text{for all $\lambda\in \Q+i\Q$.}
\end{align}
Let $\{a_i\}$ be a countable dense sequence in $A$ and let
$N=\bigcup_{i,j}N_{a_i,a_j}$.  Then $N$ is a $\mu$-null set and for
all $x\not\in N$ we know that the identities
\eqref{eq:32}-\eqref{eq:33} hold for
any elements of $\{a_i\}$. Let $S$ be the countable family of
all finite sums of elements of the form $rb_1\cdots b_n$ where $r\in
\Q+i\Q$ and $b_1,\ldots,b_n\in\{a_i\}\cup\{a_i^*\}$.  It is
straightforward to show that, for all $x\notin N$, \eqref{eq:32}-\eqref{eq:33} extend to all
of $S$.  It follows from Proposition \ref{prop:51}
that given $a\in S$ we have 
$\|\pi_x'(a)\| \leq \|\pi(a)\|\leq \|a\|$ for $\mu$-almost all $x$.
Thus, by enlarging $N$ by another countable family of null sets, we
can assume that given $x\notin N$ we have
\[
\|\pi_x'(a)\| \leq \|a\|
\]
for all $a\in S$.  Thus $\pi_x'$ is a bounded $*$-homomorphism on the
normed $*$-algebra $S$.  Therefore,
given $x\notin N$ we can extend $\pi_x'$ to a $*$-homomorphism
$\tilde{\pi}_x:A\rightarrow B(\mcal{H}(x))$. If $x\in N$ we let
$\tilde{\pi}_x$ be the trivial representation on $\mcal{H}(x)$.

Now we must show that $\{\tilde{\pi}_x(a)\}$ is also a decomposition
of $\pi(a)$ for all $a\in A$.  Suppose $a\in A$, that
$a_{i_j}\rightarrow a$, and that $\{f_n\}$ is a fundamental sequence.
Fix $n$ and $m$ and let $\phi:X\rightarrow \C$ be defined by 
$\phi(x) = (\tilde{\pi}_x(a)f_n(x)|f_m(x))$ and $\phi_j:X\rightarrow
\C$ be defined by $\phi_j(x) = (\pi'_x(a_{i_j})f_n(x)|f_m(x))$.  Then for
$x\not\in N$ we have $\phi(x) = \lim_j \phi_j(x)$.  Since the
pointwise limit of Borel functions is Borel, it follows that $\phi|_N$
is Borel.  However, $\phi(N) = 0$ so that $\phi$ must be Borel
everywhere. 
It follows that $\{\tilde{\pi}_x(a)\}$ is a Borel field
of operators that is clearly bounded by $\|a\|$.  Now, given $a\in A$
let $\{a_{i_j}\}$ be a sequence in $\{a_i\}$ converging to $a$.  It is
not difficult to use the fact that $\tilde{\pi}_x$ is an extension of
$\pi'_x$ from $S$ and the fact that each $\tilde{\pi}_x$ is a
homomorphism to check that 
\begin{align*}
&\left\| \int^\oplus_X \pi'_x(a)d\mu(x) - \int^\oplus_X \tilde{\pi}_x(a)
  d\mu(x) \right\| \\ &\leq \left\|\int^\oplus_X \pi_x'(a)d\mu(x) -
  \int^\oplus_X \pi_x'(a_{i_j}) d\mu(x)
\right\| + \left\|\int^\oplus_X \tilde{\pi}_x(a_{i_j})
  d\mu(x) - \int^\oplus_X \tilde{\pi}_x(a)d\mu(x)\right\| \\
&= \| \pi(a)-\pi(a_{i_j})\| + \left\|\int^\oplus_X
  \tilde{\pi}_x(a_{i_j}-a)d\mu(x)
\right\| 
\\ &\leq 2\|a-a_{i_j}\| \rightarrow 0. 
\end{align*}
Hence $\pi(a) = \int^\oplus_X \pi'_x(a)d\mu = \int^\oplus_X
\tilde{\pi}_x(a) d\mu$ and, as claimed, we can decompose $\pi$ so that
each $\tilde{\pi}_x$ is a homomorphism. 

Next, we show that the $\tilde{\pi}_x$ are nondegenerate almost
everywhere.  This proof is taken from \cite[Proposition
4.2.2]{invitation}.  Let $a_i$ be an approximate identity for $A$ and
observe that, because $\pi$ is nondegenerate $\pi(a_i)\rightarrow \id$
strongly. 
Let $e_j$ be a special fundamental orthogonal sequence of $X*\mfrk{H}$
and note that $e_j\in \mcal{L}^2(X*\mfrk{H},\mu)$ for all $j$.  Now,
we compute
\[
\lim_{i\rightarrow \infty}\|\pi(a_i)e_j - e_j\|^2 = 
\lim_{i\rightarrow\infty}\int_X \|\tilde{\pi}_x(a_i)e_j(x) - e_j(x)\|^2 d\mu(x) 
= 0
\]
for all $j$.  Since convergence in mean implies a subsequence
converges almost everywhere it follows that for each $j$ we can find a
subsequence $i_{k_j}$ such that
\[
\lim_{k\rightarrow \infty}\|\tilde{\pi}_x(a_{i_{k_j}})e_j(x)-e_j(x)\| = 0
\]
for all $x$ not in some null set $N_j\subset X$.  By induction we can
arrange that the $(j+1)$st subsequence is a subsequence of the $j$th.
Then consider the diagonalization $\{a_{i_{k_k}}\}$.  As long as
$x\notin \bigcup_{j} N_j$ we have 
\[
\lim_{k\rightarrow\infty}\|\tilde{\pi}_x(a_{i_{k_k}})e_j(x) - e_j(x)\| = 0
\]
for all $j$.  
Since the $\{e_j(x)\}$ form a basis for $\mcal{H}(x)$, this is enough
to show that $\tilde{\pi}_x$ is nondegenerate almost everywhere.  

Now, let $I_x$ be the ideal in $A$ such that $A(x) = A/I_x$. 
We would like to claim that $I_x\subset \ker\tilde{\pi}_x$ for all $x\in
X$.  However, we are going to have to deal with more almost everywhere
nonsense.  Recall that $I_x$ is generated by elements of the form
$\phi\cdot a$ where $a\in A$ and $\phi\in C_0(X)$ such that $\phi(x) =
0$.  Given $a\in A$, $\phi\in C_0(X)$ and 
$f\in \mcal{L}^2(X*\mfrk{H})$ we have 
\[
\int^\oplus_X \tilde{\pi}_x(\phi\cdot a)d\mu(x) f = 
\pi(\phi\cdot a)f = T_\phi \pi(a)f = T_\phi \int^\oplus_X \tilde{\pi}_x(a)
d\mu(x)f
\]
where $T_\phi$ is the diagonal operator associated to $\phi$.  
This implies that there exists a null set $N(\phi,a,f)$ such that 
\[
\tilde{\pi}_x(\phi\cdot a)f(x) = \phi(x)\tilde{\pi}_x(a)f(x)
\]
for all $x\not\in N(\phi,a,f)$.  However, if we let $\{\phi_i\}$ and
$\{a_j\}$ be countable dense subsets of $C_0(X)$ and $A$, respectively,
and $e_l$ a special orthogonal fundamental sequence of $X*\mfrk{H}$
then we can pick a single null set $N = \bigcup N(\phi_i,a_j,e_l)$
such that given $x\not\in N$ we have 
\[
\tilde{\pi}_x(\phi_i\cdot a_j)e_l(x) = \phi_i(x) \tilde{\pi}_x(a_j)e_l(x)
\]
for all $i,j,l$.  Since each $\{e_l(x)\}$ forms an orthogonal basis
for $\mcal{H}(x)$ (plus some zero vectors) it follows that given
$x\notin N$ we have $\tilde{\pi}_x(\phi_i\cdot a_j) = \phi_i(x)\tilde{\pi}_x(a_j)$ for
all $i,j$.  By continuity we get 
$\tilde{\pi}_x(\phi\cdot a) = \phi(x)\tilde{\pi}_x(a)$
for all $\phi\in C_0(X)$ and $a\in A$ as long as $x\notin N$.  Thus,
for $x\not\in N$, we get $I_x\subset \ker\tilde{\pi}_x$ and we can
assume that this is true all of the time by setting $\tilde{\pi}_x=0$
on $N$.  Furthermore, since we have only changed
$\{\tilde{\pi}_x(a)\}$ 
on a null set we still have $\pi(a) = \int^\oplus_X
\tilde{\pi}_x(a)d\mu(x)$ for all $a\in A$.  Now let $\pi_x$ be the
factorization of $\tilde{\pi}_x$ to $A(x)$ for all $x\in X$.  Then,
since $\pi_x(a(x)) = \tilde{\pi}_x(a)$ for all $x$, 
$\{\pi_x(a(x))\}$ is clearly a bounded Borel
field of operators for all $a\in A$ and we have
\[
\pi(a) = \int^\oplus_X \tilde{\pi}_x(a)d\mu(x) = \int^\oplus_X
\pi_x(a(x))d\mu(x).
\]
This yields the required decomposition of $\pi$ to the fibres $A(x)$.
Furthermore, since the factorization of a nondegenerate representation
is nondegenerate, the $\pi_x$ are nondegenerate almost everywhere.  

Finally, suppose $\rho_x$ is also a decomposition of $\pi$ as in the
statement of the proposition.  Then it follows from Lemma \ref{lem:12}
that for each $a\in A$ there exists a $\mu$-null set $N_a$ such that
for all $x\not\in N_a$ we have $\pi_x(a) = \rho_x(a)$.  Let $a_i$ be a
countable dense set in $A$ and let $N = \bigcup N_{a_i}$.  Then $N$ is
still $\mu$-null and it follows from the fact that representations of
$C^*$-algebras are continuous that $\pi_x(a) = \rho_x(a)$ for all
$a\in A$ and all $x\not\in N$.  
\end{proof}

It is a useful fact that, at least in the separable case, every
representation of a $C_0(X)$-algebra is equivalent to a
$C_0(X)$-linear one.  The following theorem is a restriction of
\cite[Theorem 8.3.2]{dixmiercstar}.

\begin{prop}
\label{prop:56}
Suppose $X$ is a second countable locally compact Hausdorff space,
$A$ is a separable $C_0(X)$-algebra and $\rho$ is a separable
representation of $A$ on $\mcal{H}_\rho$.  Then there is an analytic
Borel Hilbert bundle $X*\mfrk{H}$ and a finite measure $\mu$ on $X$ such
that $\rho$ is unitarily equivalent to a $C_0(X)$-linear
representation of $A$ on $L^2(X*\mfrk{H},\mu)$.
\end{prop}

\begin{proof}
Suppose the $C_0(X)$-action on $A$ is given by $\Phi_A$.  
Given $\rho$ as above we can extend $\rho$ to the multiplier algebra
and obtain a representation $\rho'=\overline{\rho}\circ\Phi_A$ of
$C_0(X)$ on $\mcal{H}_\rho$ which is nondegenerate because $\Phi_A$
and $\rho$ are.  It is a deep result (\cite[Theorem
E.14]{tfb2},\cite[Pages 54--55]{invitation}) that $\rho'$ is unitarily
equivalent to a representation of the form 
\[
\tilde{\rho} = (\rho_{\mu_\infty}\otimes 1_{\mcal{H}_\infty})\oplus \rho_{\mu_1} \oplus (\rho_{\mu_2}\otimes
1_{\mcal{H}_2}) \oplus \cdots 
\]
on 
\[
(L^2(X_\infty,\mu_\infty)\otimes
\mcal{H}_\infty) \oplus L^2(X_1,\mu_1) \oplus (L^2(X_2,\mu_2)\otimes
\mcal{H}_2)\oplus \cdots
\]
where each $\mu_n$ is a finite Borel measure on $X$ with $\mu_n$
disjoint from $\mu_m$ if $n\ne m$, $\rho_{\mu_n}$ is the
representation of $C_0(X)$ on $L^2(X,\mu_n)$ given by
$\rho_{\mu_n}(\phi)h(x) = \phi(x)h(x)$ for all $\phi\in C_0(X)$ and $h\in
L^2(X,\mu_n)$, and $\mcal{H}_n$ is a fixed Hilbert space of dimension
$0\leq n\leq \infty$.  Since there is no harm in replacing $\mu_n$ by
a scalar multiple of $\mu_n$ we can assume without loss of generality
that $\mu_n(X) \leq 1/2^n$.  Let $X_n=\supp \mu_n$ and observe that $X_n$ is
a Borel partition of $X$. (If we are missing any bits just throw them
into $X_1$ and give them zero measure.)  Let $X*\mfrk{H}$ be the
disjoint union $\coprod_{n=0}^{n=\infty} X_n\times\mcal{H}_n$.  It
is easy to see that $X*\mfrk{H}$ is an analytic Borel Hilbert bundle 
\cite[Example F.5]{tfb2}.  
Now let $\mu = \sum_n \mu_n$ and note that $\mu$ is a
finite Borel measure on $X$.  It is straightforward to see
\cite[Corollary F.12]{tfb2} that 
\[
L^2(X*\mfrk{H},\mu) \cong (L^2(X_\infty,\mu_\infty)\otimes
\mcal{H}_\infty) \oplus L^2(X_1,\mu_1) \oplus (L^2(X_2,\mu_2)\otimes
\mcal{H}_2)\oplus \cdots 
\]
and that the isomorphism is given by sending a section in
$L^2(X*\mfrk{H},\mu)$ to the sum of the appropriate restrictions to
each subfactor.  However, once we untangle all of the definitions it
follows that this isomorphism intertwines $\tilde{\rho}$ with the
representation $\rho_\mu$ of $C_0(X)$ on $L^2(X*\mfrk{H},\mu)$ given
by $\rho_\mu(\phi) = T_\phi$ for all $\phi\in C_0(X)$ where $T_\phi$
is the diagonal operator associated to $\phi$.  (Notice that
$\rho_{\mu_n}(\phi)$ is exactly the ``diagonal'' operator associated
to $\phi$ in $L^2(X,\mu_n)$.)

Next, let $U:\mcal{H}_\rho \rightarrow L^2(X*\mfrk{H},\mu)$ be the
unitary intertwining $\rho'$ and $\rho_\mu$ and let $\pi(a) = U
\rho(a) U^*$.  Then $\pi$ is unitarily equivalent to $\rho$ by
construction and we can compute that 
\[
\pi(\Phi_A(\phi)a) = U\overline{\rho}(\Phi_A(\phi)) U^* U\rho(a)U^*
= U\rho'(\phi)U^* \pi(a) = \rho_\mu(\phi)\pi(a) = T_\phi \pi(a).
\]
Thus $\pi$ is $C_0(X)$-linear.  
\end{proof}

\begin{remark}
\index{direct integral}
Of course, combined with Proposition \ref{prop:53}, this implies that
{\em every} separable representation of a separable $C_0(X)$-algebra
is unitarily equivalent to a direct integral of representations of the
fibres.
\end{remark}

At this point we are finally ready to define what it means to be a
covariant representation of a groupoid crossed product.  

\begin{definition}
\label{def:42}
\index{covariant representation}
Suppose $(A,G,\alpha)$ is a separable groupoid dynamical system.  A
{\em covariant representation} $(\mu,G\unit*\mfrk{H},\pi,U)$ of $(A,G,\alpha)$
consists of a unitary representation $(\mu,G\unit *\mfrk{H},U)$ of $G$
and a $C_0(G\unit)$-linear representation $\pi:A\rightarrow
B(L^2(X*\mfrk{H},\mu))$.  Furthermore, if $\{\pi_x\}$ is a
decomposition of $\pi$ and $\nu$ is the measure induced by $\mu$ we
require that there exists a $\nu$-null set $N$ such that for all
$\gamma\not\in N$ we have 
\begin{equation}
\label{eq:cov}
U_\gamma \pi_{s(\gamma)}(a) =
\pi_{r(\gamma)}(\alpha_\gamma(a))U_\gamma\quad\text{for all $a\in
  A(s(\gamma))$.}
\end{equation}
\end{definition}

\begin{remark}
\index[not]{$(\pi,U)$}
To conserve notation we will sometimes denote a covariant
representation by $(\pi,U)$ and will understand that the groupoid
representation includes a quasi-invariant measure and a Borel Hilbert
bundle.  
\end{remark}

We end this section with an example of a very important class of
representations.  Unfortunately, the example is fairly technical and
requires technology we have not introduced here.  The following
attempts to describe a coherent construction but certainly does not
include all of the relevant details.\footnote{Many thanks to Jon Brown
  and Dana Williams\index{Dana Williams} for allowing me to use their notes as a
  reference.}  One of the things
to take away from this example is that dealing with covariant
representations directly is messy.  We will develop tools in the next
section which will allow us to get around these difficulties.  

\begin{example}
\label{ex:19}
\index{left regular representation}
Suppose $(A,G,\alpha)$ is a separable groupoid dynamical system and
$\pi$ is a separable representation of $A$.  Using Proposition
\ref{prop:56} we can assume without loss of generality that there exists
a Borel Hilbert bundle $X*\mfrk{H}$ and a finite
measure $\mu$ such that $\pi$ is a $C_0(G\unit)$-linear representation
of $A$ on $L^2(X*\mfrk{H},\mu)$.  As such there is a decomposition 
\[
\pi = \int_{G\unit}^\oplus \pi_u d\mu(u)
\]
coming from Proposition \ref{prop:53}.  Next, let $[\mu]$ be the
saturation of $\mu$ and $\nu=\int_{G\unit}\lambda^u d\mu$ 
be the measure induced on $G$ by
$\mu$.  We use \cite[Theorem I.5]{tfb2} to decompose $\nu$ with
respect to the source map to get measures $\nu_u$ on $G_u$ such that 
\[
\nu(f) = \int_{G\unit}\int_G f(x)d\nu_u(x)d[\mu](u).
\]
We denote the image of $\nu_u$ under inversion as $\nu^u$.  Observe
that $\nu^u$ is a Radon measure on $G^u$.  
Now, use Example \ref{ex:35} to form the pull back bundle
$s^*(G\unit*\mfrk{H})$.  We define 
\[
\mcal{K}(u) = L^2(s^*(G\unit*\mfrk{H})|_{G^u},\nu^u)
\]
for all $u\in G\unit$.  In other words, $\mcal{K}(u)$ is the
$L^2$-space of maps from $G^u$ into $G\unit*\mfrk{H}$ such that 
$h(\gamma)\in \mcal{H}(s(\gamma))$
for all $\gamma\in G^u$.   Since each
$\nu^u$ is finite, this is a collection of separable
Hilbert spaces which we then form into the bundle $G\unit*\mfrk{K}$.  Now,
let $f_n$ be a fundamental sequence for $s^*(G\unit*\mfrk{H})$ and
$\phi_m$ a sequence of point separating functions in $C_c(G)$.  We then
define a fundamental sequence on $G\unit*\mfrk{K}$ by 
\[
g_{n,m}(u)(\gamma):= \phi_m(\gamma)f_n(\gamma).
\]
It is straightforward to show that Proposition \ref{prop:113}
implies that we can use the $g_{n,m}$ to make
$G\unit*\mfrk{K}$ into an analytic Borel Hilbert bundle. 

Let $\Delta$ be the modular function coming from $[\mu]$.  
Now, it takes a lot of work, is not at all obvious, and currently
there is no decent reference,\footnote{There is a proof of this fact
  in the personal notes of Dana Williams.\index{Dana Williams}} but we can
modify the $\nu^u$ on a null set so that 
$\gamma\cdot\nu^{s(\gamma)} = \Delta(\gamma)\nu^{r(\gamma)}$ for all
$\gamma\in G$.  
This allows us to define a unitary 
\[
L_\gamma:
L^2(s^*(G\unit*\mfrk{H})|_{G^{s(\gamma)}},\nu^{s(\gamma)}) 
\rightarrow
L^2(s^*(G\unit*\mfrk{H})|_{G^{r(\gamma)}},\nu^{r(\gamma)})
\]
by $L_\gamma h(\eta) = \Delta(\gamma)\poshalf h(\gamma\inv\eta)$ for
all $\eta\in G^{r(\gamma)}$.  This is just a souped up version of the left regular
representation from Example \ref{ex:18} and it is not hard to see that 
$L:G\rightarrow \isom(G\unit*\mfrk{K})$ such that
$L(\gamma)=(r(\gamma),L_\gamma,s(\gamma))$ defines a representation of
$G$.  Finally, we define $\tilde{\pi}:A\rightarrow
L^2(G\unit*\mfrk{K},[\mu])$ by
\[
(\tilde{\pi}(a)\xi(u))(\gamma) =
\pi_{s(\gamma)}(\alpha_\gamma\inv(a(u)))\xi(u)(\gamma)
\]
for $a\in A$ and $\xi\in L^2(G\unit*\mfrk{K},[\mu])$.  
It takes some work but one can show that $\tilde{\pi}$ is a
$C_0(G\unit)$-linear representation of $A$ and that its decomposition is
given by 
\[
\tilde{\pi}_u(a)h(\gamma) =
\pi_{s(\gamma)}(\alpha_\gamma\inv(a))h(\gamma).
\]
for $a\in A(u)$ and $h\in \mcal{K}(u)$.  Finally, given $\gamma\in G$, 
we can observe that 
\begin{align*}
(L_\gamma \tilde{\pi}_{s(\gamma)}(a)h)(\eta) &=
\Delta(\gamma)\poshalf \tilde{\pi}_{s(\gamma)}(a)h(\gamma\inv\eta) \\
&= \Delta(\gamma)\poshalf 
   \pi_{s(\eta)}(\alpha_{\gamma\inv\eta}\inv(a))h(\gamma\inv\eta) \\
&=  \pi_{s(\eta)}(\alpha_\eta\inv(\alpha_\gamma(a)))L_\gamma h(\eta) \\
&= (\tilde{\pi}_{r(\gamma)}(\alpha_\gamma(a))L_\gamma h)(\eta)
\end{align*}
for $a\in A(s(\gamma))$, $h\in \mcal{K}(s(\gamma))$ and $\eta\in
G^{r(\gamma)}$.  Thus $([\mu],G\unit*\mfrk{K},\tilde{\pi},L)$ is a
covariant representation of $(A,G,\alpha)$ called the {\em left regular
  representation}.  
\end{example}


\section{The Groupoid Crossed Product}
\label{sec:crossedprod}

We start with a couple of useful propositions that are mildly
interesting from a measure theoretic point of view.  The following proofs
were communicated to me by Dana Williams.\footnote{Thanks
  Dana!}\index{Dana Williams}

\begin{definition}
Suppose that $X$ and $Y$ are locally compact
Hausdorff spaces.  A family $\{\rho^y\}_{y\in Y}$ of Radon measures on
$X$ is called a {\em Radon family} of measures on $X$ if 
\begin{equation}
\label{eq:59}
y\mapsto \int_X f(y)d\rho^y(x)
\end{equation}
is a bounded Borel function supported on a compact set for all $f\in C_c(X)$. 
\end{definition}

\begin{example}
The most common example of a Radon family of measures will be a Haar
system on a groupoid $G$.
\end{example}

\begin{remark}
\label{rem:16}
Suppose $X$ and $Y$ are locally compact Hausdorff spaces and
$\{\rho^y\}$ is a Radon family of measures on $X$.  Then if $\mu$ is
any Radon measure on $Y$ setting
\begin{equation}
\label{eq:60}
\nu(f) := \int_Y\int_X f(x) d\rho^y(x)d\mu(y)
\end{equation}
for all $f\in C_c(X)$ 
defines a radon measure $\nu$ on $X$ called the {\em induced
  measure}.  Our goal will be to see that \eqref{eq:60} extends to a
much wider class of functions.  
\end{remark}

\begin{prop}
\label{prop:57}
Suppose $X$ and $Y$ are second countable locally compact Hausdorff
spaces and $\{\rho^y\}$ is a Radon family of measures on $X$.  
Let $\mu$ be a Radon measure on $Y$ and let
$\nu = \int \rho^y \mu(u)$ be the induced measure as in Remark
\ref{rem:16}.  
Given a positive Borel
function $f$ we define $\rho(f)$ on $Y$ by 
\begin{equation}
\label{eq:62}
\rho(f)(y) := \int_X f(x)d\rho^y(x)
\end{equation}
Then $\rho(f)$ is a positive extended real valued Borel
function on $Y$ and we have
\begin{equation}
\label{eq:41}
\int_X f(x)d\nu(x) = \int_Y \rho(f)(y)d\mu(y).
\end{equation}

Furthermore if $f$ is a $\nu$-integrable Borel function then we define
$\rho(f)$ on $Y$ by \eqref{eq:62} whenever $f$ is
$\rho^y$-integrable and $\rho(f)=0$ otherwise.  Then in this case $\rho(f)$
is a Borel function on $Y$ and \eqref{eq:41} still holds.  
\end{prop}

\begin{remark}
When convenient we will use \eqref{eq:60} instead of the more 
formal \eqref{eq:41}.  
\end{remark}

\begin{remark}
We will use the following notation.  Given a locally compact Hausdorff
space $X$ we will let $\mcal{B}(X)$ denote the set of Borel functions
on $X$ and $\mcal{B}_c^b(X)$ denote the set of bounded Borel
functions which vanish off a compact set.  We will let $\mcal{B}^+(X)$
denote the positive Borel functions and $\mcal{B}^{+,e}(X)$ denote the
set of positive extended real valued Borel functions.  
\end{remark}

We will prove Proposition \ref{prop:57} via a series of lemmas. 

\begin{lemma}
\label{lem:15}
If $f\in \mcal{B}_c^b(X)$ then $\rho(f)\in \mcal{B}_c^b(Y)$.  If $f\in
\mcal{B}^+(X)$ then $\rho(f)\in \mcal{B}^{+,e}(Y)$.  
\end{lemma}

\begin{proof}
Note that $\rho(f)$ is defined if $f\in\mcal{B}_c^b(X)$ or $f\in
\mcal{B}^+(X)$, although in the latter case the function may take
infinite values.  If $f\in \mcal{B}_c^b(X)$ then there is a $g\in
C_c^+(X)$ such that $|f|\leq g$.  Then 
\[
|\rho(f)|\leq \rho(|f|)\leq \rho(g).
\]
Since $\{\rho^y\}$ is a Radon family, $\rho(g)\in \mcal{B}_c^b(Y)$,
and it follows that $\rho(f)$ is bounded and must vanish off a compact
set.  Thus if $f\in \mcal{B}_c^b(X)$ then $\rho(f)\in
\mcal{B}_c^b(Y)$ exactly when $\rho(f)$ is Borel. 

Now let $U\subset X$ be a relatively compact open set.  Let 
\[
\mcal{A} := \{f\in \mcal{B}^b(U):\text{$\rho(f)$ is Borel}\}.
\]
Then $C_c(U)\subset\mcal{A}$.  The Dominated Convergence Theorem
implies that $\mcal{A}$ is closed under monotone sequential limits.
Then \cite[Proposition 6.2.9]{analysisnow} implies that $\mcal{A} =
\mcal{B}^b(U) \subset\mcal{B}_c^b(X)$.  Since this holds for any
relatively compact open set we see that $\rho(f)\in\mcal{B}_c^b(Y)$
for any $f\in \mcal{B}_c^b(X)$.  If $f\in \mcal{B}^+(X)$ then, since
$X$ is $\sigma$-compact, we can find $f_n\in
\mcal{B}_c^b(X)\cap\mcal{B}^+(X)$ such that $f_n\nearrow f$.  We then
see that $\rho(f)\in \mcal{B}^{+,e}(Y)$ by the Monotone Convergence
Theorem.  
\end{proof}

\begin{lemma}
\label{lem:16}
We get a measure $\bar{\nu}$ on $X$ via 
\[
\bar{\nu}(E) := \int_Y \int_X \chi_E(x)d\rho^y(x)d\mu(y)
\]
for any Borel set $E\subset X$. 
\end{lemma}

\begin{proof}
Since $\rho(\chi_E)$ is in $\mcal{B}^{+,e}(Y)$ we know $\bar{\nu}$ is
defined on all Borel sets $E$.  Clearly $\bar{\nu}(\emptyset) = 0$.
To see that $\bar{\nu}$ is countably additive requires only a couple
applications of the Monotone Convergence Theorem. 
\end{proof}

\begin{lemma}
\label{lem:17}
For all $f\in\mcal{B}^+(X)$
\begin{equation}
\label{eq:63}
\bar{\nu}(f) := \int_X f(x)d\bar{\nu}(x) = \int_Y\int_X f(x) d\rho^y(x)
d\mu(y). 
\end{equation}
In particular, $\bar{\nu}$ is a Radon measure on $X$. 
\end{lemma}

\begin{proof}
By linearity, \eqref{eq:63} holds for all nonnegative simple
functions. But if $f\in \mcal{B}^+(X)$ then there are nonnegative
simple functions $f_n\nearrow f$.  Thus \eqref{eq:63} holds for all
$f\in \mcal{B}^+(X)$ by the Monotone Convergence Theorem.  

Since $X$ is second countable, to see that $\bar{\nu}$ is a Radon
measure just requires that we demonstrate that $\bar{\nu}(K)< \infty$ for all
compact sets $K\subset X$ \cite[Theorem 2.18]{daddyrudin}.  
But we can find $g\in C_c^+(X)$ such that
$g(x) = 1$ for all $x\in K$.  Then $\bar{\nu}(K)\leq \bar{\nu}(g)<
\infty$.  
\end{proof}

\begin{proof}[Proof of Proposition \ref{prop:57}]
If $f\in\mcal{B}^+(X)$ then the assertions about $\rho(f)$ are taken
care of by Lemma \ref{lem:15}.  Furthermore, $\nu$ and
$\bar{\nu}$ are Radon measures that agree on $C_c(X)$.  Thus, $\nu
= \bar{\nu}$ which implies that \eqref{eq:60} holds as required. 

Now suppose $f$ is $\nu$-integrable.  Then $f$ is
$\bar{\nu}$-integrable and we can decompose $f$ as the sum of four
positive $\bar{\nu}$-integrable functions $f = f_1 -f_2 +i f_3 -i
f_4$.  However it now follows that, for each $i$, $\rho(f_i)(y)<\infty$
for $\mu$-almost all $y$.  Therefore, $f$ is $\rho^y$ integrable
and $\rho(f)$ is defined by
\eqref{eq:62} $\mu$-almost everywhere.  The fact that \eqref{eq:60}
holds for $f$ now follows from the fact that it holds for each $f_i$
and by linearity. 
\end{proof}

The following result is contained in  \cite[Lemma 5.2]{ramsay} or 
\cite[Lemma 4.9]{coords}.  However, both of those references make use
of measured groupoids, so a measured groupoid free proof is provided
here for convenience.  This result is particularly useful when dealing with the
fact that the covariance relation for covariant representations only
holds almost everywhere.  

\begin{definition}
\label{def:43}
Suppose $G$ is a groupoid and $N$ is a $G$-invariant subset of
$G\unit$.  Then the {\em restriction} of $G$ to $N$ is defined to be
\[
G|_N := r\inv(N) = \{\gamma\in G : r(\gamma),s(\gamma)\in N\}.
\]
\end{definition}

\begin{lemma}
\label{lem:13}
Suppose $G$ is a second countable 
locally compact Hausdorff groupoid with Haar system
$\lambda$ and $\mu$ is a finite quasi-invariant measure on $G\unit$.  Let
$\nu = \int_{G\unit} \lambda^u d\mu(u)$ and suppose $\Sigma$ is a $\nu$-conull
subset of $G$.  If $\Sigma$ is closed under multiplication then there
exists a $G$-invariant $\mu$-conull set $N\subset G\unit$ such that
$G|_N$ is conull and $G|_N \subset \Sigma$.  
\end{lemma}

\begin{proof}
Set $\Sigma_1 = \Sigma\cap\Sigma\inv = \{\gamma\in G :
\gamma,\gamma\inv\in \Sigma\}$.  Then $\Sigma_1$ is a groupoid
contained in $G$ (whose unit space may be smaller than $G\unit$) and
contains a $\nu$-conull Borel subset $B$ of $G$.  It follows from
Proposition \ref{prop:57} applied to the characteristic function of
$G\setminus B$ that if $B$ is conull
with respect to $\nu$ then there must exist a $\mu$-conull set $U$ such
that $\lambda^u(G^u\setminus B) = 0$ for all $u\in U$.  Now, let $N =
r(r\inv(U)\cap s\inv(U))$.  Clearly $N$ is a Borel set.  It is
straightforward to show that $r\inv(U)$ is $\nu$-conull and $s\inv(U)$
is $\nu\inv$-conull. However, $\nu$ is equivalent to $\nu\inv$ so that
$r\inv(U)\cap s\inv(U)$ is $\nu$-conull.  It then follows relatively
quickly that $N$ is
$\mu$-conull and $G$-invariant.  Furthermore, the
previous argument also shows that $G|_N$ is conull.  We would like to see
that $G|_N\subset \Sigma_1$.  Suppose $\gamma\in G|_N$ and let $u =
r(\gamma)$ and $v=s(\gamma)$.  Then $u,v\in U$ so that $G^u\cap B$ and
$G^v\cap B$ are $\lambda^u$-conull and $\lambda^v$-conull, respectively.
However, $\lambda$ is invariant so $\gamma\cdot(G^v\cap B)$ is
$\lambda^u$-conull.  It follows that the intersection $G^u\cap B\cap
\gamma\cdot (G^v\cap B)$ is $\lambda^u$-conull and therefore nonempty.  
Thus there
exists $\eta\in B\subset\Sigma^1$  such that $r(\eta) = s(\gamma)$ and
$\gamma\eta\in B\subset\Sigma^1$.  But $\Sigma^1$ is a subgroupoid so
therefore $\gamma\in \Sigma^1$.  
\end{proof}

This leads us, as promised, to the following proposition.
This will be extremely useful when dealing with covariant
representations because it will allow us to restrict our attention to
a conull subgroupoid (what some people call an inessential
contraction) on which the covariance relation \eqref{eq:cov} holds.  

\begin{prop}
\label{prop:55}
Suppose $(A,G,\alpha)$ is a separable groupoid dynamical system and
$(\mu,G\unit*\mfrk{H},\pi,U)$ is a covariant representation of
$(A,G,\alpha)$.  Then there exists a $G$-invariant 
$\mu$-conull set $N\subset G\unit$
such that \eqref{eq:cov} holds on all of the conull subgroupoid 
$G|_N$.
\end{prop}

\begin{proof}
Let $\Sigma$ be the set of $\gamma\in G$ such that \eqref{eq:cov}
holds and let $\nu=\int_{G\unit} \lambda^u d\mu(u)$.  
By definition $\Sigma$ is $\nu$-conull.  Since $U$ and $\alpha$ are
both homomorphisms, it is not hard to see that $\Sigma$ is closed under
multiplication.  However the result then follows from Lemma
\ref{lem:13}.
\end{proof}

The following construction is a key aspect of crossed product
theory.  It is (mostly) developed in \cite[Section 7]{renaultequiv} but
is so fundamental that it has been reproduced here.  

\begin{prop}
\label{prop:59}
\index{covariant representation!integrated form}
\index[not]{$\pi\rtimes U$}
Suppose $(A,G,\alpha)$ is a separable dynamical system and that
$(\mu,{G\unit*\mfrk{H}},\pi,U)$ is a covariant representation.   Let
$\pi = \int^\oplus_{G\unit}\pi_u d\mu$ be a decomposition of $\pi$.  Then
there is an $I$-norm decreasing, nondegenerate, $*$-representation $\pi\rtimes
U$ of $\Gamma_c(G,r^*\mcal{A})$ on $L^2(G\unit*\mfrk{H},\mu)$ given by 
\begin{equation}
\label{eq:36}
\pi\rtimes U(f)h(u) = \int_G \pi_u(f(\gamma))U_\gamma
h(s(\gamma))\Delta(\gamma)\neghalf d\lambda^u(\gamma)
\end{equation}
for $f\in \Gamma_c(G,r^*\mcal{A})$, $h\in \mcal{L}^2(G\unit*\mfrk{H},\mu)$, 
and $u\in G\unit$
where $\Delta$ is the modular function from Definition \ref{def:39}.  

Furthermore, given $h,k\in \mcal{L}^2(G\unit*\mfrk{H},\mu)$ and $f\in
\Gamma_c(G,r^*\mcal{A})$ we have
\begin{equation}
\label{eq:38}
(\pi\rtimes U(f)h,k) = \int_G (\pi_{r(\gamma)}(f(\gamma))U_\gamma h(s(\gamma)) ,
k(r(\gamma))) \Delta(\gamma)\neghalf d\nu(\gamma)
\end{equation}
where $\nu$ is the measure induced by $\mu$.  
\end{prop}

\begin{proof}
Using the quasi-invariance of $\mu$ and the Cauchy-Schwartz inequality
in $\mcal{H}(x)$ and $L^2(G,\nu)$ we have, for all $f\in
\Gamma_c(G,r^*\mcal{A})$ and $h,k\in
\mcal{L}^2(G\unit*\mfrk{H},\mu)$, 
\begin{align*}
\int_G& |(\pi_{r(\gamma)}(f(\gamma))U_\gamma
h(s(\gamma)), k(r(\gamma)))| \Delta(\gamma)\neghalf d\nu(\gamma)  \\
\leq& \int_G \|\pi_{r(\gamma)}(f(\gamma))\|\|U_\gamma
h(s(\gamma))\|\|k(r(\gamma))\|\Delta(\gamma)\neghalf d\nu(\gamma) \\
\leq& \int_G
\|f(\gamma)\|\|h(s(\gamma))\|\|k(r(\gamma))\|\Delta(\gamma)\neghalf
d\nu(\gamma) \\
\leq& \left( \int_G \|f(\gamma)\|\|h(s(\gamma))\|^2\Delta(\gamma)\inv
  d\nu(\gamma)\right)^{1/2} \\
&\ \cdot\left(\int_G \|f(\gamma)\|\|k(r(\gamma))\|^2 d\nu(\gamma)\right)^{1/2}
\\
\leq&\left( \int_{G\unit}\int_G \|f(\gamma)\|d\lambda_u(\gamma)
  \|h(u)\|^2 d\mu(u)\right)^{1/2} \\
&\ \cdot \left(\int_{G\unit}\int_G \|f(\gamma)\|d\lambda^u(\gamma) \|k(u)\|^2
  d\mu(u)\right)^{1/2} \\
\leq&\left(
  \|f\|_I\|h\|^2\right)^{1/2}\left(\|f\|_I\|k\|^2\right)^{1/2} =
\|f\|_I\|h\|\|k\|.
\end{align*}
Thus $\gamma\mapsto
(\pi_{r(\gamma)}(f(\gamma))U_\gamma h(s(\gamma)),k(r(\gamma)))\Delta(\gamma)\neghalf$
is an integrable function, and, using elementary tensors, it is 
straightforward to see that it is
Borel.  Equation \eqref{eq:38} now follows quickly from 
Proposition \ref{prop:57}. Let $R = \pi\rtimes U$.  
It then follows from the above
calculation that 
\[
|(R(f)h|k)| \leq \int_G |(\pi_{r(\gamma)}(f(\gamma))U_\gamma
h(s(\gamma)), k(r(\gamma)))| \Delta(\gamma)\neghalf d\nu(\gamma) 
\leq \|f\|_I \|h\| \|k\|.
\]
This is enough to show that $R$ is $I$-norm decreasing. 

It is clear that $R$ is linear.  We show now 
that $R$ is multiplicative.  Using Proposition
\ref{prop:55} we can assume that there is a $\mu$-null set $N$ such
that the convolution identity holds on all of $G|_N$.  If $f,g\in
\Gamma_c(G,r^*\mcal{A})$, $h\in \mcal{L}^2(G\unit*\mfrk{H},\mu)$ and 
$u\in N$ we have
\begin{align*}
R(f*g)h(u) &= \int_G \pi_u(f*g(\gamma))U_\gamma
h(s(\gamma))\Delta(\gamma)\neghalf d\lambda^u(\gamma) \\
&= \int_G \int_G \pi_u(f(\eta)\alpha_\eta(g(\eta\inv\gamma)))U_\gamma
h(s(\gamma)) \Delta(\gamma)\neghalf d\lambda^u(\eta)\lambda^u(\gamma)
\\
&= \int_G \pi_u(f(\eta))\int_G\pi_u(\alpha_\eta(g(\eta\inv\gamma)))
U_\gamma h(s(\gamma))\Delta(\gamma)\neghalf
d\lambda^u(\gamma)d\lambda^u(\eta) \\
&= \int_G
\pi_u(f(\eta))\int_G\pi_u(\alpha_\eta(g(\gamma)))U_{\eta\gamma}h(s(\gamma))
\Delta(\eta\gamma)\neghalf d\lambda^{s(\eta)}(\gamma)d\lambda^u(\eta).
\end{align*}
Recall that Theorem \ref{thm:quasi} says that we can choose $\Delta$
to be a homomorphism.  Using the fact that $U$ is a homomorphism, as
well as the covariance relation, we get
\begin{align*}
R(f*g)h(u) &= \int_G \pi_u(f(\eta))U_\eta\int_G \pi_u(g(\gamma))
  U_\gamma
h(s(\gamma))\Delta(\gamma)\neghalf d\lambda^{s(\eta)}(\gamma)
\Delta(\eta)\neghalf d\lambda^u(\eta) \\
&= \int_G \pi_u(f(\eta))U_\eta R(g)h(s(\eta))\Delta(\eta)\neghalf
d\lambda^u(\eta) \\
&= R(f)R(g)h(u).
\end{align*}

Next, we will prove that $R$ preserves the involution.  It is
straightforward to show that
$\Delta(\gamma)\neghalf d\nu(\gamma)$ is invariant under
inversion.  Using this fact we have, for appropriate $f$,$h$ and $k$,
\begin{align*}
(R(f^*)h,k) &= \int_G (\pi_{r(\gamma)}(f^*(\gamma))U_\gamma
h(s(\gamma)), k(r(\gamma))) \Delta(\gamma)\neghalf d\nu(\gamma) \\
&= \int_G (\pi_{r(\gamma)}(\alpha_\gamma(f(\gamma\inv)^*)) U_\gamma
h(s(\gamma)), k(r(\gamma))) \Delta(\gamma)\neghalf d\nu(\gamma) \\
&= \int_G (\pi_{s(\gamma)}(\alpha_\gamma\inv(f(\gamma)))^* U_\gamma^*
h(r(\gamma)) , k(s(\gamma))) \Delta(\gamma)\neghalf d\nu(\gamma) \\
&= \int_G (h(r(\gamma)) , U_\gamma
\pi_{s(\gamma)}(\alpha_\gamma\inv(f(\gamma))) k(s(\gamma)))
\Delta(\gamma)\neghalf d\nu(\gamma) \\
&= \int_G (h(r(\gamma)) , \pi_{r(\gamma)}(f(\gamma)) U_\gamma
k(s(\gamma))) \Delta(\gamma)\neghalf d\nu(\gamma) \\
&= (h,R(f)k).
\end{align*}
It follows that $R$ is a $*$-homomorphism. 

The last thing we have to prove is that $R$ is nondegenerate.  Suppose
$(R(f)h,k) = 0$ for all $f\in \Gamma_c(G,r^*\mcal{A})$ and $h\in
L^2(X*\mfrk{H},\mu)$.  Let $e_l$ be a special orthogonal fundamental
sequence in $L^2(X*\mfrk{H},\mu)$ and $a_i$ a countable dense set in
$A$.  Since $G$ is second countable it is $\sigma$-compact and we can
find a countable set of functions $\phi_n\in C_c(G)^+$ such that for all
$\gamma\in G$ there exists $\phi_n$ such that $\phi_n(\gamma) = 1$.  Now
let $\mcal{H}_0$ be the countable set of rational linear combinations
of $e_l$ and observe that 
\[
0=(R(\phi_n\otimes a_i)h,k) = \int_G \phi_n(\gamma)
(\pi_{r(\gamma)}(a_i(r(\gamma))) U_\gamma h(s(\gamma)) ,
k(r(\gamma)))\Delta(\gamma)\neghalf d\nu(\gamma)
\]
for all $n,i$ and $h\in \mcal{H}_0$.  Since this is a countable family
there exists a $\nu$-null set $N$ such that, for all  $\gamma\not\in
N$, we have 
\[
\phi_n(\gamma)(\pi_{r(\gamma)}(a_i(r(\gamma)))U_\gamma h(s(\gamma)) ,
k(r(\gamma)))
\]
for all $n,i$ and $h\in \mcal{H}_0$.  Now $0 = \nu(N) =
\int_{G\unit}\lambda^u(N) d\mu(u)$ so that there exists a $\mu$-null
set $M$ such that $\lambda^u(N) = 0$ for all  $u\not\in M$.
Furthermore, by making $M$ a little larger if necessary, we can assume
that $\pi_u$ is nondegenerate for all $u\not\in M$.  Then given
$u\not\in M$ choose some  $\gamma\in G^u$ such that $\gamma\not\in N$
and $n$ such that $\phi_n(\gamma)=1$.  Then 
\begin{equation}
\label{eq:44}
(\pi_u(a_i(u))U_\gamma h(s(\gamma)) , k(u)) = 0
\end{equation}
for all $i$ and $h\in \mcal{H}_0$.  Now, observing that
$\{e_l(s(\gamma))\}$ forms an orthogonal basis for
$\mcal{H}(s(\gamma))$ (plus zero vectors) and that $U_\gamma$ is a
unitary, it follows that $\{U_\gamma h(s(\gamma))\}_{h\in\mcal{H}_0}$
is dense in $\mcal{H}(r(\gamma))$.  Furthermore, $\{a_i(u)\}$ is dense
in $A(u)$ and $\pi_u$ is nondegenerate so $\{\pi_u(a_i(u))U_\gamma
h(s(\gamma))\}$ is dense in $\mcal{H}(r(\gamma))$.  It now follows from
\eqref{eq:44} that $k(u) = 0$.  Since this is true  $\mu$-almost
everywhere, we have $k=0$ and that $R$ is nondegenerate. 
\end{proof}

This brings us to the most important tool in the theory of
groupoid crossed products.  
The following theorem is a generalization of Renault's ``Proposition 4.2''
which we will discuss in Section \ref{sec:scalarprod}.  It is stated
and proved in \cite{renaultequiv}.  The reason it's so important is
that it allows us to avoid dealing with covariant representations directly.

\begin{theorem}[Renault's Disintegration Theorem 
{\cite[Theorem 7.12]{renaultequiv}}]
\index{disintegration theorem}
\label{thm:disintigration}
Suppose that $\mcal{H}_0$ is a dense subspace of a complex Hilbert
space $\mcal{H}$ and that $\pi$ is a homomorphism from
$\Gamma_c(G,r^*\mcal{A})$ into the algebra of linear operators on
$\mcal{H}_0$ such that 
\begin{enumerate}
\item $\spn\{\pi(f)h:f\in\Gamma_c(G,r^*\mcal{A}), h\in \mcal{H}_0\}$ is dense in
  $\mcal{H}$,
\item for each $h,k\in\mcal{H}_0$, 
\[
f\mapsto (\pi(f)h,k)
\]
is continuous in the inductive limit topology, and 
\item for each $f\in\Gamma_c(G,r^*\mcal{A})$ and all $h,k\in
  \mcal{H}_0$
\[
(\pi(f)h,k) = (h,\pi(f^*)k).
\]
\end{enumerate}
Then each $\pi(f)$ is bounded and extends to a  bounded operator
$\Pi(f)$ on $\mcal{H}$ such that $\Pi$ is an $I$-norm decreasing
$*$-representation of $\Gamma_c(G,r^*\mcal{A})$.  Furthermore, there
is a covariant representation $(\mu,G\unit*\mfrk{H},\rho,U)$ such that
$\Pi$ is equivalent to the integrated representation $\rho\rtimes U$.  
\end{theorem}

\begin{remark}
\label{rem:27}
It is worth pointing out that at this point we are deeply dependent on
separability hypotheses.  They were essential in disintegrating
representations of $C_0(X)$-algebras and they are also essential to
proving the disintegration theorem above.  This is still true if we
restrict to groupoid algebras, so that separability
assumptions will be required in that case as well. 
\end{remark}

One easy application of the disintegration theorem is the following  

\begin{corr}
\label{cor:7}
\index{covariant representation}
Suppose $(A,G,\alpha)$ is a separable dynamical system and $\pi$ is a
(nondegenerate) $*$-representation of $\Gamma_c(G,r^*\mcal{A})$ on some
Hilbert space $\mcal{H}$.  If $\pi$ is either $I$-norm decreasing or
continuous in the inductive limit topology then $\pi$ is equivalent to
the integrated form of some covariant representation. 
\end{corr}

\begin{proof}
If $\pi$ is $I$-norm continuous 
then it is continuous in the inductive
limit topology so that it suffices to address the case when $\pi$
is continuous in the inductive limit topology.  We will apply Theorem
\ref{thm:disintigration} with $\mcal{H}_0 = \mcal{H}$.  Condition (a)
holds since $\pi$ is nondegenerate.  Condition (c) holds because $\pi$
is a $*$-homomorphism.  Given $f_i\rightarrow f$
in the inductive limit topology we know that $\pi(f_i)\rightarrow
\pi(f)$.  However, this implies that $\pi(f_i)\rightarrow \pi(f)$ with
respect to the weak operator topology and it follows that (b) holds as well.
Hence, $\pi$ is equivalent to the integrated form of some covariant
representation and we are done. 
\end{proof}

This example presents a different manifestation of the left
regular representation that is much easier to deal with.  We have
developed it here in more generality than is strictly necessary for
our purposes.  As a result there are some Borel subtleties that make
things a bit more confusing.

\begin{example}
\label{ex:20}
\index{left regular representation}
Suppose $(A,G,\alpha)$ is a separable dynamical system and $\pi$ is a
nondegenerate representation of $A$ on a separable Hilbert space. Using Proposition \ref{prop:56},
assume without loss of generality that $\pi$ is a $C_0(G\unit)$-linear
representation with decomposition 
\[
\pi = \int^\oplus_{G\unit}\pi_u d\mu(u)
\]
on $L^2(G\unit*\mfrk{H},\mu)$.  Let $\nu\inv_Y = \int_{G\unit} \lambda_u d\mu(u)$
be the induced measure and $s^*(G\unit*\mfrk{H})$ the pull-back
bundle as in Example \ref{ex:35}.  
We define the {\em integrated left regular representation}
$L_\mu$, usually denoted  $L$, of $\Gamma_c(G,r^*\mcal{A})$
on $L^2(s^*(G\unit*\mfrk{H}),\nu\inv)$ by 
\[
L_\mu(f)h(\gamma) = \int_G \pi_{s(\gamma)}(\alpha_\gamma\inv(f(\eta)))
h(\eta\inv\gamma) d\lambda^{r(\gamma)}(\eta).
\]
We will show that $L_\mu$ is a nondegenerate, $I$-norm decreasing,
$*$-rep\-resent\-ation of $\Gamma_c(G,r^*\mcal{A})$.  It is relatively
straightforward to see that $L_\mu(f)h\in
\mcal{L}^2(s^*(G\unit*\mfrk{H}),\nu\inv)$ and that $L_\mu$ is linear.  Suppose
$f,g\in\Gamma_c(G,r^*\mcal{A})$ and $h,k\in
\mcal{L}^2(s^*(G\unit*\mfrk{H}),\nu\inv)$.  Then we compute
\begin{align*}
C &:=\int_{G\unit}\int_G\int_G
|(\pi_u(\alpha_\gamma\inv(f(\eta)))h(\eta\inv\gamma),k(\gamma))|
d\lambda^{r(\gamma)}(\eta)d\lambda_u(\gamma)d\mu(u) \\
&\leq \int_{G\unit}\int_G\int_G
\|f(\eta\inv)\|\|h(\eta\gamma)\|\|k(\gamma)\|
d\lambda_{r(\gamma)}(\eta)d\lambda_u(\gamma)d\mu(u) \\
&= \int_{G\unit}\int_G\int_G
\|f(\gamma\eta\inv)\|\|h(\eta)\|\|k(\gamma)\|d\lambda_u(\eta)d\lambda_u(\gamma)d\mu(u).
\end{align*}
Now, we define a measure $\beta$ on $C_c(G\times G)$ by 
\begin{equation}
\label{eq:45}
\beta(g) := \int_G\int_G g(\gamma,\eta)
d\lambda_u(\gamma)d\nu\inv(\eta) = \int_{G\unit}\int_G\int_G g(\gamma,\eta)
d\lambda_u(\gamma)d\lambda_u(\eta)d\mu(u).
\end{equation}
Using the Cauchy-Schwartz inequality with respect to $\beta$ we get
the following inequality. (Technically, need to know that
\eqref{eq:45} holds for positive Borel functions.  However this
follows from a straightforward application of  Proposition \ref{prop:57}.)
\begin{align}
\label{eq:46}
C \leq &\left( \int_{G\unit}\int_G\int_G \|f(\gamma\eta\inv)\|
  \|h(\eta)\|^2 d\lambda_u(\gamma)d\lambda_u(\eta)d\mu(u)\right)^{1/2} \\
\nonumber &\cdot \left( \int_{G\unit}\int_G\int_G \|f(\gamma\eta\inv)\|
  \|k(\gamma)\|^2 d\lambda_u(\gamma)d\lambda_u(\eta)d\mu(u)\right)^{1/2}
\end{align}
Now observe that 
\begin{align*}
\int_{G\unit}\int_G\int_G &\|f(\gamma\eta\inv)\|
  \|h(\eta)\|^2 d\lambda_{u}(\gamma)d\lambda_u(\eta)d\mu(u) \\
& =\int_G\left(\int_G \|f(\gamma\eta\inv)\|d\lambda_{s(\eta)}(\gamma)\right)
\|h(\eta)\|^2 d\nu\inv(\eta) \\
&=\int_G \left( \int_G
  \|f(\gamma)\|d\lambda_{r(\eta)}(\gamma)\right)\|h(\eta)\|^2
d\nu\inv(\eta) \\
&\leq \|f\|_I\|h\|^2.
\end{align*}
Using a similar computation for the term containing $f$ and $k$ in
\eqref{eq:46} we get
\begin{align} 
\nonumber C &=\int_{G\unit}\int_G\int_G
|(\pi_u(\alpha_\gamma\inv(f(\eta)))h(\eta\inv\gamma),k(\gamma))|
d\lambda^{r(\gamma)}(\eta)d\lambda_u(\gamma)d\mu(u) \\
 &\leq \|f\|_I\|h\|\|k\|. \label{eq:48}
\end{align}
This implies that the function 
\[
\gamma \mapsto \int_G
|(\pi_{s(\gamma)}(\alpha_\gamma\inv(f(\eta)))h(\eta\inv\gamma),k(\gamma))|
d\lambda^{r(\gamma)}(\eta)
\]
is $\nu\inv$-integrable.  Applying Proposition \ref{prop:57} to $\nu\inv$
implies that
\begin{equation}
\label{eq:47}
(L_\mu(f)h,k) = \int_{G\unit}\int_G\int_G
(\pi_u(\alpha_\gamma\inv(f(\eta)))h(\eta\inv\gamma),k(\gamma))
d\lambda^{r(\gamma)}(\eta)d\lambda_u(\gamma)d\mu(u)
\end{equation}
Hence \eqref{eq:48} also implies that $(L_\mu(f)h,k) \leq
\|f\|_I\|h\|\|k\|$ and therefore $L_\mu$ is $I$-norm decreasing.  

Next, we prove that $L_\mu$ is multiplicative by computing
\begin{align*}
L_\mu&(f*g)h(\gamma) = \int_G\int_G
\pi_{s(\gamma)}(\alpha_\gamma\inv(f(\zeta)\alpha_\zeta(g(\zeta\inv\eta))))
h(\eta\inv\gamma)d\lambda^{r(\gamma)}(\zeta)d\lambda^{r(\gamma)}(\eta) \\
&= \int_G\int_G
\pi_{s(\gamma)}(\alpha_\gamma\inv(f(\zeta)))\pi_{s(\gamma)}(\alpha_{\zeta\inv\gamma}\inv(g(\eta)))
h(\eta\inv\zeta\inv\gamma)d\lambda^{s(\zeta)}(\eta)d\lambda^{r(\gamma)}(\eta)
\\
&= \int_G
\pi_{s(\gamma)}(\alpha_\gamma\inv(f(\zeta)))(L_\mu(g)h)(\zeta\inv\gamma)
d\lambda^{r(\gamma)}(\zeta) \\
&= L_\mu(f)L_\mu(g)h(\gamma).
\end{align*}
We follow up by showing that $L_\mu$ preserves involution.  Observe that
we have to use \eqref{eq:47} again in order for this computation to hold.  
\begin{align*}
(L_\mu(f^*)h,k) &= \int_{G\unit}\int_G\int_G
(\pi_u(\alpha_{\gamma\inv\eta}(f(\eta\inv)^*))
h(\eta\inv\gamma) ,
k(\gamma))d\lambda^{r(\gamma)}(\eta)d\lambda_u(\gamma)d\mu(u) \\
&= \int_{G\unit}\int_G\int_G
(h(\eta\inv),\pi_u(\alpha_\eta(f(\eta\inv\gamma\inv)))k(\gamma))
d\lambda^u(\eta)d\lambda_u(\gamma)d\mu(u) \\
&= \int_{G\unit}\int_G\int_G
(h(\eta),\pi_{u}(\alpha_\eta\inv(f(\eta\gamma)))k(\gamma\inv)) 
d\lambda^u(\gamma)d\lambda_u(\eta)d\mu(u) \\
&= \int_G\int_G
(h(\eta),\pi_{s(\eta)}(\alpha_\eta\inv(f(\gamma)))k(\gamma\inv\eta))d\lambda^{r(\eta)}(\gamma)d\nu\inv(\eta)
\\
&= (h,L_\mu(f)k).
\end{align*}
This shows that $L_\mu$ is a $*$-homomorphism. 

Finally, we prove that $L_\mu$ is nondegenerate.  Suppose $(L_\mu(f)h,k) = 0$
for all $f\in \Gamma_c(G,r^*\mcal{A})$ and
$h\in\mcal{L}^2(s^*(G\unit*\mfrk{H}),\nu\inv)$.  It follows from
\eqref{eq:47} that 
\begin{align}
\label{eq:156}
0 &= (L_\mu(f)h,k) \\ &= \int_{G\unit}\int_G\int_G
(\pi_u(\alpha_\gamma\inv(f(\eta)))h(\eta\inv\gamma),k(\gamma))
d\lambda^{r(\gamma)}(\eta)d\lambda_u(\gamma)d\mu(u)\nonumber \\
&= \int_{G\unit}\int_G\int_G
(\pi_u(\alpha_\gamma\inv(f(\gamma\eta\inv)))
h(\eta),k(\gamma))d\lambda_u(\eta) d\lambda_u(\gamma)d\mu(u).\nonumber
\end{align}
Now, as in the proof of Proposition \ref{prop:59}, let $\{a_i\}$ be a
countable dense set in $A$ and $\{\phi_n\}$ a countable collection of
functions in $C_c(G)$ such that for all $\gamma\in G$ there exists $n$
such that $\phi_n(\gamma) = 1$.  Let $\{e_l\}$ be a special orthogonal
fundamental sequence for $G\unit*\mfrk{H}$ and let $\xi$ be an
arbitrary rational linear combination of the $e_l$.  Observe that
$\phi_n\otimes \xi(\gamma) = \phi_n(\gamma)\xi(s(\gamma))$ defines an element in
$\mcal{L}^2(s^*(G\unit*\mfrk{H}),\nu\inv)$.  By applying
\eqref{eq:156} to $\phi_m\otimes a_i$ and $\phi_n\otimes \xi$ for all
$n$,$m$,$i$, and $\xi$ we can, by taking a countable union, obtain a single
$\beta$-null set $N$ such that given $(\gamma,\eta)\not\in N$
\begin{equation}
\label{eq:49}
0 = \phi_m(\gamma\eta\inv)\phi_n(\eta)
(\pi_{s(\gamma)}(\alpha_\gamma\inv(a_i(r(\gamma))))\xi(s(\gamma)),k(\gamma))
\end{equation}
for all $n,m,i$ and $\xi$ a rational linear combination of the $e_l$.
Now, it follows from \eqref{eq:45} that we can find a $\nu\inv$-null
set $M$ such that given $\gamma\not\in M$ there exists $\eta\in G$
such that $(\eta,\gamma)\not\in N$.  Pick $n$ and $m$ so that
$\phi_m(\gamma\eta\inv) = \phi_n(\eta) = 1$.  Then \eqref{eq:49}
reduces to 
\begin{equation}
\label{eq:157}
0 = (\pi_{s(\gamma)}(\alpha_\gamma\inv(a_i(r(\gamma))))\xi(s(\gamma)),k(\gamma))
\end{equation}
for all $a_i$, $\xi$, and $\gamma\not\in M$.
Next, we can extend $M$ so
that $\pi_{r(\gamma)}$ is nondegenerate for all $\gamma\not\in M$.
Since we have made our choices so that the set
$\{\pi_{s(\gamma)}(\alpha_\gamma\inv(a_i(r(\gamma))))\xi(s(\gamma))\}$ is
dense in $\mcal{H}(s(\gamma))$, it follows from \eqref{eq:157} that
$k(\gamma) = 0$ for all $\gamma\not\in M$.  Thus $L_\mu$ is
nondegenerate. 

Note that Corollary \ref{cor:7} implies that $L_\mu$ must be equivalent to the
integrated form of some covariant representation.  
\end{example}

\begin{remark}
Readers may have observed that the computations in Example \ref{ex:20}
are very similar to those done in Proposition \ref{prop:59}.  This is
because the representation $L_\mu$ is ``close'' to being the integrated
form of the left regular representation from Example \ref{ex:19}.  
In fact, $L_\mu$ is equivalent to the integrated form of the left regular
representation.  Let $[\mu]$ be the
saturation of $\mu$ and $G\unit*\mfrk{K}$ be as in Example
\ref{ex:19}. Without going into too much
detail, the equivalence is implemented by the unitary
$U:L^2(s^*(G\unit*\mfrk{H}),\nu\inv)\rightarrow
L^2(G\unit*\mfrk{K},[\mu])$ where $Uh(u)(\gamma) = h(\gamma)$.  
We can then compute
\begin{align*}
UL_\mu(f)U^*h(u)(\gamma) &= L_\mu(f)U^*h(\gamma) \\
&= \int_G
\pi_{s(\gamma)}(\alpha_\gamma\inv(f(\eta)))U^*h(\eta\inv\gamma)d\lambda^u(\eta)\\
&= \int_G \pi_{s(\gamma)}(\alpha_\gamma\inv(
f(\eta)))h(s(\eta))(\eta\inv\gamma)\Delta(\eta)\poshalf\Delta(\eta)\neghalf
d\lambda^u(\eta)\\
&= \int_G \tilde{\pi}_u(f(\eta)) L_\eta h(s(\eta))(\gamma)
\Delta\neghalf(\eta) d\lambda^u(\eta) \\
&= \tilde{\pi}\rtimes L(f)h(u)(\gamma).
\end{align*}
Of course, we could have saved ourselves some trouble and started by
proving that $L_\mu$ is equivalent to $\tilde{\pi}\rtimes L$ and then used
Proposition \ref{prop:59}.  However, the construction of
$(\tilde{\pi},L)$ was really only outlined in Example \ref{ex:19}, and
it is not entirely obvious that $U$ is a unitary.  The
fact is that working with the integrated left regular representation 
$L_\mu$ directly is easier.  As we have said, one of the great
things about Theorem \ref{thm:disintigration} is that it frees us up
from having to work with covariant representations.    
\end{remark}

It takes some effort, but we will show that we can separate points in
$\Gamma_c(G,r^*\mcal{A})$ with these integrated regular
representations.  Unfortunately,  we will have to forward
reference Lemma \ref{lem:19}.  This isn't a problem, however, since
Lemma \ref{lem:19} is taken directly from \cite{renaultequiv}.  

\begin{lemma}
\label{lem:14}
Suppose $(A,G,\alpha)$ is a separable groupoid dynamical system.  Then
the (integrated forms of) covariant representations separate points of
$\Gamma_c(G,r^*\mcal{A})$.  
\end{lemma}

\begin{proof}
Clearly it suffices to show that if $f_0\ne 0$ in
$\Gamma_c(G,r^*\mcal{A})$ then there exists a covariant representation
such that $\pi\rtimes U(f_0) \ne 0$.  However, once we consider
Corollary \ref{cor:7} it clearly suffices to show that there exists a
(nondegenerate) $I$-norm decreasing representation $L$ such that
$L(f_0)\ne 0$. Find $\gamma_0\in G$
such that $f_0(\gamma_0)\ne 0$.  Set $u= s(\gamma_0)$ and let 
$\rho$ be a faithful representation of $A(u)$ on $\mcal{H}$.  It is easy to
see that the lift of $\rho$ to $A$, denoted $\pi$, 
is $C_0(G\unit)$-linear and that
the decomposition is given on the trivial bundle $G\unit\times
\mcal{H}$ by
\[
\pi = \int_G \pi_v d\delta_{u}(v)
\]
where $\pi_{u} = \rho$ and $\pi_v = 0$ for all $v\ne u$.
Furthermore, observe that the induced measure $\nu\inv$ is
exactly $\lambda_{u}$ in this case.  Let us form the integrated left regular
representation as in Example \ref{ex:20}.  After sorting through all
the definitions we find we have a representation $L_{\delta_{u}}$,
denoted $L$, of
$\Gamma_c(G,r^*\mcal{A})$ on $L^2(G_u,\mcal{H},\lambda_u)$ given by 
\begin{equation}
\label{eq:50}
L(f)h(\gamma) = \int_G \rho(\alpha_\gamma\inv(f(\eta)))
h(\eta\inv\gamma) d\lambda^{r(\gamma)}(\eta).
\footnote{Notice that most of the measure theoretic difficulties of Example
\ref{ex:20} disappear because $\nu\inv = \lambda_u$ and $\beta =
\lambda_u\times \lambda_u$.}
\end{equation}

Let $e_i$ be an orthonormal basis for $\mcal{H}$ and for each
$g\in \Gamma_c(G,r^*\mcal{A})$ and $\gamma\in G_u$ define 
\[
\Phi_i(g)(\gamma) = \rho(\alpha_\gamma\inv(g(\gamma)))e_i
\]
Then $\Phi_i(g)\in C_c(G_u,\mcal{H})$ and
$\Phi_i:\Gamma_c(G,r^*\mcal{A})\rightarrow
L^2(G_u,\mcal{H},\lambda_u)$.  What's more, suppose $\Phi_i(g)= 0$ for
all $i$.  Then $\rho(\alpha_\gamma\inv(g(\gamma)))e_i = 0$ almost
everywhere for all $i$.  Thus, we can find a single $\lambda_u$-null
set $N$ such that $\rho(\alpha_\gamma\inv(g(\gamma))) = 0$ for all
$\gamma\not\in N$.  Since $\rho$ is faithful this implies that
$g(\gamma) = 0$ for all $\gamma\not\in N$.   Using the fact that $g$
is continuous and $\supp \lambda_u = G_u$ we conclude that $g(\gamma)
= 0$ for all $\gamma\in G_u$.  

Next, given $g\in\Gamma_c(G,r^*\mcal{A})$ we calculate
\begin{align*}
L(f)\Phi_i(g)(\gamma) &= \int_G
\rho(\alpha_\gamma\inv(f(\eta)))\rho(\alpha_{\eta\inv\gamma}\inv(g(\eta\inv\gamma)))e_i
d\lambda^{r(\gamma)}(\eta) \\
&= \int_G
\rho(\alpha_\gamma\inv(f(\eta)\alpha_\eta(g(\eta\inv\gamma))))e_i
d\lambda^{r(\gamma)}(\eta)\\
&= \rho(\alpha_\gamma\inv(f*g(\gamma)))e_i = \Phi_i(f*g)(\gamma).
\end{align*}
Now suppose, to the contrary, that $L(f_0)= 0$.  It follows from Lemma
\ref{lem:19} that there exists a left approximate identity $\{g_j\}$ in
$\Gamma_c(G,r^*\mcal{A})$ with respect to the inductive limit
topology.  By replacing $g_j$ with $g_j^*$ we may assume that
$\{g_j\}$ is a right approximate identity.  Thus, $f_0*g_j \rightarrow
f_0$ with respect to the inductive limit topology.  However, we have
$L(f_0)\Phi_i(g_j) = \Phi_i(f_0*g_j) = 0$ for all $i$ and $j$.  It follows
from the previous paragraph that $f_0*g_j = 0$ on $G_u$.  Since
$f_0*g_j\rightarrow f_0$ uniformly, we must have $f_0(\gamma) = 0$ 
for all $\gamma\in G_u$, but
this is a contradiction since $f_0(\gamma_0) \ne 0$.  Thus
$L(f_0)\ne 0$ and we are done.
\end{proof}

\begin{remark}
We can actually make some considerable upgrades to Lemma
\ref{lem:14}.  In particular, we can identify a class of
representations of $A$ for which the integrated left regular
representation is faithful on $\Gamma_c(G,r^*\mcal{A})$.  Suppose
$\pi$ is a faithful $C_0(G\unit)$-linear representation of $A$ with
decomposition 
\[
\pi = \int_{G\unit}^\oplus \pi_u d\mu(u)
\]
on $L^2(X*\mfrk{H},\mu)$ 
such that $\pi_u$ is faithful for almost all $u$.  Observe that this
happens automatically if $A$ has Hausdorff spectrum.\footnote{It is not
  clear if many such representations exist in the non-Hausdorff case.  
  There seems to be no reason why the decomposition of a
  faithful representation should be faithful almost everywhere.}  It is
straightforward to show that for $\pi$ to be faithful $\mu$ must have
full support.  It then follows quickly that the induced measure
$\nu\inv$ has full support as well.  Let $e_i$ be a special orthogonal
fundamental sequence for $X*\mfrk{H}$ and given $g\in
\Gamma_c(G,r^*\mcal{A})$ define
\[
\Phi_i(g)(\gamma) =
\pi_{s(\gamma)}(\alpha_\gamma\inv(g(\gamma)))e_i(s(\gamma)).
\]
It is straightforward to show that $\Phi_i(g)$ is a bounded Borel
function supported on a compact set so that for each $i$ we have
$\Phi_i:\Gamma_c(G,r^*\mcal{A})\rightarrow
L^2(s^*(X*\mfrk{H}),\nu\inv)$.  Suppose $\Phi_i(g) = 0$ for all $i$.
Using the usual trick, we can find a single $\nu\inv$-null set $N$ such
that 
\[
\pi_{s(\gamma)}(\alpha_\gamma\inv(g(\gamma)))e_i(s(\gamma)) = 0
\]
for all $i$ given $\gamma\not\in N$.  Since $e_i$ is a special
orthogonal fundamental sequence, this implies that
$\pi_{s(\gamma)}(\alpha_\gamma\inv(g(\gamma)))=0$ for all
$\gamma\not\in N$.  We may thicken $N$ a bit and assume that
$\pi_{s(\gamma)}$ is faithful for all $\gamma\not\in N$ and hence
$g(\gamma) = 0$ for all $\gamma\not\in N$.  Since $g$ is continuous
and $\supp\nu\inv = G$ this implies that $g=0$.  

Now suppose $f\in \Gamma_c(G,r^*\mcal{A})$ is such that $L_\mu(f) =
0$.  Then, just as in the previous lemma,
it is a simple matter to prove that $0 = L_\mu(f)\Phi_i(g) =
\Phi_i(f*g)$ for all $g\in \Gamma_c(G,r^*\mcal{A})$.   
Use Lemma \ref{lem:19} to find a {\em right}
approximate identity $\{g_j\}$ 
with respect to the inductive limit topology.  Then $\Phi_i(f*g_j) =
0$ for all $i$ and $j$.  By the above paragraph this implies that
$f*g_j=0$ for all $j$.  Since $f*g_j\rightarrow f$ uniformly this
implies $f=0$.  Thus $L_\mu$ is faithful on
$\Gamma_c(G,r^*\mcal{A})$.  This is related to \cite[Proposition
1.11]{groupoidapproach}, which states that the left regular
representation is faithful in the scalar case as long as $\nu\inv$ has
full support. 
\end{remark}

We are, at long last, ready to define the groupoid crossed product.
This definition is slightly different than the definition given in
\cite{renaultequiv} but, as will see, because of the Disintegration
Theorem the universal norm can be obtained in any number of ways.  

\begin{definition}
\label{def:44}
\index{universal norm}
\index{crossed product}
\index[not]{$A\rtimes_\alpha G$}
Suppose $(A,G,\alpha)$ is a separable groupoid dynamical system.  We
define the {\em universal norm} on $\Gamma_c(G,r^*\mcal{A})$ by 
\begin{equation}
\|f\| := \sup\{ \|\pi\rtimes U(f)\| : \text{$(\pi,U)$ is a covariant
  representation of $(A,G,\alpha)$}\}.
\end{equation}
The completion of $\Gamma_c(G,r^*\mcal{A})$ with respect to this norm
is a $C^*$-algebra called the {\em groupoid crossed product} of $A$ by
$G$ and denoted $A\rtimes_\alpha G$.  
\end{definition}

\begin{remark}
To avoid clutter we will denote $A\rtimes_\alpha G$ by
$A\rtimes G$ whenever the chance of confusion is near zero. 
\end{remark}

\begin{remark}
There is also the notion of a reduced groupoid crossed product which
mimics the group case.  Instead of taking the universal norm to be the
supremum over all covariant representations, we take it to be the
supremum over the left regular representations, either those in Example
\ref{ex:19} or Example \ref{ex:20}.  
While reduced crossed products are interesting, they are less well
studied and we will not deal with them here.  Those who would like to
know more are referenced to \cite{jonthon}.
\end{remark}

Let us verify the claims made in Definition \ref{def:44} and begin our
exploration of the crossed product.  

\begin{prop}
\label{prop:58}
Suppose $(A,G,\alpha)$ is a separable dynamical system.  Then the
universal norm is a norm on $\Gamma_c(G,r^*\mcal{A})$ 
and is dominated by the $I$-norm.
Furthermore, the completion $A\rtimes_\alpha G$ is a separable
$C^*$-algebra.  
\end{prop}

\begin{remark}
\label{rem:14}
It follows that convergence with respect to the $I$-norm is stronger
than convergence with respect to the universal norm.  Furthermore, we
have from Proposition \ref{prop:60} that convergence with respect to
the inductive limit topology is stronger than convergence with respect
to the $I$-norm, and therefore with respect to 
the universal norm as well.  We will use this fact often.  
For example, we can extend Corollary \ref{cor:7} and conclude that
any representation of $A\rtimes G$ is equivalent to the integrated
form of a covariant representation. 
\end{remark}

\begin{proof}
Since the integrated form of a covariant representation is $I$-norm
decreasing, it is clear that $\|f\|\leq \|f\|_I$ for all $f\in
\Gamma_c(G,r^*\mcal{A})$.  Furthermore, 
given $f,g\in\Gamma_c(G,r^*\mcal{A})$ we have 
\[
\|\pi\rtimes U(f+g)\| \leq \|\pi\rtimes U(f)\|+\|\pi\rtimes U(g)\| \leq 
\|f\|+\|g\|
\]
and it follows that $\|f+g\|\leq \|f\|+\|g\|$.  The fact
that $\|\pi\rtimes U(\lambda f)\| = |\lambda|\|\pi\rtimes U(f)\|$
implies that $\|\lambda f\| = |\lambda|\|f\|$ for $\lambda\in\C$.  All
that remains to show that $\|\cdot\|$ is a norm is to prove that $\|f\|
\ne 0$ if $f\ne 0$. However, this follows immediately from Lemma
\ref{lem:14}.  

Moving on, it is straightforward to show that $\|f*g\|\leq \|f\|\|g\|$ for
$f,g\in \Gamma_c(G,r^*\mcal{A})$, and since $A\rtimes_\alpha G$ is the
completion, it follows that the norm is submultiplicative on the entire
crossed product.  Similarly we find that $\|f^*\| =\|f\|$ and  
$\|f^*f\| = \|f\|^2$ for all $f\in A\rtimes_\alpha G$. It follows that 
$A\rtimes_\alpha G$ is a $C^*$-algebra.  The last thing we need to do
is show that $A\rtimes_\alpha G$ is separable.  
Let $\{a_i\}$ be a countable dense set in
$A$.  Since $G$ is second countable we can find a countable set
$\{\phi_j\}$ which is dense in $C_c(G)$ with respect to the inductive
limit topology.  Now, given $\phi\otimes a\in C_c(G)\odot A$ suppose
$a_{i_k}\rightarrow a$ and $\phi_{j_k}\rightarrow \phi$ with respect
to the inductive limit topology.  Then 
\[
\|\phi_{j_k}(\gamma)a_{i_k}(r(\gamma)) - \phi(\gamma)a(r(\gamma))\| \leq
\|\phi_{j_k}-\phi\|_\infty \|a_{i_k}\| - \|\phi\|_\infty\|a_{i_k}-a\|
\]
Since $\{\|a_{i_k}\|\}$ is bounded, this shows that $\phi_{i_k}\otimes
a_{i_k}\rightarrow \phi\otimes a$ uniformly.  Furthermore, since the
supports of the $\phi_{j_k}$ are eventually contained in a fixed compact
set, the same is true for  $\phi_{j_k}\otimes a_{i_k}$. 
It follows quickly that the countable set $D$ of rational sums of elements of
the form $\phi_j\otimes a_i$ is
dense in $C_c(G)\odot A$ with respect to the inductive limit
topology, and hence with respect to the $I$-norm.  We conclude
from Corollary \ref{cor:3}
that $D$ is dense in $\Gamma_c(G,r^*\mcal{A})$
with respect to the $I$-norm.  Since the universal norm is bounded by
the $I$-norm, this enough to show that $A\rtimes G$ is separable.  
\end{proof}

The following identities will be quite useful when dealing with
crossed products.  Both are immediate results of the Disintegration
Theorem.  

\begin{prop}
\label{prop:61}
Suppose $(A,G,\alpha)$ is a separable groupoid dynamical system.  Then
the universal norm on $\Gamma_c(G,r^*\mcal{A})$ is also given by
\begin{align}
\label{eq:42}
\|f\| &= \sup\left\{\|\pi(f)\| : \begin{array}{l} 
 \text{$\pi$ is a (possibly degenerate) $I$-norm decreasing} \\ 
\text{$*$-representation of 
      $\Gamma_c(G,r^*\mcal{A})$}\end{array}\right\}, \\
\label{eq:43}
\|f\| &= \sup\left\{\|\pi(f)\| : \begin{array}{l}
\text{$\pi$ is a (possibly degenerate) $*$-representation
  of $\Gamma_c(G,r^*\mcal{A})$} \\ 
\text{which is continuous in the inductive limit topology} \end{array}\right\}.
\end{align}
It follows that any (possibly degenerate)
$*$-representation of $\Gamma_c(G,r^*\mcal{A})$
which is either $I$-norm decreasing or continuous with respect to the
inductive limit topology is bounded with respect to the universal norm
and extends to a representation of $A\rtimes_\alpha G$.  
\end{prop}
\begin{proof}
Let $\|\cdot\|_1$ and $\|\cdot\|_2$ be defined by \eqref{eq:42} and
\eqref{eq:43} respectively.  Since every $I$-norm decreasing
representation is continuous with respect to the inductive limit
topology, we have $\|\cdot\|_1 \leq \|\cdot\|_2$.  Furthermore, since
each covariant representation is $I$-norm decreasing, we have
$\|\cdot\| \leq \|\cdot \|_1$.  Now suppose $f\in
\Gamma_c(G,r^*\mcal{A})$ and $\pi$ is a $*$-representation which is
continuous with respect to the inductive limit topology.  Let
$\pi_{\ess}$ be the restriction to its essential subspace
$\mcal{H}_{\ess} = \cspn\{\pi(f)h: f\in\Gamma_c(G,r^*\mcal{A}), h\in\mcal{H}\}$.  Let $U$ be the
unitary map from $\mcal{H}$ onto $\mcal{H}_{\ess}\oplus
\mcal{H}_{\ess}^\perp$.  Then given
$(h,k)\in\mcal{H}_{\ess}\oplus\mcal{H}_{\ess}^\perp$ we have
\[
U\pi(f)U^*(h,k) =  U(\pi(f)h+\pi(f)k) = (\pi_{\ess}(f)h+\pi_{\ess}(f)k,0).
\]
However, given $l\in\mcal{H}$ observe that
\[
(\pi(f)k,l) = (k,\pi(f^*)l) = 0
\]
since $k\in\mcal{H}_{\ess}^\perp$.  It follows that $\pi(f)k = 0$ so
that 
\[
U\pi(f)U^*(h,k) = (\pi_{\ess}(f)h,0).
\]
It follows that $\pi$ is unitarily equivalent to the representation 
$\pi_{\ess}\oplus 0$.  Hence $\|\pi(f)\| = \|\pi_{\ess}\oplus 0(f)\| =
\|\pi_{\ess}(f)\|$ for all $f\in \Gamma_c(G,r^*\mcal{A})$.
Furthermore, it is straightforward to show that, by
construction, $\pi_{\ess}$ is a nondegenerate
$*$-representation of $\Gamma_c(G,r^*\mcal{A})$ which is continuous
with respect to the inductive limit topology.  It follows
from Corollary \ref{cor:7} that $\pi$ is equivalent to the integrated
form of a covariant representation $(\rho,U)$.  Therefore
\[
\|\pi(f)\| = \|\pi_{\ess}(f)\| = \|\rho\rtimes U(f)\| \leq \|f\|.  
\]
Since $\pi$ was
generic, it follows that $\|\cdot\|_2\leq\|\cdot\|$ and we have demonstrated
\eqref{eq:42} and \eqref{eq:43}.  The second half of the proposition
is clear from the first half.  
\end{proof}

\begin{remark}
This shows that
$A\rtimes_\alpha G$ is the ``universal enveloping algebra'' of
$\Gamma_c(G,r^*\mcal{A})$ with respect to the $I$-norm, bringing us in
line with \cite{renaultequiv}. 
\end{remark}

We end this section with a proposition that we will use liberally.  It
also hints at the author's general philosophy for working with crossed
products, which is ``stick to the inductive limit topology whenever possible.''

\begin{prop}
\label{prop:62}
Suppose $(A,G,\alpha)$ is a 
separable groupoid dynamical system, $D$ is a $C^*$-algebra, and that
$\Phi:\Gamma_c(G,r^*\mcal{A})\rightarrow D$ is
a $*$-homomorphism 
such that $\Phi$ is either continuous with respect to the inductive limit
topology or $I$-norm decreasing.  
Then $\Phi$ is norm decreasing with respect to the
universal norm and extends to a homomorphism of $A\rtimes_\alpha G$
into $D$.  
\end{prop}

\begin{proof}
Clearly we can restrict to the case where $\Phi$ is continuous with
respect to the inductive limit topology. 
Suppose $R$ is a faithful representation of $D$.  It
follows that $R \circ\Phi$ is continuous with respect to the
inductive limit topology.  Therefore \eqref{eq:43} implies that
$\|\Phi(f)\| = \|R(\Phi(f))\| \leq \|f\|$ for all
$f\in\Gamma_c(G,r^*\mcal{A})$.  Thus $\Phi$ is norm decreasing.  
The remainder of the proposition is clear. 
\end{proof}

We usually apply this proposition via the following corollary. 

\begin{corr}
\label{cor:24}
\index{inductive limit topology}
Suppose $(A,G,\alpha)$ and $(B,H,\beta)$ are separable groupoid dynamical
systems and that
$\Phi:\Gamma_c(G,r^*\mcal{A})\rightarrow\Gamma_c(H,r^*\mcal{B})$ is a $*$-homomorphism.  If
$\Phi(f_i)\rightarrow \Phi(f)$ with respect to the inductive limit
topology  whenever $f_i\rightarrow f$ with respect to the inductive
limit topology then $\Phi$ extends to a $*$-homomorphism from $A\rtimes_\alpha G$
into $B\rtimes_\beta H$. 
\end{corr}

\begin{proof}
This follows immediately from Proposition \ref{prop:62} and the fact
that convergence with respect to the inductive limit topology is
stronger than convergence with respect to the universal norm.  
\end{proof}


\chapter{Special Cases}
\label{cha:special-cases}
In this chapter we show how groupoid
crossed products generalize group crossed products in Section
\ref{sec:groupprod}, groupoid 
$C^*$-algebras in Section \ref{sec:scalarprod} and transformation
group algebras in Section \ref{sec:transform}.  In addition, we will
deal with the special case where the groupoid is actually a group
bundle in Section \ref{sec:bundleprod}, and when the action arises
from a continuous $G$-space in Section \ref{sec:transform}.
However, one
thing this chapter, and indeed the whole theory, lacks is natural
examples of groupoid dynamical systems which do not come from one of
these special cases.  The material in this section is essential for
two reasons.  First, it provides a testing ground for the theory
developed in the last chapter and many of the arguments we will use
later can be found here in simpler form.  Second, the connections made
in this chapter provide a paradigm for transporting existing theory to
the groupoid case.  

\section{Group Crossed Products}
\label{sec:groupprod}

This section deals with group crossed products as defined and
developed in \cite{tfb2}.  In particular it is assumed that the reader
is familiar with at least \cite[Chapter 1]{tfb2}.  

Groupoid crossed products are an extension of group crossed products
in two different ways.  Both are important.  The first is the
situation where the groupoid $G$ is a actually group.  Start by observing that
in this case 
the unit space of $G$ is a single point $\{e\}$ 
and that any $C^*$-algebra $A$ can be  viewed as a $C_0(\{e\})$-algebra.  
In this case $A$ only has one fibre,
the bundle associated to $A$ is just $A$, and the axioms of
a dynamical system are reduced to the following

\begin{definition}
\label{def:45}
\index{dynamical system}
Suppose $G$ is a locally compact Hausdorff group and $A$ is a
$C^*$-algebra. An action $\alpha$ of $G$ on $A$ consists of a family
of maps $\{\alpha_s\}_{s\in G}$ such that 
\begin{enumerate}
\item $\alpha_s$ is an automorphism of $A$,
\item $\alpha_{st} = \alpha_s\circ\alpha_t$ for all $s,t\in G$, and
\item $s\cdot a:= \alpha_s(a)$ defines a continuous action of $G$ on
  $A$.
\end{enumerate}
\end{definition}

Of course, this can be simplified a great deal.  

\begin{prop}
\label{prop:63}
Suppose $G$ is a locally compact Hausdorff group and $A$ is a
$C^*$-algebra.  Then $\alpha$ is an action of $G$ on $A$ if and only
if $\alpha:G\rightarrow \Aut(A)$ is a continuous homomorphism where
$\Aut(A)$ is equipped with the topology of pointwise convergence. 
\end{prop}

\begin{proof}
Suppose $\alpha$ is an action of $G$ on $A$.  It is clear that 
$\alpha:G\rightarrow \Aut(A)$ is a continuous homomorphism. Suppose
$\alpha$ is a continuous homomorphism.  Obviously the first two
conditions of Definition \ref{def:45} are satisfied.  Now if
$s_i\rightarrow s$ and $a_i\rightarrow a$ then 
\[
\|\alpha_{s_i}(a_i)-\alpha_s(a)\| \leq
\|\alpha_{s_i}(a_i-a)\| + \|\alpha_{s_i}(a) -\alpha_s(a)\|
\leq
\|a_i-a\| +\|\alpha_{s_i}(a)-\alpha_s(a)\|.
\]
Since $\alpha$ is continuous into the topology of pointwise
convergence it follows that $\alpha$ is an action of $G$ on $A$.
\end{proof}

This shows that the restriction of Definition \ref{def:33} to groups
is exactly the classical definition of a group action on a
$C^*$-algebra as given in \cite[Definition 2.6]{tfb2}.  

\begin{remark}
\label{rem:20}
\index{covariant representation}
Let us consider covariant representations.  First recall that any
group $G$ has a Haar system given by the Haar measure on $G$.
Furthermore, any 
measure on $\{e\}$ is trivially quasi-invariant and the modular function is
always given by the classical modular function of $G$ with respect to
its Haar measure.  It is also worth noting that the modular function of a
group is always {\em continuous}.  Now, according to 
Definition \ref{def:40} a
(Borel) representation of $G$ is just a Borel homomorphism from $G$
into the unitary operators on some separable Hilbert space equipped with the
Borel structure coming from the strong topology. However it follows
from \cite[Theorem D.3]{tfb2} that in this case $U$ is actually (strongly)
continuous.  Any representation $\pi$ of $A$ is trivially
$C_0(\{e\})$-linear and, since there is only one fibre, $\pi$ is its own
decomposition.  Thus, in the group case, a covariant representation of
$(A,G,\alpha)$ is given by a representation $\pi$ of $A$ and a
representation $U$ of $G$ such that 
\begin{equation}
\label{eq:52}
U_s \pi(a) = \pi(\alpha_s(a))U_s
\end{equation}
for almost all $s\in G$.  However, $U$ is continuous so that
\eqref{eq:52} holds for all $s\in G$.  Thus a covariant
representation of $(A,G,\alpha)$ in the groupoid sense is exactly a
covariant representation of $(A,G,\alpha)$ as defined in
\cite[Definition 2.10]{tfb2}.  
\end{remark}

It seems reasonable that, since the covariant representations are the
same, then the crossed products of $A$ by $G$ as a group and as a
groupoid must be the same.  Unfortunately there is the problem of the
modular function.  First, observe that $r^*A = G\times A$ and that
$\Gamma_c(G,r^*A)$ can be identified with $C_c(G,A)$.  Recall from Proposition
\ref{prop:48}
that we defined the operations on $C_c(G,A)$ to be given by 
\begin{align}
\label{eq:55}
f*g(s) &= \int_G f(t)\alpha_t(g(t\inv s))d\lambda(s),\quad\text{and} \\
\label{eq:53}
f^*(s) &= \alpha_s(f(s\inv)^*)
\end{align}
where $\lambda$ is the Haar measure on $G$.  In \cite[Section
2.3]{tfb2} the convolution operation on $C_c(G,A)$ is also defined by \eqref{eq:53}.  {\em However}, the involution operation is
defined by 
\begin{equation}
\label{eq:54}
f^*(s) = \Delta(s\inv)\alpha_s(f(s\inv)^*).
\end{equation}
Notice that this definition of involution only works in the group
case because for groupoids the modular function depends on the choice
of the quasi-invariant measure, while for groups the modular function
only depends on the Haar measure.  This is compensated by the fact that 
the integrated form of a covariant representation $(\pi,U)$ as a
groupoid dynamical system is given by \index{covariant representation!integrated form}
\[
\pi\rtimes U(f) = \int_G \pi(f(s))U_s \Delta(s)\neghalf
d\lambda(s)
\]
while the integrated form as a group dynamical system is defined in
\cite[Proposition 2.23]{tfb2} to be 
\[
\pi\rtimes' U(f) = \int_G \pi(f(s)) U_s d\lambda(s).
\]

Of course, when we refer to a group
crossed product we will use \eqref{eq:53} and the universal norm
coming from Definition \ref{def:44}.  However, it is important to see
that this group crossed product is naturally isomorphic to the one
defined in \cite[Section 2.3]{tfb2}, and to sort out the differences
described above.  

\begin{prop}
\index{crossed product!group}
\label{prop:64}
Suppose $(A,G,\alpha)$ is a separable dynamical system and that $G$ is
a group.  Let $A\rtimes_\alpha G$ be the crossed product as defined in
Definition \ref{def:44}.  Let $A\rtimes'_\alpha G$ be the crossed
product as defined in \cite[Lemma 2.27]{tfb2}.  Then the map
$\Phi:C_c(G,A)\rightarrow C_c(G,A)$ given by $\Phi(f)(s) =
\Delta(s)\neghalf f(s)$ extends to a $*$-isomorphism of $A\rtimes_\alpha G$ and
$A\rtimes'_\alpha G$.  
\end{prop}

\begin{remark}
\index{modular function}
The point of Proposition \ref{prop:64} is that for a group
dynamical system the modular function can be either included in
the formula for involution, or in the formula for the integrated form
of a covariant representation.  Either way we end up with isomorphic 
crossed products. We will view group crossed products as
special cases of groupoid crossed products.  
\end{remark}

\begin{proof}
Let $(A,G)$ and $\Phi$ be as above.  Let $C$ denote $C_c(G,A)$ given
the operations coming from \eqref{eq:55} and \eqref{eq:53} so that
$A\rtimes_\alpha G$ is the completion of $C$.  Let $D$ denote
$C_c(G,A)$ given the operations \eqref{eq:55} and \eqref{eq:54} so
that $A\rtimes'_\alpha G$ is the completion of $D$.  Let us calculate
\begin{align*}
\Phi(f)*\Phi(g)(s) &= \int_G \Phi(f)(t)\alpha_t(\Phi(g)(t\inv s))
d\lambda(t) \\
&= \int_G \Delta(t)\neghalf\Delta(t\inv s)\neghalf f(t)\alpha_t(g(t\inv s))d\lambda(t)
\\ &= \Delta(s)\neghalf f*g(s) = \Phi(f*g)(s). 
\end{align*}
We also have 
\begin{align*}
\Phi(f)^*(s) &= \Delta(s\inv)\alpha_s(\Phi(f)(s\inv)^*)
= \Delta(s\inv)\Delta(s\inv)\neghalf \alpha_s(f(s\inv)^*) \\
&= \Delta(s)\neghalf f^*(s) = \Phi(f^*)(s).
\end{align*}
Thus $\Phi$ is a $*$-homomorphism.  Furthermore, it is clear that
$\Phi$ is a bijection with inverse given by $\Phi\inv(f)(s) =
\Delta(s)\poshalf f(s)$.  

\begin{remark}
We would like to see that $\Phi$ is bounded and as such extends to the
entire crossed product.  A valid approach to this problem would be to
use the fact that the sets of covariant representations of
$(A,G,\alpha)$ as both a group and a groupoid dynamical system 
are the same and to show that $\Phi$ intertwines the two
different forms of the integrated representation.  Since the universal
norms of $A\rtimes G$ and $A\rtimes' G$ are both given by supremums
over the covariant representations, the result would follow.  However,
it is easier, and instructive, to use Proposition \ref{prop:61}.
\end{remark}

Suppose $f_i\rightarrow f$ with respect to the inductive limit
topology so that  $f_i\rightarrow f$ uniformly and $\supp f_i\subset
K$ eventually, for a fixed compact set $K$.  Since $\Delta$ is
continuous in the group case we can find $M$ such that
$\Delta(s)\neghalf < M$ for all $s\in K$.  Then we have 
\[
\|\Phi(f_i)(s)-\Phi(f)(s)\| = \Delta(s)\neghalf \|f_i(s)-f(s)\| \leq
M\|f_i-f\|_\infty
\]
and it follows quickly that $\Phi(f_i)\rightarrow \Phi(f)$ with
respect to the inductive limit topology.  
Suppose $\pi$ is  faithful representation of $A\rtimes' G$.  We know
from \cite[Corollary 2.46]{tfb2} that 
\begin{equation}
\label{eq:56}
\|f\|_{A\rtimes' G} = \sup\{ \|\pi(f)\|: \text{$\pi$ is continuous in
  the inductive limit topology}\}.
\end{equation}
It follows that $\pi$ is continuous in the inductive limit topology so
that the composition $\pi\circ\Phi$ is continuous in the inductive
limit topology.  It is clearly a $*$-homomorphism.  
Hence Proposition \ref{prop:61} and
\eqref{eq:43} imply that 
\[
\|\Phi(f)\| = \|\pi(\Phi(f))\| \leq \|f\|.
\]
Thus $\Phi$ is bounded and extends to $A\rtimes G$.  We can use
exactly the same argument in reverse to see that
$\Phi\inv$ is bounded and extends to $A\rtimes' G$.  Since the
extensions of $\Phi$ and $\Phi\inv$ are inverses on a dense set they
are inverses everywhere and $\Phi$ extends to an isomorphism. 
\end{proof}

Now, Proposition \ref{prop:64} is theoretically important but it is
not the ``right'' way to view groupoid crossed products as a
generalization of group crossed products.  The following construction
is much more useful in that regard, as we will see in Sections
\ref{sec:locally-unitary} and \ref{sec:regularity}. 

\begin{example}
\label{ex:22}
\index{crossed product!group}
Suppose $A$ is a separable $C^*$-algebra with Hausdorff spectrum $\widehat{A}$
and recall from Example \ref{ex:21} that we can view $A$ as a
$C_0(\widehat{A})$-algebra with fibres $A(\pi)=A/\ker\pi$.  Let $\mcal{A}$ be
the corresponding upper-semicontinuous bundle.  Suppose 
$\alpha$ is an action of a second countable locally compact Hausdorff
group $H$ on $A$.  Recall from \cite[Lemma 2.8]{tfb2} that there is a
jointly continuous action of $G$ on $\widehat{A}$ given by $s\cdot \pi =
\pi\circ\alpha_s\inv$ and let $G$ be the corresponding
transformation groupoid $G=H\ltimes \widehat{A}$.  Observe that the unit
space of $G$ can be identified with $\widehat{A}$ as in Example \ref{ex:4}.  

Our goal is to describe an action of $G$ on $A$. 
Given $(s,\pi)\in G$ define $\beta_{(s,\pi)}:A(s\inv\cdot
\pi)\rightarrow A(\pi)$ by 
\begin{equation}
\label{eq:57}
\beta_{(s,\pi)}(a(s\inv\cdot \pi)) = \alpha_s(a)(\pi).
\end{equation}
We would like to show that $\beta_{(s,\pi)}$ is a well defined
isomorphism.  First, suppose $a,b\in A$ such that $a(s\inv\cdot\pi) =
b(s\inv\cdot \pi)$.  However, this just says that 
\[
0=s\inv\cdot \pi(a-b) = \pi(\alpha_s(a-b))
\]
and hence $\alpha_s(a)(\pi) = \alpha_s(b)(\pi)$.  Thus
$\beta_{(s,\pi)}$ is well defined.  Next, since each $\alpha_s$ is
a $*$-homomorphism it is straightforward to show that
$\beta_{(s,\pi)}$ is as well.  Given $(e,\pi)\in G$ we have 
\[
\beta_{(e,\pi)}(a(\pi)) = \alpha_e(a)(\pi) = a(\pi)
\]
so that $\beta_{(e,\pi)} = \id$.  Furthermore, given
$(s,\pi),(t,s\inv\cdot\pi)\in G$ we have 
\begin{align*}
\beta_{(s,\pi)}\circ\beta_{(t,s\inv\cdot \pi)}(a(t\inv s\inv\cdot
\pi)) &= \beta_{(s,\pi)}(\alpha_t(a)(s\inv\cdot \pi)) = 
\alpha_s(\alpha_t(a))(\pi) \\
&= \alpha_{st}(a)(\pi) = \beta_{(st,\pi)}(a((st)\inv\cdot\pi)) \\
&= \beta_{(s,\pi)(t,s\inv\cdot \pi)}(a(t\inv s\inv \cdot\pi)).
\end{align*}
It follows that $\beta_{(s\inv,s\inv\cdot \pi)}$ is an inverse to
$\beta_{(s,\pi)}$ and that each $\beta_{(s,\pi)}$ is an isomorphism.
Furthermore, a byproduct of the above computation is that we have 
shown $\beta$ is a homomorphism.  At this point, all
we need to do to show that $\beta$ is an action of $G$ on $A$ is show
that it is continuous.  

Suppose $(s_i,\pi_i)\rightarrow (s,\pi)$ and that $a_i\rightarrow a$
in $\mcal{A}$ such that $a_i\in A(s_i\inv\cdot\pi_i)$ for all $i$ and
$a\in A(s\inv\cdot\pi)$.  Choose $b\in A$ so that $b(s\inv\cdot\pi) =
a$.  First, observe that 
\[
p(\beta_{(s_i,\pi_i)}(a_i)) = \pi_i \rightarrow \pi
=p(\beta_{(s,\pi)}(a))
\]
by assumption.  Next, observe that 
\[
\beta_{(s,\pi)}(a) = \beta_{(s,\pi)}(b(s\inv\cdot\pi)) =
\alpha_s(a)(\pi)
\]
by definition.  Furthermore, $a_i - b(s_i\inv\cdot\pi_i)\rightarrow 0$
so that by Proposition \ref{prop:35} we have $\|a_i -
b(s_i\inv\cdot\pi_i)\|\rightarrow 0$.  However, since each
$\beta_{(s_i,\pi_i)}$ is an isomorphism, it follows that 
\[
\|\beta_{(s_i,\pi_i)}(a_i) - \alpha_{s_i}(b)(\pi_i)\| = 
\|\beta_{(s_i,\pi_i)}(a_i - b(s_i\inv\cdot\pi_i))\| =
\|a_i-b(s_i\inv\cdot\pi_i)\| \rightarrow 0.
\]
Because $\alpha$ is continuous with respect to the topology of
pointwise convergence we have $\alpha_{s_i}(b)\rightarrow \alpha_s(b)$
so that 
\[
\|\alpha_{s_i}(b)(\pi_i) -\alpha_s(b)(\pi_i)\| \leq
\|\alpha_{s_i}(b)-\alpha_s(b)\| \rightarrow 0.  
\]
Combining the previous two calculations we can conclude that given
$\epsilon > 0$ eventually 
\[
\|\beta_{(s_i,\pi_i)}(a_i) - \alpha_s(b)(\pi_i)\| < \epsilon.
\]
Clearly $\alpha_s(b)(\pi_i)\rightarrow \alpha_s(b)(\pi)$ so that it follows
from Proposition \ref{prop:35} that
$\beta_{(s_i,\pi_i)}(a_i)\rightarrow \beta_{(s,\pi)}(a)$. 

Therefore we have an action $\beta$ of $G$ on $A$ and as such we can
form the crossed product $A\rtimes_\beta G$ as the completion of
$\Gamma_c(G,r^*\mcal{A})$.   We claim that $A\rtimes_\beta G$ is isomorphic
to $A\rtimes_\alpha H$.  Given $f\in \Gamma_c(G,r^*\mcal{A})$ define 
\[
\Phi(f)(s)(\pi) = f(s,\pi)
\]
for all $s\in H$ and $\pi\in \widehat{A}$.  Fix $s\in H$.  It is clear
that $\Phi(f)(s)(\pi)\in A(\pi)$ for all $\pi\in\widehat{A}$ and that
$\Phi(f)(s)$ is continuous and compactly supported.  Therefore, after
identifying $A$ with $\Gamma_0(\widehat{A},\mcal{A})$, we get an
element $\Phi(f)(s)\in A$.  Suppose $s_i\rightarrow s$ and, to the
contrary, that $\Phi(f)(s_i)\not\rightarrow \Phi(f)(s)$ in
$\Gamma_0(\widehat{A},\mcal{A})$.  This implies that there exists
$\epsilon > 0$ and, after passing to a subnet and relabeling, 
$\pi_i\in \widehat{A}$
such that 
\begin{equation}
\label{eq:58}
\|f(s_i,\pi_i)-f(s,\pi_i)\|\geq \epsilon
\end{equation}
for all $i$.  Let $K$ be the projection of $\supp f$ onto
$\widehat{A}$.  If \eqref{eq:58} is to hold we must have $\pi_i\in K$
for all $i$.  However, $K$ is compact so that we may pass to another subnet,
relabel, and find $\pi\in K$ so that $\pi_i\rightarrow \pi$.  However, we
now have, using the continuity of $f$, 
\[
f(s_i,\pi_i)-f(s,\pi_i) \rightarrow f(s,\pi)-f(s,\pi)=0.
\]
It follows from Proposition \ref{prop:35} that this contradicts
\eqref{eq:58} so that we must have had 
$\Phi(f)(s_i)\rightarrow \Phi(f)(s)$.  Since the support of $\Phi(f)$
is contained in the projection of $\supp f$ to $H$ we have $\Phi(f)\in
C_c(H,A)$.  Next, suppose $f_i\rightarrow f$ with respect to
the inductive limit topology.  Let $K$ be the compact set which
eventually contains the support of the $f_i$.  Then the projection of
$K$ to $H$ eventually contains $\supp \Phi(f_i)$.  Furthermore, we
have, for all $g\in \Gamma_c(G,r^*\mcal{A})$, 
\begin{equation}
\|\Phi(g)\|_\infty = \sup_{s\in H} \|\Phi(g)(s)\| = \sup_{s\in
  H}\sup_{\pi\in\widehat{A}}\|\Phi(g)(s)(\pi)\| = 
\sup_{(s,\pi)\in G} \|g(s,\pi)\| = \|g\|_\infty.
\end{equation}
It follows that $\Phi(f_i)\rightarrow \Phi(f)$ with respect to the
inductive limit topology.  

Next, recall that the Haar system on $G$ is given by $\lambda^\pi =
\lambda\times \delta_\pi$ where $\lambda$ is Haar measure and
$\delta_\pi$ is the Dirac delta measure at $\pi$.  We compute for $f,g\in\Gamma_c(G,r^*\mcal{A})$, 
\begin{align*}
\Phi(f*g)(s)(\pi) &= f*g(s,\pi) = \int_G f(t,\rho)
\beta_{(t,\rho)}(g((t,\rho)\inv(s,\pi))) d(\lambda\times
\delta_\pi)(t,\rho) \\
&= \int_H f(t,\pi) \beta_{(t,\pi)}(g(t\inv s, t\inv\cdot \pi)) d\lambda(t)
\\ &= \int_H \Phi(f)(t)(\pi)\beta_{(t,\pi)}(\Phi(g)(t\inv
s)(t\inv\cdot \pi)) d\lambda(t) \\
&= \int_H (\Phi(f)(t)\alpha_t(\Phi(g)(t\inv s)))(\pi)d\lambda(t) \\
&= (\Phi(f)*\Phi(g))(s)(\pi).
\end{align*}
Observe that since ``evaluation at $\pi$'' is given
by the quotient map $A\rightarrow A(\pi)$ and since this map is bounded
and linear we may move it through the integral in the last
equality of the above calculation.  We also have
\begin{align*}
\Phi(f^*)(s)(\pi) &= f^*(s,\pi) = \beta_{(s,\pi)}(f((s\inv,s\inv\cdot \pi)^*)) \\
&= \beta_{(s,\pi)}(\Phi(f)(s\inv)(s\inv\cdot \pi))^* = 
(\alpha_s(\Phi(f)(s\inv))(\pi))^* \\
&= \alpha_s(\Phi(f)(s\inv)^*)(\pi) = \Phi(f)^*(s)(\pi).
\end{align*}
It follows that $\Phi$ is a $*$-homomorphism which is continuous in the
inductive limit topology.  Therefore Proposition \ref{prop:62} implies
that $\Phi$ is bounded with respect to the universal norm and extends
to a $*$-homomorphism from $A\rtimes_\beta G$ into $A\rtimes_\alpha
H$.  We need to show
that $\ran\Phi$ is dense in $C_c(H,A)$.  Let 
$D = \Gamma_c(\widehat{A},\mcal{A})$ be viewed as a subset of $A$ and
observe that $D$ is dense.  Consider the set of sums of elementary tensors
$C_c(H)\odot D \subset C_c(H,A)$.  By viewing $C_c(G,A)$ as a set of 
sections of a trivial bundle we can use Proposition
\ref{prop:42} to show that $C_c(H)\odot D$ is dense in $C_c(H,A)$.
First, $C_c(H)\odot D$ is clearly closed under the action of
$C_0(H)$. Furthermore, given $s\in H$ choose $\phi\in C_c(H)$ such
that $\phi(s) = 1$.  Then $\phi\otimes d(s) = d$ for all $d\in D$ and
it clear that $C_c(H)\odot D$ is ``fibrewise dense.''  It follows that
$C_c(H)\odot D$ is dense in $C_c(H,A)$ with respect to the uniform
norm.  However, if $\sum_j\phi^j_i\otimes d^j_i\rightarrow f$ then we can choose
$\psi$ which is one on $\supp f$ and zero off a neighborhood of $\supp
f$.  We then have $\sum_j\psi\phi_i^j \otimes d_i^j \rightarrow f$ with respect
to the inductive limit topology.  Alternatively we could just
cite \cite[Lemma 1.87]{tfb2} and skip the previous argument.  
Now, given $\phi\otimes a\in C_c(H)\odot
D$ we define $f(s,\pi)= \phi(s)a(\pi)$ and we have chosen $\phi$ and
$a$ so that $f\in \Gamma_c(G,r^*\mcal{A})$.  Furthermore we clearly
have $\Phi(f)(s)(\pi) = \phi(s)a(\pi) = \phi\otimes a(s)(\pi)$.  It
follows that $C_c(H)\odot D\subset \ran\Phi$ and therefore $\ran\Phi$
is dense in 
$C_c(H,A)$ with respect to the inductive limit topology,
with respect to the $I$-norm, and with respect to the universal norm.

We would like to show that $\Phi$ extends to an isomorphism.  First,
observe that $\Phi$ is clearly injective on $\Gamma_c(G,r^*\mcal{A})$
so that we can define an inverse map $\Psi=\Phi\inv:\ran \Phi
\rightarrow \Gamma_c(G,r^*\mcal{A})$ given for $f\in \ran\Phi$ by
$\Psi f(s,\pi) = f(s)(\pi)$. It is
straightforward to use the fact that $\Phi$ is an injective
$*$-homomorphism to show that that $\ran\Phi$ is a $*$-subalgebra of
$A\rtimes_\alpha H$ and that $\Psi$ is a $*$-homomorphism.
Furthermore, given $f\in\ran\Phi$ we calculate
\begin{align*}
\int_G \|\Psi(f)(s,\rho)\| d(\lambda\times\delta_\pi)(s,\rho) &=
\int_H \|f(s)(\pi)\| d\lambda(s) \\
&\leq \int_H \|f(s)\| d\lambda(s) \leq \|f\|_I.
\end{align*}
It is a similar task to show that 
\[
\int_G \|\Psi(f)((s,\rho)\inv) d(\lambda\times\delta_\pi)(s,\rho)
\leq \|f\|_I
\]
and it follows that $\Psi$ is $I$-norm decreasing.\footnote{Since we
  are viewing group crossed products as special cases of groupoid
  crossed products we use the $I$-norm on $C_c(G,A)$ and not the $L^1$-norm.}
In particular, we have for all $f\in \ran\Phi$
\[
\|\Psi(f)\| \leq \|f\|_I.
\]
Since $\ran\Phi$ is dense in $C_c(H,A)$ with respect to the $I$-norm
we can extend $\Psi$ from $\ran\Phi$ to all of $C_c(H,A)$.  Notice
that $\Psi$ now maps into the crossed product $A\rtimes_\beta G$.
Furthermore, since the operations are $I$-norm continuous and $\Psi$
is a $*$-homomorphism on $\ran\Phi$ it follows that $\Psi$ is a
$I$-norm decreasing $*$-homomorphism on all of $C_c(H,A)$.  
It follows from Proposition \ref{prop:62} that 
$\Psi$ is bounded with respect to the universal norm and we can
extend it to a $*$-homomorphism from $A\rtimes_\alpha H$ into
$A\rtimes_\beta G$.  Furthermore, since $\Psi$ and  $\Phi$ are
inverses on a dense subset, they must be inverses everywhere.  It
follows that $\Phi$ is an isomorphism and $A\rtimes_\alpha H$ and
$A\rtimes_\beta G$ are isomorphic.  
\end{example}

\begin{remark}
We used some tricky arguments to show that the two crossed products
in Example \ref{ex:22} are isomorphic.  It is possible, and a useful
exercise, to show that $\Phi$ is an isomorphism by proving that given a
covariant representation $(\pi,U)$ of $(A,H,\alpha)$ then
$\pi\rtimes U\circ \Phi$ is equivalent to a covariant representation
of $(A,G,\alpha)$ and vice-versa.  To construct the covariant
representation of $(A,G,\alpha)$ use $\pi$ as the representation
of $A$.  To get a representation of $G$, view $C^*(G)$ as the
crossed product $C_0(\widehat{A})\rtimes_{\lt} H$.  Then extend
$\pi$ to a representation of $C_0(\widehat{A})$, sitting inside the
multiplier algebra of $A$, and form the covariant representation
$\overline{\pi}\rtimes U$ of $C_0(\widehat{A})\rtimes_{\lt} H\cong C^*(G)$.
Renault's Decomposition theorem then gives the desired
representation of $G$.  In order to go the other direction,
basically perform this process in reverse, obtaining the
representation of $H$ as part of a covariant decomposition of the
representation of $C_0(\widehat{A})\rtimes_{\lt} H$ coming from the
representation of $G$.  Of course, there are a lot of technicalities to
work through and this could just be taken as another example of why
it's preferable not to work directly with covariant representations.  
\end{remark}

We also show that Example \ref{ex:22} has a nice converse, further
strengthening the notion that Example \ref{ex:22} provides an
alternate method for viewing groupoid crossed products as generalizing
the group case.  

\begin{prop}
\label{prop:98}
Suppose $A$ has Hausdorff spectrum $X$ and that $(H,X)$ is a
transformation group. Let $G=H\ltimes X$ be the associated
transformation groupoid.  If $\beta$ is an action of $G$ on $A$ then 
\[
\alpha_s(a)(\pi) = \beta_{(s,\pi)}(a(s\inv\cdot \pi))
\]
defines an action of $H$ on $A$.  Furthermore $A\rtimes_\alpha H$ and
$A\rtimes_\beta G$ are isomorphic. 
\end{prop}

\begin{proof}
We view $A$ as sections of the associated bundle $\mcal{A}$.  If $a\in
A$ then it is clear enough that $\alpha_s(a)$ defines a continuous section
of $\mcal{A}$.  Now, if $\epsilon > 0$ then 
\[
\{\pi:\|\alpha_s(a)(\pi)\|\geq \epsilon \} = \{\pi:\|a(s\inv\cdot
\pi)\|\geq \epsilon\} = s\cdot \{\pi:\|a(\pi)\|\geq \epsilon\}.
\]
Since $\{\pi:\|a(\pi)\|\geq\epsilon\}$ is compact, it follows that
$\{\pi:\|\alpha_s(a)(\pi)\|\geq \epsilon\}$ is compact as well.  Thus
$\alpha_s(a)$ vanishes at infinity and $\alpha_s(a)\in A$.
Furthermore since each $\beta$ is a $*$-homomorphism it is
straightforward to show that $\alpha_s$ is a $*$-homomorphism.  Next
observe that if $e$ is the unit of $H$ then  
\[
\alpha_e(a)(\pi) = \beta_{(e,\pi)}(a(\pi)) = a(\pi).
\]
Thus $\alpha_e = \id$.  Furthermore for $s,t\in H$ we have 
\begin{align*}
\alpha_{st}(a)(\pi) &= \beta_{(st,\pi)}(a(t\inv s\inv\cdot \pi)) = 
\beta_{(s,\pi)}(\beta_{(t,s\inv\cdot \pi)}(a(t\inv s\inv\cdot \pi)))
\\
&= \beta_{(s,\pi)}(\alpha_t(a)(s\inv\cdot \pi)) =
\alpha_s(\alpha_t(a))(\pi).
\end{align*}
This shows that $\alpha_s$ is an automorphism of $A$ and that
$\alpha:H\rightarrow \Aut(A)$ is a homomorphism.  Finally, fix $a\in
A$, let $s_i\rightarrow s$ in $H$ and suppose $\alpha_{s_i}(a)$ does
not converge to $\alpha_s(a)$.   Then, by passing to a subsequence,
there exists $\epsilon > 0$ and $\pi_i \in X$ such that 
\begin{equation}
\label{eq:123}
\| \alpha_{s_i}(a)(\pi_i) - \alpha_s(a)(\pi_i) \| = 
\| \beta_{(s_i,\pi_i)}(a(s_i\inv\cdot \pi_i)) -
\beta_{(s,\pi_i)}(a(s\inv\cdot \pi_i))\| \geq \epsilon 
\end{equation}
for all $i$.  Consider the compact set $K = \{\rho : \|a(\rho)\| \geq
\epsilon/2 \}$.  If both $s_i\inv\cdot \pi_i$ and $s\inv\cdot \pi_i$ are
not in $K$ then 
\[
\| \alpha_{s_i}(a)(\pi_i) - \alpha_s(a)(\pi_i)\| \leq
\|a(s_i\inv\cdot\pi_i)\| + \|a(s\inv\cdot \pi_i)\| < \epsilon
\]
which contradicts \eqref{eq:123}.  There are two cases to consider.
If $s\inv \cdot \pi_i \in K$ infinitely often then we can pass to a
subsequence, multiply by $s$, and find $\pi\in X$ such that
$\pi_i\rightarrow \pi$.  In the other case $s_i\inv\cdot \pi_i$ is in
$K$ infinitely often.  Therefore we can pass to a subsequence,
multiply by the sequence $s_i\rightarrow s$, and find $\pi\in X$ such that
$\pi_i\rightarrow \pi$.  In either case we have
\begin{align*}
\beta_{(s,\pi_i)}(a(s\inv\cdot \pi_i))&\rightarrow
\beta_{(s,\pi)}(a(s\inv\cdot \pi)),\quad\text{and} \\
\beta_{(s_i,\pi_i)}(a(s_i\inv\cdot \pi_i))&\rightarrow 
\beta_{(s,\pi)}(a(s\inv\cdot \pi)).
\end{align*}
It follows that 
\[
\| \alpha_{s_i}(a)(\pi_i)-\alpha_s(a)(\pi_i)\|\rightarrow 0
\]
which contradicts \eqref{eq:123}.  Thus $\alpha$ is a strongly
continuous map and is an action of $H$ on $A$.  

Next, given $\pi\in X$ consider its factorization $\pi'$ to $A(\pi)$.  Then
\[
s\cdot\pi(a) = (s,s\cdot \pi)\cdot \pi(a) = 
\pi'(\beta_{(s\inv, \pi)}(a(s\cdot \pi))) 
= \pi'(\alpha_s\inv(a)(\pi)) 
= \pi(\alpha_s\inv(a))
\]
In other words, the action of $G$ on $X$ induced by $\alpha$ is exactly
the original action of $H$ on $X$ so that $G$ is the transformation
groupoid associated to both actions.  Furthermore, as in Example
\ref{ex:22}, there is an action $\beta'$ of $G$ on $A$ induced by $\alpha$.
However it is easy to see that $\beta' = \beta$.  Thus it follows from
Example \ref{ex:22} that $A\rtimes_\beta G$ is isomorphic to
$A\rtimes_\alpha H$.  
\end{proof}

\begin{remark}
\label{rem:21}
As was mentioned before, Example \ref{ex:22} and Proposition
\ref{prop:98} are an effective way of viewing groupoid crossed
products as generalizing group crossed products.  However, since we
require the algebra to have Hausdorff spectrum, this is only a partial
generalization.  In order to make this work for general $C^*$-algebras
we would need to be able to work with the transformation groupoid
associated to the action of the group on the primitive ideal space.
This would require us to expand our notion of groupoid dynamical system to
include non-Hausdorff groupoids with non-Hausdorff unit spaces that act
on bundles over non-Hausdorff spaces, and that seems difficult. 
\end{remark}


\section{Groupoid Algebras}
\label{sec:scalarprod}

In this section we explore another important special case of the
groupoid crossed product.  We start with the following observation. 

\begin{prop}
Suppose $X$ is a locally compact Hausdorff space.
Then $C_0(X)$ is a $C_0(X)$-algebra with
$C_0(X)(x)\cong \C$ for all $x\in X$.  Furthermore, the
bundle associated to $C_0(X)$ is (isomorphic to) $X\times \C$. 
\end{prop}

\begin{proof}
It is straightforward to show 
that $C_0(X)$ is a $C_0(X)$-algebra with the
action given by left multiplication.  Let $I_x$ be the ideal such that
$C_0(X)(x) = C_0(X)/I_x$ and let $J_x$ be the ideal of all
functions which vanish at $x$.  It is clear that $I_x \subset J_x$.
Let $e_l$ be an approximate unit for
$C_0(X)$.  Then given $f\in J_x$ we have, by definition, 
 $f\cdot e_l \in I_x$ for all $l$.  It follows that $f\in
I_x$. Thus $I_x = J_x$ and it is now straightforward to show that
the map $f\mapsto f(x)$ factors to an isomorphism of 
$C_0(X)(x)$ with $\C$. Let $\mcal{C}$ denote the bundle
associated to $C_0(X)$.  
We would like to show that $\mcal{C}\cong X\times \C$.  It
follows from Corollary \ref{cor:2} that it suffices to show that 
there is a $C_0(X)$-linear
isomorphism from $C_0(X)$ onto $\Gamma_0(X,X\times \C)$.  However it
is clear that $\Gamma_0(X,X\times \C)$ can be identified with 
$C_0(X)$ so that the desired
isomorphism follows. 
\end{proof}

\begin{prop}
\index{dynamical system}
Let $G$ be a second countable, locally compact Hausdorff groupoid with
a Haar system.
Then there is an action of $G$ on $C_0(G\unit)$ given by the
collection of functions
$\id_\gamma:C_0(G\unit)(s(\gamma))\rightarrow C_0(G\unit)(r(\gamma))$
such that 
\[
\id_\gamma(z) = z
\]
for all $z\in\C$ and $\gamma\in G$.  
\end{prop}

\begin{proof}
First observe that the bundle associated to $C_0(G\unit)$ as a
$C_0(G\unit)$-algebra is $G\unit\times \C$.  
Let us be a bit more formal and define
$\id_\gamma:C_0(G\unit)(s(\gamma))\rightarrow C_0(G\unit)(r(\gamma))$
by $\id_\gamma(s(\gamma),z) = (r(\gamma),z)$.
Observe that $\id_\gamma$ is an
isomorphism.  Furthermore, it is clear that $\id_{\gamma\eta} =
\id_\gamma\circ\id_\eta$ if $s(\gamma)=r(\eta)$.  The only thing left
to check is the continuity condition and this is obvious.  
\end{proof}

Thus, given any groupoid we have a natural dynamical system
associated to that groupoid.  We make the following

\begin{definition}
\index{groupoid $C^*$-algebra}
\index[not]{$C^*(G)$}
\label{def:46}
Let $G$ be a second countable locally compact Hausdorff groupoid
$G$ with Haar system $\lambda$.  We define the {\em groupoid
  $C^*$-algebra} to be the crossed product $C_0(G)\rtimes_{\id} G$ and
use the notation $C^*(G)$.
\end{definition}

\begin{remark}
While the notation doesn't reflect this, the groupoid
algebra $C^*(G)$ depends on the choice of Haar
system. When it matters the notation $C^*(G,\lambda)$ is often used.
Whether or not the groupoid $C^*$-algebra is, up to isomorphism,
independent of the choice of the Haar system is an open question.  On the
other hand, it is an immediate corollary to
Theorem \ref{thm:groupoidequiv} that up to Morita equivalence 
$C^*(G)$ is independent of the choice of Haar system. 
\end{remark} 

Of course, since the groupoid $C^*$-algebra is just a crossed product
we immediately recover the following from Proposition \ref{prop:58}. 

\begin{corr}
\label{cor:23}
Suppose $G$ is a second countable, locally compact Hausdorff groupoid
with a Haar system.  Then the
universal norm is a norm and is dominated by the $I$-norm.
Furthermore, the completion $C^*(G)$ is a separable
$C^*$-algebra.  
\end{corr}

The following remark shows how the groupoid $C^*$-algebra is an
extension of the group $C^*$-algebra.  

\begin{remark}
\label{rem:17}
Let $G$ be a second countable locally compact Hausdorff 
group.  Then the the ``groupoid'' $C^*$-algebra
associated to $G$ is just $C^*(G) = \C\rtimes_{\id} G$.  However, it
follows from Section \ref{sec:groupprod} that this is just a
group crossed product.  We then cite \cite[Example 2.33]{tfb2} to see
that this is exactly the group $C^*$-algebra associated to $G$.  
Because the group $C^*$-algebra is obtained in this way as a group crossed
product with $\C$, we often refer to 
group $C^*$-algebras as the
``scalar'' case of group crossed products.  The same terminology is
used for groupoids.  This is further justified by the fact that
$C_0(G\unit)$ is the simplest $C_0(G\unit)$-algebra possible and that
the fibres are all isomorphic to $\C$. 
\end{remark}

We make the
following straightforward observations. 

\begin{prop}
\label{prop:67}
Given a separable locally compact Hausdorff groupoid $G$ with Haar
system $\lambda$, the section algebra $\Gamma_c(G,G\unit\times \C)$ can
be identified with $C_c(G)$.  Furthermore the operations on $C_c(G)$ become 
\[
f*g(\gamma) = \int_G f(\eta)g(\eta\inv \gamma)d\lambda^{r(\gamma)}(\eta)
\quad\text{and}\quad f^*(\gamma) = \overline{f(\gamma\inv)}
\]
and the $I$-norm reduces to 
\[
\|f\|_I = \max\left\{\sup_{u\in G\unit} \int_G
  |f(\gamma)|d\lambda^u(\gamma), \sup_{u\in G\unit} \int_G |f(\gamma)|
  d\lambda_u(\gamma)\right\}.
\]
\end{prop}

\begin{proof}
In light of the identification made in Remark \ref{rem:8} it is clear
that any section $f\in \Gamma_c(G,r^*(G\unit\times \C))$ is of the
form $f(\gamma) = (r(\gamma),\tilde{f}(\gamma))$ where $\tilde{f}\in
C_c(G)$.  Thus by identifying $\tilde{f}$ with $f$ we can identify
$C_c(G)$ with $\Gamma_c(G,r^*(G\unit\times \C))$.  The rest of the
claims made in the proposition follow immediately once you recall
that the action of $G$ on $C_0(G\unit)$ is given by the identity
map. 
\end{proof}

We continue by considering the covariant
representations of $(C_0(G\unit),G,\id)$.  

\begin{prop}
\label{prop:66}
\index{covariant representation}
\index{groupoid representation}
Let $G$ be a second countable locally compact Hausdorff groupoid with
a Haar system
and $(G\unit*\mfrk{H},\mu,U)$ a unitary representation of $G$.  Then the
representation $L_\mu:C_0(G\unit)\rightarrow L^2(G\unit*\mfrk{H},\mu)$
such that $L_\mu(f)\phi(u)= f(u)\phi(u)$ for all $f\in C_0(G\unit)$
and $\phi\in\mcal{L}^2(G\unit*\mfrk{H},\mu)$ is $C_0(G\unit)$-linear and
$(G\unit*\mfrk{H},\mu,L_\mu,U)$ is a covariant representation of
$(C_0(G\unit),G,\id)$.   Furthermore, every covariant representation
of $(C_0(G\unit),G,\id)$ is of this form. 

What's more, given a unitary representation $U$ of $G$ as above,
the integrated form of $(G\unit*\mfrk{H},\mu,L_\mu,U)$, 
also denoted by $U$, is given by 
\begin{equation}
\label{eq:100}
U(f)h(\gamma) = \int_G f(\gamma)U_\gamma h(s(\gamma))
\Delta(\gamma)\neghalf d\lambda^u(\gamma)
\end{equation}
for $f\in C_c(G)$ and $h\in \mcal{L}^2(G\unit*\mfrk{H},\mu)$.
Finally, \eqref{eq:100} is characterized by 
\begin{equation}
\label{eq:101}
(U(f)h,k) = \int_G (f(\gamma)U_\gamma
h(s(\gamma)),k(r(\gamma)))\Delta(\gamma)\neghalf d\nu(\gamma)
\end{equation}
where $\nu$ is the measure induced by $\mu$. 
\end{prop}

\begin{proof}
Since $L_\mu(f)$ is just a multiplication operator it is
straightforward to show that $L_\mu$ is a representation of
$C_0(G\unit)$.  Let $\{e_l\}$ be a special orthogonal fundamental
sequence for $X*\mfrk{H}$ and view $e_l$ as an element of
$L^2(X*\mfrk{H},\mu)$.  Given
$\phi\in L^2(G\unit*\mfrk{H},\mu)$ if $(L_\mu(f)e_l,\phi) = 0$ for all
$l$ and $f\in C_c(G\unit)^+$ then 
\[
\int_{G\unit} f(u) (e_l(u),\phi(u)) d\mu(u) = 0
\]
for all $l$ and $f\in C_c(G\unit)^+$. However, it follows that
$(e_l(u),\phi(u))\ne 0$ almost everywhere.  We can choose a single
null set $N$ such that this is true for all $l$ given $u\not\in N$.
Since $\{e_l(u)\}$ is a basis for $\mcal{H}(u)$ it follows that
$\phi(u) =0$ almost everywhere.  This is enough to show
that $L_\mu$ is nondegenerate.  Finally, since $L_\mu(f)$ is exactly
the diagonal operator associated to $f$, it is clear that $L_\mu$ is
$C_0(G\unit)$-linear.  Consider its associated decomposition 
\[
L_\mu = \int^\oplus_{G\unit} \pi_u d\mu(u).
\]
Then, at least for $\mu$-almost every $\mu$, we know that $\pi_u$ is a
nondegenerate representation of $C_0(G\unit)(u) = \C$.  However, any
nondegenerate representation of $\C$ must map $1$ onto the identity
operator $\bf{1}$.  Thus, we have $\pi_u(u,z) = z
\bf{1}$ almost everywhere. 
We may as well change the decomposition of $L_\mu$ on a
null set and assume this is always true.  However, we now have, given
$z\in \C$ and $\gamma\in G$, 
\[
U_\gamma \pi_{s(\gamma)}(s(\gamma),z) = z U_\gamma
= \pi_{r(\gamma)}(\id_\gamma(s(\gamma),z)) U_\gamma.
\]
Thus $(L_\mu,U)$ is a covariant representation of
$(C_0(G\unit),G,\id)$.  

Now suppose $(\mu,G\unit*\mfrk{H},\pi,U)$ is a covariant representation of
$(C_0(G\unit),G,\id)$.  Then, by definition, $(\mu,G\unit*\mfrk{H},U)$
is a unitary representation of $G$.  All we have to do is show that
$\pi=L_\mu$.  Find a decomposition $\pi = \int^\oplus_{G\unit}\pi_ud\mu(u)$
and observe that $\pi_u$ is nondegenerate almost everywhere.  Then as
before we must have $\pi_u(u,z) = z\bf{1}$ for $\mu$-almost
all $u$.  It follows that, given $f\in C_0(X)$ and $\phi\in
\mcal{L}^2(G\unit*\mfrk{H},\mu)$,
\[
\pi(f)\phi(u) = \pi_u(f(u))\phi(u) = f(u)\phi(u) = L_\mu(f)\phi(u)
\]
$\mu$-almost everywhere.  Hence $\pi(f) = L_\mu(f)$ and we
are done.  

For the last part of the theorem suppose $U$ is a unitary
representation and $L_\mu$ is the associated representation of $C_0(G
\unit)$ with $L_\mu = \int_{G\unit}^\oplus \pi_u d\mu(u)$ as above so
that $\pi_u(z) = z \bf{1}$ for all $z\in \C$ and $u\in
G\unit$.  Then the integrated representation $L_\mu\rtimes U$, denoted
by $U$ in the statement of the proposition, is given by 
\[
U(f)\phi(u) = \int_G \pi_u(f(\gamma)) U_\gamma h(s(\gamma))
\Delta(\gamma)\neghalf \lambda^u(\gamma) = 
\int_G f(\gamma) U_\gamma h(s(\gamma))\Delta(\gamma)\neghalf
\lambda^u(\gamma).
\] 
The final part of the proposition follows from the fact that that
$L_\mu\rtimes U$ is characterized by \eqref{eq:38}.
\end{proof}

The last thing we will do is restate
Theorem \ref{thm:disintigration} for groupoids.  This ``special case''
is actually Renault's Proposition 4.2 and is used to prove the more
general Theorem \ref{thm:disintigration} \cite[Section
7]{renaultequiv}.  

\begin{theorem}[Renault's Disintegration Theorem]
\label{thm:scalardis}
\index{disintegration theorem}
Suppose that $\mcal{H}_0$ is a dense subspace of a complex Hilbert
space $\mcal{H}$ and that $u$ is a homomorphism from $C_c(G)$ into
the algebra of linear maps on $\mcal{H}_0$ such that 
\begin{enumerate}
\item $\{u(f)h:f\in C_c(G), h\in \mcal{H}_0\}$ is dense in $\mcal{H}$,
\item for each $h,k\in\mcal{H}_0$, 
\[
f\mapsto (u(f)h,k)
\]
is continuous in the inductive limit topology, and 
\item for each $f\in C_c(G)$ and all $h,k\in \mcal{H}_0$
\[
(u(f)h,k)=(h,u(f^*)k).
\]
\end{enumerate}
Then $u(f)$ is bounded and extends to a bounded operator $U(f)$ on
$\mcal{H}$ such that $U$ is an $I$-norm decreasing $*$-representation
of $C_c(G)$.  Furthermore, there is a unitary representation
$(\mu,G\unit*\mfrk{H}, V)$ such that $U$ is equivalent to the
integrated form of $V$.  
\end{theorem}

\begin{proof}
This follows immediately from Theorem \ref{thm:disintigration} once
you identify unitary representations of $G$ with covariant
representations of $(C_0(G\unit),G,\alpha)$ via Proposition
\ref{prop:66}. 
\end{proof}

From here we obtain the following corollary, which is really just a
restatement of Proposition \ref{prop:61}.   

\begin{corr}
\label{cor:5}
Suppose $G$ is a second countable locally compact Hausdorff groupoid
with a Haar system.  Then
the universal norm on $C_c(G)$ is also given by
\begin{align*}
\|f\| &= \sup\left\{\|U(f)\| : \begin{array}{l} 
 \text{$U$ is an $I$-norm decreasing} \\ 
\text{$*$-representation of
  $C_c(G)$}\end{array}\right\}, \\
\|f\| &= \sup\left\{\|U(f)\| : \begin{array}{l}
\text{$U$ is a $*$-representation of
  $C_c(G)$ which is} \\ 
\text{continuous in the inductive limit topology} \end{array}\right\}.
\end{align*}
It follows that any (possibly degenerate)
$*$-representation of $C_c(G)$
which is either $I$-norm decreasing or continuous with respect to the
inductive limit topology is bounded with respect to the universal norm
and extends to a representation of $C^*(G)$.  
\end{corr}

Now it's time for a simple example. 

\begin{example}
\label{ex:23}
Suppose $X$ is a second countable locally compact Hausdorff space and
we view $X$ as a ``cotrivial'' groupoid as in Example \ref{ex:2}.
Then $X$ has a Haar system given by the Dirac delta measures $\delta_x$
as in Example \ref{ex:24}.  Thus we can form the groupoid algebra
$C^*(X)$.  However, since integration against $\delta_x$ is
evaluation at $x$ it is easy to see that $\|f\|_I = \|f\|_\infty$.
Thus the $I$-norm is a $C^*$-norm and the completion of $C_c(X)$ with
respect to the $I$-norm is clearly $C_0(X)$.  By choosing a faithful
representation $\pi$ of $C_0(X)$ it follows from Corollary \ref{cor:5}
that 
\begin{align*}
\|f\| &= \sup\{\|\pi(f)\|:\text{$\pi$ is $I$-norm decreasing}\} \\
&\geq \|\pi(f)\| = \|f\|_I = \|f\|_\infty.
\end{align*}
It follows that in this case the universal norm is the uniform norm
and that the groupoid algebra $C^*(X)$ is just $C_0(X)$.  
\end{example}

The point of all of this is twofold.  First, groupoid $C^*$-algebras
are much easier to work with than groupoid crossed products because
the associated functions are scalar valued and it will
be useful to have the simplified operations and representation theory
written down.  The second is the following remark. 

\begin{remark}
\label{rem:31}
Groupoid $C^*$-algebras have been around for much longer than groupoid
crossed products have.  It is important to show that viewing groupoid
$C^*$-algebras as special cases of crossed products gets us back to
the older definitions.  In \cite[Chapter II]{groupoidapproach} and
\cite[Chapter 3]{coords} the groupoid $C^*$-algebra is defined to be a
universal completion of $C_c(G)$.  In the first reference we complete
with respect the norm 
\begin{equation}
\label{eq:153}
\|f\| = \sup\{\|U(f)\|:\text{$U$ is an $I$-norm decreasing
  representation}\}
\end{equation}
and in the second we complete with respect to the norm 
\begin{equation}
\label{eq:154}
\|f\| = \sup\{\|U(f)\|:\text{$U$ is a unitary representation}\}.
\end{equation}
However, it is clear once one sorts through the references
that the convolution algebra $C_c(G)$
is the same in all three cases.  In light of Proposition \ref{prop:66},
\eqref{eq:154} is exactly the universal norm and
Corollary \ref{cor:5} shows it is equal to \eqref{eq:153}.  Thus,
in both cases the completion of $C_c(G)$ gives us the same algebra as 
Definition \ref{def:46}.  
\end{remark}

Moving on, we would like to make a connection between groupoid
isomorphisms and groupoid algebras.  
The following proposition is useful and also entirely unsurprising.  

\begin{prop}
\label{prop:70}
\index{groupoid homomorphism}
Suppose $G$ and $H$ are second countable locally compact Hausdorff
groupoids with Haar systems $\lambda$ and $\beta$ respectively.  If
$\phi:G\rightarrow H$ is a groupoid isomorphism such that 
$\phi_* \lambda^u = \beta^{\phi(u)}$, that is, such that 
\[
\int_G f(\phi(\gamma)) d\lambda^u(\gamma) = \int_H
f(\gamma)d\beta^{\phi(u)}(\gamma)
\]
for all $f\in C_c(H)$, then the map
$\Phi:C_c(H)\rightarrow C_c(G)$ defined by $\Phi(f) = f\circ \phi$
extends to an isomorphism of $C^*(H)$ onto $C^*(G)$.  
\end{prop}

\begin{remark}
The fact that $\phi$ has to intertwine the Haar systems is an
annoyance.  If $G$ has a Haar system and $H$ is isomorphic to $G$ then
$H$ has a Haar system induced by the isomorphism.  However, it
has not been shown that the groupoid $C^*$-algebra is independent of
the choice of Haar measure so that we cannot assume without loss of
generality that any two Haar systems on $G$ and $H$ are intertwined.  
It's also annoying that this proposition does not immediately extend to
homomorphisms.  The problem is that if $\phi$ is simply a continuous groupoid
homomorphism then there is no reason that $f\circ \phi$ should be compactly
supported for a given $f\in C_c(G)$.  
\end{remark}

\begin{proof}
First, 
$f\circ \phi$ is clearly a continuous compactly supported function so
that $\Phi(f)$ is well defined.  It is easy to see that this map is
linear and we can check that 
\begin{align*}
\Phi(f*g)(\gamma) &= f*g(\phi(\gamma)) = 
\int_G f(\eta) g(\eta\inv\phi(\gamma)) d\beta^{r(\phi(\gamma))}(\eta)
\\
&= \int_G f(\phi(\eta)) g(\phi(\eta)\inv\phi(\gamma))
d\lambda^{r(\gamma)}(\eta) \\
&= \Phi(f)*\Phi(g)(\gamma)
\end{align*}
and
\begin{align*}
\Phi(f^*)(\gamma) &= f^*(\phi(\gamma)) =
\overline{f(\phi(\gamma)\inv)} 
= \overline{\Phi(f)(\gamma\inv)} = \Phi(f)^*(\gamma).
\end{align*}
Thus $\Phi$ is a $*$-homomorphism.  Since $\phi$ is a homeomorphism it
is easy to show that $\Phi$ is continuous in the inductive limit
topology.  It follows from Corollary \ref{cor:24} that
$\Phi$ extends to a homomorphism from $C^*(H)$ onto $C^*(G)$.
However, applying the same argument to $\phi\inv$ gives us a homomorphism
from $C^*(G)$ onto $C^*(H)$ which is clearly an inverse for $\Phi$.
Thus $C^*(G)$ and $C^*(H)$ are isomorphic.
\end{proof}

We end this section by citing a useful and interesting theorem which
makes use of the fact that, as stated in Remark \ref{rem:31}, 
viewing groupoid algebras
as crossed products brings you back to the more classical object. 

\begin{remark}
We assume that the reader is familiar with Morita equivalence,
imprimitivity bimodules, pre-Hilbert $A$-spaces and the like.  A
good reference for this material is \cite[Chapters 1,2,3]{tfb}.
\end{remark}

\begin{theorem}[Groupoid Equivalence Theorem 
{\cite[Theorem 2.8]{groupoidequiv}}]
\index{G,H-equivalence@$(G,H)$-equivalence}
\label{thm:groupoidequiv}
Suppose $G$ and $H$ are second countable locally compact Hausdorff
groupoids with Haar systems $\lambda$ and $\beta$.  Then for any
$(G,H)$-equivalence $X$, $\mcal{Z}_0 = C_c(X)$ is a
pre-$C^*(G)-C^*(H)$-imprimitivity bimodule with operations defined for
$f\in C_c(G)$, $g\in C_c(H)$, and $\phi,\psi\in C_c(X)$ by 
\begin{align}
f\cdot\phi(x) &= \int_G f(\gamma)\phi(\gamma\inv\cdot
x)d\lambda^{r(x)}(\gamma) \\
\phi\cdot g(x) &= \int_H
\label{eq:66}
\phi(x\cdot\eta)g(\eta\inv)d\beta^{s(x)}(\eta) \\
\langle \phi,\psi\rangle_{C^*(H)}(\eta) &= \int_G
\overline{\phi(\gamma\inv\cdot x)}\psi(\gamma\inv\cdot x\cdot
\eta)d\lambda^{r(x)}(\gamma) \\
\label{eq:67}
\lset{C^*(G)}\langle \phi,\psi\rangle(\gamma) &= 
\int_H \phi(\gamma\inv\cdot x\cdot \eta)\overline{\psi(x\cdot \eta)}
d\beta^{s(x)}(\eta)
\end{align}
where $x$ in \eqref{eq:66} is any element of $X$ such that $s(x) =
r(\eta)$ and $x$ in \eqref{eq:67} is any element of $X$ such that
$r(x) = r(\gamma)$.  In particular, the completion $\mcal{Z}$ of
$\mcal{Z}_0$ is a $C^*(G)-C^*(H)$-imprimitivity bimodule and
$C^*(G)$ and $C^*(H)$ are Morita
equivalent. 
\end{theorem} 

\begin{proof}[Remark]
This isn't so much a proof as it is an explanation of how this
statement of the Groupoid Equivalence Theorem can be obtained from the
statement in the reference.  Those readers not interested in chasing
down citations can skip ahead.  
Except for the explicit description of the bimodule operations, Theorem
\ref{thm:groupoidequiv} is as stated in \cite[Theorem 2.8]{groupoidequiv}.  The
operations themselves can be deduced by scanning through the two pages
following the statement of the theorem itself. 
\end{proof}

This equivalence theorem is quite nice and is extended to groupoid
crossed products in Section \ref{sec:transitive}.  


\section{Group Bundles}
\label{sec:bundleprod}
In this section we present another special case of groupoid crossed
products which will play an important role in Sections
\ref{sec:locally-unitary} and \ref{sec:crossedstab}.  In
particular, we will be interested in what happens when a groupoid
group bundle $S$ acts on a $C^*$-algebra. Suppose $S$ is a group bundle
over $S\unit$, $A$ is a $C_0(S\unit)$-algebra and $\mcal{A}$ is the
upper-semicontinuous bundle associated to $A$.  Now suppose
$\alpha$ is an action of $S$ on $A$ and consider the restriction of
$\alpha$ to $S_u$ for $u\in S\unit$.  Given $t\in S_u$, since
$r(t)=s(t) = u$, we have
$\alpha_t:A(u)\rightarrow A(u)$.  It follows from the fact
that $\alpha$ preserves the groupoid operation that
$\alpha|_{S_u}:S_u\rightarrow \Aut(A(u))$ is an automorphism. 
Furthermore the continuity
condition on $\alpha$ implies that  $\alpha|_{S_u}$ is continuous onto
the topology of pointwise convergence.  Putting this all together we
have the following

\begin{prop}
Suppose $(A,S,\alpha)$ is a separable groupoid dynamical system and
$S$ is a group bundle.  Then $(A(u),S_u,\alpha|_{S_u})$ is a group
dynamical system for each $u\in G\unit$.
\end{prop}

\begin{proof}
This follows from the above discussion and Proposition \ref{prop:63}.
\end{proof}

At this point we have a bundle of groups acting on a bundle of 
$C^*$-algebras and it
is not surprising that the resulting crossed product is a bundle of
group crossed products. However, this proposition is not as easy as it
looks and takes some effort to prove.  

\begin{prop}
\label{prop:65}
\index{crossed product!group}
Suppose $(A,S,\alpha)$ is a separable groupoid dynamical system and
$S$ is a group bundle.  Then $A\rtimes_\alpha S$ is a
$C_0(S\unit)$-algebra with the action defined for $\phi\in
C_0(S\unit)$ and $f\in\Gamma_c(S,p^*\mcal{A})$ by 
\[
\Phi(\phi)f(s)= \phi\cdot f(s) := \phi(p(s))f(s).
\]  
Furthermore the restriction map from 
$\Gamma_c(S,p^*\mcal{A})$ to $C_c(S_u,A(u))$ factors to an
isomorphism of $A\rtimes_\alpha S(u)$ onto
$A(u)\rtimes_{\alpha|_{S_u}} S_u$.  
\end{prop}

\begin{remark}
\index{crossed product!bundle}
\index{upper-semicontinuous!cstar-bundle@$C^*$-bundle}
It follows that if $S$ is a group bundle acting on $A$ then there is
an upper-semicontinuous bundle associated to the crossed product
$A\rtimes S$ such that, up to isomorphism, the fibres are the group
crossed products $A(u)\rtimes S_u$. For this reason we will sometimes
refer to $A\rtimes S$ as a crossed product bundle, a bundle crossed
product, or occasionally, a crossed bundle product.   
\end{remark}

Before we prove Proposition \ref{prop:65} we need a series of lemmas.  

\begin{lemma}[{\cite[Lemma 8.3]{tfb2}}]
\label{lem:18}
Suppose $A$ is a $C^*$-algebra and $T\in M(A)$ is such that $T(ab) =
aT(b)$ for all $a,b\in A$.  Then $T\in ZM(A)$.  
\end{lemma}

\begin{proof}
Let  $S\in M(A)$ and $a,b\in A$.  Then 
\[
ST(ab) = S(aT(b)) = S(a)T(b) = T(S(a)b) = TS(ab).
\]
This suffices as $A^2$ is dense in $A$.
\end{proof}

\begin{lemma}
\label{lem:20}
\index{upper-semicontinuous}
The uniform limit of a net of upper-semicontinuous functions is
upper-semicontinuous.
\end{lemma}

\begin{proof}
This is straightforward once you pick the ``right'' definition of
upper semicontinuous.  Recall that $f:X\rightarrow \R$ is 
upper-semicontinuous if and only if $\{x:f(x)<a\}$ is open for all
$a\in\R$.  We claim this is equivalent to the condition that given
$x_0\in X$ and $\epsilon > 0$ there exists an open neighborhood $U$ of
$x_0$ such that $x\in U$ implies $f(x)\leq f(x_0)+\epsilon$.  Let's
start by proving the forward direction.  Given $x_0\in X$ and
$\epsilon > 0$ the set $\{x:f(x)<f(x_0)+\epsilon\}$ is open by
assumption and obviously gives us the desired neighborhood of
$x_0$.  Now let's tackle the reverse direction.  Suppose $a\in \R$ and
$f(x_0)< a$.  Then let $0< \epsilon < a-f(x_0)$ and find a
neighborhood $U$ of $x_0$ as above.  Then if $x\in U$ we have 
\[
f(x) \leq f(x_0)+\epsilon < a
\]
This is enough to show that $\{x:f(x)<a\}$ is open. 

Now suppose $f_i\rightarrow f$ uniformly and $f_i$ is
upper-semicontinuous for all $i$.  Then given $\epsilon >0$ and
$x_0\in X$ choose $I$ such that $\|f_I-f\|_\infty < \epsilon/3$ and
$U$ such that $x\in U$ implies that $f_I(x)\leq f_I(x_0)+\epsilon/3$.
We then have 
\[
f(x) = f(x)-f_I(x)+f_I(x) \leq \frac{2\epsilon}{3} + f_I(x_0)
= \frac{2\epsilon}{3} + f(x_0) + f_I(x_0)-f(x_0) \leq f(x_0) + \epsilon. \qedhere
\]
\end{proof}

The following lemma is interesting in its own right.  Notice that it
holds for any groupoid, not just group bundles.   

\begin{lemma}
\label{lem:21}
Suppose $(A,G,\alpha)$ is a separable dynamical system and $G$ has
Haar system $\lambda$.  Then the function 
\[
u\mapsto \int_G \|f(\gamma)\| d\lambda^u(\gamma)
\]
is upper-semicontinuous for all $f\in\Gamma_c(G,r^*\mcal{A})$.  
\end{lemma}

\begin{proof}
Given $f\in\Gamma_c(G,r^*\mcal{A})$ define $\lambda(f):G\unit\rightarrow \R$ by 
\[
\lambda(f)(u) = \int_G \|f(\gamma)\| d\lambda^u(\gamma).
\]
If $\phi\otimes a$ is an elementary tensor in $C_c(G)\odot A$ then 
\[
\lambda(\phi\otimes a)(u) = \int_G |\phi(\gamma)|d\lambda^u(\gamma) \|a(u)\|.
\]
It follows from the continuity of the Haar system and the fact
that $u\mapsto \|a(u)\|$ is upper-semicontinuous that
$\lambda(\phi\otimes a)$ is upper semicontinuous.  Now suppose
$f\in\Gamma_c(G,r^*\mcal{A})$.  Then there exists a set of
elementary tensors $\{\phi_i^j\otimes a_i^j\}$ such that $k_i = \sum_j
\phi_i^j\otimes a_i^j$ converges to $f$ with 
respect to the inductive limit topology and
therefore with respect to the $I$-norm.  Now, $\lambda(k_i) =
\sum_j\lambda(\phi_i^j\otimes a_i^j)$, and it is straightforward to show
that sums of upper-semicontinuous functions are upper-semicontinuous.
Hence, $\lambda(k_i)$ is upper-semicontinuous.  It now follows quickly
from the computation 
\begin{align*}
\left|\int_G \|f(\gamma)\|d\lambda^u(\gamma) - \int_G
  \|g(\gamma)\|d\lambda^u(\gamma)\right| & \leq
\int_G |\,\|f(\gamma)\|-\|g(\gamma)\|\,| d\lambda^u(\gamma) \\
&\leq \int_G \|f(\gamma)-g(\gamma)\| d\lambda^u(\gamma)\\
&\leq \|f-g\|_I
\end{align*}
that $\lambda(k_i)\rightarrow \lambda(f)$ uniformly and the
result follows from Lemma \ref{lem:20}.
\end{proof}

\begin{proof}[Proof of Proposition \ref{prop:65}]
Given $\phi\in C_0(S\unit)$ and $f\in \Gamma_c(S,p^*\mcal{A})$ define
$\Phi(\phi)f$ as in the statement of the proposition.  It is easy to see
that $\Phi(\phi)f\in \Gamma_c(S,p^*\mcal{A})$ and that $\Phi(\phi)$ is
linear as a function on $\Gamma_c(S,p^*\mcal{A})$.  Next observe that
given $\phi,\psi\in C_0(S\unit)$ and $f\in\Gamma_c(S,p^*\mcal{A})$ we
have 
\begin{equation}
\label{eq:61}
\Phi(\phi)\Phi(\psi)f(s) = \phi(p(s))\Phi(\psi)f(s) =
\phi(p(s))\psi(p(s)) f(s) = \Phi(\phi\psi)f(s).
\end{equation}
Thus $\Phi$ preserves the multiplication on $C_0(S\unit)$.  We need to extend
$\Phi(\phi)$ to an element of the multiplier algebra. The following
slick proof is thanks to Dana Williams.\index{Dana Williams}  
Recall that multipliers on
$A\rtimes S$ are adjointable operators on $A\rtimes S$ where we view
$A\rtimes S$ as a right $A\rtimes S$-module with the operations
\begin{align*}
f\cdot g &= f*g, & \langle f,g \rangle &= f^**g.
\end{align*}
We start by showing that, for all $f,g\in\Gamma_c(S,p^*\mcal{A})$,
\[
\Phi(\phi)(f*g)(s) = \phi(p(s))\int_G f(t)\alpha_t(g(t\inv s))
d\beta^{p(s)}(t) = (\Phi(\phi)f)*g(s).
\]
Thus $\Phi(\phi)$ is $A\rtimes S$-linear on
$\Gamma_c(S,p^*\mcal{A})$.  Next we show that
\begin{equation}
\label{eq:35}
\langle \Phi(\phi)f, g \rangle = 
\langle f, \Phi(\overline{\phi})g \rangle. 
\end{equation}
We compute for $s\in S$ and $u= p(s)$
\begin{align*}
\langle \Phi(\phi)f,g\rangle(s) &= (\phi\cdot f)^**g(s) = 
\int_S \alpha_t(\phi\cdot f(t\inv)^*)\alpha_t(g(t\inv s)) d\beta^u(t)
\\
&= \int_S \alpha_t(\overline{\phi(u)} f(t\inv)^*g(t\inv s))d\beta^u(t)
= \int_S \alpha_t(f(t\inv)^*\overline{\phi}\cdot g(t\inv
s))d\beta^u(t) \\
&= f^* * \overline{\phi}\cdot g(s) = \langle f,
\Phi(\overline{\phi})g\rangle.
\end{align*}
Now extend $\Phi$ to the unitization $C_0(S\unit)^1$ by setting
$\Phi(\phi+\lambda 1)f = \Phi(\phi)f + \lambda f$.  It is an easy
computation to show that \eqref{eq:61} and \eqref{eq:35} extend to all of 
$C_0(S\unit)^1$. Now
suppose $\phi\in C_0(G\unit)$ and $f\in\Gamma_c(S,p^*\mcal{A})$.  It
follows from the $C^*$-identity that $\|\langle f,f\rangle\| =
\|f\|^2$.  If we want to show that $\|\Phi(\phi) f\|\leq
\|\phi\|_\infty\|f\|$ then, because the norm respects the ordering in a
$C^*$-algebra, it suffices to show that as $C^*$-algebra elements
$\langle \phi\cdot f, \phi\cdot f \rangle \leq \|\phi\|_\infty^2 \langle f,
f\rangle$.  However, using \eqref{eq:35} we find that this is
equivalent to showing 
\begin{align*}
0 &\leq \|\phi\|_\infty^2 \langle f, f \rangle - \langle
\Phi(\phi)f,\Phi(\phi)f\rangle \\
&= \langle \|\phi\|_\infty^2 f, f\rangle + \langle
\Phi(\overline{\phi})\Phi(\phi)f,f\rangle \\
&= \langle \Phi(\|\phi\|_\infty^2 1 - \overline{\phi}\phi)f,f\rangle.
\end{align*}
Since general $C^*$-algebra theory assures us that $\|\phi\|_\infty^2 1 -
\overline{\phi}\phi$ is positive in $C_0(S\unit)^1$ it follows that
there is some  $\xi \in C_0(S\unit)^1$ such that $\xi^*\xi =
\|\phi\|_\infty^2 1 - \overline{\phi}\phi$.  We now compute
\begin{align*}
\langle \Phi(\|\phi\|_\infty^2 1 - \overline{\phi}\phi)f,f\rangle  &= 
\langle \Phi(\xi^*)\Phi(\xi) f, f\rangle  \\
&= \langle \Phi(\xi)f,\Phi(\xi)f\rangle  \\
&= (\Phi(\xi)f)^**(\Phi(\xi)f) \geq 0
\end{align*}
It follows that $\Phi(\phi)$ is a bounded linear operator on
$\Gamma_c(S,p^*\mcal{A})$ with norm less than $\|\phi\|_\infty$ and as
such extends to an operator on $A\rtimes_\alpha S$.  Since
$\Phi(\phi)$ is $A\rtimes S$-linear on a dense subset, it is $A\rtimes
S$-linear everywhere.  Furthermore
\eqref{eq:35} implies that $\Phi(\phi)$ is adjointable with
$\Phi(\phi)^* = \Phi(\overline{\phi})$.  Hence $\Phi(\phi)\in
M(A\rtimes_\alpha S)$.  It also follows from \eqref{eq:61} and
\eqref{eq:35} that $\Phi$ is a $*$-homomorphism. We would like to show
that $\Phi$ maps into the center.  Using Lemma \ref{lem:18},
it suffices to show that $\Phi(\phi)(f*g) = f*\Phi(\phi)g$ for
all $f,g\in\Gamma_c(S,p^*\mcal{A})$. We compute 
for $s\in S$ and $u = p(s)$
\begin{align*}
\Phi(\phi)(f*g) &= \phi\cdot (f*g)(s) = \int_S \phi(u) f(t)
\alpha_t(g(t\inv s)) d\beta^u(t) \\
&= \int_S f(t)\alpha_t(\phi\cdot g(t\inv s))d\beta^u(t) = f*\Phi(\phi)g(s).
\end{align*}
The last thing we need to do to show that $A\rtimes S$ is a
$C_0(S\unit)$-algebra is show that the set $\Phi(C_0(S\unit))\cdot
A\rtimes S$ is dense in $A\rtimes S$.  However, given
$f\in\Gamma_c(S,p^*\mcal{A})$ we can find $\phi\in C_c(S\unit)$ such
that $\phi$ is one on $p(\supp f)$ and in this case $\Phi(\phi)f =
f$.  It follows immediately that $\Phi$ is nondegenerate.  

Now we tackle the second part of the assertion.  Fix $u\in S\unit$
and define $R:\Gamma_c(S,p^*\mcal{A})\rightarrow C_c(S_u,A(u))$ by
restriction.  It is clear that $R$ is a well-defined linear map and it
is trivial to show that $R$ is a $*$-homomorphism.  Furthermore $R$ is
uniform norm decreasing and it is straightforward to show that it is
continuous with respect to the inductive limit topology.  Thus
Corollary \ref{cor:24} implies that 
$R$ extends to a map $R:A\rtimes_\alpha S\mapsto A(u)\rtimes_{\alpha|_{S_u}} S_u$.
Now, given
$\phi\in C_c(S_u)$ and $a\in A(u)$ find $b\in A$ such that $b(u) =
a$.  Next, since $S$ is second countable, we can use the Tietze
Extension Theorem to extend $\phi$ so that $\phi\in C_c(S)$.  However
it is clear that $R(\phi\otimes b) = \phi\otimes a$ and therefore $\ran R$
contains the elementary tensors.  Treating $C_c(S_u,A(u))$ as sections of
the trivial bundle it follows from Corollary \ref{cor:3} that $\ran R$
is dense in $C_c(S_u,A(u))$.  Since the range is dense, $R$ must be onto.  

Let 
\[
I_u =\cspn \{ \phi\cdot f : \phi\in C_0(S\unit), \phi(x)=0,
f\in\Gamma_c(S,p^*\mcal{A})\}
\]
and recall that, by definition, $A\rtimes S(u) = A\rtimes S/I_u$.  We
would like to show that $I_u = \ker R$. This next part of the proof is
inspired by \cite[Lemma 2.3]{locunitarystab}.  If $\phi(u) = 0$ then for all
$s\in S_u$ we have $\phi\cdot f(s) = \phi(u)f(s) = 0$.  It follows
that $R(\phi\cdot f) =0$ and that $I_u \subset \ker R$.  Next, suppose
that $f\in \Gamma_c(S,p^*\mcal{A})$ and $R(f) = 0$.  Since
$s\mapsto \|f(s)\|$ is upper-semicontinuous, given $\epsilon >
0$ we can find an open neighborhood $U$ of $S_u$ such that
$\|f(s)\|<\epsilon$ for all $s\in U$.  Choose $\phi\in C_c(S\unit)$
such that $0\leq \phi\leq 1$, 
$\phi(u) = 0$, and $\phi$ is one on $p(\supp f\setminus U)$.
Then $\|f-\phi\cdot f\|_\infty <\epsilon$ and $\supp(\phi\cdot
f)\subset \supp f$.  This is enough to show, after a straightforward
argument, that $f\in I_u$.  Now, let
$\pi$ be a representation of $A\rtimes S$ such that $\ker\pi = I_u$.
It follows from the above that if $f,g\in\Gamma_c(S,p^*\mcal{A})$
such that $R(f)-R(g) = 0$ then we have $f-g\in\ker\pi$.  Hence, we
can define a representation $\rho$ of $C_c(S_u,A(u))$ by $\rho(R(f)) =
\pi(f)$.  It is easy to show, using the fact that $\pi$ and $R$
preserve the operations, that $\rho$ does as well.  Furthermore, given
$f\in\Gamma_c(S,p^*\mcal{A})$ for any $\phi\in C_0(S\unit)^+$ such that
$\phi(u)=1$ we have 
\begin{equation}
\|\rho(R(f))\| = \|\pi(f)\|=\|\pi(\phi\cdot f)\| \leq \|\phi\cdot f\|
\leq 
\|\phi \cdot f\|_I \label{eq:64}
\end{equation}
Let $B= \{\phi\in C_0(S\unit)^+: \phi(u) = 1\}$, $M = \inf_{\phi\in
  B} \|\phi \cdot f\|_I$ and observe that \eqref{eq:64} implies
$\|\rho(R(f))\| \leq M$.  We make the following claim. 

\begin{claim}
Given $f$ and $M$ as above we have
\begin{equation}
M = \|R(f)\|_I = \max\left\{ \int_S \|f(s)\| d\beta^u(s), \int_S
  \|f(s\inv)\| d\beta^u(s)\right\}.
\end{equation}
\end{claim}

\begin{proof}[Proof of Claim.]
First observe that, almost by definition, $\|R(f)\|_I \leq M$.  
Suppose $\epsilon > 0$.  It follows from Lemma \ref{lem:21} that the function
\[
v\mapsto \int_S \|f(s)\|d\beta^v(s)
\]
is upper-semicontinuous.  As such there exists a relatively compact 
neighborhood $U$ of $u$ such that
$v\in U$ implies 
\[
\int_S \|f(s)\|d\beta^v(s) \leq \int_S \|f(s)\|d\beta^u(s) + \epsilon
\leq \|R(f)\|_I+\epsilon.
\]
Since $s\mapsto f(s\inv)$ defines an element of
$\Gamma_c(S,p^*\mcal{A})$ we can use Lemma \ref{lem:21} again to see
that the function 
\[
v\mapsto \int_S \|f(s\inv)\|d\beta^v(s)
\]
is upper-semicontinuous.  As such there exists a relatively compact 
neighborhood $V$ of $u$
such that $v\in V$ implies 
\[
\int_S \|f(s\inv)\|d\beta^v(s) \leq \int_S \|f(s\inv)\|d\beta^u(s) +
\epsilon \leq \|R(f)\|_I+\epsilon.
\]
Choose $\phi\in C_c(S\unit)$ such that
$0\leq\phi\leq 1$, $\phi(u)=1$ and $\phi(v) = 0$ for all $v\not\in
U\cap V$.  Now if $v\in U\cap V$
we have 
\begin{align*}
\phi(v) \int_S \|f(s)\|d\beta^v(s) &\leq \phi(v)(\|R(f)\|_I+\epsilon)
\leq  \|R(f)\|_I +\epsilon,\quad\text{and} \\
\phi(v) \int_S \|f(s\inv)\|d\beta^v(s) &\leq
\phi(v)(\|R(f)\|_I+\epsilon) 
\leq \|R(f)\|_I+\epsilon.
\end{align*}
If $v\not\in U\cap V$ we have 
\begin{align*}
\phi(v) \int_S \|f(s)\|d\beta^v(s) &= 0 \leq \|R(f)\|_I+\epsilon,\quad\text{and} \\
\phi(v) \int_S \|f(s\inv)\|d\beta^v(s) &=0 \leq \|R(f)\|_I+\epsilon.
\end{align*}
In any case, we certainly have 
\begin{align*}
\|\phi\cdot f\|_I &= \max\left\{\sup_{v\in S\unit} \phi(v)\int_S
  \|f(s)\|d\beta^v(s), \sup_{v\in\in S\unit} \phi(v)\int_S
  \|f(s\inv)\|d\beta^v(s)\right\} \\
&\leq \|R(f)\|_I+\epsilon.
\end{align*}
Since we were able to find such a
$\phi$ for any $\epsilon$ it follows that $M \leq \|R(f)\|_I$.  The
claim follows. 
\end{proof}

At this point we have shown that $\|\rho(R(f))\|\leq \|R(f)\|_I$ for any
$f\in \Gamma_c(S,p^*\mcal{A})$ and as such $\rho$ is bounded on
$C_c(S_u,A(u))$ with respect to the $I$-norm.  We have already
asserted that $\rho$ is a $*$-homomorphism.
Thus, we can use Proposition \ref{prop:61} to extend $\rho$ to the
entire group crossed product $A(u)\rtimes S_u$.  Furthermore since
$\rho\circ R= \pi$ on a dense subset, this identity extends to all of
$A\rtimes S$.  It follows that $\ker R \subset \ker\pi = I_u$.  This
shows that $\ker R=I_u$ and hence restriction factors to an
isomorphism of $A\rtimes S(u)$ onto $A(u)\rtimes S_u$.
\end{proof}

\begin{remark}
\index{irreducible representation}
One important consequence of Proposition \ref{prop:65} is that the
irreducible representations of $A\rtimes S$ are well behaved.  To
elaborate, Proposition \ref{prop:36} and the second part of
Proposition \ref{prop:65} says that, as a set,
the spectrum $(A\rtimes S)\sidehat$ can be identified with the
disjoint union $\coprod_{u\in S_u} (A(u)\rtimes S_u)\sidehat$.  In
other words, every irreducible representation of the crossed product
bundle $A\rtimes S$ is lifted from an irreducible covariant
representation of the group crossed product $A(u)\rtimes S_u$ for some
$u\in S\unit$.  This will be an important theme in Section
\ref{sec:crossedstab}. 
\end{remark}

We end this section with a restriction of Proposition \ref{prop:65} to
the scalar case.  This is the proposition used in Section
\ref{sec:duality} to give the total space of the dual bundle a
topology.  

\begin{prop}
\label{prop:groupbundalg}
If $S$ is a locally compact Hausdorff group bundle with a Haar system
then $C^*(S)$ is a $C_0(S\unit)$-algebra.  Furthermore the restriction
map from $C_c(S)$ onto $C_c(S_u)$ factors to an isomorphism from
$C^*(S)(u)$ onto the group $C^*$-algebra $C^*(S_u)$ for all $u\in S\unit$.  
\end{prop}

\begin{proof}
The groupoid algebra $C^*(S)$ is defined to be the crossed
product $C_0(S\unit)\rtimes_{\id} S$.  By Proposition \ref{prop:65} the
crossed product is a $C_0(S\unit)$-algebra and the restriction map
factors to an isomorphism from $C^*(S)(u)$ onto the crossed product
${C_0(S\unit)(u)\rtimes_{\id} S_u}$.  However, 
$C_0(S\unit)(u) = \C$ and $\id|_{S_u}$
is still the identity action so that $\C\rtimes S_u$ is equal to
$C^*(S_u)$ by definition.  The fact that $C^*(S_u)$ is also the group algebra
associated to $S_u$ follows from Remark \ref{rem:17}. 
\end{proof}


\section{Transformation Groupoid Algebras}
\label{sec:transform}
Our last special case of groupoid crossed products comes from the notion
of a groupoid action introduced in Section \ref{sec:actions}.
This section is particularly important because it also introduces a very
natural groupoid action associated to any crossed product.  
We start with the following definition.

\begin{definition}
\index{transformation groupoid!$C^*$-algebra}
\index{G-space@$G$-space}
\index[not]{$C^*(G,X)$}
\label{def:47}
Suppose a second countable locally compact Hausdorff groupoid $G$ with a Haar
system acts continuously on a second countable locally compact
Hausdorff space $X$. 
The {\em transformation groupoid $C^*$-algebra} of $(G,X)$ is defined
to be the groupoid $C^*$-algebra of the transformation groupoid $G\ltimes
X$.  Furthermore, we use the notation $C^*(G,X) := C^*(G\ltimes X)$.
When $G$ acts on the right of $X$ the transformation groupoid
$C^*$-algebra is defined in an analogous way and we use the notation
$C^*(X,G) := C^*(X\rtimes G)$.  
\end{definition}

This definition makes sense because the transformation
groupoid has a Haar system whenever $G$ does, by Proposition \ref{prop:14}.

\begin{remark}
Definition \ref{def:47} is slightly specious in that the notation
$C^*(G,X)$ isn't any simpler than $C^*(G\ltimes X)$.  However, it does
connect back to the notation used for transformation group
$C^*$-algebras.  Furthermore, although we won't address it here,
$C^*(G,X)$ simplifies to the transformation group $C^*$-algebra when
$G$ is a group, and of course is also the same as the $C^*$-algebra of
the transformation group groupoid in this case. 
\end{remark}

The ultimate goal of this section will be to show that if $G$ acts on
$X$ then $G$ acts on $C_0(X)$ and that the groupoid crossed product in
this case is the same as the transformation groupoid $C^*$-algebra.
However, before we can prove this result we need the following lemma.
Basically, we need a tool to deal with the topology on the
upper-semicontinuous bundles we will encounter later on.  

\begin{lemma}
\label{lem:26}
Suppose $X$, $Y$  and $Z$ are locally compact Hausdorff spaces and that
$\sigma:Y\rightarrow X$ and $\tau:Z\rightarrow X$ are continuous
surjections. Let $Z*Y = \{(z,y)\in Z\times Y : \tau(z) =
\sigma(y)\}$.  Then the map $\iota:C_0(Z*Y)\rightarrow \tau^*(C_0(Y))$
such that $\iota(f)(z)(y) = f(z,y)$ is an isometric isomorphism. 
Furthermore $\iota(f)$ is compactly supported if $f$ is and $\iota$
preserves convergence with respect to the inductive limit topology.  
\end{lemma}

\begin{proof}
Let $\mcal{C}$ be the upper-semicontinuous bundle associated to
$C_0(Y)$ as a $C_0(X)$-algebra and recall that $C_0(Y)(x) =
C_0(\sigma\inv(x))$.  
Define $\iota: C_c(Z*Y)\rightarrow \Gamma_c(Z,\tau^*\mcal{C})$
by $\iota(f)(z)(y) = f(z,y)$.  
It is clear that $\iota(f)(z)$ is a continuous compactly supported function on
$\sigma\inv(\tau(z))$.  Therefore 
$\iota(f)(z) \in C_0(Y)(\tau(z))$ and $\iota$ is a well defined
section of $\tau^*\mcal{C}$.  It is also clear
that $\iota(f)$ has compact support.  Notice that this also verifies the claim
that $\iota(f)$ is compactly supported if $f$ is.  We need to see that
$\iota(f)$ is continuous.  We start by demonstrating this in a simpler
case. Suppose $g\in C_c(Z)$, $h\in C_c(Y)$ and define $g\otimes
h(z,y) = g(z)h(y)$ for all $(z,y)\in Z*Y$. Let
$z_i \rightarrow z$.  We would like to show that
$\iota(g\otimes h)(z_i)\rightarrow \iota(g\otimes h)(z)$ in
$\mcal{C}$.  Since $h\in C_c(X)$ we can view $h$ as a continuous
section of $\mcal{C}$ with $h(x) = h|_{\sigma\inv(x)}$ for all $x\in
X$ and it follows that $h(\tau(z_i))\rightarrow h(\tau(z))$
in $\mcal{C}$ since $\tau$ is continuous.  Furthermore,
scalar multiplication is continuous and $g(z_i)\rightarrow
g(z)$ so we also have 
\[
g(z_i)h(\tau(z_i))\rightarrow g(z)h(\tau(z))
\]
in $\mcal{C}$.  However, it is clear with some thought that
$\iota(g\otimes h)(w) = g(w)h(\tau(w)) = g(w)h|_{\sigma\inv(\tau(w))}$ for
any $w\in Z$ and therefore $\iota(g\otimes h)\in
\Gamma_c(Z,\tau^*\mcal{C})$.  

Now suppose we have $f\in C_c(Z*Y)$.  
Since $Z*Y$ is closed in $Z\times Y$
we can use Lemma \ref{lem:8} to extend $f$ to the product space.  Next
find $g_i^j\in C_c(Z)$ and $h_i^j\in C_c(Y)$ such that 
$k_i = \sum_j g_i^j\otimes h_i^j
\rightarrow f$ uniformly.  Now let $z_i\rightarrow z$ and
observe that $\tau(z_i)\rightarrow \tau(z)$.  We will show that
$\iota(f)(z_i)\rightarrow \iota(f)(z)$ using Proposition
\ref{prop:35}.  Let $\epsilon > 0$ and choose $I$ such that
$\|k_I -f \|_\infty < \epsilon$.  Since sums of continuous functions
are continuous,
$\iota(k_I)(z_i)\rightarrow \iota(k_I)(z)$ by the previous paragraph.  
Furthermore given any
$w\in Z$ we have 
\[
\|\iota(k_I)(w) -\iota(f)(w)\|_\infty  = 
\sup_{y\in \sigma\inv(\tau(w))} |k_I(w,y) - f(w,y)| 
\leq \|k_I - f\|_\infty < \epsilon.
\]
Since this is true for all $z_i$ and $z$ as well it follows
from the last part of Proposition \ref{prop:35} that
$\iota(f)(z_i)\rightarrow \iota(f)(z)$.  Thus $\iota(f)\in
\Gamma_c(Z,r^*\mcal{C})$.  

Suppose $f\in C_c(Z*Y)$.  We calculate
\begin{align*}
\|\iota(f)\|_\infty &= \sup_{z\in Z} \|\iota(f)(z)\|_\infty
\\ &= \sup_{z\in Z}\sup_{y\in \sigma\inv(\tau(z))} |f(z,y)| =
\|f\|_\infty.
\end{align*}
Hence $\iota$ is isometric.  Next we show that
$\ran \iota$ is dense in $\tau^*C_0(X) = \Gamma_0(Z,\tau^*\mcal{C})$.  
If $f\in C_c(Z*Y)$ and $g\in C_0(Z)$ then observe
that if we define $h(z,y) = g(z)f(z,y)$ then $h\in
C_c(Z*Y)$ and
\[
(g\cdot\iota(f)(z))(y) = g(z)f(z,y) = \iota(h)(z)(y).
\]
It follows that $\ran \iota$ is closed under the action of $C_0(Z)$.
Next, fix $z\in G$ and $f\in C_c(\sigma\inv(\tau(z)))$.  
Find $g\in C_c(Y)$ such that $g|_{\sigma\inv(\tau(z))} = f$ and 
$h\in C_c(Z)$ such that $h(z) = 1$.  Then clearly
$\iota(h\otimes g)(z) = f$ so that $\ran\iota$ is dense
fibrewise.  It follows from Proposition \ref{prop:42} that $\ran\iota$
is dense in $\Gamma_0(Z,\tau^*\mcal{C})$.  Since
$\iota:C_c(Z*Y)\rightarrow \tau^*C_0(Y)$ is an isometry mapping onto a
dense set we can extend $\iota$ to an isomorphism
$\iota:C_0(Z*Y) \rightarrow \tau^*C_0(Y)$.  Furthermore, if $f\in
C_0(Z*Y)$ then we can find a sequence in $C_c(Z*Y)$ such that
$f_i\rightarrow f$ uniformly.  However, it follows that
$\iota(f_i)\rightarrow \iota(f)$ uniformly.  Thus $\iota(f_i)(z)\rightarrow
\iota(f)(z)$.  This convergence takes place in
$C_0(Y)(\tau(z)) = C_0(\sigma\inv(\tau(z)))$ and is therefore
convergence with respect to the uniform norm.  Again this implies that
for $y\in \sigma\inv(\tau(z))$ we have 
\[
\iota(f_i)(z)(y) = f_i(z,y)\rightarrow \iota(f)(z)(y).
\]
However we also have $f_i(z,y)\rightarrow f(z,y)$.  Thus $\iota$ has
the desired form on all of $C_0(Z*Y)$.  

Finally, observe that if $f_i\rightarrow f$ with respect to the
inductive limit topology then $\iota(f_i)\rightarrow \iota(f)$
uniformly.  Furthermore if the supports of $f_i$ are eventually
contained in $K$ then the supports of $\iota(f_i)$ are eventually
contained in the projection of $K$ onto $Z$.  Thus
$\iota(f_i)\rightarrow \iota(f)$ with respect to the inductive limit
topology. 
\end{proof}

\begin{prop}
\label{prop:68}
\index{crossed product}
Suppose a second countable, locally compact Hausdorff groupoid $G$ acts on a
second countable, locally compact Hausdorff space $X$.  Then $C_0(X)$ is a
$C_0(G\unit)$-algebra and there is an action of $G$ on $C_0(X)$ given
by the maps $\lt_\gamma:C_0(r_X\inv(s(\gamma))) \rightarrow
C_0(r_X\inv(r(\gamma)))$ where 
\begin{equation}
\lt_\gamma(f)(y) = f(\gamma\inv\cdot y)
\end{equation}
for all $f\in C_0(r_X\inv(s(\gamma)))$ and $y\in r_X\inv(r(\gamma))$.
Furthermore, the groupoid crossed product $C_0(X)\rtimes_{\lt} G$ is
isomorphic to $C^*(G\ltimes X)$ and the isomorphism is given on
$C_c(G\ltimes X)$ by $\iota(f)(\gamma)(x) = f(\gamma,x)$.  
\end{prop}

\begin{proof}
The fact that $C_0(X)$ is a $C_0(G\unit)$-algebra is really just 
Example \ref{ex:15} using the map $r_X:X\rightarrow
G\unit$.  The action is given by $\Phi(f)(g)(x) = f(r(x))g(x)$ and it
is straightforward to use the Stone-Weierstrass theorem to show that
this makes $C_0(X)$ into a $C_0(G\unit)$-algebra.  Let $I_u$ be the
ideal such that $C_0(X)(u)$ is defined to be $C_0(X)/I_u$ and 
$I_{r\inv(u)}$ the ideal
of functions which are zero on $r_X\inv(u)$.  It is clear from the
definition of $I_u$ that $I_u\subset I_{r\inv(u)}$.  However, given
$f\in I_{r\inv(u)}$ and $\epsilon > 0$ since $|f|$ is continuous
the set $O=\{x\in X: |f|<\epsilon\}$ is open.  Choose $g\in
C_c(G\unit)$ so that $g(u) = 0$ and $g$ is one on $\supp f \setminus
U$.  It's easy to check that $\|\Phi(g)f - f\|_\infty < \epsilon$. It 
follows quickly that $f\in I_u$ and $I_{r\inv(u)} = I_u$.  However,
$I_{r\inv(u)}$ is the ideal of functions which are zero on some
closed set and it follows from classical theory that
$C_0(X)/I_{r\inv(u)} \cong C_0(r\inv(u))$.  Thus, we can view the
upper-semicontinuous bundle $\mcal{C}$ associated to $C_0(X)$ as
having fibres $C_0(r_X\inv(u))$.  
We deal with the fact that the topology on $\mcal{C}$
is mysterious by using Proposition \ref{prop:44}.  

Use Lemma \ref{lem:26} to observe that there are two isomorphisms 
\begin{align*}
\iota_s&: C_0(G*X)\rightarrow s^*C_0(X), & \iota_r&:C_0(G\ltimes
X)\rightarrow r^*C_0(X)
\end{align*}
where
\begin{align*}
G*X &= \{(\gamma,x)\in G\times X : s(\gamma) = r(x)\},\quad\text{and}
\\
G\rtimes X &= \{(\gamma,x)\in G\times X : r(\gamma) = r(x)\}.
\end{align*}
We use these isomorphisms to, perhaps
foolishly, identify
$C_0(G*X)$ with $s^*C_0(X)$ and $C_0(G\ltimes X)$ with $r^*C_0(X)$.
We define a map $\lt: s^*C_0(X)\rightarrow r^*C_0(X)$ by 
\[
\lt(f)(\gamma,x) = f(\gamma,\gamma\inv\cdot x).
\]
It is clear that $\lt$ is a $*$-homomorphism and it is easy to
construct an inverse for $\lt$ using ``right translation.''
Furthermore, if $\phi\in C_0(G)$ we have
\[
\phi\cdot \lt(f)(\gamma,x) = \phi(\gamma) f(\gamma,\gamma\inv\cdot x)
= \lt(\phi\cdot f)(\gamma,x).
\]
Thus $\lt$ is a $C_0(G)$-linear $*$-isomorphism.  As such it defines a
family of isomorphisms $\lt_\gamma:C_0(X)(s(\gamma))\rightarrow
C_0(X)(r(\gamma))$.  Given $f\in C_0(X)(s(\gamma)) =
C_0(r_X\inv(s(\gamma)))$ extend $f$ to a function $g\in C_0(X)$ and
pick $h\in C_c(G)$ such that $h(\gamma) = 1$.  By construction
$h\otimes g(\gamma) = f$ so that we have, by the definition of $\lt_\gamma$,
\[
\lt_\gamma(f) = \lt_\gamma(h\otimes g(\gamma)) = 
\lt(h\otimes g)(\gamma).
\]
But now we can compute
\[
\lt_\gamma(f)(x) = \lt(h\otimes g)(\gamma,x) =
h(\gamma)g(\gamma\inv\cdot x) = f(\gamma\inv\cdot x).
\]
Thus the action $\lt$ has the desired form.  Furthermore, it is now
easy to verify that $\lt_\gamma\circ \lt_\eta = \lt_{\gamma\eta}$.  It
follows from Proposition \ref{prop:44} that $(C_0(X),G,\lt)$ is a
groupoid dynamical system.  

Now consider the restriction of the isomorphism $\iota_r$ to 
$\iota:C_c(G\ltimes X)\rightarrow 
\Gamma_c(G,r^*\mcal{C})$. We already used the fact that $\iota$
preserves the pointwise operations.  Now we show that it preserves the
convolution and involution as well.  We compute
\begin{align*}
\iota(f)*\iota(g)(\gamma)(x) &= 
\int_G
\iota(f)(\eta)(x)\lt_\eta(\iota(g)(\eta\inv\gamma))(x)d\lambda^{r(\gamma)}(\eta)
\\
&= \int_G f(\gamma,x)g(\eta\inv\gamma,\eta\inv\cdot x)
d\lambda^{r(\gamma)}(\eta) \\
&= \int_G f(\gamma,x)g((\eta,x)\inv(\gamma,x))d\lambda^{r(x)}(\eta) \\
&= f*g(\gamma,x) = \iota(f*g)(\gamma)(x)
\end{align*}
and
\begin{align*}
\iota(f)^*(\gamma)(x) &= \lt_\gamma(\iota(f)(\gamma\inv)^*)(x) = 
\overline{f(\gamma\inv,\gamma\inv\cdot x)} \\
&= \overline{f((\gamma,x)\inv)} = f^*(\gamma,x) \\
&= \iota(f^*)(\gamma)(x).
\end{align*}
Thus $\iota$ is a $*$-homomorphism with respect to the crossed product
operations as well.  Furthermore we proved in Lemma \ref{lem:26} that
$\iota$ is continuous with respect to the inductive limit topology.
It then follows from Corollary \ref{cor:24} 
that $\iota$ extends to a homomorphism $I:C^*(G\ltimes
X)\rightarrow C_0(X)\rtimes_{\lt} G$. 
We would like to see that this map is an isomorphism. Let $j:\ran\iota
\rightarrow C_c(G\ltimes X)$ be the inverse of $\iota$ so that
$j(f)(\gamma,x) = f(\gamma)(x)$.  Then $j$ is clearly onto, injective,
and a $*$-homomorphism.  Given $x\in X$ we have 
\begin{align*}
\int_G |j(f)(\gamma,x)| d\lambda^{r(x)}(\gamma) &= 
\int_G |f(\gamma)(x)| d\lambda^{r(x)}(\gamma) \\
&\leq \int_G \|f(\gamma)\|_\infty d\lambda^{r(x)}(\gamma) \leq \|f\|_I.
\end{align*}
Similarly we calculate that 
\[
\int_G |j(f)(\gamma\inv,\gamma\inv\cdot x)| d\lambda^{r(x)}(\gamma)
\leq
\int_G \|f(\gamma\inv)\|_\infty d\lambda^{r(x)}(\gamma) \leq \|f\|_I.
\]
Once we recall the definition of the Haar system on $G\ltimes
X$ it follows immediately that $j$ is $I$-norm decreasing.  
It is straightforward to show, using the fact that $i_r$ is an
isomorphism, that $\ran\iota$ is dense in $\Gamma_c(G,r^*\mcal{C})$
with respect to the inductive limit topology and with respect to the
$I$-norm. Now, extend $j$ to $\Gamma_c(G,r^*\mcal{C})$.  Then 
Proposition \ref{prop:62} implies that $j$ extends to a
$*$-homomorphism $J:C_0(X)\rtimes_{\lt} G\rightarrow C^*(G\ltimes X)$.
However, $J$ and $I$ are inverses on a dense set so they must be
inverses everywhere.  It follows that $I$ is an isomorphism of the
transformation groupoid $C^*$-algebra onto the groupoid crossed
product.  
\end{proof}

There is a ``converse'' to Proposition \ref{prop:68} in that we
can identify those crossed products which arise from transformation
groupoids.  However, the following piece of this construction is
necessary for what we will do in Chapter \ref{cha:fine-structure} 
and is very interesting in its own right. 

\begin{prop}
\label{prop:69}
\index{G-space@$G$-space!induced from a dynamical system}
If $(A,G,\alpha)$ is a groupoid dynamical system then there is a
continuous action of $G$ on $\widehat{A}$ given by $\gamma\cdot \pi=
\pi\circ \alpha_\gamma\inv$.  
\end{prop}

\begin{proof}
Since $A$ is a $C_0(G\unit)$-algebra it follows from Proposition
\ref{prop:36} that there is a continuous map $r:\widehat{A}\rightarrow
G\unit$.  Furthermore, as in Proposition \ref{prop:36}, we view
$\widehat{A}$ as being fibred over $G\unit$ so that if
$\pi\in\widehat{A}$ with $r(\pi) = u$ then $I_u\subset \ker\pi$ and 
we can factor $\pi$ to a
representation $\pi'$ of $A(u)$.  Oftentimes we will not distinguish
between $\pi$ and $\pi'$, but it is important to do so
now.  Given $\gamma\in G$ we know
$\alpha_\gamma:A(s(\gamma))\rightarrow A(r(\gamma))$ so that if
$r(\pi) = s(\gamma)$ we can define $\gamma\cdot \pi\in \widehat{A}$ by
\[
\gamma\cdot \pi(a) = \pi'(\alpha_\gamma\inv(a(r(\gamma)))).
\]
Of course when we factor $\gamma\cdot \pi$ to $A(r(\gamma))$ we get
$(\gamma\cdot \pi)' = \pi'\circ\alpha_\gamma\inv$ as desired.  
Given $u\in G\unit$ we know  $\alpha_u = \id$ so it is clear that
$r(\pi)\cdot \pi = \pi$.  Furthermore if $\gamma$ and $\eta$ are
composable then 
\[
\gamma\cdot(\eta\cdot \pi)(a) =
\pi'(\alpha_\eta\inv(\alpha_\gamma\inv(a(r(\gamma))))) 
= (\gamma\eta)\cdot \pi(a).
\]
All that is left is to check that the action is continuous.  Suppose
$\gamma_i \rightarrow \gamma$ and  $\pi_i\rightarrow \pi$ such that
$s(\gamma_i) = r(\pi_i)$ for all $i$ and $s(\gamma) = r(\pi)$.
Let $O_J = \{\pi\in\widehat{A}: J\not\subset \ker\pi\}$ be an open set in
$\widehat{A}$ containing $\gamma\cdot \pi$.  Recall that every open set is
of this form \cite[Corollary A.28]{tfb}.  Suppose, to the contrary, 
that $\gamma_i\cdot
\pi_i$ is not eventually in $O_J$.  By passing to a subnet and
relabeling  we can assume $\gamma_i\cdot\pi_i\not\in O_J$ for all
$i$.  Let $a\in J$ and choose $b\in A$ such that $b(s(\gamma)) =
\alpha_\gamma\inv(a(r(\gamma)))$.  Since the action is continuous,
$\alpha_{\gamma_i}\inv(a(r(\gamma_i))) \rightarrow b(s(\gamma))$.
Since the norm is upper-semicontinuous, the set
$\{ a\in \mcal{A} : \|a\| < \epsilon\}$ is open for all $\epsilon >
0$.  Because $\alpha_{\gamma_i}\inv(a(r(\gamma_i))) - b(s(\gamma_i))
\rightarrow 0$ we eventually have
$\|\alpha_{\gamma_i}\inv(a(r(\gamma_i))) - b(s(\gamma_i))\| <
\epsilon$ for all $\epsilon >0$.  Hence
\[
\|\alpha_{\gamma_i}\inv(a(r(\gamma_i))) - b(s(\gamma_i))\| \rightarrow
0.
\]
Next, $\gamma_i\cdot \pi_i\not\in O_J$ for all $i$ so we have 
\[
\gamma_i\cdot \pi_i(a) = \pi'(\alpha_{\gamma_i}\inv(a(r(\gamma_i)))) =
0
\]
for all $i$.  Hence
\begin{align}
\nonumber
\|\pi_i(b)\| &= \|\pi'_i(b(s(\gamma_i)))\| = \|\pi'_i(b(s(\gamma_i))
- \alpha_{\gamma_i}\inv(a(r(\gamma_i))))\| \\
&\leq \|b(s(\gamma_i))-\alpha_{\gamma_i}\inv(a(r(\gamma_i)))\|
\rightarrow 0. \label{eq:65}
\end{align}
It is shown in \cite[Lemma A.30]{tfb} that the map
$\pi\mapsto\|\pi(b)\|$ is lower-semicontinuous on $\widehat{A}$.  
In other words given $\epsilon \geq 0$ the set 
\[
\{ \rho\in \widehat{A} : \|\rho(b)\|\leq \epsilon \}
\]
is closed.  Now, given $\epsilon > 0$ equation \eqref{eq:65} implies
that eventually
$\pi_i\in\{\rho\in\widehat{A}:\|\rho(b)\|\leq\epsilon\}$.  Since this
set is closed and $\pi_i\rightarrow \pi$ we know that $\|\pi(b)\|\leq
\epsilon$.  This is true for all $\epsilon > 0$ so that 
\[
0 = \pi(b) = \pi'(b(s(\gamma))) =
\pi'(\alpha_\gamma\inv(a(r(\gamma)))) = \gamma\cdot \pi(a).
\]
This is a contradiction since $a\in J$ was arbitrary and we assumed
that $\gamma\cdot \pi\in O_J$.  It follows that the action of $G$ on
$\widehat{A}$ is continuous.\footnote{The author finds it amusing that
 this proof uses both upper and lower semicontinuity.\index{upper-semicontinuous}}
\end{proof}

This allows us to prove that every groupoid action on an {\em abelian}
$C^*$-algebra comes from a groupoid action on a space, which provides
the promised ``converse'' to Proposition \ref{prop:68}. 

\begin{prop}
\index{dynamical system}
Suppose $X$ is a locally compact Hausdorff space and
$(C_0(X),G,\alpha)$ is a separable dynamical system.  Then there is
an action of $G$ on $X$ such that $\alpha $ is given by left
translation.  Consequently $C_0(X)\rtimes_\alpha G \cong C^*(G,X)$.  
\end{prop}

\begin{proof}
We know from Proposition \ref{prop:69} that there is a continuous
action of $G$ on $C_0(X)\sidehat$.  We can identify $C_0(X)\sidehat$
and $X$ via the Gelfand transform which takes $x\in X$ to ``evaluation
at x'', denoted $\ev_x$, and $\pi\in C_0(X)\sidehat$ to the element
$\hat{\pi}\in X$ determined by $f(\hat{\pi}) = \pi(f)$ for $f\in
C_0(X)$.  Thus, using Proposition \ref{prop:69}, 
we can view $G$ as acting on $X$ via the formula
\begin{equation}
\gamma\cdot x = (\gamma\cdot \ev_x)\sidehat = (\ev_x\circ
\alpha_\gamma\inv)\sidehat.
\end{equation}
This action of $G$ on $X$ induces an action of $G$ on $C_0(X)$ via left
translation as in Proposition \ref{prop:68}.  That is, given $f\in
C_0(X)(s(\gamma)) = C_0(r_X\inv(s(\gamma)))$ we define for $x\in
r_X\inv(r(\gamma))$, 
\[
\lt_\gamma (f)(x) = f(\gamma\inv\cdot x).
\]
However, we compute that 
\begin{align*}
\lt_\gamma (f)(x) &= f(\gamma\inv\cdot x) =
f((\ev_x\circ\alpha_\gamma)\sidehat) \\
&= \ev_x(\alpha_\gamma(f)) = \alpha_\gamma(f)(x)
\end{align*}
Thus $\lt_\gamma =\alpha_\gamma$ and $\alpha$ is given by left
translation.  It follows from Proposition \ref{prop:68} that $C^*(G,X)
\cong C_0(X)\rtimes_\alpha G$.  
\end{proof}

\subsection{Stone Von-Neumann Theorem}
\label{sec:neumann}
While we are on the subject of groupoid actions and groupoid spaces we
may as well discuss one of the most basic groupoid actions, the action
of $G$ on itself by translation.  Recall from Example
\ref{ex:11} that the action is defined by $\gamma\cdot
\eta = \gamma\eta$.  This is a very well behaved example and we can
give a nice description of its associated transformation groupoid
$C^*$-algebra.  This theorem is a generalization of the Stone
Von-Neumann theorem for group crossed products \cite[Theorem
4.24]{tfb2}.  It is slightly surprising that it can be extended to
groupoids with such generality. 

\begin{theorem}
\label{thm:stonevn}
\index{New Result}
Suppose $G$ is a second countable locally compact Hausdorff groupoid
with a Haar system and let $G$ act on the left of itself by
translation.  Then there is an isomorphism $\Phi$ from $C^*(G,G)$ onto
$\mcal{K}(\mcal{Z})$ where $\mcal{Z}$ is the completion of the pre-Hilbert
$C_0(G\unit)$-module $\mcal{Z}_0=C_c(G)$ given the operations 
\begin{align*}
z\cdot \phi(\gamma) &= z(\gamma)\phi(s(\gamma)) &
\langle w,z\rangle(u) &= \int_G
\overline{w(\gamma)}z(\gamma)d\lambda_u(\gamma).
\end{align*}
Furthermore, given $f\in C_c(G\ltimes G)$ and $z\in Z_0$ we have 
\begin{equation}
\Phi(f)z(\gamma) = \int_G
f(\eta,\gamma)z(\eta\inv\gamma)d\lambda^{r(\gamma)}(\eta).
\end{equation}
\end{theorem}

\begin{proof}
First, consider the fact that if we let $G\unit$ act on the right of
$G$ by the trivial action then $G$ is a strong, principal, right
$G\unit$-space.  As such we can form the imprimitivity groupoid $H =
G^{G\unit}$.  Since the action is trivial, we have 
\[
H = G*G = \{(\gamma,\eta)\in G\times G: s(\gamma) = s(\eta)\}
\]
and the operations are just $(\gamma,\eta)(\eta,\zeta) = (\gamma,\zeta)$
and $(\gamma,\eta)\inv = (\eta,\gamma)$.  It follows from Proposition
\ref{prop:19} that $H$ acts on the left of $G$ and that with this
action $G$ is a $(H,G\unit)$-equivalence.  Furthermore, once we
untangle the definitions, it follows that the action of $H$ on $G$ is defined by
$(\gamma,\eta)\cdot \eta = \gamma$.  Now $G$ has a Haar system by
assumption, and since the restriction of the range and source maps to
$G\unit$ are obviously open it follows from Proposition \ref{prop:21}
that $H$ has a Haar system $\mu$ defined by 
\[
\int_H f(\eta,\gamma)d\mu^{\zeta}(\eta,\gamma) = \int_G f(\zeta,\gamma)
d\lambda_{s(\zeta)}(\gamma).
\]
Furthermore, we saw in Example \ref{ex:23} that $C^*(G\unit)$ is just
$C_0(G\unit)$ and the Haar system is given by the Dirac delta measures
$\{\delta_u\}$.  

Since both $G\unit$ and $H$ have Haar systems we can use Theorem
\ref{thm:groupoidequiv} to view $\mcal{Z}_0 = C_c(G)$ as a
pre-$C^*(H)-C_0(G\unit)$-imprimitivity bimodule.  In particular 
$\mcal{Z}_0$ is
a pre-Hilbert $C_0(G\unit)$-module and we can compute that the
operations are given for $z\in C_c(G)$ and $\phi\in C_0(G\unit)$ by
\begin{align*}
z\cdot \phi(\gamma) &= \int_{G\unit} z(\gamma\cdot u)\phi(u\inv)
d\delta_{s(\gamma)}(u) \\
&= z(\gamma)\phi(s(\gamma)),
\end{align*}
and by picking $\gamma\in G$ such that $s(\gamma) = u$,
\begin{align*}
\langle z,w \rangle_{C_0(G\unit)}(u) &= \int_H
\overline{z((\eta,\zeta)\inv \cdot \gamma)} w((\eta,\zeta)\inv \cdot
\gamma \cdot u) d\mu^{\gamma}(\eta,\zeta) \\
&= \int_G
\overline{z((\zeta,\gamma)\cdot\gamma)}w((\zeta,\gamma)\cdot\gamma)
d\lambda_{s(\gamma)}(\zeta) \\
&= \int_G \overline{z(\zeta)}w(\zeta) d\lambda_u(\zeta).
\end{align*}
Thus $\mcal{Z}_0$ has the appropriate operations and defines the desired
Hilbert $C_0(G\unit)$-module $\mcal{Z}$.  Furthermore, since $\mcal{Z}$ is an
$C^*(H)-C_0(G\unit)$-imprimitivity bimodule it follows from
\cite[Proposition 3.8]{tfb} that $C^*(H)$ is isomorphic to the compact
operators $\mcal{K}(\mcal{Z})$ and that the isomorphism is given by
$\Psi(f)(z) = f\cdot z$.  Thus we can compute for $f\in C_c(H)$ and
$z\in C_c(G)$ that 
\begin{align*}
\Psi(f)(z)(\gamma) &= f\cdot z(\gamma) = \int_H
f(\eta,\zeta)z((\eta,\zeta)\inv\cdot \gamma) d\mu^\gamma(\eta,\zeta)
\\
&= \int_G f(\gamma,\zeta) z(\zeta) d\lambda_{s(\gamma)}(\zeta).
\end{align*}

Now, consider the transformation groupoid $G\ltimes G$.  Define
$\phi:H\rightarrow G\ltimes G$ by $\phi(\gamma,\eta) =
(\gamma\eta\inv,\gamma)$.  This map is clearly continuous and has a
continuous inverse given by $(\gamma,\eta)\mapsto
(\eta,\gamma\inv\eta)$.  Furthermore, we can check that 
\begin{align*}
\phi((\gamma,\eta)(\eta,\xi))&= \phi(\gamma,\xi)   =
(\gamma\xi\inv,\gamma) \\
&= (\gamma\eta\inv,\gamma)(\eta\xi\inv,\eta) =
\phi(\gamma,\eta)\phi(\eta,\xi),\quad\text{and} \\
\phi((\gamma,\eta)\inv) &= \phi(\eta,\gamma) = (\eta\gamma\inv,\eta) \\
&= (\gamma\eta\inv,\gamma)\inv = \phi(\gamma,\eta)\inv.
\end{align*}
Thus $\phi$ is a groupoid isomorphism.  Now recall that the Haar
system on $G\ltimes G$ is defined by $\beta^\gamma =
\lambda^{r(\gamma)}\times \delta_\gamma$.  We can then check that, for
$f\in C_c(G\ltimes G)$ 
\begin{align*}
\int_{H} f(\phi(\zeta,\eta)) d\mu^{\gamma}(\zeta,\eta) &= 
\int_G f(\phi(\gamma,\eta))d\lambda_{s(\gamma)}(\eta)
= \int_G f(\gamma\eta\inv,\gamma)d\lambda_{s(\gamma)}(\eta) \\
&= \int_G f(\eta,\gamma)d\lambda^{r(\gamma)}(\eta) 
= \int_{G\ltimes G} f(\eta,\zeta)d\beta^{\gamma}(\eta,\zeta).
\end{align*}
Thus $\phi$ intertwines the Haar systems on $G\ltimes G$ and $H$.  It
follows from Proposition \ref{prop:70} that the map
$\Upsilon:C^*(G)\rightarrow C^*(H)$ defined for $f\in C_c(G)$ by
$\Upsilon(f)= f\circ \phi$ is an isomorphism.  But now we can define
an isomorphism $\Phi: C^*(G)\rightarrow \mcal{K}(\mcal{Z})$ by $\Phi=\Psi
\circ \Upsilon$.  Furthermore, we can compute for $f\in C_c(G\ltimes G)$ and
$z\in C_c(G)$ that 
\begin{align*}
\Phi(f)z(\gamma) &= \Psi(f\circ\phi)z(\gamma) =
\int_G f\circ\phi(\gamma,\zeta)z(\zeta)d\lambda_{s(\gamma)}(\zeta) \\
&= \int_G f(\gamma\zeta\inv,\gamma)z(\zeta)d\lambda_{s(\gamma)}(\zeta)
= \int_G f(\gamma\zeta,\gamma)z(\zeta\inv)d\lambda^{s(\gamma)}(\zeta)
\\
&= \int_G
f(\zeta,\gamma)z(\zeta\inv\gamma)d\lambda^{r(\gamma)}(\zeta).
\qedhere
\end{align*}
\end{proof}

Theorem \ref{thm:stonevn} is not just a generalization of the
Stone-von Neumann theorem in the thematic sense.  The classical
Stone-von Neumann theorem can be obtained as a corollary.  

\begin{corr}
\label{cor:8}
Suppose $G$ is a second countable locally compact group.  Let $M$ be
the representation of $C_0(G)$ on $L^2(G)$ given by multiplication and
let $L$ be the left regular representation.  Then $M\rtimes
L$ is a faithful representation of $C_0(G)\rtimes_{\lt} G$ onto
$\mcal{K}(L^2(G))$.  
\end{corr}

\begin{proof}
Let us first use Proposition \ref{prop:68} to identify
$C_0(G)\rtimes_{\lt} G$ with $C^*(G,G) = C^*(G\ltimes G)$.  Once we
perform this identification the representation $M\rtimes L$
becomes, for $f\in C_c(G\ltimes G)$ and $\phi\in \mcal{L}^2(G)$, 
\[
M\rtimes L(f)\phi(s) = 
\int_G M(f(t))L_t\phi(s)\Delta(t)\neghalf dt = \int_G f(t,s)
\phi(t\inv s)\Delta(t)\neghalf dt.
\]
Now, since $G$ is a group, $G\unit$ consists of a single point and
$G\ltimes G = G\times G$.  We can use Theorem \ref{thm:stonevn} to
construct an isomorphism $\Phi$ of $C^*(G\ltimes G)$ with 
$\mcal{K}(\mcal{Z})$ where
$\mcal{Z}$ is the completion of the pre-Hilbert $C_0(G\unit)$-module
$\mcal{Z}_0=C_c(G)$.  However $C_0(G\unit) = \C$ so that $\mcal{Z}$ is a Hilbert
space.  Furthermore, given $\phi,\psi\in C_c(G)$ we have 
\[
\langle \phi,\psi \rangle = \int_G \overline{\phi(s)}\psi(s)
d\lambda_e(s) = \int_G \overline{\phi(s)}\psi(s) \Delta(s)ds.
\]
Now define $U:\mcal{Z}_0\rightarrow L^2(G)$ by $U\phi(s) =
\phi(s)\Delta(s)\neghalf$.  It is straightforward to show that $U$
extends to a unitary from $\mcal{Z}$ onto $L^2(G)$ and thus $U\Phi U^*$ is an
isomorphism from $C^*(G,G)$ onto $\mcal{K}(L^2(G))$.  Furthermore,
given $f\in C_c(G\ltimes G)$ and $\phi\in C_c(G)$ we have 
\begin{align*}
U\Phi(f)U^*\phi(s) &= \int_G f(t,s)U^*\phi(t\inv s)\Delta\neghalf(s)dt
\\
&= \int_G f(t,s)\phi(t\inv s)\Delta\neghalf(t)dt = M\rtimes
L(f)\phi(s).
\end{align*}
It follows that $M\ltimes L = U^*\Phi(f)U$ so that $M$ is a
faithful representation onto $\mcal{K}(L^2(G))$.  
\end{proof}


\chapter{Basic Constructions}
\label{cha:basic}
Half of this chapter consists of technical tools and the other half of
capstone results.  In particular, Section \ref{sec:transitive} is a
mixture of both.  We introduce some major results for groupoid
crossed products in this section, but we also apply them and show that
a transitive groupoid crossed product is Morita equivalent to a group
crossed product.  Section \ref{sec:ideals} is entirely
technical but is essential for Section \ref{sec:regularity}.  In
Section \ref{sec:unitary} we introduce the notion of a unitary action
and prove that for unitary actions the crossed product reduces to a
tensor product.  This result is to be expected.  What's more
interesting is the definition of locally unitary actions given in
Section \ref{sec:locally-unitary}.  There we show that locally
unitary actions give rise to principal bundles and that the exterior
equivalence class of the action is characterized by the cohomological
invariant of the associated bundle. Furthermore, we show that locally
unitary actions can be constructed from any principal bundle.  

\section{Transitive Groupoid Crossed Products}
\label{sec:transitive}

In this section, we will compute the representations of crossed
products by transitive groupoids.  This work will be essential in
Section \ref{sec:regularity} and also generalizes the latter half of
\cite{groupoidequiv}.  In order to accomplish our goal we will have to
introduce the notion of an equivalence bundle for groupoid dynamical
systems.  
These objects are intimidating
because they are very complex.  However, as we will see, 
in practice they are fairly easy to work with.  Now, we won't do
anything other than write down the definition and cite the major
theorem for equivalences.  The reader is referred to
\cite{renaultequiv} for the whole story.  It is also mildly
helpful to keep in mind that this material is meant to be a
generalization of groupoid equivalences, which were introduced in 
Section \ref{sec:equivalence}.

\begin{remark}
In what follows we use the following notation.  Given two bundles
${p:E\rightarrow X}$ and $q:F\rightarrow X$ we define 
\[
E*F = \{(e,f)\in E\times F : p(e) = q(f) \}.
\]
\end{remark}

\begin{definition}
\label{def:54}
\index{dynamical system equivalence}
Suppose $(A,G,\alpha)$ and $(B,H,\beta)$ are two separable groupoid
dynamical systems with associated upper-semicontinuous bundles
$\mcal{A}$ and $\mcal{B}$.  An {\em equivalence} between the dynamical systems
$(A,G,\alpha)$ and $(B,H,\beta)$ is an upper-semicontinuous Banach
bundle\footnote{The symbol $\erune$
  is pronounced ``humpf'' as in ``Humpf-Humpf-a-Dumpfer'' \cite{seuss}.} $p\suberune:\erune\rightarrow X$ 
over a $(G,H)$-equivalence $X$ together with
$A(r(x))-B(s(x))$-imprimitivity bimodule structures on each fibre 
$\erune_x = p\suberune\inv(x)$ and commuting {\em strongly} continuous
actions of $G$ and
$H$ on the left and right, respectively, of $\erune$ such that the
following additional properties are satisfied for $e,f\in \erune$,
$a\in\mcal{A}$, $b\in\mcal{B}$, $\gamma\in G$ and $\eta\in H$: 
\begin{enumerate}
\item (Continuity)  The maps induced by the imprimitivity bimodule
  inner products from $\erune*\erune\rightarrow \mcal{A}$ and $\erune*\erune
  \rightarrow \mcal{B}$ are continuous, as are the maps
  $\mcal{A}*\erune\rightarrow \erune$ and $\erune*\mcal{B}\rightarrow
  \erune$ induced by the imprimitivity bimodule actions. 
\item (Equivariance) The bundle map $p\suberune$ is equivariant with
  respect to the groupoid actions.  In other words,
  $p\suberune(\gamma\cdot e) = \gamma\cdot p\suberune(e)$ and
  $p\suberune(e\cdot \eta) = p\suberune(e)\cdot \eta$.  
\item (Compatibility) The groupoid actions are compatible with the
  imprimitivity bimodule structure: 
\begin{align*}
\lset{A}\langle \gamma\cdot e,\gamma\cdot f\rangle &=
\alpha_\gamma(\lset{A}\langle e,f\rangle) &
\langle e \cdot \eta, f\cdot \eta \rangle_B &= \beta_\eta\inv(\langle
e,f \rangle_A) \\
\gamma\cdot(a\cdot e) &= \alpha_\gamma(a)\cdot(\gamma\cdot e) &
(e\cdot b)\cdot\eta &= (e\cdot \eta)\cdot \beta_\eta\inv(b).
\end{align*}
\item (Invariance) The $G$-action commutes with the $B$-action on
  $\erune$ and the $H$-action commutes with the $A$-action.  That is,
  $\gamma\cdot (e\cdot b) = (\gamma\cdot e)\cdot b$ and $(a\cdot
  e)\cdot \gamma = a\cdot (e\cdot\gamma)$.  
\end{enumerate}
\end{definition}

The only reason anyone is inspired to try and form an equivalence bundle is
that, in the same vein as Theorem \ref{thm:groupoidequiv}, we get the
following useful theorem. This result is can be found, in French, 
in \cite[Corollaire 5.4]{frenchrenault} or, in English, in
\cite[Theorem 5.5]{renaultequiv}.

\begin{theorem}[Renault's Equivalence Theorem]
\label{thm:renaultequiv}
Suppose $G$ and $H$ are second countable locally compact Hausdorff
groupoids with Haar systems $\lambda_G$ and $\lambda_H$ respectively.
Furthermore, suppose that $(A,G,\alpha)$ and $(B,H,\beta)$ are
separable dynamical systems and $p\suberune:\erune\rightarrow X$ is an
equivalence between $(A,G,\alpha)$ and $(B,H,\beta)$.  Then $\mcal{Z}_0 =
\Gamma_c(X,\erune)$ becomes a pre-$A\rtimes_\alpha G-B\rtimes_\beta
H$-imprimitivity bimodule with respect to the following operations
for $f\in \Gamma_c(G,r^*\mcal{A})$, $g\in\Gamma_c(H,r^*\mcal{B})$ and
$z,w\in \Gamma_c(X,\erune)$: 
\begin{align}
\label{eq:74}
f\cdot z(x) &:= \int_G f(\gamma)\cdot(\gamma\cdot z(\gamma\inv\cdot
x)) d\lambda_G^{r(x)}(\gamma) \\
\label{eq:75}
z\cdot g(x) &:= \int_H
(z(x\cdot\eta)\cdot\eta\inv)\cdot\beta_\eta(g(\eta\inv))
d\lambda_H^{s(x)}(\eta) \\
\label{eq:76}
\llangle z,w \rrangle_{B\rtimes_\beta H}(\eta) &:= \int_G \langle z(\gamma\inv\cdot x
\cdot \eta\inv),w(\gamma\inv\cdot x)\cdot \eta\inv\rangle_{B(r(\eta))}
d\lambda_G^{r(x)}(\gamma) \\
\label{eq:73}
\lset{A\rtimes_\alpha G}\llangle z,w\rrangle(\gamma) &:= \int_H
\lset{A(r(\gamma))}\langle z(\gamma\cdot x\cdot \eta),\gamma\cdot w(x\cdot
\eta)\rangle d\lambda_H^{s(x)}(\eta)
\end{align}
where $x$ in \eqref{eq:76} is any element of $X$ such that $s(x) =
s(\eta)$ and $x$ in \eqref{eq:73} is any element of $X$ such that
$r(x) = s(\gamma)$.  In particular, the completion $\mcal{Z}$ of
$\mcal{Z}_0$ is a $A\rtimes_\alpha G-B\rtimes_\beta H$-imprimitivity
bimodule and $A\rtimes_\alpha G$ and $B\rtimes_\beta H$ are Morita
equivalent. 
\end{theorem}

This theorem will be essential to our study of induced representations.
We can also give a straightforward application in the case of
transitive groupoid actions.  We will start by
citing a deep theorem by Ramsay.  The following is a slightly
trimmed down transcription of \cite[Theorem 2.1]{groupoiddichotomy}.

\begin{theorem}[Mackey-Glimm Dichotomy]
\label{thm:glimmdich}
\index{Mackey-Glimm Dichotomy}
Let $G$ be a second countable locally compact Hausdorff groupoid.
Then the following are equivalent:
\begin{enumerate}
\item For each $u\in G\unit$, the map $\gamma\cdot S_u\mapsto
  r(\gamma)$ of $G_u/S_u$ to the orbit
  $G\cdot u$ is a homeomorphism.   
\item $G\unit/G$ is a $T_0$-space.  
\item Each orbit is locally closed in $G\unit$.  
\item The quotient topology on $G\unit/G$ generates the quotient Borel
  structure. 
\item The quotient Borel structure on $G\unit/G$ is countably
  separated.  
\item $G\unit/G$ is almost Hausdorff.\footnote{See Definition
    \ref{def:almosthauss}.}
\item $G\unit/G$ is a standard Borel space. 
\item The quotient map $\pi:G\unit\rightarrow G\unit/G$ has a Borel
  section.
\end{enumerate}
\end{theorem}

\begin{proof}[Remark.]
This isn't so much a proof as it is an explanation of how this
statement of the Mackey-Glimm Dichotomy can be obtained from the
statement in the reference.  Those readers not interested in chasing
down citations can skip ahead.  
The result in \cite{groupoiddichotomy} is for Polish groupoids.  Since
second countable locally compact Hausdorff spaces are completely metrizable the
result clearly applies.  Furthermore, in order to use the full power
of \cite[Theorem 2.1]{groupoiddichotomy} we need to know that the orbit
groupoid $R$ is a $F_\sigma$ subset of $G\unit\times G\unit$.
However, $G$ is second countable locally compact Hausdorff so that it
is the countable union of compact sets, say $\{K_i\}_{i=1}^\infty$.  
Consider the
canonical map $\pi:G\rightarrow R_P$.  Since $\pi$ is continuous,
$\pi(K_i)$ is compact (with respect to the product topology) and
therefore closed in $G\unit\times G\unit$.  Thus $R=
\bigcup_{i=1}^\infty \pi(K_i)$ is $F_\sigma$.  Therefore all fourteen (!)
conditions listed in \cite{groupoiddichotomy} are equivalent and the
conditions stated in the theorem are a subset of those. 
\end{proof}  

Theorem \ref{thm:glimmdich} is known as the Mackey-Glimm Dichotomy
because it says that either a number of nice conditions hold or none
of them do. We will primarily be interested in the case where
$G\unit/G$ is $T_0$, but will spend some time on the ``other side'' of
the dichotomy as well.  As we will see in Section
\ref{sec:groupstab}, things don't always work as smoothly as
they do when Theorem \ref{thm:glimmdich} is satisfied.  

The next proposition makes more sense if you realize that given
$\gamma\in G$ conjugation by $\gamma$ defines an isomorphism
from $S_{s(\gamma)}$ onto $S_{r(\gamma)}$.  Thus, in a
transitive groupoid all of the stabilizer subgroups must be
isomorphic.  

\begin{prop}[{\cite[Example 2.2]{groupoidequiv}}]
\label{prop:78}
Suppose $G$ is a transitive second countable locally compact Hausdorff
groupoid.  Then given $u\in G\unit$ the space $X = G_u$ is
a $(G,S_u)$-equivalence with respect to the natural actions of $G$ and $S_u$. 
\end{prop}

\begin{proof}
Let $G$ act on $X=G_u$ by left translation and $S_u$ act on $X$ by
right translation.  It is easy to show that $X$ is a free continuous left
$G$-space and a free continuous right $S_u$-space and that the actions
commute.  Now suppose $\gamma_i\rightarrow \gamma$ in $X$ and
$\eta_i\gamma_i\rightarrow \zeta$ in $X$.  Then
$\eta_i\rightarrow \zeta\gamma\inv$ so that by Proposition
\ref{prop:18} the action of $G$ on $X$ is proper. The same argument
shows that the action of $S_u$ on $X$ is proper as well.  
Furthermore since the source map $s_X$ on $X$ maps onto a single point it
must be open.  We also note that if $\gamma,\eta\in X$ then $\gamma$
and $\eta\inv$ are composable and that $(\gamma\eta\inv)\cdot \eta =
\gamma$.  Thus the action of $G$ on $X$ is transitive so that $s_X$
factors to a bijection.  Next, let $r_X:X\rightarrow G\unit$ be the
restriction of the range map to $X$.  Given $\gamma,\eta\in X$, 
if $r(\gamma) = r(\eta)$ then
$\gamma\inv\eta\in S_u$ so that $\gamma\cdot \gamma\inv\eta = \eta$.
This shows that $r_X$ factors to a bijection from $X/S_u$ onto
$G\unit$.  All that is left is to show that $r_X$ is open.  This is
not necessarily true in the nonseparable case, \cite[Example
2.2]{groupoidequiv}.  However, we assumed that
$G$ is second countable and, since it is
transitive, the quotient space $G\unit/G$ is trivial.  Clearly $G\unit/G$ is
$T_0$ so that we may use Theorem \ref{thm:glimmdich} to conclude that
the map $\rho:G_u/S_u\rightarrow G\cdot u = G\unit$ such that
$\rho(\gamma\cdot S_u) = r(\gamma)$ is a homeomorphism.
However, $r_X$ is just the composition of $\rho$ with the quotient map
from $G_u$ onto $G_u/S_u$.   Since both of these maps are open $r_X$
is open as well. 
\end{proof}

\begin{remark}
It is worth noting that you do not need to use the Mackey-Glimm
Dichotomy to prove Proposition \ref{prop:78}.  A more elementary proof
can be found in \cite[Theorem 2.2A,2.2B]{groupoidequiv}.  This argument
makes use of the Baire Category Theorem and is quite interesting.  
\end{remark}

We can now use Proposition \ref{prop:78} to identify the Morita
equivalence class of transitive groupoid crossed products.   As before,
observe that because $G$ is transitive each of the stabilizer subgroups
must be isomorphic.  Furthermore, the $\alpha_\gamma$ guarantee that all
of the fibres of $\mcal{A}$ are isomorphic as well.  

\begin{theorem}
\label{thm:transprod}
\index{groupoid!transtive}
\index{New Result}
Suppose $(A,G,\alpha)$ is a separable dynamical system and that $G$ is
transitive with Haar system $\lambda$.  Fix $u\in G\unit$ and
let $\beta$ be Haar measure on $S_u$.  Then 
$\mcal{X}_0 = C_c(G_u,A(u))$ becomes a
pre-$A\rtimes_\alpha G-A(u)\rtimes_{\alpha|_{S_u}} S_u$-imprimitivity
bimodule with respect to the following operations for $f\in
\Gamma_c(G,r^*\mcal{A})$, $g\in C_c(S_u,A(u))$ and $z,w\in
\mcal{X}_0$:
\begin{align}
\label{eq:77}
f\cdot z(\gamma) &= \int_G
\alpha_\gamma\inv(f(\eta))z(\eta\inv\gamma)d\lambda^{r(\gamma)}(\eta) \\
\label{eq:78}
z\cdot g(\gamma) &= \int_{S_u}
\alpha_s(z(\gamma s)g(s\inv))d\beta(s) \\
\label{eq:79}
\llangle z,w \rrangle_{A(u)\rtimes S_u}(s) &= 
\int_G z(\eta\inv)^*\alpha_s(w(\eta\inv s))d\lambda^u(\eta)
\\ \label{eq:80}
\lset{A\rtimes G}\llangle z,w \rrangle(\gamma) &= \int_{S_u}
\alpha_{\gamma\zeta s}(z(\gamma\zeta s)w(\zeta s)^*)d\beta(s) 
\end{align}
where $\zeta$ in \eqref{eq:80} is any element of $G_u$ such that
$r(\zeta)= s(\gamma)$.  
In particular, the completion $\mcal{X}$ of $\mcal{X}_0$ is an
$A\rtimes_\alpha G-A(u)\rtimes_{\alpha|_{S_u}}S_u$-imprimitivity
bimodule and $A\rtimes_\alpha G$ is Morita equivalent to 
$A(u)\rtimes_{\alpha|_{S_u}} S_u$. 
\end{theorem}

\begin{proof}
Suppose $(A,G,\alpha)$ is as in the statement of the theorem and that $u\in
G\unit$.  First, let $X=G_u$ and recall that $X$ is a $(G,S_u)$-equivalence with
respect to the natural actions of $G$ and $S_u$ on $X$ by Proposition
\ref{prop:78}.  Next, consider the trivial bundle $\erune = X\times A(u)$.
This is clearly a Banach bundle and the space of compactly supported 
sections can be identified with $\mcal{X}_0 = C_c(G_u,A(u))$.  
Given $\gamma \in X$ we need to equip $\erune_\gamma \cong A(u)$ with an
$A(r(\gamma))-A(u)$-imprimitivity bimodule structure.  However
$\alpha_\gamma:A(u)\rightarrow A(r(\gamma))$ is an isomorphism so that
we can use \cite[Example 3.14]{tfb} to give $\erune_\gamma$ the bimodule
structure with respect to the following operations for $a,b\in A(u)$
and $c\in A(r(\gamma))$:
\begin{align*}
(\gamma,a)\cdot b &= (\gamma,ab) & c\cdot (\gamma,a) &=
(\gamma,\alpha_\gamma\inv(c)b) \\
\langle (\gamma,a),(\gamma,b)\rangle_{A(u)} &= a^*b & 
\lset{A(r(\gamma))}\langle (\gamma,a),(\gamma,b)\rangle 
&= \alpha_\gamma(ab^*)
\end{align*}
Next, we define the range and source maps, as well as the $G$ and
$S_u$ actions, on $\erune$ in the following manner 
for $(\gamma,a)\in\erune$, $\eta\in G$ and $s\in S_u$
\begin{align*}
r\suberune(\gamma,a) &= r(\gamma) & s\suberune(\gamma,a) &= u  \\
\eta\cdot(\gamma,a) &= (\eta\gamma,a) & (\gamma,a)\cdot s &= (\gamma
s, \alpha_s\inv(a)).
\end{align*}
It is relatively obvious that these are continuous actions.  Since
$s\suberune$ maps onto a single point it must be open.  Furthermore,
since $X$ is a $(G,S_u)$-equivalence 
we know that $r_X$ must be open.  However,
$r\suberune$ is just the composition of $r_X$ with the projection map
from $\erune$ onto the first factor and therefore must be open as
well.  Thus the actions of $G$ and $S_u$ on $\erune$ are strongly
continuous in the sense of Definition \ref{def:11}.  The
computation 
\begin{align*}
(\eta\cdot(\gamma,a))\cdot s &= (\eta\gamma,a)\cdot s = 
(\eta\gamma s,\alpha_s\inv(a)) \\
&= \eta\cdot(\gamma s,\alpha_s\inv(a)) = \eta((\gamma,a)\cdot s)
\end{align*}
shows that the actions commute.  Now we must verify the different
requirements for equivalence.  The continuity condition is clear once
you recall that multiplication on $\mcal{A}$ and the action $\alpha$
are both continuous.  Furthermore, it is very easy to show that the bundle map
$p\suberune:\erune\rightarrow X$ is equivariant.  Next, we must show
that the groupoid actions are compatible with the imprimitivity
actions.  Fix $\gamma\in X$, $a,b\in A(u)$, $c\in A(r(\gamma))$, $s\in S_u$
and $\eta\in G$.  We can now perform the following computations
without too much difficulty:
\begin{align*}
\lset{A(r(\eta))}\langle \eta\cdot (\gamma,a),\eta\cdot(\gamma,b)\rangle &=
\alpha_{\eta\gamma}(ab^*) \\ &= \alpha_\eta(\lset{A(r(\gamma))}\langle
(\gamma,a),(\gamma,b)\rangle), \\
\langle (\gamma,a)\cdot s, (\gamma,b)\cdot s\rangle_{A(u)} &=
\langle (\gamma s, \alpha_s\inv(a)),(\gamma
s,\alpha_s\inv(b)\rangle_{A(u)}
=  \alpha_s\inv(a^*b)\\ & =
\alpha_s\inv(\langle(\gamma,a),(\gamma,b)\rangle_{A(u)}), \\
\eta\cdot (c\cdot(\gamma,a)) &= \eta\cdot (\gamma,\alpha_\gamma\inv(c)a)
= (\eta\gamma,\alpha_{\eta\gamma}\inv(\alpha_\eta(c))a) \\ &=
\alpha_\eta(c)\cdot(\eta\cdot (\gamma,a)), \\
((\gamma,a)\cdot b)\cdot s &= (\gamma s, \alpha_s\inv(ab))
= (\gamma s, \alpha_s\inv(a)\alpha_s\inv(b)) \\ &=
((\gamma,a)\cdot s)\cdot \alpha_s\inv(b).
\end{align*}
The last thing to check are the invariance conditions.  We compute
\[
\eta\cdot((\gamma,a)\cdot b) = (\eta\gamma, ab) =
(\eta\cdot(\gamma,a))\cdot b
\]
and 
\begin{align*}
(c\cdot(\gamma,a))\cdot s &= (\gamma,\alpha_\gamma\inv(c)a)\cdot s = 
(\gamma s, \alpha_s\inv(\alpha_\gamma\inv(c)a)) \\
&= (\gamma s, \alpha_{\gamma s}\inv(c)\alpha_s\inv(a)) =
c\cdot (\gamma s, \alpha_s\inv(a)) \\
&= c\cdot((\gamma,a)\cdot s)
\end{align*}
We have now verified all of the conditions for $\erune$ to be an
equivalence between $(A,G,\alpha)$ and $(A(u),S_u,\alpha|_{S_u})$.  We
use Renault's Equivalence Theorem to conclude
that $\mcal{X}_0$ completes into an $A\rtimes G-A(u)\rtimes
S_u$-imprimitivity bimodule.  What's more we can use \eqref{eq:74}
through \eqref{eq:73} to compute the operations on $\mcal{X}_0$
explicitly.  For instance, fix $f\in\Gamma_c(G,r^*\mcal{A})$, $g\in
C_c(S_u,A(u))$ and $z,w\in \mcal{X}_0$.
Observe that if $z(\eta\inv\gamma)=(\eta\inv\gamma, a)$ then
$\eta\cdot z(\eta\inv\gamma) = (\gamma, a)$ so that $f(\eta)\cdot
(\eta\cdot z(\eta\inv\gamma)) = (\gamma,
\alpha_\gamma\inv(f(\eta))a)$.  Making the usual identification of $a$ with
$z(\eta\inv\gamma)$ we get
\begin{align*}
f\cdot z(\gamma) &= \int_G f(\eta)\cdot(\eta\cdot z(\eta\inv\cdot
\gamma))d\lambda^{r(\gamma)}(\eta) \\
&= \int_G \alpha_\gamma\inv(f(\eta))
z(\eta\inv\gamma)d\lambda^{r(\gamma)}(\eta).  
\end{align*}
In a similar fashion we obtain
\begin{align*}
z\cdot g(\gamma) &= \int_{S_u} (z(\gamma\cdot s)\cdot s\inv)\cdot
\alpha_s(g(s\inv))d\beta(s) \\
&= \int_{S_u} \alpha_s(z(\gamma s)g(s\inv))d\beta(s).
\end{align*}
Next, in \eqref{eq:76} we may as well pick $x = u$ so that we have 
\begin{align*} 
\llangle z,w \rrangle_{A(u)\rtimes S_u}(s) &= 
\int_G \langle z(\eta\inv\cdot u \cdot s\inv), 
w(\eta\inv \cdot u)\cdot s\inv\rangle_{A(u)} d\lambda^u(\eta)
\\ &= \int_G z(\eta\inv s\inv)^* \alpha_s(w(\eta\inv))
d\lambda^u(\eta) \\
&= \int_G z(\eta\inv)^* \alpha_s(w(\eta\inv s)) d\lambda^u(\eta).
\end{align*}
Finally, given $\gamma\in G$ as in the definition of \eqref{eq:73} we
choose some $\zeta\in X$ such that $r(\zeta) = s(\gamma)$, and then observe
\begin{align*}
\lset{A\rtimes G} \llangle z,w \rrangle(\gamma) &= 
\int_{S_u} \lset{A(r(\gamma))}\langle z(\gamma\cdot \zeta\cdot
s),\gamma\cdot w(\zeta\cdot s)\rangle d\beta(s) \\
&= \int_{S_u} \alpha_{\gamma\zeta s}(z(\gamma\zeta s)w(\zeta s)^*) d\beta(s)
\end{align*}
Thus, all of the operations on $\mcal{X}_0$ have the correct form and
we are done. 
\end{proof}

\begin{remark}
As we noted in the beginning of the section, equivalence bundles seem
frightening because there are many different operations.  However, in
Theorem \ref{thm:transprod} all of the actions had natural
definitions and working with them posed no real difficulty.  This is
often the case with equivalence bundles.  
\end{remark}

We also take this opportunity, while we are on the topic of
equivalence bundles, to introduce the following technical lemma.  
It may not seem important but approximate identities like this are
very useful for proving nondegeneracy conditions.  They can also be
quite fiddly which is why we only reference this result.  

\begin{remark}
Recall from Definition \ref{def:55} that given a groupoid $G$ a
neighborhood of $G\unit$ is called {\em conditionally compact} if for
every compact $K\subset G\unit$ its intersection with $r\inv(K)$ is
compact.  
\end{remark}

\begin{lemma}[{\cite[Proposition 6.7]{renaultequiv}}]
\label{lem:19}
\index{approximate identity}
Suppose $(A,G,\alpha)$ and $(B,H,\beta)$ are separable groupoid dynamical systems and
$p\suberune:\erune\rightarrow X$ 
is an equivalence between $(A,G,\alpha)$ and $(B,H,\beta)$.
Let $\{a_l\}_{l\in\Lambda}$ be an approximate identity for $A$.  Then for
each 4-tuple $(K,U,l,\epsilon)$ consisting of a compact subset
$K\subset G\unit$, a conditionally compact neighborhood $U$ of
$G\unit$ in $G$, $l\in\Lambda$ and $\epsilon >0$ there is an $e =
e_{(K,U,l,\epsilon)}\in \Gamma_c(G,r^*\mcal{A})$ such that 
\begin{enumerate}
\item $\supp e \subset U$, 
\item $\int_G \|e(\gamma)\|d\lambda^u(\gamma) \leq 4$ for all $u\in K$
  and 
\item $\left\| \int_G e(\gamma)d\lambda^u(\gamma)-a_l(u)\right\|<\epsilon$
    for all $u\in K$.  
\end{enumerate}
Furthermore $\{e_{(K,U,l,\epsilon)}\}$ directed by increasing $K$ and
  $l$ and decreasing $U$ and $\epsilon$ is an approximate identity with
  respect to the inductive limit topology for both the left action of
  $\Gamma_c(G,r^*\mcal{A})$ on itself and the left action of
  $\Gamma_c(G,r^*\mcal{A})$ on $\Gamma_c(X,\erune)$.  
Specifically, $e_{(K,U,l,\epsilon)}*g\rightarrow g$ with respect to the inductive
  limit topology for all $g\in\Gamma_c(G,r^*\mcal{A})$ and
  $e_{(K,U,l,\epsilon)}\cdot z \rightarrow z$ with respect to the
  inductive limit topology for all $z\in \Gamma_c(X,\erune)$.  
\end{lemma} 

\begin{proof}[Remark]
This is just an explanation of how to extract this statement of
the lemma from the reference.  Those readers not concerned with
chasing down citations can skip ahead.  
Lemma \ref{lem:19} is really a result of both \cite[Lemma 6.6, Lemma
6.7]{renaultequiv}. Specifically, \cite[Lemma 6.6]{renaultequiv} 
states that the net $\{e_{(K,l,U,\epsilon)}\}$ given
above is an approximate identity in the inductive limit topology.
Then such a net is constructed in the proof of \cite[Lemma
6.7]{renaultequiv}.  The construction of this approximate identity is
fairly complex and occupies the whole of \cite[Section 6]{renaultequiv}.
\end{proof}


\section{Invariant Ideals}
\label{sec:ideals}
This is mainly a technical section. We will investigate certain nice
ideals of crossed products and show that they behave in a reasonable
manner.  This material is structured along the lines of \cite[Section
3.4]{tfb2}.  We start by considering what happens when we cut down
upper-semicontinuous bundles.  

\begin{definition}
\label{def:48}
\index[not]{$A(Y)$}
Let $A$ be a $C_0(X)$-algebra, $\mcal{A}$ its associated
upper-semicontinuous bundle, and $Y$ a locally compact subset of
$X$.  We define the {\em restriction} of $A$ to $Y$ to be 
\[
A(Y) := \Gamma_0(Y,\mcal{A}).
\]
\end{definition}

\begin{prop}
\label{prop:81}
Suppose $A$ is a $C_0(X)$-algebra, $\mcal{A}$ its associated
upper-semi\-continu\-ous bundle, and $Y$ is a nonempty locally compact
subset of $X$.  Then $A(Y)$ is a $C_0(Y)$-algebra with $A(Y)(y) =
A(y)$ for all $y\in Y$.  Furthermore the upper-semicontinuous bundle
associated to $A(Y)$ is $\mcal{A}|_Y = p\inv(Y)$.  
\end{prop}

\begin{proof}
Let $\mcal{B} = p\inv(Y)$.  It
follows immediately from the fact that $\mcal{A}$ is an
upper-semicontinuous $C^*$-bundle that $\mcal{B}$ is as well.  Observe
that, by definition, $A(Y)$ is the space of sections of $\mcal{B}$
which vanish at infinity.  It follows from Example \ref{ex:16}
that $A(Y)$ is a $C_0(Y)$-algebra.  Now, once we make the usual
identification of $A(Y)(y)$ with the fibre $\mcal{B}_y$,
we have $A(Y)(y) = \mcal{B}_y = \mcal{A}_y = A(y)$.  Finally, still
using this identification, it is vacuously true that $\mcal{B}$ has the
unique topology making $A(Y)$ into the space of sections of $\mcal{B}$
which vanish at infinity.  Thus, by definition, $\mcal{B}$ is the
upper-semicontinuous bundle associated to $A(Y)$.  
\end{proof}

Now, the notation $A(Y)$ is going to be subject to the same abuse that
we subject $a(x)$ to.  Specifically, the next proposition shows that
for closed sets $A(Y)$ is (isomorphic to) a quotient of $A$, bringing us back to
definition of the similarly denoted $A(x)$.  

\begin{remark}
\label{rem:7}
In the same vein as Definition \ref{def:29}, if $C$ is a closed subset
of $X$ and $J_C$ is the ideal of functions in $C_0(X)$ which vanish on
$C$ then we will denote the ideal $\overline{\Phi_A(J_C)\cdot A}$ by
$I_C$.  
\end{remark}

\begin{prop}
\label{prop:71}
Suppose $A$ is a $C_0(X)$-algebra and $\mcal{A}$ is its associated
upper-semicontinu\-ous bundle.  Let $U$ be an open subset of $X$ and $C
= X\setminus U$.  Then there are $*$-homomorphisms
$\iota:A(U)\rightarrow A$ and $\rho:A\rightarrow A(C)$ where $\iota$
is the inclusion of $\Gamma_0(U,\mcal{A})$ into $\Gamma_0(X,\mcal{A})$
and $\rho$ is the restriction map from $\Gamma_0(X,\mcal{A})$ into
$\Gamma_0(C,\mcal{A})$.  Furthermore, the sequence 
\[
\begin{CD} 
0 @>>> A(U) @>\iota>> A @>\rho>> A(C) @>>>0 
\end{CD}
\]
is short exact.  Finally, $\ran \iota = I_C$ so that $A(U)$ is
isomorphic to $I_C$ and $A(C)$ is isomorphic to the quotient $A/I_C$.  
\end{prop}

\begin{proof}
We will identify $A$ with
the sections of $\mcal{A}$ which vanish at infinity.  
Recall that, by definition, $A(U)$ and $A(C)$ are the sections on $U$ and
$C$ which vanish at infinity.  Our first goal is to see that $\iota$
and $\rho$ are well defined.  Suppose $f\in A(U)$ and let $\iota(f)$
be the extension of $f$ to $X$ by letting $\iota(f)(x) = 0_x$ for all
$x\not\in U$.  We claim that $\iota(f)$ is continuous on $X$.  Suppose
$x_i\rightarrow x$ and $x\in U$.  Then, eventually, $x_i\in U$ so that
$\iota(f)(x_i) = f(x_i)\rightarrow f(x) = \iota(f)(x)$.  Now suppose
$x\not\in U$.  As we have done before, we will pass to a subnet and
show that a sub-subnet converges to $\iota(f)(x)=0_x$.  If
$x_i\not\in U$ infinitely often then pass to a subnet and assume this
is always true.  It follows from condition (d) of
Definition \ref{def:26} that $\iota(f)(x_i) = 0_{x_i}\rightarrow
0_x$.  If $x_i\in U$ eventually then pass to a subnet and assume this
is always true.  If $\epsilon > 0$ then, since $f$ vanishes at
infinity, $K=\{x\in U : \|f(x)\|\geq \epsilon\}$ is a compact subset
of $U$ and is therefore
closed in $X$.  Thus, $X\setminus K$ is an open neighborhood of
$x$ and therefore eventually contains $x_i$.  However, this implies that
eventually $\|f(x_i)\|<\epsilon$.  Since $\epsilon >0$ was arbitrary we
must have $\|f(x_i)\|\rightarrow 0$ and it follows 
that $f(x_i)\rightarrow 0_x$.  Thus $\iota(f)$ is continuous. This
also shows that $\iota(f)$ vanishes at infinity in $X$, since the
set $K$ above is compact as a subset of $X$.  Finally,
it is clear that $\iota$ is a $*$-homomorphism under the pointwise
operations and that $\iota$ is isometric with respect to the uniform
norm.  It follows that $\iota$ is injective and an
isomorphism onto its range.  

Clearly $\rho(f)$ is a continuous
section on $C$.  Furthermore if $K = {\{x:\|f(x)\|\geq \epsilon\}}$ for
some  $\epsilon > 0$ then $K\cap C$ is a closed subset of a compact
set and is therefore compact.  It follows that $\rho(f)$ vanishes at
infinity.  Obviously, $\rho$ is a $*$-homomorphism with respect
to the pointwise operations and  $\|\rho(f)\|_\infty \leq
\|f\|_\infty$.  Now consider $\ran\rho$.  We would like to show that
$\ran\rho$ is closed under the $C_0(C)$-action.  Well, $\ran\rho$
is closed and if $\phi_i\rightarrow \phi$ uniformly in $C_0(C)$ then
$\phi_i\cdot a \rightarrow \phi\cdot a$ uniformly in $A(C)$.  It follows that it
suffices to show that $\ran\rho$ is closed under the action of
$C_c(C)$.  Given $\phi\in C_c(C)$ use Lemma \ref{lem:8} to extend
$\phi$ to $C_c(X)$.  We then observe that $\rho(\phi\cdot f) =
\phi\cdot\rho(f)$ for all $f\in A$.  Hence $\ran\rho$ is closed under
the action of $C_0(C)$.  Now, given $x\in C$ and $a\in A(C)(x)=A(x)$
we pick $f\in A$ such that $f(x) = a$.  Then $\rho(f)(x) = f(x) = a$.
Thus $\ran\rho$ is fibrewise dense.  It follows from Proposition
\ref{prop:42} that $\ran\rho$ is dense in $A(C)$.  Since $\ran\rho$ is
closed it follows that $\rho$ is surjective.  

All we have to do now is show that $\ran\iota = \ker\rho = I_C$.
If $f\in A(U)$ then $\rho(\iota(f))(x) = \iota(f)(x) = 0_x$ for all
$x\in C=X\setminus U$.  Thus $\ran\iota \subset \ker\rho$.  Now
suppose that $f\in \ker\rho$ so that $f$ is zero on $C$.  (Note: The
following argument is slightly easier if we are in the second
countable case.)  Given
$\epsilon > 0$ we know that $K=\{x\in X: \|f(x)\|\geq \epsilon/2\}$ is
compact. Using the fact that $K$ is a compact set which is disjoint
from $C$ we can find a relatively compact open neighborhood $U$ of
$K$ which is disjoint from $C$.  Use Urysohn's Lemma to find a
function $\phi\in C_c(\overline{U})$ such that $\phi$ is one $K$ and
zero on $\overline{U}\setminus U$.  Now extend $\phi$ to
$\overline{U}\cup C$ by letting $\phi$ be zero on $C$.  If
$x_i\rightarrow x$ and $x\in U$ then $x_i$ is eventually in $U$ and
$\phi(x_i)\rightarrow \phi(x)$ by construction.  Suppose
$x_i\rightarrow x$ and $x\notin U$.  Then if $x_i\not\in U$
infinitely often pass to a subnet and assume this is always true.  It
follows that $\phi(x_i) = 0\rightarrow 0 = \phi(x)$. If eventually
$x_i\in U$ pass to a subnet and assume that this is always true.  But
then $x\in \overline{U}$ and $\phi(x_i)\rightarrow \phi(x)$ by
assumption.  Thus, our extension $\phi$ is a continuous, compactly
supported function on $C\cup \overline{U}$.  We can now use Lemma
\ref{lem:8} to extend $\phi$ to a function in $C_c(X)$.  Since $\phi$
was constructed to be zero on $C$ we have $\phi\cdot f \in I_C$.
Furthermore it is easy to check that $\|\phi\cdot f - f\|_\infty <
\epsilon$.  Since $\epsilon > 0$ was arbitrary it follows that $f\in
I_C$ and thus $\ker\rho\subset I_C$. Finally, given $\phi\in C_c(X)$
such that $\phi$ is zero on $C$ and $f\in A$ observe that $\phi\cdot
f$ is zero on $C$.  Hence, if $\epsilon > 0$ then $K=\{x\in
X:\|\phi\cdot f(x)\|\geq \epsilon\} \subset U$.  It follows that the
restriction of $\phi\cdot f$ to $U$, denoted $g$, vanishes at
infinity.  Since $\iota(g) = \phi\cdot f$ it follows that
$I_C\subset \ran\iota$.  Thus $\ran\iota = \ker\rho = I_C$ and the
rest of the proposition follows. 
\end{proof}

As a corollary to the above we find that the spectra of these
``restriction'' bundles are well behaved. 

\begin{corr}
\label{cor:9}
Suppose $A$ is a $C_0(X)$-algebra, $U$ is an open subset of $X$ and $C
= X\setminus U$.  Let $\sigma:\widehat{A}\rightarrow X$
be the associated map on the spectrum.  Then $\sigma\inv(U) \cong
(A(U))\sidehat$ and $\sigma\inv(C)\cong (A(C))\sidehat$.  Furthermore,
if we view both $A(U)$ and $\sigma\inv(U)$ as the disjoint union
$\coprod_{x\in U} A(x)\sidehat$ then this identification is given
by the identity.  The corresponding statement holds for $C$.  
\end{corr}

\begin{proof}
The map $\iota$ is an isomorphism of $A(U)$ onto 
$I_C$ and as such we may identify 
$\widehat{I_C}$ and $(A(U))\sidehat$ 
via the map $\phi_1(\pi)=\pi\circ\iota$.  
Furthermore,
since $I_C$ is an ideal in $A$ we can, and do, identify the set
$\{\pi\in\widehat{A}:I_C\not\subset \ker\pi\}$ with $\widehat{I_C}$
\cite[Proposition A.26]{tfb} via the map $\phi_2(\pi)=\pi|_{I_C}$.    
We would like to show that $\widehat{I_C} =
\sigma\inv(U)$.  If $\pi\in\sigma\inv(U)$ then, by definition,
$I_x \subset \ker\pi$ for some $x\in U$ so that $\pi$ factors to an
irreducible representation $\pi'$ of $A(x)$.  Since $\pi'$
is irreducible it must be non-zero, 
so there exists $a\in A(x)$ such that $\pi(x) \ne 0$.  However, if we
choose $b\in A$ such that $b(x) = a$ and  $\phi\in C_c(U)$ such that
$\phi(x) = 1$ then $\phi\cdot b\in I_C$ and
\[
\pi(\phi\cdot b) = \pi'(\phi(x)b(x)) = \pi(a) \ne 0.
\]
Thus $\pi\in \widehat{I_C}$.  Next suppose $\pi\in \widehat{I_C}$.
By construction, $I_{\sigma(\pi)}\subset\ker \pi$.  If
$\sigma(\pi)\in C$ then $I_C\subset I_{\sigma(\pi)} \subset\ker\pi$
but this is a clear contradiction. 
At this point we can form the homeomorphism  $\phi= \phi_2\circ\phi_1$ mapping
$\sigma\inv(U)$ onto $A(U)$ defined via $\phi(\pi) =
\pi|_{I_C}\circ\iota$.  
Fix $x\in U$ and suppose $\pi\in A(x)\sidehat$.  Let $\pi'$
be the lift of $\pi$ to $A$ and $\pi''$ the lift of $\pi$ to $A(U)$. 
Observe that 
\[
\phi(\pi')(a) = \pi'(\iota(a)) = \pi(\iota(a)(x)) = \pi(a(x)) = \pi''(a).
\]
It follows
that if we identify both $A(U)$ and $\sigma\inv(U)$ with $\coprod_{x\in
  U}A(x)\sidehat$ then $\phi$ is given by the identity map.  

Next, observe that the restriction map 
$\rho$ factors to an isomorphism $\bar{\rho}$ of $A/I_C$ onto $A(C)$
and thus we can define a homeomorphism from $(A/I_C)\sidehat$ onto
$(A(C))\sidehat$ by $\phi_1(\pi) =  \pi\circ \bar{\rho}\inv$.  
Furthermore we can, and do, identify the set
$\{\pi\in\widehat{A} : {I_C\subset \ker\pi}\}$ with $(A/I_C)\sidehat$ 
\cite[Proposition A.28]{tfb} via the map $\phi_2(\pi) = \bar{\pi}$
where $\bar{\pi}$ is the factorization of $\pi$ to $A/I_C$.  
We would like to show that
$(A/I_C)\sidehat=\sigma\inv(C)$.  Suppose  $\pi\in \widehat{A}$ such that
$I_C\subset\ker\pi$ and suppose $x=\sigma(\pi)\not\in C$.  By definition $\pi$
factors to an irreducible representation $\pi'$ of $A(x)$ and
therefore there exists
$a\in A(x)$ such that $\pi'(a) \ne 0$.  Now find $b\in A$ such that
$b(x)=a$ and $\phi\in C_c(U)$ such that $\phi(x) = 1$.  Then
$\phi\cdot b \in I_C$ so that
\[
0=\pi(\phi\cdot b) = \pi'(\phi(x)b(x)) = \pi'(a)
\]
which is a contradiction.  Now suppose $\pi\in\widehat{A}$ such that
$\sigma(\pi)\in C$.  Then $I_C\subset I_x \subset\pi$ and we are
done. Thus we can define a homeomorphism from 
$\sigma\inv(C)$ onto $A(C)\sidehat$
via $\phi = \phi_1\circ\phi_2$.  Furthermore given $\pi\in
\sigma\inv(C)$ we have $\phi(\pi) = \bar{\pi}\circ \bar{\rho}\inv$.  Now
fix $x\in C$ and suppose $\pi\in A(x)\sidehat$.  Let $\pi'$ be the
lift of $\pi$ to $A$ and $\pi''$ the lift of $\pi$ to $A(C)$.  
If $q:A\rightarrow A/I_C$ is the quotient map and $a\in A$ then
\[
\phi(\pi')(\bar{\rho}(q(a))) = \overline{\pi'}(q(a)) = \pi'(a) =
\pi(a(x)) = \pi(\rho(a)(x)) = \pi''(\rho(a)) = \pi''(\bar{\rho}(q(a))).
\]
Since $\bar{\rho}\circ q$ is surjective this implies that $\phi(\pi')
= \pi''$.  Thus we can view $\phi$ as
the identity map.
\end{proof}

Now we will extend some of this theory to groupoid crossed products.
In particular we need to be able to cut down groupoids by restricting
the unit space.   

\begin{definition}
\index{groupoid}
Suppose $G$ is a locally compact Hausdorff groupoid and $Y$ is a
locally compact $G$-invariant subset of $G\unit$. Then we define the {\em
  restriction} of $G$ to $Y$ to be $G|_Y := r\inv(Y)=s\inv(Y)$.  
\end{definition}

\begin{prop}
\label{prop:72}
\index{Haar system}
Let $G$ be a locally compact Hausdorff groupoid with Haar system $\lambda$
and $Y$ a locally compact $G$-invariant subset of $G\unit$.
Then $G|_Y$ is a locally compact Hausdorff subgroupoid of
$G$.  Furthermore the restriction of the Haar system to $G|_Y$ is a
Haar system.  We will always equip $G|_Y$ with this Haar system.  
\end{prop}

\begin{proof}
Since $Y$ is locally compact in $G\unit$ it must be the intersection of an
open set and a closed set \cite[Lemma 1.25,1.26]{tfb2}.  It follows
that $G|_Y$ is locally compact and, since $Y$
is $G$-invariant, that $G|_Y$ is a subgroupoid of $G$ such that $(G|_Y)^u =
G^u$ for all $u\in Y$.  Consider the
restriction of the Haar system to $G|_Y$.  Because
$G|_Y$ is the restriction of $G$ to a $G$-invariant subset, $\supp
\lambda^u = (G|_Y)^u = G^u$ for all $u\in Y$.  Given a compactly
supported function $f\in C_c(G|_Y)$ we can extend $f$ to $G$ using
Lemma \ref{lem:8} and then
the continuity of the Haar system follows.  Finally,
left-invariance is immediate from the fact that $\lambda$ is a Haar
system for $G$. 
\end{proof}

Of course, we would like to be able to restrict a groupoid action to
$G|_Y$, but this means restricting $A$ as well.  The result is the following

\begin{prop}
Suppose $(A,G,\alpha)$ is a groupoid dynamical system and $Y$ is a
locally compact $G$-invariant subset of $G\unit$.  Then the
restriction of $\alpha$ to $G|_Y$ defines an
action of $G|_Y$ on $A(Y)$. 
\end{prop}

\begin{proof}
The fact that $G|_Y$ is a groupoid with a Haar system follows from
Proposition \ref{prop:72} since $Y$ is $G$-invariant.
Since $A(Y)(y)
= A(y)$, it is easy to see that $\alpha$ satisfies the first two
conditions of an action of $G|_Y$ on $A(Y)$.  On the other hand, since the
bundle associated to $A(Y)$ is just the restriction of the bundle
associated to $A$, the continuity condition is easy to verify as well. 
\end{proof}

Now that we can restrict actions we will show, similar to
Proposition \ref{prop:71}, that the crossed product respects the
restrictions.  In particular, we will have to work with the following
object.  

\begin{definition}
\label{def:49}
\index[not]{$\Ex(U)$}
Given a separable groupoid dynamical system $(A,G,\alpha)$ and an open 
$G$-invariant
subset $U$ of $G\unit$ let $\Ex(U)$ be the closure in $A\rtimes_\alpha
G$ of the set 
\[
\{f\in \Gamma_c(G,r^*\mcal{A}) : \supp f \subset G|_U\}.
\]
\end{definition}

There are more elementary ways to prove the next result, however the
following proof, inspired by \cite[Lemma 3.3.1]{lisathesa}, is pretty slick.  

\begin{prop}
\label{prop:73}
Given a separable groupoid dynamical system $(A,G,\alpha)$ and an open
$G$-invariant subset $U$ of $G\unit$ then the inclusion map
\[
\iota:\Gamma_c(G|_U,r^*\mcal{A})\rightarrow \Gamma_c(G,r^*\mcal{A})
\]
extends to an isomorphism of $A(U)\rtimes_{\alpha} G|_U$ onto
$\Ex(U)$.  Furthermore, $\Ex(U)$ is an ideal in $A\rtimes G$.  
\end{prop}

\begin{proof}
First, observe that $G|_U$ is an open subset of $G$ so that we can apply
Proposition \ref{prop:71} to $r^*A$ and see that $\iota(f)$ is
continuous for all $f\in \Gamma_c(G|_U,r^*\mcal{A})$.
However, $\supp\iota(f) = \supp f$ and it follows that $\iota(f)\in
\Ex(U)$.  Since the action of $G|_U$ on $A(U)$ is just the restriction
of the action of $G$ on $A$, it is straightforward to show that $\iota$
is a $*$-homomorphism.  Furthermore, it is clear that $\iota$ is
continuous with respect to the inductive limit topology.  Next, observe that if
$f,g\in \Gamma_c(G,r^*\mcal{A})$ such that $\supp f\subset G|_U$ then for all
$\gamma\not\in G|_U$ we have
\begin{equation}
\label{eq:51}
f*g(\gamma) = \int_G f(\eta)\alpha_\eta(g(\eta\inv\gamma))
d\lambda^{r(\gamma)}(\eta) = 0
\end{equation}
since $r(\eta)=r(\gamma)\in U$ implies $\eta\in G|_U$.  Similarly we
find that $g*f(\gamma) = 0$ for all $\gamma\not\in G|_U$.  
This shows that $g*f,f*g\in \Ex(U)$ and it is
enough to imply that $\Ex(U)$ is an ideal in $A\rtimes G$.
Furthermore, if $f\in \Gamma_c(G,r^*\mcal{A})$ such that $\supp f
\subset G|_U$ then we can clearly view $f$ as a compactly supported
section on $G|_U$ and in this case $\iota(f) = f$.  Thus $\iota$ maps
onto a dense subset of $\Ex(U)$.  

Now, it follows from Corollary \ref{cor:24} that $\iota$ is
bounded and extends to a $*$-homomorphism on $A(U)\rtimes G|_U$.
Furthermore, it is clear that $\ran \iota = \Ex(U)$.  
Next, suppose $R$ is a faithful representation of $A(U)\rtimes G|_U$ on a
separable Hilbert space $\mcal{H}$. 
Let $\mcal{H}_0 = \spn\{R(f)h:
f\in\Gamma_c(G|_U,r^*\mcal{A}),h\in\mcal{H}\}$ and observe that
$\mcal{H}_0$ is dense in $\mcal{H}$.  If  
$f\in\ran\iota$ and $g\in \Gamma_c(G, r^*\mcal{A})$ then we know from
\eqref{eq:51} that $f*g(\gamma)
= 0$ unless $r(\gamma)\in r(\supp f) \subset U$.  Thus $r(\supp f*g)
\subset U$ so that $\supp f*g \subset G|_U$.  In particular, we
can view $f*g$ as a function in $\Gamma_c(G|_U,r^*\mcal{A})$.  We
define a representation of $\Gamma_c(G,r^*\mcal{A})$ on
$\mcal{H}_0$ via 
\begin{equation}
\label{eq:68}
T(f)\sum_{i=1}^n R(g_i)h_i = \sum_{i=1}^n R(f*g_i)h_i.  
\end{equation}
Of course, we need to check that $T$ is well defined.  It will suffice
to show that if $\sum_i R(g_i)h_i = 0$ then $T(f) = 0$ for all $f\in
\Gamma_c(G,r^*\mcal{A})$.  Let $\{e_\kappa\}\subset
\Gamma_c(G|_U,r^*\mcal{A})$ be the approximate identity
from Lemma \ref{lem:19} so that $e_\kappa*g_i\rightarrow g_i$ with respect
to the inductive limit topology for all $i$.  We now have,
\begin{align*}
\sum_{i=1}^n R(f*g_i)h_i &= \sum_{i=1}^n R(f*\lim_\kappa (e_\kappa*g_i)) \\
&= \sum_{i=1}^n R(\lim_\kappa(f*e_\kappa)*g_i) h_i  \\
&= \lim_\kappa R(f*e_\kappa) \sum_{i=1}^n R(g_i)h_i = 0
\end{align*}
where each limit, except for the last, is taken in the inductive limit
topology.  Thus $T$ is well defined and it is easy to see that it is a
homomorphism into the algebra of linear operators on $\mcal{H}_0$.
Next, observe that if $f\in \Gamma_c(G|_U,r^*\mcal{A})$ then
$T(\iota(f)) = R(f)$.  We
now verify the conditions of Theorem \ref{thm:disintigration}.
First, we have
\[
T(\iota(e_k))R(f)h = R(e_k*f)h\rightarrow R(f)h
\]
for all $f\in\Gamma_c(G|_U,r^*\mcal{A})$ and $h\in\mcal{H}$.  
Since $R$ is nondegenerate
this suffices to show that the set
$\spn\{T(f)k:f\in\Gamma_c(G,r^*\mcal{A}),k\in\mcal{H}_0\}$ is dense in
$\mcal{H}$.  In order to verify the continuity condition it clearly
suffices to show that 
\[
f\mapsto (T(f)R(g)h,R(k)l)
\]
is continuous with respect to the inductive limit topology for all
$g,k\in \Gamma_c(G|_U,r^*\mcal{A})$ and $h,l\in \mcal{H}$.  However,
$(T(f)R(g),R(k)l) = (R(f*g)h,R(k)l)$ so the continuity follows from
the fact that $R$ and convolution are both continuous with respect to the
inductive limit topology.  Lastly, for the third condition it will
suffice to check that 
\[
(T(f)R(g)h,R(k)l) = (R(g)h,T(f^*)R(k)l)
\]
for all $g,h$ and $k,l$ as before.  We can compute
\begin{align*}
(T(f)R(g)h,R(k)l) &= (R(f*g)h,R(k)l) = (h,R((f*g)^**k)l) \\
&= (h,R(g^**f^**k)l) = (R(g)h,R(f^**k)l) \\
&= (R(g),T(f^*)R(k)l). 
\end{align*}
Notice that we are being a little schizophrenic about which algebra
the convolution and involution are actually occurring in.  However,
since the Haar system and action of $(A(U),G,\alpha|_U)$ are the
restrictions of the Haar system and action of $(A,G,\alpha)$, the
convolution and involution formulas for both algebras are the same.
In any case, by Theorem \ref{thm:disintigration}, $T$ is bounded with respect
to the universal norm and
extends to a representation of $A\rtimes
G$.  Furthermore, since $R = T\circ\iota$ on
$\Gamma_c(G|_U,r^*\mcal{A})$ this identity holds in general.  Thus,
given $f\in A(U)\rtimes G|_U$ we have 
\[
\|f\| = \|R(f)\| = \|T(\iota(f))\| \leq \|\iota(f)\|.
\]
It follows that $\iota$ is isometric and we are done.  
\end{proof}

The complement to Proposition \ref{prop:73} is the following 

\begin{prop}
\label{prop:74}
Suppose $(A,G,\alpha)$ is a separable groupoid dynamical system and $C$ is a
closed $G$-invariant subset of $G\unit$.  Then the restriction map
\[
\rho:\Gamma_c(G,r^*\mcal{A})\rightarrow \Gamma_c(G|_C,r^*\mcal{A})
\]
extends to a surjective homomorphism from $A\rtimes G$ onto
$A(C)\rtimes G|_C$. Furthermore, $\rho(\Gamma_c(G,r^*\mcal{A}))$ is
dense in $\Gamma_c(G|_C,r^*\mcal{A})$ with respect to the inductive
limit topology.  
\end{prop}

\begin{proof}
Since $C$ is closed, $\rho(f)$ is compactly supported in $G|_C$ for all
$f\in \Gamma_c(G,r^*\mcal{A})$.  Thus $\rho$ is well defined and it is
straightforward to see that it is a $*$-homomorphism which is
continuous with respect to the inductive limit topology.  It follows
from Proposition \ref{prop:62} that $\rho$ 
extends to a $*$-homomorphism on $A\rtimes G$.
It follows from Proposition \ref{prop:71} that the restriction map
from $C_0(G)$ onto $C_0(G|_C)$ is surjective.  Given $\phi\in
C_0(G|_C)$ use the aforementioned surjectivity to extend $\phi$ to a
function $\tilde{\phi}\in C_0(G)$ and 
observe that $\phi\cdot\rho(f) = \rho(\tilde{\phi}\cdot
f)$.  Thus, $\ran\rho$ is closed under the action of $C_0(G|_C)$.
Next, given $\gamma \in G|_C$ and $a\in A(r(\gamma))$ choose 
$b\in A$ such that $b(r(\gamma)) = a$ and $f\in C_c(G)$ such that
$f(\gamma) = 1$.  Then $\rho(f\cdot b)(\gamma) = f(\gamma)b(r(\gamma))
= a$.  Thus $\ran\rho$ is fibrewise dense and therefore dense in the
uniform norm by Proposition \ref{prop:42}.  
What's more, we can perform the standard trick of
multiplying a uniformly converging sequence by an appropriately
supported function in $C_c(G)$ to see that $\ran\rho$ is dense in the
inductive limit topology.  Hence $\ran\rho$ is dense in $A(C)\rtimes
G|_C$ and therefore $\rho$ is surjective. 
\end{proof}

Now we can put everything together to get a nice result mimicking
\cite[Proposition 3.19]{tfb2}.  

\begin{theorem}
\label{thm:invtideal}
Let $(A,G,\alpha)$ be a separable groupoid dynamical system, $U$ an
open $G$-invariant subset of $G\unit$, and $C$ the closed
$G$-invariant set $G\unit\setminus U$.  Then inclusion and restriction
extend to $*$-homomorphisms $\iota:A(U)\rtimes G|_U\rightarrow
A\rtimes G$ and $\rho:A\rtimes G\rightarrow A(C)\rtimes G|_C$,
respectively.  Furthermore, the following sequence is short exact
\[
\begin{CD} 
0 @>>> A(U)\rtimes G|_U @>\iota>> A\rtimes G @>\rho>> A(C)\rtimes G|_C @>>>0 
\end{CD}
\]
and $\ran\iota = \ker\rho = \Ex(U)$ so that $A(C)\rtimes G|_C$ is
isomorphic to the quotient space $A\rtimes G/\Ex(U)$.  
\end{theorem}

\begin{proof}
Proposition \ref{prop:73} shows that $\iota$ is well defined and
injective, and Proposition \ref{prop:74} shows that $\rho$ is well
defined and surjective.  All that is left is to show that
$\ran\iota = \ker\rho = \Ex(U)$.  We have already shown that
$\ran\iota = \Ex(U)$.  Furthermore, it is clear that given $f\in
\Gamma_c(G|_U,r^*\mcal{A})$ we have $\rho\circ\iota(f) =
0$.  It follows that $\ran\iota \subset \ker\rho$.  Thus we are
reduced to proving that $\ker\rho \subset \Ex(U)$.  

Let $R$ be a representation of $A\rtimes G$ such that $\ker R =
\Ex(U)$.  Now, suppose $f,g\in \Gamma_c(G,r^*\mcal{A})$ such that
$\rho(f)=\rho(g)$.  Unfortunately, 
just because $f-g$ is zero on $G|_C$ doesn't mean
$f-g$ is supported on $G|_U$.  However, consider $K = \{\gamma\in G :
\|f(x)-g(x)\| \geq \epsilon\}$.  Since $K$ is a closed subset of $\supp
(f-g)$ it must be compact.  Furthermore, since $K$ is compact and
disjoint from $G|_C$ we can find some relatively compact neighborhood $V$
of $K$ such that $K\subset V \subset G|_U$.  Now choose $\phi\in
C_c(G)$ such that $\phi$ is one on $K$ and zero off $V$.  Then
$\phi\cdot (f-g)$ is supported inside $G|_U$ so that
$\phi\cdot(f-g)\in \Ex(U)$.  However, it is easy to see that
$\|\phi\cdot(f-g)-(f-g)\|_\infty < \epsilon$.  Since $\supp\phi\cdot
(f-g)\subset \supp f-g$ we can use this construction to find a
sequence in $\Ex(U)$ which converges to $f-g$ in the inductive limit
topology.  It follows that $f-g\in \Ex(U)= \ker R$.  Thus the
representation $T$ of $\Gamma_c(G|_C,r^*\mcal{A})$ given by 
\[
T(\rho(f)) = R(f)
\]
is well defined.  Furthermore, since $R$ and $\rho$ are
$*$-homomorphisms, it follows that $T$ is as well.  
We would like to see that $T$ is $I$-norm
decreasing.  Suppose $f\in \Gamma_c(G,r^*\mcal{A})$ and fix
$\epsilon > 0$.  Since Lemma \ref{lem:21} implies that $u\mapsto
\int_G \|f(\gamma)\|d\lambda^u(\gamma)$ is upper-semicontinuous we can
find for each $v\in r(\supp f)\cap C$ some relatively compact 
open set $O_v$ such that
$w\in O_v$ implies 
\begin{equation}
\int_G \|f(\gamma )\|d\lambda^w(\gamma) \leq \int_G \|f(\gamma)\|
d\lambda^v(\gamma) + \epsilon \leq \|\rho(f)\|_I + \epsilon.
\end{equation}
By considering the continuous compactly supported function
$\gamma\mapsto f(\gamma\inv)$ we can, in the same fashion, 
also find for each $v\in r(\supp
f)\cap C$ some relatively compact 
open set $V_v$ such that $w\in V_v$ implies 
\begin{equation}
\int_G \|f(\gamma)\|d\lambda_w(\gamma) \leq \int_G \|f(\gamma)\|
d\lambda_v(\gamma) +\epsilon \leq \|\rho(f)\|_I + \epsilon.
\end{equation}
Since $\{O_v\cap V_v\}$ is an open cover of the compact set
$r(\supp f)\cap C$, there exists some finite subcover $\{O_{v_i}\cap
V_{v_i}\}$.  Let $O = \bigcup_{i} O_{v_i}\cap V_{v_i}$ and observe
that, because the union is finite, $O$ is still relatively
compact.  Now choose $\phi\in C_c(G\unit)$ such that $\phi$ is one on
$r(\supp f)\cap C$, zero off $O$, and $0\leq \phi\leq 1$, and define
$g\in \Gamma_c(G,r^*\mcal{A})$ by $g(\gamma) =
\phi(r(\gamma))f(\gamma)$.  If $v\in O$ then, by construction,
\begin{align*}
\phi(v) \int_G \|f(\gamma)\| d\lambda^v(\gamma) &\leq 
\|\rho(f)\|_I + \epsilon,\quad\text{and} \\
\phi(v) \int_G \|f(\gamma)\| d\lambda_v(\gamma) &\leq
\|\rho(f)\|_I + \epsilon.
\end{align*}
Furthermore, if $v\not\in O$ then 
\begin{align*}
\phi(v) \int_G \|f(\gamma)\| d\lambda^v(\gamma) & = 0 \leq 
\|\rho(f)\|_I + \epsilon,\quad\text{and} \\
\phi(v) \int_G \|f(\gamma)\| d\lambda_v(\gamma) & = 0 \leq
\|\rho(f)\|_I + \epsilon.
\end{align*}
It follows that $\|g\|_I \leq \|\rho(f)\|_I + \epsilon$.  However,
$g- f$ is zero on $C$ by construction so that 
\[
\|T(\rho(f))\| = \|R(f)\| = \|R(g)\| \leq \|g\| \leq \|g\|_I \leq
\|\rho(f)\|_I + \epsilon.
\]
Since $\epsilon$ was chosen arbitrarily, this implies that 
$T$ is $I$-norm decreasing.  
Since $T$ is an $I$-norm decreasing $*$-representation, it
follows that $T$ extends to a representation of $A(C)\rtimes G|_C$.
Finally, since the identity $T\circ \rho = R$ holds on a dense
subset it holds everywhere.  Thus $\ker\rho \subset \ker R = \Ex(U)$
and we are done.  
\end{proof}

\begin{remark}
This section serves as an excellent demonstration of the fact that
kernels are not well behaved with respect to completions.  In
Proposition \ref{prop:72} it is clear that $\iota$ is injective on a
dense subalgebra but this does not imply that its extension to the
crossed product is injective.  We had to put in the extra effort to
show that it was isometric.  In Proposition \ref{prop:74} it is easy
to show that those elements in $\Gamma_c(G,r^*\mcal{A})$ for which
$\rho(f) = 0$ are contained in $\Ex(U)$.  However, $\ker \rho$ is {\em
  not} the completion of $\ker \rho \cap \Gamma_c(G,r^*\mcal{A})$.  The
solution in this case was to work with representations since they are
determined by their action $\Gamma_c(G,r^*\mcal{A})$. 
\end{remark}

Of course, if we restrict Theorem \ref{thm:invtideal} to the group
bundle case then there is even more structure to worry about.  In
particular, we want to see that restriction preserves the bundle
structure of the crossed product.  

\begin{prop}
\label{prop:88}
\index{group bundle}
Suppose $(A,S,\alpha)$ is a separable dynamical system, $S$ is a group
bundle, and that $U$ is
an open subset of $S\unit$.  Then $A\rtimes_\alpha S(U)$ and
$A(U)\rtimes_\alpha S|_U$ are isomorphic as $C_0(U)$-algebras.  
Similarly if $C$ is a closed subset of
$S\unit$ then $A\rtimes_\alpha S(C)$ and $A(C)\rtimes_\alpha S|_C$ are 
isomorphic as $C_0(C)$-algebras.  
\end{prop}

\begin{proof}
Let $U$ be open in $S\unit$ and let $C = S\unit\setminus U$.  
Recall that $A\rtimes S$ is a $C_0(S\unit)$-algebra.  It follows from
Proposition \ref{prop:71} that $A\rtimes S(U)$ is isomorphic to the
ideal 
\begin{align*}
I_U &= \cspn\{\phi\cdot f: \phi\in C_0(S\unit), f\in A\rtimes S,
\phi(C) = 0\} \\
&= \cspn\{\phi\cdot f : \phi\in
C_0(S\unit),f\in\Gamma_c(S,p^*\mcal{A}), \phi(C) = 0\}
\end{align*}
via the inclusion map $\iota_1:A\rtimes S(U)\rightarrow A\rtimes S$.  
We claim that $I_U = \Ex(U)$.  Recall that $\Ex(U)$ is the closure of
$\Gamma_c(S|_U,p^*\mcal{A})$ inside $A\rtimes S$.  Given $\phi\in
C_0(G\unit)$ and $f\in\Gamma_c(S,p^*\mcal{A})$ such that $\phi(C) = 0$ we
would like to show $\phi\cdot f \in \Ex(U)$.  Let $K =
\{v\in S\unit:|\phi(v)|\geq \epsilon\}$ and observe that $K$ is disjoint from
$C$. Thus we can find a function $\psi\in C_c(S\unit)^+$ such that
$\psi$ is one on $K$ and $\psi$ is zero off a neighborhood $V\subset
U$ of $K$.  Then $\supp
\psi\phi\cdot f \subset p\inv(V)\subset S|_U$ so that $\psi\phi\cdot
f\in \Ex(U)$.  Furthermore we constructed $\psi$ so that
$\|\psi\phi\cdot f - \phi\cdot f\|_\infty < \epsilon$.  Since we also
have $\supp\psi\phi\cdot f\subset \supp f$ we can use this
construction to find a sequence in $\Ex(U)$ which converges to
$\phi\cdot f$ in the inductive limit topology.  Hence $\phi\cdot f\in
\Ex(U)$ and it follows that $I_C\subset\Ex(U)$.  
Next suppose $f\in \Gamma_c(S|_U,p^*\mcal{A})$.  
Then $\supp f$ is a compact set which is
disjoint from $C$.  Let $\phi\in C_c(S\unit)^+$ be one on $K$ and zero
on $C$.  Then $f=\phi\cdot f \in I_C$ and it follows that
$\Ex(U)\subset I_C$.  Now, we also know that the inclusion map
$\iota_2:\Gamma_c(S|_U,p^*\mcal{A})\rightarrow \Gamma_c(S,p^*\mcal{A})$
extends to an isomorphism of $A(U)\rtimes S|_U$ with $\Ex(U)$.
Consider the isomorphism $\iota_2\inv\circ\iota_1$ from $A(U)\rtimes S|_U$
onto $A\rtimes S(U)$.  If $\phi\in C_0(U)$ and $f\in
\Gamma_c(S|_U,p^*\mcal{A})$ then $
\iota_2\inv\circ\iota_1(\phi\cdot f)(u) = \iota_2\inv(\phi\cdot
f)(u)$.\footnote{Recall that, because the quotient map is given by
  restriction on sections, $g(u)$ is the restriction of $g$ to $S_u$.}
Now $\iota_2\inv(\phi\cdot f)(u)$ is just the restriction of
$\phi\cdot f$ to $C_c(S_u,A(u))$ and therefore 
\[
\iota_2\inv\circ\iota_1(\phi\cdot f)(u) = \phi(u)f(u) = \phi\cdot
(\iota_2\inv\circ\iota_1(f))(u).
\]
Thus $\iota_2\inv\circ\iota_1$ is $C_0(U)$-linear. 

Next, it follows from Proposition \ref{prop:73} that
the restriction map factors to an isomorphism $\bar{\rho}_1$ of 
$A\rtimes S/\Ex(U)$ onto $A(C)\rtimes S|_C$.  We also know from
Proposition \ref{prop:71} that the
restriction map factors to an isomorphism $\bar{\rho}_2$ of $A\rtimes
S/I_C$ onto $A\rtimes S(C)$.  Since $I_C=\Ex(U)$ we may form the
isomorphism $\bar{\rho}_2 \circ \bar{\rho}_1\inv$ of $A(C)\rtimes
S|_C$ onto $A\rtimes S(C)$.   
Suppose $f\in \Gamma_c(S_C,p^*\mcal{A})$ and $\phi\in
C_0(C)$.  Choose $a\in A\rtimes S$ so that $\rho_1(a) = f$.  Applying
Proposition \ref{prop:71} to $C_0(S\unit)$ it follows that the
restriction map from $C_0(S\unit)$ to $C_0(C)$ is surjective.  In
particular we can extend $\phi$ to an element to of $C_0(S\unit)$.  
Clearly $\rho_1(\phi\cdot a) = \phi\cdot f$.
Thus
\[
\bar{\rho}_2(\bar{\rho}_1\inv(\phi\cdot f))(u) = \
\rho_2(\phi\cdot a)(u) = \phi(u)a(u)
= \phi(u) \rho_2(a)(u) = \phi\cdot
\bar{\rho}_2(\bar{\rho}_1\inv(f))(u)
\]
and therefore the isomorphism $\bar{\rho}_2\circ\bar{\rho}_1\inv$ is
$C_0(C)$-linear. 
\end{proof}

\begin{corr}
\label{cor:11}
Suppose $(A,S,\alpha)$ is a separable dynamical system, $S$ is a group
bundle, and that $U$ is
an open subset of $S\unit$.  Then $(A\rtimes_\alpha S(U))\sidehat \cong
(A(U)\rtimes_\alpha S|_U)\sidehat$.  Furthermore, if we view both of
these sets as the disjoint union $\coprod_{u\in U}(A(u)\rtimes
S_u)\sidehat$ then the identification is given by the identity.  
Corresponding statements hold if $C$ is a closed subset of $S\unit$.
\end{corr}

\begin{proof}
Suppose $U$ is an open subset of $S\unit$ and let $C = S\unit\setminus
U$.  Let $\sigma:(A\rtimes S)\sidehat\rightarrow S\unit$ be the map arising
from $A\rtimes S$ as a $C_0(S\unit)$-algebra.  
It follows from Corollary \ref{cor:9} that we can identify $(A\rtimes
S(U))\sidehat$ with the set $\sigma\inv(U)$ and that, if we view both
of these sets as $\coprod_{u\in U}(A(u)\rtimes S_u)\sidehat$ then the
identification is given by the identity.  Thus it suffices to show
that we can identify $\sigma\inv(U)$ and $A(U)\rtimes S|_U$ in the
appropriate manner.  Now, the inclusion map $\iota$ extends to an
isomorphism from $A(U)\rtimes S|_U$ to 
$\Ex(U)$ so that we can form the homeomorphism $\phi_1:\Ex(U)\sidehat\rightarrow
(A(U)\rtimes S_U)\sidehat$ given by 
$\phi_1(\pi) = \pi\circ\iota$.  Furthermore, we
can identify the set $P=\{\pi\in (A\rtimes
  S)\sidehat : \pi(\Ex(U))\ne 0\}$ with $\Ex(U)\sidehat$ 
via the map $\phi_2:P\rightarrow \Ex(U)\sidehat$ given by 
$\phi_2(\pi) = \pi|_{\Ex(U)}$.  
However $\Ex(U) = I_C$ and it was shown  
in the proof of Corollary \ref{cor:9} that in this case 
$P = \sigma\inv(U)$.  Thus
$\phi = \phi_1\circ \phi_2$ is a homeomorphism from $\sigma\inv(U)$
onto $A(U)\rtimes S|_U$.  Now suppose $\pi$ is a representation of
$A(u)\rtimes S_u$ and let $\pi'$ be its lift to $A\rtimes S$ and
$\pi''$ its lift to $A(U)\rtimes S|_U$.  
Furthermore, recall that the quotient map from $A\rtimes S$ to
$A(u)\rtimes S_u$ is given by restriction on $\Gamma_c(S,p^*\mcal{A})$.  Then 
for $f\in \Gamma_c(S_U,p^*\mcal{A})$, we have 
\[
\phi(\pi')(f) = \pi'(\iota(f)) = \pi(\iota(f)|_{S_u}) = \pi(f|_{S_u})
= \pi''(f)
\]
Thus $\phi(\pi')=\pi''$
and if we identify $\sigma\inv(U)$ and $A(U)\rtimes S|_U$ 
with the disjoint union
$\coprod_{u\in U} (A(u)\rtimes S_u)\sidehat$ then $\phi$ is given by
the identity map. 

Moving on, we can also use Corollary \ref{cor:9} to reduce to
the problem of showing that $A(C)\rtimes S|_C$ can be identified with
$\sigma\inv(C)$ in the appropriate fashion.  Recall that the
restriction map $\rho$ factors to an isomorphism of $A\rtimes
S/\Ex(U)$ with $A(C)\rtimes S|_C$ so that we can build a
homeomorphism $\phi_1:(A\rtimes S/\Ex(U))\sidehat \rightarrow
(A(C)\rtimes S|_C)\sidehat$ by letting $\phi_1(\pi) = \pi\circ
\overline{\rho}\inv$.  Then we identify $Q=\{\pi\in (A\rtimes
S)\sidehat: \Ex(U)\subset\ker\pi\}$ with $(A\rtimes S/\Ex(U))\sidehat$
via the map $\phi_2(\pi) = \bar{\pi}$ where $\bar{\pi}$ is the
factorization of $\pi$ to the quotient.  Well $\Ex(U)=I_C$ and we
showed in the proof of Corollary \ref{cor:9} that in this case
$Q=\sigma\inv(C)$.
Thus the desired homeomorphism is $\phi = \phi_1\circ\phi_2$. 
Fix $u\in C$ and $\pi\in (A(u)\rtimes S_u)\sidehat$. Let $\pi'$ be the
lift of $\pi$ to $A\rtimes S$ and $\pi''$ the lift to $A(C)\rtimes
S_C$. Then given $f\in \Gamma_c(S,p^*\mcal{A})$ we have 
\[
\phi(\pi')(\rho(f)) = \pi'(f) = \pi(f|_{S_u}) = \pi(\rho(f)|_{S_u}) = 
\pi''(\rho(f)).
\]
Since $\rho$ is surjective, this shows that
$\phi(\pi')=\pi''$.  Thus
we can view $\phi$ as the identity map. 
\end{proof}


\section{Unitary Actions}
\label{sec:unitary}

In this section we will discuss what it means for a groupoid to act
trivially.  The main goal will be to show that if the action is
trivial then the crossed product reduces to a tensor product.  As with
group crossed products, trivial actions are going to be defined by
unitaries.  

\begin{definition}
\label{def:51}
\index{dynamical system!unitary}
Suppose $S$ is a locally compact Hausdorff 
groupoid group bundle and $A$ is a
$C_0(S\unit)$-algebra.  Then a {\em unitary action} of $S$ on $A$ is
defined to be a collection $\{u_s\}_{s\in S}$ such that 
\begin{enumerate}
\item $u_s \in UM(A(p(s)))$ for all $s\in S$,
\item $u_{st} = u_su_t$ whenever $p(s)=p(t)$, and 
\item $s\cdot a := u_s a$ defines a (strongly) continuous action of $S$ on the
  associated upper-semi\-continu\-ous bundle $\mcal{A}$.
\end{enumerate}
The triple $(A,S,u)$ is called a unitary dynamical system.  
\end{definition}

\begin{remark}
We show in Section \ref{sec:locally-unitary} that if $u$ is a unitary action
of $S$ on $A$ then the restriction of $u$ to $S_v$ for $v\in S\unit$
gives a unitary action of $S_v$ on $A(v)$ in the sense of
\cite[Definition 2.70]{tfb2}.  Thus, Definition \ref{def:51} is really
just a ``bundled'' version of the notion of a unitary action of
a group on a $C^*$-algebra.  
\end{remark}

As with groupoid dynamical systems there is an ``unbundled''
definition.  However, we first take this opportunity to present some
of the basic facts about multipliers of $C_0(X)$-algebras. 

\begin{lemma}[{\cite[Lemma 2]{yinglee}}]
\label{lem:22}
Suppose that $A$ is a $C_0(X)$-algebra and that $m\in M(A)$.  Then for
each $x\in X$ there exists a multiplier $m(x)\in M(A(x))$ such that
$m(x)(a(x)) = m(a)(x)$.  Conversely,  if we are given $m_x\in M(A(x))$
for all $x\in X$ and if for each $a\in A$ there are elements $b,c\in
A$ such that for all $x\in X$
\begin{equation}
\label{eq:69}
b(x) = m_xa(x)\quad\text{and}\quad c(x)=m_x^*a(x), 
\end{equation}
then 
\begin{enumerate}
\item there is a $m\in M(A)$ such that $m(x) = m_x$ for all $x\in X$,
  and 
\item $\sup_{x\in X}\|m(x)\| = \|m\| <\infty$.
\end{enumerate}
\end{lemma}

\begin{remark}
\label{rem:18}
Condition \eqref{eq:69} is equivalent to requiring that for all $a\in
A$ there are elements $b,b'\in A$ such that $b(x)=m_xa(x)$ and
$b'(x)=a(x)m_x$ for all $x\in X$. 
\end{remark}

\begin{proof}
Since any ideal $I$ in $A$ is also an ideal in $M(A)$ any multiplier
$m$ defines a multiplier $m_I$ of $A/I$ via $m_I(a+I)=m(a)+I$.  The
first portion of the lemma follows immediately.  Now suppose we are
given $m_x$ as in the statement of the lemma and define a map
$m:A\rightarrow A$ by $m(a)(x) = m_xa(x)$.  The map $m$ is well
defined by assumption and it is easy to show that $m$ is
an adjointable $A$-linear operator on $A$ 
with adjoint $m^*(a) = m_x^*a(x)$.  Since we view
multipliers as the adjointable $A$-linear operators on $A_A$ we have
established part (a).  

Let $L = \sup_x\|m(x)\|$.  Then, viewing $A$ as
$\Gamma_0(X,r^*\mcal{A})$, we have 
\begin{align*}
\|m(a)\| &= \sup_{x} \|m(a)(x)\| \\
&\leq \sup_x \|m(x)\|\|a(x)\| \\
&\leq \sup_x \|m(x)\| \|a\|.
\end{align*}

Thus $\|m\|\leq L$.  Fix $\epsilon > 0$ and $x\in X$.  We can find
$b\in A(x)$ of norm one such that
$\|m(x)b\|\geq\|m(x)\|-\epsilon$.  Since the norm on $A(x)$ is the
quotient norm there is an $a\in A$ with $a(x)=b$ and $\|a\|\leq
1+\epsilon$.  But then $\|m(a)\|\geq \|m(a)(x)\| \geq
\|m(x)\|-\epsilon$ and it follows that 
\[
\|m\|\geq \frac{\|m(x)\|-\epsilon}{1+\epsilon}
\]
Since $\epsilon$ was arbitrary $\|m\|\geq \|m(x)\|$ for all $x\in X$
and $L\leq \|m\|$. 
\end{proof}

Now we can present an equivalent form of Definition \ref{def:51}.  

\begin{prop}
\label{prop:76}
Suppose $(A,S,u)$ is a unitary dynamical system.  Then there is an
element $u\in UM(p^*A)$ such that $u(s) = u_s$ for all $s\in S$.
Conversely, if we have $u\in UM(p^*A)$ then there are elements $u_s\in
UM(A(p(s)))$ for all $s\in S$ and if $u_{st} = u_s u_t$ whenever
$p(s)=p(t)$ then $\{u_s\}$ defines a unitary action of $S$ on $A$.
\end{prop}

\begin{proof}
Suppose $(A,G,u)$ is a unitary action and $f\in p^*A$.  We need
to show that 
\begin{align*}
h(s) &:= u_s f(s),&
g(s) &:= u_s^*f(s)
\end{align*}
define elements of $p^*A$.  The continuity of $h$ is obvious from
condition (c) of Definition~\ref{def:51}.  Suppose $s_i\rightarrow s$
and $a_i\rightarrow a$.  First, observe that condition (b) of
Definition~\ref{def:51} guarantees that $u_{s\inv} = u_s\inv = u_s^*$
for all $s\in S$.  Now, we know $s_i\inv\rightarrow s\inv$ and
therefore
\[
u_{s_i\inv}a_i = u_{s_i}^*a_i \rightarrow u_s^*a = u_{s\inv}a.
\]
It follows immediately that $g$ is continuous as well.  Furthermore, 
\[
\|h(s)\| = \|u_s f(s)\| = \|f(s)\| = \|u_s^* f(s)\| = \|g(s)\|
\]
so that both $h$ and $g$ must vanish at infinity because $f$ does.
Thus $h,g\in p^*A$.  Hence Lemma \ref{lem:22} implies that there is
a multiplier $u$ such that $u(f)(s) = u_sf(s)$ for all $s\in S$.
Since each $u_s$ is a unitary, it is clear that $u$ must be a
unitary.  

Next, suppose we are given $u\in UM(p^*A)$.  Then, via Lemma
\ref{lem:22}, we know there exists multipliers $u_s$ such that
$u_s(f(s)) = u(f)(s)$.  However, since $u$ is a unitary each $u_s$
must be as well.  Furthermore, condition (b) of Definition \ref{def:51}
holds by assumption.  All that is left is to show the action is
continuous.  Suppose $s_i\rightarrow s$ and $a_i\rightarrow a$ such
that $p(s_i) =
p(a_i)$ and $p(s) = p(a)$.  Choose $f\in p^*A$ such that $f(s) = a$.
Then $u(f)\in p^*A$ and $u(f)(s) = u_s(f(s)) = u_s a$.  Furthermore,
since $f(s_i)-a_i\rightarrow 0$, we have 
\[
\|u_{s_i}a_i - u(f)(s_i)\| = \|a_i-f(s_i)\| \rightarrow 0.
\]
It follows from Proposition \ref{prop:37} that the action is
continuous. 
\end{proof}

Next, given a unitary dynamical system we can form an associated
groupoid dynamical system as follows. 

\begin{remark}
\index[not]{$\Ad u$}
Suppose $A$ is a $C^*$-algebra and $u\in UM(A)$.  Then $u$ defines an
automorphism on $A$ via conjugation.  This automorphism is denoted $\Ad
u$ and is given by $\Ad u(a) = u a u^*$. 
\end{remark}

\begin{prop}
\index{dynamical system}
Suppose $(A,S,u)$ is a unitary dynamical system.  
Then the collection $\{\Ad
u_s\}_{s\in S}$ defines a groupoid action of $S$ on $A$. 
\end{prop}

\begin{proof}
Given a unitary action let $u$ be the corresponding element of
$UM(p^*A)$ guaranteed by Proposition \ref{prop:76}.  Then define
$\Ad u:p^*A\rightarrow p^*A$ by $\Ad u(f) = ufu^*$.  Clearly $\Ad u$ is a
$C_0(S\unit)$-linear automorphism of $p^*(A)$.  As in Proposition
\ref{prop:44} there exists isomorphisms $(\Ad u)_s :
A(p(s))\rightarrow A(p(s))$ for all $s\in S$.  Furthermore, these
isomorphisms are given by 
\[
(\Ad u)_s(f(s)) = \Ad_u(f)(s) = ufu^*(s) = u_s f(s)u^*_s.
\]
Thus $(\Ad u)_s = \Ad u_s$.  Finally if $p(s) = p(t)$ then 
\[
\Ad u_s\circ \Ad u_t (a) = u_su_t a u_t^* u_s^* = \Ad u_{st}(a).
\]
It now follows from Proposition \ref{prop:44} that $\Ad u$ is an
action of $S$ on $A$. 
\end{proof}

This allows us to define a special class of groupoid actions which we
will eventually see are the aforementioned 
``trivial'' dynamical systems. 

\begin{definition}
\label{def:52}
Suppose $S$ is a group bundle and $A$ is a $C_0(S\unit)$-algebra.
Then a dynamical system $(A,S,\alpha)$ is said to be {\em unitary} or
{\em unitarily implemented} if there exists a unitary action $u$ of
$S$ on $A$ such that $\alpha = \Ad u$.  
\end{definition}

At this point we need to make a brief detour through the notion of
equivalent actions.  The following construction will
play the role of isomorphism for dynamical systems. 

\begin{definition}
\label{def:53}
\index{exterior equivalent}
Suppose $G$ is a locally compact Hausdorff 
groupoid and $A$ is a $C_0(X)$-algebra.  Furthermore,
suppose $\alpha$ and $\beta$ are actions of $G$ on $A$.  Then we say
that $\alpha$ and $\beta$ are {\em exterior equivalent} if there is a
collection $\{u_\gamma\}_{\gamma\in G}$ such that 
\begin{enumerate}
\item $u_\gamma \in UM(A(r(\gamma)))$ for all $\gamma\in G$, 
\item $u_{\gamma\eta} = u_\gamma\overline{\alpha}_\gamma(u_\eta)$ for
  all $\gamma,\eta\in G$ such that $s(\gamma)=r(\eta)$,\footnote{Here
    $\overline{\alpha}_\gamma$ denotes the canonical extension of
    $\alpha_\gamma$ to the multiplier algebra.} 
\item the map $(\gamma,a)\mapsto u_\gamma a$ is jointly continuous on
  the set $\{(\gamma,a)\in G\times\mcal{A}:{r(\gamma) = q(a)}\}$, and
\item $\beta_\gamma = \Ad u_\gamma \circ \alpha_\gamma$ for all
  $\gamma\in G$. 
\end{enumerate}
\end{definition}

The following lemma expands on condition (b) above and gives us a
formula for inversion.  

\begin{lemma}
\label{lem:23}
Suppose that $(A,G,\alpha)$ and $(A,G,\beta)$ are exterior equivalent
dynamical systems and that the equivalence is implemented by
$\{u_\gamma\}$.  Then
\begin{enumerate}
\item $u_w = \id$ for all $w\in G\unit$, and 
\item $u_{\gamma\inv} = \overline{\alpha}_\gamma\inv(u_\gamma^*)$ for
  all $\gamma\in G$.
\end{enumerate}
\end{lemma}

\begin{proof}
If $w\in G\unit$ then $\alpha_w = \id$ so that 
\[
u_w = u_{w^2} = u_w \overline{\alpha}_w(u_w) = u_w^2.
\]
Since $u_w$ is a unitary this implies $u_w = \id$.  Now, given
$\gamma\in G$ we have 
\[
\id = u_{\gamma\inv\gamma} = u_{\gamma\inv}
\overline{\alpha}_{\gamma\inv}(u_\gamma).
\]
After recalling that $\alpha_{\gamma\inv} = \alpha_\gamma\inv$ and
that $\overline{\alpha}_\gamma\inv(u_\gamma)$ is a unitary we conclude
\[
\overline{\alpha}_\gamma\inv(u_\gamma^*) = 
\overline{\alpha}_\gamma\inv(u_\gamma)^* = u_{\gamma\inv}.\qedhere
\]
\end{proof}

As before, we present an alternate definition which removes the bundle theory. 

\begin{prop}
\label{prop:77}
Suppose $\alpha$ and $\beta$ are exterior equivalent actions of the
locally compact groupoid $G$ on the $C_0(G\unit)$-algebra $A$ with
the collection $\{u_\gamma\}$ implementing the equivalence.  Then
there is an element $u\in UM(r^*A)$ such that $u(f)(\gamma) = u_\gamma
f(\gamma)$ for all $f\in r^*A$ and $\gamma\in G$.  

Conversely, if $u\in UM(r^*A)$ then there are $u_\gamma\in
UM(A(r(\gamma)))$ for all $\gamma\in G$.  If $u_{\gamma\eta} =
u_\gamma \overline{\alpha}_\gamma(u_\eta)$ whenever $s(\gamma) =
r(\eta)$ and $\beta_\gamma = \Ad u_\gamma\circ\alpha_\gamma$ for all
$\gamma\in G$ then $\alpha$ and $\beta$ are exterior equivalent. 
\end{prop}

\begin{proof}
Suppose $\alpha$, $\beta$ and $\{u_\gamma\}$ are as in the first part of the
proposition and let $\mcal{A}$ be the upper-semicontinuous bundle
associated to $A$.  Given $f\in r^*A$ we must show that
\begin{align*}
g(\gamma) &:= u_\gamma f(\gamma), &
h(\gamma) &:= u_\gamma^* f(\gamma)
\end{align*}
define elements of $r^*A$.  It is clear from condition (c) of
Definition \ref{def:53} that $g$ defines a continuous section of
$r^*\mcal{A}$. Showing that $h$ is continuous takes a little more
work.  Suppose $\gamma_i\rightarrow \gamma$ in $G$ and $a_i\rightarrow
a$ in $\mcal{A}$ such that $r(\gamma_i) = q(a_i)$ for all $i$ and
$r(\gamma) = q(a)$.  It follows that $\gamma_i\inv\rightarrow
\gamma\inv$ and $a_i^*\rightarrow a^*$.  Furthermore, since $\alpha$
is a continuous action, we have $\alpha_{\gamma_i}\inv(a_i^*)
\rightarrow \alpha_\gamma\inv(a)$.  It follows from the continuity of
$\{u_\gamma\}$ that 
\begin{equation}
\label{eq:70}
u_{\gamma_i\inv} \alpha_{\gamma_i}\inv(a_i^*) \rightarrow
u_{\gamma\inv}\alpha_\gamma\inv(a^*).
\end{equation}
Applying Lemma \ref{lem:23} we conclude that 
\begin{equation}
\label{eq:71}
\alpha_{\gamma_i}\inv(u^*_{\gamma_i}a_i^*) \rightarrow
\alpha_\gamma\inv(u_\gamma^* a^*).
\end{equation}
If we apply the continuity of $\alpha$ with respect to
$\gamma_i\rightarrow \gamma$ and \eqref{eq:71} we obtain 
\begin{equation}
u_{\gamma_i}^* a_i^* \rightarrow u_\gamma^*a^*
\end{equation}
and therefore $a_iu_{\gamma_i}\rightarrow au_\gamma$.  It
follows immediately that $h$ is also a continuous section.
Furthermore, since $u_\gamma$ is unitary for all $\gamma\in G$, we have 
\[
\|g(\gamma)\| = \|u_\gamma f(\gamma)\| = \|f(\gamma)\| =
\|f(\gamma)u_\gamma\| = \|h(\gamma)\|.
\]
Because $f$ vanishes at infinity this is enough to show that $g,h\in
r^*A$.  Therefore we can conclude from Lemma \ref{lem:22} that there
exists $u\in UM(r^*A)$ which has the required form. 

Now suppose we are given $u\in UM(r^*A)$ with the properties listed in
the second half of the proposition.  It follows from Lemma
\ref{lem:22} that there are $u_\gamma\in UM(A(r(\gamma)))$ for all
$\gamma\in G$. Furthermore, by assumption, the only condition of
Definition \ref{def:53} which isn't satisfied is condition (c). 
Suppose $\gamma_i\rightarrow \gamma$ in $G$ and $a_i\rightarrow a$ in
$\mcal{A}$ such that $r(\gamma_i) = q(a_i)$ and $r(\gamma)
= q(a)$.  Choose $f\in r^*A$ such that $f(\gamma) = a$. 
Then $u(f)\in r^*A$ and $u(f)(\gamma) = u_\gamma(f(\gamma)) = u_\gamma
a$.  Furthermore,
since $f(\gamma_i)-a_i\rightarrow 0$, we have 
\[
\|u_{\gamma_i}a_i - u(f)(\gamma_i)\| = \|a_i-f(\gamma_i)\| \rightarrow 0.
\]
The continuity condition now follows from Proposition \ref{prop:37}.
\end{proof}

The most important fact about exterior equivalent
actions is the following 

\begin{prop}
\label{prop:84}
Suppose $(A,G,\alpha)$ and $(A,G,\beta)$ are exterior equivalent
separable groupoid dynamical systems with the equivalence implemented
by $\{u_\gamma\}$.  Then the map
$\phi:\Gamma_c(G,r^*\mcal{A})\rightarrow \Gamma_c(G,r^*\mcal{A})$
defined by 
\begin{equation}
\label{eq:72}
\phi(f)(\gamma) = f(\gamma)u_\gamma^*
\end{equation}
for all $\gamma\in G$ extends to an isomorphism from $A\rtimes_\alpha
G$ onto $A\rtimes_\beta G$.  
\end{prop}

\begin{proof}
Let $(A,G,\alpha)$ and $(A,G,\beta)$ be exterior equivalent dynamical
systems with the equivalence implemented by $\{u_\gamma\}$.  Use
Proposition \ref{prop:77} to find $u\in UM(r^*A)$ such that
$uf(\gamma) = u_\gamma f(\gamma)$ for all $\gamma$.  Given $f\in
\Gamma_c(G,r^*\mcal{A})$ view $f$ as an element of $r^*A$ and define
$\phi(f) = fu^*$.  It is clear from the construction of $u$ that
$\phi$ is also given by \eqref{eq:72}.  Furthermore $fu^*\in r^*A$ so
that $\phi(f)$ is a continuous section.  However, it follows
from \eqref{eq:72} that $\phi(f)$ is compactly supported as well.

Obviously, $\phi$ is linear.  We would like to show
that it is a $*$-homomorphism.  Given
$f,g\in\Gamma_c(G,r^*\mcal{A})$ we have 
\begin{align*}
\phi(f)*\phi(g)(\gamma) &= \int_G f(\eta)u_\eta^* \beta_\eta(g(\eta\inv\gamma)
u_{\eta\inv\gamma}^*) d\lambda^{r(\gamma)}(\eta) \\
&= \int_G f(\eta)u_\eta^* u_\eta
\alpha_\eta(g(\eta\inv\gamma))\overline{\alpha}_\eta(u_{\eta\inv}^*
\overline{\alpha}_\eta\inv(u_\gamma))^*u_\eta^*
d\lambda^{r(\gamma)}(\eta) \\
&= \int_G f(\eta)\alpha_\eta(g(\eta\inv\gamma))u_\gamma^* 
\overline{\alpha}_\eta(u_{\eta\inv}^*)u_\eta^*
d\lambda^{r(\gamma)}(\eta).
\end{align*}
Applying Lemma \ref{lem:23} to our calculation we obtain
\begin{align*}
\phi(f)*\phi(g)(\gamma) &= \int_G f(\eta)\alpha_\eta(g(\eta\inv
\gamma)) u_\gamma^* u_\eta u_\eta^* d\lambda^{r(\gamma)}(\eta) \\
&= f*g(\gamma)u_\gamma^* = \phi(f*g)(\gamma).
\end{align*}
We can also use Lemma \ref{lem:23} to show that 
\begin{align*}
\phi(f)^*(\gamma) &= \beta_\gamma(f(\gamma\inv)u_{\gamma\inv}^*)^* = 
u_\gamma \alpha_\gamma(u_{\gamma\inv} f(\gamma\inv)^*)u_\gamma^* \\
&= u_\gamma
\overline{\alpha}_\gamma(u_{\gamma\inv})\alpha_\gamma(f(\gamma\inv)^*)u_\gamma^*
\\
&= u_\gamma u_\gamma^* f^*(\gamma)u_\gamma^* = \phi(f^*)(\gamma).
\end{align*}
Thus $\phi$ is a $*$-homomorphism. 

Next, observe that given $f\in \Gamma_c(G,r^*\mcal{A})$ we have 
\[
\|\phi(f)(\gamma)\| = \|f(\gamma)u_\gamma^*\| = \|f(\gamma)\|.
\]
It follows quickly that $\phi$ is continuous with respect to the
inductive limit topology and therefore Corollary \ref{cor:24} implies
that $\phi$ extends to a $*$-homomorphism from $A\rtimes_\alpha G$
into $A\rtimes_\beta G$.  
We can define an inverse $\psi$ for $\phi$ on
$\Gamma_c(G,r^*\mcal{A})$ by $\psi(f)(\gamma) = f(\gamma)u_\gamma$.  An
argument nearly identical to the above shows that $\psi$ extends to a
$*$-homomorphism on $A\rtimes_\beta G$.  Since $\phi$ and $\psi$ are
inverses on a dense subset they are inverses on the entire algebra and
$\phi$ is an isomorphism. 
\end{proof}

Moving on, our statement that the unitary actions are
``trivial'' dynamical systems is supported by the next lemma.  However,
let us first introduce an action which is as trivial as possible.  

\begin{example}
\label{ex:25}
\index{dynamical system}
Suppose $S$ is a locally compact group bundle and $A$ is a
$C_0(S\unit)$-algebra.  Consider the identity map
$\id:p^*A\rightarrow p^*A$.  This isomorphism is clearly
$C_0(S\unit)$-linear.  Furthermore, $\id_s:A(p(s))\rightarrow A(p(s))$
is the identity map for all $s\in S$.  Therefore $\id_{st} =
\id_s\circ \id_t$ and Proposition \ref{prop:44} implies that 
the collection of identity maps
$\id_s:A(p(s))\rightarrow A(p(s))$ defines an action of $S$ on $A$. 
Observe that group bundles are the only groupoids which can act
trivially in this way.  If the source and range map are not equal then
$s^*A$ is not equal to $r^*A$ and we cannot use the identity map to induce
a groupoid action.  
\end{example}

\begin{lemma}
\label{lem:24}
If $(A,S,\alpha)$ is a unitary dynamical system then it is exterior
equivalent to the trivial system $(A,S,\id)$.  
\end{lemma}

\begin{proof}
Suppose $\alpha$ is implemented by the unitaries $\{u_s\}$.  We claim
that $\{u_s\}$ also implements an equivalence between $\id$ and
$\alpha$.  Condition (a) of Definition \ref{def:53} holds by
assumption, as does condition (c).  Next observe that if $p(s)=p(t)$ then
\[
u_{st} = u_s u_t = u_s \overline{\id}_s(u_t)
\]
so that condition (b) holds as well.  Finally we check that 
\[
\alpha_s(a) = \Ad u_s(a) = \Ad u_s\circ \id_s(a).\qedhere
\]
\end{proof}

\begin{remark}
The curious reader may wonder why we have only defined unitary actions
for a special class of groupoids.  Unitary actions should
always be equivalent to the trivial action and, as stated in Example
\ref{ex:25}, the trivial action only makes sense for group bundles.
Thus, it only makes sense to define unitary actions for group bundles.
\end{remark}

\subsection{Tensor Products}
\label{sec:tensor}
\index{tensor product}

We want to show that crossed products of unitary dynamical systems are
tensor products.  However, we are working with fibred
objects so we need to use a tensor product which respects the bundle
structure on the algebras.  It is assumed that the reader is familiar
with the basics of $C^*$-algebraic tensor products.  In particular we
will cite \cite[Appendix B]{tfb} frequently. In fact, just for
reference, we reproduce the following  

\begin{prop}[{\cite[Theorem B.27]{tfb}}]
\label{prop:85}
Suppose $A$ and $B$ are $C^*$-algebras.  Then there are nondegenerate
homomorphisms $\iota_A:A\rightarrow M(A\otimes_{\max} B)$ and
$\iota_B:B\rightarrow M(A\otimes_{\max} B)$ such that 
\begin{enumerate}
\item $\iota_A(a)\iota_B(b) = \iota_B(b)\iota_A(a) = a\otimes b$ for
  $a\in A$ and $b\in B$, 
\item if $\phi$ and $\psi$ are representations of $A$ and $B$ with
  commuting ranges then there is a representation
  $\phi\otimes_{\max}\psi$ of $A\otimes_{\max} B$ such that 
\[
\phi\otimes_{\max} \psi(\iota_A(a)\iota_B(b)) = \phi(a)\psi(b)
\]
for $a\in A$ and $b\in B$, 
\item $A\otimes_{\max} B = \cspn\{\iota_A(a)\iota_B(b):a\in A, b\in
  B\}$.
\end{enumerate}
If $D$ is a $C^*$-algebra and $j_A:A\rightarrow M(D)$ and
$j_B:B\rightarrow M(D)$ are homomorphisms satisfying the analogues of
(a),(b) and (c) then there is an isomorphism $\theta$ of $A\otimes_{\max}
B$ onto $D$ such that $\phi(a\otimes b) = j_A(a)j_B(b).$
\end{prop}

The proper notion of a ``fibred'' tensor product is the balanced
tensor product, defined below.  

\begin{definition}
\label{def:50}
\index{tensor product!balanced}
\index[not]{$A\otimes_{C_0(X)} B$}
Suppose $A$ and $B$ are $C_0(X)$-algebras and let $A\otimes_{\max} B$ denote
the (maximal) tensor product of $A$ and $B$.  The {\em balancing
  ideal} $I_X$ is the ideal in $A\otimes_{\max} B$ generated by 
\[
\{ f\cdot a \otimes b - a\otimes f\cdot b : f\in C_0(X), a\in A, b\in
B\}.
\]
The {\em balanced tensor product} $A\otimes_{C_0(X)} B$ is defined to
be the quotient $A\otimes_{\max} B/I_X$.  
\end{definition}

\begin{remark}
If we view two $C^*$-algebras $A$ and $B$ as being
$C_0(\{\text{pt}\})$-algebras then the balanced tensor product is just the
usual maximal tensor product.
\end{remark}

\begin{remark}
Our tensor products will generally be maximal tensor products
\cite[Appendix B]{tfb}.    However, most of the time we will be
working with nuclear $C^*$-algebras so that we will not have
to make this distinction.  
\end{remark}

Moving on, one of the key facts about the balanced tensor product is
that, at least for nice $C^*$-algebras, its spectrum is the fibre
product of the spectra of its components.  This
proposition is a reproduction of \cite[Lemma 1.1]{pullback}.

\begin{prop}
\label{prop:83}
Suppose $A$ and $B$ are separable $C_0(X)$-algebras and that either $A$ or $B$
is nuclear.  Define the bundle product of $\widehat{A}$ with
$\widehat{B}$ to be
\[
\widehat{A}\times_X\widehat{B} := \{(\pi,\rho)\in
\widehat{A}\times\widehat{B}: \sigma_A(\pi) = \sigma_B(\rho)\}.
\]
\begin{enumerate}
\item The map $(\pi,\rho)\mapsto \pi\otimes_\sigma \rho$ induces a
  homeomorphism $\Phi$ of $\widehat{A}\times_X\widehat{B}$ onto its
  range in $(A\otimes_{C_0(X)}B)\sidehat$.  
\item If either $A$ or $B$ is GCR then this homeomorphism is
  surjective.  
\end{enumerate}
\end{prop}

\begin{proof}
Since at least one of $A$ or $B$ is nuclear there is a unique tensor
product $A\otimes B$.  Furthermore, we cite \cite[Theorem B.45]{tfb}
to see that the map $(\pi,\rho)\mapsto \pi\otimes_\sigma \rho$ induces a
homeomorphism $\Phi$ of $\widehat{A}\times \widehat{B}$ onto its range
in $(A\otimes B)\sidehat$ and is surjective if either $A$ or $B$ is
GCR.  Since $A\otimes_{C_0(X)} B$ is a quotient
of $A\otimes B$ by the balancing ideal $I$, we can identify
$(A\otimes_{C_0(X)} B)\sidehat$ with the closed set 
$\{R\in (A\otimes B)\sidehat : I\subset \ker R\}$.  If
$\pi\in \widehat{A}$ and $\rho\in\widehat{B}$ then $\pi\otimes_\sigma \rho(I)
= 0$ if and only if $\pi(\phi\cdot a)\rho(b) = \pi(a)\rho(\phi\cdot
b)$ for all $a\in A$, $b\in B$ and $\phi\in C_0(X)$. Let $x = \sigma_A(\pi)$
so that $\pi$ factors to a representation $\overline{\pi}$ 
of $A(x)$ and $y= \sigma_B(\rho)$
so that $\rho$ factors to a representation $\overline{\rho}$ 
of $B(y)$.  Then
$\pi(\phi\cdot a)\rho(b) = \pi(a)\rho(\phi\cdot b)$ if and only if 
\begin{equation}
\label{eq:97}
\phi(x)\overline{\pi}(a(x))\overline{\rho}(b(y)) = 
\phi(y)\overline{\pi}(a(x))\overline{\rho}(b(y))
\end{equation}
for all $\phi\in C_0(X)$, $a\in A$, and $b\in B$.  However
\eqref{eq:97} holds if and only if $x = y$.  Thus the restriction
of $\Phi$ to the closed set $\widehat{A}\times_X \widehat{B}$  maps
onto $(A\otimes_{C_0(X)} B)\sidehat$.
\end{proof}

Recall that in Definition \ref{def:32} we defined the pull back of a
$C^*$-algebra to be the section algebra of the pull back of the
associated bundle.  This is not how the pull back is classically
defined. However, we now have the tools to prove 
the following proposition, which brings us back to the usual
definition.

\begin{prop}[{\cite[Proposition 1.3]{pullback}}]
\index{pull back}
\label{prop:75}
Suppose $X$ and $Y$ are locally compact Hausdorff spaces, $A$ is a
$C_0(X)$-algebra and $\tau:Y\rightarrow X$ is a continuous
surjection.  Then the pull back algebra $\tau^*A$ is isomorphic to the
balanced tensor product $C_0(Y)\otimes_{C_0(X)} A$. 
\end{prop}

\begin{proof}
First, recall that if $\tau:Y\rightarrow X$ is a continuous surjection then
we can view $C_0(Y)$ as a $C_0(X)$-algebra as in Example \ref{ex:15}.  
Let $\iota:C_0(Y)\odot A\rightarrow \Gamma_0(Y,\tau^*\mcal{A})$ be
such that $\iota(f\otimes a)(y) = f(y)a(\tau(y))$.\footnote{In other
  words, view $f\otimes a$ as an elementary tensor in $\tau^*A$.}  It is
clear that this defines a continuous section and $\iota(f\otimes a)$
vanishes at infinity because $f$ does and $a$ is bounded.  
It is straightforward to
show that $\iota$ is a $*$-homomorphism and it follows quickly from
Proposition \ref{prop:46} that $\iota$ maps onto a dense subset.
Finally, because $\iota$ is a homomorphism, pulling back the uniform
norm on $\Gamma_0(Y,\tau^*\mcal{A})$ defines a $C^*$-seminorm on
$C_0(Y)\odot A$ by $\|f\otimes a\|_\iota = \|\iota(f\otimes
a)\|_\infty$.  However this implies that $\|\iota(f\otimes a)\|_\infty
\leq \|f\otimes a\|_{\max}$ and that $\iota$ extends to a representation
of $C_0(Y)\otimes A$.  Furthermore we clearly have 
\[
\iota(\phi\cdot f\otimes a)(y) = \phi(\tau(y))f(y)a(\tau(y)) = 
\iota(f\otimes \phi\cdot a)(y)
\]
for all $y\in Y$.  Hence $\iota$ vanishes on the balancing ideal and
factors to a homomorphism on the balanced tensor product
$C_0(Y)\otimes_{C_0(X)} A$, which we also denote by $\iota$.  

Now suppose $R$ is an irreducible representation of $C_0(Y)\otimes_{C_0(X)}
A$.  Since $C_0(Y)$ is abelian it is both nuclear and GCR so that by
Proposition \ref{prop:83} there exists $y\in Y$ and
$\pi\in\widehat{A}$ such that $R = \ev_y\otimes\sigma \pi$ where $\ev_y$ is
the evaluation representation.  Furthermore we
must have $\sigma(\pi) = \tau(y)$ so that $\pi$ factors to a representation
$\overline{\pi}$ of $A(\tau(y))$.  If $\pi$ acts on $\mcal{H}$ then 
$\ev_y\otimes_\sigma \pi$ acts on $\C\otimes \mcal{H}$, which we
identify with $\mcal{H}$, via $\ev_y \otimes_\sigma \pi(f\otimes a) =
f(y)\pi(a)$.  We may now compute for $\sum_if_i\otimes a_i \in C_0(Y)\odot A$
\begin{align*}
\left\|\ev_y\otimes_\sigma \pi\left(\sum_i f_i\otimes a_i\right)\right\| &= 
\left\|\sum_{i}f_i(y)\pi(a_i)\right\| = \left\|\pi\left(\sum_if_i(y)a_i
\right)\right\| \\
&=\left\|\overline{\pi}\left(\sum_i f_i(y)a_i(\tau(y))\right)\right\|
= \left\|\overline{\pi}\left(\sum_i f_i\otimes a_i(y)\right)\right\|\\
&\leq \left\|\sum_i f_i\otimes a_i(y)\right\| \leq
\left\|\iota\left(\sum_i f_i\otimes a_i\right)\right\|_\infty.
\end{align*}
Since this is true for every irreducible representation of
$C_0(Y)\otimes_{C_0(X)} A$ we conclude $\left\|\sum_i f_i\otimes a_i\right\|\leq
\left\|\iota\left(\sum_i f_i\otimes a_i\right)\right\|_\infty$.  It follows that $\iota$ is isometric
on a dense subset and therefore must be isometric everywhere.  Hence
$\iota$ is an isomorphism of $C_0(Y)\otimes_{C_0(X)} A$ with
$\tau^*A$. 
\end{proof}

The last thing we want to show about balanced tensor products is that the
resulting algebra is still a $C_0(X)$-algebra.  This is, perhaps,
unsurprising considering how balanced tensor products work for
modules.  Of course, it follows from Proposition \ref{prop:83} that,
at least for nice algebras, the spectrum of the balanced tensor
product is the bundle product of the spectra of its components.  The
$C_0(X)$-algebra structure then follows from Theorem
\ref{thm:c0xalgs}.  However, this construction can be done in greater
generality.  

\begin{prop}
Suppose $A$ and $B$ are $C_0(X)$-algebras.  Then $A\otimes_{C_0(X)} B$
is a $C_0(X)$-algebra with the action characterized by 
\[
\phi\cdot (a\otimes b) := (\phi\cdot a)\otimes b = a\otimes(\phi\cdot
b)
\]
\end{prop}

\begin{proof}
Suppose $\Phi_A$ and $\Phi_B$ implement the $C_0(X)$-algebra structure
on $A$ and $B$ respectively.  Then $\Phi_A$ and $\Phi_B$ are
nondegenerate homomorphisms into $ZM(A)$ and $ZM(B)$ respectively.
Let $\iota_A$ and $\iota_B$ be the nondegenerate homomorphisms from
Proposition \ref{prop:85} and let $\pi:A\otimes B\rightarrow
A\otimes_{C_0(X)} B$ be the quotient map.  Consider the map
\[
\Psi = \overline{\pi}\circ \overline{\iota}_A\circ \Phi_A.
\]
This is certainly a homomorphism and one can check, given $\phi\in
C_0(X)$, $a\in A$ and $b\in B$, that 
\begin{align}
\label{eq:104}
\Psi(\phi)(a\otimes b) &=
\pi(\overline{\iota}_A(\Phi_A(\phi))\iota_A(a)\iota_B(b)) 
= \pi(\iota_A(\phi\cdot a)\iota_B(b))\\  &= (\phi\cdot a)\otimes b =
a\otimes (\phi\cdot b). \nonumber
\end{align}
Thus $\Psi$ has the desired form.  Now we need to see that it maps
into the center of the multiplier algebra.  Using Lemma \ref{lem:18}
and linearity it will suffice to show that 
\[
\Psi(\phi)((a\otimes b)(c\otimes d)) = (a\otimes
b)(\Psi(\phi)(c\otimes d))
\]
for all $a,c\in A$ and $b,d\in B$.  However, by \eqref{eq:104} we
have 
\[
\Psi(\phi)(ac\otimes bd) = (\phi\cdot ac)\otimes bd = 
(a\otimes b)((\phi\cdot c)\otimes d) = (a\otimes
b)(\Psi(\phi)(c\otimes d)).
\]
Thus $\Psi$ maps into the center of the multiplier algebra.  The last
thing we need to show is that $\Psi$ is nondegenerate.  It will
suffice to show that elements of the form $\phi\cdot a \otimes b$ span
a dense subset.  Since elements of the form $\phi\cdot a$
span a dense subset in $A$, we are done.  
\end{proof}

\begin{remark}
We won't use this fact directly, and so do not prove it here, but it
is clear enough that $(A\otimes_{C_0(X)} B)(x) = A(x)\otimes B(x)$ for
all $x\in X$.
\end{remark}

Let us get back to the matter at hand.  
We now are able to prove the main theorem concerning unitary actions,
which is that they have trivial crossed products.  

\begin{theorem}
\index{New Result}
\label{thm:unitary}
\index{dynamical system!unitary}
\index{tensor product!balanced}
Suppose $(A,S,\alpha)$ is a separable unitary dynamical system with
$\alpha$ implemented by $u$.  Then
there is a $C_0(S\unit)$-linear isomorphism $\phi:{C^*(S)\otimes_{C_0(S\unit)}
A}\rightarrow {A\rtimes_\alpha S}$ which is characterized for $a\in
A$ and $f\in C_c(S)$ by 
\begin{equation}
\phi(f\otimes a)(s) = f(s)a(p(s))u_s^*
\end{equation}
\end{theorem}

\begin{remark}
Since $\phi$ is $C_0(S\unit)$-linear it factors, by Proposition
\ref{prop:43}, to isomorphisms $\phi_u:C^*(S_u)\otimes A(u) \rightarrow
A(u)\rtimes S_u$.  It is not difficult to check that these are the
usual isomorphisms that arise from unitary actions 
\cite[Lemma 2.73]{tfb2}.
\end{remark}

\begin{proof}
First, let $\beta$ be a Haar system for $S$ and 
consider the trivial action $\id$ of $S$ on $A$.  Given $a\in A$
and $f\in C_c(S)$ define $\iota(f\otimes a)(s) := f(s)a(p(s))$.  In
other words, view $f\otimes a$ as an ``elementary tensor'' in
$\Gamma_c(S,p^*\mcal{A})$ in the sense of Proposition \ref{prop:45}.
It follows that $\iota(f\otimes a)\in\Gamma_c(S,p^*\mcal{A})$.  Extend $\iota$
to the algebraic tensor product $C_c(S)\odot A$ by linearity so that
$\iota:C_c(S)\odot A\rightarrow \Gamma_c(G,r^*\mcal{A})$. Observe that
$\ran \iota$ is
dense with respect to the inductive limit topology.    We would
like to show that $\iota$ is a $*$-homomorphism.  It is enough to do
the calculations on elementary tensors.  For $f,g\in C_c(S)$ and
$a,b\in A$ we have 
\begin{align*}
\iota(f\otimes a)*\iota(g\otimes b)(s) &= \int_S f(t)a(p(t))g(t\inv
s)b(p(t))d\beta^{p(s)}(t) \\
&= \int_S f(t)g(t\inv s)d\beta^{p(s)}(t) ab(p(s)) \\
&= \iota(f*g\otimes ab)(s).
\end{align*}
and
\[
\iota(f\otimes a)^*(s) = (f(s\inv)a(p(s)))^* =
\overline{f(s\inv)}a^*(p(s)) = \iota(f^*\otimes a^*)(s).
\]
Thus $\iota$ is a $*$-homomorphism. 

Now we check that $\iota$ is bounded.  Suppose 
$(S\unit*\mfrk{H},\mu,\pi,U)$ is a covariant
representation of $(A,S,\id)$.  Then $U$ is a groupoid representation
of $S$ and we can form the integrated representation as in Proposition
\ref{prop:66} which we also denote by $U$.  Let $\pi =
\int_{S\unit}^\oplus \pi_u d\mu(u)$ be a decomposition of $\pi$.
Since $(\pi,U)$ is covariant we must have, for all $a\in A$ and almost
every $s\in S$,
\begin{equation}
\label{eq:99}
\pi_{p(s)}(a(p(s)))U_s = U_s\pi_{p(s)}(a(p(s))).
\end{equation}
However, we can now compute for $f\in C_c(G)$ and 
$h\in \mcal{L}^2(S\unit*\mfrk{H},\mu)$ that 
\begin{align*}
(\pi(a)U(f))h(u) &= \pi_u(a(u))U(f)h(u) \\
&= \int_S \pi_u(a(u))f(s)U(s)h(u)\Delta(s)\neghalf d\beta^u(s) \\
&= \int_S f(s)U(s)\pi(a(u))h(u)\Delta(s)\neghalf d\beta^u(s) \\
&= (U(f)\pi(a))h(u).
\end{align*}
We can extend this by continuity to all $f\in C^*(S)$ and conclude
that $\pi$ and $U$ are commuting representations of $A$ and $C^*(S)$.
It follows from Proposition \ref{prop:83} that there exists a
representation $U\otimes\pi$ on $C^*(S)\otimes A$ such that
\[
U\otimes \pi(f\otimes a) = U(f)\pi(a).
\]
Given $f\in C_c(G)$ and $a\in A$ we check that 
\begin{align}
\label{eq:103}
\pi\rtimes U(\iota(f\otimes a))h(u) &= \int_S
\pi_u(f(s)a(u))U_sh(u)\Delta(s)\neghalf d\beta^u(s) \\ \nonumber
&= \pi_u(a(u))\int_S f(s)U_sh(u)\Delta(s)\neghalf d\beta^u(s) \\ \nonumber
&= \pi(a)U(f)h(u) = U\otimes \pi(f\otimes a)h(u).
\end{align}
Using linearity, we conclude that $\pi\rtimes U(\iota(\xi)) = U\otimes \pi(\xi)$
for all $\xi \in C_c(S)\odot A$.  Thus, given $\xi\in C_c(S)\odot A$, 
\[
\|\pi\rtimes U(\iota(\xi))\| = \|U\otimes \pi(\xi)\|
\leq \|\xi\|.
\]
Since this is true for all covariant representations $(\pi,U)$, it
follows that $\iota$ is bounded and extends to a homomorphism on
$C^*(S)\otimes A$. Furthermore, since the range of $\iota$ is dense, it
must be surjective. What's more, given
$\phi\in C_0(S\unit)$, $f\in C_c(G)$ and $a\in A$ we have 
\[
\iota(\phi\cdot f\otimes a)(s) = \phi(p(s))f(s)a(p(s)) = 
\iota(f\otimes\phi\cdot a)(s).
\]
It follows by continuity and linearity that $\iota$ factors through
the balancing ideal and induces a homomorphism
$\hat{\iota}:C^*(S)\otimes_{C_0(X)} A\rightarrow A\rtimes S$.  

We would like to show that $\hat{\iota}$ is isometric.  Suppose $R$ is
a faithful representation of $C^*(S)\otimes_{C_0(X)} A$ and let
$\overline{R}$ be its lift to $C^*(S)\otimes A$.  It follows
\cite[Corollary B.22]{tfb} that there are commuting representations $\pi$
and $U$ of $A$ and $C^*(S)$ such that $\overline{R}=U\otimes \pi$.
Furthermore, since $U\otimes \pi$ contains the balancing ideal, a
quick computation shows that $U(\phi\cdot f)\pi(a) =
U(f)\pi(\phi\cdot a)$ for all $\phi \in C_0(S\unit)$, $f\in C^*(S)$,
and $a\in A$.  Now, without loss of generality, we can use Theorem
\ref{thm:scalardis} to assume that $U$ is the integrated form of some
groupoid representation $(S\unit*\mfrk{H},\mu,U)$.  Furthermore we
have for all $\phi\in C_0(S\unit)$, $a\in A$, and $f\in C^*(s)$
\begin{align*}
\pi(\phi\cdot a)U(f)h(u) &= U(\phi\cdot f)\pi(a)h(u) \\
&= \int_S \phi(u)f(s)U_s \pi(a)h(u)\Delta(s)\neghalf d\beta^u(s) \\
&= \phi(u)U(f)\pi(a)h(u)\\
&= T_\phi \pi(a)U(f)h(u)
\end{align*}
where $T_\phi$ is the diagonal operator associated to $\phi$. Since $U$
is nondegenerate this implies that $\pi$ is $C_0(X)$-linear.  Let $\pi =
\int_{S\unit}^\oplus \pi_ud\mu(u)$ be the decomposition of $\pi$ and
let $\nu$ be the measure on $S$ induced by $\mu$.  All we need to do 
to prove that $(\pi,U)$ is a covariant representation of $(A,S,\id)$
is to verify the covariance relation.  In other words we need to show
that \eqref{eq:99} holds $\nu$-almost everywhere.  Let $\{a_i\}$ be a
countable dense subset in $A$ and $e_l$ a special orthogonal
fundamental sequence for $S\unit* \mfrk{H}$.  Since the ranges of
$\pi$ and $U$ commute, we have for all $i,l,k$ and $f\in C_c(G)$
\begin{align*}
0 =& (\pi(a_i)U(f)e_l,e_k) - (U(f)\pi(a_i)e_l,e_k) \\
=& \int_S (f(s)\pi_{p(s)}(a_i(p(s)))U_s
e_l(p(s)),e_k(p(s)))\Delta(s)\neghalf d\nu(s) \\
&- 
 \int_S
 (f(s)U_s\pi_{p(s)}(a_i(p(s)))e_l(p(s)),e_k(p(s)))\Delta(s)\neghalf
 d\nu(s) \\
=& \int_S f(s)((\pi_{p(s)}(a_i(p(s)))U_s -
U_s\pi_{p(s)}(a_i(p(s))))e_l(p(s)),e_k(p(s))) d\nu(s).
\end{align*}
This holds for all $f\in C_c(G)$ so that we may conclude for each $i,l$
and $k$ there exists a $\nu$-null set $N_{i,l,k}$ such that 
\begin{equation}
\label{eq:102}
((\pi_{p(s)}(a_i(p(s)))U_s -
U_s\pi_{p(s)}(a_i(p(s))))e_l(p(s)),e_k(p(s))) = 0
\end{equation}
for all $s\not\in N_{i,l,k}$.  However if we let $N = \bigcup_{i,l,k}
N_{i,l,k}$ then $N$ is still a $\nu$-null set and for each $s\not\in
N$ \eqref{eq:102} holds for all $i,l$ and $k$.  Since $\{e_l(p(s))\}$
is a basis (plus zero vectors) for each $p(s)$ this implies that for
$s\not\in N$ we have 
\[
\pi_{p(s)}(a_i(p(s)))U_s = U_s \pi_{p(s)}(a_i(p(s)))
\]
for all $i$.  It follows from the fact that $\{a_i\}$ is dense in $A$
that this holds for all $a\in A$.  Thus $(\pi, U)$ is  covariant
representation of $(A,G,\id)$.  Furthermore, we can reuse the
computation in \eqref{eq:103} to show
that in this case $\pi\rtimes U\circ \iota =
\pi\otimes U$.  Given $\xi\in C^*(S)\otimes A$ let $\xi'$ be its
image in $C^*(S)\otimes_{C_0(X)} A$.  We then have
\[
\|\xi'\| = \|R(\xi')\| = \|U\otimes \pi(\xi)\| = 
\|\pi\rtimes U(\iota(\xi))\| \leq \|\iota(\xi)\| = \|\hat{\iota}(\xi')\|.  
\]
It follows that $\hat{\iota}$ is isometric and is therefore an
isomorphism.  

To finish the proof, observe that because of Proposition \ref{prop:84} and
Lemma \ref{lem:24}, the map $\psi:A\rtimes_{\id} S\rightarrow
A\rtimes_\alpha S$ given by $\psi(f)(s) = f(s)u_s^*$ is an
isomorphism.  Thus $\phi = \psi\circ \hat{\iota}$ is an isomorphism from
$C^*(G)\otimes_{C_0(X)} A$ onto $A\rtimes_\alpha S$ and given $f\in
C_c(S)$ and $a\in A$ we have 
\[
\phi(f\otimes a)(s) = \iota(f\otimes a)(s)u_s^* =
f(s)a(p(s))u_s^*.
\]
It follows quickly that $\phi$ is a $C_0(S\unit)$-linear isomorphism
and we are done.
\end{proof}


\section{Locally Unitary Actions}
\label{sec:locally-unitary}
Now that we have developed the theory of unitary actions
we can modify Definition \ref{def:51} and introduce some new concepts.
The basic idea is that we weaken the continuity condition and see what
kind of structure we have left.  
This material is a generalization of \cite{locunitary}.   

\begin{definition}
\label{def:57}
\index{dynamical system!pointwise unitary}
Suppose $S$ is a group bundle and $A$ is a $C_0(S\unit)$-algebra.
Then a dynamical system $(A,G,\alpha)$ is said to be {\em pointwise
  unitary} if $\alpha|_{S_u}$ is unitarily implemented for each $u\in
S\unit$.  
\end{definition}

Notice that in a pointwise unitary dynamical system $\alpha_s$ is
still given by conjugation by a unitary for all $s\in S$.  
What we have done is restrict the continuity of the unitaries to just
the fibres.  Of course, this should bring us back to the usual notion
of a unitary group action, which we show in the following
proposition.  

\begin{remark}
Given a $C^*$-algebra $A$ a function $f:X\rightarrow M(A)$ is said to
be strictly continuous if $x\mapsto f(x)a$ is continuous for all $a\in
A$.
\end{remark}

\begin{prop}
\label{prop:86}
Suppose $G$ is a locally compact group and $A$ is a $C^*$-algebra.
Then a map $u:G\rightarrow UM(A)$ is strictly continuous if and only
if the function $(s,a)\mapsto u_s a$ is continuous on $G\times A$.  In
particular, a unitary action of $G$ on $A$ is given by a strictly
continuous homomorphism $u:G\rightarrow UM(A)$. 
\end{prop}

\begin{proof}
The reverse direction is clear since strict continuity is weaker than
joint continuity.  Now suppose $u$ is strictly continuous,
$s_i\rightarrow s$, and $a_i\rightarrow a$.  Then 
\[
\|u_{s_i}a_i - u_s a\| \leq \|u_{s_i} (a_i - a)\| + \|u_{s_i}a-u_s a\|
\leq \|a_i-a\| + \|u_{s_i} a - u_s a\| \rightarrow 0
\]
It follows that $u_{s_i}a_i\rightarrow u_s a$.  
The rest of the proposition follows.  
\end{proof}

The problem with pointwise unitary actions is that the unitaries tell
you very little about the total space of the bundle structure.  It
will turn out to be more interesting if we pick a point ``between''
unitary and pointwise unitary. 

\begin{definition}
\label{def:56}
\index{group bundle}
\index{dynamical system!locally unitary}
\index{locally unitary|see{dynamical system}}
Suppose $S$ is a group bundle and $A$ is a $C_0(S\unit)$-algebra. A
dynamical system $(A,S,\alpha)$ is said to be {\em locally unitary} if
there is an open cover $\{U_i\}_{i\in I}$ of $S\unit$ such that
$(A(U_i),S|_{U_i},\alpha|_{S_{U_i}})$ is unitarily implemented for all
$i\in I$.  
\end{definition}

Note that if $S$ is a group bundle then every set in $S\unit$ is
$S$-invariant so the above definition makes sense.  
Our goal will be to analyze the exterior equivalence classes of
abelian locally unitary
actions on $C^*$-algebras with Hausdorff spectrum.  In particular, the
rest of the $C^*$-algebras in this section will have Hausdorff
spectrum and we will view them as $C_0(\widehat{A})$-algebras as in
Example \ref{ex:21}. This will allow us to identify the 
spectrum of the crossed product for unitary
actions. 

\begin{prop}
\label{prop:87}
Suppose $S$ is an abelian, second countable, locally compact Hausdorff
continuously varying group bundle, that $A$ is a $C^*$-algebra with
Hausdorff spectrum $S\unit$ and that $(A,S,\alpha)$ 
is a unitary dynamical system.  Let $\{u_s\}$
be the unitaries implementing $\alpha$ and for all $v\in S\unit$ let
$\pi_v$ be the unique (up to equivalence) irreducible representation
of $A(v)$. Define, for $\omega\in \widehat{S}$, 
\[
\omega\overline{\pi}_{\hat{p}(\omega)}(u)(s) :=
\omega(s)\overline{\pi}_{\hat{p}(\omega)}(u_s).
\]
Then the map $\phi :\widehat{S} \rightarrow
(A\rtimes_\alpha S)\sidehat$ given by $\phi(\omega) = \pi_{\hat{p}(\omega)}\rtimes
\omega \overline{\pi}_{\hat{p}(\omega)}(u)$ is a bundle homeomorphism.
\end{prop}

\begin{remark}
Since we can't use $u$ to denote both unitaries and units we will
temporarily use $x$ to denote elements of $S\unit$.  
\end{remark}

\begin{proof}
Let $(A,S,\alpha)$ and $u$ be as above. 
It follows from Theorem \ref{thm:unitary} that the map
$\psi:C^*(S)\otimes_{C_0(S\unit)} A\rightarrow A\rtimes S$ characterized
by $\psi(a\otimes f)(s) = f(s)a u_s^*$ is an isomorphism.  
Hence, there is a homeomorphism $\phi_1:(C^*(S)\otimes
A)\sidehat\rightarrow (A\rtimes_{C_0(S\unit)} S)\sidehat$ such that
$\phi_1(R) = R\circ\psi\inv$.  Next, recall that we identify the dual
group bundle $\widehat{S}$ with $C^*(S)\sidehat$.
Since $C^*(S)$ is an abelian $C^*$-algebra, and is therefore GCR and
nuclear, it follows from Proposition \ref{prop:83} that
$\phi_2:\widehat{S}\times_{S\unit}\widehat{A}\rightarrow(C^*(S)\otimes_{C_0(S\unit)}
A)\sidehat$ given by $\phi_2(\omega,\pi) = \omega\otimes_\sigma \pi$ 
is a homeomorphism.  Recall that if $\pi$ is a representation on $\mcal{H}$
then $\omega\otimes_\sigma \pi$ is a representation on $\C\otimes \mcal{H}$,
which we will of course identify with $\mcal{H}$, characterized by 
$\omega\otimes_\sigma \pi(f\otimes a) = \omega(f)\pi(a)$.  Moving on, since
$\widehat{A}= S\unit$ we can define another homeomorphism
$\phi_3:\widehat{S}\rightarrow \widehat{S}\times_{S\unit} \widehat{A}$
by $\phi_3(\omega)= (\omega,\pi_{\hat{p}(\omega)})$.  
Let $\phi = \phi_1\circ\phi_2\circ\phi_3$ and observe that
$\phi:\widehat{S}\rightarrow (A\rtimes S)\sidehat$ is a
homeomorphism.  Furthermore given $\omega\in \widehat{S}$ we have
$\phi(\omega) = \omega\otimes_\sigma \pi_{\hat{p}(\omega)}\circ\psi\inv$.  

Now, fix $x\in S\unit$ and $\omega \in \widehat{S}_x$ and define the map
$U:S_x\rightarrow U(\mcal{H})$ by $U_s = 
\omega(s)\overline{\pi}_x(u_s)$.  Since $u$ is a continuous action,
and since $\omega$ is continuous, it follows quickly that $U$ is a
unitary representation of $S_x$.  Furthermore we can compute for $a\in
A(x)$ and $s\in S_x$ that 
\[
U_s\pi_x(a) = \omega(s)\pi_x(u_sa) = \omega(s)\pi_x(u_s a
u_s^*u_s) = \pi_x(\alpha_s(a))U_s.
\]
Thus $(\pi_x,U)$ is a covariant representation of
$(A(x),S_x,\alpha)$.\footnote{We described what covariant
representations looked like for group dynamical systems in Remark
\ref{rem:20}.}  As such we can form the integrated representation
$\pi_x\rtimes U$.   Recall that $A\rtimes S$ is a
$C_0(S\unit)$-algebra and that the restriction map $\rho$ 
factors to an isomorphism 
between $A\rtimes S(x)$ and $A(x)\rtimes S_x$.   Using the
restriction map to view  $\pi_x\rtimes U$ as a representation of
$A\rtimes S$ we claim that $\pi_x\rtimes U = \phi(\omega)$.
It will suffice to show that given an elementary tensor $f\otimes a$
then $\pi_x\rtimes U(\psi(f\otimes a)) = \omega\otimes_\sigma \pi_x(f\otimes a)$.  
We compute, observing that the modular function is one since
$S$ is abelian, 
\begin{align*}
\pi_x\rtimes U(\psi(f\otimes a))h
&= \int_S \pi_x(f(s)a(x)u_s^*)\omega(s) \overline{\pi}_x(u_s)h d\beta^x(s)
\\
&= \int_S f(s)\omega(s)d\beta^x(s) \pi_x(a(x))h \\
&= \omega(f)\pi_x(a)h = (\omega\otimes_\sigma \pi_x)(f\otimes a)h.
\end{align*}
Thus $\phi(\omega) = \pi_x\rtimes U$ and since $U$ is just an abbreviated
notation for $\omega\overline{\pi}_x(u)$ we are almost done.  All that
is left is to show that $\phi$ preserves the fibres, but if
$\omega$ is a representation of $S_x$ then $\phi(\omega) =
\pi_x\rtimes U$ is clearly a representation of the fibre $A(x)\rtimes
S_x$.  
\end{proof}

Let us see what this implies in the weaker pointwise unitary
case.  

\begin{corr}
\label{cor:10}
\index{dynamical system!pointwise unitary}
Suppose $S$ is an abelian, second countable, locally compact Hausdorff
continuously varying group bundle, that $A$ is a $C^*$-algebra with
Hausdorff spectrum $S\unit$ and that $(A,S,\alpha)$ 
is a pointwise unitary dynamical system.  Then we can view
$(A\rtimes_\alpha S)\sidehat$ as a topological bundle over $S\unit$ and
for all $x\in S\unit$ the fibre over $x$ is isomorphic to $\widehat{S}_x$
\end{corr}

\begin{proof}
Recall $A\rtimes S$ is a $C_0(S\unit)$-algebra with fibres
$A(x)\rtimes S_x$. By definition, 
$\alpha$ is unitarily implemented on $S_x$ for all $x$.  Thus
we can apply Proposition \ref{prop:87} to each fibre and
conclude that $(A(x)\rtimes S_x)\sidehat$ is homeomorphic to
$\widehat{S}_x$.  If we view $(A\rtimes
S)\sidehat$ as a topological bundle with fibres $(A(x)\rtimes
S_x)\sidehat$ then the result follows.
\end{proof}

So if $\alpha$ is pointwise unitary then 
$(A\rtimes S)\sidehat$ is fibrewise isomorphic to the dual bundle
$\widehat{S}$.  A good question is to ask when is $(A\rtimes S)\sidehat$
a principal $\widehat{S}$-bundle in the sense of Section
\ref{sec:principal}.  Of course, the answer is hidden in the
title of this section.

\begin{theorem}
\label{thm:locunit}
\index{New Result}
\index{dynamical system!locally unitary}
\index{principal S-bundle@principal $S$-bundle}
Suppose $S$ is an abelian, second countable, locally compact Hausdorff,
continuously varying group bundle, that $A$ is a $C^*$-algebra with
Hausdorff spectrum $S\unit$ and that $(A,S,\alpha)$ 
is a locally unitary dynamical system.  Let $u^i$ implement $\alpha$
on $S|_{U_i}$ where $\{U_i\}$ is an open cover of $S\unit$ and let
$q:(A\rtimes_\alpha S)\sidehat\rightarrow S\unit$ be the bundle map.  
Then for each $i$ the map $\psi_i:\hat{p}_i\inv(U_i)\rightarrow
q\inv(U_i)$ such that 
\begin{equation}
\label{eq:105}
\psi_i(\omega) = \pi_{\hat{p}(\omega)}\rtimes
\omega\overline{\pi}_{\hat{p}(\omega)}(u^i)
\end{equation}
is a homeomorphism and the map $\gamma_{ij}$ such that 
\begin{equation}
\gamma_{ij}(p(s))(s) =
\overline{\pi}_{p(s)}((u_s^i)^*u_s^j)
\end{equation}
defines a continuous section of $\widehat{S}$.  Furthermore these maps
make $(A\rtimes S)\sidehat$ into a principal $\widehat{S}$-bundle with 
trivialization $(\mcal{U},\psi\inv,\gamma)$.  
\end{theorem}

\begin{proof}
Let $(A,G,\alpha)$ be as in the statement of the theorem.  Let
$\{u^i\}$ implement $\alpha$ on $S|_{U_i}$ where $U_i$ is an element of some
open cover $\mcal{U}$.  Given an open set $U\in\mcal{U}$ we identify
each of $(A(U)\rtimes S|_{U})\sidehat$, 
$(A\rtimes S(U))\sidehat$ and $q\inv(U)$ with the
disjoint union $\coprod_{x\in U}(A(x)\rtimes S_x)\sidehat$.
Furthermore Corollary \ref{cor:9} and Corollary
\ref{cor:11} imply that this identification respects the topologies on
all three spaces. 
In a similar fashion we identify each of $(C^*(S)(U))\sidehat$,
$C^*(S|_U)\sidehat$ and $\hat{p}\inv(U)$ with the disjoint union
$\coprod_{x\in U}\widehat{S}_x$ and again observe that this
identification preserves the topologies.  

Now, fix $U_i\in\mcal{U}$.  By assumption 
$\alpha|_{S|_{U_i}}$, denoted $\alpha$ whenever possible, is
unitarily implemented by $\{u^i\}$ and as such Proposition
\ref{prop:87} implies that the map $\psi_i:(S|_{U_i})\sidehat
\rightarrow  (A(U_i)\rtimes S|_{U_i})\sidehat$ defined via
\eqref{eq:105} is a homeomorphism.  However, under the
identifications made in the previous paragraph, we can view $\psi_i$ as
a map from $\hat{p}\inv(U_i)$ onto $q\inv(U_i)$.  Furthermore, $q\circ
\psi_i= \hat{p}$ since $\psi_i$ is a bundle isomorphism.  We
define the trivializing maps on $(A\rtimes S)\sidehat$ to be $\phi_i =
\psi_i\inv$.  What's more, since $(A\rtimes S)\sidehat$ is locally
homeomorphic to a locally compact Hausdorff space, we can conclude that
$(A\rtimes S)\sidehat$ is locally compact Hausdorff.  

Next, suppose $U_i,U_j\in\mcal{U}$ and for each $x\in
U_{ij}$ let $\pi_x$ be the (unique) irreducible representation of
$A(x)$.  On $A(x)\rtimes S_x$ both $u^i$ and $u^j$ implement $\alpha$
so that we compute, for $s\in S_x$ and $a\in A(x)$,
\[
\overline{\pi}_x((u_s^i)^*u_s^j)\pi_x(a) = 
\pi_x((u_s^i)^* u_s^j a) = \pi_x(\alpha_s\inv(\alpha_s(a))(u_s^i)^*u_s^j)
= \pi_x(a)\overline{\pi}_x((u_s^i)^*u_s^j).
\]
Since $\pi_x$ is irreducible it follows \cite[Lemma A.1]{tfb} that
$\gamma_{ij}(x)(s):= \overline{\pi}_x((u_s^i)^*u_s^j)$ is a scalar.
Since $u_s^i$ and $u_s^j$ are unitaries $\gamma_{ij}(x)(s)$ must be a
unitary as well and therefore has modulus one.  Next, observe that 
\begin{align*}
\gamma_{ij}(x)(st) &= \overline{\pi}_x((u_{st}^i)^*u_{st}^j) = 
\overline{\pi}_x((u_t^i)^*)\overline{\pi}_x((u_s^i)^*u_s^j)\overline{\pi}_x(u_t^j)
\\ &= \overline{\pi}_x((u_t^i)^*)\gamma_{ij}(x)(s)\overline{\pi}_x(u_t^j)
= \gamma_{ij}(x)(s)\gamma_{ij}(x)(t)
\end{align*}
Thus $\gamma_{ij}(x)$ is a homomorphism on $S_x$.  Finally if
$s_l\rightarrow s$ then, given $a\in A(x)$ and $h\in H$, we use the
continuity of $u^i$ and $u^j$ to conclude
\[
\gamma_{ij}(x)(s_l)\pi_x(a)h = \pi_x((u_{s_l}^i)^*u_{s_l}^j a)h 
\rightarrow \pi_x((u_s^i)^*u_s^j a)h = \gamma_{ij}(x)(s)\pi_x(a)h.
\]
Of course, this implies that $\gamma_{ij}(x)$ is continuous so that
$\gamma_{ij}(x)$ is a character on $S_x$ and the map $\gamma_{ij}$ is
a section of $\widehat{S}$ on $U_{ij}$.  Next, we compute for
$\omega\in \hat{p}\inv(U_{ij})$
\begin{equation}
\label{eq:106}
\phi_i\circ\phi_j\inv(\omega) = \psi_i\inv\circ\psi_j(\omega) = 
\psi_i(\pi_{\hat{p}(\omega)}\rtimes
\omega\overline{\pi}_{\hat{p}(\omega)}(u^j)).
\end{equation}
Given $s\in S_{\hat{p}(\omega)}$ we have 
\[
\overline{\pi}_{\hat{p}(\omega)}(u_s^j)  = 
\overline{\pi}_{\hat{p}(\omega)}(u_s^i)\overline{\pi}_{\hat{p}(\omega)}((u_s^i)^*u_s^j)
=
\gamma_{ij}(\hat{p}(\omega))(s)\overline{\pi}_{\hat{p}(\omega)}(u_s^i).
\]
Applying this to \eqref{eq:106} we obtain
\begin{equation}
\label{eq:107}
\phi_i\circ\phi_j\inv(\omega) = 
\psi_i\inv(\pi_{\hat{p}(\omega)}\rtimes
(\omega\gamma_{ij}(\hat{p}(\omega))\overline{\pi}_{\hat{p}(\omega)}(u^i)))
= \omega\gamma_{ij}(\hat{p}(\omega))
\end{equation}
Thus \eqref{eq:107} shows that the $\gamma_{ij}$ are transition
functions for the $\phi_i$.  Furthermore, suppose
$u_l\rightarrow u$.  Then
\[
\gamma_{ij}(u_l) = 
\phi_i\circ\phi_j\inv(u_l) \rightarrow 
\phi_i\circ\phi_j(u) = \gamma_{ij}(u).
\]
This suffices to show that $\gamma_{ij}$ is continuous.  It now follows
that the trivialization $(\mcal{U},\phi,\gamma)$ makes $(A\rtimes
S)\sidehat$ into a principal $\widehat{S}$-bundle.
\end{proof}

Of course this is little more than a curiosity unless we can use the
principal bundle structure to tell us something about the action
$\alpha$.  Fortunately, we can do just that.  

\begin{theorem}
\label{thm:unique}
\index{New Result}
\index{dynamical system!locally unitary}
\index{principal S-bundle@principal $S$-bundle}
Suppose $S$ is an abelian, second countable, 
locally compact Hausdorff, continuously
varying group bundle and that $A$ has Hausdorff spectrum $S\unit$.
Two locally unitary actions $(A,S,\alpha)$ and $(A,S,\beta)$ are
exterior equivalent if and only if
$(A\rtimes_\alpha S)\sidehat$ and $(A\rtimes_\beta S)\sidehat$ are
isomorphic as $\widehat{S}$-bundles.  
\end{theorem}

\begin{proof}
Suppose $\alpha$ and $\beta$ are equivalent locally unitary actions
and the equivalence is implemented by the collection $\{u_s\}$.  It
follows from Proposition \ref{prop:84} that the map
$\phi:A\rtimes_\alpha S\rightarrow A\rtimes_\beta S$ defined for $f\in
\Gamma_c(S,p^*\mcal{A})$ by $\phi(f)(s) = f(s)u_s^*$ is an
isomorphism.  As such it induces a homeomorphism $\Phi:(A\rtimes_\beta
S)\sidehat\rightarrow (A\rtimes_\alpha S)\sidehat$ via the map
$\Phi(\pi) = \pi\circ\phi$.  Furthermore, $\phi$ is
$C_0(S\unit)$-linear so that $\phi$ factors to an isomorphism on each
of the fibres.  This implies that, if $I^\alpha_x$ is the ideal such
that $(A\rtimes_\alpha S)/I^\alpha_x = A\rtimes_\alpha S(x)$ and
$I^\beta_x$ is the corresponding ideal for $\beta$, then we must have
$\phi(I^\alpha_x) = I^\beta_x$.  Thus if $I^\beta_x\subset
\ker \pi$ then $I^\alpha_x \subset \ker \pi\circ\phi$, and if $q_\alpha$ is
the bundle map on $(A\rtimes_\alpha S)\sidehat$ and $q_\beta$ is the
corresponding map for $\beta$ then we must have, by definition, 
$q_\beta(\pi) = q_\alpha(\Phi(\pi))$.  Therefore $\Phi$ is a bundle
isomorphism.  

Next, let us establish some notation.  Since $\alpha$ and $\beta$ are
both locally trivial we may as well pass to some common refinement and
assume that there exists an open cover $\mcal{U}$ of $S\unit$ such
that on $S|_{U_i}$ the unitary actions $v^i$ and $w^i$ implement
$\alpha$ and $\beta$, respectively.  Let $\phi_i$ and $\psi_i$ be the
trivializing maps induced by $v^i$ and $w^i$, respectively.
Furthermore given $x\in S\unit$ let $\pi_x$ be the (unique) irreducible
representation of $A(x)$ associated to $x$.  Now fix
$U_i\in \mcal{U}$ and $x\in U_i$.  In order to conserve notation we will drop
the $i$'s on the $v^i$ and $w^i$ unless they are needed. 
Recall that $\beta_s = \Ad u_s\circ \alpha_s$ so that we can compute
for $s\in S_x$
\begin{align*}
u_s^* w_s v_s^* a &= u_s^*\beta_s(\alpha_s\inv(a))w_sv_s^* \\
&= \Ad(u_s^*)\circ \beta_s \circ\alpha_s\inv(a)u_s^*w_sv_s^* \\
&= \Ad(u_s^*)\circ\Ad(u_s) \circ \alpha_s \circ\alpha_s\inv(a) u_s^*w_s v_s^*
\\
&= a u_s^* w_s v_s^*.
\end{align*}
It follows that $\beta_i(x)(s) :=
\overline{\pi}_{x}(u_s^*w_sv_s^*)$ commutes with
$\pi_{x}(A(x))$. Since $\pi_{x}$ is irreducible this implies
that $\beta_i(x)(s)$ must be a scalar.  Since $u_s,w_s$ and $v_s$
are all unitaries $\beta_i(x)(s)$ must have modulus one.
Furthermore, it is straightforward to show that the continuity
conditions on $u, v$ and $w$ all conspire to make $\beta_i(x)$
continuous on $S_{x}$.  Lastly we compute
\begin{align*}
\beta_i(x)(st) &= \overline{\pi}_x(u_{st}^*w_{st}v_{st}^*) = 
\overline{\pi}_x(\overline{\alpha}_s(u_t^*)u_s^*w_s w_t v_t^* v_s^*) \\
&= \overline{\pi}_x(v_s u_t^* v_s^* u_s^*w_s v_s^* v_s w_t v_t^* v_s^*)
\\
&= \overline{\pi}_x(v_s u_t^* v_s^*)\beta_i(x)(s)\overline{\pi}_x(v_s
w_t v_t^* v_s^*) \\
&= \beta_i(x)(s) \overline{\pi}_x(v_s u_t^* w_t v_t^* v_s^*) \\
&= \beta_i(x)(s) \overline{\pi}_x(v_s)
\beta_i(x)(t)\overline{\pi}_x(v_s^*) \\
&= \beta_i(x)(s)\beta_i(x)(t).
\end{align*}
Since $\beta_i(x)$ is a continuous $\T$-valued homomorphism, it is an
element of $\widehat{S}_x$ and thus $\beta_i$ is a section of $\widehat{S}$
on $U_i$. Now suppose $\omega \in \widehat{S}_x$. Then we compute for
$ f\in \Gamma_c(S,p^*\mcal{A})$
\begin{align}
\label{eq:109}
\pi_x\rtimes (\omega\overline{\pi}_x(w))(\phi(f)) &= 
\int_S \pi_x(\phi(f)(s)) \omega(s) \overline{\pi}_x(w_s)d\beta^x(s)
\\ \nonumber
&= \int_S \pi_x(f(s))\omega(s)\overline{\pi}_x(u_s^*w_s)d\beta^x(x) \\ \nonumber
&= \int_S \pi_x(f(s))\omega(s)\beta_i(x)(s)
\overline{\pi}_x(v_s)d\beta^x(s) \\\nonumber
&= \pi_x\rtimes (\omega\beta_i(x) \overline{\pi}_x(v))(f).
\end{align}
It follows from \eqref{eq:109} that 
\begin{align*}
\phi_i\circ \Phi \circ \psi_i\inv(\omega) &= 
\phi_i(\pi_x\rtimes (\omega\overline{\pi}_x(w))\circ \phi) \\
&= \phi_i(\pi_x \rtimes (\omega\beta_i(x) \overline{\pi}_x(v))) \\
&= \omega \beta_i(x).
\end{align*}
Therefore $\beta_i$ implements $\Phi$ on trivializations.
Furthermore since $\phi_i$, $\Phi$, and $\psi_i$ are all continuous, it
is now straightforward to show that $\beta_i$ is a continuous
section.  Thus $(\mcal{U},\Phi,\beta)$ is an $\widehat{S}$-bundle
isomorphism of $(A\rtimes_\beta S)\sidehat$ onto $(A\rtimes_\alpha
S)\sidehat$.  

Now suppose that $(\mcal{U},\Phi,\beta)$ is an $\widehat{S}$-bundle
isomorphism of $(A\rtimes_\alpha S)\sidehat$ onto $(A\rtimes_\beta
S)\sidehat$.  Let $w^i$ and $v^i$ implement $\alpha$ and $\beta$,
respectively. Notice that
$\mcal{U}$ must be a common refinement of 
the local trivializing cover for $\alpha$ and $\beta$ so that we may
as well assume $w^i$ and $v^i$ are defined on $\mcal{U}$.
Fix $U_i\in \mcal{U}$ and $x\in U_i$.  
For each $s\in S_x$ we define a unitary $u_s\in
UM(A(x))$ by 
\begin{equation}
\label{eq:110}
u_s := \beta_i(x)(s)w_s^i(v_s^i)^*.
\end{equation}
We need to show that \eqref{eq:110} doesn't depend on the choice of
$U_i$. So suppose $x\in U_j$ as well.  (Notice we are going to have to
keep track of the $i$ and $j$ for a while.)  Let $\gamma_{ij}$ and
$\eta_{ij}$ be the transition maps for $(A\rtimes_{\alpha} S)\sidehat$
and $(A\rtimes_\beta S)\sidehat$ respectively.  Recall that we can 
view $x$ as an element of $\widehat{S}_x$. 
Using the general theory of principal bundles, we obtain
\begin{align*}
\beta_i(x)\gamma_{ij}(x) &=
\psi_i(\Phi(\phi_i\inv(\phi_i(\phi_j\inv(x))))) \\
&= \psi_i(\psi_j\inv(\psi_j(\Phi(\phi_j\inv(x))))) \\
&= \eta_{ij}(x)\beta_j(x).
\end{align*}
We use this fact to compute
\begin{align}
\label{eq:111}
\beta_i(x)(s)\overline{\pi}_x(w_s^i(v_s^i)^*) &= 
\beta_i(x)(s) \overline{\pi}_x(w_s^i(v_s^i)^*v_s^j(v_s^j)^*) \\
\nonumber &=
\beta_i(x)(s)\overline{\pi}_x(w_s^i)\overline{\pi}_x((v_s^i)^*v_s^j)\overline{\pi}_x((v_s^j)^*)
\\ \nonumber
&= \beta_i(x)(s)\gamma_{ij}(x)(s)\overline{\pi}_x(w_s^i(v_s^j)^*) \\\nonumber
&= \beta_j(x)(s)\eta_{ij}(x)(s)\overline{\pi}_x(w_s^i(v_s^j)^*) \\\nonumber
&= \beta_j(x)(s)\overline{\pi}_x(w_s^i((w_s^i)^*w_s^j)(v_s^j)^*) \\\nonumber
&= \beta_j(x)(s)\overline{\pi}_x(w_s^j(v_s^j)^*).
\end{align}

\begin{remark}
Since $A$ has Hausdorff spectrum, each fibre $A(x)$ is simple
\cite[Lemma 5.1]{tfb}.  Thus, the only proper closed ideal is trivial, and
$\pi_x$ must be a faithful representation.  It is then straightforward
to calculate that the extension $\overline{\pi}_x$ is also faithful. 
\end{remark}

Hence \eqref{eq:111} implies that 
\[
\beta_i(x)(s)w_s^i(v_s^i)^* = \beta_j(x)(s)w_s^j(v_s^j)^*
\]
and that $u_s$ is well defined.  We now show that the $u_s$ implement
an equivalence between $\alpha$ and $\beta$.  Observe that the first
condition of Definition \ref{def:53} is satisfied by construction.  
Dropping the $i$'s again for convenience, we compute for $s,t\in S_x$
\begin{align*}
\overline{\pi}(u_{st}) &=
\beta_i(x)(st)\overline{\pi}_x(w_{st}v_{st}^*) \\
&=
\beta_i(x)(s)\beta_i(x)(t)\overline{\pi}_x(w_sv_s^*)\overline{\pi}_x(v_sw_t
v_t^* v_s^*) \\
&= \overline{\pi}_x(u_s v_s u_t v_s^*) \\
&= \overline{\pi}_x(u_s \overline{\alpha}_s(u_t)).
\end{align*}
Again using the fact that $\pi$ is faithful, this is sufficient to
verify the second condition of equivalence.  The continuity condition
is straightforward to prove using the fact that the actions $u$ and
$v$ are continuous, as well as the fact that $\beta_i$ is a continuous
section.  For the last condition observe that for $a\in A(x)$
\begin{align*}
\pi(\Ad u_s(\alpha_s(a))) &= \pi(u_s v_s a v_s^* u_s^*) \\
&= \beta_i(x)(s)\overline{\beta_i(x)(s)} \pi(w_s v_s^* v_s a v_s^* v_s
w_s^*) \\
&= \pi(w_s a w_s^*) = \pi(\beta_s(a)).
\end{align*}
Therefore $\Ad u_s \circ \alpha_s = \beta_s$ and $\{u_s\}$ implements
an equivalence between $\alpha$ and $\beta$. 
\end{proof}

Of course, this leads to the following corollary.  

\begin{corr}
\index{cohomology}
A locally unitary action of a continuously varying 
abelian group bundle $S$ on a
$C^*$-algebra with Hausdorff spectrum $S\unit$ is determined, up to
exterior equivalence, by the associated cohomological invariant of
$(A\rtimes S)\sidehat$ as a principal $\widehat{S}$-bundle.
Furthermore, the isomorphism class of $A\rtimes S$ is characterized by
this invariant. 
\end{corr}

\begin{proof}
This corollary just puts together Theorem \ref{thm:unique} with
Theorem \ref{prop:principcohom} and Proposition \ref{prop:84}.
\end{proof}

\begin{remark}
Since non-exterior equivalent actions can still have isomorphic
crossed products, it is possible that $A\rtimes_\alpha S$ can be
isomorphic to $A\rtimes_\beta S$ even though their invariants are
different.  
\end{remark}

The last piece of the puzzle will be to prove that locally unitary
actions are about as abundant as they can be. In other words we will
show that every
principal bundle can be obtained through a locally unitary action.  

\begin{theorem}
\label{thm:exist}
\index{New Result}
\index{principal S-bundle@principal $S$-bundle}
\index{dynamical system!locally unitary}
\index{transformation groupoid}
Suppose $S$ is an abelian, second countable, locally compact Hausdorff,
continuously varying group bundle and $q:X\rightarrow S\unit$ is a 
principal $S$-bundle.
Then the transformation groupoid $C^*$-algebra $C^*(S,X)$ has
Hausdorff spectrum $S\unit$ and the dual action of $\widehat{S}$ on
$C^*(S,X)$ defined for $\omega\in \widehat{S}_u$ 
on $C_c(S_u\times q\inv(u))$ by 
\begin{equation}
\label{eq:112}
\widehat{\lt}_\omega(f)(s,x) = \omega(s)f(s,x)
\end{equation}
is locally unitary.  Furthermore, $(C^*(S,X)\rtimes
\widehat{S})\sidehat$ and $X$ are isomorphic $S$-bundles. 
\end{theorem}

Of course, before we can get down to the details of Theorem
\ref{thm:exist} we have a number of things to check.  

\begin{lemma}
\label{lem:25}
Suppose $S$ is a continuously varying abelian group bundle which acts
continuously on a locally compact Hausdorff space $X$.  Let $q$ be the
range map on $X$ and define $X_u := q\inv(u)$ for all $u\in S\unit$.
Then $C^*(S,X)$ is a $C_0(S\unit)$-algebra and for each $u\in S\unit$
restriction factors to an isomorphism from $C^*(S,X)(u)$ onto
$C^*(S_u,X_u)$.  
\end{lemma}

\begin{remark}
Of course, now that we have a space called $X$ floating around we can
no longer denote elements of $S\unit$ by $x$.  Unfortunately, $u$ is
still going to conflict with our unitary notation, but we will make
the best of it.  To make matters worse, we are going to need to use
$\omega$ to denote characters, which will look ugly when paired with
the $w$'s that turn up. It is moments like this which make the author
wish the alphabet were longer. \cite{seuss} 
\end{remark}

\begin{proof}
Let $S\ltimes X$ be the transformation group\-oid associated to $S$
and $X$ and recall that by
definition $C^*(S,X)$ is equal to $C^*(S\ltimes X)$.  Furthermore, recall that
Proposition \ref{prop:68} states that the map $\Phi:C^*(S,X)\rightarrow
C_0(X)\rtimes_{\lt} S$ defined on $C_c(S\ltimes X)$ by $\Phi(f)(s)(x)
= f(s,x)$ is an isomorphism.  Let $\mcal{C}$ be the upper
semicontinuous bundle associated to $C_0(X)$ and recall from Example
\ref{ex:15} that $\mcal{C}$ has fibres $C_0(X)(u) = C_0(X_u)$.  
Proposition \ref{prop:65} says that
$C_0(X)\rtimes S$ is a $C_0(S\unit)$-algebra with the action defined
for $\phi\in C_0(S\unit)$ and $f\in \Gamma_c(S,p^*\mcal{C})$ by 
\[
\phi\cdot f(s)(x)= \phi(p(s))f(s)(x).
\]
Therefore we can use the isomorphism $\Phi$ to give $C^*(S,X)$ a
$C_0(S\unit)$-algebra structure defined via 
\[
\phi\cdot f(s,x) = \Phi\inv(\phi\cdot \Phi(f))(s,x) = 
\phi(p(s))\Phi(f)(s)(x) = \phi(p(s))f(s,x).
\]
Now fix $u\in S\unit$.  By construction $\Phi$ is $C_0(X)$-linear so
that $\Phi$ factors to an isomorphism
$\bar{\Phi}:C^*(S,X)(u)\rightarrow C_0(X)\rtimes S(u)$.   
The restriction map
$\rho$ factors to an isomorphism $\bar{\rho}$ of $C_0(X)\rtimes S(u)$ with
$C_0(X_u)\rtimes S_u$.  Furthermore since the action of $S_u$ on
$C_0(X_u)$ is still given by left translation, there is another
isomorphism $\Psi:C_0(X_u)\rtimes S_u\rightarrow C^*(S_u,X_u)$
defined in the same manner as $\Phi$.   Thus, we get an isomorphism
$\Psi\circ\bar{\rho}\circ \bar{\Phi}$ of $C^*(S,X)(u)$ onto
$C^*(S_u,X_u)$.  We would like to see that this isomorphism is given
by restriction for $f\in C_c(S\ltimes X)$.  Let $\sigma$ be the restriction
map from $C_c(S\ltimes X)$ onto $C_c(S_u\times X_u)$ and (foolishly)
let $q$ denote both the quotient map from $C^*(S,X)$ onto $C^*(S,X)(u)$
and the quotient map from $C_0(X)\rtimes S$ onto $C_0(X)\rtimes S(u)$.  We
then have for $f\in C_c(S\ltimes X)$
\begin{align*}
\Psi\circ \bar{\rho}\circ \bar{\Phi}(q(f))(s,x) &= 
\bar{\rho}(\bar{\Phi}(q(f)))(s)(x) = 
\bar{\rho}(q(\Phi(f)))(s)(x) \\
&= \rho(\Phi(f))(s)(x) = \Phi(f)(s)(x) \\
&= f(s,x) = \sigma(f)(s,x).
\end{align*}
Therefore, the isomorphism from $C^*(S,X)(u)$ onto $C^*(S_u,X_u)$ is
just the factorization of the restriction map on $C_c(S\ltimes X)$.  
\end{proof} 

\begin{remark}
Propositions like the one above are really just notational trickery.
Philosophically, $C^*(S,X)$ and $C_0(X)\rtimes S$ are (basically)
completions of the same function algebra and should be treated as the
same object.  
\end{remark}

Now, in order for Theorem \ref{thm:exist} to work we need to know that
$C^*(S,X)$ has Hausdorff spectrum.  

\begin{prop}
\label{prop:90}
Suppose $S$ is a continuously varying abelian group bundle and that
$q:X\rightarrow S\unit$ is a principal $S$-bundle.  
Then the transformation groupoid
algebra $C^*(S,X)$ has Hausdorff spectrum $S\unit$. Furthermore the
fibre of $C^*(S,X)$ over $u\in S\unit$ 
is $C^*(S_u,X_u)$ where $X_u = q\inv(u)$.  
\end{prop}

\begin{proof}
Using Lemma \ref{lem:25} we conclude that $C^*(S,X)$ is a
$C_0(S\unit)$-algebra and that restriction factors to an isomorphism of
$C^*(S,X)(u)$ with $C^*(S_u,X_u)$.  Since $C^*(S,X)$ is a
$C^*(S\unit)$-algebra there is a continuous surjection $q$ of
$C^*(S,X)\sidehat$ onto $S\unit$.  Furthermore we identify $q\inv(u)$
with $C^*(S_u,X_u)\sidehat$ in the usual fashion.  
Next, let $\phi:X_u\rightarrow S_u$ be the restriction of one of the
trivializing maps for $X$ to $X_u$.  Since $\phi$ is a homeomorphism
we can pull back the group structure from $S_u$ to $X_u$ and turn
$\phi$ into a group isomorphism.  Furthermore it follows from
Proposition \ref{prop:27} that $\phi$ is equivariant with respect to
the action of $S_u$ on $X_u$.  Therefore
\[
s\cdot \phi\inv(t) = \phi\inv(st)
\]
so that if we identify $X_u$ with $S_u$ then the action of $S_u$ on
$X_u$ becomes the action of $S_u$ on itself by translation.  In other
words $C^*(S_u,X_u)$ is isomorphic to $C^*(S_u,S_u)$.  It follows from the von
Neumann Theorem, Corollary \ref{cor:8}, that $C^*(S_u,S_u)$ is
isomorphic to the compact operators on some separable Hilbert space.
Hence $C^*(S_u,S_u)$, and therefore $C^*(S,X)(u)$, has a unique
irreducible representation.  It follows that the map $q$ is
injective.  

All that remains is to show that $q$ is open, or equivalently,
closed.  Suppose $C$ is a closed subset of $C^*(S,X)\sidehat$.  Then
there is some ideal $I$ such that $C = \{\pi\in C^*(S,X)\sidehat :
I\subset \ker \pi\}$. Let $D = \{ u\in S\unit : I\subset I_u\}$
where the ideal $I_u$ in $C^*(S,X)$ is given by 
\[
I_u = \cspn\{\phi\cdot f:
\phi\in C_0(S\unit), f\in {C_c(S\ltimes X)}, \phi(u) = 0\}.
\]  
We claim
that $D = q(C)$.  If $u\in D$ and $\pi = q\inv(u)$ then $\pi$ factors
to a faithful representation of $C^*(S_u,X_u)$ so that $\ker \pi =
I_u$.  Thus $I \subset I_u = \ker \pi$ and $\pi \in C$.  Conversely if
$\pi \in C$ and $u=q(\pi)$ then $\pi$ factors to a faithful
representation of $C^*(S_u,X_u)$ so that $I_u = \ker \pi$.  It follows
that $I \subset \ker\pi = I_u$ and $u\in D$. All that is left is to show
that $D$ is closed.  Suppose $u_i\rightarrow u$ in $S\unit$ and
$u_i\in D$ for all $i$.  Then since $I\subset I_{u_i}$ for all $i$ we
have $f(u_i) = 0$ for all $f\in I$.  However, $f$ is continuous 
when viewed as a function on $S\unit$ so that $f(u) = 0$.  Thus $f\in
I_u$ and $u\in D$.  
\end{proof}

Next, we show that there is a dual action of $\widehat{S}$ on
$C^*(S,X)$ induced by left translation.  Since it isn't much harder, 
we actually prove this result in greater generality. Unfortunately,
we can't just jump right in.  Verifying the continuity condition will
take work.  In particular we have to deal with the topology associated
to the crossed bundle product $A\rtimes S$.  

\begin{lemma}
\label{lem:39}
Suppose $(A,S,\alpha)$ is a separable dynamical system and that 
$S$ is an abelian group bundle.  Let $\mcal{A}$ be the
upper-semicontinuous bundle associated to $A$, and define
\[
\widehat{S}*S = \{(\omega,s)\in \widehat{S}\times S : \hat{p}(\omega)
= p(s)\},
\]
and let $p:\widehat{S}*S \rightarrow S\unit$ be given by $p(\omega,s)
= p(s)$.  Then there is a map
$\iota:\Gamma_c(\widehat{S}*S,p^*\mcal{A})\rightarrow \hat{p}^*
(A\rtimes S)$ such that 
\[
\iota(f)(\omega)(s) = f(\omega, s).
\]
Furthermore $\iota$ is continuous with respect to the inductive limit
topology and the range of $\iota$ is dense.  
\end{lemma}

\begin{proof}
First observe that, since $f$ is continuous and compactly supported, 
$\iota(f)(\omega)$ will be a continuous, compactly supported
function from $S_{\hat{p}(\omega)}$ into $A(\hat{p}(\omega))$.  
Thus $\iota(f)(\omega)\in C_c(S_{\hat{p}(\omega)},
A(\hat{p}(\omega)))$ and  $\iota(f)$ is a section of
$\hat{p}^*(A\rtimes S)$.  Furthermore, it is clear that $\iota(f)$ has
compact support.  We would like to show that $\iota(f)$ is
continuous as a function into $\mcal{E}$ where $\mcal{E}$ is the
upper-semicontinuous bundle associated to $A\rtimes S$.  

We start out with a simpler function.  Suppose $g\in
C_c(\widehat{S})$, $h\in C_c(S)$ and $a\in A$.  Define $g\otimes
h\otimes a$ on $\widehat{S}*S$  by $g\otimes h \otimes a(\omega,s) =
g(\omega)h(s)a(p(s))$.  It is clear that $g\otimes h\otimes a\in
\Gamma_c(\widehat{S}*S,p^*\mcal{A})$.  Furthermore if we view
$h\otimes a$ as an element of $\Gamma_c(S,p^*\mcal{A})$ then
$\iota(g\otimes h\otimes a)(\omega) = g(\omega) (h\otimes a)(\hat{p}(\omega))$
where $(h\otimes a)(\hat{p}(\omega))$ is just the restriction of
$h\otimes a$ to $S_{\hat{p}(\omega)}$.  Since $h\otimes a$ defines a
continuous section of $\mcal{E}$ it is easy to see that
$\iota(g\otimes h\otimes a)$ is a continuous function from
$\widehat{S}$ into $\mcal{E}$.  Thus $\iota(g\otimes h\otimes a) \in
\Gamma_c(\widehat{S},\hat{p}^*\mcal{E})$.  

We now show $\iota$ preserves convergence with respect to the
inductive limit topology.  Suppose $f_i\rightarrow f$ uniformly in
$\Gamma_c(\widehat{S}*S,p^*\mcal{A})$ and that eventually the supports
are contained in some fixed compact set $K$.  Clearly the supports of
$\iota(f)$ are eventually contained in the projection of $K$ to
$\widehat{S}$.  Fix $\epsilon > 0$ and let $M$ be an upper bound for
$\{\beta^u(L)\}$ where $L$ is the projection of $K$ to $S$.  Then
eventually $\|f_i-f\|_\infty < \epsilon/L$.  Thus for large $i$ we
have, for $\omega \in \widehat{S}_u$ and making use of the fact that
$S_u$ is abelian so the $I$-norm on $C_c(S_u,A(u))$ only has one term, 
\begin{align*}
\|\iota(f_i)(\omega)-\iota(f)(\omega)\|&\leq \|\iota(f_i)(\omega) -
\iota(f)(\omega)\|_I \\
&= \int_S \|f_i(\omega,s) - f(\omega,s)\| d\beta^u(s) \\
&\leq \|f_i-f\|_\infty L < \epsilon 
\end{align*}
Thus $\iota(f_i)\rightarrow \iota(f)$ uniformly and hence with
respect to the inductive limit topology.  

Now suppose $f\in \Gamma_c(\widehat{S}*S,p^*\mcal{A})$ and that
$\omega_i\rightarrow \omega$ in $\widehat{S}$. Fix $\epsilon > 0$ and
let $U$ and $V$ be relatively compact neighborhoods of the projections
of $\supp f$ to $\widehat{S}$ and $S$ respectively.  
Since $f$ is a
section of a pull back bundle we can use Proposition \ref{prop:46} to
find $\{g_i\}_{i=1}^N\in C_c(\widehat{S}*S)$ and $\{a_i\}_{i=1}^N\in A$
such that $\|f-\sum_i g_i\otimes
a_i\|_\infty < \epsilon/2$.  For each $1\leq i \leq N$ extend $g_i$ to all of
$C_c(\widehat{S}\times S)$ and choose $h_i^j\in C_c(\widehat{S})$ and
$k_i^j\in C_c(S)$ such that $\|g_i-\sum_jh_i^j\otimes k_i^j\|_\infty < \epsilon/(2N\|a_i\|)$.
It then follows from some simple computations that 
\begin{align*}
\left\|f - \sum_{i=1}^N\sum_j h_i^j\otimes k_i^j \otimes a_i\right\|_\infty \leq&
\left\|f- \sum_{i=1}^N g_i \otimes a_i\right\|_\infty \\ &+ 
\left\|\sum_{i=1}^N g_i \otimes a_i - \sum_{i=1}^N \sum_j h_i^j\otimes
  k_i^j \otimes a_i \right\|_\infty \\
\leq & \epsilon/2 +
\sum_{i=1}^N\|a_i\|\left\|g_i-\sum_j h_i^j\otimes k_i^j\right\|_\infty
< \epsilon. 
\end{align*}
Furthermore, we can multiply the $h_i^j$ and $k_i^j$ by 
functions which vanish off
$U$ and $V$ respectively so that $\supp h_i^j\otimes k_i^j\otimes a \subset
\overline{U}\times \overline{V}$.  This construction shows that
sums of elements of the form $h\otimes k\otimes a$ for $h\in
C_c(\widehat{S})$, $k\in C_c(S)$ and $a\in A$ are dense in
$\Gamma_c(\widehat{S}*S, p^*\mcal{A})$ with respect to the inductive
limit topology.  

At last we can show that $\iota(f)$ is continuous for 
$f\in \Gamma_c(\widehat{S}*S,p^*S)$.  Let $g_i = \sum_k h_i^k\otimes
k_i^k\otimes a_i^k$ 
be a sequence converging to $f$ in the inductive limit topology as
above.  Let $\omega_j\rightarrow \omega$ and fix $\epsilon > 0$.  For
some very large $I$ we have $\|\iota(f)-\iota(g_I)\|_\infty <
\epsilon$.  In particular
$\|\iota(f)(\omega)-\iota(g_I)(\omega)\|<\epsilon$ and
$\|\iota(f)(\omega_j)-\iota(g_I)(\omega_j)\|<\epsilon$ for all $i$.
Since sums of continuous functions are continuous, 
$\iota(g_I)(\omega_j)\rightarrow \iota(g_I)(\omega)$ and it follows
from the last part of Proposition \ref{prop:35} that $\iota(f)$ is
continuous. 

Thus $\iota$ maps $\Gamma_c(\widehat{S}*S,p^*\mcal{A})$ into
$\hat{p}^*(A\rtimes S)$.  The last thing we need to do is verify that $\iota$
has dense range.  Given $\phi\in C_0(\widehat{S})$ and $f\in
\Gamma_c(\widehat{S}*S,p^*\mcal{A})$ we can define a new function
$h\in\Gamma_c(\widehat{S}*S,p^*\mcal{A})$ by $h(\omega,s) =
\phi(\omega)h(\omega,s)$.  It is easy enough to see that
$\phi\cdot\iota(f) = \iota(h)$.  Thus $\ran\iota$ is closed under the
$C_0(\widehat{S})$ action.  Now fix $\omega\in \widehat{S}_u$.  
Given $a\in A$ and $f\in C_c(S_u)$ extend $f$ to a function $h\in
C_c(S)$ and choose $g\in C_c(\widehat{S})$ so
that $g(\omega) = 1$.  Then clearly $\iota(g\otimes f\otimes
a)(\omega) = h|_{S_u}\otimes a = f\otimes a$.  Thus $\ran\iota$
contains sums of elementary tensors in $A(u)\rtimes S_u$.  
It follows from Proposition \ref{prop:42} that $\ran\iota$ is
dense.
\end{proof}

The following corollary isn't necessary to build the dual action, but
it will be needed in the proof of Theorem \ref{thm:exist} so we
include it here.  

\begin{corr}
\label{cor:12}
Suppose $S$ is an abelian, second countable, locally compact Hausdorff
continuously varying group bundle and $q:X\rightarrow S\unit$ is a principal
$S$-bundle. Define 
\[
\widehat{S}*S*X := \{(\omega,s,x)\in \widehat{S}\times S\times X :
\hat{p}(\omega) = p(s) = q(x)\}
\]
Then there is a map $\iota:C_c(\widehat{S}*S*X)\rightarrow
\hat{p}^*C^*(S,X)$ such that 
\[
\iota(f)(\omega)(s,x) = f(\omega,s,x).
\]
Furthermore $\iota$ is continuous with respect to the inductive limit
topology and the range of $\iota$ is dense. 
\end{corr}

\begin{proof}
Let $\mcal{C}$ be the upper-semicontinuous bundle associated to
$C_0(X)$.  Since both algebras can be viewed as completions of
$C_c(S\ltimes X)$ we will use Proposition \ref{prop:68} to identify
$C_0(X)\rtimes S$ with $C^*(S,X)$.  In particular, we will not
distinguish between functions in $C_c(S\ltimes X)$ and their
corresponding functions in $\Gamma_c(S,p^*\mcal{C})$.  Consider the map
$\iota_1:\Gamma_c(\widehat{S}*S,p^*\mcal{C})\rightarrow
\hat{p}^*C^*(S,X)$ given by 
\[
\iota_1(f)(\omega)(s,x) = \iota_1(\omega)(s)(x) := f(\omega,s)(x).  
\]
It follows from Lemma \ref{lem:39} that this map is
continuous with respect to the inductive limit topology and its range
is dense in $\hat{p}^*C^*(S,X)$.  Now consider the map
$\iota_2:C_c(\widehat{S}*S*X)\rightarrow
\Gamma_c(\widehat{S}*S,p^*\mcal{C})$ given by $\iota_2(f)(\omega,s)(x)
= f(\omega,s,x)$. It follows from Lemma \ref{lem:26} that $\iota_2$ is
surjective and preserves the inductive limit topology.  Thus
the map $\iota = \iota_2\circ\iota_1$ has the correct form and all the
right properties. 
\end{proof}

Now we can finally tackle the dual action construction.  This will
provide the last tool we need to demonstrate Theorem \ref{thm:exist}.  

\begin{prop}
\label{prop:89}
\index{dual action}
\index{dual action|see{Takai duality}}
\index{dual bundle}
Suppose $(A,S,\alpha)$ is a separable dynamical system and that $S$ is
an abelian group bundle.  Then for each $\omega\in \widehat{S}$ there is
an automorphism $\hat{\alpha}_\omega$ on $A\rtimes S(\hat{p}(\omega))$
defined for $f\in C_c(S_{\hat{p}(\omega)},A(\hat{p}(\omega)))$ by 
\[
\hat{\alpha}_\omega(f)(s) = \overline{\omega(s)}f(s).
\]
With this action $(A\rtimes S, \widehat{S},\hat{\alpha})$ is a
dynamical system. 
\end{prop}

\begin{proof}
Recall that $A\rtimes S$ is a $C_0(S\unit)$-algebra with restriction
factoring to an isomorphism of $A\rtimes S(u)$ with $A(u)\rtimes
S_u$. Fix $u\in S\unit$, $\omega\in \widehat{S}_u$ and define
$\alpha_\omega$ as above.  Since $\omega(s)$ is unimodular for all
$s\in S$, it follows from a simple calculation 
that $\alpha_\omega$ is continuous with respect to the
inductive limit topology.  Furthermore, we have for $f,g\in
C_c(S_u,A(u))$
\begin{align*}
\alpha_\omega(f)*\alpha_\omega(g)(s) &= \int_S \overline{\omega(t)}
f(t) \alpha_t(\overline{\omega(t\inv s)} g(t\inv s)) d\beta^u(t) \\
&= \overline{\omega(s)} f*g(s) = \alpha_\omega(f*g)(s),\quad\text{and} \\
\alpha_\omega(f)^*(s) &= \alpha_s(\omega(s\inv)f(s\inv)^*) \\
&= \overline{\omega(s)}f^*(s) = \alpha_\omega(f^*)(s).
\end{align*}
Thus $\alpha_\omega$ is a $*$-homomorphism and it follows from
Corollary \ref{cor:24} that $\alpha_\omega$ extends to a map from
$A(u)\rtimes S_u$ into itself.  Furthermore observe that
$\alpha_u = \id$ and 
\[
\alpha_\omega(\alpha_\chi(f))(s) = \overline{\omega(s)\chi(s)}f(s) = 
\alpha_{\omega\chi}(f)(s).
\]
It follows that, one, $\alpha_\omega$ is an automorphism, and, two, that
$\alpha$ respects the groupoid structure.  All that remains is to
verify the continuity condition.  

Let $\mcal{E}$ be the bundle associated to $A\rtimes S$ and suppose 
$\omega_i\rightarrow \omega$ in $\widehat{S}$ and $f_i\rightarrow f$
in $\mcal{E}$.  Now choose $g\in \hat{p}^*(A\rtimes S)$ such that
$g(\omega) = f$.  It follows from Lemma \ref{lem:25} that we can view
$\Gamma_c(\widehat{S}*S,p^*\mcal{A})$ as lying inside
$\hat{p}^*(A\rtimes S)$ and can
choose $h\in \Gamma_c(\widehat{S}*S,p^*\mcal{A})$ such that
$\|h-g\|_\infty < \epsilon/2$.  Define $\alpha(h)(\omega,s) =
\overline{\omega(s)} h(\omega,s)$.  It is clear that $\alpha(h)\in
\Gamma_c(\widehat{S}*S,p^*\mcal{A})$.  It is also easy to see that
$\alpha(h)(\omega) = \alpha_\omega(h(\omega))$.  Thus
\[
\|\alpha(h)(\omega)-\alpha_\omega(f)\| =
\|\alpha_\omega(h(\omega)-g(\omega))\| < \epsilon/2 < \epsilon.
\]
Next, since $g(\omega_i)\rightarrow g(\omega) =f$ and $f_i\rightarrow
f$ we have $\|g(\omega_i)-f_i\|\rightarrow 0$.  Therefore, eventually,
we have 
\begin{align*}
\|\alpha(h)(\omega_i)-\alpha_{\omega_i}(f_i)\| &\leq
\|\alpha_{\omega_i}(h(\omega_i)-g(\omega_i))\| +
\|\alpha_{\omega_i}(g(\omega_i)-f_i)\| \\
&\leq \epsilon/2 + \|g(\omega_i)-f_i\| < \epsilon.
\end{align*}
It follows from Proposition \ref{prop:35} that
$\alpha_{\omega_i}(f_i)\rightarrow \alpha_\omega(f)$.  
\end{proof}

\begin{remark}
\index{Takai duality}
The action from Proposition \ref{prop:89} is a generalization of the
usual Takai dual action for abelian groups \cite[Section 7]{tfb2}.  In
particular there is a Takai Duality Theorem for abelian group bundles
which states that $(A\rtimes_\alpha S)\rtimes_{\hat{\alpha}}
\widehat{S}$ is isomorphic to $A\otimes_{C_0(X)}\mcal{K}(\mcal{Z})$
where $\mcal{Z}$ is the Hilbert $C_0(S\unit)$-bundle from 
Theorem \ref{thm:stonevn}.  However, this theorem is really nothing
more than a ``bundled'' version of the theorem for groups.  Although
it takes some work, the proof boils down to observing that the
isomorphism given in \cite[Section 7]{tfb2} respects the total space
structure.  Since there is no interesting ``groupoid'' component to
this result, it has been omitted.  
\end{remark}

We are now ready to prove our existence theorem. 

\begin{proof}[Proof of Theorem \ref{thm:exist}]
We have shown in Proposition \ref{prop:90} 
that $C^*(S,X)$ has Hausdorff spectrum $S\unit$ and that
restriction factors to an isomorphism of $C^*(S,X)(u)$ with
$C^*(S_u,X_u)$ where $X_u = q\inv(u)$.  Furthermore we showed in
Proposition \ref{prop:89} that
there is an action of $\widehat{S}$ which, after making the usual
identification of $C^*(S,X)$ with $C_0(X)\rtimes_{\lt} S$, is given by 
\[
\widehat{\lt}_\omega(f)(s,x) = \overline{\omega(s)}f(s,x)
\]
for $f\in C_c(S_u\times X_u)$.  We need to show that $\widehat{\lt}$
is locally unitary. Let $\mcal{U}$ be a trivializing cover of $X$ and
let $\phi_i$ be the local trivializations.  
Fix $U_i\in\mcal{U}$.  Then for all 
$w\in U_i$, $\omega \in \widehat{S}_w$ and $f\in C_c(S_w\times X_w)$
define 
\begin{equation}
\label{eq:108}
u_\omega f(s,x) := \overline{\omega(\phi_i(x))}f(s,x).
\end{equation}
Viewing $w$ as the trivial element\footnote{It even looks
  like it should be a character! $(\omega \approx w)$} in $\widehat{S}_w$ we clearly have $u_w =
\id$. Furthermore, given $\omega,\chi\in \widehat{S}_w$, 
\[
u_{\omega\chi}f(s,x)=
\omega(\phi_i(x))\chi(\phi_i(x))f(s,x) = u_\omega u_\chi
f(s,x).
\]
Thus $u$ is a homomorphism on $S_w$.  Next we will show that $u$ is
adjointable.  Recall that we equip $C^*(S_w,X_w)$ with the inner
product $\langle f,g\rangle = f^**g$.  For all $f,g\in C_c(S_w\times
X_w)$ we have 
\begin{align*}
\langle u_\omega f,g\rangle(s,x) &= (u_\omega f)^* * g(s,x) \\
&= \int_S \overline{u_\omega f(t\inv,t\inv\cdot x)}g(t\inv s, t\inv\cdot
x) d\beta^w(t) \\
&= \int_S \omega(\phi_i(t\inv\cdot x)) \overline{f(t\inv,t\inv\cdot x)}g(t\inv
s, t\inv\cdot x)d\beta^w(t) \\
&= \int_S \overline{f(t\inv,t\inv\cdot x)}u_{\omega\inv}g(t\inv s,
t\inv\cdot x) d\beta^w(t) \\
&= f^**(u_{\omega\inv}g)(s,x) = \langle f,u_{\omega\inv}g\rangle(s,x).
\end{align*}
This shows that $u_\omega$ is adjointable on $C_c(S_w,X_w)$ and we can
also observe that 
\begin{align*}
\|u_\omega f\|^2 &= \|\langle u_\omega f, u_\omega f\rangle\| = 
\|\langle f, u_{\omega\inv}u_\omega f\rangle \| \\
&= \|\langle f, f\rangle\| = \|f\|^2.
\end{align*}
Thus $u_\omega$ is isometric on $C_c(S_w,X_w)$ and as such it can be
extended to an operator on $C^*(S_w,X_w)$.  Since $u_\omega$ is
adjointable on a dense subspace with $u_\omega^* = u_{\omega\inv} =
u_\omega\inv$, we know $u_\omega$ is a unitary multiplier on
$C^*(S_w,E_w)$.  Consider the collection
$\{u_\omega\}_{\omega\in\hat{p}\inv(U_i)}$.  All that remains for
$\{u_\omega\}$ to define a unitary action of $\hat{p}\inv(U_i)$ on
$C^*(S,X)(U_i)$ is continuity.  

Let $\mcal{E}$ be the bundle associated to $C^*(S,X)$ and fix
$\epsilon > 0$.  Suppose $\omega_j\rightarrow \omega$ in $\hat{p}\inv(U_i)$ and
$f_j\rightarrow f$ in $\mcal{E}|_{U_i}$.  Choose $g\in \hat{p}^*
C^*(S,X)$ such that $g(\omega)=f$.   Using Corollary \ref{cor:12} we can find a
continuous, compactly supported function $h$ on $\widehat{S}*S*X$ such that
$\|\iota(h)-g\|_\infty < \epsilon/2$.  Consider the open set $O =
\hat{p}\inv(U_i)*p\inv(U_i)*q\inv(U_i)$ in $\widehat{S}*S*X$. We
define a new function $k\in C(O)$ by 
\[
k(\chi,s,x) = \psi(p(s))\overline{\chi(\phi_i(x))}h(\chi,s,x)
\]
where $\psi\in C_c(U_i)$ is some function which is one on a
neighborhood of $\hat{p}(\omega)$. Now $k$ is clearly compactly
supported with $\supp k \subset
\hat{p}\inv(\supp\psi)*p\inv(\supp\psi)*q\inv(\supp\psi)\subset O$.
Therefore we can, and do, extend $k$ by zero to all of $\widehat{S}*S*X$.  Next
we observe the following facts.  First, 
\[
\iota(k)(\omega)(s,x) =
\overline{\omega(\phi_i(x))}h(\omega,s,x) 
= u_\omega \iota(h)(\omega)(s,x).
\]
In a similar fashion we see that eventually $\iota(k)(\omega_j) =
u_{\omega_j} \iota(h)(\omega_j$).  Second, observe that 
\[
\|\iota(k)(\omega) - u_\omega f\| = \|u_\omega(\iota(h)(\omega) -
g(\omega))\| = \|\iota(h)(\omega) - g(\omega)\| < \epsilon/2.
\]
Furthermore $f_i\rightarrow f$ and $g(\omega_i)\rightarrow g(\omega) =
f$ so that $\|f_i-g(\omega_i)\|\rightarrow 0$.  Thus, eventually,
we have 
\begin{align*}
\|\iota(k)(\omega_i) - u_{\omega_i}f_i\| &\leq
\|u_{\omega_i}(\iota(h)(\omega_i) - g(\omega_i)) \| +
\|u_{\omega_i}(g(\omega_i) - f_i)\| \\
&\leq \epsilon/2 + \|g(\omega_i)-f_i\| < \epsilon.
\end{align*}
Finally we observe that $\iota(k)(\omega_i)\rightarrow
\iota(k)(\omega)$ since $\iota(k)$ is a continuous section.  It
follows that $u_{\omega_i}f_i\rightarrow u_\omega f$ and that
$\{u_\omega\}$ defines a unitary action of $\hat{p}\inv(U_i)$ on
$C^*(S,X)(U_i)$.  

Next, we show that $u$ implements $\widehat{\lt}$ on
$\hat{p}\inv(U_i)$.  Suppose $u\in U_i$ and $\omega\in
\widehat{S}_u$.  Then for $f\in C_c(S_u\times X_u)$ we have 
\begin{align*}
u_\omega f u_\omega^*(s,x) &= \overline{\omega(\phi_i(x))} (u_\omega
f^*)^*(s,x) = \overline{\omega(\phi_i(x)) u_\omega
  f^*(s\inv,s\inv\cdot x)} \\
&= \overline{\omega(\phi_i(x))} \omega(\phi_i(s\inv\cdot x)) f(s,x) = 
\overline{\omega(\phi_i(x))}\overline{\omega(s)}\omega(\phi_i(x))
f(s,x) \\
&= \overline{\omega(s)}f(s,x) = \widehat{\lt}f(s,x)
\end{align*}
where we have used the fact that $\phi_i$ is equivariant with respect to
the action of $S$ on $X$.  Thus  $\widehat{\lt}$ is unitarily
implemented on $\hat{p}\inv(U_i)$.  Since we performed this
construction for each element
of the cover $\mcal{U}$ it follows that $\widehat{\lt}$ is locally
unitary.  

Consider $Y = (C^*(S,X)\rtimes \widehat{S})\sidehat$.  Now, $Y$ is
a principal $\doubledual{S}$-bundle and in light of Theorem
\ref{thm:duality} a principal $S$-bundle as well.  We would
like to show that $Y$ is isomorphic to $X$.  
Using Theorem \ref{prop:principcohom} it suffices to
show that $X$ and $Y$ have the same cohomological invariant.  Let
$\gamma_{ij}$ be the transition functions for $X$ with respect to the
trivializing maps $\phi_i$.  Let $\eta_{ij}$ be the transition
functions for $Y$ and recall from Theorem \ref{thm:locunit} that for
$v\in U_{ij}$ and $\omega \in \widehat{S}_v$ 
\[
\eta_{ij}(v)(\omega) = \overline{\pi}_v((u_\omega^i)^*u_\omega^j)
\]
where
$\pi_v$ is the unique irreducible representation of
$A(v)$ and $u_\omega^i, u_\omega^j$ are the unitaries
constructed above.  Recall that 
$\phi_i \circ\phi_j\inv(s)  = \gamma_{ij}(p(s)) s$
so that, letting $x=\phi_j\inv(s)$, we have 
$\phi_i(x) = \gamma_{ij}(q(x))\phi_j(x)$.
We now compute for $f\in C_c(S_v,X_v)$ 
\begin{align*}
((u_\omega^i)^*u_\omega^j f)(s,x) &=
\overline{\omega\inv(\phi_i(x))\omega(\phi_j(x))}f(s,x) \\
&=
\omega(\gamma_{ij}(v)\phi_j(x))\overline{\omega(\phi_j(x))}f(s,x)\\
&= \omega(\gamma_{ij}(v))f(s,x) \\
&= \hat{\gamma}_{ij}(v)(\omega) f(s,x)
\end{align*}
where $\hat{\gamma}_{ij}(v)$ denotes the image of $\gamma_{ij}(v)$ in
the double dual.  Therefore $(u_\omega^i)^* u_\omega^j =
\hat{\gamma}_{ij}(v)(\omega)$ and, since $\pi_v$ is faithful,
$\eta_{ij}(v)(\omega) = \hat{\gamma}_{ij}(v)(\omega)$.  Thus, once we
identify $S$ with $\doubledual{S}$, the cohomological invariants of
$X$ and $Y$ are identical.  
\end{proof}

\begin{example}
\index{locally sigma trivial@locally $\sigma$-trivial}
Theorem \ref{thm:exist} says that any
principal $S$-bundle $X$ gives rise to a locally unitary action of
$\widehat{S}$ on $C^*(S,X)$ and in particular this holds for locally
$\sigma$-trivial bundles.  Thus Examples \ref{ex:13} and \ref{ex:26}
yield examples of locally unitary actions.  
\end{example}

\begin{remark}
It is worth describing, at least briefly, how this material
generalizes \cite{locunitary}.  
Suppose $H$ is an abelian group and $A$ has
Hausdorff spectrum $X$.  If $\alpha$ is an action of $H$ on $A$
then, as in Example \ref{ex:22}, we can form the transformation
groupoid $H\ltimes X$ and there is an action $\beta$ of $H\ltimes X$ on
$A$.  Furthermore we have $A\rtimes_\alpha X \cong A\rtimes_\beta
(H\ltimes X)$.  Without getting into the details, $\alpha$ is locally
unitary according \cite{locunitary} if, for each $\pi\in X$, there is an
open neighborhood $U$ of $\pi$ and a strictly continuous map $u:H\rightarrow
M(A)$ such that for each $\rho\in U$, $\bar{\rho}\circ u$ is a
representation of $H$ on $\mcal{H}_\rho$ which implements $\alpha$.  In
particular this implies that $\rho\circ \alpha_s\inv =
\bar{\rho}(u_s)\rho \bar{\rho}(u_s^*)$ is equivalent to $\rho$.
Thus the action of $H$ on $X$ induced by $\alpha$ is trivial and the
transformation groupoid $H\ltimes X$ is the trivial
group bundle.  What's more, since $A$ has Hausdorff spectrum, 
it is not hard to show
that $u_s(x)$ implements $\beta_{(s,x)}$ on $A(x)$ and that $\beta$ is
unitarily implemented by $\{u_s(x)\}$ on $H\times U$.  Thus $\beta$ is
a locally unitary action of $H\times S$ on $A$.  Now,
the dual of $H\times X$ is $\widehat{H}\times X$
and we have, according to Theorem \ref{thm:locunit}, that
$(A\rtimes_\alpha X)\sidehat \cong
(A\rtimes_\beta(H\times X))\sidehat$ is a principal $\widehat{H}\times
X$ bundle.  However, it follows from Proposition \ref{prop:91} that this
implies $(A\rtimes_\alpha X)\sidehat$ is a principal
$\widehat{H}$-bundle.  From here it is straightforward to see how the results of
this section are related to those in \cite{locunitary}.
\end{remark}


\chapter{Fine Structure of Groupoid Crossed Products}
\label{cha:fine-structure}
In this chapter we present the main results of the thesis.  An
important aspect of the proof is the induction process detailed in
Section \ref{sec:indreps}.  This allows us to induce representations
from the crossed product by any closed subgroupoid to the whole
crossed product.  We will use this in Section \ref{sec:regularity} to
show that when the orbits are $T_0$ every irreducible representation
of the crossed product is induced from a stabilizer
subgroup.  We use this result in Section \ref{sec:crossedstab} to
identify the spectrum of $A\rtimes G$ with a quotient of the spectrum
of $A\rtimes S$, where $S$ is the stabilizer subgroupoid of $G$.  

\section{Induction}
\label{sec:indreps}

A key tool in our study of the representation theory of crossed
products will be the ability to induce representations from closed
subgroupoids.  This notion has been around for groups since \cite{frobenius},
although this particular section is more closely modeled after
\cite[Chapter 5]{tfb2}.  Our starting point will be the imprimitivity
groupoid from Section \ref{sec:equivalence}.  Specifically,
we begin by describing an action which arises from coupling
the imprimitivity groupoid with a dynamical system.

\begin{prop}
\label{prop:79}
\index{dynamical system}
\index{imprimitivity groupoid}
\index{pull back}
Let $(A,G,\alpha)$ be a separable groupoid dynamical system and suppose
$H$ is a closed subgroupoid of $G$ with a Haar system.  Let $X =
s\inv(H\unit)$ and let $G^H$ be the associated imprimitivity
groupoid.  Define $\rho:X/H\rightarrow G\unit$ via $\rho(\gamma\cdot H) =
r(\gamma)$ and let $\rho^*A$ be the pull back algebra. Then the collection
$\sigma = \{\sigma_{[\gamma,\eta]}\}_{[\gamma,\eta]\in G^H}$ where 
$\sigma_{[\gamma,\eta]}:A(r(\eta))\rightarrow A(r(\gamma))$ is defined
for $a\in A(r(\gamma))$ by
\begin{equation}
\label{eq:81}
\sigma_{[\gamma,\eta]}(a) = \alpha_{\gamma\eta\inv}(a)
\end{equation}
defines an action of $G^H$ on $\rho^*A$. 
\end{prop}

\begin{proof}
It follows from Proposition \ref{prop:45} that the pull back $\rho^*A$
is a $C_0(X/H)$-algebra with fibres $A(\gamma\cdot H) =
A(r(\gamma))$.  Recall that we identify $(G^H)\unit$ with $X/H$ so
that we may also view $\rho^*A$ as a $C_0((G^H)\unit)$-algebra.  Next, we
want to show that $\sigma_{[\gamma,\eta]}$ is independent of the
choice of representatives $\gamma$ and $\eta$.  If
$[\gamma,\eta] = [\gamma',\eta']$ then there exists $\zeta\in H$ such
that $\gamma \zeta= \gamma'$ and $\eta\zeta = \eta'$.  However this
implies that $\gamma\eta\inv = \gamma'(\eta')\inv$ and that
$\sigma_{[\gamma,\eta]}$ is well defined. 

Moving on, it is clear
that $\sigma_{[\gamma,\eta]} = \alpha_\gamma\circ \alpha_\eta\inv$ is
an isomorphism of $A(s([\gamma,\eta])) = A(r(\eta))$ onto
$A(r([\gamma,\eta])) = A(r(\gamma))$.  
Next, if $[\gamma,\eta],[\eta,\zeta]\in G^H$ and
$a\in A(r(\zeta))$ then 
\[
\sigma_{[\gamma,\eta]}\circ\sigma_{[\eta,\zeta]}(a) =
\alpha_{\gamma\eta\inv}(\alpha_{\eta\zeta\inv}(a)) =
\alpha_{\gamma\zeta\inv}(a) = \sigma_{[\gamma,\zeta]}(a).
\]
Thus the action respects the groupoid operations.  Lastly we have to
show that if $[\gamma_i,\eta_i]\rightarrow [\gamma,\eta]$ in $G^H$ and
$a_i\rightarrow a$ in $\mcal{A}$ such that $p(a_i) = r(\eta_i)$ for
all $i$ and $p(a) = r(\eta)$ then 
\begin{equation}
\label{eq:82}
\sigma_{[\gamma_i,\eta_i]}(a_i) = \alpha_{\gamma_i\eta_i\inv}(a_i)
\rightarrow \alpha_{\gamma\eta\inv}(a) = \sigma_{[\gamma,\eta]}(a).
\end{equation}
After passing to a subnet it will suffice to show that a sub-subnet
converges.  However, we can pass to a subnet, relabel and choose new
representatives so that $\gamma_i\rightarrow \gamma$ and
$\eta_i\rightarrow \eta$.  However, \eqref{eq:82} now holds
because $\alpha$ is a continuous action. 
\end{proof}

\begin{remark}
The crossed product $\rho^*A\rtimes_\sigma G^H$ is the completion
of $\Gamma_c(G^H,r^*(\rho^*\mcal{A}))$.  As in Remark \ref{rem:8}, 
elements of this function
algebra can be viewed as continuous, compactly supported maps from
$G^H$ into $\mcal{A}$ such that $f([\gamma,\eta])\in A(r(\gamma))$ for
all $[\gamma,\eta]\in G^H$. 
\end{remark}

Now we use the equivalence theorem to build an imprimitivity
bimodule.  This shows that up to Morita equivalence 
we really don't get anything new from
$\sigma$, and that it is equivalent to the restriction of $\alpha$ to
$H$.  

\begin{prop}
\label{prop:80}
\index{dynamical system equivalence}
Suppose $(A,G,\alpha)$ is a separable dynamical system and let
$\lambda$ be a Haar system for $G$.  Furthermore, suppose $H$ is a
closed subgroupoid of $G$ with Haar system $\lambda_H$.  Let
$X=s\inv(H\unit)$, $G^H$ be the imprimitivity groupoid, and $\sigma$
be the action of $G^H$ on $\rho^*A$.  Then $\mcal{Z}_0 =
\Gamma_c(X,s^*\mcal{A})$ becomes a pre-$\rho^*A\rtimes_\sigma
G^H-A(H\unit)\rtimes_{\alpha|_H} H$-imprimitivity bimodule with
respect to the following actions for $f\in \Gamma_c(G^H,r^*(\rho^*\mcal{A}))$,
$g\in\Gamma_c(H,r^*\mcal{A})$, and $z,w\in \mcal{Z}_0$:
\begin{align}
\label{eq:83}
f\cdot z(\gamma) &= \int_G
\alpha_\gamma\inv(f([\gamma,\eta]))z(\eta)d\lambda_{s(\gamma)}(\eta)
\\
\label{eq:84}
z\cdot g(\gamma) &= \int_H
\alpha_\eta(z(\gamma\eta)g(\eta\inv))d\lambda_H^{s(\gamma)}(\eta) \\
\label{eq:85}
\llangle z,w \rrangle_{A\rtimes H}(\eta) &= \int_G
z(\xi\eta\inv)^* \alpha_\eta(w(\xi))d\lambda_{s(\eta)}(\xi)
\\ \label{eq:86}
\lset{\rho^*A\rtimes G^H}\llangle z,w \rrangle([\gamma,\eta]) &= \int_H
\alpha_{\gamma\xi}(z(\gamma \xi)w(\eta\xi)^*)d\lambda_H^{s(\gamma)}(\xi)
\end{align}
The completion $\mcal{Z}_H^G$ of $\mcal{Z}_0$ is a
$\rho^*A\rtimes_\sigma G^H-A(H\unit)\rtimes_{\alpha|_H} H$-imprimitivity
bimodule and $\rho^*A\rtimes_\sigma G^H$ and $A(H\unit)\rtimes_\alpha H$ are
Morita equivalent. 
\end{prop}

\begin{proof}
Let $A,G,\alpha,$ and $H$ be as above and suppose $\mcal{A}$ is the
upper-semicontinuous bundle associated to $A$. 
Let $X=s\inv(H\unit)$ be the canonical
$(G^H,H)$-equivalence, let $s_X$ be the
restriction of the source map to $X$, and let $r_X:X\rightarrow X/H$ be the
quotient map.  
Consider the pull back $\erune = s_X^*(\mcal{A})$.  This is clearly an
upper semicontinuous bundle and we define 
$\mcal{Z}_0=\Gamma_c(X,s^*\mcal{A})$.  We will construct an
equivalence from $\erune$.  Observe that given $\gamma\in X$ we have
$\erune_\gamma = \{\gamma\}\times A(s(\gamma))$, 
$A(s_X(\gamma))=A(s(\gamma))$ and
$\rho^*A(r_X(\gamma)) = \rho^*A(\gamma\cdot H) = A(r(\gamma))$.  Thus we can equip
$\erune_\gamma$ with the $A(r(\gamma))-A(s(\gamma))$-imprimitivity
bimodule structure coming from the isomorphism $\alpha_\gamma$
\cite[Example 3.14]{tfb}.  Specifically, given $a,b\in
A(s(\gamma))$ and $c\in A(r(\gamma))$ we have 
\begin{align*}
c\cdot (\gamma,a) &= (\gamma,\alpha_\gamma\inv(c)a), & (\gamma,a)\cdot b = (\gamma,ab), \\
\lset{A(r(\gamma))}\langle (\gamma,a),(\gamma,b) \rangle &= \alpha_\gamma(ab^*),  &
\langle (\gamma,a),(\gamma,b) \rangle_{A(s(\gamma))} = a^*b.
\end{align*}

Next, let $p\suberune:\erune\rightarrow X$ be the bundle map and 
define $r\suberune:\erune\rightarrow X/H$ so that
$r\suberune(\gamma,a) = \gamma\cdot H$ and $s\suberune:\erune\rightarrow
H\unit$ so that $s(\gamma,a) = s(\gamma)$.  These maps are clearly
continuous.  Furthermore $r\suberune$ is just the the composition of
$p\suberune$ and $r_X$ so that $r\suberune$ is open.  Similarly $s\suberune
= s_X\circ p\suberune$ so $s\suberune$ is open as well.  (Recall that
$s_X$ is open because $X$ is a saturated closed set in $G$.)  Now we
define actions of $G^H$ and $H$ on $\erune$ for $[\eta,\zeta]\in
G^H$, $\xi\in H$, and $(\gamma,a)\in \erune$ by 
\begin{align}
\label{eq:87}
(\gamma,a)\cdot \xi &:= (\gamma\xi,\alpha_\xi\inv(a)) \\
\label{eq:88}
[\eta,\zeta]\cdot (\gamma,a) &:= ([\eta,\zeta]\cdot \gamma, a)  = 
(\eta\zeta\inv\gamma, a).
\end{align}
The second equality in \eqref{eq:88} follows from the fact that
$[\eta,\zeta]\cdot \gamma = \eta\delta$ where $\delta$
is the unique element of $H$ such that
$\gamma = \zeta\delta$.  Notice that \eqref{eq:88} is well defined
because $a\in A(s(\gamma)) = A(s(\eta\zeta\inv\gamma))$.
Furthermore it is easy to show that the action of $G^H$ on $\erune$ is
continuous and respects the groupoid structure.  Now consider
\eqref{eq:87}.  It is also easy to show that this defines an action of
$H$ on $\erune$ which is continuous because $\alpha$ is continuous.
Thus it follows that $\erune$ is a strong left $G^H$-space and a strong
right $H$-space.  Finally, we recall that the actions of $G^H$ and $H$
on $X$ commute so that 
\begin{align*}
([\eta,\zeta]\cdot(\gamma,a))\cdot \xi &=
(([\eta,\zeta]\cdot \gamma)\cdot \xi,\alpha_\xi\inv(a)) \\
&= ([\eta,\zeta]\cdot(\gamma\cdot \xi),\alpha_\xi\inv(a)) = 
[\eta,\zeta]\cdot((\gamma,a)\cdot \xi).
\end{align*}
Hence the actions on $\erune$ commute. 

At this point we need to verify the equivalence conditions on
$\erune$.  The continuity condition follows in a straightforward
manner from the fact that the operations on $\mcal{A}$ are continuous
and the fact that $\alpha$ is a continuous action.  Next we need to
show that $p\suberune$ is equivariant.  However, it is clear from
\eqref{eq:87} and \eqref{eq:88} that this is the case.  The third
condition to verify is compatibility.  Suppose $[\eta,\zeta]\in G^H$,
$\xi\in H$, $(\gamma,a),(\gamma,b)\in \erune$, $c\in A(s(\gamma))$,
and $d\in A(r(\gamma))$.  Then we compute
\begin{align*}
\lset{A(r(\eta))}\langle
[\eta,\zeta]\cdot(\gamma,a),[\eta,\zeta]\cdot (\gamma,b)\rangle &=
\lset{A(r(\eta))}\langle
(\eta\zeta\inv\gamma,a),(\eta\zeta\inv\gamma,b)\rangle \\
&= \alpha_{\eta\zeta\inv\gamma}(ab^*) =
\alpha_{\eta\zeta\inv}(\alpha_\gamma(ab^*)) \\
&=
\sigma_{[\eta,\zeta]}(\lset{A(r(\gamma))}\langle(\gamma,a),(\gamma,b)\rangle),
\\
\langle (\gamma,a)\cdot \xi,(\gamma,b)\cdot \xi\rangle_{A(s(\xi))} &= 
\langle
(\gamma\xi,\alpha_\xi\inv(a)),(\gamma\xi,\alpha_\xi\inv(b))\rangle_{A(s(\xi))}
= \alpha_\xi\inv(a^*b) \\
&= \alpha_\xi\inv(\langle (\gamma,a),(\gamma,b)\rangle_{A(s(\gamma))}),
\\
[\eta,\zeta]\cdot(d\cdot(\gamma,a)) &=
[\eta,\zeta]\cdot(\gamma,\alpha_\gamma\inv(d)a) =
(\eta\zeta\inv\gamma,\alpha_\gamma\inv(d)a) \\
&=
(\eta\zeta\inv\gamma,\alpha_{\gamma\inv\zeta\eta\inv}(\alpha_{\eta\zeta\inv}(d))a)
\\ &=
(\eta\zeta\inv\gamma,\alpha_{\eta\zeta\inv\gamma}\inv(\sigma_{[\eta,\zeta]}(d))a)
\\
&= \sigma_{[\eta,\zeta]}(d)\cdot ([\eta,\zeta]\cdot(\gamma,a)), \\
((\gamma,a)\cdot c)\cdot \xi &= (\gamma\xi, \alpha_\xi\inv(ac)) =
(\gamma\xi,\alpha_\xi\inv(a)\alpha_\xi\inv(c)) \\
&= ((\gamma,a)\cdot \xi)\cdot \alpha_\xi\inv(c).
\end{align*}
This shows that the operations are compatible and all that is left is
to verify the invariance condition.  Once again, we calculate 
\begin{align*}
[\eta,\zeta]\cdot((\gamma,a)\cdot c) &= ([\eta,\zeta]\cdot\gamma,ac)
= ([\eta,\zeta]\cdot(\gamma,a))\cdot c, \\
(d\cdot(\gamma,a))\cdot \xi &= (\gamma,\alpha_\gamma\inv(d)a)\cdot \xi
= (\gamma\xi,\alpha_\xi\inv(\alpha_\gamma\inv(d)a)) \\
&= (\gamma\xi,\alpha_{\gamma\xi}\inv(d)\alpha_\xi\inv(a)) = 
d\cdot (\gamma\xi,\alpha_\xi\inv(a)) \\
&= d\cdot((\gamma,a)\cdot \xi).
\end{align*}
At this point we have shown that $\erune$ is an equivalence between
$(A(H\unit),H,\alpha|_H)$ and $(\rho^*A,G^H,\sigma)$.  We can apply Theorem
\ref{thm:renaultequiv} to conclude that $\mcal{Z}_0$ completes to the
desired imprimitivity bimodule.  What's more we can use \eqref{eq:74}
through \eqref{eq:73} to compute the bimodule operations.  First
recall that $G^H$ has a Haar system $\mu$ defined by \eqref{eq:2}.  Fix $f\in
\Gamma_c(G^H,r^*(\rho^*\mcal{A}))$, $g\in \Gamma_c(H,r^*\mcal{A})$ and $z,w\in
\Gamma_c(X,s^*\mcal{A})$.  Let $z([\eta,\gamma]\cdot\gamma) = z(\eta)
= (\eta,a)$ and recall that $[\gamma,\eta]\cdot(\eta,a) =
([\gamma,\eta]\cdot \eta,a) = (\gamma,a)$.  Furthermore $f([\gamma,\eta])\cdot
(\gamma,a) = (\gamma,\alpha_\gamma\inv(f([\gamma,\eta]))a)$.  Making the
usual identification of $a$ with $z(\eta)$ we then have 
\begin{align*}
f\cdot z(\gamma) &= \int_{G^H} f([\zeta,\eta])\cdot([\zeta,\eta]\cdot
z([\zeta,\eta]\inv\cdot \gamma)) d\mu^{\gamma\cdot H}([\zeta,\eta]) \\
&= \int_G f([\gamma,\eta])\cdot([\gamma,\eta]\cdot
z([\eta,\gamma]\cdot \gamma)) d\lambda_{s(\gamma)}(\eta) \\
&= \int_G
\alpha_{\gamma}\inv(f([\gamma,\eta]))z(\eta)d\lambda_{s(\gamma)}(\eta).
\end{align*}
Similarly we compute 
\begin{align*}
z\cdot g(\gamma) &= \int_H (z(\gamma\cdot \eta)\cdot \eta\inv)\cdot
\alpha_\eta(g(\eta\inv))d\lambda_H^{s(\gamma)}(\eta) \\
&= \int_H
\alpha_\eta(z(\gamma\eta))\alpha_\eta(g(\eta\inv))d\lambda_H^{s(\gamma)}(\eta)
\\
&= \int_H
\alpha_\eta(z(\gamma\eta)g(\eta\inv))d\lambda_H^{s(\gamma)}(\eta).
\end{align*}
Next, given $\gamma\in H$, in \eqref{eq:76} we are allowed to choose any
$\delta\in X$ such that $s(\delta) = s(\gamma)$.  However we may as well
choose $\delta = s(\gamma)$ so that 
\begin{align*}
\llangle z,&w\rrangle_{A\rtimes H}(\gamma) \\ &= 
\int_{G^H}\langle z([\zeta,\eta]\inv \cdot s(\gamma) \cdot \gamma\inv),
w([\zeta,\eta]\inv\cdot s(\gamma))\cdot
\gamma\inv\rangle_{A(r(\gamma))}d\mu^{s(\gamma)\cdot H}([\zeta,\eta])
\\
&= \int_G z([\eta,s(\gamma)]\cdot \gamma\inv)^*
\alpha_\gamma(w([\eta,s(\gamma)]\cdot s(\gamma)))
d\lambda_{s(\gamma)}(\eta) \\
&= \int_G
z(\eta\gamma\inv)^*\alpha_\gamma(w(\eta))d\lambda_{s(\gamma)}(\eta).
\end{align*}
Finally, given $[\gamma,\eta]\in G^H$ in \eqref{eq:73} we are allowed
to choose any $\delta\in X$ such that $r_X(\delta) = s([\gamma,\eta]) =
\eta\cdot H$.  Therefore we may as well choose $\delta = \eta$ so that 
\begin{align*}
\lset{\rho^*A\rtimes_\sigma G^H}\llangle z,w\rrangle([\gamma,\eta]) &= 
\int_G \lset{A(r(\gamma))}\langle z([\gamma,\eta]\cdot \eta \cdot
\xi), [\gamma,\eta]\cdot w(\eta\cdot \xi)\rangle
d\lambda^{s(\eta)}(\xi) \\
&= \int_G \lset{A(r(\gamma))}\langle z(\gamma\xi),w(\eta\xi)\rangle
d\lambda^{s(\gamma)}(\xi) \\
&= \int_G
\alpha_{\gamma\xi}(z(\gamma\xi)^*w(\eta\xi))d\lambda^{s(\gamma)}(\xi). \qedhere
\end{align*}
\end{proof}

\begin{remark}
\index{groupoid!transitive}
Suppose $G$ is a transitive groupoid and $u\in G\unit$.  It is worth
pointing out that Theorem \ref{thm:transprod} becomes a special case
of Proposition \ref{prop:80} after we identify $G$ with $G^{S_u}$ via
Proposition \ref{prop:20}.  We will actually use this fact indirectly 
in Section \ref{sec:regularity}. 
\end{remark}

\begin{remark}
\index{crossed product!group}
It is not obvious, but 
Proposition \ref{prop:80} reduces to Green's Imprimitivity Theorem
\cite[Theorem 4.22]{tfb2} when $(A,G,\alpha)$ is a group dynamical
system and $H$ is a subgroup of $G$.  We
will sketch this construction without going into detail.  
First, observe that $X = G$.  It is
not difficult to show that the imprimitivity groupoid $G^H$ is
isomorphic to the transformation groupoid $G\ltimes G/H$ 
associated to the left action of $G$ on $G/H$.  As a bundle, $A$ has a
single fibre so that $\rho^*A = C_0(G/H,A)$.  Proposition
\ref{prop:79} implies that there is an action $\sigma$ of $G\ltimes
G/H$ on $C_0(G/H,A)$.  Now, define $\bar{\sigma} = \lt\otimes \alpha$
and let $G$ act on $C_0(G/H,A)$ via $\bar{\sigma}$.  It is
straightforward to show, using Example \ref{ex:22} as a guide, that
$C_0(G/H,A)\rtimes_\sigma (G/H\ltimes G)$ is naturally isomorphic to
$C_0(G/H,A)\rtimes_\sigma G$ via the map $\phi(f)(s)([t]) = f(s,[t])$
for $f\in C_c(G,C_c(G/H,A))$.  Thus Proposition \ref{prop:80} implies
that $C_0(G/H,A)\rtimes_\sigma G$ is Morita equivalent to
$A\rtimes_\alpha H$, just as in Green's Imprimitivity Theorem.  While
the imprimitivity algebras from these theorems are not the same, they
are related.  Since $A$ has a single fibre, $\mcal{Z}_0 =
C_c(G,A)$.  Let $\mcal{X}_0 = C_c(G,A)$ be the pre-imprimitivity
bundle coming from Green's Imprimitivity Theorem and define $\psi:\mcal{Z}_0
\rightarrow \mcal{X}_0$ such that $\psi(f)(s) = \alpha_s(f(s))$.
It is not difficult, but requires some lengthy computations, to
show that $\psi$ is a bijection which preserves all of the operations
defined in Proposition \ref{prop:80} and \cite[Theorem 4.22]{tfb2}. 
\end{remark}

The whole purpose of building the imprimitivity bimodule
$\mcal{Z}_H^G$ was so that we can mimic \cite[Section 5.1]{tfb2}
and \cite{rieffelinduce} to create an induction process for representations. 

\begin{remark}
\label{rem:19}
It is assumed that the reader is familiar with the material from
\cite[Section 2.4]{tfb}.  We will be making use of the
induced representation construction described there. In particular,
given a separable groupoid dynamical system $(A,G,\alpha)$, a closed
subgroupoid $H$ with a Haar system, and a representation $\pi$ of
$A\rtimes_\alpha H$ on $\mcal{H}$ we will define the Hilbert space
$\mcal{Z}_H^G\otimes_{A\rtimes H} \mcal{H}$ to be the completion of
the vector space tensor product $\mcal{Z}_H^G\odot \mcal{H}$ with
respect to the inner product characterized by 
\begin{equation}
(z\otimes h, w\otimes k) := (\pi(\llangle w,z\rrangle_{A\rtimes
  H})h,k).
\end{equation}
\end{remark}

In any case, our ultimate goal will be to prove the following 

\begin{theorem}
\label{thm:induce}
\index{induced representation}
\index{New Result}
\index[not]{$\Ind_H^G \pi$}
Suppose $(A,G,\alpha)$ is a separable 
groupoid dynamical system and that $H$ is
a closed subgroupoid of $G$ with a Haar system.  Then given a
representation $\pi$ of $A(H\unit)\rtimes_\alpha H$ on $\mcal{H}$  
we may form the induced
representation $\Ind_H^G \pi$ of $A\rtimes_\alpha G$ on
$\mcal{Z}_H^G\otimes_{A\rtimes H}\mcal{H}$ which is defined for
$f\in\Gamma_c(G,r^*\mcal{A})$, $z\in\mcal{Z}_0$ and $h\in\mcal{H}$ by
\[
\Ind_H^G\pi(f)(z\otimes h) = f\cdot z\otimes h
\]
where 
\begin{equation}
\label{eq:ind}
f\cdot z(\gamma) = \int_G
\alpha_\gamma\inv(f(\eta))z(\eta\inv\gamma)d\lambda^{r(\gamma)}(\eta).
\end{equation}
\end{theorem}

We start by proving that we can let $A\rtimes G$ act nondegenerately
as adjointable linear operators on $\mcal{Z}_G^H$.  Actually,
considering the remarks at the end of \cite[Section 3.3]{tfb}, this
gets us most of the way there.  

\begin{prop}
\label{prop:82}
Let $(A,G,\alpha)$ be a separable groupoid dynamical system, $H$ a closed
subgroupoid of $G$ with a Haar system and $\mcal{Z}_H^G$ be the
associated imprimitivity bimodule.  Let $\lambda$ be a Haar system for
$G$, $X=s\inv(H\unit)$, and $\mcal{Z}_0 = \Gamma_c(X,s^*\mcal{A})$.
Then there is a nondegenerate homomorphism $\phi:A\rtimes_\alpha
G\rightarrow \mcal{L}(\mcal{Z}_H^G)$ such that for $f\in
\Gamma_c(G,r^*\mcal{A})$ and $z\in \mcal{Z}_0$
\begin{equation}
\label{eq:90}
\phi(f)z(\gamma) = \int_G
\alpha_\gamma\inv(f(\eta))z(\eta\inv\gamma)d\lambda^{r(\gamma)}(\eta). 
\end{equation}
\end{prop}

\begin{proof}
We will construct $\phi$ by showing that $A\rtimes G$ sits
inside the multiplier algebra of $\rho^*A\rtimes_\sigma G^H$ where $G^H$ is
the imprimitivity groupoid and $\sigma$ is the associated action of
$G^H$ on $\rho^*A$.  Given $f\in \Gamma_c(G,r^*\mcal{A})$ and $g\in
\Gamma_c(G^H,r^*(\rho^*\mcal{A}))$ define 
\begin{equation}
\label{eq:89}
M_f(g)([\gamma,\eta]) = \int_G
f(\xi)\alpha_{\xi}(g([\xi\inv\gamma,\eta]))d\lambda^{r(\gamma)}(\xi).
\end{equation}
We start by proving that $M_f(g)$ is a continuous compactly supported
section.  This argument is nearly the same as the argument that
convolution is well defined.  Consider the function
\[
\kappa(\xi,[\gamma,\eta])=
f(\xi)\alpha_{\xi}(g([\xi\inv\gamma,\eta]))
\]
on $G*G^H = \{(\xi,[\gamma,\eta])\in G\times G^H : r(\xi) =
r(\gamma)\}$.  Suppose $\xi_i\rightarrow \xi$ in $G$,
$[\gamma_i,\eta_i]\rightarrow [\gamma,\eta]$ in $G^H$, that $r(\xi_i)
= r(\gamma_i)$ for all $i$ and $r(\xi) = r(\gamma)$.  Then, after
passing to a subsequence, choosing new representatives and relabeling, we
may assume that $\gamma_i\rightarrow \gamma$ and $\eta_i\rightarrow
\eta$.  However it follows immediately that
$[\xi_i\inv\gamma_i,\eta_i]\rightarrow [\xi\inv\gamma,\eta]$.  From
here, one just observes that $\alpha$ is continuous to conclude that
$\kappa$ is a continuous function.  Furthermore, suppose the
sequence $(\xi_i,[\gamma_i,\eta_i])$ is in $\supp \kappa$.  Then we
must have $\xi_i\in \supp f$ and $[\xi_i\inv\gamma_i,\eta_i]\in \supp
g$ for all $i$.  Since both $f$ and $g$ are compactly supported we can
find subsequences which converge.  However, after passing to another
subsequence and choosing new representatives $\gamma_i$ and $\eta_i$
we can assume that there exists $\xi,\gamma,\eta\in G$ such that
$\xi_i\rightarrow \xi$, $\xi_i\inv\gamma_i\rightarrow \gamma$ and
$\eta_i\rightarrow \eta$.  However it is now clear that
$(\xi_i,[\gamma_i,\eta_i])$ has a subsequence which converges to
$(\xi,[\gamma\xi,\eta])$.  Thus $\kappa$ is a compactly supported
continuous function so that $\kappa\in \Gamma_c(G*G^H,\bar{r}^*\mcal{A})$,
where in this case $\bar{r}(\xi,[\gamma,\eta]) = r(\xi)$.   

Now, given $\kappa\in \Gamma_c(G*G^H,\bar{r}^*\mcal{A})$ we 
wish to show that the function 
\begin{equation}
\label{eq:92}
[\gamma,\eta]\mapsto \int_G
\kappa(\xi,[\gamma,\eta])d\lambda^{r(\gamma)}(\xi)
\end{equation}
is continuous and compactly supported.  As in Lemma \ref{lem:9},
if $\kappa_i\rightarrow \kappa$ with respect to the inductive limit
topology and \eqref{eq:92} is continuous for each $\kappa_i$ then it
is continuous for $\kappa$ as well.  Since sums of elementary tensors are
dense in $\Gamma_c(G*G^H,\bar{r}^*\mcal{A})$, and since sums of
continuous functions are continuous, we may as well assume that
$\kappa = h\otimes a$ for $h\in C_c(G*G^H)$ and $a\in A$.  
It is not difficult to
see that $G*G^H$ is closed in $G\times G^H$.  As a result we can
extend $h$ using Lemma \ref{lem:8} to a compactly supported
function on $G\times G^H$.  Since sums of functions of the form $k\otimes
l(\xi,[\gamma,\eta]) = k(\xi)l([\gamma,\eta])$ are dense in
$C_c(G\times G^H)$ we can, as above, assume that  
$h = k\otimes l$ for $k\in C_c(G)$ and $l\in C_c(G^H)$.  However in
this case 
\[
\int_G (k\otimes l)\otimes
a(\xi,[\gamma,\eta])d\lambda^{r(\gamma)}(\xi) = 
l([\gamma,\eta])a(r(\gamma))\int_G k(\xi)d\lambda^{r(\gamma)}(\xi)
\]
and thus for $\kappa = k\otimes l \otimes a$ 
\eqref{eq:92} is clearly a continuous and compactly supported
function on $G^H$.  It
now follows that $M_f g\in \Gamma_c(G^H,r^*(\rho^*\mcal{A}))$.  

Next, suppose that the sequence $f_i\rightarrow f$ with respect to the inductive
limit topology in $\Gamma_c(G,r^*\mcal{A})$ and $g_i\rightarrow g$
with respect to the inductive limit topology in
$\Gamma_c(G^H,r^*(\rho^*\mcal{A}))$.  Let $K$ be a compact set in $G$
eventually containing $\supp f_i$ and $L$ a compact set in $G^H$
eventually containing $\supp g_i$.  We compute
\begin{align*}
\|&M_{f_i}g_i([\gamma,\eta])-M_f g([\gamma,\eta])\| \\
&\leq 
\int_G \|f_i(\xi)\alpha_{\xi}(g_i([\xi\inv\gamma,\eta])) - 
f(\xi)\alpha_{\xi}(g([\xi\inv\gamma,\eta]))\|
d\lambda^{r(\gamma)}(\xi) \\
&\leq \int_G \|f_i(\xi)-f(\xi)\|\|g_i([\xi\inv\gamma,\eta])\| + 
\|f(\xi)\|\|g_i([\xi\inv\gamma,\eta])-g([\xi\inv\gamma,\eta])\|
d\lambda^{r(\gamma)}(\xi) \\
&\leq (\|f_i-f\|_\infty \|g_i\|_\infty +
\|f\|_\infty\|g_i-g\|_\infty)\lambda^{r(\gamma)}(K)
\end{align*}
Since $\{\|g_i\|_\infty\}$ and $\{\lambda^u(K)\}$ are bounded this shows that
$M_{f_i}g_i\rightarrow M_f g$ uniformly.  Furthermore it is
straightforward to show that $\{[\gamma\xi\inv,\eta]\in G^H : \xi\in K,
[\gamma,\eta]\in L\}$ is a compact set which eventually contains
$\supp M_{f_i}g_i$.  Thus $M$ is jointly continuous with respect to
the inductive limit topology.  

In order to prove that $M_f$ defines a multiplier we need to show that
it extends to an adjointable linear operator when we view 
$\rho^*A\rtimes_\sigma G^H$ as a
right $\rho^*A\rtimes_\sigma G^H$-module in the usual fashion.  
First, it is clear that $M_f$ is linear.  Next we show that $M_f$
preserves the left action of $\rho^*A\rtimes_\sigma G^H$ on
$\Gamma_c(G^H,r^*(\rho^*\mcal{A}))$.   Let $\mu$ be the Haar
system on $G^H$ from Proposition \ref{prop:21}.  Using the left invariance
of Haar measure we compute for $f\in \Gamma_c(G,r^*\mcal{A})$ and
$g,h\in\Gamma_c(G^H,r^*(\rho^*\mcal{A}))$ 
\begin{align*}
M_f(g*h)&([\gamma,\eta]) \\
&= \int_G\int_{G^H}
f(\xi)\alpha_\xi(g([\delta,\zeta])\sigma_{[\delta,\zeta]}(h([\delta,\zeta]\inv[\xi\inv\gamma,\eta])))
d\mu^{\xi\inv\gamma\cdot H}([\delta,\zeta])d\lambda^{r(\gamma)}(\xi)
\\
&=\int_G\int_G
f(\xi)\alpha_\xi(g([\xi\inv\gamma,\zeta])\alpha_{\xi\inv\gamma\zeta\inv}(h([\zeta,\eta])))
d\lambda_{s(\gamma)}(\zeta)d\lambda^{r(\gamma)}(\xi) \\
&= \int_G\int_G
f(\xi)\alpha_{\xi}(g([\xi\inv\gamma,\zeta]))\alpha_{\gamma\zeta\inv}(h([\zeta,\eta]))d\lambda^{r(\gamma)}(\xi)d\lambda_{s(\gamma)}(\zeta)
\\
&= \int_{G^H}
M_fg([\delta,\zeta])\sigma_{[\delta,\zeta]}(h([\delta,\zeta]\inv[\gamma,\eta]))
d\mu^{\gamma\cdot H}([\delta,\zeta]) \\
&= (M_fg)*h([\gamma,\eta]).
\end{align*}
Next, we show that $M$ is
adjointable on $\Gamma_c(G^H,r^*(\rho^*\mcal{A}))$ with adjoint $M_{f^*}$ by computing
\begin{align*}
(M_f g)^*&*h([\gamma,\eta]) = \int_{G^H}
(M_fg)^*([\zeta,\xi])\sigma_{[\zeta,\xi]}(h([\zeta,\xi]\inv[\gamma,\eta]))d\mu^{\gamma\cdot
  H}([\zeta,\xi]) \\
&= \int_G \sigma_{[\gamma,\xi]}((M_f g([\xi,\gamma]))^*h([\xi,\eta]))
d\lambda_{s(\gamma)}(\xi) \\
&= \int_G\int_G
\alpha_{\gamma\xi\inv}(\alpha_\zeta(g([\zeta\inv\xi,\gamma]))^*f(\zeta)^*h([\xi,\eta]))d\lambda^{r(\xi)}(\zeta)d\lambda_{s(\gamma)}(\xi)
\\
&= \int_G\int_G
\alpha_{\gamma\xi\zeta}(g([\zeta\inv\xi\inv,\gamma])^*)\alpha_{\gamma\xi}(f(\zeta)^*
h([\xi\inv,\eta]))d\lambda^{s(\xi)}(\zeta)d\lambda^{s(\gamma)}(\xi) \\
&= \int_G \int_G
\alpha_{\gamma\zeta}(g([\zeta\inv,\gamma])^*)\alpha_{\gamma\xi}(f(\xi\inv\zeta)^*h([\xi\inv,\eta]))d\lambda^{s(\gamma)}(\zeta)d\lambda^{s(\gamma)}(\xi)
\\
&= \int_G \int_G
\alpha_{\gamma\zeta}(g([\zeta\inv,\gamma])^*\alpha_{\zeta\inv\xi}(f(\xi\inv\zeta)^*h([\xi\inv,\eta])))d\lambda^{s(\gamma)}(\xi)d\lambda^{s(\gamma)}(\zeta)
\\
&= \int_G \int_G \alpha_{\gamma\zeta}(g([\zeta\inv,\gamma])^*
\alpha_\xi(f(\xi\inv)^*)\alpha_\xi(h([\xi\inv\zeta\inv,\eta])))d\lambda^{s(\zeta)}(\xi)d\lambda^{s(\gamma)}(\zeta)
\\
&= \int_G \int_G
\alpha_{\gamma\zeta\inv}(g([\zeta,\gamma])^*f^*(\xi)\alpha_\xi(h([\xi\inv\zeta,\eta])))d\lambda^{r(\zeta)}(\xi)d\lambda_{s(\gamma)}(\zeta)
\\
&= \int_G \sigma_{[\gamma,\zeta]}(g([\zeta,\gamma])^*
M_{f^*}(h)([\zeta,\eta])) d\lambda_{s(\gamma)}(\zeta) \\
&= \int_{G^H} \sigma_{[\xi,\zeta]}(g([\xi,\zeta]\inv)^*
M_{f^*}(h)([\xi,\zeta]\inv[\gamma,\eta])) d\mu^{\gamma\cdot
  H}([\xi,\zeta]) \\
&= g^* * (M_{f^*}h)([\gamma,\eta]).
\end{align*}
Finally, we prove that $M$ preserves convolution on
$\Gamma_c(G,r^*\mcal{A})$ by computing for
$f,g\in\Gamma_c(G,r^*\mcal{A})$ and $h\in\Gamma_c(G^H,r^*(\rho^*\mcal{A}))$
that 
\begin{align*}
M_{f*g}h([\gamma,\eta]) &= 
\int_G\int_G
f(\delta)\alpha_\delta(g(\delta\inv\xi))\alpha_{\xi}(h([\xi\inv\gamma,\eta]))d\lambda^{r(\gamma)}(\delta)d\lambda^{r(\gamma)}(\xi)
\\
&= \int_G \int_G
f(\delta)\alpha_\delta(g(\delta\inv\xi)\alpha_{\delta\inv\xi}(h([\xi\inv\gamma,\eta]))d\lambda^{r(\gamma)}(\xi)d\lambda^{r(\gamma)}(\delta)
\\
&= \int_G \int_G f(\delta)
\alpha_\delta(g(\xi)\alpha_{\xi}(h([\xi\inv\delta\inv\gamma,\eta])))
d\lambda^{s(\delta)}(\xi)d\lambda^{r(\gamma)}(\delta) \\
&= \int_G
f(\delta)\alpha_\delta(M_gh([\delta\inv\gamma,\eta]))d\lambda^{r(\gamma)}(\delta)
\\
&= M_fM_gh([\gamma,\eta]).
\end{align*}

Next we show that $M$ is nondegenerate.  Specifically we will show
that elements of the form $M_f (g)$ are dense in
$\Gamma_c(G^H,r^*(\rho^*\mcal{A}))$ with respect to the inductive 
limit topology.  Let $a_l$ be an approximate identity for
$A$ and let $e_{(K,U,l,\epsilon)}$ be the approximate identity coming
from Lemma \ref{lem:19}.   We will use $\kappa$ to denote a generic
4-tuple $(K,U,l,\epsilon)$.  We would like to show that given
$g\in\Gamma_c(G^H,r^*(\rho^*\mcal{A}))$ we have $M_{e_\kappa}g\rightarrow g$ with
respect to the inductive limit topology.  Fix $\epsilon_1> 0$ and
let $L= \supp g$.  We make the following claim.

\begin{claim}
There exists a conditionally compact
neighborhood $U_1$ such that $\xi\in U_1$ implies 
\begin{equation}
\label{eq:91}
\|\alpha_\xi(g([\xi\inv\gamma,\eta])) - g([\gamma,\eta])\| <
\epsilon_1
\end{equation}
for all $[\gamma,\eta]\in G^H$ such that $r(\gamma) = r(\xi)$.
\end{claim}

\begin{proof}[Proof of Claim.] 
Suppose the claim does not hold.  
Fix some conditionally compact neighborhood $W$.  Then for any
conditionally compact neighborhood $U\subset W$ there exists $\xi_U\in
U$ and $[\gamma_U,\eta_U]\in G^H$ such that 
\[
\|\alpha_{\xi_U}(g([\xi_U\inv\gamma_U,\eta_U]))-g([\gamma_U,\eta_U])\|\geq
\epsilon_1.
\]
However, for this to hold one of the terms must be nonzero so that we
must have, recalling that $W$ is a neighborhood of $G\unit$ which
contains $U$,  
\begin{equation}
\label{eq:96}
[\gamma_U,\eta_U]\in \widetilde{L}=\{[\xi\gamma,\eta]:\xi\in W,
[\gamma,\eta]\in L\ \text{and}\ s(\xi) = r(\gamma)\}.
\end{equation}
Suppose $\{[\xi_i\inv\gamma_i,\eta_i]\}$ is contained in
$\widetilde{L}$.  Then, by passing to a subsequence, relabeling and choosing new
representatives we can use the fact that $L$ is compact to find
$\gamma,\eta\in G$ such that $\gamma_i\rightarrow \gamma$ and
$\eta_i\rightarrow \eta$.  Since $\rho$ is continuous it follows that 
$K_1=\rho(r(L))$ is
compact and by assumption contains $\{s(\xi_i)\}$.  Since $W$ is
conditionally compact the set $W\cap s\inv(K_1)$ is compact
and also contains $\{\xi_i\}$.  Now we can pass to
another subnet and find $\xi$ such that $\xi_i\rightarrow \xi$.  It
follows immediately that $\widetilde{L}$ is compact.  Thus, ordering
$\{[\gamma_U,\eta_U]\}$ by decreasing $U$, we can pass to a subnet
(twice actually),
relabel and find new representatives such that there exists
$\gamma,\eta\in G$ with $\gamma_U\rightarrow \gamma$ and
$\eta_U\rightarrow \eta$. 
Next, observe that we have $r(\xi_U)\in \rho(r(\widetilde{L}))$ for all $U$ so
that $\{\xi_U\}$ is contained in the compact set $U\cap
r\inv(\rho(r(\widetilde{L})))$.  
Therefore we can pass to yet another subnet and find
$\xi$ such that $\xi_U\rightarrow \xi$.  However, by construction,
$\xi\in U$ for any conditionally compact neighborhood of $G\unit$.  It
follows that $\xi\in G\unit$.  Hence 
\[
\alpha_{\xi_U}(g([\xi_U\inv\gamma_U,\eta_U]))\rightarrow
g([\gamma,\eta]).
\]
However, it follows that eventually 
\[
\|\alpha_{\xi_U}(g([\xi_U\inv\gamma_U,\eta_U]))-g([\gamma_U,\eta_U])\|
< \epsilon_1 
\]
which is a contradiction. 
\end{proof}

Unfortunately we need another claim before we can tackle
nondegeneracy.  

\begin{claim}
There exists $l_1$ such that $l\geq l_1$ implies 
\begin{equation}
\label{eq:93}
\|a_l(r(\gamma))g([\gamma,\eta]) - g([\gamma,\eta])\| < \epsilon_1 \quad\text{for all $[\gamma,\eta]\in G^H$.}
\end{equation}
\end{claim}

\begin{proof}[Proof of Claim.]
It clearly suffices to verify this identity on $L$.  
Since $a_l$ factors to an approximate identity on each fibre we have
$a_l(r(\gamma))g([\gamma,\eta])\rightarrow g([\gamma,\eta])$ for each
$[\gamma,\eta]\in G^H$.  We use the fact
that the norm is upper-semicontinuous to choose for each
$[\gamma,\eta]\in L$ some neighborhood $O_{[\gamma,\eta]}$ of
$[\gamma,\eta]$ and some $b_{[\gamma,\eta]}\in \{a_l\}$ such that 
\[
\|b_{[\gamma,\eta]}(r(\xi))g([\xi,\zeta]) - g([\xi,\zeta])\| <
\frac{\epsilon_1}{3}
\]
for all $[\xi,\zeta]\in O_{[\gamma,\eta]}$.  Since $L$ is compact we
can find some finite subcover $\{O_i\}$.  Let $\phi_i\in C_c(G^H)$
be a partition of unity with respect to $\{O_i\}_{i=1}^N$ so that
$\supp\phi_i\subset O_i$ and $\sum \phi_i([\gamma,\eta]) = 1$ if
$[\gamma,\eta]\in L$.  Define $h \in \Gamma_c(G^H,r^*(\rho^*\mcal{A}))$ by 
\[
h = \sum_{i=1}^N \phi_i \otimes b_{[\gamma_i,\eta_i]}.
\]
Then by construction, for all $[\xi,\zeta]\in L$ we have 
\begin{align}
\nonumber
\|h([\xi,\zeta])g([\xi,\zeta]) - g([\xi,\zeta])\| 
&\leq \sum_{i=1}^\infty
\phi_i([\xi,\zeta])\|b_{[\gamma_i,\eta_i]}(r(\xi))g([\xi,\zeta])-
  g([\xi,\zeta])\| \\
\label{eq:94}
&< \frac{\epsilon_1}{3}.
\end{align}
Moving on, we can find $l_1$ such that if $l \geq l_1$ then 
\[
\|a_lb_{[\gamma_i,\eta_i]} - b_{[\gamma_i,\eta_i]}\| < \frac{\epsilon_1}{3\|g\|_\infty}
\]
for all $i$.  Using the fact
that passing to fibres is norm contractive this implies that
for $[\xi,\zeta]\in L$ we have
\begin{align}
\nonumber
\|a_l(r(\xi))h([\xi,\zeta]) - h([\xi,\zeta])\| &\leq \sum_{i=1}^N
\phi_i([\zeta,\xi])\|a_l(r(\xi))b_{[\gamma_i,\eta_i]}(r(\xi)) -
b_{[\gamma_i,\eta_i]}(r(\xi))\| \\
\label{eq:95}
&< \frac{\epsilon_1}{3\|g\|_\infty}.
\end{align}
Therefore, using \eqref{eq:94} and \eqref{eq:95} and the fact that
$\|a_l\|\leq 1$ we compute for $l\geq l_1$ and $[\xi,\zeta]\in L$
\begin{align*}
\|a_l(r(\xi))g([\xi,\zeta]) - g([\xi,\zeta])\| \leq&
\|a_l(r(\xi))(g([\xi,\zeta]) - h([\xi,\zeta])g([\xi,\zeta]))\|\\
&+\|(a_l(r(\xi))h([\xi,\zeta]) - h([\xi,\zeta]))g([\xi,\zeta])\|\\ 
&+\|h([\xi,\zeta])g([\xi,\zeta]) - g([\xi,\zeta])\| \\
<& \|a_l(r(\xi))\|\frac{\epsilon_1}{3} + \|g\|_\infty
\frac{\epsilon_1}{3\|g\|_\infty} + \frac{\epsilon_1}{3} \leq
\epsilon_1. \qedhere
\end{align*}
\end{proof}

Now suppose we are given $\epsilon_0 >0$ and let $\epsilon_1 =
\epsilon_0/(5+\|g\|_\infty)$.  Choose $U_1$ and $l_1$ for $\epsilon_1$
as above.  Then given $e =
e_{(K,U,l,\epsilon)}$ with $K_1=\rho(r(L))\subset K$, $U\subset U_1$, $l_1\leq 1$
and $\epsilon < \epsilon_1$ we compute for
$[\gamma,\eta]\in G^H$
\begin{align*}
\|M_e g([\gamma,\eta]) - g([\gamma,\eta])\| \leq& \left\|\int_G
  e(\xi)(\alpha_\xi(g([\xi\inv\gamma,\eta]))-g([\gamma,\eta]))d\lambda^{r(\gamma)}(\xi)\right\|
\\
&+ \left\|\left(\int_G
e(\xi)d\lambda^{r(\gamma)}(\xi)-a_l(r(\gamma))\right)g([\gamma,\eta])\right\|\\ 
&+ \|a_l(r(\gamma))g([\gamma,\eta]) - g([\gamma,\eta])\| \\
<& \int_U
\|e(\xi)\|\|\alpha_\xi(g([\xi\inv\gamma,\eta]))-g([\gamma,\eta])\|d\lambda^{r(\gamma)}(\xi)
\\
&+ \epsilon\|g([\gamma,\eta])\| + \epsilon_1 \\
<& 4\epsilon_1 + \epsilon_1\|g\|_\infty + \epsilon_1 = \epsilon_0.
\end{align*}
Thus $M_{e_\kappa} g\rightarrow g$ uniformly.  Furthermore, if
$\kappa =(K,U,l,\epsilon)$ such that $U\subset U_1$ 
then, considering the fact that $\supp
e_{\kappa}\subset U$ we have $M_{e_\kappa}g([\gamma,\eta]) \ne 0$
only if 
\[
[\gamma,\eta]\in \widetilde{L} = \{[\xi\gamma,\eta]:\xi\in U_1, [\gamma,\eta]\in
L,\ \text{and}\ s(\xi) = r(\gamma)\}
\]
However this set has the same form as \eqref{eq:96} and we proved that
$\widetilde{L}$ was compact there.  Thus we eventually have 
$\supp M_{e_{\kappa}} g\subset \widetilde{L}$ so that $M_{e_\kappa}g \rightarrow
g$ with respect to the inductive limit topology.  This of course
implies that elements of the form $M_f g$ are dense in
$\Gamma_c(G^H,r^*(\rho^*\mcal{A}))$ with respect to the inductive
limit topology.  

This nondegeneracy argument was a bear but it will be crucial in what
follows.  We would like to show that $M_f$ extends to an adjointable
operator on $\rho^*A\rtimes_\sigma G^H$ and that $\|M_f\|\leq \|f\|$ so that
$M$ extends to $A\rtimes G$.  Well, suppose $\tau$ is a state on
$\rho^*A\rtimes_\sigma G^H$ and, following the usual GNS construction, define a
pre-inner product on $\rho^*A\rtimes_\sigma G^H$ via 
\[
(g,h)_\tau := \tau(g^**h).
\]
Let $\mcal{H}_\tau$ denote the resulting Hilbert space and
$\mcal{H}_0$ the image of $\Gamma_c(G^H,r^*(\rho^*\mcal{A}))$ in
$\mcal{H}_\tau$.  Observe that $\mcal{H}_0$ is a dense subspace.  Now
given $f\in \Gamma_c(G,r^*\mcal{A})$ we would like to define an
operator $\pi(f)$ on $\mcal{H}_0$ by 
\[
\pi(f)g = M_f g.
\]
Of course, we need to see that this factors correctly.  Suppose $g\in
\Gamma_c(G^H,r^*(\rho^*\mcal{A}))$  is such that $(g,h)_\tau = 0$ for all
$h\in \Gamma_c(G^H,r^*(\rho^*\mcal{A}))$.  Then 
\[
(\pi(f)g,h)_\tau = \tau((M_f g)^**h) = \tau(g^**(M_{f^*}h)) =
(g,\pi(f^*)h)_\tau = 0
\]
for all $h\in \Gamma_c(G^H,r^*(\rho^*\mcal{A}))$ so that $\pi(f)g = 0$.  Thus
$\pi$ is well defined and obviously defines a linear operator on
$\mcal{H}_0$.  Furthermore, it follows from the fact that $M$ is
linear in $f$ and that $M_{f*g} =M_f M_g$ that $\pi$ is a homomorphism
from $\Gamma_c(G,r^*\mcal{A})$ into the algebra of linear operators on
$\mcal{H}_0$.  We will now show that we can apply Theorem
\ref{thm:disintigration}.  Since elements of the form $M_f(g)$ are
dense in the inductive limit topology it is clear that elements of the
form $\pi(f)h$ are dense in $\mcal{H}_\tau$.  Furthermore, we have
shown that $M_f g$ is jointly continuous with respect to the inductive
limit topology.  Therefore, if we fix
$g,h\in\Gamma_c(G^H,r^*(\rho^*\mcal{A}))$ and let $f_i\rightarrow f$ with
respect to the inductive limit topology we must have 
\[
(\pi(f_i)g,h)_\tau = \tau((M_{f_i}g)^**h)\rightarrow \tau((M_fg)^**h)
= (\pi(f)g,h)_\tau.
\]
Finally it is clear from the fact that $M_f$ is adjointable on
$\Gamma_c(G^H, r^*(\rho^*\mcal{A}))$ with adjoint $M_{f^*}$ that 
\[
(\pi(f)g,h)_\tau = \tau((M_f g)^**h) = \tau(g^**(M_{f^*}h)) =
(g,\pi(f^*)h)_\tau.
\]
Thus we may apply Theorem \ref{thm:disintigration} and conclude that
$\pi$ extends to representation of $A\rtimes G$.  In particular this
implies that for $f\in\Gamma_c(G,r^*\mcal{A})$ and
$g\in\Gamma_c(G^H,r^*(\rho^*\mcal{A}))$ we have 
\begin{align*}
\tau((M_f g)^* *(M_f g)) &= \|\pi(f)g\|^2_\tau \leq \|f\|^2\|g\|^2_\tau = 
\|f\|^2 \tau(g^* * g) \\
& \leq \|f\|^2\|g^* * g\| = \|f\|^2\|g\|^2.
\end{align*}
However, $\tau$ is an arbitrary state on $\rho^*A\rtimes_\sigma G^H$ so that
by choosing $\tau$  such that $\tau((M_f g)^* * (M_f g)) = \|M_f
g\|^2$ \cite[Lemma A.3]{tfb} we must have 
\begin{equation}
\label{eq:98}
\|M_f g\| \leq \|f\| \|g\|.
\end{equation}
Therefore $M_f$ is bounded and as such extends to a linear operator on
$\rho^*A\rtimes_\sigma G^H$.  However, $M_f$ is $\rho^*A\rtimes_\sigma G^H$-linear
and adjointable on a dense subspace so that this must be true in
general.  It follows that $M_f$ defines a multiplier on
$\rho^*A\rtimes_\sigma G^H$.  Furthermore it is clear that $f\mapsto M_f$ is
linear and we have already verified that it preserves convolution on
$\Gamma_c(G,r^*\mcal{A})$.  In addition we calculated that the adjoint
of $M_f$ is $M_{f^*}$ so that $M$ preserves involution as well.  Since 
\eqref{eq:98} implies that $\|M_f\|\leq  \|f\|$ it follows that $M$
extends to a $*$-homomorphism from $A\rtimes G$ into
$M(\rho^*A\rtimes_\sigma G^H)$.  

At this point we are essentially done since the space of multipliers
on $\rho^*A\rtimes_\sigma G^H$ can be identified with
$\mcal{L}(\mcal{Z}_H^G)$.  Specifically let $\mcal{Z}_H^G$ be the
imprimitivity bimodule associated to $G$ and $H$.  Recall that since
$\mcal{Z}_H^G$ is a Hilbert $\rho^*A\rtimes_\sigma G^H$-module every
element of $\mcal{Z}$ is of the form $g\cdot z$ for $g\in
\rho^*A\rtimes_\sigma G^H$ and $z\in Z$ \cite[Proposition 2.31]{tfb}.
As such, given $f\in A\rtimes G$ we can define $\phi(f)$ on $\mcal{Z}$
by setting
\[
\phi(f)(g\cdot z) = M_f g\cdot z
\]
whenever $g\in \rho^*A\rtimes_\sigma G^H$ and $z\in \mcal{Z}_H^G$.  
It is clear that $\phi(f)$ defines  linear operator on $Z$.
Next, we compute for $f\in
A\rtimes G$, $g,h\in \rho^*A\rtimes_\sigma G^H$ and $z,w\in \mcal{Z}$ that 
\begin{align*}
\llangle \phi(f)g\cdot z, h\cdot w\rrangle_{A\rtimes H} &= 
\llangle g\cdot z, (M_f g)^* h \cdot w\rrangle_{A\rtimes H} = 
\llangle z, g^* M_{f^*} h\cdot w\rrangle_{A\rtimes H} \\
&= \llangle g\cdot z, M_{f^*} h\cdot w\rrangle_{A\rtimes H} = 
\llangle g\cdot z, \phi(f^*)h\cdot w\rrangle_{A\rtimes H}
\end{align*}
Thus $\phi(f)$ is adjointable with adjoint $\phi(f^*)$ and it follows
that $\phi:A\rtimes G\rightarrow \mcal{L}(\mcal{Z})$.  Furthermore
$\phi$ preserves involution and it is easy to see that $\phi$ is
linear.  We now calculate 
\begin{align*}
\phi(f*g)h\cdot z = M_{f*g} h\cdot z = M_f M_g h\cdot z =
\phi(f)\phi(g)h\cdot z.
\end{align*}
Thus $\phi$ is a $*$-homomorphism.  Finally, since elements of the
form $M_f g$ are dense in $\rho^*A\rtimes_\sigma G^H$,
it follows that elements of
the form $\phi(f)g\cdot z$ are dense in $\mcal{Z}$.  Hence $\phi$ is
nondegenerate.  To finish we calculate for $f\in
\Gamma_c(G,r^*\mcal{A})$, $g\in\Gamma_c(G^H,r^*(\rho^*\mcal{A}))$ and
$z\in\mcal{Z}_0$ 
\begin{align*}
\phi(f)g\cdot z(\gamma) &= \int_G \alpha_\gamma\inv(M_f
g([\gamma,\eta]))z(\eta) d\lambda_{s(\gamma)}(\eta) \\
&= \int_G\int_G
\alpha_\gamma\inv(f(\xi)\alpha_\xi(g[\xi\inv\gamma,\eta])) z(\eta)
d\lambda^{r(\gamma)}(\xi)d\lambda_{s(\gamma)}(\eta) \\
&= \int_G \alpha_\gamma\inv(f(\xi))\int_G
\alpha\inv_{\xi\inv\gamma}(g([\xi\inv\gamma,\eta]))z(\eta)d\lambda_{s(\gamma)}(\eta)d\lambda^{r(\gamma)}(\xi)
\\
&= \int_G \alpha_\gamma\inv(f(\xi)) g\cdot
z(\xi\inv\gamma)d\lambda^{r(\gamma)}(\xi).
\end{align*}
Thus $\phi$ is given by \eqref{eq:90} on elements of the form $g\cdot
z$ when $g\in \Gamma_c(G^H,r^*(\rho^*\mcal{A}))$ and $z\in\mcal{Z}_0$.
However, it follows from Lemma \ref{lem:19} that there is a net in 
$\Gamma_c(G^H,r^*(\rho^*\mcal{A}))$ which is an approximate identity in
the inductive limit topology with respect to the left action of
$\Gamma_c(G^H,r^*(\rho^*\mcal{A}))$ on $\mcal{Z}_0$.  This implies that
elements of the form $g\cdot z$ with $g\in\Gamma_c(G^H,r^*(\rho^*\mcal{A}))$
and $z\in\mcal{Z}_0$ are dense in $\mcal{Z}_0$ with respect
to the inductive limit topology.  Fix $f\in \Gamma_c(G,r^*\mcal{A})$
and given $z\in\mcal{Z}_0$ define 
$\hat{z}(\xi,\gamma) = \alpha_\gamma\inv(f(\xi))z(\xi\inv\gamma)$ and
$\bar{z}(\gamma) = \int_G \hat{z}(\xi,\gamma) d\lambda^{r(\gamma)}(\xi)$.  If
$z_i\rightarrow z$ with respect to the inductive limit topology then
it is straightforward to show $\hat{z}_i\rightarrow \hat{z}$ with
respect to the inductive limit topology and in turn that $\bar{z}_i
\rightarrow \bar{z}$ with respect to the inductive limit topology.  
However, it now follows
in a straightforward fashion that \eqref{eq:90}
holds in general on $\mcal{Z}_0$.   
\end{proof}

The upshot of all of this is that we can now prove the desired
induction theorem. 

\begin{proof}[Proof of Theorem \ref{thm:induce}]
Using $\phi$ from Proposition \ref{prop:82} we can let $A\rtimes G$
act nondegenerately as adjointable operators on the Hilbert
$A\rtimes_\alpha H$-module $\mcal{Z}_H^G$.  However \cite[Proposition
2.66]{tfb} then implies that there is a nondegenerate induced 
representation $\mcal{Z}-\Ind(\pi)$ of $A\rtimes G$ acting on elementary
tensors in $\mcal{Z}\otimes_{A\rtimes H}\mcal{H}$ by 
\[
\mcal{Z}-\Ind (\pi)(f)(z\otimes h) = \phi(f)z\otimes h.
\]
Thus we can define $\Ind_H^G \pi = \mcal{Z}-\Ind(\pi)$ and
\eqref{eq:90} shows that $\Ind_H^G\pi$ has the desired action on
$\mcal{Z}\otimes_{A\rtimes H}\mcal{H}$.  
\end{proof}

\begin{example}
\label{ex:32}
\index{stabilizer subgroup}
Suppose $(A,G,\alpha)$ is a separable dynamical system.  Then for any
$u\in G\unit$ the stabilizer subgroup $S_u$ is a closed subgroup with
a Haar system.  Therefore we can induce representations from
$A(u)\rtimes S_u$ to $A\rtimes G$.  
\end{example}

\begin{example}
\index{left regular representation!induced}
Suppose $(A,G,\alpha)$ is a separable dynamical system and $\pi$ is a
representation of $A$ on $\mcal{H}$.  Without loss of generality
assume that there exists a Borel Hilbert bundle $G\unit*\mfrk{H}$ and
measure $\mu$ so that $\mcal{H} = L^2(G\unit*\mfrk{H},\mu)$ and $\pi$
decomposes as $\pi = \int_{G\unit}^\oplus \pi_u d\mu(u)$.  Consider
that $G\unit$ is a closed subgroupoid of $G$ with Haar system given by
the $\delta$ measures.  We can then form the induced representation
$\Ind_{G\unit}^G \pi$.  Observe that in this situation $A\rtimes
G\unit$ is the completion of $\Gamma_c(G\unit,\mcal{A})$ and that the
$I$-norm on $A\rtimes G\unit$ is just the uniform norm.  Since the
uniform norm is a $C^*$-norm, the enveloping algebra is just the
uniform norm 
completion of $\Gamma_c(G\unit,\mcal{A})$.  In other words $A\rtimes
G\unit = A$.  The right hand 
operations on $\mcal{Z}_{G\unit}^G$
then simplify to 
\begin{align*}
z\cdot g(\gamma) &= z(\gamma)g(s(\gamma)) \\
\llangle z,w\rrangle_A(u) &= \int_G z(\eta)^*w(\eta)d\lambda_{s(\gamma)}(\eta).
\end{align*}
Thus $\mcal{Z}\otimes_{A} \mcal{H}$ is the completion of $\mcal{Z}_0 =
\Gamma_c(G,s^*\mcal{A})$ with respect to the inner product 
\begin{align*}
(f\otimes h, g\otimes k) &:= \int_G (\pi(g(\eta)^*f(\eta))h,k)
d\lambda_{s(\gamma)}(\eta)  \\
&= \int_G (\pi_{s(\eta)}(f(\eta)^*f(\eta))h(s(\eta)),k(s(\eta))) d\nu\inv(\eta)
\\ &= \int_G (\pi_{s(\eta)}(f(\eta))h(s(\eta)),\pi_{s(\eta)}(g(\eta))h(s(\eta)))d\nu\inv(\eta).
\end{align*}
where $\nu\inv = \int_{G\unit} \lambda_u d\mu(u)$.  Furthermore the
action of $\Ind \pi$ on $\mcal{Z}\otimes_{A} \mcal{H}$ is determined
by 
\[
\Ind \pi(f)(g\otimes h) = f\cdot g \otimes h
\]
where 
\[
f\cdot g(\gamma) = \int_G
\alpha_\gamma\inv(f(\eta))g(\eta\inv\gamma)d\lambda^{r(\gamma)}(\eta). 
\]
Now recall that in Example \ref{ex:20}
we defined the left regular representation $L$ of $\pi$ to act on
$L^2(s^*(G\unit*\mfrk{H}),\nu\inv)$ via
\[
L(f)h(\gamma) = \int_G \pi_{s(\gamma)}(\alpha_\gamma\inv(f(\eta)))
h(\eta\inv\gamma)d\lambda^{r(\gamma)}(\eta).
\]
Without going into the details, it is (more or less) 
straightforward to show that the map
$U:\mcal{Z}_0\odot \mcal{H}\rightarrow
\mcal{L}^2(s^*(G\unit*\mfrk{H}),\mu)$ characterized by
\[
U(f\otimes h)(\gamma) = \pi_{s(\gamma)}(f(\gamma))h(s(\gamma))
\]
extends to a unitary map $U:\mcal{Z}\otimes_A \mcal{H}\rightarrow
L^2(s^*(G\unit*\mfrk{H}),\nu\inv)$.  Furthermore 
\begin{align*}
L(f)U(g\otimes h)(\gamma) &= \int_G
\pi_{s(\gamma)}(\alpha_\gamma\inv(f(\eta))) U(g\otimes
h)(\eta\inv\gamma) d\lambda^{r(\gamma)}(\eta) \\
&= \int_G
\pi_{s(\gamma)}(\alpha_\gamma\inv(f(\eta))g(\eta\inv\gamma))h(s(\gamma))d\lambda^{r(\gamma)}(\eta)
\\
&= \pi_{s(\gamma)}(f\cdot g(\gamma))h(s(\gamma)) = U(f\cdot g\otimes
h)(\gamma) \\
&= U\Ind\pi(f)(g\otimes h)(\gamma).
\end{align*}
Thus the left regular representation associated to $\pi$ is equivalent
to the induced representation $\Ind_{G\unit}^G \pi$.
\end{example}


\section{Stabilizers and $T_0$ Orbits}
\label{sec:regularity}

As we noted in Example \ref{ex:32}, given a groupoid crossed product
$A\rtimes G$ and $u\in G\unit$ we can induce representations of
$A(u)\rtimes S_u$ to $A\rtimes G$.  Since $A(u)\rtimes S_u$ is a group
crossed product its representation theory is relatively well
understood.  In this section we will find conditions so that every irreducible
representation of $A\rtimes G$ can be obtained in this fashion.  In
particular we will consider the ``nice'' groupoids for which the
conditions of the Mackey-Glimm dichotomy hold.  

\begin{remark}
Recall that $G$ acts on $G\unit$ by left translation.  We denote the
image of $u$ in $G\unit/G$ by $G\cdot u$.   However this notation is
also used for the orbit of $u$ in $G\unit$.  We will regularly confuse
the two and place the burden of deciding which interpretation to
use on the reader. 
\end{remark}

The key to this section will be to reduce to the case where the orbit
space is Hausdorff, because in this case we get the following result.  

\begin{prop}
\label{prop:92}
\index{upper-semicontinuous!cstar-bundle@$C^*$-bundle}
Suppose $(A,G,\alpha)$ is a separable groupoid dynamical system and
$G\unit/G$ is Hausdorff.  Then $A\rtimes_\alpha G$ is a
$C_0(G\unit/G)$-algebra with the action $\Phi$ defined for $\phi\in
C_0(G\unit/G)$ and $f\in \Gamma_c(G,r^*\mcal{A})$ by 
\[
\Phi(\phi)f(\gamma) = \phi(G\cdot r(\gamma))f(\gamma).
\]
Furthermore, restriction factors to an isomorphism of $A\rtimes
G(G\cdot u)$ onto the fibre ${A(G\cdot u)\rtimes G|_{G\cdot u}}$. 
\end{prop}

\begin{proof}
First recall that $G\unit/G$ is always locally compact so that in this
case $G\unit/G$ is a second countable locally compact Hausdorff
space.  Suppose $\Phi$ is defined as above.  It is clear that $\Phi$ is
at least linear in $\phi$ and $f$.  Furthermore we check for
$f,g\in\Gamma_c(G,r^*\mcal{A})$ that 
\begin{align}
\label{eq:113}
(\Phi(\phi)f)^* * g(\gamma) &= \int_G
\alpha_\eta((\Phi(\phi)f(\eta\inv))^*
g(\eta\inv\gamma))d\lambda^{r(\gamma)}(\eta) \\ \nonumber
&= \int_G \alpha_\eta(f(\eta\inv)^*\overline{\phi(G\cdot
  s(\eta))}g(\eta\inv\gamma)) d\lambda^{r(\gamma)}(\eta) \\ \nonumber
&= \int_G
\alpha_\eta(f(\eta\inv)^*\Phi(\overline{\phi})g(\eta\inv\gamma))d\lambda^{r(\gamma)}(\eta)
\\ \nonumber
&= f^**(\Phi(\overline{\phi})g)(\gamma).
\end{align}
This shows that $\Phi(\phi)$ is adjointable on
$\Gamma_c(G,r^*\mcal{A})$.  Next we check that 
\begin{align}
\label{eq:114}
\Phi(\phi)(f*g)(\gamma) &= \phi(G\cdot r(\gamma)) \int_G
f(\eta)\alpha_\eta(g(\eta\inv\gamma))d\lambda^{r(\gamma)}(\eta) \\ \nonumber
&= \int_G 
\Phi(\phi)f(\eta)\alpha_\eta(g(\eta\inv\gamma))d\lambda^{r(\gamma)}(\eta)
\\ \nonumber
&= (\Phi(\phi)f)*g.
\end{align}
This shows that $\Phi(\phi)$ is linear with respect to the left action
of $\Gamma_c(G,r^*\mcal{A})$ on itself.  Now let $C_0(G\unit/G)^1$ be the
unitization of $C_0(G\unit/G)$ and extend $\Phi$ to $C_0(G\unit/G)^1$ by
setting $\Phi(\phi + \lambda 1)f = \Phi(\phi)f +\lambda f$.  It is a
simple matter to show that \eqref{eq:113} and \eqref{eq:114} extend to
$C_0(G\unit/G)^1$.  Let
$\langle f,g \rangle = f^* * g$ be the usual inner product on
$A\rtimes G$ as an $A\rtimes G$-module.  We would like to show that
$\|\Phi(\phi)f\| \leq \|\phi\|_\infty \|f\|$ for all 
$f\in\Gamma_c(G,r^*\mcal{A})$. It will suffice to show that 
\[
\|\phi\|_\infty^2 \langle
f,f\rangle - \langle \Phi(\phi)f,\Phi(\phi)f\rangle \geq 0
\]
as elements of $A\rtimes G$.  However, using the fact that $\Phi$ is
adjointable on $C_0(G\unit/G)^1$, this amounts to showing
\begin{equation}
\label{eq:115}
\langle \Phi(\|\phi\|_\infty^2 1 - \overline{\phi} \phi)f, f\rangle \geq 0.
\end{equation}
All elements of the form $\|\phi\|_\infty^2 1 - \overline{\phi} \phi$ are
positive in $C_0(G\unit/G)^1$ so there exists $\xi\in C_0(G\unit/G)^1$
such that $\|\phi\|_\infty^2 1 - \overline{\phi}\phi = \xi^* \xi$.  
Therefore we have 
\[
\langle \Phi(\|\phi\|_\infty^2 1 - \overline{\phi}\phi)f,f\rangle
 = \langle \Phi(\xi^*\xi) f,f\rangle = \langle
\Phi(\xi)f,\Phi(\xi)f\rangle \geq 0.
\]
It follows that $\Phi(\phi)$ is a bounded operator on
$\Gamma_c(G,r^*\mcal{A})$ with norm less than $\|\phi\|_\infty$.  Thus
$\Phi(\phi)$ extends to an operator on $A\rtimes_\alpha G$.  Furthermore
\eqref{eq:113} and \eqref{eq:114} imply that $\Phi(\phi)$ is linear
with respect to the action of $A\rtimes S$ on itself and that
$\Phi(\phi)$ is adjointable with adjoint $\Phi(\phi)^* =
\Phi(\overline{\phi})$.  Hence $\Phi(\phi)\in M(A\rtimes G)$.  We have
already shown that $\Phi$ preserves the involution on $C_0(G\unit/G)$
and the computation 
\[
\Phi(\phi)\Phi(\psi)f(\gamma) = \phi(G\cdot r(\gamma))\psi(G\cdot
r(\gamma))f(\gamma) = \Phi(\phi\psi)f(\gamma)
\]
shows that it preserves multiplication as well. Thus $\Phi$ is a
$*$-homomorphism.   In order to show that $\Phi$ maps into the center
it will suffice to show, using Lemma \ref{lem:18}, that $\Phi(f*g)=
f*\Phi(g)$ for all $f,g\in\Gamma_c(G,r^*\mcal{A})$.  However, observe
that $r(\eta\inv\gamma) = \eta\inv\cdot r(\gamma)$ so that 
\begin{align*}
\Phi(\phi)(f*g)(\gamma) &= \phi(G\cdot r(\gamma)) \int_G f(\eta)
\alpha_\eta(g(\eta\inv \gamma))d\lambda^{r(\gamma)}(\eta) \\
&= \int_G f(\eta)\alpha_\eta(\phi(G\cdot
r(\eta\inv\gamma))g(\eta\inv\gamma)) d\lambda^{r(\gamma)}(\eta) \\
&= f*\Phi(\phi)g.
\end{align*}
Thus $\Phi(f)\in ZM(A\rtimes G)$.  The last thing we need to do 
is to show that
the set $\Phi(C_0(G\unit/G))\cdot A\rtimes S$ is dense in $A\rtimes
S$.  However given $f\in\Gamma_c(G,r^*\mcal{A})$ the image of $r(\supp
f)$ in $G\unit/G$, denoted $K= G\cdot r(\supp f)$, is compact.  Therefore
we can find $\phi\in C_c(G\unit/G)$ which is one on $K$ and zero off
some neighborhood of $K$.  It is easy to see that in this case
$\Phi(\phi)f = f$ and it follows immediately that $\Phi$ is
nondegenerate.  

Next we identify the fibers of $A\rtimes G$.  Fix $u\in G\unit$.  
Since $G\unit/G$
is Hausdorff, $G\cdot u$ is closed in $G\unit$.
Let $O = G\unit \setminus G\cdot u$.  It is clear that $G\cdot u$ and
$O$ are both $G$-invariant so that we may apply Theorem
\ref{thm:invtideal} to conclude that the restriction map $\rho$ factors to an
isomorphism of $A\rtimes G/\ker\rho$ with $A(G\cdot u)\rtimes G|_{G\cdot
  u}$.  Now define 
\[
I_u = \cspn \{\Phi(\phi)f : \phi\in C_0(G\unit/G), f\in
\Gamma_c(G,r^*\mcal{A}), \phi(G\cdot u) = 0\}.
\]
Since, by definition, $A\rtimes G(G\cdot u) = A\rtimes G/I_u$ it will
suffice to show that $I_u = \ker\rho$.  If $\phi(G\cdot u) = 0$ then
$\Phi(\phi)f(\gamma) = 0$ for all $\gamma\in G|_{G\cdot u}$ so that we
must have $\Phi(\phi)f \in \ker\rho$.  It follows that $I_u
\subset\ker\rho$.  On the other hand we also know from Theorem
\ref{thm:invtideal} that $\ker \rho = \Ex(O)$ where $\Ex(O)$ is the
ideal generated by those functions $f\in\Gamma_c(G,r^*\mcal{A})$ such
that $\supp f \subset G|_O$.  Now let $q:G\unit\rightarrow G\unit/G$
be the quotient map. Given $f\in\Gamma_c(G,r^*\mcal{A})$ such that 
$\supp f\subset G|_O$ we must have $q(r(\supp f))$ disjoint from
$G\cdot u$ in $G\unit/G$.  Since $q(r(\supp f))$ is compact we can
find some $\phi\in C_c(G\unit/G)$ such that $\phi$ is one on
$q(r(\supp f))$ and $\phi(G\cdot u) = 0$.  It follows that
$\Phi(\phi)f = f \in I_u$.  Thus $\ker\rho = \Ex(O)\subset I_u$ and we
are done. 
\end{proof}

The reason that this is a useful result is that we know a lot about
the fibres of $A\rtimes G$ when $G\unit/G$ is Hausdorff.  

\begin{corr}
\index{groupoid!transitive}
\label{cor:13}
Suppose $(A,G,\alpha)$ is a separable dynamical system and that the
orbit space 
$G\unit/G$ is Hausdorff.  Given $u\in G\unit$ the fibre $A\rtimes
G(G\cdot u)$ is Morita equivalent to $A(u)\rtimes S_u$.  
\end{corr}

\begin{proof}
Since $G\unit/G$ is Hausdorff $A\rtimes G$ is a $C_0(G\unit/G)$-algebra with
fibres 
\[
A\rtimes G(G\cdot u) \cong A(G\cdot u)\rtimes G|_{G\cdot u}.
\]
However, $G|_{G\cdot u}$ is a {\em transitive} groupoid so the result
follows from Theorem \ref{thm:transprod}. 
\end{proof}

Thus, in the case where $G\unit/G$ is Hausdorff, every irreducible
representation is lifted from a fibre $A\rtimes G(G\cdot u)$, 
and every irreducible representation of $A\rtimes G(G\cdot u)$ comes from
an irreducible representation of $A(u)\rtimes S_u$.  We will show that
this two stage description is nothing more than the usual induction
process.  

\begin{prop}
\label{prop:93}
Suppose $(A,G,\alpha)$ is a separable dynamical system and that the
orbit space
$G\unit/G$ is Hausdorff.  Then every irreducible representation of
$A\rtimes_\alpha G$ is equivalent to one of the form $\Ind_{S_u}^G R$
where $u\in G\unit$ and $R$ is an irreducible representation of
$A(u)\rtimes_\alpha S_u$. 
\end{prop}

We start with a remark and a useful lemma.  

\begin{remark}
Suppose we have an $A-B$-imprimitivity bimodule $\mcal{X}$ and a
representation $\pi$ of $B$.  The Rieffel induction process 
\cite[Proposition 2.66]{tfb} yields an induced
representation $\mcal{X}-\Ind \pi$ of $A$.  It is assumed that the
reader is familiar with this process.  If not they may wish to use
\cite[Section 2.4]{tfb} as a reference.  Furthermore, the Rieffel
correspondence provides a very strong link between the ideal structure
and representation theory of $A$ and the ideal structure and
representation theory of $B$.  This material can be found in
\cite[Section 3.3]{tfb}. 
\end{remark}

\begin{lemma}
\label{lem:27}
Suppose $(A,G,\alpha)$ is a separable dynamical system and that the
orbit space 
$G\unit/G$ is Hausdorff.  Given $u\in G\unit$ let $\rho:A\rtimes
G\rightarrow A(G\cdot u)\rtimes G|_{G\cdot u}$ be the extension of the
restriction map on $\Gamma_c(G,r^*\mcal{A})$.  Furthermore, let
$\mcal{X}$ be the $A(G\cdot u)\rtimes G|_{G\cdot u}-A(u)\rtimes S_u$
imprimitivity bimodule from Theorem \ref{thm:transprod}.  If 
$R$ is a representation of $A(u)\rtimes S_u$ then $\Ind_{S_u}^G R =
\mcal{X}-\Ind(R)\circ \rho$.  
\end{lemma}

\begin{proof}
First let us establish some notation.  Let $\beta$ be Haar measure on
$S_u$.  Recall that $\mcal{X}$ is the
completion of the pre-$A(G\cdot u)\rtimes G|_{G\cdot u}-A(u)\rtimes
S_u$-imprimitivity bimodule $\mcal{X}_0 = C_c(G_u, A(u))$ and that the
left hand operations on $\mcal{X}_0$ are given by 
\begin{align*}
z\cdot g(\gamma) &= \int_{S_u} \alpha_s(z(\gamma s)g(s\inv))d\beta(s),
\\
\llangle z,w\rrangle_{A(u)\rtimes S_u}(s) &= \int_G
z(\eta\inv)^*\alpha_s(w(\eta\inv s))d\lambda^u(\eta).
\end{align*}
Next let $X = s\inv((S_u)\unit) = G_u$ and recall that the imprimitivity
bimodule $\mcal{Z}_{S_u}^G$ is the completion of $\mcal{Z}_0 =
\Gamma_c(X,s^*\mcal{A}) = C_c(G_u,A(u))$ and carries the left hand
actions 
\begin{align*}
z\cdot g(\gamma) &= \int_{S_u} \alpha_s(z(\gamma s)g(s\inv))d\beta(s)
\\
\llangle z,w \rrangle_{A(u)\rtimes S_u}(s) &= \int_G z(\eta
s\inv)^*\alpha_s(w(\eta))d\lambda_u(\eta), \\
&= \int_G z(\eta\inv)^*\alpha_s(w(\eta\inv s))d\lambda^u(\eta).
\end{align*}
It is a happy fact that $\mcal{Z}_{S_u}^G$ and $\mcal{X}$ are
obviously equal as right
Hilbert $A(u)\rtimes S_u$-modules.  

Suppose $R$ is a representation of $A(u)\rtimes S_u$ on
$\mcal{H}$.  Recall from Theorem \ref{thm:induce} that $\Ind_{S_u}^G
R$ acts on $\mcal{K} = \mcal{Z}_{S_u}^G \otimes_{A(u)\rtimes S_u}
\mcal{H}$ via $\Ind R(f)(z\otimes h) = f\cdot z\otimes h$ where, given
$f\in\Gamma_c(G,r^*\mcal{A})$ and $z\in C_c(G_u,A(u))$,
\begin{equation}
\label{eq:116}
f\cdot z(\gamma) = \int_G
\alpha_\gamma\inv(f(\eta))z(\eta\inv\gamma)d\lambda^{r(\gamma)}(\eta).
\end{equation}
However, since $\mcal{Z}_{S_u}^G = \mcal{X}$ as right Hilbert
$A(u)\rtimes S_u$-modules, the representation $\mcal{X}-\Ind(R)$ also
acts on $\mcal{K}$.  Furthermore the action is given by
$\mcal{X}-\Ind(S)f(z\otimes h) = f\odot z\otimes h$ where 
$f\odot z$ is the left module action of $\mcal{X}$ and is given
for $f\in \Gamma_c(G|_{G\cdot u},r^*\mcal{A})$ and $z\in
C_c(G_u,A(u))$ by 
\begin{equation}
\label{eq:117}
f\odot z(\gamma) = \int_G
\alpha_\gamma\inv(f(\eta))z(\eta\inv\gamma)d\lambda^{r(\gamma)}(\eta)
\end{equation}
However \eqref{eq:116} and \eqref{eq:117} are basically the same
action and since $\rho$ is  an extension of the restriction map it
is clear that for $f\in\Gamma_c(G,r^*\mcal{A})$ we have $\Ind_{S_u}^G
R(f) = \mcal{X}-\Ind(R)(\rho(f))$ and this extends to the entire
crossed product by continuity. 
\end{proof}

Actually, with this result at our disposal we are mostly done.  

\begin{proof}[Proof of Proposition \ref{prop:93}]
Suppose $(A,G,\alpha)$ is a separable dynamical system and that the
orbit space
$G\unit/G$ is Hausdorff.  By Proposition \ref{prop:92} $A\rtimes
G$ is a $C_0(G\unit/G)$-algebra.  It then follows from Proposition
\ref{prop:36} that  any irreducible
representation $T$ is of the form $T = L\circ \rho$ where $u\in
G\unit$, $L$ is an
irreducible representation of $A(G\cdot u)\rtimes G|_{G\cdot u}$ and
$\rho$ is the canonical extension of the restriction map on
$\Gamma_c(G,r^*\mcal{A})$.  However $A(G\cdot u)\rtimes G|_{G\cdot u}$
is Morita equivalent to $A(u)\rtimes S_u$ by Corollary \ref{cor:13}.
Let $\mcal{X}$ be the bimodule implementing the equivalence and let
$\widetilde{\mcal{X}}$ be its ``inverse'' bimodule.  Set
$R=\widetilde{\mcal{X}}-\Ind L$.  It follows from 
\cite[Theorem 3.29]{tfb} that $\mcal{X}-\Ind R$ is naturally equivalent to
$L$ and from \cite[Corollary 3.32]{tfb} that $R$ is an 
irreducible representation.  However, it follows that the representations $T =
L\circ \rho$ and $\Ind_{S_u}^G R = \mcal{X}-\Ind(R)\circ \rho$ are also 
equivalent and we are done. 
\end{proof}

The reason we separated out Lemma \ref{lem:27} is it allows us to
easily prove the following 

\begin{prop}
\label{prop:94}
Suppose $(A,G,\alpha)$ is a separable dynamical system and that the
orbit space 
$G\unit/G$ is Hausdorff.  If $R$ is an irreducible representation of
$A(u)\rtimes S_u$ then $\Ind_{S_u}^G R$ is irreducible.  Furthermore
if $L$ and $R$ are both irreducible representations of $A(u)\rtimes
S_u$ and $\Ind_{S_u}^G L$ is equivalent to $\Ind_{S_u}^G R$ then $L$
is equivalent to $R$. 
\end{prop}

\begin{proof}
Given $(A,G,\alpha)$ and $R$ as above let $\rho:A\rtimes G\rightarrow
A(G\cdot u)\rtimes G|_{G\cdot u}$ be the extension of the restriction
map.  Furthermore let $\mcal{X}$ be the $A(G\cdot u)\rtimes G|_{G\cdot
  u}-A(u)\rtimes S_u$-imprimitivity bimodule from Theorem
\ref{thm:transprod}.  Since the Rieffel correspondence preserves
irreducibility $\mcal{X}-\Ind(R)$ is an irreducible representation of
$A(G\cdot u)\rtimes G|_{G\cdot u}$.  Hence $\Ind_{S_u}^G R =
\mcal{X}-\Ind(R)\circ \rho$ must be irreducible.  

Next suppose we are given two irreducible representations $L$ and $R$
of ${A(u)\rtimes S_u}$ and suppose $\Ind_{S_u}^G L$ and $\Ind_{S_u}^G R$
are equivalent.  Since $\rho$ is surjective it follows from Lemma
\ref{lem:27} that $\mcal{X}-\Ind(L)$ and $\mcal{X}-\Ind(R)$ are
equivalent.  Since $\mcal{X}$ is an imprimitivity bimodule this
implies $L$ is equivalent to $R$. 
\end{proof}

\index{Mackey-Glimm Dichotomy}
We are going to extend Proposition \ref{prop:93} to groupoids
which satisfy the Mackey-Glimm dichotomy.  For our purposes the most
useful condition of Theorem \ref{thm:glimmdich} will be the fact that
the orbit space is almost Hausdorff.  

\begin{definition}
\index{almost Hausdorff}
\label{def:almosthauss}
A, not necessarily Hausdorff, locally compact space $X$ is said to be
{\em almost Hausdorff} if each locally compact subspace $V$ contains
a relatively open nonempty Hausdorff subset. 
\end{definition}

The key fact we will use about almost Hausdorff spaces is the
following proposition, which we cite without proof. Those readers
unfamiliar with ordinals are referenced to \cite[Chapter 6]{hrbacek}.  

\begin{prop}[{\cite[Lemma 6.3]{tfb2}}]
\label{prop:95}
Suppose $X$ is a, not necessarily Hausdorff, locally compact space.
Then the following are equivalent. 
\begin{enumerate}
\item $X$ is almost Hausdorff
\item Every nonempty closed subspace of $X$ has a relatively open
  nonempty Hausdorff subspace.  
\item Every closed subspace of $X$ has a dense relatively open
  Hausdorff subspace. 
\item There is an ordinal $\gamma$ and open sets $\{U_\alpha :\alpha
  \leq \gamma\}$ such that 
\begin{enumerate}
\item $\alpha < \beta \leq \gamma$ implies that $U_\alpha \subsetneq
  U_\beta$,
\item $\alpha < \gamma$ implies that $U_{\alpha+1}\setminus
  U_{\alpha}$ is a dense Hausdorff subspace of $X\setminus U_\alpha$, 
\item if $\delta\leq \gamma$ is a limit ordinal then 
\[
U_\delta = \bigcup_{\alpha < \delta} U_\alpha , 
\]
\item $U_0 = \emptyset$ and $U_\gamma = X$.
\end{enumerate}
\item Every subspace of $X$ has a relatively open dense Hausdorff
  subspace. 
\end{enumerate}
\end{prop}

The main reason we care about Proposition
\ref{prop:95} is that condition (d) will allow us to build the
following object.

\begin{definition}
\label{def:58}
A {\em composition series} in a $C^*$-algebra $A$ is a family
$\{I_\alpha\}_{\alpha\in \Lambda}$ of ideals $I_\alpha$ indexed by a
segment $\Lambda$ of ordinals $0\leq \alpha \leq \gamma$ such that 
\begin{enumerate}
\item $I_0 = \{0\}$ and $I_\gamma = A$, 
\item $\alpha < \beta \leq \gamma$ implies $I_\alpha \subsetneq
  I_\beta$ and
\item if $\delta\leq \gamma$ is a limit ordinal then 
\[
I_\delta = \overline{\bigcup_{\alpha < \delta}I_\alpha}.
\]
\end{enumerate}
\end{definition}

Continuing the chain, the main reason we care about Definition
\ref{def:58} is the following 

\begin{lemma}[{\cite[Lemma 8.13]{tfb2}}]
\label{lem:28}
Suppose $\{I_\alpha\}_{\alpha\in\Lambda}$ is a composition series for
a $C^*$-algebra $A$.  Then every irreducible representation $\pi$ of $A$
lives on a subquotient $I_{\alpha+1}/I_\alpha$ for some $\alpha$.  In
other words there is an irreducible representation $\rho$ of
$I_{\alpha+1}/I_\alpha$ such that $\pi$ is equal to the canonical
extension of the lift of $\rho$ to $I_{\alpha+1}$.  
\end{lemma}

\begin{proof}
Let $\pi$ be an irreducible representation of $A$.  Let
$S=\{\alpha\leq \gamma: I_\alpha\not\subset \ker\pi\}$.  If $\beta =
\min S$ is a limit ordinal, then it follows from part (c) of
Definition \ref{def:58} that $I_\beta \subset \ker\pi$.  However, this
contradicts the fact that $\beta\in S$.  Thus $\beta$ has an immediate
predecessor $\alpha$.  Let $\rho$ be the factorization of
$\pi|_{I_\beta}$ to $I_\beta/I_\alpha$.  Then clearly the lift of
$\rho$ to $I_\beta$ is $\pi|_{I_\beta}$ and since $I_\beta\not\subset
\ker\pi$ the extension of $\pi|_{I_\beta}$ to $A$ is $\pi$. 
\end{proof}

At this point it may be clear where we are going.  If
$G\unit/G$ is almost Hausdorff then we will build a composition series
of crossed products where the orbit space associated to the subquotients is
Hausdorff.  This will allow us to use Proposition \ref{prop:93} to
prove the following theorem.  

\begin{theorem}
\index{New Result}
\label{thm:ehregularity}
Suppose $(A,G,\alpha)$ is a separable groupoid dynamical system and that
$G/G\unit$ is a $T_0$ space.  Then every irreducible representation of
$A\rtimes_\alpha G$ is equivalent to one of the form $\Ind_{S_u}^G R$
where $u\in G\unit$ and $R$ is an irreducible representation of
$A(u)\rtimes_\alpha S_u$. 
\end{theorem}

As before we start with a utility lemma. 

\begin{lemma}
\label{lem:29}
Suppose $(A,G,\alpha)$ is a separable dynamical system and that
$U\subset V \subset G\unit$ are open $G$-invariant sets.  Then we may
identify $A(V\setminus U)\rtimes G|_{V\setminus U}$ with the
subquotient 
$\Ex(V)/\Ex(U)$.  Furthermore if $u\in U$ and $R$ is a representation of
${A(u)\rtimes S_u}$ then the canonical extension of
$\Ind_{S_u}^{G|_{V\setminus U}} R$ to $A\rtimes G$ is equal to $\Ind_{S_u}^G R$. 
\end{lemma}

\begin{proof}
First recall that we equip $G|_V$ with the restriction of the Haar system
from $G$.  Furthermore we equip $(G|_V)|_{V\setminus U} =
G|_{V\setminus U}$ with the restriction of the Haar system coming from
$G|_V$, and therefore from $G$.  Since $V$ is an open $G$-invariant
set we use Proposition \ref{prop:73} to identify $A(V)\rtimes G|_V$
with the ideal $\Ex(V)$ via the inclusion map $\iota$.  Furthermore we also
identify $A(U)\rtimes G|_U$ with $\Ex(U)$.  Since $U\subset V$ any
function which is supported in $G|_U$ must be supported in $G|_V$ as
well and therefore $\Ex(U)\subset \Ex(V)$.  
Now $U$ is also an open $G|_V$ invariant
subset of $V$ so that we can also identify ${A(U)\rtimes G|_U}$ with its
image $\Ex'(U)$ in $A(V)\rtimes G|_V$.  We would like to see that
$\iota(\Ex'(U)) = \Ex(U)$.  However, this is obvious since both $\Ex(U)$ and
$\iota(\Ex'(U))$ are the completion of
$\Gamma_c(G|_U,r^*\mcal{A})$ inside $A\rtimes G$.  It now follows that
$\iota$ factors to an isomorphism of $A(V)\rtimes G|_V/\Ex'(U)$ onto
the subquotient $\Ex(V)/\Ex(U)$.   Furthermore Theorem
\ref{thm:invtideal} implies that the restriction map $\rho$ factors to an
isomorphism of $A(V)\rtimes G|_V/\Ex'(U)$ with $A(V\setminus U)\rtimes
G|_{V\setminus U}$. Thus $A(V\setminus U)\rtimes G|_{V\setminus U}$ is
isomorphic to the subquotient $\Ex(V)/\Ex(U)$.  Therefore, 
given a representation $T$ of $A(V\setminus U)\rtimes G|_{V\setminus
  U}$ we can take its lift $T\circ \rho$ to $A(V)\rtimes G|_V$ and
then extend $T\circ \rho\circ \iota\inv$ from $\Ex(V)$ to a
representation of $A\rtimes G$.  
Of course when we are
working with elements of $\Gamma_c(G|_V,r^*\mcal{A})$ the $\rho$ and
$\iota\inv$ maps basically disappear so that we will usually not be
this precise about viewing $A(V\setminus U)\rtimes G|_{V\setminus U}$
as a subquotient.  

So suppose $R$ is a representation of $A(u)\rtimes S_u$ on $\mcal{H}$ 
for $u\in V\setminus U$.  Recall that, as in the proof of
Lemma \ref{lem:27},  $\mcal{Z}_{S_u}^G$ is
the completion of $C_c(G_u,A(u))$ with respect to the following left operations 
\begin{align}
\label{eq:118}
z\cdot g(\gamma)&= \int_G \alpha(z(\gamma s)g(s\inv))d\beta(s), \\
\label{eq:119}
\llangle z,w\rrangle_{A(u)\rtimes S_u}(s) &= \int_G
z(\eta\inv)^*\alpha_s(w(\eta\inv s))d\lambda^u(\eta).
\end{align}
Furthermore $\Ind_{S_u}^G R$ acts on $\mcal{K} =
Z_{S_u}^G\otimes_{A(u)\rtimes S_u} \mcal{H}$ via $\Ind_{S_u}^G
R(f)(z\otimes h) = f\cdot z \otimes h$ 
where, given $f\in\Gamma_c(G,r^*\mcal{A})$, we define 
\begin{equation}
\label{eq:120}
f\cdot z(\gamma) = \int_G
\alpha_\gamma\inv(f(\eta))g(\eta\inv\gamma)d\lambda^{r(\gamma)}(\eta).
\end{equation}
Since the Haar system on $G|_{V\setminus U}$ is just the
restriction of the Haar system of $G$,
$\mcal{Z}_{S_u}^{G|_{V\setminus U}}$ is also the completion of
$C_c(G_u,A(u))$ with respect to the operations \eqref{eq:118} and
\eqref{eq:119}.  Hence $\Ind_{S_u}^{G|_{V\setminus U}} R$ also acts on
$\mcal{K}$ and the action is given by $\Ind_{S_u}^{G|_{V\setminus
U}} R(f)(z\otimes h) = f\cdot z\otimes h$ where $f\cdot z$ is
defined via \eqref{eq:120} on $\Gamma_c(G|_{V\setminus U},r^*\mcal{A})$.
At this point it is clear that for $f\in\Gamma_c(G|_V,r^*\mcal{A})$ we
have 
\[
\Ind_{S_u}^G R(f) = \Ind_{S_u}^{G|_{V\setminus U}}R (\rho(\iota\inv(f))).
\]
It follows that $\Ind_{S_u}^G R$ agrees with
$\Ind_{S_u}^{G|_{V\setminus U}}R\circ \rho\circ\iota\inv$ on $\Ex(V)$.
Hence $\Ind_{S_u}^G R$ is equal to the unique extension of 
$\Ind_{S_u}^{G|_{V\setminus U}} R\circ\rho\circ\iota\inv$ to $A\rtimes
G$ and we are done.  
\end{proof}

This get us most of the way there since it shows that the process of
lifting representations from a subquotient and induction are
compatible.  We can now prove the main result of this section.  

\begin{proof}[Proof of Theorem \ref{thm:ehregularity}]
Since $G$ is a second countable, locally compact Hausdorff groupoid, the
fact that $G\unit/G$ is $T_0$ implies, by Theorem \ref{thm:glimmdich},
that $G\unit/G$ is almost Hausdorff.  Therefore there are open sets
$\{V_\beta\}_{0\leq \beta\leq \gamma}$ in $G\unit/G$ satisfying properties
(i)-(iv) of Proposition \ref{prop:95}.  Let $q:G\unit\rightarrow
G\unit/G$ be the quotient map and $U_\beta = q\inv(V_\beta)$ for all
$0\leq \beta\leq \gamma$.  Then each $U_\beta$ is an open $G$-invariant
subset and we define $I_\beta = \Ex(U_\beta)$.  Since $U_0 =
\emptyset$ and $U_{\gamma} = G\unit$ we must have $I_0 = \{0\}$ and
$I_\gamma = A\rtimes G$.  Furthermore if $\delta< \beta\leq \gamma$
then $U_\delta \subsetneq U_\beta$.  Thus any function supported in
$G|_{U_\delta}$ must be supported in $G|_{U_\beta}$ as well so we must
have $I_\delta \subset I_\beta$.  Since $U_\delta \ne U_\beta$ it is
easy to build a function supported on $U_\beta$ and not $U_\delta$.
Thus $I_\delta \ne I_\beta$.  Finally suppose $\delta\leq \gamma$
is a limit ordinal and $f\in \Gamma_c(G|_{U_\delta},r^*\mcal{A})$.
Because $r(\supp f)\subset U_\delta = \bigcup_{\beta<\delta}U_\beta$ the
collection $\{U_\beta\}_{\beta<\delta}$ is an open cover of $r(\supp
f)$.  Since $r(\supp f)$ is compact there must be a finite subcover
and since the $U_\beta$ are nested this implies that there exists
$\beta'<\delta$ such that $r(\supp f)\subset U_{\beta'}$.  Hence $f\in
\Gamma_c(G|_{U_{\beta'}},r^*\mcal{A})\subset I_{\beta'}$.  It follows
that $I_\delta \subset \overline{\bigcup_{\beta < \delta}
  I_\beta}$. The other inclusion is trivial so that we have
\[
I_\delta = \overline{\bigcup_{\beta < \delta} I_\beta}. 
\]
Thus $\{I_\beta\}$ is a composition series for $A\rtimes G$.  

Now suppose $L$ is a irreducible representation of $A\rtimes G$.  
Lemma \ref{lem:28} implies that there exists $\beta$ such that $L$
lives on $I_{\beta+1}/I_\beta$.  In other words, there is an
irreducible representation $T$ of $I_{\beta+1}/I_\beta$ such that $L$ is the
unique extension of the lift of $T$.  Next, Lemma \ref{lem:29} implies
that we can identify $I_{\beta+1}/I_\beta$ with
$A(U_{\beta+1}\setminus U_\beta)\rtimes G|_{U_{\beta+1}\setminus
  U_\beta}$.  Furthermore $(U_{\beta+1}\setminus
U_\beta)/G = V_{\beta+1}\setminus V_\beta$ is Hausdorff so that by
Proposition \ref{prop:93} there exists $u\in U_{\beta+1}\setminus
U_\beta$ and an irreducible representation $R$ of $A(u)\rtimes S_u$
such that $T$ is equivalent to $R'=\Ind_{S_u}^{G|_{U_{\beta+1}\setminus
    U_\beta}} R$.  Hence the extension of $T$ to $A\rtimes G$, which
is $L$, is equivalent to the extension of $R'$ to $A\rtimes G$, which
is $\Ind_{S_u}^G R$ by Lemma \ref{lem:29}.  
\end{proof}

As before, we separated out Lemma \ref{lem:29} so that we could prove
the following 

\begin{prop}
\label{prop:96}
Suppose $(A,G,\alpha)$ is a separable dynamical system and the orbit space
$G\unit/G$ is $T_0$. If $R$ is an irreducible representation
of $A(u)\rtimes S_u$ then $\Ind_{S_u}^G R$ is irreducible.
Furthermore, if $R$ and $L$ are both irreducible representations of
$A(u)\rtimes S_u$ and $\Ind_{S_u}^G R$ is equivalent to $\Ind_{S_u}^G
L$ then $R$ is equivalent to $L$. 
\end{prop}

\begin{proof}
Suppose $R$ is an irreducible representation of $A(u)\rtimes S_u$.
Using Theorem \ref{thm:glimmdich} $G\unit/G$ must be locally
Hausdorff.  Let $\{V_\beta\}$ be as in Proposition \ref{prop:95}.
Consider $\Gamma = \{\beta\leq \gamma: G\cdot u \in V_\beta\}$.  If
$\delta = \min\Gamma$ is a limit ordinal then $G\cdot u \in
\bigcup_{\beta < \delta} V_\beta$.  However this implies $G\cdot u \in
V_\beta$ for some $\beta < \delta$.  This is a contradiction.  It
follows that $\delta$ has an immediate predecessor $\sigma$ and
$G\cdot u\in V_{\delta}\setminus V_\sigma$.  Let $q:G\unit\rightarrow
G\unit/G$ be the quotient map, $U_\delta = q\inv(V_\delta)$ and
$U_\sigma = q\inv(V_\sigma)$.  Then $u\in U_\delta\setminus U_\sigma$
and since $V_\delta\setminus V_\sigma$ is Hausdorff we can use
Proposition \ref{prop:94} to conclude that
$R'=\Ind_{S_u}^{G|_{U_\delta\setminus U_{\sigma}}} R$ is irreducible.
However it follows that the extension of $R'$ to $A\rtimes G$ is irreducible
and by Lemma \ref{lem:29} this is exactly $\Ind_{S_u}^G R$.  

Now suppose $R$ and $L$ are both irreducible representations on
$A(u)\rtimes S_u$ and that $\Ind_{S_u}^G R$ is equivalent to
$\Ind_{S_u}^G L$.  This implies that their factorizations to the
subquotient $A(U_\delta\setminus U_\sigma)\rtimes
G|_{U_\delta\setminus U_\sigma}$ are equivalent.  It follows from
Lemma \ref{lem:29} that these factorizations are 
$\Ind_{S_u}^{G|_{U_\delta\setminus U_{\sigma}}} R$ and
$\Ind_{S_u}^{G|_{U_\delta\setminus U_{\sigma}}} L$ is respectively.
Hence, Proposition \ref{prop:94} implies that $R$ and $L$ must be
equivalent.  
\end{proof}

\begin{remark}
\label{rem:22}
It would be tempting, in light of Theorem \ref{thm:ehregularity}, to say that
every representation of $A\rtimes G$ is induced from a stabilizer.
Unfortunately, this notion has a conflicting definition in
\cite[Definition 8.10]{tfb2}.  The problem lies in the meaning of the
word stabilizer.  In \cite{tfb2} the stabilizers are the stabilizer
subgroups with respect to the action of $G$ on $\Prim A$.  In Theorem
\ref{thm:ehregularity} the stabilizers are with respect to the action
of $G$ on its unit space, which may be larger. 
Of course, when $A$ has Hausdorff spectrum equal to $G\unit$
these two notions match up.  Furthermore when $A$ has Hausdorff
spectrum it is not hard to show, using Example \ref{ex:22}, that
Theorem \ref{thm:ehregularity} generalizes \cite[Theorem 8.16]{tfb2}.  
It's also worth pointing out that this is the only way we can view the
results of this section as generalizing the group case.  Theorem
\ref{thm:ehregularity} is trivial if we take the naive approach and
treat groups as groupoids with a single unit.  
\end{remark}

\begin{remark}
\index{Mackey-Glimm Dichotomy}
Generalizing this result to groupoids which do not satisfy the
Mackey-Glimm Dichotomy is difficult.  For group crossed products the
result is known as the Gootman-Rosenberg-Sauvageot (GRS) theorem.  The
method of attack was developed by Sauvageot in \cite{sauvageot1,sauvageot2} and
the complete solution was given by Gootman and Rosenberg in \cite{goodmanrosen}.
The result is also proved in \cite[Chapter 9]{tfb2}.  
For groupoid $C^*$-algebras,  
the corresponding result is proved in \cite{geneffhan} and for general 
groupoid crossed products the question is still open. 
\end{remark}


\section{Crossed Products with Abelian Isotropy}
\label{sec:crossedstab}

Theorem \ref{thm:ehregularity} is a nice enough result, but if we want to
study the fine structure of $A\rtimes G$ we need to consider more than
just individual representations.  This next proposition adds a
topological component to the results of the last section.  

\begin{prop}
\label{prop:97}
Suppose $(A,G,\alpha)$ is a separable dynamical system.  Furthermore,
suppose that the isotropy subgroupoid $S$ varies continuously and that
$G\unit/G$ is a $T_0$ space.  Then $\Phi:(A\rtimes S)\sidehat
\rightarrow (A\rtimes G)\sidehat$ given by $\Phi(R) = \Ind_S^G R$ is
a continuous surjection. 
\end{prop}

Recall that $A\rtimes S$ is a $C_0(G\unit)$-algebra and that
restriction factors to an isomorphism of $A\rtimes S(u)$ with
$A(u)\rtimes S_u$. The main difficulty is to show that induction
respects this fibering.  

\begin{lemma}
\label{lem:30}
Suppose $(A,G,\alpha)$ is a separable dynamical system and that the
isotropy subgroupoid $S$ varies continuously.  Given $u\in G\unit$ and
a representation $R$ of $A(u)\rtimes S_u$ let $\rho:A\rtimes
S\rightarrow A(u)\rtimes S_u$ be given on $\Gamma_c(S,p^*\mcal{A})$ by
restriction.  Then $\Ind_{S_u}^G R$ is equivalent to $\Ind_S^G (R\circ \rho)$.  
\end{lemma}

\begin{proof}
First, since $S$ is continuously varying, it is a closed subgroupoid of
$G$ with its own Haar system, which we call $\beta$. 
Let $u\in G\unit$ and $R$ be a representation of $A(u)\rtimes S_u$ on
$\mcal{H}$.  Consider the right Hilbert
$A(u)\rtimes S_u$-module $\mcal{Z}_{S_u}^G$ associated to
$\Ind_{S_u}^G R$.  It follows from Proposition \ref{prop:80} that
$\mcal{Z}_{S_u}^G$ is the completion of $C_c(G_u,A(u))$ with respect
to the right actions 
\begin{align*}
z\cdot g(\gamma) &= \int_{S_u} \alpha_s(z(\gamma
s)g(s\inv))d\beta^u(s) \\
\llangle z,w \rrangle_{A(u)\rtimes S_u}(s) &= \int_G
z(\eta s\inv)^*\alpha_s(w(\eta))d\lambda_u(\eta).
\end{align*}
Furthermore, $\Ind_{S_u}^G R$ acts on $\mcal{Z}_{S_u}^G
\otimes_{A(u)\rtimes S_u}\mcal{H}$ which, as in Remark \ref{rem:19}, is
the completion of
$C_c(G_u,A(u))\odot \mcal{H}$ with respect to the inner product
characterized by
\[
(f\otimes h,g\otimes k) = (R(\llangle g,f\rrangle_{A(u)\rtimes
  S_u})h,k).
\]
Now consider the right Hilbert $A\rtimes S$-module $\mcal{Z}_S^G$
associated to $\Ind_S^G (R\circ \rho)$.  
It follows from Proposition \ref{prop:80}
that $\mcal{Z}_S^G$ is the completion of $\Gamma_c(G,s^*\mcal{A})$ with
respect to the operations 
\begin{align*}
z\cdot g(\gamma) &= \int_S \alpha_s(z(\gamma s)g(s\inv))
d\beta^{s(\gamma)}(s) \\
\llangle z,w \rrangle_{A\rtimes S}(s) &= \int_G z(\eta s\inv)^*
\alpha_s(w(\eta))d\lambda_{p(s)}(\eta).
\end{align*}
Furthermore $\Ind_S^G (R\circ\rho)$ acts on 
$\mcal{Z}_{S}^G\otimes_{A\rtimes S} \mcal{H}$ which is the completion
of $\Gamma_c(G,s^*\mcal{A})\odot \mcal{H}$ with respect to the inner
product characterized by 
\[
(f\otimes h, g\otimes k) = (R(\rho(\llangle g,f\rrangle_{A\rtimes
  S}))h,k)
\]

We would like to define a unitary map $U:\mcal{Z}_S^G\otimes
\mcal{H}\rightarrow \mcal{Z}_{S_x}^G\otimes \mcal{H}$.  Start by
letting $\pi:\Gamma_c(G,s^*\mcal{A})\rightarrow C_c(G_u,A(u))$ be
given by restriction.  Since $\pi$ is clearly linear we can define
$U:\Gamma_c(G,s^*\mcal{A})\odot \mcal{H}\rightarrow C_c(G_u,A(u))\odot
\mcal{H}$ on elementary tensors by $U(f\otimes h) = \pi(f)\otimes
h$.  It is clear enough that given $f,g\in \Gamma_c(G,s^*\mcal{A})$ we have 
$\rho(\llangle f,g\rrangle_{A\rtimes S}) = \llangle \pi(f),\pi(g)
\rrangle_{A(u)\rtimes S_u}$ so that 
\begin{align*}
(f\otimes h, g\otimes k) &= (R(\rho(\llangle g,f\rrangle_{A\rtimes
  S}))h,k) \\
&=(R(\llangle \pi(g),\pi(f)\rrangle_{A(u)\rtimes
  S_u})h,k) \\ 
&= (\pi(f)\otimes h, \pi(g)\otimes k)\\
&= (U(f\otimes h),U(f\otimes k)).
\end{align*}
Thus $U$ is isometric on $\Gamma_c(G,s^*\mcal{A})\odot\mcal{H}$ and we
can extend it to an isometry on $\mcal{Z}_S^G\otimes \mcal{H}$.  We
need to show that $U$ is surjective.  Suppose $f\in C_c(G_u)$ and
$a\in A(u)$.  Choose $b\in A$ such that $b(u) = a$ and extend $f$
to a function $g\in C_c(G)$.  Then $g\otimes b\in
\Gamma_c(G,s^*\mcal{A})$ and given $\gamma\in G_u$ we clearly have
$\pi(g\otimes b)(\gamma) = f(\gamma)a = f\otimes a(\gamma)$.  Thus
$\ran\pi$ contains all of the elementary tensors in $C_c(G_u,A(u))$
and as such is dense in the inductive limit topology.  Now suppose
$z_i\rightarrow z$ with respect to the inductive limit topology in
$C_c(G_u,A(u))$.  Let $K$ be some compact set which eventually
contains the supports of the $z_i$.  Then 
\begin{align*}
\|\llangle z_i, z_i\rrangle&_{A(u)\rtimes S_u}(s) - \llangle z,z
  \rrangle_{A(u)\rtimes S_u}(s)\| \\
&\leq 
\int_G \|z_i(\eta s\inv)^*\alpha_s(z_i(\eta)) - z(\eta
s\inv)^*\alpha_s(z(\eta))\| d\lambda_u(\eta) \\
&\leq \int_G \|z_i(\eta s\inv)-z(\eta s\inv)\|\|z_i(\eta)\| + 
\|z(\eta s\inv)\| \|z_i(\eta)-z(\eta)\|d\lambda_u(\eta) \\
&\leq \|z_i-z\|_\infty(\|z_i\|_\infty+\|z\|_\infty)\lambda^u(K)
\end{align*}
Since $\{\|z_i\|_\infty\}$ is bounded this shows that $\llangle
z_i,z_i\rrangle \rightarrow \llangle z,z\rrangle$ uniformly.
Furthermore $\supp \llangle z_i,z_i\rrangle$ is eventually contained
in $K\inv K$, which is compact.  Thus $\llangle
z_i,z_i\rrangle\rightarrow \llangle z,z\rrangle$ with respect to the
inductive limit topology and hence 
\[
R(\llangle z_i,z_i\rrangle_{A(u)\rtimes S_u})\rightarrow
R(\llangle z,z\rrangle_{A(u)\rtimes S_u}).
\]
Now suppose $z\in C_c(G_u,A(u))$ and $h\in \mcal{H}$.   Choose $z_i =
\pi(w_i)$ such that $z_i\rightarrow z$ with respect to the inductive
limit topology. Then, by the above, since $z_i-z\rightarrow 0$ with
respect to the inductive limit topology, 
\[
\|U(w_i\otimes h)-z\otimes h\|^2 = \|(z_i-z)\otimes h\|^2
= (R(\llangle z_i-z,z_i-z\rrangle_{A(u)\rtimes S_u})h,h) \rightarrow 0.
\]
It follows that $\ran U$ is dense in
$\mcal{Z}_{S_u}^G\otimes\mcal{H}$ and that $U$ is a unitary.  

The last step is to show that $U$ intertwines $\Ind_S^G (R\circ
\rho)$ and $\Ind_{S_u}^G R$.  According to Theorem \ref{thm:induce}
$\Ind_{S_u}^R$ acts via $\Ind_{S_u}^G R(f)(z\otimes h) = f\cdot
z\otimes h$ where 
\begin{equation}
\label{eq:121}
f\cdot z(\gamma)=\int_G
\alpha_\gamma\inv(f(\eta))z(\eta\inv\gamma)d\lambda^{r(\gamma)}(\eta).
\end{equation}
and that $\Ind_S^G (R\circ \rho)$ also acts via
$\Ind_S^G(R\circ\rho)(f)(z\otimes h)= f \cdot
z \otimes h$ where $f\cdot z$ is still given by \eqref{eq:121}.  It is
clear that for $z\in\Gamma_c(G,s^*\mcal{A})$ we have $f\cdot
\pi(z)(\gamma) = \pi(f\cdot z)(\gamma)$.  Thus on
$\Gamma_c(G,s^*\mcal{A})\odot \mcal{H}$
\[
\Ind_S^G(R\circ \rho)(f)U(z\otimes h) = f\cdot \pi(z)\otimes h = 
\pi(f\cdot z)\otimes h = U\Ind_{S_u}^GR(f)(z\otimes h).
\]
This suffices to show that $U$ intertwines $\Ind_S^G(R\circ \rho)$ and
$\Ind_{S_u}^G R$.  
\end{proof}

\begin{remark}
\label{rem:25}
In light of how natural the unitary intertwining $\Ind_{S_u}^G R$ and
$\Ind_S^G(R\circ \rho)$ is we will often confuse the two.
Furthermore, since every irreducible representation of $A\rtimes S$ is
lifted from a fibre via restriction we will feel free to use the
notation $\Ind_S^G R$ even when $R$ is an irreducible representation of
$A(u)\rtimes S_u$.  Furthermore we will interpret $\Ind_S^G R$ as
either $\Ind_{S_u}^G R$ or as $\Ind_S^G(R\circ \rho)$ as we see
fit. We trust the reader will forgive the author for these abuses.  
\end{remark}

The advantage of viewing the induction as occurring on $S$ is that
induction from a fixed algebra is a continuous process.  

\begin{proof}[Proof of Proposition \ref{prop:97}]
Since $S$ is a continuously varying group bundle we know $A\rtimes S$
exists and that every irreducible representation is lifted from a
fibre.  In particular every irreducible representation is of the form
$R\circ \rho$ where $R$ is an irreducible representation of
$A(u)\rtimes S_u$ for some $u$ and $\rho$ is the canonical extension of
the restriction map.  Since $G\unit/G$ is $T_0$ we may use Proposition
\ref{prop:96} to conclude that $\Ind_{S_u}^G R$, and hence
$\Ind_S^G(R\circ \rho)$ is irreducible.  Thus $\Phi(R) = \Ind_S^G R$
is a well defined map from $(A\rtimes S)\sidehat$ into $(A\rtimes
G)\sidehat$.  Furthermore Theorem \ref{thm:ehregularity} tells us that
every irreducible representation is (equivalent to one) 
of the form $\Ind_{S_u}^G R$ so that $\Phi$ is surjective.  
Finally, we show that $\Phi$ is continuous.  This 
follows from the general theory of Rieffel induction.  In particular 
\cite[Corollary 3.35]{tfb} and the definition of $\Ind_{S}^G$ 
implies that the map 
\[
\ker R \mapsto \Ind_{S_u}^G \ker R = \ker \Ind_S^G R
\]
is continuous.  Since the topology on the spectrum of a $C^*$-algebra
is inherited from the Jacobson topology on the space of primitive
ideals it is straightforward to show that $\Phi$ must be continuous. 
\end{proof}

\subsection{Groupoid Actions}

At this point we begin to transition over to the use of abelian
isotropy subgroups, an assumption that will stick with us.  The reason
is because this assumption will allow us to identify the equivalence
classes determined by $\Phi$.  It is worth pointing out that the odd
result here or there may extend to the nonabelian case.  

\begin{prop}[{\cite[Lemma 4.1]{ctgIII}}]
Suppose $G$ is  a second countable locally compact groupoid and that
the isotropy subgroupoid $S$ is abelian and varies continuously.  
Then there is a continuous $S$-invariant homomorphism $\omega$ from
$G$ to $\R^+$ such that for all $f\in C_c(S)$
\begin{equation}
\label{eq:122}
\int_S f(s)d\beta^{r(\gamma)}(s) = \omega(\gamma)\int_S f(\gamma s
\gamma\inv)d\beta^{s(\gamma)}(s).
\end{equation}
\end{prop}

\begin{proof}
Let $\beta$ be a Haar system for $S$.  
Given $\gamma\in G$ consider the map $\phi^\gamma:S_{s(\gamma)}\rightarrow
S_{r(\gamma)}$ defined by $\phi^\gamma(s) = \gamma s \gamma\inv$.  It is
clear that $\phi^\gamma$ is a group isomorphism so that we can push forward
the Haar measure $\beta^{s(\gamma)}$ to a Haar measure on
$S_{r(\gamma)}$ defined for $f\in C_c(S)$ by 
\[
\phi^\gamma_* \beta^{s(\gamma)}(f) = \int_S f(\gamma
s\gamma\inv)d\beta^{s(\gamma)}(s).
\]
However, Haar measure is unique up to a scalar multiple so there
exists $\omega(\gamma) \in \R^+$ such that $\beta^{r(\gamma)} =
\omega(\gamma) \phi^\gamma_*\beta^{s(\gamma)}$.  It is clear that
$\omega $ is the map we are looking for.  Furthermore, it is easy to
show that if $\gamma$ and $\eta$ are composable then
$\phi^{\gamma\eta} = \phi^\gamma\circ \phi^\eta$ so that 
\[
\phi^{\gamma\eta}_*\beta^{s(\eta)} = \phi^\gamma_*\phi^\eta_*
\beta^{s(\eta)}.
\]
It follows that $\omega(\gamma\eta) = \omega(\gamma)\omega(\eta)$.
Finally, if $s\in S_{s(\gamma)}$ then, since the stabilizers are
abelian, we have $\phi^{\gamma s} = \phi^\gamma$ and $\omega$ is
$S$-invariant on the left.  A similar argument shows that it is
invariant on the right. 
 
Now we show that $\omega$ is continuous.  This portion of the proof is
taken from \cite[Lemma 4.1]{ctgIII}. Suppose to the contrary that
there exists $\gamma_n\rightarrow \gamma_0$ such that
$|\omega(\gamma_n)-\omega(\gamma_0)| \geq \epsilon > 0$ for all $n$.
We can certainly choose $f\in C_c(S)$ such that $\int_S f
d\beta^{r(\gamma_0)} = 1$.  Thus $\int_S f d\beta^{r(\gamma_n)}$ is
eventually nonzero.  We claim that we may as well assume that
$s(\gamma_n) \ne s(\gamma_0)$ for all $n > 0$.  If not, then we can
pass to a subsequence, relabel, and assume that $s(\gamma_n) =
s(\gamma_0) = u$ for all $n$.  Now suppose $\delta > 0$ and that there
exists $s_n$ such that 
\begin{equation}
\label{eq:124}
|f(\gamma_n s_n \gamma_n\inv) - f(\gamma_0s_n \gamma_0\inv)| \geq
\delta
\end{equation}
for all $n> 0$.  For \eqref{eq:124} to hold we must have either
have $\gamma_n s_n \gamma_n\inv\in \supp f$ infinitely often or
$\gamma_0 s_n \gamma_0\in \supp f$ infinitely often.  Either way we
can pass to a subnet and find $s$ such that $s_n\rightarrow s$.
However we then have 
\begin{align*}
f(\gamma_n s_n \gamma_n\inv) &\rightarrow f(\gamma_0 s \gamma_0\inv) \\
f(\gamma_0 s_n \gamma_0\inv) &\rightarrow f(\gamma_0 s \gamma_0\inv)
\end{align*}
which contradicts \eqref{eq:124}.  This shows that $f\circ
\phi^{\gamma_n}\rightarrow f\circ \phi^{\gamma_0}$ uniformly on
$C_c(S_u)$.  Next, let $U$ be some compact neighborhood of
$\gamma_0$.  Then eventually $\gamma_n \in U$ and $\supp
(f\circ\phi^{\gamma_n})$ is contained in the compact set 
\[
\{\gamma\inv s \gamma : \gamma\in U, s\in \supp f, s(\gamma) = p(s), r(\gamma)=
u\}.
\]
Thus $f\circ \phi^{\gamma_n}\rightarrow f\circ\phi^{\gamma_0}$ with
respect to the inductive limit topology and 
\[
\phi^{\gamma_n}_* \beta^u(f) = \int_S f(\gamma_n s
\gamma_n\inv)d\beta^u(s) \rightarrow \phi^{\gamma_0}_* \beta^u(f) = \
\int_S f(\gamma_0 s \gamma_0\inv) d\beta^u(s).
\]
Therefore 
\begin{equation}
\label{eq:125}
\omega(\gamma_n)\inv = (\beta^{r(\gamma_n)}(f))\inv\phi^{\gamma_n}_*\beta^{s(\gamma_n)}(f)
\rightarrow \omega(\gamma_0)\inv = (\beta^{r(\gamma_0)}(f))\inv\phi^{\gamma_n}_*\beta^{s(\gamma_0)}(f)
\end{equation}
which leads to a contradiction. 

This proves our claim so that we may assume $s(\gamma_n)\ne
s(\gamma_0)$ for all $n > 0$.  By passing to a subsequence and
relabeling we can also assume that $s(\gamma_n) \ne s(\gamma_m)$ for
all $n\ne m$.  Then $C = p\inv(\{s(\gamma_n)\}_{n=0}^\infty)$ is
closed in $S$ and we can define $\iota$ on $C$ by $\iota(s) = n$  if
and only if $p(s) = s(\gamma_n)$.  Then it is straightforward to show
that the function 
\[
F_0(s)= f(\gamma_{\iota(s)} s \gamma_{\iota(s)}\inv)
\]
is continuous and compactly supported on $C$.  Therefore we can find
an extension $F\in C_c(S)$.  However, we then have 
\[
\phi^{\gamma_n}_* \beta^{s(\gamma_n)}(f) = \int_S
F(s)d\beta^{s(\gamma_n)}(s) \rightarrow \phi^{\gamma_0}_*
\beta^{s(\gamma_0)}(f) = \int_S F(s) d\beta^{s(\gamma_0)}(s)
\]
and we obtain a contradiction just as in \eqref{eq:125}. 
\end{proof}

We can now perform a construction which is, in many ways, interesting in
its own right, even though we will only make use of it indirectly.  

\begin{prop}
\label{prop:99}
\index{dynamical system}
Suppose $(A,G,\alpha)$ is a separable dynamical system and that the
isotropy subgroupoid $S$ is abelian and varies continuously.  Then
there is an action of $G$ on $A\rtimes_\alpha S$ defined by the
collection $\{\delta_\gamma\}_{\gamma\in G}$ where, for $f\in
C_c(S_{s(\gamma)},A(s(\gamma)))$, 
\[
\delta_\gamma(f)(s) = \omega(\gamma)\inv\alpha_\gamma(f(\gamma\inv
s\gamma)).
\]
\end{prop}

\begin{proof}
First, recall that $A\rtimes S$ is a $C_0(G\unit)$-algebra so that it
makes sense to define a groupoid action.  Furthermore, as usual we will
use the restriction map to identify the fibres with $A(u)\rtimes
S_u$.  It is easy to see that $\delta_\gamma$ maps
$C_c(S_{s(\gamma)},A(s(\gamma)))$ into
$C_c(S_{r(\gamma)},A(r(\gamma)))$ and that $\delta_\gamma$ is
continuous with respect to the inductive limit topology (on both
algebras).  We will show it is a $*$-homomorphism.  Given $f,g\in
C_c(S_{s(\gamma)},A(s(\gamma)))$ we have 
\begin{align*}
\delta_\gamma(f*g)(s) &= \omega(\gamma)\inv
\alpha_\gamma(f*g(\gamma\inv s \gamma)) \\
&= \int_S \omega(\gamma)\inv
\alpha_\gamma(f(t)\alpha_t(g(t\inv\gamma\inv s
\gamma)))d\beta^{s(\gamma)}(\eta) \\
&= \int_S \omega(\gamma)^{-2}\alpha_\gamma(f(\gamma\inv t
\gamma)\alpha_{\gamma\inv t \gamma}(g(\gamma\inv t\inv s
\gamma)))d\beta^{r(\gamma)}(t) \\
&= \int_S \delta_\gamma(f)(t)\alpha_t(\delta_\gamma(g)(t\inv
s))d\beta^{r(\gamma)}(t) \\
&= \delta_\gamma(f)*\delta_\gamma(g)(s),
\end{align*}
as well as 
\begin{align*}
\delta_\gamma(f^*)(s) &= \omega(\gamma)\inv
\alpha_\gamma(f^*(\gamma\inv s \gamma)) \\
&= \omega(\gamma)\inv \alpha_\gamma(\alpha_{\gamma\inv s
  \gamma}(f(\gamma\inv s\inv \gamma)^*)) \\
&= \alpha_s(\omega(\gamma)\inv \alpha_\gamma(f(\gamma\inv s \inv
\gamma))^*) \\
&= \alpha_s(\delta_\gamma(f)(s\inv)^*) = \delta_\gamma(f)^*(s).
\end{align*}
Since $\delta_\gamma$ is a $*$-homomorphism which is continuous with
respect to the inductive limit topology Corollary \ref{cor:24}
shows that it is bounded and extends to $A(u)\rtimes S_u$.  
Next we observe that $\delta_u = \id$ and
that for composable $\gamma$ and $\eta$ 
\begin{align*}
\delta_\gamma(\delta_\eta(f))(s) &= \omega(\gamma)\inv
\omega(\eta)\inv \alpha_{\gamma}(\alpha_\eta(f(\eta\inv\gamma\inv s
\gamma \eta))) \\
&= \omega(\gamma\eta)\inv \alpha_{\gamma\eta}(f((\gamma\eta)\inv s
(\gamma\eta))) \\
&= \delta_{\gamma\eta}(f)(s).
\end{align*}
This not only shows that $\delta_\gamma$ is an isomorphism with
inverse $\delta_{\gamma\inv}$, it also shows that $\delta$ preserves
the groupoid operations.  All we need to do now is show that the
action is continuous. 

Let $\mcal{E}$ be the upper-semicontinuous bundle associated to
$A\rtimes S$.  Suppose
$\gamma_n\rightarrow \gamma_0$ and that $a_n\rightarrow a_0$ in
$\mcal{E}$ such that $s(\gamma_n) = p(a_n) = u_n$ for all $n \geq 0$.
Fix $\epsilon > 0$ and let $v_n = r(\gamma_n)$ for all $n\geq 0$. 
First, choose $b \in A\rtimes S$ such that
$b(u_0) = a_0$.  Next, using the fact that $\Gamma_c(S,p^*\mcal{A})$ is dense in
$A\rtimes S$, we can choose $f\in \Gamma_c(S,p^*\mcal{A})$ such that
$\|f - b\| < \epsilon/2$.  In particular 
\[
\|f(u) - b(u) \| < \epsilon/2
\]
for all $u\in U$ where we recall that $f(u)$ denotes the restriction of
$f$ to $S_u$.  We make the following claim. 

\begin{claim}
If $f\in \Gamma_c(S,p^*\mcal{A})$ and $\gamma_n\rightarrow
\gamma_0$ as above then $\delta_{\gamma_n}(f(u_n)) \rightarrow
\delta_{\gamma_0}(f(u_0))$ in $\mcal{E}$.  
\end{claim}
\begin{proof}[Proof of Claim.]

As in the proof that $\omega$ is continuous, we first suppose that $v_n
= v_0$ infinitely often.  Then we can pass to a
subsequence, relabel, and assume $v_n = v_0$ for all $n$.  Suppose
that we can pass to a subsequence such that for each $n> 0$ there
exists $s_n\in S_{v_0}$ such that 
\begin{equation}
\label{eq:126}
\|\delta_{\gamma_n}(f(u_n))(s_n) - \delta_{\gamma_0}(f(u_0))(s_n)\| \geq
\epsilon > 0
\end{equation}
for all $n>0$.  If this is to hold then we either must have 
$\gamma_n\inv s_n \gamma_n \in \supp f$ infinitely often or $\gamma_0
\inv s_n \gamma_0 \in \supp f$ infinitely often.  In either case we may
pass to a subsequence, multiply by the appropriate groupoid elements, and
find $s_0$ such that $s_n\rightarrow s_0$.  However we then have 
\begin{align*}
f(\gamma_n s_n \gamma_n\inv)&\rightarrow f(\gamma_0 s_0 \gamma_0\inv),
\quad\text{and}\\
f(\gamma_0 s_n \gamma_0\inv)&\rightarrow f(\gamma_0 s_0 \gamma_0\inv).
\end{align*}
Since $\alpha$ and $\omega$ are continuous, 
it follows that $\delta_{\gamma_n}(f(u_n))(s_n)$ and
$\delta_{\gamma_0}(f(u_0))(s_n)$ both converge to
$\delta_{\gamma_0}(f(u_0))(s_0)$.  This contradicts \eqref{eq:126}.  As a
result $\delta_{\gamma_n}(f(u_n))$ must converge to
$\delta_{\gamma_0}(f(u_0))$ uniformly.  Let $U$ be a compact neighborhood
of $\gamma_0$.  Eventually $\gamma_n \in U$ and therefore eventually
$\supp \delta_{\gamma_n}(f(u_n))$ is contained in the compact set 
\[
\{ \gamma\inv s \gamma  : \gamma \in U, s\in \supp f, s(\gamma) =
p(s),r(\gamma) = v_0\}.
\]
Thus $\delta_{\gamma_n}(f(u_n))\rightarrow \delta_{\gamma_0}(f(u_0))$
with respect to the inductive limit topology and thus in
$A(v_0)\rtimes S_{v_0}\subset \mcal{E}$.  

Next, suppose that we may remove an initial segment and assume that
$v_n \ne v_0$ for all $n > 0$.  Furthermore, we may pass to a
subsequence, relabel,
and assume that $v_n \ne v_m$ for all $n \ne m$.  Then $C =
p\inv(\{v_n\}_{n=0}^\infty)$ is closed in $S$ and we can define
$\iota$ on $C$ by $\iota(s) = n$ if and only if $p(s) = v_n$.  We
would like to show that 
\[
F_0(s) = \delta_{\gamma_{\iota(s)}}(f)(s) =
\omega(\gamma_{\iota(s)})\inv
\alpha_{\gamma_{\iota(s)}}(f(\gamma_{\iota(s)}\inv s
\gamma_{\iota(s)}))
\]
defines a compactly supported continuous function on $C$.  Suppose
$s_i\rightarrow s$.  If $\iota(s_i)$ is eventually constant then the
convergence of $F_0(s_i)\rightarrow F_0(s)$ is easy.  However, it is
also easy if $\iota(s_i)\rightarrow\infty$ because in this case we
just use the fact that $\gamma_i\rightarrow \gamma_0$ and that all the
different components are continuous.  Furthermore it is
straightforward to check that $F_0$ has support contained in the
compact set 
\[
\{\eta s \eta\inv : \eta\in \{\gamma_n\}_{n=0}^\infty, s\in \supp f,
s(\eta) = p(s) \}.
\]
Now, $K=\{v_n\}_{n=0}^\infty$ is an $S$-invariant closed subset of
$G\unit$ so that we may use Proposition \ref{prop:88} to conclude that
$A\rtimes S(K)$ is isomorphic to $A(K)\rtimes C$.  Thus, since 
the restriction map is a surjective homomorphism of
$A\rtimes S$ onto $A\rtimes S(K)$ by Proposition \ref{prop:71}, and since $F_0
\in\Gamma_c(C,p^*\mcal{A})$, there must be some $F\in A\rtimes S$ such
that $F(v_n) = F_0(v_n)$ for all $n \geq 0$.  In particular, since $F$
is a continuous section of $\mcal{E}$, we have $F_0(v_n)\rightarrow
F_0(v_0)$.  However we clearly constructed $F_0$ so that $F_0(v_n) =
\delta_{\gamma_n}(f(u_n))$ for all $n\geq 0$.  This proves our claim. 
\end{proof}

Thus $\delta_{\gamma_n}(f(u_n))\rightarrow \delta_{\gamma_0}(f(u_0))$.
Furthermore we have 
\[
\|\delta_{\gamma_0}(f(u_0)) - \delta_{\gamma_0}(a_0)\|=
\|f(u_0)-b(u_0)\| < \epsilon/2 < \epsilon.
\]
by construction.  Since both $a_n \rightarrow a_0$ and
$b(u_n)\rightarrow a_0$ it follows that $\|a_n-b(u_n)\| \rightarrow 0$ so that
eventually 
\[
\|\delta_{\gamma_n}(f(u_n)) - \delta_{\gamma_n}(a_n)\|\leq 
\|f(u_n) - b(u_n)\| + \|b(u_n)-a_n\| < \epsilon.
\]
It now follows from the last part of Proposition \ref{prop:35} that
$\delta_{\gamma_n}(a_n)\rightarrow \delta_{\gamma_0}(a_0)$. Hence the
action is continuous and we are done. 
\end{proof}

We get the following immediate and important corollary from
Proposition \ref{prop:69}.  Looking ahead, this corollary lays the
foundation for our identification of the equivalence classes
determined by $\Phi$.  

\begin{corr}
\label{cor:14}
Suppose $(A,G,\alpha)$ is a separable dynamical system and that the
stabilizer subgroupoid $S$ is abelian and continuously varying.  Then
the action $\delta$ induces an action of $G$ on $(A\rtimes S)\sidehat$
such that $\gamma \cdot R = R \circ \delta_\gamma\inv$ for all
$\gamma\in G$ and $R\in (A(s(\gamma))\rtimes S_{s(\gamma)})\sidehat$. 
\end{corr}

Of course, we would like to find a covariant decomposition for the
above action. 

\begin{prop}
\label{prop:100}
Suppose $(A,G,\alpha)$ is a separable dynamical system and that the
isotropy subgroupoid $S$ is abelian and continuously varying.  If $R =
\pi\rtimes U$ is a representation of $A(u)\rtimes S_u$ then
$\gamma\cdot R = \rho\rtimes V$ where 
\begin{equation}
\rho(a) = \pi(\alpha_\gamma\inv(a)),\quad\text{and}\quad
V_s = U_{\gamma\inv s\gamma}.
\end{equation}
\end{prop}

\begin{proof}
Suppose we are given $(A,G,\alpha)$ and $R$ as above with $\gamma\in
G$ such that $s(\gamma) = u$.  Let $v =
r(\gamma)$. Since
$A(v)\rtimes S_{v}$ is a group crossed product we have
a lot of technology at our disposal.  In particular we know that there
must be a covariant representation $(\rho,V)$ such that $\gamma\cdot
R= \rho\rtimes V$.  It follows from
\cite[Proposition 2.34]{tfb2} that 
$\rho = \overline{R}\circ \iota_{A(v)}$ and $V =
\overline{R}\circ\iota_{S_{v}}$ where $\iota_{A(v)}$ and
  $\iota_{S_{v}}$ are the canonical maps given by 
\begin{align*}
\iota_{A(v)}(a)f(s) &= af(s), &
\iota_{S_{v}}(s)f(t) &= \alpha_s(f(s\inv t)).
\end{align*}
We compute that 
\begin{align*}
V_s (\gamma\cdot R)(f)h &= \gamma\cdot R(\iota_{S_{v}}(s)f)h =
R(\delta_\gamma\inv(\iota_{S_v}(s)f))h \\
&= \int_S \pi(\delta_\gamma\inv(\iota_{S_v}(s)f)(t))U_t h d\beta^u(t) \\
&= \int_S \omega(\gamma)\pi(\alpha_\gamma\inv(\iota_{S_v}(s)f(\gamma
t\gamma\inv))) U_t h d\beta^u(t) \\
&= \int_S \omega(\gamma)\pi(\alpha_{\gamma\inv s}(f(s\inv \gamma t
\gamma\inv))) U_t hd\beta^u(t) \\
&= \int_S \omega(\gamma)U_{\gamma\inv s
  \gamma}\pi(\alpha_\gamma\inv(f(\gamma(\gamma\inv s\inv
\gamma)t\gamma\inv))) U_{\gamma\inv s\inv \gamma t} h d\beta^u(t).
\end{align*}
Using the fact that $\beta^u$ is left invariant we obtain
\begin{align*}
V_s(\gamma\cdot R)(f)h &= U_{\gamma\inv s \gamma} \int_S
\omega(\gamma) \pi(\alpha_\gamma\inv(f(\gamma t \gamma\inv))) U_t h
d\beta^u(t)  \\
&= U_{\gamma\inv s \gamma} \int_S \pi(\delta_\gamma\inv(f)(t))U_t h
d\beta^u(t) \\
&= U_{\gamma\inv s\gamma} R(\delta_\gamma\inv(f)) h 
= U_{\gamma\inv s \gamma}(\gamma\cdot R)(f)h.
\end{align*}
Since $\gamma\cdot R$ is nondegenerate this shows that $V_s =
U_{\gamma\inv s\gamma}$.  In the same manner we compute
\begin{align*}
\rho(a) (\gamma\cdot R)(f) h&= \gamma\cdot R(\iota_{A(v)}(a)f)h =
R(\delta_\gamma\inv(\iota_{A(v)}(a)f))h \\
&= \int_S \pi(\delta_\gamma\inv(\iota_{A(v)}(a)f(t)))U_t h d\beta^u(t)
\\
&= \int_S \omega(\gamma)\pi(\alpha_\gamma\inv(af(\gamma
t\gamma\inv)))U_t h d\beta^u(t) \\
&= \pi(\alpha_\gamma\inv(a)) \int_S
\omega(\gamma)\inv\pi(\alpha_\gamma\inv(f(\gamma t \gamma \inv)))U_t
hd\beta^v(t) \\
&= \pi(\alpha_\gamma\inv(a))R(\delta_\gamma\inv(f))h =
\pi(\alpha_\gamma\inv(a))(\gamma\cdot R)(f)h.
\end{align*}
Once again using the fact that $\gamma\cdot R$ is nondegenerate, we
conclude that $\rho(a) = \pi(\alpha_\gamma\inv(a))$.  
\end{proof}

\subsection{Equivalent Representations}

Now it is time to explore the structure of representations induced
from the stabilizers.
We need to find better ways of writing them and in particular will
find a couple of very nice equivalent representations.  This material
is at least inspired by the work done in \cite{ctgIII} when it doesn't
copy it directly.  In order to proceed we
need the following 

\begin{lemma}[{\cite[Lemma 2.1]{ctgIII}}]
\label{lem:32}
Let $G$ be a second countable locally compact Hausdorff groupoid. 
Suppose $u\in G\unit$, that $A$ is an abelian subgroup of $S_u$ and
that $\beta$ is a Haar measure on $A$.  Then the following hold. 
\begin{enumerate}
\item The formula
\[
Q(f)([\gamma]) = \int_A f(\gamma s)d\beta(s)
\]
defines a surjection from $C_c(G)$ onto $C_c(G_u/A)$.  
\item There is a non-negative, bounded, continuous function $b$ on
  $G_u$ such that for any compact set $K\subset G_u$ the support of
  $b$ and $KA$ have compact intersection and for all $\gamma\in G_u$
\begin{equation}
\label{eq:127}
\int_A b(\gamma s)d\beta(s) = 1.
\end{equation}
\item There is a Radon measure $\sigma$ with full support on $G_u/A$
  such that 
\begin{equation}
\label{eq:129}
\int_G f(\gamma)d\lambda_u(\gamma) = \int_{G_u/A}\int_A f(\gamma
s)d\beta(s)d\sigma([\gamma]).
\end{equation}
\end{enumerate}
\end{lemma}

\begin{proof}
This proof is taken (almost) verbatim from \cite{ctgIII}. The
properness of the $A$-action implies that $G_u/A$ is locally compact
Hausdorff and that $Q$ takes values in
$C_c(G_u/A)$.  The existence of a function $b'$ satisfying the
requirements of (2) with the exception of \eqref{eq:127} follows from
\cite[Lemme 1]{bourbaki}.  Now, the rest of (2) follows by normalizing
$b'$ and then the rest of (1) follows from (2).  Part (3) will follow
(except for the support statement) if we can show that the equation 
\[
\sigma(Q(f)) = \int_G f(\gamma)d\lambda_u(\gamma)
\]
yields a well-defined, positive linear functional on $C_c(G_u/A)$.
But this amounts to showing that given $f\in C_c(G)$ such that
\begin{equation}
\label{eq:128}
\int_A f(\gamma s)d\beta(s) = 0
\end{equation}
for all $\gamma\in G_u$ then $\sigma(Q(f))=0$.  However, if
\eqref{eq:128} holds then for any $h\in C_c(G)$,
\begin{align*}
\int_A h*f(s)d\beta(s) &= \int_A \int_G
h(\gamma)f(\gamma\inv s)d\lambda^u(\gamma) d\beta(s) \\
&= \int_G h(\gamma) \int_A f(\gamma\inv s)d\beta(s) d\lambda^u(\gamma)
= 0
\end{align*}
On the other hand,
\begin{align}
\int_A h*f(s)d\beta(s) &= \int_A \int_G h(\gamma)f(\gamma\inv
s)d\lambda^u(\gamma) d\beta(s) \label{eq:130} \\\nonumber
&= \int_G \int_A h(s\gamma) f(\gamma\inv)d\beta(s)d\lambda^u(\gamma)
\\\nonumber
&= \int_G \left(\int_A (\bar{h})^*(\gamma\inv s)d\beta(s)\right)
f(\gamma\inv) d\lambda^u(\gamma) 
\end{align}
where we replaced $s$ by $s\inv$ in the final equality and used the
fact that $A$ is abelian and hence unimodular.   Now consider $K =
\supp f \cap G_u$.  By part (2) $\supp b \cap KA$ is compact and therefore we
can use the Tietze Extension Theorem to extend $b$ from $\supp b\cap
KA$ to a function $d \in C_c(G)^+$.  If we let $h = d^*$ then, whenever
$\gamma\inv \in \supp f\cap G_u$, 
\[
\int_A (\bar{h})^*(\gamma\inv s)d\beta(s)=
\int_A b(\gamma\inv s) d\beta(s)= 1.
\]
It now follows from \eqref{eq:130} that $\sigma(Q(f)) = 0$.  Thus our
linear functional is well defined and the Radon measure $\sigma$
exists.  

Next we need to show that $\supp \sigma = G_u/A$.  Suppose $O$ is an
open neighborhood of $[\eta]$ in $G_u/A$.  We must show $\sigma(O)>
0$. Choose $f\in C_c(G_u/A)$ such that $0\leq f \leq 1$, $f([\eta]) =
1$, and $\supp f \subset O$.  Then we have 
\begin{align*}
\sigma(O) &\geq \int_{G_u/A} f([\gamma]) d\sigma([\gamma]) =
\int_{G_u/A}\int_A f([\gamma s])b(\gamma s)d\beta(s)d\sigma([\gamma]) \\
&= \int_G f([\gamma])b(\gamma)d\lambda_u(\gamma).
\end{align*}
Now $\gamma\mapsto f([\gamma])b(\gamma)$ is a continuous function.
Furthermore since $\int b(\eta s)d\beta(s) = 1$ there must be some
$s\in A$ such that $b(\eta s) > 0$.  However we also have $f([\eta s])
= f([\eta]) = 1 > 0$.  Since $\supp \lambda_u = G_u$ and the integrand
is continuous and nonzero on $G_u$, it follows that 
\[
\int_G f([\gamma])b(\gamma)d\lambda_u(\gamma) > 0
\]
and hence $\sigma(O)> 0$.  
\end{proof}

\begin{remark}
\label{rem:23}
\index[not]{$\sigma^u$}
Suppose $(A,G,\alpha)$ is a separable dynamical system and that the
stabilizer subgroupoid $S$ is abelian.  For all $u\in S\unit$ let
$\beta^u$ be a Haar measure on $S_u$.  Using Lemma \ref{lem:32} for each
$u\in G\unit$ there exists a Radon measure $\sigma^u$ with full support on
$G_u/S_u$ such that 
\[
\int_G f(\gamma)d\lambda_u(\gamma) = \int_{G_u/S_u}\int_S f(\gamma s)
d\beta^u(s) d\sigma^u([\gamma]).
\]
For the rest of this section whenever we have $(A,G,\alpha)$ and $S$
as above we will let $\sigma = \{\sigma^u\}$ be defined in this way.  
It is worth
mentioning that if the $\beta^u$ form a Haar system for $S$ then 
the $\sigma^u$ form a Haar system on $R_Q$ \cite[Lemma 4.2]{ctgIII}.  
\end{remark}

\begin{lemma}
\label{lem:33}
Suppose $(A,G,\alpha)$ is a separable dynamical system and that the
stabilizer subgroupoid $S$ is abelian.  Given
$u\in G\unit$ let $R=\pi\rtimes U$ be a representation of $A(u)\rtimes
S_u$ which acts on a separable Hilbert Space $\mcal{H}$.  Let
$\mcal{V}$ be the set of Borel functions $\phi:G_u\rightarrow
\mcal{H}$ such that \footnote{Since $\mcal{H}$ is
  separable we don't have to worry about the measurably
  considerations described in \cite[Appendix I.4]{tfb2}.}
\begin{equation}
\label{eq:134}
\phi(\gamma s) = U_s^* \phi(\gamma)
\end{equation}
for all $\gamma \in G_u$ and $s\in S_u$.  Define
\[
\mcal{L}^2_U(G_u,\mcal{H},\sigma^u):=\left\{ \phi\in \mcal{V}:
\int_{G_u/A_u} \|\phi(\gamma)\|^2d\sigma^u([\gamma])<\infty\right\}
\]
and let $L^2_U(G_u,\mcal{H},\sigma^u)$ be the quotient of
$\mcal{L}^2_U(G_u,\mcal{H},\sigma^u)$ where we identify functions
which agree $\lambda_u$-almost everywhere.  If $\phi,\psi\in
\mcal{L}^2_U(G_u,\mcal{H},\sigma^u)$ then 
\begin{equation}
\label{eq:135}
(\phi|\psi):= \int_{G_u/S_u} (\phi(\gamma),\psi(\gamma))
d\sigma^u([\gamma])
\end{equation}
defines an inner product which makes
$L^2_U(G_u,\mcal{H},\sigma^u)$ into a Hilbert space. 
\end{lemma}

\begin{proof}
It is clear that $L^2_U(G_u,\mcal{H},\sigma^u)$ is at least a
vector space.  The usual Cauchy-Schwartz considerations will show that
\eqref{eq:135} is integrable and 
\eqref{eq:134} guarantees that \eqref{eq:135} is well defined on
$\mcal{L}^2_U(G_u,\mcal{H})$.  
Furthermore, it is easy to see that \eqref{eq:135}
defines a sesqui-linear form.   
Let $b$ be as in part (b) of Lemma \ref{lem:32} for $S_u$.  
Suppose $\phi\in \mcal{L}^2_U(G_u,\mcal{H})$ 
and suppose $\phi$ is zero $\lambda_u$-almost everywhere. Then 
\begin{align*}
\|\phi\|^2  &= \int_{G_u/S_u} \|\phi(\gamma)\|^2 d\sigma^u([\gamma]) \\
&= \int_{G_u/S_u}\int_S b(\gamma s)\|\phi(\gamma
s)\|^2d\beta^u(s)d\sigma^u([\gamma]) \\
&= \int_{G} b(\gamma)\|\phi(\gamma)\|^2 d\lambda_u(\gamma)= 0.
\end{align*}
This suffices to show that \eqref{eq:135} is well defined on
$L^2(G_u,\mcal{H})$.  Now suppose $\|\phi\|=0$.  Then in particular
$\|\phi(\gamma)\| = 0$ for all $[\gamma]\not\in N$ where $N$ is some
$\sigma^u$-null set.  This implies that $\|\phi(\gamma)\|=0$ for all
$\gamma \in NS_u$.  However it follows from \eqref{eq:129} that $NS_u$
is $\lambda_u$-null.  Hence \eqref{eq:135} is positive definite on
$L^2_U(G_u,\mcal{H})$.  

The last thing to show is that $L^2_U(G_u,\mcal{H})$ is
complete.  This portion of the proof is inspired by \cite[Page
290]{tfb2}.  Suppose $\phi_n$ is a Cauchy sequence.  We can pass to a
subsequence, relabel and assume that 
\[
\|\phi_{n+1} - \phi_n\| < \frac{1}{2^n}
\]
for all $n$.  We define the
following extended real valued functions on $G_u$ by 
\begin{align*}
z_n(\gamma) &= \sum_{i=1}^n \|\phi_{i+1}(\gamma)-\phi_i(\gamma)\|, \\
z(\gamma) &= \sum_{i=1}^\infty \|\phi_{i+1}(\gamma)-\phi_i(\gamma)\|.
\end{align*}
Of course, $z_n$ is constant on $S_u$ orbits and factors to a Borel
map on $G_u/S_u$.  Using the triangle inequality in
$L^2(G_u/S_u,\sigma^u)$ we find 
\[
\|z_n\| \leq \sum_{i=1}^n \left(\int_{G_u/S_u}
\|\phi_{i+1}(\gamma)-\phi_i(\gamma)\|^2d\sigma^u([\gamma])\right)\poshalf =
\sum_{i=1}^n \|\phi_{i+1}-\phi_i\| \leq 1
\]
Since $\|z_n\|^2 = \int_{G_u/S_u} z_n(\gamma)^2d\sigma^u([\gamma])$ it follows from
the Monotone Convergence Theorem that 
\[
\|z\|^2 = \int_{G_u/S_u} z(\gamma)^2d\sigma^u([\gamma]) \leq 1.
\]
Hence, there is a $\sigma^u$-null set $N$ such that $[\gamma]\not\in
N$ implies that $z(\gamma)< \infty$.  In particular we can lift $N$ to
$G_u$ and get a $\lambda_u$-null set $NS_u$ such that $\gamma\not\in
NS_u$ implies 
\begin{equation}
\label{eq:31}
\sum_{i=1}^\infty \phi_{i+1}(\gamma)-\phi_i(\gamma)
\end{equation}
is absolutely convergent.  Thus \eqref{eq:31} converges to some
$\phi'(\gamma)\in \mcal{H}$ for all $\gamma\not\in NS$.  Furthermore 
\[
\phi'(\gamma) = \lim_{n\rightarrow \infty} \sum_{i=1}^n
\phi_{i+1}(\gamma) - \phi_i(\gamma) = \lim_{n\rightarrow\infty}
\phi_{n+1}(\gamma) - \phi_1(\gamma)
\]
Thus $\phi(\gamma) := \phi'(\gamma)-\phi_1(\gamma)$ satisfies
\begin{equation}
\label{eq:143}
\phi(\gamma) = \lim_{n\rightarrow\infty} \phi_n(\gamma)
\end{equation}
for all $\gamma\not\in NS_u$.  Hence $\phi_n\rightarrow \phi$ almost
everywhere and $\phi$ is a Borel function.  Now let $\phi$ be zero off $NS_u$.
Then, using \eqref{eq:143} and the fact that $NS_u$ is saturated we
find that 
\[
\phi(\gamma s) = U_s^* \phi(\gamma)
\]
for all $\gamma\in G_u$ and $s\in S_u$.  Next, given $\epsilon > 0$
there exists $M$ such that $\|\phi_n-\phi_m\| < \epsilon$ for all
$n,m\geq M$.  If $\gamma \not\in NS_u$ then 
\[
\|\phi(\gamma)-\phi_i(\gamma)\| = \lim_{n\rightarrow
  \infty}\|\phi_n(\gamma)-\phi_i(\gamma)\|.
\]
Thus, if $k\geq M$, Fatou's Lemma implies that 
\[
\|\phi-\phi_k\|^2 \leq \liminf_{n\rightarrow\infty}
\|\phi_n-\phi_k\|^2 \leq \epsilon^2
\]
Furthermore we have 
\begin{align*}
\|\phi(\gamma)\|^2 &\leq
(\|\phi(\gamma)-\phi_k(\gamma)\|+\|\phi_k(\gamma)\|)^2 \\
&\leq 3\|\phi(\gamma)-\phi_k(\gamma)\|^2 + 3\|\phi_k(\gamma)\|^2
\end{align*}
so that 
\[
\int_{G_u/S_u} \|\phi(\gamma)\|^2 d\sigma^u([\gamma]) \leq 3\|\phi-\phi_k\|^2 +
3\|\phi_k\|^2 < \infty.
\]
Thus $\phi\in\mcal{L}_U^2(G_u,\mcal{H},\sigma^u)$,
$\phi_n\rightarrow \phi$ in $L_U^2(G_u,\mcal{H},\sigma^u)$ and, to
quote \cite{tfb2}, ``this completes the proof of completeness.''
\end{proof}

The whole point of building this Hilbert space is so that we can use
it to define a representation.  

\begin{lemma}
\label{lem:31}
Suppose $(A,G,\alpha)$ is a separable dynamical system and that the stabilizer
subgroupoid $S$ is abelian.  Given $u\in G\unit$ let
$R = \pi\rtimes U$ be a representation of $A(u)\rtimes S_u$ which acts
on a separable Hilbert space 
$\mcal{H}$. 
Then $\Ind_{S_u}^G R$ is equivalent to the representation $T^R$ on
$L^2_U(G_u,\mcal{H},\sigma^u)$ defined for $f\in
\Gamma_c(G,r^*\mcal{A})$ and $\phi\in\mcal{L}^2_U(G_u,\mcal{H},\sigma^u)$ by 
\begin{equation}
\label{eq:131}
T^R(f)\phi(\gamma) = \int_G \pi(\alpha_\gamma\inv(f(\gamma \eta)))
\phi(\eta\inv) d\lambda^u(\eta).
\end{equation}
\end{lemma}

\begin{proof}
Recall from Theorem \ref{thm:induce} that $\Ind R$ acts on
$\mcal{Z}_{S_u}^G\otimes_{A(u)\rtimes S_u}\mcal{H}$ and we have ${\Ind
R(f)(z\otimes h)} = f\cdot z \otimes h$ where $f\cdot z$ is given by
\eqref{eq:ind}.  Now define $V:C_c(G_u,A(u))\odot \mcal{H}\rightarrow
\mcal{L}_U^2(G_u,\mcal{H})$ by 
\begin{equation}
\label{eq:132}
V(z\otimes h)(\gamma) = \int_S U_s \pi(z(\gamma s))hd\beta^u(s).
\end{equation}
It is straightforward to show that the integrand in \eqref{eq:132} is
jointly continuous in $\gamma$ and $s$.  However, this implies that it
is Borel on the product space.  The fact that $V(z\otimes h)$ is Borel
now follows from Fubini's Theorem (for vector integration).  
Furthermore, given $s\in S_u$ we have 
\[
V(z\otimes h)(\gamma s)  = \int_S U_t \pi(z(\gamma st))hd\beta^u(t) =
\int_S U_{s\inv}U_t \pi(z(\gamma t))hd\beta^u(t) = 
U_s^* V(z\otimes h)(\gamma).
\]
Finally, observe that $V(z\otimes h)$ is supported on the (compact) image of
$\supp z$ in $G_u/S_u$ so that 
\[
\int_{G_u/S_u} \|V(z\otimes h)(\gamma)\|^2 d\sigma^u([\gamma]) < \infty.
\] 
Thus $V(z\otimes h)$ maps $C_c(G_u,A(u))\odot \mcal{H}$ into
$\mcal{L}^2_U(G_u,\mcal{H})$.  Next, we compute 
\begin{align*}
(z\otimes h&, w\otimes k) = (R(\llangle w, z\rrangle_{A(u)\rtimes
  S_u})h,k) \\
&= \int _S (\pi(\llangle w,z\rrangle_{A(u)\rtimes S_u}(s))U_s
h,k)d\beta^u(s) \\
&= \int_S \int_G (\pi(w(\gamma s\inv)^*\alpha_s(z(\gamma)))U_s
h,k)d\lambda_u(\gamma) d\beta^u(x) \\
&= \int_S\int_{G_u/S_u}\int_S (\pi(w(\gamma t s\inv))^*
\pi(\alpha_{s}(z(\gamma
t)))U_{s}h,k)d\beta^u(t)d\sigma^u([\gamma])d\beta^u(s)\\
&= \int_{G_u/S_u}\int_S\int_S (\pi(w(\gamma s\inv))^*
\pi(\alpha_{st}(z(\gamma t)))
U_{st}h,k)d\beta^u(s)d\beta^u(t)d\sigma^u([\gamma])
\end{align*}
where we used the fact that $S_u$ is abelian to right translate by
$t$.  Continuing the computation
\begin{align*}
(z\otimes h, w\otimes k)
&= \int_{G_u/S_u}\int_S\int_S (U_s U_t \pi(z(\gamma t))h,\pi(w(\gamma
s\inv))k)d\beta^u(s)d\beta^u(t)d\sigma^u([\gamma]) \\
&= \int_{G_u/S_u}\int_S\int_S (U_t \pi(z(\gamma t))h , U_s
\pi(w(\gamma s)) k) d\beta^u(s)d\beta^u(t)d\sigma^u([\gamma])\\
&= \int_{G_u/S_u}(V(z\otimes h)(\gamma),V(w\otimes
k)(\gamma))d\sigma^u([\gamma]) \\
&= (V(z\otimes h),V(w\otimes k))
\end{align*}
where we again used the fact that $S_u$ is abelian and hence
unimodular.  
Thus $V$ is an isometry and extends to a map from
$\mcal{Z}_{S_u}^G\otimes_{A(u)\rtimes S_u} \mcal{H}$ into $L^2_U(G_u,\mcal{H})$.  
We will show that it is surjective.  Suppose
$\phi\in\mcal{L}_U^2(G_u,\mcal{H})$ is such that $(V(z\otimes
h),\phi) = 0$ for all $z\otimes h\in C_c(G_u,A(u))\odot
\mcal{H}$. It will suffice to show $\phi$ is zero $\lambda_u$-almost
everywhere.  We have 
\begin{align}
\label{eq:158}
0 &= (V(z\otimes h),\phi) = \int_{G_u/S_u} (V(z\otimes
h)(\gamma),\phi(\gamma))d\sigma^u([\gamma]) \\ \nonumber
&= \int_{G_u/S_u}\int_S (U_s \pi(z(\gamma s))h,
\phi(\gamma))d\beta^u(s)d\sigma^u([\gamma]) \\ \nonumber
&= \int_{G_u/S_u}\int_S (\pi(z(\gamma s))h, \phi(\gamma
s))d\beta^u(s)d\sigma^u([\gamma]) \\ \nonumber
&= \int_G (((\pi\circ z)\otimes h)(\gamma),\phi(\gamma))d\lambda_u(\gamma).
\end{align}
where $(\pi\circ z)\otimes h$ denotes the function $\gamma\mapsto
\pi(z(\gamma))h$.  Now, $\phi$ is probably not an element of
$L^2(G_u,\mcal{H})$.  We get around this using the following trick. 
Suppose $K\subset G_u$ is compact.  Let $\phi|_K$ be the
function obtained by letting $\phi$ be zero off $K$, and let $g\in C_c(G_u)$
be one on $K$.  Then 
\[
F([\gamma]) = \int_S g(\gamma s) d\beta^u(s)
\]
defines an element of $C_c(G_u/S_u)$.  We observe that 
\begin{align*}
\int_G \|\phi|_K(\gamma)\|^2 d\lambda_u(\gamma) &\leq
\int_G g(\gamma)\|\phi(\gamma)\|^2 d\lambda_u(\gamma) \\
&= \int_{G_u/H_u} \|\phi(\gamma)\|^2 \int_{S_u} g(\gamma s)
d\beta^u(s)d\sigma^{u}([\gamma]) \\
&\leq \|\phi\|^2 \|F\|_\infty.
\end{align*}
Thus $\phi|_K\in L^2(G_u,\mcal{H})$.  Given $z\in C_c(G_u,A(u))$ such
that $\supp z \subset K$ we conclude from \eqref{eq:158} that 
\[
0 = \int_G (((\pi\circ z)\otimes
h)(\gamma),\phi(\gamma))d\lambda_u(\gamma) = ((\pi\circ z|_K)\otimes
h, \phi|_{K})_{L^2(K,\mcal{H},\lambda_u)}.
\]
Hence $\phi|_K$ will be zero $\lambda_u$-almost everywhere if we can show
elements of the form $(\pi\circ z)\otimes h$ span a dense set in
$L^2(K,\mcal{H},\lambda_u)$.  However, we can restrict ourselves even further
and show that elements of the form 
\[
 f\otimes \pi(a) h = ((f\otimes a)\circ \pi)\otimes h
\]
span a dense set, where $a\in A$, $h\in\mcal{H}$ and $f\in C_c(K)$.
Recall that $L^2(K,\mcal{H})\cong
L^2(K)\otimes \mcal{H}$ \cite[Example 2.62]{tfb2} and that
elementary tensors span a dense set in $L^2(K,\mcal{H})$.   The result
now follows quickly once we recall that
$C_c(G_u)$ is dense in $L^2(G_u)$ because $\supp \lambda_u = G_u$ 
and $\pi(A(u))\mcal{H}$ is
dense in $\mcal{H}$ because $\pi$ is nondegenerate.   
Thus $\phi|_K$ is zero $\lambda_u$-almost
everywhere for each compact set $K\subset G_u$.  
Since $G_u$ is $\sigma$-compact this
implies that $\phi$ is zero $\lambda_u$-almost everywhere.  

The fact that $V$ is a unitary implies that there is a
representation $T^R$ of $A\rtimes G$ defined by $T^R(f) = V\Ind_{S_u}^G
R(f)V^*$.  We would like to see that $T^R$ is given by \eqref{eq:131}.
This follows from the following computation for $f\in
\Gamma_c(G,r^*\mcal{A})$ and $\phi\in
\mcal{L}^2_U(G_u,\mcal{H})$.
\begin{align*}
T^R(f)V(z\otimes h)(\gamma) &= V\Ind_{S_u}^G R(f)(z\otimes h)(\gamma)
=
V(f\cdot z\otimes h)(\gamma) \\
&=\int_S U_s\pi(f\cdot z(\gamma s))hd\beta^u(s) \\
&= \int_S \int_G U_s \pi(\alpha\inv_{\gamma s}(f(\eta))z(\eta\inv
\gamma s))hd\lambda^{r(\gamma)}(\eta)d\beta^u(s) \\
&= \int_G \int_S U_s \pi(\alpha\inv_{\gamma
  s}(f(\gamma\eta))z(\eta\inv s))h d\beta^u(s) d\lambda^u(\eta) \\
&= \int_G \int_S \pi(\alpha_\gamma\inv(f(\gamma\eta)))U_s
\pi(z(\eta\inv s))hd\beta^u(s)d\lambda^u(\eta) \\
&= \int_G \pi(\alpha_\gamma\inv(f(\gamma \eta)))V(z\otimes
h)(\eta\inv)d\lambda^u(\eta).\qedhere
\end{align*}
\end{proof}

This representation is sort of ``halfway'' to where we want to be.  In
particular we would like to work with a more standard Hilbert space.
The following observation will be crucial for this.  

\begin{remark}
\label{rem:26}
Suppose $G$ is a second countable groupoid and fix $u\in G\unit$.  The
fact that $G_u$ is second
countable implies that we can find a Borel cross section
$c:G_u/S_u\rightarrow G_u$ for the quotient map
\cite[Theorem 3.4.1]{invitation}.
Furthermore we can then define a Borel map $\delta:G_u\rightarrow S_u$
such that $\gamma = c([\gamma])\delta(\gamma)$.  We will make use of
these maps in what follows. One of the key
properties about $\delta$ that we will need is that
\[
\delta(\gamma s) = c([\gamma s])\inv (\gamma s) = c([\gamma])\inv
\gamma s = \delta(\gamma) s.
\]
\end{remark}

Using the Borel cross section $\delta$ we can transform $T^R$ into a
representation which acts on $L^2(G_u/S_u,\mcal{H},\sigma^u)$.  

\begin{prop}
\label{prop:102}
Suppose $(A,G,\alpha)$ is a separable dynamical system and that the
stabilizer subgroupoid $S$ is abelian.  Given
$u\in G\unit$ let $R=\pi\rtimes U$ be a separable representation of $A(u)\rtimes
S_u$ which acts on $\mcal{H}$.  Then $T^R$, and hence $\Ind_{S_u}^G
R$, is equivalent to the
representation $N^R$ given on $L^2(G_u/S_u,\mcal{H},\sigma^u)$ by 
\begin{equation}
\label{eq:141}
N^R(f)(\phi)([\gamma]) = \int_G
U_{\delta(\gamma)}\pi(\alpha_\gamma\inv(f(\eta)))U_{\delta(\eta\inv\gamma)}^*
\phi([\eta\inv\gamma])d\lambda^{r(\gamma)}(\eta).
\end{equation}
\end{prop}
 
\begin{proof}
We start by defining a map
$W$ on $L^2_U(G_u,\mcal{H},\sigma^u)$ by $W(\phi)([\gamma]) =
\phi(c([\gamma]))$.  Then 
\begin{align*}
(W(\phi),W(\psi)) &= \int_{G_u/S_u}
(W(\phi)([\gamma]),W(\psi)([\gamma])) d\sigma^u([\gamma]) \\
&= \int_{G_u/S_u}
(\phi(c([\gamma])),\psi(c([\gamma])))d\sigma^u([\gamma]) \\
&=
\int_{G_u/S_u}(U_{\delta(\gamma)}\phi(\gamma),U_{\delta(\gamma)}\psi(\gamma))d\sigma^u([\gamma]) 
\\
&= (\phi,\psi).
\end{align*}
Thus $W$ maps into $L^2(G_u/S_u,\mcal{H},\sigma^u)$ and is in fact
isometric.  Furthermore we can define an inverse by $W\inv(\phi)(\gamma)
:= U_{\delta(\gamma)}\phi([\gamma])$ so that $W$ is actually a unitary
map.  We now define a representation on $A\rtimes G$ by $N^R(f) =
WT^R(f)W^*$ and we show that $N^R$ has the desired form by computing, for
$f\in\Gamma_c(G,r^*\mcal{A})$
\begin{align*}
N^R(f)W\phi([\gamma]) &= WT(f)\phi([\gamma]) = T(f)\phi(c([\gamma]))
\\
&= \int_G \pi(\alpha_{c([\gamma])}\inv(f(c([\gamma])
\eta)))\phi(\eta\inv)d\lambda^u(\eta) \\
&= \int_G
\pi(\alpha_{c([\gamma])}\inv(f(\gamma\eta)))\phi(\eta\inv\delta(\gamma)\inv)d\lambda^u(\eta)
\\
&= \int_G
\pi(\alpha_{\delta(\gamma)\gamma\inv}(f(\gamma\eta)))
\phi(c([\eta\inv])\delta(\eta\inv)\delta(\gamma)\inv) d\lambda^u(\eta) \\
&= \int_G
U_{\delta(\gamma)}\pi(\alpha_\gamma\inv(f(\gamma\eta)))U_{\delta(\gamma)}^*
U_{\delta(\gamma)}U_{\delta(\eta\inv)}^* \phi(c([\eta\inv]))d\lambda_u(\eta) \\
&= \int_G U_{\delta(\gamma)}
\pi(\alpha_\gamma\inv(f(\gamma\eta)))U_{\delta(\eta\inv)}^*
W\phi([\eta\inv]) d\lambda^u(\eta) \\
&= \int_G U_{\delta(\gamma)}\pi(\alpha_\gamma\inv(f(\eta)))U_{\delta(\eta\inv\gamma)}W\phi([\eta\inv\gamma])d\lambda^{r(\gamma)}(\eta). \qedhere
\end{align*}
\end{proof}

\begin{remark}
\label{rem:24}
Before we state the next proposition we need to use some more measure
theoretic trickery.  Observe that the range map $r$ factors to
a continuous bijection $\tilde{r}$ between $G_u/S_u$ and $G\cdot u$.  Since
$G_u/S_u$ and $G\cdot u$ are both second countable locally compact
Hausdorff spaces, we cite Souslin's Theorem \cite[Theorem
3.2.3]{invitation} to conclude that $\tilde{r}$ is a Borel
isomorphism.  We use $\tilde{r}$ to push the measure $\sigma^u$
forward to a measure $\sigma^u_*$ on $G\cdot u$.  
It is clear that by identifying
$L^2(G_u/S_u,\mcal{H},\sigma^u)$ and ${L^2(G\cdot u,\mcal{H},\sigma^u_*)}$
via $\tilde{r}$ we can view $N^R$ as a representation on the latter
space.  It's easy to see that in this case its action is given by 
\[
N^R(f)(\phi)(\gamma\cdot u) = \int_G
U_{\delta(\gamma)}\pi(\alpha_\gamma\inv(f(\eta)))U_{\delta(\eta\inv\gamma)}^*
\phi(\eta\inv\gamma\cdot u)d\lambda^{r(\gamma)}(\eta).
\]
Since this identification is fairly natural, we won't make too much of
a fuss about it.  
\end{remark}

The reason we went through the effort to build $N^R$ is that, as the
next lemma demonstrates, it interfaces nicely with the multiplication
representation of $C^b(G\cdot u)$ on $L^2(G\cdot u,\mcal{H})$.  We
will be able to take advantage of this later on.  

\begin{lemma}
\label{lem:34}
Suppose $(A,G,\alpha)$ is a separable dynamical system and that the
stabilizer subgroupoid $S$ is abelian.  Let
$u\in G\unit$ and $R=\pi\rtimes U$ be a representation of $A(u)\rtimes
S_u$.  Consider the representation of $C_0(G\unit)$ on $L^2(G\cdot u,
\mcal{H}, \sigma^u_*)$ defined via
\[
N^u(f)\phi(v) = f(v)\phi(v).
\]
Furthermore, given $f\in C_0(G\unit)$ and $g\in \Gamma_c(G,r^*\mcal{A})$
define $f\cdot g(\gamma) := f(r(\gamma))g(\gamma)$.  Then
$N^u(f)N^R(G) = N^R(f\cdot g)$ for all $f\in C_0(G\unit)$ and
$g\in\Gamma_c(G,r^*\mcal{A})$.  
\end{lemma}

\begin{proof}
We discussed $\sigma_*^u$ and how to view $N^R$ as acting on
$L^2(G\cdot u, \mcal{H},\sigma^u_*)$ in Remark \ref{rem:24}.  The
representation $N^u$ is nothing more than the restriction map from
$C_0(G\unit)$ to $C^b(G\cdot u)$ composed with the usual multiplication
representation of $C^b(G\cdot u)$ on $L^2(G\cdot u,\mcal{H})$.  It is
also easy to see that if $f$ and $g$ are as above then 
$f\cdot g\in \Gamma_c(G,r^*\mcal{A})$.  For the
last statement we just compute 
\begin{align*}
N^u(f)N^R(g)\phi(\gamma\cdot u) &= f(r(\gamma))N^R(g)\phi(\gamma\cdot
u) \\
&= \int_G
f(r(\gamma))U_{\delta(\gamma)}\pi(\alpha_\gamma\inv(g(\eta)))U^*_{\delta(\eta\inv\gamma)}\phi(\eta\inv\gamma\cdot
u)d\lambda^{r(\gamma)}(\eta) \\
&= \int_G
U_{\delta(\gamma)}\pi(\alpha_\gamma\inv(f\cdot
g(\eta)))U^*_{\delta(\eta\inv\gamma)}\phi(\eta\inv\gamma\cdot
u)d\lambda^{r(\gamma)}(\eta) \\
&= N^R(f\cdot g)\phi(\gamma\cdot u).\qedhere
\end{align*}
\end{proof}

The final technical aspect we need to deal with before proving
something more interesting is to demonstrate the relationship between
$\omega$ and the $\sigma^u$.  

\begin{lemma}
\label{lem:36}
Suppose $G$ is a second countable locally compact Hausdorff groupoid
and that the stabilizer groupoid $S$ is abelian and varies
continuously.  Given $u\in G\unit$ and $\gamma\in G_u$ we have,
letting $v=\gamma\cdot u$ and $\phi$ be a Borel function on $G_v/S_v$,
\begin{equation}
\label{eq:144}
\int_{G_v/S_v}\omega(\gamma)\phi([\eta\gamma])d\sigma^v([\eta]) = 
\int_{G_u/S_u}\phi([\eta])d\sigma^u([\eta]).
\end{equation}
\end{lemma} 

\begin{proof}
Since $\sigma^u$ and $\sigma^v$ are Radon
measures, it suffices to verify \eqref{eq:144} for $f\in
C_c(G_u/S_u)$.  Using the fact that $\eta
s\gamma = \eta \gamma (\gamma\inv s \gamma)$ is not difficult to show
that the map $[\eta]\mapsto [\eta\gamma]$ defines a homeomorphism from
$G_v/S_v$ onto $G_u/S_u$.  Hence $[\eta]\mapsto f([\eta\gamma])$ is
continuous and compactly supported.  
We let $b$ be as in Lemma \ref{lem:32} and compute 
\begin{align*}
\int_{G_u/S_u}f([\eta])d\sigma^u([\eta]) &= \int_{G_u/S_u}\int_S
f([\eta s])b(\eta s)d\beta^u(s) d\sigma^u([\eta]) \\
&= \int_{G} f([\eta])b(\eta)d\lambda_u(\eta) 
= \int_{G} f([\eta\gamma])b(\eta\gamma)d\lambda_v(\eta) \\
&= \int_{G_v/S_v}\int_S f([\eta s\gamma])b(\eta s \gamma )d\beta^v(s)d\sigma^v([\eta])\\
&= \int_{G_v/S_v}\int_S \omega(\gamma)f([\eta\gamma s])b(\eta\gamma
s)d\beta^u(s) d\sigma^v([\eta]) \\
&= \int_{G_v/S_v} \omega(\gamma)f([\eta\gamma])
d\sigma^v([\eta]).\qedhere
\end{align*}
\end{proof}

We can now prove the following proposition, which tells us that the
equivalence classes on $S$ induced by $\Phi$ are exactly the orbits of
the $G$ action.  

\begin{prop}
\label{prop:104}
Suppose $(A,G,\alpha)$ is a separable dynamical system and that the
isotropy subgroupoid $S$ is abelian and continuously varying.  Let
$u\in G\unit$ and $R$ be an irreducible representation of $A(u)\rtimes S_u$ on a
separable Hilbert space $\mcal{H}$. Then
$\Phi(R)$ is equivalent to $\Phi(\gamma\cdot R)$ for all
$\gamma\in G_u$.   Furthermore if $G\unit/G$ is $T_0$ and $L$ and $R$
are irreducible representations of $A(u)\rtimes S_u$ and $A(v)\rtimes
S_v$, respectively, 
then $\Phi(L)$ is equivalent to $\Phi(R)$ if and only if there
exists $\gamma\in G$ such that $\gamma\cdot L$ is equivalent to $R$.  
\end{prop}

\begin{proof}
Let $v = \gamma\cdot u$.  Suppose 
$R = \pi\rtimes U$ and $\gamma\cdot R = \rho\rtimes V$ as in
Proposition \ref{prop:100}.  Recall that $\Phi(R) = \Ind_{S_u}^G R$ is
equivalent to the representation $T^R$ on
$L^2_U(G_u,\mcal{H},\sigma^u)$ and $\Phi(\gamma\cdot R) =
\Ind_{S_v}^G\gamma\cdot R$ is equivalent to the representation
$T^{\gamma\cdot R}$ on $L^2_V(G_v,\mcal{H},\sigma^v)$.  We define a
map $W$ on $\mcal{L}^2_U(G_u,\mcal{H},\sigma^u)$ by 
\[
W(\phi)(\eta) = \omega(\gamma)\poshalf f(\eta\gamma)\quad
\text{for all $\eta\in G_v$.}
\]
Clearly $W(\phi)$ is Borel, and we compute for $s\in S_v$
\begin{align*}
W(\phi)(\eta s) &= \omega(\gamma)\poshalf \phi(\eta s \gamma) = 
\omega(\gamma)\poshalf \phi(\eta \gamma (\gamma\inv s \gamma)) \\
&= \omega(\gamma)\poshalf U_{\gamma\inv s \gamma}^* \phi(\eta\gamma) 
= V_s^* W(\phi)(\eta).
\end{align*}
Furthermore, we use Lemma \ref{lem:36} to conclude that 
\begin{align*}
(W(\phi),W(\psi)) &= \int_{G_v/S_v} (W(\phi)(\eta),W(\psi)(\eta))
d\sigma^v([\eta]) \\ 
&= \int_{G_v/S_v} \omega(\gamma)
(\phi(\eta\gamma),\psi(\eta\gamma))d\sigma^v([\eta]) \\
&= \int_{G_u/S_u} (\phi(\eta),\psi(\eta)) d\sigma^u([\eta]) =
(\phi,\psi).
\end{align*}
This calculation proves two things.  First, that $W(\phi)$ is in
$\mcal{L}_V^2(G_v,\mcal{H},\sigma^v)$ and, second, that $W$ is
isometric.  Since $W$ has an obvious inverse it must be a unitary
map.  

Next we show $W$ intertwines $T^R$ and $T^{\gamma\cdot R}$.  We
see for $f\in\Gamma_c(G,r^*\mcal{A})$ that 
\begin{align*}
WT^R(f)\phi(\eta) &= \omega(\gamma)\poshalf T^R(f)\phi(\eta\gamma)
\\
&= \int_G \omega(\gamma)\poshalf
\pi(\alpha_{\eta\gamma}\inv(f(\eta\gamma\zeta)))\phi(\zeta\inv)d\lambda^u(\zeta)
\\
&=\int_G \omega(\gamma)\poshalf
\pi(\alpha_\gamma\inv(\alpha_\eta\inv(f(\eta\zeta))))\phi(\zeta\inv\gamma)d\lambda^v(\zeta)
\\
&= \int_G
\rho(\alpha_\eta\inv(f(\eta\zeta)))W\phi(\zeta\inv)d\lambda^v(\zeta)
\\
&= T^{\gamma\cdot R}(f)W\phi(\eta).
\end{align*}

Moving on, suppose $G\unit/G$ is $T_0$ and that we are given two
irreducible representations $L$ and $R$ of $A(u)\rtimes S_u$ and
$A(v)\rtimes S_v$, respectively.  Suppose $\Phi(R)$ is equivalent to
$\Phi(L)$. Proposition \ref{prop:102} implies that $N^L$ is equivalent
to $N^R$.  Let $U$ be the intertwining unitary and let $N^u$ and $N^v$
be as in Lemma \ref{lem:34}.  Then, given $f\in C_0(G\unit)$ and $g\in
\Gamma_c(G,r^*\mcal{A})$, we have 
\begin{align*}
U N^v(f)N^R(g)h &= UN^R(f\cdot g)h = N^L(f\cdot g)Uh \\
&= N^u(f)N^L(g)Uh = N^u(f)UN^R(g)h.
\end{align*}
Since $N^R$ is nondegenerate this implies that $N^v$ is unitarily
equivalent to $N^u$.  However, if $G\cdot u\cap G\cdot v = \emptyset$
then \cite[Lemma 4.15]{primtrangroup} says that $N^v$ and $N^u$ can
have no equivalent subrepresentations. Therefore we must have $G\cdot
u = G\cdot v$. So let $\gamma\in G$ be such that $\gamma\cdot u = v$.
Then $\gamma\cdot L$ and $R$ are both irreducible representation of
$A(v)\rtimes S_v$ and we assumed that that $\Phi(R)$ is
equivalent to $\Phi(L)$ which is in turn  equivalent to
$\Phi(\gamma\cdot L)$ by the above.  
It now follows from Proposition \ref{prop:96}
that $\gamma\cdot L$ is equivalent to $R$.  
\end{proof}

\subsection{Restricting Representations}

Now that we know which representations have the same image under
$\Phi$ it is time to try and show that $\Phi$ is open.  The key
construction is a restriction process from $A\rtimes G$ to $A\rtimes
S$.  This is defined using the following map.  

\begin{prop}
\label{prop:101}
Suppose $(A,G,\alpha)$ is a separable dynamical system and the
stabilizer subgroupoid $S$ is abelian and continuously varying.  Then
there is a nondegenerate homomorphism $M:A\rtimes S\rightarrow
M(A\rtimes G)$ such that 
\begin{equation}
M(f)g(\gamma) = \int_S
f(s)\alpha_s(g(s\inv\gamma))d\beta^{r(\gamma)}(s)
\end{equation}
for $f\in \Gamma_c(S,p^*\mcal{A})$ and $g\in
\Gamma_c(G,r^*\mcal{A})$.  
\end{prop}

\begin{proof}
Since $M(f)g$ is basically given by convolution, it is straightforward
to show that $M(f)g\in\Gamma_c(G,r^*\mcal{A})$ and we will not detail
a proof here.  Instead, we show that
$M_f$ is adjointable and $A\rtimes G$-linear on
$\Gamma_c(G,r^*\mcal{A})$.  First
\begin{align*}
M(f)(g*h)(\gamma) &= \int_S\int_G
f(s)\alpha_s(g(\eta)\alpha_\eta(h(\eta\inv s\inv\gamma)))
d\lambda^{r(\gamma)}(\eta)d\beta^{r(\gamma)}(s) \\
&= \int_S \int_G f(s) \alpha_s(g(\eta)) \alpha_{s\eta}(h(\eta\inv
s\inv \gamma)) d\lambda^{r(\gamma)}(\eta) d\beta^{r(\gamma)}(s) \\
&= \int_S \int_G f(s) \alpha_s(g(s\inv\eta)) \alpha_{\eta}(h(\eta\inv
\gamma)) d\lambda^{r(\gamma)}(\eta) d\beta^{r(\gamma)}(s) \\
&= \int_G M(f)g(\eta) \alpha_\eta(h(\eta\inv\gamma))
d\lambda^{r(\gamma)}(\eta) \\
&= (M(f)g)*h(\gamma).
\end{align*}
Next we compute 
\begin{align*}
(M&(f)g)^**h(\gamma) = \int_G
\alpha_\eta(M(f)g(\eta\inv)^*h(\eta\inv\gamma))d\lambda^{r(\gamma)}(\eta)
\\
&= \int_G \int_S \alpha_\eta(\alpha_s(g(s\inv\eta\inv)^*)f(s)^*
h(\eta\inv \gamma))d\beta^{s(\eta)}(s)d\lambda^{r(\gamma)}(\eta) \\
&= \int_G \int_S \omega(\eta\inv) \alpha_\eta(\alpha_{\eta\inv s\eta}(
g(\eta\inv s\inv)^*)f(\eta\inv s \eta)^* h(\eta\inv\gamma))
d\beta^{r(\gamma)}(s)d\lambda^{r(\gamma)}(\eta) \\
&= \int_S \int_G \omega(\eta\inv s)\alpha_{s\inv\eta}(\alpha_{\eta\inv
  s\eta}(g(\eta\inv)^*)f(\eta\inv s \eta)^*h(\eta\inv s
\gamma))d\lambda^{r(\gamma)}(\eta)d\beta^{r(\gamma)}(s) \\
&= \int_G \int_S\omega(\eta\inv) \alpha_\eta(g(\eta\inv)^*)
\alpha_{s\inv \eta}(f(\eta\inv s \eta)^*h(\eta\inv s
\gamma))d\beta^{r(\gamma)}(\eta)d\lambda^{r(\gamma)}(\eta) \\
&= \int_G\int_S g^*(\eta) \alpha_{\eta s\inv}(f(s)^*h(s\eta\inv
\gamma))d\beta^{s(\eta)}(s)d\lambda^{r(\gamma)}(\eta) \\
&= \int_G\int_S g^*(\eta) \alpha_{\eta} (f^*(s)\alpha_s(h(s\inv
\eta\inv \gamma)))d\beta^{s(\eta)}(s)d\lambda^{r(\gamma)}(\eta) \\
&= g^* *(M(f^*)h)(\gamma)
\end{align*}
where we have used the fact that $\omega$ is $S$-invariant and $S_u$
is unimodular.  Finally we show that $M$ preserves convolution on
$\Gamma_c(S,p^*\mcal{A})$ by calculating
\begin{align*}
M(f*g)h(\gamma) &= \int_S\int_S f(t)\alpha_t(g(t\inv s))
\alpha_s(h(s\inv \gamma))d\beta^{r(\gamma)}(t)d\beta^{r(\gamma)}(s) \\
&= \int_S f(t) \alpha_t(g(s))\alpha_{ts}(h(s\inv
t\inv\gamma))d\beta^{r(\gamma)}(s) d\beta^{r(\gamma)}(t) \\
&= \int_S f(t) \alpha_t(g(s)\alpha_s(s\inv t\inv\gamma))
d\beta^{r(\gamma)}(s) d\beta^{r(\gamma)}(t) \\
&= M(f)M(g)h(\gamma).
\end{align*}

Moving on, we show that elements of the form $M(f)g$ are dense in
$\Gamma_c(G,r^*\mcal{A})$ with respect to the inductive limit
topology.  As in Section \ref{sec:indreps}, this argument will be a
lengthy one.  Fix $\epsilon > 0$ and suppose $g\in
\Gamma_c(G,r^*\mcal{A})$.  Let $K = r(\supp g)$ and choose some fixed
open neighborhood $U$ of $K$ in $S$.  We make the
following claim.
\begin{claim}
There is a relatively compact open neighborhood $O$ of $K$ in $S$
such that $O\subset U$ and for all $\gamma\in G$ and $s\in O$
\begin{equation}
\label{eq:136}
\|\alpha_s(g(s\inv \gamma)) - g(\gamma)\| < \epsilon/2.
\end{equation}
\end{claim}
\begin{proof}[Proof of Claim.]
Suppose not.  Then for every relatively compact neighborhood $W\subset
U$ of $K$ there exists $\gamma_W\in G$ and $s_W \in W$ such that 
\begin{equation}
\label{eq:137}
\|\alpha_{s_W}(g(s_W\inv \gamma_W)) - g(\gamma_W)\| \geq \epsilon/2.
\end{equation}
When we order $W$ by reverse inclusion the sets $\{\gamma_W\}$ and
$\{s_W\}$ form nets in $G$ and $S$ respectively.  In order for
\eqref{eq:137} to hold we must have either $s_W\inv \gamma_W \in\supp g$
or $\gamma_W\in\supp g$ for each $W$.  In either case we have
$r(\gamma_W)\in K$ and, since $W$ is a neighborhood of $K$, 
$\gamma_W \in W\supp g\subset \overline{U}\supp g$.  
Furthermore, $s_W\in W\subset \overline{U}$ for all $W$.  
Since $\overline{U}$ and
$\overline{U}\supp g$ are compact, we can pass to a subnet, twice,
relabel, and find $s\in S$ and $\gamma\in G$ such that 
 $s_W\rightarrow s$ and $\gamma_W\rightarrow
\gamma$.  However, since $s_W$ is eventually in every neighborhood of
$K$ we must have $s\in K \subset G\unit$.  This implies that $s_W\inv
\gamma_W \rightarrow \gamma_W$.  However, using the continuity of the
action, this contradicts \eqref{eq:137}.  
\end{proof}

Let $O$ be the open set from above and let $f\in C_c(S)^+$ such that
$\supp f \subset O$ and $f$ is one on $K$.  Then the function
\[
u\mapsto c(u):=\int_S f(s)d\beta^u(s)
\]
is continuous and nonzero on $K$.  In particular,
the function $1/c$ is continuous on $K$ so that we may use the Tietze
Extension Theorem to find $\tilde{c}\in C_c(G\unit)$
which extends $1/c$ off $K$.  Then we can replace $f$ by
$(\tilde{c}\circ p)f$ and assume without loss of generality that 
\[
\int_S f(s)\beta^u(s) = 1
\]
for all $u\in K$.  Next, let $\{a_l\}$ be an approximate identity for $A$.
We make the following claim. 
\begin{claim}
There exists $l_0$ such that 
\begin{equation}
\label{eq:138}
\|a_{l_0}(r(\gamma))\alpha_s(g(s\inv \gamma)) - \alpha_s(g(s\inv
\gamma))\| < \epsilon /2 
\end{equation}
for all $s\in \supp f$ and $\gamma \in G$.  
\end{claim}
\begin{proof}[Proof of Claim.]
Suppose not.  Then for each $l$ there exists $\gamma_l\in G$ and
$s_l\in\supp f$ such
that
\begin{equation}
\label{eq:139}
\|a_l(r(\gamma_l))\alpha_{s_l}(g(s_l\inv \gamma_l)) -
\alpha_{s_l}(g(s_l\inv\gamma_l))\| \geq \epsilon /2.
\end{equation}
However, in order for \eqref{eq:139} to hold we must have
$s_l\inv\gamma_l\in \supp g$ for all $l$. But then $\gamma_l \in
(\supp f)\inv\supp g$.  Since both this set and $\supp f$ are
compact we can pass to two subnets, relabel, and find $\gamma\in G$
and $s\in S$ such that
$\gamma_l\rightarrow \gamma$ and $s_l\rightarrow s$.  However we now
have $\alpha_{s_l}(g(s_l\inv\gamma_l))\rightarrow
\alpha_s(g(s\inv\gamma))$.  Choose $b\in A$ such that
$b(r(\gamma)) = \alpha_s(g(s\inv \gamma))$.  Then $a_l
b \rightarrow b$ and $b(r(\gamma_l))\rightarrow b(r(\gamma))$.
Since $\alpha_{s_l}(g(s_l\inv\gamma_l))\rightarrow b(r(\gamma))$, 
we must have 
\[
\|\alpha_{s_l}(g(s_l\inv\gamma_l))-b(r(\gamma_l))\| \rightarrow 0.
\]
Putting everything together, it follows that, eventually,
\begin{align*}
\|a_l(r(\gamma_l))\alpha_{s_l}(g(s_l\inv\gamma_l)) -
\alpha_{s_l}(g(s_l\inv\gamma_l))\| \leq&
\|a_l(r(\gamma_l))\alpha_{s_l}(g(s_l\inv\gamma_l)) -
a_l(r(\gamma_l))b(r(\gamma_l)) \| \\
&+\|a_l(r(\gamma_l))b(r(\gamma_l)) + b(r(\gamma_l))\|\\
&+\|b(r(\gamma_l)) - \alpha_{s_l}(g(s_l\inv\gamma_l))\| \\
\leq& 2\|\alpha_{s_l}(g(s_l\inv \gamma_l)) -
b(r(\gamma_l))\| + \|a_l b - b\| \\
<&\epsilon /2
\end{align*}
and this contradicts \eqref{eq:139}.
\end{proof}

Now consider $f\otimes a_{l_0}\in \Gamma_c(S,p^*\mcal{A})$.  First
observe that $\supp f\otimes a_{l_0} \subset U$ and that $U$ was
chosen independently of $\epsilon$.  Next, given $\gamma\in G$ if
$r(\gamma)\not\in K$ then $g(s\gamma) = 0$ for all $s\in
S_{r(\gamma)}$ so that in particular 
\[
M(f\otimes a_{l_0})g(\gamma) - g(\gamma) = \int_S f(s)a_{l_0}(r(\gamma))
\alpha_s(g(s\inv \gamma))d\beta^{r(\gamma)}(s) = 0.
\]
If $r(\gamma)\in K$ then 
\begin{align*}
\|M(f\otimes a_{l_0})&g(\gamma) - g(\gamma)\| \\
=& \left\| \int_S
  f(s)a_{l_0}(r(\gamma))\alpha_s(g(s\inv\gamma))d\beta^{r(\gamma)}(s)
  - \int_S f(s)d\beta^{r(\gamma)}(s) g(\gamma)\right\| \\
\leq& \int_S
  f(s)\|a_{l_0}(r(\gamma))\alpha_s(g(s\inv\gamma))-g(\gamma)\|d\beta^{r(\gamma)}(s)
  \\
\leq& \int_S f(s)
\|a_{l_0}(r(\gamma))\alpha_s(g(s\inv\gamma))-\alpha_s(g(s\inv\gamma))\|
d\beta^{r(\gamma)}(s) \\
&+ \int_S f(s) \|\alpha_s(g(s\inv\gamma)) - g(\gamma)\|
d\beta^{r(\gamma)}(s) \\
&< \epsilon/2 + \epsilon/2 = \epsilon.
\end{align*}
Hence $\|M(f\otimes a_{l_0})g - g\|_\infty < \epsilon$.  This suffices
to show that elements of the form $M(f)g$ are dense in
$\Gamma_c(G,r^*\mcal{A})$ with respect to the inductive limit
topology.

Next, we want to show that $M(f)$ is bounded so that it extends to a
multiplier on $A\rtimes G$. Let $\rho$ be a state on $A\rtimes G$ and define
an inner product on $\Gamma_c(G,r^*\mcal{A})$ via 
\[
(f,g)_\rho = \rho(\langle f, g\rangle)
\]
where we give $A\rtimes G$ its usual inner-product as an $A\rtimes G$-module.
This is clearly sesqui-linear and is positive because states
are positive.  Let $\mcal{H}_\rho$ be the Hilbert space completion of
$\Gamma_c(G,r^*\mcal{A})$ with respect to this pre-inner product.
Just as in the proof of Proposition \ref{prop:82}, we
would like to show that we can apply Theorem \ref{thm:disintigration}
when $\mcal{H}_0$ is the image of $\Gamma_c(G,r^*\mcal{A})$ in
$\mcal{H}_\rho$.  Define $\pi$ on $\mcal{H}_0$ by 
\[
\pi(f)g = M(f)g
\]
for $f\in \Gamma_c(S,p^*\mcal{A})$ and $g\in
\Gamma_c(G,r^*\mcal{A})$.  If $(g,h)_\rho = 0$ for all $h\in
\Gamma_c(G,r^*\mcal{A})$ then 
\[
(\pi(f)g,h)_\rho = \rho((M(f)g)^**h) = \rho(g^**M(f^*)h) = (g,M(f^*)h)_\rho
= 0.
\]
Thus, $\pi$ is well defined and it is
clear that $\pi$ is a homomorphism from $\Gamma_c(S,p^*\mcal{A})$
to the algebra of operators on $\mcal{H}_0$.  Next, elements of the
form $M(f)g$ are dense in $\Gamma_c(G,r^*\mcal{A})$ with respect to
the inductive limit topology and therefore with respect to the norm
topology as well.  It follows that elements of the form $\pi(f)g$ 
are dense in $\mcal{H}_\rho$.
Fix $g,h\in \Gamma_c(G,r^*\mcal{A})$.  We would like to see that
$f\mapsto (\pi(f)g,h)_\rho$ is continuous with respect to the inductive
limit topology.  It suffices to see that the map $f\mapsto M(f)g$ is
continuous with respect to the inductive limit topology.  Suppose
$f_i\rightarrow f$ uniformly and $\supp f_i\subset K$ for some
compact set $K$.  Then 
\begin{align*}
\|M(f_i)g(s) - M(f)g(s)\| &\leq \int_S
\|f_i(s)-f(s)\|\|\alpha_s(g(s\inv\gamma))\|d\beta^{r(\gamma)}(s) \\
&\leq \int_S \|f_i(s)-f(s)\|\|g\|_\infty d\beta^{r(\gamma)}(s) \\
&\leq \|f_i-f\|_\infty\|g\|_\infty\beta^{r(\gamma)}(K)
\end{align*}
Since $\{\beta^u(K)\}$ is bounded this shows that
$M(f_i)g\rightarrow M(f_i)g$ uniformly.  Furthermore it is easy to
see that $\supp M(f_i) g \subset K\supp g$ so that $M(f_i)g\rightarrow
M(f)g$ with respect to the inductive limit topology.  Finally, the
fact that $(\pi(f)g,h)_\rho = (g,\pi(f^*)h)_\rho$ follows immediately from the
fact that $(M(f)g)^**h = g^**(M(f^*)h)$.  Thus, it follows from Theorem
\ref{thm:disintigration} that $\pi$ extends to a representation of
$A\rtimes G$.  In particular we have 
\[
\rho(\langle M(f)g,M(f)g\rangle) = (\pi(f)g,\pi(f)g)_\rho \leq \|f\|^2
(g,g)_\rho \leq \|f\|^2 \|g\|^2.
\]
By choosing $\rho$ so that $\rho(\langle M(f)g,M(f)g\rangle) =
\|M(f)g\|^2$ we conclude that $\|M(f)g\|\leq \|f\|\|g\|$.  Thus $M(f)$
is bounded and is adjointable with adjoint $M(f^*)$ 
on a dense subspace and therefore
extends to a multiplier on $A\rtimes G$.  Furthermore $\|M(f)\|\leq
\|f\|$ so that $M$ also extends to all of $A\rtimes S$.  We
have shown that $M$ is a homomorphism on a dense subspace so it is a
homomorphism everywhere.  Finally, the fact that $M$ is nondegenerate
follows from the fact that elements of the form $M(f)g$ are dense in
the inductive limit topology. 
\end{proof}

The point is that nondegenerate maps into multiplier algebras yield
continuous restriction processes through the usual general nonsense.  

\begin{corr}
\label{cor:15}
Suppose $(A,G,\alpha)$ is a separable dynamical system and that the
stabilizer subgroupoid $S$ is abelian and continuously varying.  
Then there exists a restriction map $\Res_M:\mcal{I}(A\rtimes G)\rightarrow
\mcal{I}(A\rtimes S)$ such that $\Res_M$ is continuous and is
characterized by $\Res_M(\ker \pi) = \ker \overline{\pi}\circ M$ for all
representations $\pi$ of $A\rtimes G$. 
\end{corr}

\begin{proof}
We proved in Proposition \ref{prop:101} that $M:A\rtimes S\rightarrow
M(A\rtimes G)$ is a nondegenerate homomorphism.  Therefore it follows
as in the latter half of \cite[Page 61]{tfb} that $M$ defines the
required continuous restriction map.  
\end{proof}

One might be tempted into thinking that restriction and induction are
dual, or inverse, in some sense.  The next lemma shows that this is not
the case.  

\begin{lemma}
\label{lem:35}
Suppose $(A,G,\alpha)$ is a separable dynamical system and that the
stabilizer subgroupoid $S$ is abelian and continuously varying.  Then
given $u\in G\unit$ and an irreducible representation $R=\pi\rtimes U$ of 
$A(u)\rtimes S_u$ we have 
\begin{equation}
\label{eq:140}
\Res_M \ker \Ind_{S_u}^G R = \bigcap_{\gamma\in G_u} \ker(\gamma\cdot
R).
\end{equation}
\end{lemma}

\begin{proof}
We know from Proposition \ref{prop:102} that $\Ind_{S_u}^G R$ is
equivalent to $N^R$ acting on $L^2(G_u/S_u,\mcal{H},\sigma^u)$ via
\eqref{eq:141}.  Let $Q = \overline{N^R}\circ M$ so that 
\[
\Res_M \ker \Ind_{S_u}^G R = \Res_M \ker N^R = \ker Q.
\]
Given $\gamma\in G_u$ recall that we can decompose $\gamma\cdot R$ as 
$\pi^\gamma\rtimes U^\gamma$ where $\pi^\gamma$ and $U^\gamma$ are
given in Proposition \ref{prop:100}.  Furthermore, we will leave it to
the reader to keep track of when we are treating $R$ as a representation
of $A(u)\rtimes S_u$ or of $A\rtimes S$.  

Given $f\in A\rtimes S$ it is straightforward to show that the
collection $\{c([\gamma])\cdot R(f)\}$ is a Borel
field of operators on the trivial bundle $G_u/S_u\times \mcal{H}$.
Use Proposition \ref{prop:41} to form the direct integral representation
$\int_{G_u/S_u}^\oplus c([\gamma])\cdot R\,
d\sigma^u([\gamma])$.  \index{direct integral} We can then compute for
$f\in \Gamma_c(S,p^*\mcal{A})$, $g\in \Gamma_c(G,r^*\mcal{A})$ and
$\phi \in \mcal{L}^2(G_u/S_u,\mcal{H},\sigma^u)$ that 
\begin{align*}
Q&(f)N^R(g)\phi([\gamma]) = N^R(M(f)g)\phi([\gamma]) \\
&=\int_G
U_{\delta(\gamma)}\pi(\alpha_\gamma\inv(M(f)g(\eta)))U^*_{\delta(\eta\inv\gamma)}\phi([\eta\inv\gamma])
d\lambda^{r(\gamma)}(\eta) \\
&= \int_G \int_S
U_{\delta(\gamma)}\pi(\alpha_\gamma\inv(f(s)\alpha_s(g(s\inv\eta))))U^*_{\delta(\eta\inv\gamma)}\phi([\eta\inv\gamma])d\beta^{r(\gamma)}(s)
d\lambda^{r(\gamma)}(\eta) \\
&= \int_S \int_G 
U_{\delta(\gamma)}\pi(\alpha_\gamma\inv(f(s)\alpha_s(g(\eta))))U^*_{\delta(\eta\inv
  s\inv \gamma)}\phi([\eta\inv s\inv
\gamma])d\lambda^{r(\gamma)}(\eta)d\beta^{r(\gamma)}(s) \\
&= \int_S \int_G
U_{\delta(\gamma)}\pi(\alpha_\gamma\inv(f(s)\alpha_s(g(\eta))))U^*_{\gamma\inv
  s\inv \gamma} U^*_{\delta(\eta\inv\gamma)}\phi([\eta\inv \gamma ])d\lambda^{r(\gamma)}(\eta) d\beta^{r(\gamma)}(s)\\
&=\int_S\int_G
U_{\delta(\gamma)}\pi(\alpha_\gamma\inv(f(s)))U_{\gamma\inv
  s\gamma}\pi(\alpha_\gamma\inv(g(\eta)))U_{\delta(\eta\inv \gamma)}^*
\phi([\eta\inv\gamma])d\lambda^{r(\gamma)}(\eta)d\beta^{r(\gamma)}(s)
\\
&= 
\int_S \int_G
\pi(\alpha_{\delta(\gamma)\gamma\inv}(f(s)))U_{\delta(\gamma)\gamma\inv
  s \gamma\delta(\gamma)\inv} U_{\delta(\gamma)}
\pi(\alpha_\gamma\inv(g(\eta))) \ldots\\
&\hspace{3in} \ldots U_{\delta(\eta\inv\gamma)}^*\phi([\eta\inv\gamma])
d\lambda^{r(\gamma)}(\eta)d\beta^{r(\gamma)}(s)
\\
&= \int_S \pi(\alpha_{c([\gamma])}\inv(f(s)))U_{c([\gamma])\inv s
  c([\gamma])} N^R(g)\phi([\gamma]) d\beta^{r(\gamma)}(s) \\
&= \pi^{c([\gamma])}\rtimes U^{c([\gamma])}(f)N^R(g)\phi([\gamma])
= c([\gamma])\cdot R(f) (N^R(g)\phi([\gamma]))\\
&=\int_{G_u/S_u}^\oplus 
c([\eta])\cdot R(f)\, d\sigma^u([\eta]) N^R(g)\phi([\gamma]).
\end{align*}
Since $N^R$ is nondegenerate, this implies that $Q = \int_{G_u/S_u}^\oplus
c([\gamma])\cdot R\, d\sigma^u([\gamma])$ and that
\begin{equation}
\label{eq:133}
Q(f)\phi([\gamma]) = (c([\gamma])\cdot R)(f)\phi([\gamma])
\end{equation}
for all $f\in A\rtimes S$ and $\phi\in \mcal{L}^2(G_u/S_u,\mcal{H})$.  
Now suppose $f\in A\rtimes S$ and $Q(f) =
0$.  Let $g_i\in C_c(G_u/S_u)$ be a countable set of functions which
separate points and let $h_j$ be a countable basis for $\mcal{H}$.  
Then for each $g_i$ and $h_j$ we have 
\begin{equation}
\label{eq:142}
(c([\gamma])\cdot R)(f)(g_i\otimes h_j)([\gamma]) = 
g_i([\gamma])(c([\gamma])\cdot R)(f)h_j = 0
\end{equation}
for all $[\gamma]\not\in N_{ij}$ where $N_{ij}$ is a $\sigma^u$-null
set.  Let $N = \bigcup_{ij} N_{ij}$ and observe that given
$[\gamma]\not\in N$ \eqref{eq:142} holds for all $i$ and $j$.  In
particular, we can pick $g_i$ so that $g_i([\gamma])\ne 0$ and conclude
that $(c([\gamma])\cdot R)(f) = 0$.  Thus $(c([\gamma])\cdot R)(f) =
0$ for all $[\gamma]\not\in N$.  Now consider $NS_u$.  Since $N$ is a
$\sigma^u$-null set it follows from \eqref{eq:129} that $NS_u$ is
$\lambda_u$-null.  We conclude that $(c([\gamma])\cdot R)(f) = 0$ for
almost every $\gamma\in G_u$.  

Next, suppose $s\in S_u$ and $s\cdot R = \pi^s\rtimes
U^s$.  We have $U^s_t = U_{s\inv t s} = U_t$ and $\pi^s = \pi\circ
\alpha_s\inv$.  Then for each $f\in C_c(S_u,A(u))$ we compute 
\[
\pi^s \rtimes U^s(f) = \int_S \pi(\alpha_s\inv(f(t)))U_t d\beta^u(t) = 
\int_S U_s^* \pi(f(s)) U_s U_t d\beta^u(t) = U_s^* \pi\rtimes U(f) U_s.
\]
Hence $s\cdot R$ and $R$ are unitarily equivalent.  In particular
$\gamma\cdot R = c([\gamma])\cdot (\delta(\gamma)\cdot R) \cong
c([\gamma]) \cdot R$ and the previous paragraph implies that
$\gamma\cdot R(f) = 0$ for $\lambda_u$-almost all $\gamma$.  Since $G$
acts continuously on $(A\rtimes S)\sidehat$, the map $\gamma\mapsto
\gamma\cdot R(f)$ is continuous.  Furthermore $\supp \lambda_u = G_u$ and
$\gamma\cdot R(f) = 0$ for $\lambda_u$-almost every $\gamma\in G_u$ so
that we
must have $\gamma\cdot R(f) = 0$ for all $\gamma\in G_u$.  Hence $\ker
Q \subset \bigcap_{\gamma\in G_u} \ker (\gamma\cdot R)$.  The other
inclusion is obvious since if $f\in \ker (\gamma\cdot R)$ for all
$\gamma \in G_u$ then for any $\phi\in L^2(G_u/S_u,\mcal{H})$ we have
\[
Q(f)\phi([\gamma]) = c([\gamma])\cdot R(f)\phi([\gamma]) = 0
\]
for all $[\gamma]$.  Thus $Q(f)=0$ and $\ker Q = \bigcap_{\gamma\in
  G_u} \ker(\gamma\cdot R)$.
\end{proof}

\begin{remark}
\index{direct integral}
Since $Q$ is a representation of a $C_0(G\unit)$-algebra 
it must have a decomposition.  
We showed in the proof of Lemma \ref{lem:35} that $Q$ decomposes as 
\[
Q = \int_{G_u/S_u}^\oplus c([\gamma])\cdot R\,d\sigma^u([\gamma]).
\]
where we view $L^2(G_u/S_u,\mcal{H},\sigma^u)$ as sections of the
trivial bundle.  Furthermore, modulo some $\sigma$-finite concerns, 
we could use the fact that $G_u/S_u$ is
Borel isomorphic to $G\cdot u\subset G\unit$ to view this as a
decomposition over $G\unit$.
\end{remark}

\subsection{The Main Result}

We now have all but one of the tools we need to prove the main result.  
In order to make
proper use of our restriction map we need this useful technical lemma.  

\begin{lemma}[{\cite[Lemma 8.38]{tfb2}}]
\label{lem:37}
Suppose $A$ is a $C^*$-algebra.  For each closed set $F\subset \Prim
A$ let $I(F)$ be the corresponding ideal in $\mcal{I}(A)$.  Then a net
$\{I(F_j)\}$ converges to $I(F)$ in $\mcal{I}(A)$ if and only if
given $P\in F$ there is a subnet $\{I(F_{j_k})\}$ and $P_k \in
F_{j_k}$ such that $P_k\rightarrow P$ in $\Prim A$. 
\end{lemma}

\begin{proof}
Suppose that $I(F_j)\rightarrow I(F)$ in $\mcal{I}(A)$ and that $P\in
F$.  Let $U$ be a neighborhood of $P$ in $\Prim A$ and let $J =
I(\Prim A \setminus U)$ be the ideal corresponding to the complement of
$U$.  Then $I(F)\not\supset J$ and therefore
\[
\mcal{O}_J = \{I\in\mcal{I}(A):I\not\supset J\}
\]
is a neighborhood of $I(F)$.  Thus there is a $j_0$ such that $j\geq
j_0$ implies that $I(F_j)\in \mcal{O}_J$.  In particular, if $j\geq
j_0$ then $U\cap F_j \neq \emptyset$.  Next, if we define 
\[
M := \{ (U,j):\text{$U$ is a neighborhood of $P$ and $U\cap F_j\neq
  \emptyset$}\}.
\]
then $M$ is directed by decreasing $U$ and increasing $j$.
Observe that $\{I(F_{U,j})\}$ is a subnet of $\{I(F_j)\}$.  
For each $m=(U,j)\in M$ we can pick $P_m \in
F_j\cap U$.  Then $\{P_m\}$ converges to $P$ as required.  

For the converse, suppose that $\{I(F_j)\}$ has the property given in
the lemma and that $I(F_j)\not\rightarrow I(F)$.  After passing to a
subnet, and relabeling, we can assume that there is an open set
$U\subset \Prim A$ such that $U\cap F \neq \emptyset$ and such that
$F_j\cap U = \emptyset$ for all $j$.  But if $P\in F\cap U$ then we
can pass to a subnet, relabel, and find $P_j\in F_j$ such that
$P_j\rightarrow P$. Then $P_j$ must eventually be in $U$ which is a
contradiction.  
\end{proof}

We have now acquired everything we need to identify the spectrum of
$A\rtimes G$.  

\begin{theorem}
\index{New Result}
\label{thm:crossedstab}
Suppose $(A,G,\alpha)$ is a separable dynamical system and that the
isotropy subgroupoid $S$ is abelian and has continuously varying
stabilizers.  If  $G\unit/G$ is $T_0$ then $\Phi:(A\rtimes
S)\sidehat\rightarrow (A\rtimes G)\sidehat$ defined by $\Phi(R) =
\Ind_S^GR$ is open and factors to a homeomorphism from $(A\rtimes
S)\sidehat/G$ onto $(A\rtimes G)\sidehat$. 
\end{theorem}

\begin{proof}
It follows from Proposition \ref{prop:97} that $\Phi$ is a continuous
surjection and from Proposition \ref{prop:104} that $\Phi$ factors to
a bijection on $(A\rtimes S)\sidehat/G$. All that remains is to show
that $\Phi$ is open.  We will use Proposition \ref{prop:9}.  Suppose
$\Phi(R_i)\rightarrow \Phi(R)$ so that, almost by definition,
$\ker \Phi(R_i)\rightarrow \ker \Phi(R)$.   Using Corollary
\ref{cor:15} we know that $\Res_M$ is continuous and therefore 
\[
\Res_M \ker \Phi(R_i) = \Res_M \ker \Ind_S^G R_i \rightarrow
\Res_M \ker \Phi(R) = \Res_M\ker \Ind_S^G R.
\]
Let $u= \sigma(R)$ and $u_i = \sigma(R_i)$ for all $l$ where $\sigma:(A\rtimes
S)\sidehat \rightarrow S\unit$ is the natural map.  
Using the identifications made in Remark \ref{rem:25} and Lemma
\ref{lem:35} we have 
\begin{align*}
\Res_M\ker \Ind_S^G R &= \bigcap_{\gamma\in G_u} \ker(\gamma\cdot R),\quad\text{and}\\
\Res_M\ker \Ind_S^G R_i  &= \bigcap_{\gamma\in
  G_{u_i}}\ker(\gamma\cdot R_i)
\end{align*}
for all $i$.  It follows from the definition of the Jacobson topology 
that the closed sets
associated to $\Res_M\ker \Ind_S^G R$  and $\Res_M\ker \Ind_S^G R_i$
are 
\begin{align*}
F &= \overline{\{ \ker \gamma\cdot R : \gamma\in G_u\}},\quad\text{and} \\
F_i &= \overline{\{\ker\gamma\cdot R_i : \gamma \in G_{u_i}\}},
\end{align*}
respectively.  Since $\ker R\in F$ it follows from Lemma \ref{lem:37}
that, after passing to a subnet and relabeling, there exists $P_i \in
F_i$ such that $P_i\rightarrow \ker R$.  

Let $\mcal{U}=\{U\}$ be a neighborhood basis of $\ker R$.  
For each $U\in\mcal{U}$ there exists $i_0$ such that $i\geq i_0$ implies that
$P_i\in U$.  We let 
\[
M:= \{ (U,i) : U\in\mcal{U}, P_i\in U\}.
\]
and direct $M$ by decreasing $U$ and increasing $i$.  Then $M$ is a
subnet of $i$ such that $P_{(U,i)}\in U$ for all $(U,i)\in M$.  Given
$(U,i)\in M$ since $U$ is an open set
containing $P_i$ there exists
$\gamma_{(U,i)}\in G_{u_i}$ such that $\ker \gamma_{(U,i)}\cdot R_i \in
U$.  Now, given any $U_0\in \mcal{U}$, choose $i_0$ so
that $P_{i_0}\in U$ and $(U_0,i_0)\in M$.  If $(U,i)\in M$
such that $(U_0,i_0)\leq (U,i)$ then $\ker\gamma_{(U,i)}\cdot R_i \in
U\subset U_0$.  Thus 
\[
\ker \gamma_{(U,i)}\cdot R_i\rightarrow \ker R.
\]
However this implies that any ideal not contained in $\ker R$ is
eventually not contained in $\ker \gamma_{(U,i)}\cdot R_i$.  Thus, by
definition, $\gamma_{(U,i)}\cdot R_i\rightarrow R$.  This suffices to
show that $\Phi$ is open. 
\end{proof}

\begin{remark}
\label{rem:29}
If there is a problem with Theorem \ref{thm:crossedstab} it is that
$(A\rtimes S)\sidehat$ can be just as mysterious as $(A\rtimes
G)\sidehat$.  As we will see, and have seen, there are times when
$(A\rtimes S)\sidehat$ can be analyzed, but in general it is
difficult.  For instance, if $A$ has Hausdorff spectrum (and is
separable) then each
fibre $A(u)$ can be identified with the compacts and in this case 
$A(u)\rtimes S_u$ is relatively well understood \cite[Section
7.3]{tfb2} and in particular is isomorphic to 
$C^*(S_u,\bar{\omega}_u)$ where $[\omega_u]$ is the Mackey obstruction for
$\alpha|_{S_u}$.  However, even if the stabilizers vary continuously,
the collection $\{\omega_u\}$ may be poorly behaved and identifying
the total space topology of $(A\rtimes S)\sidehat$ may be difficult. 
On the bright side, $A\rtimes S$ is a bundle product and there are
times when we can say something about it.  For instance, consider the
scalar case.  Then $(A\rtimes S)\sidehat$ becomes $\widehat{S}$ 
which is much simpler.
Or consider the case where $\alpha$ is ``locally unitary on the
stabilizers.''  Then Theorem \ref{thm:locunit} tells us that
$(A\rtimes S)\sidehat$ is a principal bundle.  In particular, it is
determined up to isomorphism by a cohomology class and thus
$(A\rtimes G)\sidehat$ has a nice cohomological invariant.  
\end{remark}

The following corollary is immediate and interesting enough to be
worth writing down.  We will explore further applications of Theorem
\ref{thm:crossedstab} in the next chapter.  

\begin{corr}
Suppose $(A,G,\alpha)$ is a separable dynamical system, that $G$ is a
principal groupoid, and that $G\unit/G$ is $T_0$.  Then $(A\rtimes
G)\sidehat$ is isomorphic to $\widehat{A}/G$. 
\end{corr}

\begin{proof}
Since $G$ is principal it clearly has continuously varying abelian
isotropy.  In fact the isotropy subgroupoid is just $G\unit$.
Furthermore, we have $A\rtimes G\unit = A$ and the result now
follows from Theorem \ref{thm:crossedstab}.
\end{proof}


\chapter{Examples and Applications}
\label{cha:examples}
In this chapter we present a number of applications of Theorem
\ref{thm:crossedstab}.  Section \ref{sec:groupstab} contains a
strengthening of the main result in the scalar case.  The scalar case
is particularly interesting because we identify the spectrum of
$C^*(G)$ as a quotient of the much better understood $\widehat{S}$.
In Section \ref{sec:redux} we apply these results to
transformation groupoids and transformation groupoid algebras.  This
allows us to present a couple of interesting examples and as well as
connect these results back to the existing theory.  Finally, in Section
\ref{sec:haussdorff} we give a partial analysis of when a
groupoid $C^*$-algebra has Hausdorff spectrum.  

\section{Groupoid Algebras with Abelian Isotropy}
\label{sec:groupstab}
We would like to address the concerns made in Remark \ref{rem:29} and
show that, at least in the scalar case, we can come up with a much
more concrete result.  In particular, the topology on $\widehat{S}$ is
well understood so that the following corollary to Theorem
\ref{thm:crossedstab} gives us a very useful identification of the
topology of $C^*(G)\sidehat$.

\begin{corr}
\label{cor:16}
Let $G$ be a second countable locally compact Hausdorff groupoid with
a Haar system.  Furthermore, suppose the stabilizer subgroupoid $S$ is
abelian and varies continuously.  If $G\unit/G$ is $T_0$ then the map
$\omega\mapsto \Ind_{S}^G\omega$ is open and factors to a
homeomorphism of $\widehat{S}/G$ onto $C^*(G)\sidehat$. 
\end{corr}

\begin{proof}
If $G$ is as above then we define the groupoid $C^*$-algebra to be
$C^*(G) = C_0(G\unit)\rtimes G$.  Furthermore, $C^*(S)\sidehat$ is
exactly the dual $\widehat{S}$ of $S$ as described in Section
\ref{sec:duality}.  Under these conditions it follows
immediately from Theorem \ref{thm:crossedstab} that the induction map
is open and factors to the desired homeomorphism.  
\end{proof}

We can actually use the machinery developed in Chapter
\ref{cha:fine-structure} to do
better though.  We would like to remove the assumption that $G\unit/G$
is $T_0$.  However, we have the following proposition to consider.

\begin{prop}
\label{prop:103}
Let $G$ be a second countable, locally compact Hausdorff groupoid with
a Haar system. Furthermore, suppose the stabilizer subgroupoid $S$ is
abelian and varies continuously.  Then the following are equivalent:
\begin{enumerate}
\item $G\unit/G$ is $T_0$.
\item $C^*(G)$ is GCR.
\item $C^*(G)$ is Type I.
\end{enumerate}
\end{prop}

\begin{proof}
Since abelian groups are GCR, the fact that (a) is equivalent to (b) follows
from \cite[Theorem 1.1]{ccrgca}.  The fact that (b) and (c) are
equivalent for separable $C^*$-algebras 
is a well known result \cite[Theorem 9.1]{dixmiercstar}.  
\end{proof}

This shows that if we are to deal with the ``non-$T_0$'' case then we
are going to have to work with non-Type I algebras.  This means
working with primitive ideals instead of the
spectrum.  Before we begin, let us consider the action of $G$ on $\widehat{S}$. 

\begin{corr}
\label{cor:17}
\index{G-space@$G$-space}
Suppose $G$ is a second countable, locally compact Hausdorff groupoid
with a Haar system and that the stabilizer subgroupoid $S$ is abelian
and continuously varying.  Then then there is a strongly continuous 
action of $G$ on
$\widehat{S}$ given for $\gamma\in G$ and $\omega\in \widehat{S}$ by 
\begin{equation}
\label{eq:145}
\gamma\cdot \omega(s) = \omega(\gamma\inv s \gamma).
\end{equation}
\end{corr}

\begin{proof}
Of course, the existence of such an action is shown in Corollary
\ref{cor:14}.  Furthermore the action is strongly continuous in this
case because the structure map for $\widehat{S}$, namely $\hat{p}$, is
open since $\widehat{S}$ is continuously varying.  Finally we use
Proposition \ref{prop:100} to see that the action is given by
\eqref{eq:145}.
\end{proof}

\begin{remark}
Lemma \ref{lem:30} still holds even if $G\unit/G$ is not $T_0$ and
in particular we will continue to make the identifications of Remark
\ref{rem:25}.  We will regularly confuse $\Ind_S^G \omega$ and
$\Ind_{S_u}^G \omega$.  
\end{remark}

Now, if $G\unit/G$ is not $T_0$ then we cannot use Proposition
\ref{prop:96} to conclude that $\Ind_S^G \omega$ is irreducible if
$\omega\in \widehat{S}$. It is the main result in \cite{irredreps} that
every representation of a second countable groupoid induced from an
irreducible representation of a stability group is irreducible, even
when the stabilizers are non-abelian.  However, in the abelian case
this result is much closer to the surface and we will give an account
here.  

\begin{prop}[{\cite[Lemma 2.5]{ctgIII}}]
\label{prop:112}
Let $G$ be a second countable locally compact Hausdorff groupoid with
a Haar system and suppose the stabilizer subgroupoid $S$ is abelian.
Then $\Ind_{S_u}^G\omega$ is irreducible for all 
$\omega\in\widehat{S}_u$. 
\end{prop}

\begin{proof}
By Proposition \ref{prop:102}, it suffices to
show that the representation $N^\omega$ is irreducible.  First, since
$\omega$ acts on $\C$, $N^\omega$ acts on
$L^2(G_u/S_u,\sigma^u)$.  
Furthermore, given $f\in C_c(S_u)$ we can compute that 
\begin{align*}
N^\omega(f)\phi([\gamma]) &= \int_G \omega(\delta(\gamma))
\id_\gamma\inv(f(\eta))\overline{\omega(\delta(\eta\inv\gamma))}
\phi([\eta\inv\gamma])d\lambda^{r(\gamma)}(\eta)\\
&= \int_G \overline{\omega(\delta(\eta\inv\gamma)\delta(\gamma)\inv)}
f(\eta)\phi([\eta\inv\gamma])d\lambda^{r(\gamma)}(\eta)
\end{align*}
where $\delta$ is as in Remark \ref{rem:26}.  Furthermore, as in
Remark \ref{rem:24}, we use the range map to identify $G_u/S_u$ and
$G\cdot u$ as Borel spaces, push the measure $\sigma^u$ forward to a
measure $\sigma^u_*$ on $G\cdot u$, and view $N^\omega$ as acting on
$L^2(G\cdot u,\sigma^u_*)$ via 
\[
N^\omega(f)\phi(\gamma\cdot u) = \int_G
\overline{\omega(\delta(\eta\inv\gamma)\delta(\gamma)\inv)} f(\eta)
\phi(\eta\inv\gamma\cdot u)d\lambda^{r(\gamma)}(\eta).
\]  
It is straightforward to show that
$\overline{\omega(\delta(\eta\inv\gamma)\delta(\gamma)\inv)}$ only depends
on $v = \gamma \cdot u$ and $\eta$.  We write $\theta(\eta,v)$ for
the corresponding Borel function.  In particular, this allows us to
rewrite $N^\omega$ as 
\begin{equation}
\label{eq:146}
N^\omega(f)\phi(v) = \int_G \theta(\eta,v) f(\eta)\phi(\eta\inv \cdot
v)d\lambda^v(\eta).
\end{equation}
Let $N^u$ be the representation of $C_0(G\unit)$ on $L^2(G\cdot
u,\sigma^u_*)$ defined in Lemma \ref{lem:34} and recall that
$N^\omega(f\cdot g) = N^u(f)N^\omega(g)$ for all $f\in C_0(G\unit)$
and $g\in C_c(G)$.  Now suppose $T\in N^u(C_0(G\unit))''$ and $P$ is a
projection commuting with $N^\omega(C^*(G))$.  Then in particular 
\[
P N^u(f) N^\omega(g)h = P N^\omega(f\cdot g)h = N^\omega(f\cdot g)Ph =
N^u(f)PN^\omega(g)h.
\]
Since $N^\omega$ is nondegenerate, this implies that $P$ is in the
commutant of $N^u(C_0(G\unit))$ so that $P$ commutes with $T$.  Since
$C_c(G\unit)$ separates points of $G\cdot u$, the von Neumann algebra
$N^u(C_0(G\unit))''$ is a maximal abelian subalgebra of operators on
$L^2(G\cdot u,\sigma_*^u)$.  It follows that any 
projection commuting with $N^\omega(C^*(G))$ must be of the form
$N^u(f)$ where $f = \chi_E$ is the characteristic function of some set
$E\subset G\cdot u$ and where we have extended $N^u$ to $L^\infty(G)$ in
the obvious fashion.  It is easy to see that $N^u(f)$ still commutes
with every $N^\omega(g)$.  Thus for each $g\in C_c(G)$ we have 
\begin{align*}
N^u(f)N^\omega(g)\phi(v) &= f(v) \int_G
\theta(\eta,v)g(\eta)\phi(\eta\inv\cdot v)d\lambda^v(\eta) \\
&= \int_G \theta(\eta,v)g(\eta)f(\eta\inv\cdot v)\phi(\eta\inv\cdot
v)d\lambda^v(\eta) = N^\omega(g)N^u(f)\phi(v)
\end{align*}
for $\sigma_*^u$-almost all $v$.  This suffices to show that $f$ is
constant almost everywhere on $G\cdot u$ and that the projection
$N^\omega(f)$ is a multiple of the identity.  Therefore, the only
projections commuting with $N^\omega(C^*(G))$ are the scalars and $N^\omega$
is irreducible. 
\end{proof}

Thus, even when $G\unit/G$ is not $T_0$, we can still induce
representations from $\widehat{S}$ to elements in the spectrum of
$C^*(G)$.  

\begin{prop}
\label{prop:105}
Let $G$ be a second countable, locally compact Hausdorff groupoid with
a Haar system.  Furthermore suppose the isotropy subgroupoid $S$ is
continuously varying and abelian.  Then $\Phi:\widehat{S}\rightarrow
C^*(G)\sidehat$ defined by $\Phi(\omega) = \Ind_S^G\omega$ is
continuous and open.  
\end{prop}

\begin{proof}
It follows from Proposition \ref{prop:112} that $\Phi$ maps into
$C^*(G)\sidehat$.  Furthermore, just as in Proposition \ref{prop:97},
the continuity of $\Phi$ follows from the general theory of Rieffel
induction.  All that is left to do is show $\Phi$ is open.  This proof
is almost exactly the same as the openness calculation in the proof of
Theorem \ref{thm:crossedstab}.  Suppose $\Ind \omega_i \rightarrow \Ind
\omega$ in $C^*(G)\sidehat$.  Since $\Res_M$ is continuous, it follows
that 
\[
I_i = \Res_M \ker \Ind_S^G \omega_i \rightarrow I= \Res_M \ker \Ind_S^G \omega.
\]
Now, the spectrum of $C^*(S)$ is Hausdorff so that we can identify
$\widehat{S}$ and $\Prim C^*(S)$.  In particular, under this
identification  Lemma \ref{lem:35} tells us that 
\begin{align*}
I = \bigcap_{\gamma\in G_{u}} \gamma\cdot \omega,\quad\text{and}\quad I_i = \bigcap_{\gamma\in
  G_{u_i}} \gamma \cdot \omega_i
\end{align*}
for all $i$.  Hence, the closed set associated to $I$ is
$\overline{G\cdot \omega}$ and the closed set associated to $I_i$ is
$\overline{G\cdot \omega_i}$ for all $i$.  
Since $\omega\in \overline{G\cdot \omega}$ it follows from Lemma \ref{lem:37}
that, after passing to a subnet and relabeling, there exists $\chi_i \in
\overline{G\cdot \omega_i}$ such that $\chi_i\rightarrow \omega$.  

Let $u=\hat{p}(\omega)$, $u_i=\hat{p}(\omega_i)$ for all $i$, and 
$\mcal{U}$ be a neighborhood basis of $\omega$.  
For each $U\in\mcal{U}$ there exists an $i_0$ such that $i\geq i_0$ implies that
$\chi_i\in U$.  We let 
\[
M:= \{ (U,i) : U\in\mcal{U}, \chi_i\in U\}.
\]
and direct $M$ by decreasing $U$ and increasing $i$.  Then $M$ is a
subnet of $i$ such that $\chi_i\in U$ for all $(U,i)\in M$.  Given
$(U,i)\in M$ since $U$ is an open set
containing $\chi_i$ there exists
$\gamma_{(U,i)}\in G_{u_i}$ such that $\gamma_{(U,i)}\cdot \omega_i \in
U$.  Now, given any $U_0\in \mcal{U}$ choose $i_0$ so
that $\chi_{i_0}\in U$ and $(U_0,i_0)\in M$.  If $(U,i)\in M$
such that $(U_0,i_0)\leq (U,i)$ then $\gamma_{(U,i)}\cdot \omega_i \in
U\subset U_0$.  Thus 
\[
\gamma_{(U,i)}\cdot \omega_i\rightarrow \omega.
\]
This suffices to show that $\Phi$ is open. 
\end{proof}

Of course, as we said before, 
we need to be working with primitive ideals so we will make the
switch now.  

\begin{corr}
\label{cor:19}
Let $G$ be a second countable, locally compact Hausdorff groupoid with
a Haar system.  Furthermore suppose the isotropy subgroupoid $S$ is
continuously varying and abelian.  Then $\Psi:\widehat{S}\rightarrow
\Prim C^*(G)$ defined by $\Psi(\omega) = \ker \Ind_S^G\omega$ is
continuous and open.   
\end{corr}

\begin{proof}
Of course, $\Psi$ is just $\Phi$ composed with the map $\pi\mapsto
\ker \pi$.  Since both maps are continuous and open, their composition
must be also.  
\end{proof}

We would like to factor $\Psi$ to a homeomorphism and to do that 
we will need to get a
handle on the equivalence relation determined by $\Psi$.  

\begin{lemma}
\label{lem:38}
Let $G$ be a second countable, locally compact Hausdorff groupoid with
a Haar system.  Furthermore suppose the isotropy subgroupoid $S$ is
continuously varying and abelian.  Then $\Psi(\omega) = \Psi(\chi)$ if
and only if $\overline{G\cdot \omega} = \overline{G\cdot \chi}$. 
\end{lemma}

\begin{proof}
If $\ker \Ind_S^G \chi = \ker \Ind_S^G \omega$ then we must have 
\[
\Res_M \ker \Ind_S^G \chi = \Res_M \ker \Ind_S^G \omega.
\]
However, it now follows from Lemma \ref{lem:35}, after identifying
$\widehat{S}$ and $\Prim C^*(S)$, that 
\[
\bigcap_{\gamma\in G_{\hat{p}(\omega)}}\gamma\cdot \omega = 
\bigcap_{\gamma\in G_{\hat{p}(\chi)}}\gamma\cdot \chi.
\]
This implies that the closed sets in $\widehat{S}$ associated
to these ideals must be the same.  Hence $\overline{G\cdot \omega} =
\overline{G\cdot \chi}$.  

For the reverse direction suppose that $\overline{G\cdot \omega} =
\overline{G\cdot \chi}$.  This implies that there exists $\gamma_i$
such that $\chi = \lim_i \gamma_i\cdot \omega$.  It follows
from Proposition \ref{prop:104} that $\Ind_S^G \omega$ is equivalent
to $\Ind_S^G \gamma\cdot \omega$ for all $\gamma$.  Thus $\Psi$ is
$G$-invariant and, since $\Psi$ is continuous, we get
$\Psi(\gamma_i\cdot \omega) = \Psi(\omega)\rightarrow \Psi(\chi)$.
Thus $\Psi(\chi) \in \overline{\{\Psi(\omega)\}}$ and, by definition
of the hull-kernel topology, $\Psi(\omega)\subset\Psi(\chi)$.
Reversing the roles of $\omega$ and $\chi$ above will yield the other
inclusion. 
\end{proof}

Technically the next definition uses induction for ideals, which we
haven't actually introduced.  All the reader needs to know is that
it is characterized by the formula $\Ind_H^G \ker \pi = \ker \Ind_H^G
\pi$.  

\begin{definition}
\index{EH-regular}
Let $G$ be a second countable locally compact Hausdorff groupoid
with a Haar system.  We say that $G$ is {\em EH-regular} if every
primitive ideal is induced from an isotropy subgroup.  That is, given
$P\in \Prim C^*(G)$ there exists $u\in G\unit$ and $Q\in\Prim
C^*(S_u)$ such that $P = \Ind_{S_u}^G P$. 
\end{definition}

Of course, we have already met a large class of groupoids which are
EH-regular.  

\begin{corr}
\label{cor:18}
If $G$ is a second countable locally compact Hausdorff groupoid and
$G\unit/G$ is $T_0$ then $G$ is EH-regular. 
\end{corr}

\begin{proof}
Theorem \ref{thm:ehregularity} tells us that every irreducible
representation is induced from a stabilizer, which of course implies
that every primitive ideal is induced from a stabilizer. 
\end{proof}

Of course, the whole point is to get away from the $T_0$ case so we
cite the following result.  

\begin{theorem}[{\cite[Theorem 2.1]{geneffhan}}]
\label{thm:geneffhan}
Assume that $G$ is a second countable, locally compact Hausdorff
groupoid with a Haar system.  If $G$ is amenable then every
primitive ideal is induced from a stability group.  In other words, $G$
is EH-regular. 
\end{theorem}

\begin{remark}
\label{rem:28}
If this theorem leaves something to be desired it is that, as we saw in
Section \ref{sec:amenable}, groupoid amenability is not a
transparent condition.  It is worth noting that not all
principal groupoids are amenable so that in particular not all
groupoids with abelian, continuously varying stabilizer are amenable.
Thus, the amenability condition in Theorem \ref{thm:scalarstab} below is not
superfluous.  
\end{remark}

This theorem allows us to give a strengthening of Theorem
\ref{thm:crossedstab} in the scalar case, however we need to, briefly,
introduce a new construction. 

\begin{definition}
\label{def:59}
\index[not]{$(X)^{T_0}$}
\index{Tzeroization@$T_0$-ization}
If $X$ is a topological space, then the {\em $T_0$-ization} of $X$ is
the quotient space $(X)^{T_0}:=X/\sim$ where $\sim$ is the equivalence
relation on $X$ defined by $x\sim y$ if $\overline{\{x\}} =
\overline{\{y\}}$.  We equip $(X)^{T_0}$ with the quotient topology. 
\end{definition}

We will not use the following lemma and do not provide a proof, 
but it sheds some light on the
$T_0$-ization definition. 

\begin{lemma}[{\cite[Lemma 6.10]{tfb2}}]
If $X$ is a topological space then $(X)^{T_0}$ is a $T_0$ space.  If
$Y$ is any $T_0$ topological space and if $f:X\rightarrow Y$ is
continuous then there is a continuous map $f':(X)^{T_0}\rightarrow Y$
such that $f = f' \circ q$ where $q:X\rightarrow (X)^{T_0}$ is the
quotient map. 
\end{lemma}

Now we have enough technology to prove the main result of this
section.  

\begin{theorem}
\label{thm:scalarstab}
\index{New Result}
\index{groupoid $C^*$-algebra}
Suppose $G$ is a second countable locally compact Hausdorff groupoid
with a Haar system and that the stabilizer subgroupoid $S$ is abelian
and continuously varying.  If $G$ is EH-regular, and in
particular if $G$ is amenable or $G\unit/G$ is $T_0$, then the map
$\Psi:\widehat{S}\rightarrow \Prim C^*(G)$ such that $\Psi(\omega) =
\ker\Ind_S^G \omega$ factors to a homeomorphism
of $\Prim C^*(G)$ with $(\widehat{S}/G)^{T_0}$.  
\end{theorem}

\begin{proof}
It follows from Corollary \ref{cor:19} that $\Psi$ is continuous and
open.  Furthermore, once we identify $\widehat{S}$ with $\Prim
C^*(S)$, it is clear that $\Psi$ is surjective if $G$ is EH-regular, 
which, by Theorem \ref{thm:geneffhan}, occurs whenever $G$ is amenable,
or, by Corollary \ref{cor:18}, when $G\unit/G$ is $T_0$.
Finally, it is straightforward to show that $\overline{G\cdot \omega}
= \overline{G\cdot \chi}$ in $\widehat{S}$ if and only if
$\overline{\{G\cdot \omega\}} = \overline{\{G\cdot \chi\}}$ in
$\widehat{S}/G$.  Thus it follows 
from Lemma \ref{lem:38} that the factorization of
$\Psi$ to $(\widehat{S}/G)^{T_0}$ is injective and therefore a
homeomorphism. 
\end{proof}

\begin{remark}
If $G\unit/G$ is $T_0$ then $G$ is EH-regular by
Corollary \ref{cor:18}, and $C^*(G)$ is Type I by Proposition
\ref{prop:103}.  In particular $\pi\mapsto \ker \pi$ is a
homeomorphism of $C^*(S)\sidehat$ onto $\Prim C^*(S)$ so that in this
case Theorem \ref{thm:scalarstab} reduces to Corollary \ref{cor:16}.
\end{remark}

As in Section \ref{sec:crossedstab} we get the following
corollary. 

\begin{corr}
If $G$ is a second countable, locally compact, EH-regular,
principal groupoid with a Haar system then $\Prim C^*(G)$ is
homeomorphic to $(G\unit/G)^{T_0}$. 
\end{corr}

\begin{proof}
Since $G$ is principal it has trivial, and hence continuously varying, abelian
stabilizer $S = G\unit$ and the result follows from Theorem
\ref{thm:scalarstab}.
\end{proof}

\begin{remark}
While requiring the stabilizer subgroupoid $S$ to have a Haar system
is natural from a groupoid point of view, it is a strong assumption.
In particular, as will see in Section \ref{sec:redux}, it
is this assumption which prevents Theorem \ref{thm:scalarstab} from
completely generalizing \cite[Theorem 8.39]{tfb2}.  However, removing
this hypothesis is a serious challenge in that, unless $S$ is
continuously varying, $C^*(S)$ and $\widehat{S}$ do not exist as we
have defined them.  
\end{remark}


\section{Transformation Groupoids Redux}
\label{sec:redux}

The purpose of this section is threefold.  First, we would like to
apply Theorem \ref{thm:scalarstab} to groupoid actions and restate the
results in terms of the transformation 
groupoid algebra.  Second, we would like to
show that Theorem \ref{thm:scalarstab} is a partial generalization of
the results in \cite[Section 8.3]{tfb2}.  Finally, we will present two
examples which show how we can use this theory.  
The first task is straightforward.  Recall from Proposition
\ref{prop:15} that the stabilizers and orbit space of a groupoid
action appear naturally as the stabilizers and orbit space of the
transformation groupoid.  Proving statements about transformation
groupoids often requires little more than rewording corresponding
statements for groupoids.  For instance EH-regularity is remolded into
the following

\begin{definition}
\label{def:60}
\index{EH-regular}
Let $G$ be a second countable, locally compact Hausdorff groupoid with
a Haar system acting
continuously on a second countable, locally compact Hausdorff space $X$.  We say
that $(G,X)$ is {\em EH-regular} if every primitive ideal is induced
from an isotropy subgroup.  That is, given $P\in \Prim C^*(G,X)$ then
there exists $x\in X$ and $Q\in\Prim C^*(G_x)$ such that $P =
\Ind_{G_x}^G Q$.  
\end{definition}

Observe that the stabilizer $G_x$ is just the stabilizer of the
transformation groupoid $G\ltimes X$ at $x$ and Definition \ref{def:60} is
exactly the same as requiring that $G\ltimes X$ be EH-regular. With
this in mind the following corollary is unsurprising. 

\begin{corr}
\label{cor:20}
\index{transformation groupoid!$C^*$-algebra}
Suppose $G$ is a second countable, locally compact Hausdorff groupoid
with a Haar system acting continuously on a second countable, locally
compact Hausdorff space $X$.  Furthermore, suppose that the stabilizer
group bundle of the action $S$ is abelian and varies continuously.
Then there is an action of $G$ on $\widehat{S}$ defined for $\gamma\in
G$ and $\omega\in \widehat{S}$ by 
\begin{equation}
\label{eq:155}
\gamma\cdot \omega(s) = \omega(\gamma\inv s \gamma).
\end{equation}
Furthermore, if
$(G,X)$ is EH-regular, which holds if $(G,X)$ is
amenable or $X/G$ is $T_0$, then
the map $\Psi:\widehat{S}\rightarrow \Prim C^*(G,X)$ such that
$\Psi(\omega) = \ker\Ind_S^{G\ltimes X} \omega$ factors to a
homeomorphism of $\Prim C^*(G,X)$ with $(\widehat{S}/G)^{T_0}$.  
\end{corr}

\begin{proof}
Since $G$ and $X$ are second countable, the transformation groupoid is
as well.  Since $G$ has a Haar system the transformation groupoid
does.  Furthermore the stabilizer group bundle of the action is
defined to be the stabilizer subgroupoid of $G\ltimes X$ and is
abelian and continuously varying by assumption.  As we noted above, the
condition that $(G,X)$ is EH-regular is equivalent to requiring that
$G\ltimes X$ is EH-regular.  In Definition \ref{def:16}
we say that $(G,X)$ is amenable if and only if $G\ltimes X$ is and 
since $G$ and $G\ltimes X$ have the same action on $X$ it is clear
that $X/G=X/(G\ltimes X)$.  At
this point we have everything we need to apply Theorem
\ref{thm:scalarstab}. The only thing that is not clear is what we mean
by $\widehat{S}/G$.  Observe that by composing $\hat{p}$ with the
range map $r:X\rightarrow G\unit$ we get a range map on
$\widehat{S}$.  Suppose $\omega \in S$, $\gamma\in G$ and $s(\gamma)=
r(\hat{p}(\omega))$.  Let $x=\hat{p}(\omega)$.  Then $s(\gamma) =
r(x)$ so that $(\gamma,\gamma\cdot x)\in G\ltimes X$.  Furthermore, 
$s(\gamma,\gamma\cdot x) = x = \hat{p}(\omega)$ so that we can let
$(\gamma,\gamma\cdot x)$ act on $\omega$.  We may as
well {\em define} 
\[
\gamma\cdot \omega := (\gamma,\gamma\cdot \hat{p}(\omega))\cdot\omega.
\]  
Then, given $s\in S_x$ we have 
\[
(\gamma,\gamma\cdot x)\cdot \omega(s) = 
\omega((\gamma,\gamma\cdot x)\inv (s,x) (\gamma,\gamma\cdot x)) 
= \omega(\gamma\inv s\gamma)
\]
where we are being a little sloppy about distinguishing between $s$
and $(s,x)$.  Thus the action of $G$ on $\omega$ is given by 
\eqref{eq:155}.  It is now straightforward to show
that with this action $\widehat{S}$ is a continuous $G$-space.
Furthermore if $\omega = \gamma \cdot \chi$ then $\omega =
(\gamma,\gamma\cdot \hat{p}(\chi))\cdot \chi$ by definition.  Conversely, if
$\omega = (\gamma,x)\cdot \chi$ then, observing that
$\gamma\inv\cdot x = \hat{p}(\chi)$, we have $\omega =\gamma\cdot \chi$.
Thus $G$ and $G\ltimes X$ have the same orbits in $S$ so that $S/G =
S/G\ltimes X$.   The corollary now follows from Theorem
\ref{thm:scalarstab}. 
\end{proof}

With this we have a convenient restatement of Theorem \ref{thm:scalarstab}
which doesn't (directly) use the transformation groupoid.  We will use
this result to explore some examples later on in this section.  First,
though, we would like to show that Corollary \ref{cor:20} is a partial
generalization of the known results for group actions.  First, the
relevant theorems from \cite{tfb2} have been reproduced below. 

\begin{theorem}[Gootman-Rosenberg-Sauvageot {\cite[Theorem 8.21]{tfb2}}]
\label{thm:grsthm}
Suppose that $(A,G,\alpha)$ is a separable {\em group} dynamical
system with $G$ amenable.  Then $(A,G,\alpha)$ is EH-regular.  
\end{theorem}

\begin{remark}
In particular, once one sorts out what EH-regularity means for
group dynamical systems, this implies that any second countable,
abelian transformation group $(G,X)$ is EH-regular.
\end{remark}

\begin{theorem}[{\cite[Theorem 8.39]{tfb2}}]
\label{thm:groupstab}
Let $(G,X)$ be a locally compact transformation group with $G$
abelian.  Then $\Phi:X\times \widehat{G}\rightarrow \Prim C^*(G)$ such
that $\Phi(x,\omega) = \Ind_{G_x}^G(\omega|_{G_x})$,
where $G_x$ is the stabilizer at $x$, is continuous and open.
Furthermore $\Phi$ factors through $X\times \widehat{G}/\sim$ where
$(x,\omega)\sim(x,\chi)$ if $\overline{G\cdot x} = \overline{G\cdot y}$
and $\chi\overline{\omega} \in G_x^\perp$, and defines a homeomorphism
onto its range.  If $(G,X)$ is EH-regular, which is automatic if
$(G,X)$ is second countable by the GRS-theorem, then $\Phi$ defines a
homeomorphism of $X\times \widehat{G}/\sim$ onto $\Prim C^*(G,X)$.
\end{theorem}

\begin{proof}[Remark.]
Those readers who are careful about their references will notice some
minor discrepancies.  The main difference is that
the result in \cite{tfb2} is stated in
terms of the crossed product $C_0(X)\rtimes_{\lt} G$.  Of course, we
saw in Section \ref{sec:transform} that this is isomorphic to
$C^*(G,X)$ and we will not distinguish between the two here.  
\end{proof}

Before we begin our analysis in earnest let us make two remarks.

\begin{remark}
First, let us consider the problem of separability.  
As we noted in Remark \ref{rem:27}, groupoid
crossed products and groupoid algebras are heavily dependent on
separability hypothesis.  We would not expect to be able to reproduce
Theorem \ref{thm:groupstab} in the nonseparable case using groupoids.  
Second, let us consider amenability.  As we noted in
Remark \ref{rem:28}, the amenability hypothesis in Theorem
\ref{thm:scalarstab}, and hence Corollary \ref{cor:20}, are not
specious.  On the other hand, while Theorem \ref{thm:grsthm} does have
an amenability hypothesis this assumption will disappear in Theorem
\ref{thm:groupstab} because abelian groups are always amenable.  
\end{remark}

\begin{remark}
We should also mention the most important difference between Theorem
\ref{thm:groupstab} and Corollary \ref{cor:20}.  Suppose a group $G$ acts
on a space $X$. If $G$ is abelian then all of the stabilizers are abelian.
However, they certainly don't have to vary continuously.  In
particular Theorem \ref{thm:groupstab} holds for non-continuously
varying stabilizers.  On the other hand, even in the transformation
group case, Corollary \ref{cor:20} doesn't
make sense if $S$ doesn't vary continuously.  
This is an unsatisfactory aspect of the current theory and it
is an open question if/how it can be addressed. Let us finish by
pointing out that Corollary \ref{cor:20} does have its uses.  
In Example \ref{ex:27} we present an action of a {\em
  nonabelian} group with continuously varying abelian stabilizers.
This action can be studied using Corollary \ref{cor:20} but is outside
the scope of Theorem \ref{thm:groupstab}.  
\end{remark}

So, we cannot fully reproduce Theorem \ref{thm:groupstab}
with our current theory.  Our goal will be to show that given a
transformation group $G$ such that $G$ is abelian {\em and} has
continuously varying stabilizers then Theorem \ref{thm:groupstab} and
Corollary \ref{cor:20} say the same thing.  

\begin{prop}
Let $G$ be a second countable, locally compact Hausdorff abelian
group acting on a second countable, locally compact
Hausdorff space $X$.  Furthermore suppose that the stabilizers vary
continuously in $G$.  Then $X\times \widehat{G}/\sim$ is naturally
homeomorphic to $(\widehat{S}/G)^{T_0}$ so that Theorem
\ref{thm:groupstab} and Corollary \ref{cor:20} have the same
conclusion. 
\end{prop}

\begin{proof}
First, observe that if $G$ is abelian then all of the stabilizers are
abelian.  Thus $(G,X)$ satisfies both the requirements of Theorem
\ref{thm:groupstab} and Corollary \ref{cor:20}.  Furthermore, the
GRS-theorem implies that $(G,X)$ is EH-regular so that this condition
is satisfied for both Theorem \ref{thm:groupstab} and Corollary
\ref{cor:20}.  All that is left is to show that
$X\times\widehat{G}/\sim$ and $(\widehat{S}/G)^{T_0}$ are naturally
isomorphic so that in this case each theorem can be obtained from the
other. 

Define $\rho:X\times\widehat{G}\rightarrow \widehat{S}$ by
$\rho(x,\omega) = \omega|_{S_x}$.  It is straightforward to use
Proposition \ref{prop:33} to show that $\rho$ is continuous.  It is a
classical result \cite{rudinfourier} that the dual of a subgroup is
isomorphic to a quotient of the dual of the full group via restriction.  Hence
characters on subgroups can
be extended to characters on the full group and $\rho$
is a surjection.  Define
$\equiv$ on $X\times \widehat{G}$ by $(x,\omega) \equiv (y,\chi)$ if and
only if $x=y$ and $\chi\overline{\omega} \in G_x^\perp$. It is clear
from \cite{rudinfourier} 
that $\rho(x,\omega) = \rho(y,\chi)$ if and only if
$(x,\omega)\equiv(y,\chi)$.  Thus, if we can show $\rho$ is open, it
will follow that it factors to a homeomorphism from
$X\times\widehat{G}/\equiv$ onto $\widehat{S}$.  Suppose
$\rho(x_i,\omega_i)\rightarrow \rho(x,\omega)$.  Since
$\hat{p}(\rho(y,\chi)) = y$ for all $(y,\chi)\in X\times\widehat{G}$
we clearly have $x_i\rightarrow x$.  Recall from \cite[Lemma 2.35]{tfb2}
that,  after identifying $C^*(G,X)$ with $C_0(X)\rtimes G$, there is a
natural map
$\iota:C^*(G)\rightarrow M(C^*(G,X))$.  This map induces a restriction
map $\Res_G$ from representations of $C^*(G,X)$ to $C^*(G)$.  In
particular, as with all such restriction maps, it is a continuous
process and we have  
\[
\Res_G \Ind_S^{G\ltimes X} \omega_i|_{G_{x_i}} \rightarrow \Res_G
\Ind_S^{G\ltimes X} \omega|_{G_x}.
\]
Basically, what is going on is that $\Ind_S^{G\ltimes X} \chi|_{G_y}$
is a representation of $C^*(G,X)\cong C_0(X)\rtimes G$ for all
$(y,\chi)\in X\times\widehat{G}$, and as such
must be the integrated form of some covariant
representation $(\pi,U)$.  The restriction map $\Res_G$ gives us
the unitary part $U$.  However, \cite[Corollary 5.6]{tfb2} says
that this unitary part is equivalent to
$\Ind_{G_y}^G \chi|_{G_y}$.  Thus, putting it all together, plus a
little more, we have 
\begin{equation}
\label{eq:147}
\ker\Ind_{G_{x_i}}^G \omega_i|_{G_{x_i}} \rightarrow \ker\Ind_{G_x}^G
\omega|_{G_x}.
\end{equation}
It follows from \cite[Proposition 5.14]{tfb2} that given
$(y,\omega)\in X\times\widehat{G}$, the closed set in $\widehat{G}$
associated to $\ker\Ind_{G_{y}}^G \chi|_{G_y}$ is $\chi G_y^\perp$.
At this point we can use \eqref{eq:147} and
Lemma \ref{lem:37} to pass to a subnet,
relabel, and find $\sigma_i\in \omega_i G_{x_i}^\perp$ such that
$\sigma_i\rightarrow \omega$.  It follows that
$(\sigma_i,x_i)\rightarrow (\omega,x)$.  Since we clearly have
$(\sigma_i,x_i)\equiv (\omega_i,x_i)$ for all $i$ this suffices to
show that $\phi$ is open.  

Using $\phi$ we can transport the action of $G$ on $\widehat{S}$ to
$X\times \widehat{G}/\equiv$.  It is easy to show, using the fact
that $G$ is abelian, that this action is given by 
\begin{equation}
\label{eq:148}
s\cdot [x,\omega] = [s\cdot x, \omega].
\end{equation}
We would like to show that $((X\times \widehat{G}/\!\equiv)/G)^{T_0} =
X\times \widehat{G}/\sim$.  In other words, we want to see that the
equivalence relation induced by the iterated quotient is exactly
$\sim$.  Observe that, almost by definition, $(x,\omega)$ and $(y,\chi)$ will
be identified in ${((X\times \widehat{G}/\equiv)/G)^{T_0}}$ if and only
if $\overline{G\cdot[x,\omega]} = \overline{G\cdot [y,\chi]}$.  
Suppose $(x,\omega)\sim (y,\chi)$.  Since $\overline{G\cdot y} =
\overline{G\cdot x}$ we must have $s_i$ such that $s_i\cdot x
\rightarrow y$.  It follows from \eqref{eq:148} that
$[y,\omega]\in \overline{G\cdot [x,\omega]}$.  However, $[y,\omega] =
[y,\chi]$ so that we must have $\overline{G\cdot [x,\omega]} =
\overline{G\cdot [y,\chi]}$.  Next, consider the opposite
direction.  If $\overline{G\cdot[x,\omega]}=\overline{G\cdot[y,\chi]}$
then there exists a sequence $s_i\in G$ such that 
\[
[s_i\cdot x, \omega]\rightarrow [y, \chi].
\]
It is straightforward to show, using the fact that $\rho$ is open,
that the quotient map $X\times \widehat{G}\rightarrow X\times
\widehat{G}/\equiv$ is open.  Use this fact 
to pass to a subnet, relabel, and find
$\sigma_i\in \widehat{G}$ such that $(s_i\cdot x, \omega)\equiv
(s_i\cdot x, \sigma_i)$ for all $i$ and 
\[
(s_i\cdot x, \sigma_i)\rightarrow (y,\chi).
\]
First, observe that this implies that $\overline{G\cdot x} =
\overline{G\cdot y}$.  Furthermore, we have
$\sigma_i\overline{\omega}\in G_x^\perp$ for all $i$.  Since this set
is closed and $\sigma_i\rightarrow \chi$ we must also have
$\chi\overline{\omega} \in G_x^\perp$.  Thus $(x,\omega)\sim
(y,\chi)$.  Hence, the equivalence relation induced by the iterated
crossed product is exactly $\sim$ and we have 
\[
(\widehat{S}/G)^{T_0} \cong ((X\times\widehat{G}/\equiv)/G)^{T_0} =
X\times \widehat{G}/\sim.  \qedhere
\]
\end{proof}

\subsection{Examples}

The first example we present is of a non-abelian group action with abelian,
continuously varying stabilizers. Unfortunately, it is also a
transitive action.  As such it really falls under the purview of
\cite[Section 3]{groupoidequiv}.  However, it is elegant enough that we
will include it here for reference.  

\begin{example}
\index{groupoid!transitive}
Let $G = \R\setminus \{0\}\times \R$ be the $ax+b$ group and let $G$ act
on $X=\R$ by evaluation.  That is to say, we have $(ax+b)\cdot r =
ar+b$.  This is a classic continuous group action.  Now, the $ax+b$ group is
non-abelian, but it turns out that the stabilizers are
abelian.  Given $r\in X$ we have 
\[
S_r= \{ax+b\in G: ar+b = r\} = \{ ax+r(1-a):a\in\R\setminus\{0\}\},
\]
and it is straightforward 
to show that $S_r$ is isomorphic to the multiplicative group of $\R$.  
Suppose $r_i\rightarrow r$ in $X$ and
$a\in \R\setminus \{0\}$.  Then $ax+r_i(1-a)\rightarrow ax+r(1-a)$ in
$G$ and this suffices to show that the stabilizers vary continuously.  
Next, the orbit space $X/G$ is a single point since $G$ acts
transitively and is trivially $T_0$.  Thus we can use Corollary
\ref{cor:20} to identify $\Prim C^*(G,X)$ with
$(\widehat{S}/ G)^{T_0}$.  Fix $r\in X$ and consider the continuous map
$\phi:\widehat{S}_r\rightarrow \widehat{S}/G$ defined by $\phi(\omega)
= G\cdot \omega$.  If $\omega,\chi\in \widehat{S}_r$ and $(ax+b)\cdot \omega =
\chi$ then we must have $(ax+b)\in S_r$.  However, $S_r$ is abelian
and therefore $\chi =(ax+b)\cdot \omega = \omega$.  Thus $\phi$ is
injective.  Next suppose $\omega \in \widehat{S}_l$ for some $l\in X$.  Then
there exists $(ax+b)$ such that $(ax+b)\cdot l = r$.  Thus 
$(ax+b)\cdot \omega \in \widehat{S}_r$ and $\phi((ax+b)\cdot \omega) =
G\cdot ((ax+b)\cdot \omega) = G\cdot \omega$.  Hence $\phi$ is
surjective as well. Now suppose $G\cdot \omega_i\rightarrow G\cdot
\omega$.  
We can pass to a subnet, relabel, and choose new
representative so that $\omega_i\rightarrow \omega$.  Since the
transformation groupoid $G\times X$ is
second countable and transitive we can cite \cite[Theorem
2.2A,2.2B]{groupoidequiv} to conclude that the restriction of the
source map on $G\times X$ to $(G\times X)^r$ is open.  Since
$\hat{p}(\omega_i)\rightarrow \hat{p}(\omega)$ we can pass to another
subsequence, relabel, and find $(a_ix+b_i,r)\in G\times X$ and
$(ax+b,r)\in G\times X$
such that $(a_ix+b_i,r)\rightarrow (ax+b,r)$, $(a_ix+b_i)\inv\cdot r =
\hat{p}(\omega_i)$ for all $i$ and $(ax+b)\inv\cdot r =
\hat{p}(\omega)$.  Thus we have  
\[
(a_ix+b_i)\cdot \omega_i \rightarrow (ax+b,r)\cdot \omega.
\]
But each $(a_ix+b_i)\cdot \omega_i$ and $(ax+b)\cdot \omega$ is
in $S^r$ and they map to $G\cdot \omega_i$ and $G\cdot\omega$,
respectively.  This suffices to show that $\phi\inv$ is a
continuous map.  The 
upshot is that $\widehat{S}/G$ is homeomorphic to
$\widehat{S}_r$.   Since $\widehat{S}_r$ is already $T_0$, taking the
$T_0$-ization doesn't do anything and we have the following chain of
identifications 
\[
\Prim C^*(G,X) \cong (\widehat{S}/G)^{T_0} \cong \widehat{S}_r \cong
\widehat{\R^\times} \cong \R^\times
\]
where $\R^\times$ is the multiplicative group of $\R$ and we have used
the fact that this group is self dual. 
\end{example}

As we noted, this example is unsatisfactory since $G$ acts
transitively.  What's more, we actually used an important result from
\cite{groupoidequiv} at one point so we would have been better off
using \cite[Theorem 3.1]{groupoidequiv} from the start.  This next
example is much better in that it requires the full power of Corollary
\ref{cor:20}.  

\begin{example}
\label{ex:27}
Let $G = SO(3,\R)$ and $X = \R^3\setminus\{(0,0,0)\}$.  Let $G$ act on
$X$ by rotation.  Once again this is a classic continuous group
action.  It is clear that $G$ is not abelian.  However, it does have
abelian isotropy.  Given a vector $v\in X$ it's easy to see that $S_v$ is
the set of rotations about the line described by $v$.  In particular, 
if we let $v$
be the first vector of an orthogonal basis for $\R^3$ then it is
straightforward to show that $S_v$ is the set of matrices in $SO(3)$
which fix the first coordinate.  This is isomorphic to $SO(2)$
which is itself isomorphic to the circle group and is therefore
abelian.  It is a little bit more complicated to see that the
stabilizers vary continuously.  Suppose $v_i\rightarrow v\in X$ and
that $S$ is a rotation about $v$.  If $S$ is rotation by $\theta$ then
the goal will be to show that the rotations about $v_i$ by $\theta$,
say $S_i$, converge to $S$.  This is intuitively clear.  Now, it takes
some computation, but one can show that the matrix $S^\theta_w$
which rotates $\theta$ degrees around a vector $w=(x,y,z)$ is given by 
\[
\frac{1}{L^2}
\left[\begin{array}{c c c}
x^2 + (y^2+z^2)\cos\theta & xy(1-\cos \theta)-zL\sin\theta &
xz(1-\cos\theta)+yL\sin\theta \\
xy(1-\cos\theta)+zL\sin\theta & y^2 + (x^2+z^2)\cos\theta & yz(1-\cos
\theta) - xL\sin\theta \\
xz(1-\cos\theta)-yL\sin\theta & yz(1-\cos\theta)+xL\sin\theta &
z^2+(x^2+y^2)\cos\theta \end{array}\right]
\]
where $L = \sqrt{x^2+y^2+z^2}$.  Observe that
$S^\theta_w$ varies continuously with respect to $w$ and hence 
$S_i\rightarrow S$.
This shows that the stabilizer bundle is
continuously varying.  Next, consider $X/G$.  With a little thought one can
convince oneself that this space is homeomorphic to the open
half-line.  Thus $X/G$ is obviously $T_0$.  We may now use Corollary
\ref{cor:20} to identify $\Prim C^*(G,X)$ with
$(\widehat{S}/G)^{T_0}$.  Without, going into the details we will cap
this example by examining $(\widehat{S}/G)^{T_0}$.  
We already observed that $S_v$ is isomorphic
to $\T$.  It is not particularly difficult, considering the matrix
formula above, to show that $S$ is
isomorphic to $X\times \T$.  Thus the dual bundle $\widehat{S}$ is
isomorphic to $X\times \Z$.  Next, let us consider the action of $G$
on $X\times \T$.  Suppose $S\in S_v$ is a rotation by $\theta$ around
$v$.  Given $U\in G$ a computation shows that $USU^*$
is just rotation about $Uv$ by $\theta$.  Thus, the action of $G$ on
$X\times \T$ is given by $U\cdot (v,\theta) = (Uv,\theta)$.
Using this fact, it is straightforward to show that the action of $G$
on $X\times \Z$ is given by $U\cdot (v,z) = (Uv,z)$.  In particular,
the quotient $X\times Z/G$ is isomorphic to $(0,\infty)\times \Z$.
Since this space is clearly $T_0$, taking the $T_0$-ization doesn't do
anything and we have 
\[
\Prim C^*(G,X) \cong (0,\infty)\times \Z.
\]
\end{example}

\begin{remark}
It would be nice to apply the results of this section to genuine
groupoid actions.  Unfortunately the field is laking in naturally
defined groupoids with interesting continuously varying isotropy. 
\end{remark}


\section{Groupoid Algebras with Hausdorff Spectrum}
\label{sec:haussdorff}
We finish this chapter, and the thesis, with an in depth examination
of which conditions imply that groupoid $C^*$-algebras have 
Hausdorff spectrum.  At
the risk of spoiling the punchline, it turns out that we don't find
any good conditions.  In particular, we are trying to generalize
\cite{tghs}, which states that, for abelian transformation groups with
$T_0$-orbit space, the
spectrum of $C^*(G,X)$ is Hausdorff if and only if the stabilizers
vary continuously and the orbit space is Hausdorff.  It turns out that
the naive generalization of
\cite{tghs} doesn't work and the situation is more complex.  Of
course,  this is in many ways more interesting than if the
straightforward generalization would have held.  There
is a notion that ``everything which is true for transformation groups
is true for groupoids'' and this provides a situation where such wishful
thinking fails.  

For now, let us drop down a couple of ``levels'' and consider the problem of
when the spectrum of $C^*(G)$ is $T_0$.  Suppose we are working
with a groupoid that has continuously varying abelian stabilizers.  
If $G\unit/G$ is $T_0$ then it follows from Proposition
\ref{prop:103} that $C^*(G)$ is Type I.  Hence its spectrum is
isomorphic to its primitive ideal space and must be $T_0$.
Furthermore, if $G\unit/G$ is $T_0$ it follows from
Corollary \ref{cor:16} that $C^*(G)\sidehat$ is homeomorphic to
$\widehat{S}/G$.  Thus we can make the interesting deduction that
$\widehat{S}/G$ is $T_0$ whenever $G\unit/G$ is $T_0$.  There is
actually a direct proof for this which can be generalized to the $T_1$
case. 

\begin{prop}
\label{prop:106}
Let $G$ be a second countable, locally compact Hausdorff groupoid
with a Haar system and continuously varying abelian stabilizers $S$.
If $G\unit/G$ is $T_0$ (resp. $T_1$) then $S/G$ and $\widehat{S}/G$
are $T_0$ (resp. $T_1$) as well.  
\end{prop}

\begin{proof}
First, recall that $G$ acts on $S$ by conjugation so that $\gamma\cdot
s = \gamma s \gamma\inv$.  
Suppose we are given $[s],[t]\in S/G$ such that $[s]\ne[t]$.  Let
$\tilde{p}$ denote the factorization of $p$ to $S/G$ and set
$[u]=\tilde{p}([s])$ and $[v] = \tilde{p}([t])$.  Suppose $[u]\ne
[v]$.  If $G\unit/G$ is $T_0$ then we can find an open set $O$ in
$G\unit/G$ containing either $[u]$ or $[v]$ and not the other.
Clearly $\tilde{p}\inv(O)$ is an open set containing either $[s]$ or
$[t]$ and not the other.  If
$G\unit/G$ is $T_1$ then we can find open sets $U,V$ such that $[u]\in
U$, $[v]\in V$ and $[u]\not\in V$, $[v]\not\in U$.  Then clearly
$\tilde{p}\inv(U)$ contains $[s]$ and not $[t]$ and $\tilde{p}\inv(U)$
contains $[t]$ and not $[s]$.  Thus, when $[u]\ne [v]$ we can
separate $[s]$ and $[t]$ to the same degree that we can separate $[u]$
and $[v]$.  

Next, suppose $[u] = [v]$.  We can assume without loss of generality
that $s,t\in S_u$.  However, since $[s]\ne [t]$, we must have $s\ne t$.  Let
$q:S\rightarrow S/G$ be the quotient map and recall that it is open.
Fix a neighborhood $U$ of $s$.  If $t\not\in G\cdot U$ then
$[t]\not\in q(U)$ and $q(U)$ separates $[s]$ from $[t]$.
Now suppose that $t\in G\cdot U$ for all neighborhoods $U$ of $s$.
Then for each $U$ there exists $\gamma_U\in G$ and $s_U\in U$ such that
$s_U = \gamma_U\cdot t$.  If we direct $s_U$ by decreasing $U$ then it
is clear that $s_U\rightarrow s$.  This implies that 
\begin{equation}
\label{eq:149}
\gamma_U\cdot u = r(\gamma_U) = p(s_U) \rightarrow u.
\end{equation}
Since $G\unit/G$ is (at least) $T_0$ we can use Theorem
\ref{thm:glimmdich} to conclude that $[\gamma]\mapsto r(\gamma)$ is a
homeomorphism from $G_u/S_u$ to $[u]$.  It follows from \eqref{eq:149}
that $[\gamma_U]\rightarrow [u]$ in $G_u/S_u$.  However, the quotient
map from $G_u$ onto $G_u/S_u$ is open so that we can pass to a subnet,
relabel, and choose $r_U\in S_u$ such that $\gamma_U r_U \rightarrow
u$.  Hence $\gamma_Ur_U\cdot t \rightarrow u\cdot t = t$.  But $S_u$
is abelian so that $r_U\cdot t = t$ for all $U$.  Therefore we also
have $\gamma_Ur_U\cdot t = \gamma_U\cdot t = s_U \rightarrow s$.  But
then $s = t$, which is a contradiction.  It follows that we must have
been able to separate $[s]$ from $[t]$.  This argument is completely
symmetric so that we can also find an open set around $[t]$ which does
not contain $[s]$ (even if $G\unit/G$ is only $T_0$).  It now follows
that $S/G$ is $T_0$ (resp. $T_1$) if $G\unit/G$ is $T_0$
(resp. $T_1$).  The argument for $\widehat{S}$ is exactly the same
and we end up with the same result.  
\end{proof}

This gives us the following provocative corollary. 

\begin{corr}
\label{cor:21}
Let $G$ be a second countable, locally compact Hausdorff groupoid
with a Haar system and abelian continuously varying stabilizers $S$.
Then $C^*(G)\sidehat$ is $T_0$ if $G\unit/G$ is $T_0$ and
$C^*(G)\sidehat$ is $T_1$ if $G\unit/G$ is $T_1$. 
\end{corr}

\begin{proof}
As long as $G\unit/G$ is at least $T_0$ we can use Corollary
\ref{cor:16} to identify $C^*(G)\sidehat$ with $\widehat{S}/G$.
The result now follows from Proposition \ref{prop:106}.
\end{proof}

If we could extend this result to the $T_2$ (i.e. Hausdorff) case then
we would have made great progress in identifying when the spectrum of
$C^*(G)$ is Hausdorff.  In particular we would have generalized one
direction of the main result in \cite{tghs} which states that,
for abelian transformation groups, the transformation group algebra
has Hausdorff spectrum if and only if $X/G$ is Hausdorff and the
stabilizers vary continuously.  Interestingly enough, it turns out
that Proposition \ref{prop:106}, and hence Corollary \ref{cor:21},
doesn't extend to the Hausdorff case.  In order to build our
counterexample we will have to use Green's famous example from
\cite{tgasos}. 

\begin{example}
\label{ex:28}
The space $X\subset\R^3$ will consist of countably many orbits, with
the points $x_0=(0,0,0)$ and $x_n = (2^{-2n},0,0)$ for $n\in\N$
as a family of representatives.  The action of $\R$ on $X$ is
described by defining maps $\phi_n :\R\rightarrow X$ such that
$\phi_n(s) = s\cdot x_n$.  In particular we let
\[
\phi_0(s) = (0,s,0) 
\]
and for $n\geq 1$
\[
\phi_n(s) = \begin{cases} (2^{-2n},s,0) & s \leq n \\
(2^{-2n}-(s-n)2^{-2n-1}, n\cos(\pi(s-n)),n\sin(\pi(s-n))) & n < s <
n+1 \\
(2^{-2n-1}, s-1-2n,0) & s \geq n+1.
\end{cases}
\]
This is a well known example of a continuous, free
group action that is not proper.  It is straightforward to observe that
the orbit space $X/\R$ is homeomorphic to the subset
$\{x_n\}_{n=0}^\infty$ of $\R^3$.  We will also make use of the
restriction of this action to an action of $\Z$ on the subset
\[
Y = \{\phi_n(m) : n\in \N, m\in \Z\}.
\]
In particular, the restriction to an action of $\Z$ on $X$
clearly yields a continuous action, as does a further restriction
to the action of $\Z$ on the $\Z$-invariant subset $Y$. 
\end{example}

Next, we will build an example of a groupoid $G$ with continuously
varying stabilizers such that $G\unit/G$
is Hausdorff and both $S/G$ and $\widehat{S}/G$ are not Hausdorff.  In
particular, since $G\unit/G$ is definitely $T_0$, it will follow that
$C^*(G)\sidehat$ is not Hausdorff.  

\begin{example}
\label{ex:29}
Let $X$ be as in Example \ref{ex:28} and let $G= \R\rtimes_\phi \R$ be
the semidirect product where we define $\phi(r)(s) := s e^r$.  Note that
$\phi(r)$ is clearly a continuous automorphism of $\R$.  We also have
\[
\phi(r+t)(s) = se^{r+t} = s e^r e^t - \phi(r)(\phi(t)(s)).
\]
Thus $\phi$ is a homomorphism into the automorphism group and it
follows that the semidirect product is a well defined, second countable,
locally compact Hausdorff group.  The group operations are given by 
\begin{align*}
(r,s)(u,v) &= (r+\phi(s)(u),s+v) = (r+ue^s, s+v), \\
(r,s)\inv &= (\phi(-s)(-r),-s) = (-re^{-s},-s).
\end{align*}

Next, let the second factor of $G$ act on $X$ as in Example
\ref{ex:28}.  In other words, let $(r,s)\cdot x := s\cdot x$ where
$s\cdot x$ is defined via the $\phi_n$.  It is straightforward to show
that this is a continuous group action.  It follows that the
transformation groupoid $G\ltimes X$ is a second countable, locally
compact Hausdorff groupoid with a Haar system.  Furthermore, the
stabilizer group at $x$ is $S_x = \{(s,0):s\in \R\}$ for all $x\in
X$.  Since the stabilizers are constant, they must vary continuously,
both in $G$ and in $G\ltimes X$.
Furthermore
\[
(s,0)(t,0) = (s + e^0 t, 0) = (s+t,0)
\]
and this clearly implies that $S_x$ is abelian.  Thus $G\ltimes X$ has
continuously varying abelian stabilizers.  Furthermore, the action of
$G\ltimes X$ on $X$ has the same orbits as the action of $G$ on $X$
which in turn has the same orbits as the action of $\R$ on $X$.  In
particular $X/G$ is homeomorphic to $\{x_i\}_{i=0}^\infty$ which is
clearly Hausdorff.  

Let $S$ be the stabilizer subgroupoid of $G\ltimes H$ and $s_n =
\{((e^{-2n-1},0),(2^{-2n},0,0))\}$.  Then $\{s_n\}\subset S$ and clearly
$s_n \rightarrow s = ((0,0),(0,0,0))$.  Consider 
\[
\gamma_n = ((0,2n+1),(2^{-2n-1},0,0))\quad\text{for all $n$.}
\]
One can compute that 
\begin{align*}
s(\gamma_n ) =
(2^{-2n},0,0), \quad\text{and}\quad r(\gamma_n) = (2^{-2n-1},0,0).
\end{align*}
In particular if we
let $\gamma_n$ act on $s_n$ then we obtain 
\begin{align*}
\gamma_n\cdot s_n &= ((0,2n+1)(e^{-2n-1},0)(0,-2n-1),(2^{-2n-1},0,0))
\\
&= ((1,0),(2^{-2n-1},0,0)).
\end{align*}
Therefore we have $\gamma_n \cdot s_n \rightarrow t =
((1,0),(0,0,0))$.  This, of course, implies that $[s_n]$ converges to
both $[s]$ and $[t]$ in $S/G\ltimes X$.  If $\gamma\cdot s
= t$ then we would have $r(\gamma)=s(\gamma)=(0,0,0)$ so
that $\gamma \in S_{(0,0,0)}$.  In particular, $\gamma =
((h,0),(0,0,0))$ for some $h\in \R$.  But if this is the case then it
is easy to compute that $t=\gamma\cdot s = s$.  This is a contradiction
so that we must have $[s]\ne [t]$.  Hence $[s_n]$ has two distinct
limits in $S/G\ltimes X$.  

Next, we show that $\widehat{S}/G\ltimes X$ is not Hausdorff.  First,
however, we have to compute the dual.  Since the stabilizers are
constant in $G$ it follows that $S$ must be a trivial group
bundle.  In particular, $S$ is isomorphic to $\R\times X$ via the map
$((s,0),x)\mapsto (s,x)$.  Thus we can identify $\widehat{S}$ with
$\widehat{\R}\times X\cong \R\times X$ where we recall that
$\hat{s}\in\R$ acts as a character on $\R$ via $\hat{s}(t) =
e^{ist}$.  There is an action of $G\ltimes X$ on $\widehat{S}$
given by $\gamma\cdot \omega(s) = \omega(\gamma\inv\cdot s)$.  We can
calculate that in our example 
\begin{align*}
((s,t),x)\cdot (\hat{r},(s,t)\inv\cdot x)(q,x) &=
(\hat{r},-t\cdot  x)(((s,t),x)\inv\cdot (q,x)) \\
&= (\hat{r},-t\cdot x)(qe^{-t},-t\cdot x) \\
&= e^{irqe^{-t}} = (\widehat{re^{-t}},x)(q,x).
\end{align*}
It follows that 
\begin{equation}
\label{eq:150}
((s,t),x)\cdot (\hat{r},-t\cdot x) = (\widehat{re^{-t}},x)
\end{equation}
First we observe the following from \eqref{eq:150}.
Suppose $((s,t),x)\cdot (\hat{r},x) = (\hat{q},x)$.  Then we must have
$(s,t)\in S_x$ and therefore $t = 0$.  But then 
\[
(\hat{q},x) = ((s,t),x)\cdot (\hat{r},x) = (\widehat{re^0},x ).
\]
In particular, the action of $G\ltimes X$ is trivial when restricted
to a fixed fibre in $\widehat{S}$.  Moving on, let $\gamma_n\in
G\ltimes X$ be as
above and $\omega_n = (\hat{1},(2^{-2n},0,0))$ for all $n$.  It is
clear that $\omega_n \rightarrow \omega = (\hat{1},(0,0,0))$.
However, it follows from \eqref{eq:150} that 
\[
\gamma_n\cdot \omega_n = ((e^{-2n-1})\sidehat,(2^{-2n-1},0,0)).
\]
Thus $\gamma_n\cdot \omega_n \rightarrow \chi = (\hat{0},(0,0,0))$ and
$[\omega_n]$ converges to both $[\omega]$ and $[\chi]$ in
$\widehat{S}/G\ltimes X$.  Furthermore, since $\omega$ and $\chi$ are
distinct elements of a single fibre, we must have $[\omega]\ne
[\chi]$.  
\end{example}

This example shows that 
Proposition \ref{prop:106} does not extend to the $T_2$ case and
in particular we cannot use it to determine when the spectrum of a
groupoid is Hausdorff or not.  However, Example \ref{ex:29} was
constructed to behave poorly, and there are large classes of groupoids
for which Proposition \ref{prop:106} does extend.  We would like to
find an additional hypothesis which will allow us to make this
extension.  Consider the following

\begin{prop}
\label{prop:107}
Let $G$ be a second countable, 
locally compact Hausdorff groupoid with a Haar system and continuously
varying abelian stabilizers $S$.  Then the action of $G$ on $S$
factors to an action of the orbit groupoid $R$ on $S$ which is strongly continuous when
we give $R$ the quotient topology.  Similarly the the
action of $G$ on $\widehat{S}$ factors to an action of $R$ on
$\widehat{S}$ which is strongly continuous when we give $R$ the
quotient topology. 
\end{prop}

\begin{proof}
Let $\pi:G\rightarrow R$ be the canonical map.  
Define an action of the orbit groupoid 
$R$ on $S$ by factoring the action of $G$ through
$\pi$ and setting 
\begin{equation}
\pi(\gamma)\cdot s := \gamma\cdot s.
\end{equation}
whenever $s(\gamma) = p(s)$.  We need to show that this action is well
defined.  Given $\gamma,\eta\in G$ such that
$\pi(\gamma)=\pi(\eta)$ and $\pi(\gamma)$ and $\pi(\eta)$ act on
$s$, we have $\gamma\inv\eta \in S_{p(s)}$.  In particular, since
$S_{p(s)}$ is abelian, we obtain 
\[
\gamma \cdot s = \gamma\cdot (\gamma\inv\eta \cdot s) = \eta\cdot s.
\]
Hence the action is well defined.  It is straightforward, using the
fact that $\pi$ is a homomorphism, to show that the action respects
the groupoid operations.  Furthermore, the structure map for $S$ is
open by assumption.  We would like to show that this action is
continuous when we give $R$ the quotient topology.  
Suppose $\pi(\gamma_i)\rightarrow \pi(\gamma)$ in $R_Q$
and $s_i\rightarrow s$ such that $p(s_i) = s(\gamma_i)$ for all $i$
and $p(s) = s(\gamma)$.  Then, citing Proposition \ref{prop:11}, we
can pass to a subnet, relabel, and choose new representatives so that
$\gamma_i\rightarrow \gamma$.  However, this implies $\gamma_i\cdot
s_i\rightarrow \gamma\cdot s$ and we are essentially done. 

Next, define an action of $R$ on $\widehat{S}$ via 
\[
(u,v)\cdot \omega(s) = \omega((v,u)\cdot s).
\]
It is straightforward to show that this action respects the groupoid
operations, and since $\widehat{S}$ varies continuously, the structure
map $\hat{p}$ is open.  All we need to do is show that the action is
continuous when we give $R$ the quotient topology.   
Suppose $\omega_i\rightarrow \omega$ and
$(u_i,v_i)\rightarrow (u,v)$ in $R_Q$ such that
$\hat{p}(\omega_i) = u_i$ and $\hat{p}(\omega) = u$.  Next suppose
$s_i\rightarrow s$ such that $p(s_i) = v_i$ and $p(s) = v$.  Then,
using Proposition \ref{prop:33} and the continuity of the action of
$R_Q$ on $S$, we have
\[
\omega_i((u_i,v_i)\cdot s_i)\rightarrow \omega((u,v)\cdot s).
\]
This suffices to show that $(u_i,v_i)\cdot \omega_i\rightarrow
(u,v)\cdot \omega$.  
\end{proof}

We will see in Example \ref{ex:30} that Proposition \ref{prop:107}
doesn't hold if we use the product topology instead of the quotient
topology.  It turns out that this fact is an obstruction to
generalizing Proposition \ref{prop:106} to the Hausdorff case.  Recall
that we use $R_P$ to denote the orbit groupoid equipped with the
restriction of the product topology, while $R_Q$ denotes the orbit
groupoid with the quotient topology. 

\begin{prop}
\label{prop:108}
Let $G$ be a second countable, locally compact Hausdorff groupoid
with a Haar system and continuously varying abelian stabilizers.
Furthermore, suppose that the action of $R_P$ on $S$ is continuous.
Then $S/G$ is Hausdorff if $G\unit/G$ is Hausdorff.  
Similarly, if the action of
$R_P$ on $\widehat{S}$ is continuous and $G\unit/G$ is Hausdorff then
$\widehat{S}/G$ is Hausdorff. 
\end{prop}

\begin{proof}
Suppose $[s_i]\rightarrow [s]$ and $[s_i]\rightarrow [t]$ in
$S/G$. Let $\tilde{p}$ be the factorization of $p$ to $S/G$ and set
$[u] = \tilde{p}([s])$ and $[v]=\tilde{p}([t])$.  Suppose $[u]\ne
[v]$.  Then, because $G\unit/G$ is Hausdorff, we can find disjoint
open sets $U$ and $V$ which separate $[u]$ and $[v]$.  Hence
$\tilde{p}\inv(U)$ and $\tilde{p}\inv(V)$ are disjoint open sets which
separate $[s]$ and $[t]$.  However, $[s_i]$ must eventually be in both
of these sets, which is a contradiction. It follows that $[u] = [v]$.
We may as well assume that $s,t\in S_u$. Since the quotient map
$S\rightarrow S/G$ is open, we can pass to a subsequence, twice,
relabel, choose new representatives, and then find $\gamma_i\in G$
such that $s_i\rightarrow s$ and $\gamma_i\cdot s_i \rightarrow t$.
Let $\pi(\gamma_i) = (u_i,v_i)$.  Then $u_i = p(\gamma_i\cdot
s_i)\rightarrow u$ and $v_i = p(s_i) \rightarrow u$.  Hence
$(u_i,v_i)\rightarrow (u,u)$ in $R_P$.  Since the action of $R_P$ is
continuous, we have $(u_i,v_i)\cdot s_i\rightarrow (u,u)\cdot s = s$.
However, the action of $R$ is just the factorization of the action of
$G$ so that $(u_i,v_i)\cdot s_i = \gamma_i\cdot s_i \rightarrow t$.
Hence $t = s$ and $S/G$ is Hausdorff.  The corresponding proof for
$\widehat{S}$ is exactly the same. 
\end{proof}

We can combine this fact with our identification of the
spectrum to obtain the following corollary. 

\begin{corr}
\label{cor:22}
Let $G$ be a second countable locally compact, Hausdorff groupoid
with a Haar system and abelian continuously varying stabilizers $S$.  If
$G\unit/G$ is Hausdorff and the action of $R_P$ on $\widehat{S}$ is
continuous then $C^*(G)$ has Hausdorff spectrum. 
\end{corr}

\begin{proof}
Since $G\unit/G$ is definitely $T_0$ we can apply Corollary
\ref{cor:16}.  The result now follows from Proposition \ref{prop:108}.
\end{proof}

Of course, this isn't very useful unless we can prove
that there are interesting groupoids for which the action of $R_P$ on
$\widehat{S}$ is continuous.  We start our search by finding a number
of equivalent conditions. 

\begin{prop}
\label{prop:109}
Let $G$ be a second countable, locally compact Hausdorff groupoid
with a Haar system and continuously varying abelian stabilizers $S$.
Then the following are equivalent. 
\begin{enumerate}
\item The action of $R_P$ on $S$ is continuous. 
\item The action of $R_P$ on $\widehat{S}$ is continuous. 
\item If $\{s_i\}\subset S$, $s\in S$ and $\{\gamma_i\}\subset G$ such that
  $s_i\rightarrow s$, $s(\gamma_i) = p(s_i)$ and 
  $r(\gamma_i)\rightarrow p(s)$ then $\gamma_i\cdot s_i\rightarrow s$.
\item If $\{\omega_i\}\subset \widehat{S}$, $\omega\in \widehat{S}$
  and $\{\gamma_i\}\subset G$ such that $\omega_i\rightarrow \omega$,
  $s(\gamma_i) = \hat{p}(\omega_i)$ and $r(\gamma_i)\rightarrow
  \hat{p}(\omega)$ then $\gamma_i\cdot\omega_i\rightarrow \omega$. 
\item The map $S\rightarrow G\unit*S/G = \{(u,[s]) : [p(s)]=[u]\}$
  given by $s\mapsto (p(s),[s])$ is a homeomorphism. 
\item The map $\widehat{S}\rightarrow G\unit*\widehat{S}/G =
  \{(u,[\omega]) : [\hat{p}(\omega)] = [u]\}$ given by $\omega\mapsto
  (\hat{p}(\omega),[\omega])$ is a homeomorphism. 
\end{enumerate}
\end{prop}

From a certain point of view, what the last two conditions in
Proposition \ref{prop:109} are saying is that the topology on $S$ and
$\widehat{S}$ is somehow ``constant'' over $G$ orbits. 

\begin{proof}
We start by proving that (a),(c) and (e) are equivalent.  First we
show (a) implies (e).  Given (a) let $\phi:S\rightarrow G\unit*S/G$ be
given by $\phi(p(s),[s])$.  It is clear that $\phi$ is continuous.
Furthermore, if $\phi(s) = \phi(t)$ then $p(s) = p(t)$ and $[s]=[t]$.
But the action of $G$ is trivial when restricted to a single fibre of
$S$ so that $s = t$.  Next, if we have $(u,[s])\in G\unit*S/G$ then
$[u] = [p(s)]$.  In particular, there exists $\gamma\in G$ such that
$r(\gamma) = u$ and $s(\gamma) = p(s)$.  Then $\phi(\gamma\cdot s) =
(u,[\gamma\cdot s]) = (u,[s])$.  Thus $\phi$ is onto.  It
is easy to see that in general $\phi\inv(u,[s]) = (u,p(s))\cdot s$.
Next, suppose $(u_i,[s_i])\rightarrow (u,[s])$.  Pass to a subnet,
relabel, choose new representatives, and assume $s_i\rightarrow s$.
Since $(u_i,p(s_i))\rightarrow (u,p(s))$ in $R_P$, we have
$(u_i,p(s_i))\cdot s_i \rightarrow (u,p(s))\cdot s$.  This suffices to
show that the inverse map is continuous.  
Next, we show that (e) implies (c).  Suppose $s_i\rightarrow s$,
$s(\gamma_i) = p(s_i)$ and $r(\gamma_i)\rightarrow p(s)$.  Then we
must have 
\[
(r(\gamma_i),[s_i])\rightarrow (p(s),[s])
\]
in $G\unit*S/G$.  Since $\phi$ is a homeomorphism,
\[
\phi\inv(r(\gamma_i),[s_i]) = (r(\gamma_i),p(s_i))\cdot s_i
\rightarrow \phi\inv(p(s),[s]) = s.
\]
However, since the action of $R$ is the factorization of the action of
$G$, we have $(r(\gamma_i),p(s_i))\cdot s_i = \gamma_i\cdot s_i$ and
the result follows. 

Finally, let us show that (c) implies (a).  Suppose
$(u_i,v_i)\rightarrow (u,v)$ in $R_P$ and $s_i\rightarrow s$ such that
$p(s_i) = v_i$ and $p(s) = v$.  Fix $\gamma_i$ and $\gamma$ in $G$
such that $\pi(\gamma_i) = (u_i,v_i)$ for all $i$ and $\pi(\gamma) =
(u,v)$.  Since the range map on $G$ is open, we
can pass to a subnet, relabel, and find $\eta_i$ such that 
$r(\eta_i)= u_i$ and $\eta_i\rightarrow \gamma_i$.  But then
$s(\eta_i\inv\gamma_i) = p(s_i)$ for all $i$ and 
\[
r(\eta_i\inv\gamma_i) = s(\eta_i) \rightarrow s(\gamma) = p(s). 
\]
It follows from part (c) that $\eta_i\inv\gamma_i\cdot s\rightarrow
s$.  Since $\eta_i\rightarrow \gamma$ we must have $\gamma_i\cdot s
\rightarrow \gamma\cdot s$ and that the action of $R_P$ is
continuous.  

The proof that (b), (d) and (f) are equivalent is almost exactly
the same and we will not reproduce it here.  We will end by showing
that (c) and (d) are equivalent.  Suppose (c) holds and that
$\omega_i$, $\omega$ and $\gamma$ are as in the statement of (d).
Then given $s_i\rightarrow s$ such that $p(s_i) = r(\gamma_i)$ for all
$i$ and $p(s) = \hat{p}(\omega)$ we have 
\[
\gamma_i\cdot \omega_i(s_i) = \omega_i(\gamma_i\inv\cdot s_i).
\]
Now $s(\gamma_i\inv) = p(s_i)$ and $r(\gamma_i\inv) =
\hat{p}(\omega_i) \rightarrow \hat{p}(\omega)=p(s)$ so that we may apply
part (c) to conclude that $\gamma_i\inv\cdot s\rightarrow s$.  It follows that
$\gamma_i \cdot \omega_i(s_i)\rightarrow \omega(s)$.  Hence
$\gamma_i\cdot \omega_i\rightarrow \omega$.  Thus (c) implies (d).
Now suppose (d) holds.  We can replace $S$ by $\widehat{S}$ in the
above argument to conclude that given $s_i$, $s$ and $\gamma$ as in (c)
then $\gamma_i\cdot \hat{s}_i \rightarrow \hat{s}$ in
$\doubledual{S}$.  
However, $s\mapsto \hat{s}$ is an isomorphism by Theorem
\ref{thm:duality} and it is easy to see that $\eta \cdot \hat{t} =
(\eta\cdot t)\sidehat$ so that (c) follows. 
\end{proof}

Of course, this proposition allows us to relax the conditions of
Proposition \ref{prop:108} in the obvious way.  Furthermore, now that
we have all of these equivalent conditions it is easy for us to see
that there are some fairly wide classes of groupoids for which the
action of $R_P$ is continuous.  

\begin{definition}
Suppose $G$ is a locally compact Hausdorff groupoid which acts on a
locally compact Hausdorff space $X$.  A set $A$ in $X$ is {\em
  wandering} if 
\[
\{\gamma\in G : \gamma\cdot A \cap A \ne \emptyset\}
\]
is compact.  A locally compact Hausdorff groupoid $G$ is {\em Cartan} if
every point in the unit space of $G$ has a neighborhood which is
wandering. 
\end{definition}

\begin{prop}
\label{prop:110}
Let $G$ be a second countable, locally compact Hausdorff groupoid
with a Haar system and continuously varying abelian stabilizers  $S$.
Then the action of $R_P$ on $S$ is continuous if any of the following
are true. 
\begin{enumerate}
\item $G$ is proper. 
\item $G$ is transitive. 
\item $G$ is Cartan. 
\item $G$ is principal. 
\item $G$ is the transformation groupoid of an abelian group action. 
\end{enumerate}
\end{prop}

\begin{proof}
We know that the topology on $R_Q$ is finer than the topology on
$R_P$.  Suppose $\pi:G\rightarrow R_P$ is open.  If $O\subset
R_Q$ is open then $\pi(\pi\inv(O))=O$ is open in $R_P$.  
Hence $R_P = R_Q$.  Since the action of $R_Q$ is always
continuous by Proposition \ref{prop:107}, it clearly suffices to show
that $\pi$ is open.  If $G$ is transitive then this follows from
\cite[Theorem 2.2B]{groupoidequiv}.  Suppose $G$ is proper and that
we have $\gamma_i\in G$ and $\gamma\in G$ such that 
$\pi(\gamma_i)\rightarrow \pi(\gamma)$.  Then in particular
$s(\gamma_i)\rightarrow s(\gamma)$ and $\gamma_i\cdot
s(\gamma_i)\rightarrow \gamma\cdot r(\gamma)$.  Since the action of
$G$ on its unit space is proper, by definition, we can then pass to a
subnet, relabel, and find $\eta\in G$ so that $\gamma_i\rightarrow
\eta$.  It follows
quickly that $\pi(\gamma) = \pi(\eta)$ and therefore $\gamma\eta\inv \in
S$.  Using the fact that $p$ is open we can pass to another subnet,
relabel, and find $s_i\in S$ such that $s_i = r(\gamma_i)$ for all $i$
and $s_i\rightarrow \gamma\eta\inv$. But then $\pi(s_i\gamma_i) =
\pi(\gamma_i)$ for all $i$ and $s_i\gamma_i\rightarrow \gamma$.  It
follows that $\pi$ is open.  

Moving on, suppose $G$ is principal.  Then the action of $R_P$ on $S =
G\unit$ is trivially continuous.  Next, suppose $G$ is Cartan and we
have $\gamma_i\in G$, $s_i\in S$ and $s\in S$ such that
$s_i\rightarrow s$, $s(\gamma_i) = p(s_i)$ and $r(\gamma_i)\rightarrow
p(s)$.  Let $W$ be a wandering neighborhood of $p(s)$.  Then
$s(\gamma_i)\rightarrow p(s)$ so that $s(\gamma_i)$ is eventually in
$W$.  However $r(\gamma_i)\rightarrow p(s)$ as well so that eventually
$r(\gamma_i) \in W$.  But we can then pass to a subnet, relabel, and assume
that $s(\gamma_i),r(\gamma_i)\in W$ for all $i$.  This implies
that 
\[
\gamma_i \in \{\gamma\in G : \gamma\cdot W \cap W \ne \emptyset\}
\]
for all $i$.  Since this set is compact we may pass to a subnet,
relabel, and find $\gamma\in G$ such that $\gamma_i\rightarrow
\gamma$. It follows quickly that 
$s(\gamma) = r(\gamma) = p(s)$ so that $\gamma\in S_{p(s)}$ and
therefore 
\[
\gamma_i\cdot s_i \rightarrow \gamma\cdot s = s.
\]

Finally, suppose $G = H\ltimes X$ where $H$ is an abelian group which
acts continuously on $X$.  Suppose $(s_i,x_i)\in S$, $(s,x)\in S$ and
$(t_i,y_i)\in G$ such that $s(t_i,y_i) = t_i\inv\cdot y_i = x_i$ for
all $i$ and $y_i\rightarrow x$.  Observe that 
\[
(t_i,y_i)\cdot (s_i,x_i) = (t_is_it_i\inv, y_i) = (s_i,y_i).
\]
Since $(s_i,y_i)\rightarrow (s,x)$ we are done.  
\end{proof}

At this point we would like to show that given a second countable
locally compact Hausdorff groupoid $G$ with  a Haar system and abelian,
continuously varying stabilizers $S$ such that $G\unit/G$ and
$\widehat{S}/G$ are Hausdorff then the action of $R_P$ on $S$ is
continuous.  This would prove that if $G$ has
abelian, continuously varying stabilizers and $G\unit/G$ is Hausdorff
then $C^*(G)$ has Hausdorff spectrum if and only if the action of
$R_P$ on $S$ is continuous.  Unfortunately, as we demonstrate in the
following examples, this is not true.  

\begin{example}
\label{ex:30}
Let $Y$ be as in Example \ref{ex:28} and recall that $\Z$ acts on $Y$
via the $\phi_n$ also defined there.  Let $G = \Q_D\rtimes_\phi\Z$ be
the semidirect product where $\Q_D$ denotes the rationals equipped
with the discrete topology.  Furthermore, we define
\begin{equation}
\phi(n)(r) = r2^n
\end{equation}
for all $n\in \Z$ and $r\in \Q$.  It is clear that $\phi(r)$ is an
automorphism of $\Q_D$.  Furthermore, it is easy to show that $\phi$
is a homomorphism from $\Z$ into the automorphism group of $\Q_D$.
Thus $G$ is a locally compact Hausdorff group which is second
countable because it is, in fact, countable.  Let the second factor of
$G$ act on $Y$ as in Example \ref{ex:28}.  In other words, let
\[
(r,n)\cdot y = n\cdot y
\]
where $n\cdot y$ is defined via the $\phi_n$.  As in Example
\ref{ex:29} this gives us a continuous group action of $G$ on $Y$.
Thus the transformation group $G\ltimes Y$ is a second countable,
locally compact Hausdorff groupoid with a Haar system.  Again as in
Example \ref{ex:29}, we can see that the stabilizer subgroups are
given by $S_x = \{(r,0):r\in\Q\}$ and that they must vary continuously since they
are constant.  Furthermore, it is easy to see that they are abelian so
that $G\ltimes Y$ has continuously varying abelian stabilizers.  It
will be important for us to observe that $S$ is isomorphic to $\Q_D
\times X$ via the isomorphism 
\[
((q,0),x)\mapsto (q,x).
\]
In fact, we will often just drop the extra zero and confuse
stabilizers of $G$ with stabilizers of $G\ltimes Y$.  Finally,
$\{x_n\}_{n=0}^\infty$ forms a set of representatives for the orbit
space and it is not difficult to show that $Y/G$ is actually
homeomorphic to $\{x_n\}_{n=0}^\infty$ and is therefore Hausdorff.
Next, given
$((q,n),y)\in G\ltimes Y$ and $(r,x)\in S$ we have 
\begin{equation}
\label{eq:151}
((q,n),y)\cdot (r,x) = (r2^n,y).
\end{equation}

We would like to show that $S/G\ltimes Y$ is Hausdorff.  Suppose
$[s_i]\rightarrow [s]$ and $[s_i]\rightarrow [t]$.  in $S/G\ltimes
Y$.  Since $X/G$ is Hausdorff we can perform the usual trick to see that we
must have $\tilde{p}([s]) = \tilde{p}([t]) = [u]$.  In fact, we may as
well assume that $p(s) = p(t) = u$.  In this case that means $s =
(r,u)$ and $t = (q,u)$ for $r,q\in\Q$.  Now, we can pass to subnets
and lift, twice, choose new representatives, and find $\gamma_i\in
G\ltimes Y$ so that $s_i\rightarrow s$ and $\gamma_i\cdot s_i
\rightarrow t$.  Suppose $s_i = (r_i,x_i)$ and $\gamma_i =
((p_i,n_i),y_i)$.  Then it follows from \eqref{eq:151} that
$\gamma_i\cdot s_i = (r_i 2^{n_i},y_i)$.  Hence $r_i\rightarrow r$ and
$r_i2^{n_i}\rightarrow q$.  However, we gave $\Q_D$ the discrete
topology so that, eventually,
\[
q = 2^{n_i}r_i = 2^{n_i}r.
\]
Now, if $r = 0$ then we have $q=0$ so that $s=t$.  If $r\ne 0$ we know
that eventually $n_i = n =\log_2(q/r)$.  We may as well pass to a
subnet and assume this is always true.  But then $n_i \cdot x_i
\rightarrow n\cdot x$.  However, we also have $n_i\cdot x_i =
\gamma_i\cdot x_i = y_i \rightarrow x$.  Thus $n\cdot x = x$.  But the
action of $\Z$ is free so that we must have $n=0$.  Thus $\log_2(q/r)
= 0$ and $q=r$.  It follows that $s=t$ and that $S/G$ is Hausdorff.  

We have shown that $G\ltimes Y$ is a second countable, locally compact
Hausdorff groupoid with a Haar system and continuously varying abelian
stabilizers.  Furthermore, both $Y/G$ and $S/G\ltimes Y$ are
Hausdorff.  However, we will show that the action of $R_P$ on $S$ is
not continuous.  Consider $s_n = (1,(2^{-2n},0,0))$ for all $n$.  Then
it is clear that $s_n \rightarrow s = (1,(0,0,0))$ in $S$.  Let
$\gamma_n = ((0,2n+1),(2^{-2n-1},0,0))$.  Then we compute that 
\begin{align*}
s(\gamma_n) = (2^{-2n},0,0),\quad\text{and}\quad r(\gamma_n) =(2^{-2n-1},0,0).
\end{align*}
Thus $s(\gamma_n) = p(s_n)$ and $r(\gamma_n)\rightarrow p(s)$.
However, 
\[
\gamma_n\cdot s_n = (2^{2n+1},(2^{-2n-1},0,0))
\]
and this sequence doesn't converge to anything.  
It follows from Proposition \ref{prop:109} that
the action of $R_P$ on $S$ is not continuous. 
\end{example}

This is not quite the example we are looking for.  We also want to
know that $\widehat{S}/G$ is Hausdorff so that the spectrum of
$C^*(G)$ is as well.  This may very well be true in Example
\ref{ex:30}, but in order to compute it we would need to work with the
dual of the discrete rationals, which is ugly indeed.  Furthermore,
there doesn't seem to be any inherent reason why $\widehat{S}/G$ should
be Hausdorff whenever $S/G$ is.  However, we can be tricky and form a
``dualized'' version of Example \ref{ex:30}. 

\begin{example}
\label{ex:31}
Let $Y,\Z$ and $\phi_n$ be as in Example \ref{ex:28}.  Now let $H =
\widehat{\Q_D}\rtimes_\psi \Z$ be the semidirect product of $\Z$ by
the dual of $\Q_D$.  We define 
\[
\psi(n)(\omega)(q) = \omega(\phi(-n)(q)) = \omega(2^{-n}q)
\]
for $n\in\Z$, $\omega\in \widehat{\Q_D}$ and $q\in \Q$.  Recall that,
since $\Q_D$ is discrete, the topology on $\widehat{\Q_D}$ is just the topology
of pointwise convergence.  If $\omega_i\rightarrow \omega$ in
$\widehat{\Q_D}$  then 
\[
\psi(n)(\omega_i)(q) = \omega_i(2^{-n}q)\rightarrow \omega(2^{-n}q) =
\psi(n)(\omega)(q).
\]
Thus $\psi(n)$ is a continuous function.  Furthermore, 
$\psi(0)$ is the identity and 
\[
\psi(n+m)(\omega)(q) = \omega(2^{-m}2^{-n}q) =
\psi(n)(\psi(m)(\omega))(q).
\]
Thus $\psi(n+m)=\psi(n)\circ\psi(m)$.  It now follows that each
$\psi(n)$ is an automorphism of $\widehat{\Q_D}$ and that $\phi$ is a
homomorphism.  Thus the semidirect product is well defined and $H$ is
a locally compact Hausdorff group.  Furthermore, the topology on
$\widehat{\Q}_D$ is the smallest topology such that all the point
evaluations $\ev_q$ are continuous.  In particular, if we fix a
countable basis $\{U_i\}$ of $\T$ then the collection
$\{\ev_q\inv(U_i)\}$ ranging over all $q\in\Q$ and $i$ forms a countable
sub-basis for the topology.  Hence 
$\widehat{\Q_D}$ is second countable and therefore 
$H$ is second countable.  Just as in Example
\ref{ex:30} let the second factor of $H$ act on $Y$ so that
$(\omega,n)\cdot y= n\cdot y$.  Then, as usual, given $y\in Y$ the
stabilizer subgroup is $T_y = \{(\omega,0)\in H : \omega\in
\widehat{\Q_D}\}$. In particular, $T_y$ abelian and the stabilizers
are continuously varying.  In fact, we can identify the stabilizer
subgroupoid $T$ of $H\ltimes Y$ with $\widehat{\Q_D}\times Y$.  Thus the
dual bundle is, using Pontryagin duality, $\widehat{T} = \Q_D\times Y$.  
It's no accident that $\widehat{T}$ is isomorphic to the bundle $S$
from Example \ref{ex:30}.  We are going to show that the
action of $H\ltimes Y$ on $\widehat{T}$ is nearly the same as the
action of $G\ltimes Y$ on $S$.  Given $(q,x)\in \widehat{T}$ and
$((\omega, n),y)\in H\ltimes Y$ such that $-n\cdot y = x$ we have 
\begin{align*}
((\omega,n),y)\cdot (q,x)(\rho,y) &= (q,x)(((\omega,n),y)\inv\cdot
(\rho,y)) \\
&= (q,x)(\psi(-n)(\omega\inv \rho \omega),x) \\
&= \phi(-n)(\rho)(q) = \rho(2^n q) \\
&= (2^nq,y)(\rho,y)\\
\end{align*}
This implies that 
\begin{equation}
\label{eq:152}
((\omega,n),y)\cdot (q,x) = (2^n q, y).
\end{equation}
Observe that \eqref{eq:151} implies that $[(q,x)] = [(r,y)]$ in
$S/G\rtimes Y$ if and only if $y = n\cdot x$ and $q = 2^n r$ for some $n\in \Z$.
However, \eqref{eq:152} implies the same thing about
$\widehat{T}/H\rtimes Y$.  In particular both spaces have the same
orbits so that the quotient spaces $S/G\ltimes Y$ and
$\widehat{T}/H\ltimes Y$ are identical.  What's more, 
we showed that $S/G\ltimes Y$ was Hausdorff in Example \ref{ex:30}.  

Thus we have demonstrated that $H\rtimes Y$ is a second countable,
locally compact Hausdorff groupoid with a Haar system and continuously
varying abelian stabilizer.  Furthermore, the orbit space $Y/H$ is
still homeomorphic to $\{x_n\}$ and is Hausdorff.  We also showed that
$\widehat{T}/H\ltimes Y$ is Hausdorff.  We will now see that the
action of $R_P$ on $\widehat{T}$ is not continuous.  We will do this
by first observing that $G\rtimes Y$ and $H\rtimes Y$ both have unit
space $Y$, they both act in the same way, and therefore
the orbit groupoid $R$ is the same in each case.  Furthermore, since
$R_P$ inherits its topology from $Y\times Y$, it
follows that $R_P$ is the same for both $G\ltimes Y$ and $H\ltimes
Y$.  It now follows from \eqref{eq:151} and \eqref{eq:152} that $R_P$
acts on $S$ and $\widehat{T}$ in the same manner.  Thus, since $R_P$ does
not act continuously on $S$, it does not act continuously on
$\widehat{T}$.  Furthermore Proposition \ref{prop:109} now
implies that the action of $R_P$ on $T$ must not be continuous either.
\end{example}

At this point we are about done since we have constructed a groupoid
with continuously varying abelian stabilizer such that $G\unit/G$ and
$\widehat{S}/G$ are Hausdorff, but the action of $R_P$ isn't continuous.
In particular, in this case the spectrum of $C^*(G)$ is Hausdorff.
This shows that whatever condition is equivalent to assuming Hausdorff
spectrum is weaker than requiring $R_P$ to act continuously on the
stabilizers.  

\begin{remark}
One of the implications of all of this is that, as a hypothesis for
groupoids, continuously varying abelian stabilizers does not play the
same role as abelian does for groups.  Furthermore, 
it is notable that Examples \ref{ex:30} and \ref{ex:31} are
transformation {\em group} actions.  So, in some sense, the notion of
groupoids being ``generalized transformation groups'' holds true in
that we didn't have to leave the transformation group setting to find
counter examples.  
\end{remark}


\appendix

\backmatter

\cleardoublepage
\bibliographystyle{amsalpha}
\addcontentsline{toc}{chapter}{\bibname}
\bibliography{}

\cleardoublepage
\addcontentsline{toc}{chapter}{Notation and Symbol Index}
\printindex[not]

\newpage
\addcontentsline{toc}{chapter}{Index}
\printindex

\end{document}